\documentclass[mitthesis,11pt]{report}
\usepackage{macros/baththesis}

\usepackage{amsmath,amssymb,amsthm}
\usepackage{amscd}
\usepackage{array}
\usepackage{bm}
\usepackage{caption}
\usepackage[inline]{enumitem}
\usepackage{lscape}
\usepackage{multicol}
\usepackage{multirow}
\usepackage{relsize}
\usepackage{stmaryrd}
\usepackage[explicit]{titlesec}
\usepackage{theoremref}
\usepackage{url}
\usepackage{verbatim}
\usepackage[usenames,dvipsnames,table]{xcolor}

\usepackage{tikz}
\usepackage{pgffor}
\usetikzlibrary{decorations.markings,decorations.pathreplacing}

\newcommand{\parse}[1]{\pgfmathparse{#1}\pgfmathresult}

\newcommand{\DynkinScale}           {1}            
\newcommand{\DynkinFontSize}        {\scriptsize}  
\newcommand{\DynkinRadius}          {1.6mm}        
\newcommand{\DynkinStep}            {5.0mm}        
\newcommand{\DynkinXDotSize}        {0.7mm}        
\newcommand{\DynkinDotsScale}       {0.5}          
\newcommand{\DynkinDotsSize}        {0.4mm}        
\newcommand{\DynkinDotsOffset}      {1.2mm}        
\newcommand{\DynkinDoubleLineScale} {0.45}         
\newcommand{\DynkinTripleLineScale} {0.9}          
\newcommand{\DynkinArrowSize}       {1.4mm}        
\newcommand{\DynkinArrowPosition}   {0.57}         
\newcommand{\DynkinInvolutionOffset}{0.09mm}       

\newcommand{\dynkinwrapper}[2][]{              
  \begin{tikzpicture}[                           
    scale=\DynkinScale,
    baseline=\DynkinStep*#1,
    every node/.style={font=\upshape\DynkinFontSize},
    every label/.style={text height=1ex, text depth=0ex}
  ]
    #2
  \end{tikzpicture}
}
\newcommand{\dynkin}[3][-0.2]{                 
  \dynkinwrapper[#1]{#2}
}

\newcommand{\DynkinBDot}[3][]                  
{
  \pgfmathsetmacro{\PosX}{#2}
  \pgfmathsetmacro{\PosY}{#3}
  \node[
    draw,
    circle,
    fill=black,
    minimum size=\DynkinRadius,
    inner sep=0mm,
    label=#1
  ]
  at (\DynkinStep*\PosX,\DynkinStep*\PosY)
     {};
}
\newcommand{\DynkinWDot}[3][]                  
{
  \pgfmathsetmacro{\PosX}{#2}
  \pgfmathsetmacro{\PosY}{#3}
  \node[
    draw,
    circle,
    fill=white,
    minimum size=\DynkinRadius,
    inner sep=0mm,
    label=#1
  ]
  at (\DynkinStep*\PosX,\DynkinStep*\PosY)
     {};
}
\newcommand{\DynkinXDot}[3][]                  
{
  \pgfmathsetmacro{\PosX}{#2}
  \pgfmathsetmacro{\PosY}{#3}
  \node[
    minimum size=\DynkinRadius,
    inner sep=0mm,
    label=#1
  ]
  at (\DynkinStep*\PosX,\DynkinStep*\PosY)
     {};
  \draw[thick]
    (\PosX*\DynkinStep-\DynkinXDotSize,\PosY*\DynkinStep-\DynkinXDotSize) --
    (\PosX*\DynkinStep+\DynkinXDotSize,\PosY*\DynkinStep+\DynkinXDotSize);
  \draw[thick]
    (\PosX*\DynkinStep-\DynkinXDotSize,\PosY*\DynkinStep+\DynkinXDotSize) --
    (\PosX*\DynkinStep+\DynkinXDotSize,\PosY*\DynkinStep-\DynkinXDotSize);
}

\newcommand{\DynkinLine}[4]                    
{
  \pgfmathsetmacro{\StartX}{#1}
  \pgfmathsetmacro{\StartY}{#2}
  \pgfmathsetmacro{\EndX}{#3}
  \pgfmathsetmacro{\EndY}{#4}
  \draw[thin]
  (\DynkinStep*\StartX,\DynkinStep*\StartY) -- (\DynkinStep*\EndX,\DynkinStep*\EndY);
}
\newcommand{\DynkinDoubleLine}[4]              
{
  \pgfmathsetmacro{\StartX}{#1}
  \pgfmathsetmacro{\StartY}{#2}
  \pgfmathsetmacro{\EndX}{#3}
  \pgfmathsetmacro{\EndY}{#4}
  \draw[
    decoration={
      markings,
      mark=at position \DynkinArrowPosition with {
        \draw[thick] (0:0mm) -- +(+135:\DynkinArrowSize);
        \draw[thick] (0:0mm) -- +(-135:\DynkinArrowSize);
      },
    },
    double distance=\DynkinRadius*\DynkinDoubleLineScale,
    postaction={decorate}
  ]
  (\DynkinStep*\StartX,\DynkinStep*\StartY) -- (\DynkinStep*\EndX,\DynkinStep*\EndY);
}
\newcommand{\DynkinTripleLine}[4]              
{
  \pgfmathsetmacro{\StartX}{#1}
  \pgfmathsetmacro{\StartY}{#2}
  \pgfmathsetmacro{\EndX}{#3}
  \pgfmathsetmacro{\EndY}{#4}
  \draw[
    decoration={
      markings,
      mark=at position \DynkinArrowPosition with {
        \draw[thick] (0:0mm) -- +(+135:\DynkinArrowSize);
        \draw[thick] (0:0mm) -- +(-135:\DynkinArrowSize);
      },
    },
    double distance=\DynkinRadius*\DynkinTripleLineScale,
    postaction={decorate}
  ]
    (\DynkinStep*\StartX,\DynkinStep*\StartY) -- (\DynkinStep*\EndX,\DynkinStep*\EndY);
  \draw[thin]
    (\DynkinStep*\StartX,\DynkinStep*\StartY) -- (\DynkinStep*\EndX,\DynkinStep*\EndY);
}
\newcommand{\DynkinEllipsis}[4]                
{
  \pgfmathsetmacro{\StartX}{#1}
  \pgfmathsetmacro{\StartY}{#2}
  \pgfmathsetmacro{\EndX}{#3}
  \pgfmathsetmacro{\EndY}{#4}
  \foreach \x in {-1,0,1}                      
    \node[                                       
      circle,                                    
      fill=black,
      minimum size=\DynkinDotsSize,
      inner sep=0mm
    ]
    at (\DynkinStep*\StartX+\DynkinStep+\x*\DynkinDotsOffset,\DynkinStep*\StartY)
      {};
}
\newcommand{\DynkinDots}[4]                    
{
  \pgfmathsetmacro{\StartX}{#1}
  \pgfmathsetmacro{\StartY}{#2}
  \pgfmathsetmacro{\EndX}{#3}
  \pgfmathsetmacro{\EndY}{#4}
  \draw[thin]
    (\DynkinStep*\StartX,\DynkinStep*\StartY)
    -- (\DynkinStep*\StartX+\DynkinStep*\DynkinDotsScale,\DynkinStep*\StartY);
  \draw[thin]
    (\DynkinStep*\EndX-\DynkinStep*\DynkinDotsScale,\DynkinStep*\EndY)
    -- (\DynkinStep*\EndX,\DynkinStep*\EndY);
  \DynkinEllipsis{#1}{#2}{#3}{#4}
}
\newcommand{\DynkinInvolution}[4]              
{
  \pgfmathsetmacro{\StartX}{#1}
  \pgfmathsetmacro{\StartY}{#2}
  \pgfmathsetmacro{\EndX}{#3}
  \pgfmathsetmacro{\EndY}{#4}
  \draw[thin,<->]
    (\DynkinStep*\StartX,\DynkinStep*\StartY+\DynkinStep*\DynkinInvolutionOffset)
    -- (\DynkinStep*\EndX,\DynkinStep*\EndY-\DynkinStep*\DynkinInvolutionOffset);
}
\newcommand{\DynkinConnector}[4]{
  \pgfmathsetmacro{\StartX}{#1}
  \pgfmathsetmacro{\StartY}{#2}
  \pgfmathsetmacro{\EndX}{#3}
  \pgfmathsetmacro{\EndY}{#4}
  \draw[thin,->]
    (\DynkinStep*\StartX,\DynkinStep*\StartY) -- (\DynkinStep*\EndX,\DynkinStep*\EndY);
}

\newcommand{\DynkinLabel}[3]                   
{
  \pgfmathsetmacro{\PosX}{#2}
  \pgfmathsetmacro{\PosY}{#3}
  \node[label=#1]
  at (\DynkinStep*\PosX,\DynkinStep*\PosY)
     {};
}
\newcommand{\DynkinWt}[2]{                     
  \ifx  \relax#1\relax ~                         
  \else                \parse{int({#1}[#2])}
  \fi
}
\newcommand{\DynkinUnderbrace}[5]              
{
  \pgfmathsetmacro{\StartX}{#2}
  \pgfmathsetmacro{\StartY}{#3}
  \pgfmathsetmacro{\EndX}{#4}
  \pgfmathsetmacro{\EndY}{#5}
  \draw[
    thick,
    decoration={
      brace,
      mirror
    },
    decorate
  ]
  (\DynkinStep*\StartX-\DynkinStep*0.2,\DynkinStep*\StartY)
    -- (\DynkinStep*\EndX+\DynkinStep*0.2,\DynkinStep*\EndY);
  \path (\DynkinStep*\StartX,\DynkinStep*\StartY-\DynkinStep*0.15)
    -- node[inner sep=0mm,label=below:#1] {}
    (\DynkinStep*\EndX,\DynkinStep*\EndY-\DynkinStep*0.15);
}
\newcommand{\DynkinCases}[5]                   
{                                                
  \pgfmathsetmacro{\StartX}{#2}
  \pgfmathsetmacro{\StartY}{#3}
  \pgfmathsetmacro{\EndX}{#4}
  \pgfmathsetmacro{\EndY}{#5}
  \draw[
    thick,
    decoration={
      brace,
      amplitude=4pt,
      aspect=#1
    },
    decorate
  ]
  (\DynkinStep*\StartX+\DynkinStep*0.2,\DynkinStep*\StartY-\DynkinStep*0.2)
    -- (\DynkinStep*\EndX+\DynkinStep*0.2,\DynkinStep*\EndY+\DynkinStep*0.2) node{};
}
\newcommand{\dynkinnameoffset}{}
\newcommand{\dynkinname}[3][]{            
  \underset{ \raisebox{
    \ifx  \relax#1\relax \dynkinnameoffset em
    \else                #1em
    \fi }{ \text{\normalsize$#3$} } }{ #2 } }

\usepackage{macros/readarray}

\newcommand{\dynkinA}[5][]{ \dynkin{
  \DynkinLine{1 }{0}{#2  }{0};
  \DynkinDots{#2}{0}{#2+2}{0};
  \foreach \x in {1,...,#2} { \DynkinWDot[\DynkinWt{#1}{\x-1}]{\x}{0}; }
  \ifthenelse{#3=0}{
    \DynkinLine{#2+2}{0}{#2+#4+1}{0};
    \foreach \x in {1,...,#4} { \DynkinWDot[\DynkinWt{#1}{#2+\x-1}]{#2+1+\x}{0}; }
  }{
    \DynkinLine{#2   +2}{0}{#2+#3   +1}{0};
    \DynkinDots{#2+#3+1}{0}{#2+#3   +3}{0};
    \DynkinLine{#2+#3+3}{0}{#2+#3+#4+2}{0};
    \foreach \x in {1,...,#3} { \DynkinWDot[\DynkinWt{#1}{#2   +\x-1}]{#2   +\x+1}{0}; }
    \foreach \x in {1,...,#4} { \DynkinWDot[\DynkinWt{#1}{#2+#3+\x-1}]{#2+#3+\x+2}{0}; }
  }
}{#5} }
\newcommand{\dynkinAp}[5][]{ \dynkin{
  \pgfmathsetmacro{\E}{#4-1}
  \DynkinLine{1 }{0}{#2  }{0};
  \DynkinDots{#2}{0}{#2+2}{0};
  \foreach \x in {1,...,#2} { \DynkinWDot[\DynkinWt{#1}{\x-1}]{\x}{0}; }
  \ifthenelse{#3=0}{
    \DynkinLine{#2+2}{0}{#2+#4+1}{0};
    \foreach \x in {1,...,\E} { \DynkinWDot[\DynkinWt{#1}{#2+\x-1}]{#2+1+\x}{0}; }
    \DynkinXDot[\DynkinWt{#1}{#2+#4-1}]{#2+#4+1}{0};
  }{
    \DynkinLine{#2   +2}{0}{#2+#3   +1}{0};
    \DynkinDots{#2+#3+1}{0}{#2+#3   +3}{0};
    \DynkinLine{#2+#3+3}{0}{#2+#3+#4+2}{0};
    \foreach \x in {1,...,#3} { \DynkinWDot[\DynkinWt{#1}{#2   +\x-1}]{#2   +\x+1}{0}; }
    \foreach \x in {1,...,\E} { \DynkinWDot[\DynkinWt{#1}{#2+#3+\x-1}]{#2+#3+\x+2}{0}; }
    \DynkinXDot[\DynkinWt{#1}{#2+#3+#4-1}]{#2+#3+#4+2}{0};
  }
}{#5} }
\newcommand{\dynkinApp}[5][]{ \dynkin{
  \pgfmathsetmacro{\E}{#4-2}
  \DynkinLine{1 }{0}{#2  }{0};
  \DynkinDots{#2}{0}{#2+2}{0};
  \foreach \x in {1,...,#2} { \DynkinWDot[\DynkinWt{#1}{\x-1}]{\x}{0}; }
  \ifthenelse{#3=0}{
    \DynkinLine{#2+2}{0}{#2+#4+1}{0};
    \foreach \x in {1,...,\E} { \DynkinWDot[\DynkinWt{#1}{#2+\x-1}]{#2+\x+1}{0}; }
    \DynkinXDot[\DynkinWt{#1}{#2+#4-2}]{#2+#4  }{0};
    \DynkinWDot[\DynkinWt{#1}{#2+#4-1}]{#2+#4+1}{0};
  }{
    \DynkinLine{#2   +2}{0}{#2+#3   +1}{0};
    \DynkinDots{#2+#3+1}{0}{#2+#3   +3}{0};
    \DynkinLine{#2+#3+3}{0}{#2+#3+#4+2}{0};
    \foreach \x in {1,...,#3} { \DynkinWDot[\DynkinWt{#1}{#2   +\x-1}]{#2   +\x+1}{0}; }
    \foreach \x in {1,...,\E} { \DynkinWDot[\DynkinWt{#1}{#2+#3+\x-1}]{#2+#3+\x+2}{0}; }
    \DynkinXDot[\DynkinWt{#1}{#2+#3+#4-2}]{#2+#3+#4  }{0};
    \DynkinWDot[\DynkinWt{#1}{#2+#3+#4-1}]{#2+#3+#4+1}{0};
  }
}{#5} }

\newcommand{\dynkinApk}[4][]{ \dynkin{
  \pgfmathsetmacro{\E}{#3}
  \DynkinLine{1 }{0}{#2  }{0};
  \DynkinDots{#2}{0}{#2+2}{0};
  \foreach \x in {1,...,#2} { \DynkinWDot[\DynkinWt{#1}{\x-1}]{\x}{0}; }
  \DynkinLine{#2+2}{0}{#2+4}{0};
  \DynkinDots{#2+4}{0}{#2+6}{0};
  \DynkinWDot[\DynkinWt{#1}{#2  }]{#2+2}{0};
  \DynkinXDot[\DynkinWt{#1}{#2+1}]{#2+3}{0};
  \DynkinWDot[\DynkinWt{#1}{#2+2}]{#2+4}{0};
  \DynkinLine{#2+6}{0}{#2+5+#3}{0};
  \foreach \x in {1,...,\E} { \DynkinWDot[\DynkinWt{#1}{#2+2+\x}]{#2+5+\x}{0}; }
}{#4} }

\newcommand{\dynkinASU}[5][]{ \dynkin[-0.22]{
  \DynkinLine{1 }{ 0.75}{#2  }{ 0.75};
  \DynkinLine{1 }{-0.75}{#2  }{-0.75};
  \DynkinDots{#2}{ 0.75}{#2+2}{ 0.75};
  \DynkinDots{#2}{-0.75}{#2+2}{-0.75};
  \foreach \x in {1,...,#2} {
    \DynkinWDot[      \DynkinWt{#1}{         \x-1}]{\x}{ 0.75};
    \DynkinWDot[below:\DynkinWt{#1}{#2+#3+#4+\x-1}]{\x}{-0.75};
  }
  \ifthenelse{#3=0}{
    \DynkinLine{#2   +2}{ 0.75}{#2+#4+1}{ 0.75};
    \DynkinLine{#2   +2}{-0.75}{#2+#4+1}{-0.75};
    \DynkinLine{#2+#4+1}{ 0.75}{#2+#4+2}{ 0   };
    \DynkinLine{#2+#4+1}{-0.75}{#2+#4+2}{ 0   };
    \foreach \x in {1,...,#4} {
      \DynkinWDot[      \DynkinWt{#1}{         #2+\x-1}]{#2+1+\x}{ 0.75};
      \DynkinWDot[below:\DynkinWt{#1}{#2+#3+#4+#2+\x-1}]{#2+1+\x}{-0.75};
    }
    \DynkinWDot[right:\DynkinWt{#1}{#2+#3+#4+#2+#3+#4}]{#2+#4+2}{0};
  }{
    \DynkinLine{#2      +2}{ 0.75}{#2+#3   +1}{ 0.75};
    \DynkinLine{#2      +2}{-0.75}{#2+#3   +1}{-0.75};
    \DynkinDots{#2+#3   +1}{ 0.75}{#2+#3   +3}{ 0.75};
    \DynkinDots{#2+#3   +1}{-0.75}{#2+#3   +3}{-0.75};
    \DynkinLine{#2+#3   +3}{ 0.75}{#2+#3+#4+2}{ 0.75};
    \DynkinLine{#2+#3   +3}{-0.75}{#2+#3+#4+2}{-0.75};
    \DynkinLine{#2+#3+#4+2}{ 0.75}{#2+#3+#4+3}{ 0   };
    \DynkinLine{#2+#3+#4+2}{-0.75}{#2+#3+#4+3}{ 0   };
    \foreach \x in {1,...,#3} {
      \DynkinWDot[      \DynkinWt{#1}{         #2+\x-1}]{#2+1+\x}{ 0.75};
      \DynkinWDot[below:\DynkinWt{#1}{#2+#3+#4+#2+\x-1}]{#2+1+\x}{-0.75};
    }
    \foreach \x in {1,...,#4} {
      \DynkinWDot[      \DynkinWt{#1}{         #2+#3+\x-1}]{#2+#3+2+\x}{ 0.75};
      \DynkinWDot[below:\DynkinWt{#1}{#2+#3+#4+#2+#3+\x-1}]{#2+#3+2+\x}{-0.75};
    }
    \DynkinWDot[right:\DynkinWt{#1}{#2+#3+#4+#2+#3+#4}]{#2+#3+#4+3}{0};
  }
}{#5} }
\newcommand{\dynkinASUp}[5][]{ \dynkin[-0.22]{
  \DynkinLine{1 }{ 0.75}{#2  }{ 0.75};
  \DynkinLine{1 }{-0.75}{#2  }{-0.75};
  \DynkinDots{#2}{ 0.75}{#2+2}{ 0.75};
  \DynkinDots{#2}{-0.75}{#2+2}{-0.75};
  \foreach \x in {1,...,#2} {
    \DynkinWDot[      \DynkinWt{#1}{         \x-1}]{\x}{ 0.75};
    \DynkinWDot[below:\DynkinWt{#1}{#2+#3+#4+\x-1}]{\x}{-0.75};
  }
  \ifthenelse{#3=0}{
    \DynkinLine{#2   +2}{ 0.75}{#2+#4+1}{ 0.75};
    \DynkinLine{#2   +2}{-0.75}{#2+#4+1}{-0.75};
    \DynkinLine{#2+#4+1}{ 0.75}{#2+#4+2}{ 0   };
    \DynkinLine{#2+#4+1}{-0.75}{#2+#4+2}{ 0   };
    \foreach \x in {1,...,#4} {
      \DynkinWDot[      \DynkinWt{#1}{         #2+\x-1}]{#2+1+\x}{ 0.75};
      \DynkinWDot[below:\DynkinWt{#1}{#2+#3+#4+#2+\x-1}]{#2+1+\x}{-0.75};
    }
    \DynkinXDot[right:\DynkinWt{#1}{#2+#3+#4+#2+#3+#4}]{#2+#4+2}{0};
  }{
    \DynkinLine{#2      +2}{ 0.75}{#2+#3   +1}{ 0.75};
    \DynkinLine{#2      +2}{-0.75}{#2+#3   +1}{-0.75};
    \DynkinDots{#2+#3   +1}{ 0.75}{#2+#3   +3}{ 0.75};
    \DynkinDots{#2+#3   +1}{-0.75}{#2+#3   +3}{-0.75};
    \DynkinLine{#2+#3   +3}{ 0.75}{#2+#3+#4+2}{ 0.75};
    \DynkinLine{#2+#3   +3}{-0.75}{#2+#3+#4+2}{-0.75};
    \DynkinLine{#2+#3+#4+2}{ 0.75}{#2+#3+#4+3}{ 0   };
    \DynkinLine{#2+#3+#4+2}{-0.75}{#2+#3+#4+3}{ 0   };
    \foreach \x in {1,...,#3} {
      \DynkinWDot[      \DynkinWt{#1}{         #2+\x-1}]{#2+1+\x}{ 0.75};
      \DynkinWDot[below:\DynkinWt{#1}{#2+#3+#4+#2+\x-1}]{#2+1+\x}{-0.75};
    }
    \foreach \x in {1,...,#4} {
      \DynkinWDot[      \DynkinWt{#1}{         #2+#3+\x-1}]{#2+#3+2+\x}{ 0.75};
      \DynkinWDot[below:\DynkinWt{#1}{#2+#3+#4+#2+#3+\x-1}]{#2+#3+2+\x}{-0.75};
    }
    \DynkinXDot[right:\DynkinWt{#1}{#2+#3+#4+#2+#3+#4}]{#2+#3+#4+3}{0};
  }
}{#5} }

\newcommand{\dynkinAA}[5][]{ \dynkin[-0.22]{
  \DynkinLine{1 }{ 0.75}{#2  }{ 0.75};
  \DynkinLine{1 }{-0.75}{#2  }{-0.75};
  \DynkinDots{#2}{ 0.75}{#2+2}{ 0.75};
  \DynkinDots{#2}{-0.75}{#2+2}{-0.75};
  \foreach \x in {1,...,#2} {
    \DynkinWDot[      \DynkinWt{#1}{         \x-1}]{\x}{ 0.75};
    \DynkinWDot[below:\DynkinWt{#1}{#2+#3+#4+\x-1}]{\x}{-0.75};
  }
  \ifthenelse{#3=0}{
    \DynkinLine{#2+2}{ 0.75}{#2+#4+1}{ 0.75};
    \DynkinLine{#2+2}{-0.75}{#2+#4+1}{-0.75};
    \foreach \x in {1,...,#4} {
      \DynkinWDot[      \DynkinWt{#1}{         #2+\x-1}]{#2+1+\x}{ 0.75};
      \DynkinWDot[below:\DynkinWt{#1}{#2+#3+#4+#2+\x-1}]{#2+1+\x}{-0.75};
    }
  }{
    \DynkinLine{#2   +2}{ 0.75}{#2+#3   +1}{ 0.75};
    \DynkinLine{#2   +2}{-0.75}{#2+#3   +1}{-0.75};
    \DynkinDots{#2+#3+1}{ 0.75}{#2+#3   +3}{ 0.75};
    \DynkinDots{#2+#3+1}{-0.75}{#2+#3   +3}{-0.75};
    \DynkinLine{#2+#3+3}{ 0.75}{#2+#3+#4+2}{ 0.75};
    \DynkinLine{#2+#3+3}{-0.75}{#2+#3+#4+2}{-0.75};
    \foreach \x in {1,...,#3} {
      \DynkinWDot[      \DynkinWt{#1}{         #2+\x-1}]{#2+1+\x}{ 0.75};
      \DynkinWDot[below:\DynkinWt{#1}{#2+#3+#4+#2+\x-1}]{#2+1+\x}{-0.75};
    }
    \foreach \x in {1,...,#4} {
      \DynkinWDot[      \DynkinWt{#1}{         #2+#3+\x-1}]{#2+#3+2+\x}{ 0.75};
      \DynkinWDot[below:\DynkinWt{#1}{#2+#3+#4+#2+#3+\x-1}]{#2+#3+2+\x}{-0.75};
    }
  }
}{#5} }
\newcommand{\dynkinAAk}[8][]{ \dynkin[-0.22]{
  \DynkinLine{1 }{ 0.75}{#2  }{ 0.75};
  \DynkinDots{#2}{ 0.75}{#2+2}{ 0.75};
  \foreach \x in {1,...,#2} { \DynkinWDot[\DynkinWt{#1}{\x-1}]{\x}{ 0.75}; }
  \ifthenelse{#3=0}{
    \DynkinLine{#2+2}{ 0.75}{#2+#4+1}{ 0.75};
    \foreach \x in {1,...,#4} { \DynkinWDot[\DynkinWt{#1}{#2+\x-1}]{#2+1+\x}{ 0.75}; }
  }{
    \DynkinLine{#2   +2}{ 0.75}{#2+#3   +1}{ 0.75};
    \DynkinDots{#2+#3+1}{ 0.75}{#2+#3   +3}{ 0.75};
    \DynkinLine{#2+#3+3}{ 0.75}{#2+#3+#4+2}{ 0.75};
    \foreach \x in {1,...,#3} { \DynkinWDot[\DynkinWt{#1}{#2   +\x-1}]{#2   +\x+1}{ 0.75}; }
    \foreach \x in {1,...,#4} { \DynkinWDot[\DynkinWt{#1}{#2+#3+\x-1}]{#2+#3+\x+2}{ 0.75}; }
  }
  \pgfmathsetmacro{\O}{#2+#3+#4}
  \DynkinLine{   2}{-0.75}{#5+1}{-0.75};
  \DynkinDots{#5+1}{-0.75}{#5+3}{-0.75};
  \foreach \x in {1,...,#5} { \DynkinWDot[below:\DynkinWt{#1}{\O+\x-1}]{\x+1}{-0.75}; }
  \ifthenelse{#6=0}{
    \DynkinLine{#5+3}{-0.75}{#5+#7+2}{-0.75};
    \foreach \x in {1,...,#7} { \DynkinWDot[below:\DynkinWt{#1}{\O+#5+\x-1}]{#5+2+\x}{-0.75}; }
  }{
    \DynkinLine{#5   +3}{-0.75}{#5+#6   +2}{-0.75};
    \DynkinDots{#5+#6+2}{-0.75}{#5+#6   +4}{-0.75};
    \DynkinLine{#5+#6+4}{-0.75}{#5+#6+#7+3}{-0.75};
    \foreach \x in {1,...,#6} {
      \DynkinWDot[below:\DynkinWt{#1}{\O+#5   +\x-1}]{#5   +\x+2}{-0.75}; }
    \foreach \x in {1,...,#7} {
      \DynkinWDot[below:\DynkinWt{#1}{\O+#5+#6+\x-1}]{#5+#6+\x+3}{-0.75}; }
  } 
}{#8} }
\newcommand{\dynkinAAp}[5][]{ \dynkin[-0.22]{
  \pgfmathsetmacro{\E}{#4-1}
  \DynkinLine{1 }{ 0.75}{#2  }{ 0.75};
  \DynkinLine{1 }{-0.75}{#2  }{-0.75};
  \DynkinDots{#2}{ 0.75}{#2+2}{ 0.75};
  \DynkinDots{#2}{-0.75}{#2+2}{-0.75};
  \foreach \x in {1,...,#2} {
    \DynkinWDot[      \DynkinWt{#1}{         \x-1}]{\x}{ 0.75};
    \DynkinWDot[below:\DynkinWt{#1}{#2+#3+#4+\x-1}]{\x}{-0.75};
  }
  \ifthenelse{#3=0}{
    \DynkinLine{#2+2}{ 0.75}{#2+#4+1}{ 0.75};
    \DynkinLine{#2+2}{-0.75}{#2+#4+1}{-0.75};
    \foreach \x in {1,...,\E} {
      \DynkinWDot[      \DynkinWt{#1}{      #2+\x-1}]{#2+1+\x}{ 0.75};
      \DynkinWDot[below:\DynkinWt{#1}{#2+#4+#2+\x-1}]{#2+1+\x}{-0.75};
    }
    \DynkinXDot[      \DynkinWt{#1}{      #2+#4-1}]{#2+#4+1}{ 0.75};
    \DynkinXDot[below:\DynkinWt{#1}{#2+#4+#2+#4-1}]{#2+#4+1}{-0.75};
  }{
    \DynkinLine{#2   +2}{ 0.75}{#2+#3   +1}{ 0.75};
    \DynkinLine{#2   +2}{-0.75}{#2+#3   +1}{-0.75};
    \DynkinDots{#2+#3+1}{ 0.75}{#2+#3   +3}{ 0.75};
    \DynkinDots{#2+#3+1}{-0.75}{#2+#3   +3}{-0.75};
    \DynkinLine{#2+#3+3}{ 0.75}{#2+#3+#4+2}{ 0.5};
    \DynkinLine{#2+#3+3}{-0.75}{#2+#3+#4+2}{-0.5};
    \foreach \x in {1,...,#3} {
      \DynkinWDot[      \DynkinWt{#1}{         #2+\x-1}]{#2+1+\x}{ 0.75};
      \DynkinWDot[below:\DynkinWt{#1}{#2+#3+#4+#2+\x-1}]{#2+1+\x}{-0.75};
    }
    \foreach \x in {1,...,\E} {
      \DynkinWDot[      \DynkinWt{#1}{         #2+#3+\x-1}]{#2+#3+2+\x}{ 0.75};
      \DynkinWDot[below:\DynkinWt{#1}{#2+#3+#4+#2+#3+\x-1}]{#2+#3+2+\x}{-0.75};
    }
    \DynkinXDot[      \DynkinWt{#1}{         #2+#3+#4-1}]{#2+#3+#4+2}{ 0.75};
    \DynkinXDot[below:\DynkinWt{#1}{#2+#3+#4+#2+#3+#4-1}]{#2+#3+#4+2}{-0.75};
  }
}{#5} }
\newcommand{\dynkinAApk}[8][]{ \dynkin[-0.22]{
  \pgfmathsetmacro{\E}{#4-1}
  \DynkinLine{1 }{ 0.75}{#2  }{ 0.75};
  \DynkinDots{#2}{ 0.75}{#2+2}{ 0.75};
  \foreach \x in {1,...,#2} { \DynkinWDot[\DynkinWt{#1}{\x-1}]{\x}{ 0.75}; }
  \ifthenelse{#3=0}{
    \DynkinLine{#2+2}{ 0.75}{#2+#4+1}{ 0.75};
    \foreach \x in {1,...,\E} { \DynkinWDot[\DynkinWt{#1}{#2+\x-1}]{#2+1+\x}{ 0.75}; }
    \DynkinXDot[\DynkinWt{#1}{#2+#4-1}]{#2+#4+1}{ 0.75};
  }{
    \DynkinLine{#2   +2}{ 0.75}{#2+#3   +1}{ 0.75};
    \DynkinDots{#2+#3+1}{ 0.75}{#2+#3   +3}{ 0.75};
    \DynkinLine{#2+#3+3}{ 0.75}{#2+#3+#4+2}{ 0.75};
    \foreach \x in {1,...,#3} { \DynkinWDot[\DynkinWt{#1}{#2   +\x-1}]{#2   +\x+1}{ 0.75}; }
    \foreach \x in {1,...,\E} { \DynkinWDot[\DynkinWt{#1}{#2+#3+\x-1}]{#2+#3+\x+2}{ 0.75}; }
    \DynkinXDot[\DynkinWt{#1}{#2+#3+#4-1}]{#2+#3+#4+2}{ 0.75};
  }
  \pgfmathsetmacro{\E}{#7-1}
  \pgfmathsetmacro{\O}{#2+#3+#4}
  \DynkinLine{   2}{-0.75}{#5+1}{-0.75};
  \DynkinDots{#5+1}{-0.75}{#5+3}{-0.75};
  \foreach \x in {1,...,#5} { \DynkinWDot[\DynkinWt{#1}{\O+\x-1}]{\x+1}{-0.75}; }
  \ifthenelse{#6=0}{
    \DynkinLine{#5+3}{-0.75}{#5+#7+2}{-0.75};
    \foreach \x in {1,...,\E} { \DynkinWDot[\DynkinWt{#1}{\O+#5+\x-1}]{#5+2+\x}{-0.75}; }
    \DynkinXDot[\DynkinWt{#1}{\O+#5+#7-1}]{#5+#7+2}{-0.75};
  }{
    \DynkinLine{#5   +3}{-0.75}{#5+#6   +2}{-0.75};
    \DynkinDots{#5+#6+2}{-0.75}{#5+#6   +4}{-0.75};
    \DynkinLine{#5+#6+4}{-0.75}{#5+#6+#7+3}{-0.75};
    \foreach \x in {1,...,#6} { \DynkinWDot[\DynkinWt{#1}{\O+#5   +\x-1}]{#5   +\x+2}{-0.75}; }
    \foreach \x in {1,...,\E} { \DynkinWDot[\DynkinWt{#1}{\O+#5+#6+\x-1}]{#5+#6+\x+3}{-0.75}; }
    \DynkinXDot[\DynkinWt{#1}{\O+#5+#6+#7-1}]{#5+#6+#7+3}{-0.75};
  } 
}{#8} }

\newcommand{\dynkinB}[5][]{ \dynkin{
  \pgfmathsetmacro{\E}{#4-1}
  \DynkinLine{1 }{0}{#2  }{0};
  \DynkinDots{#2}{0}{#2+2}{0};
  \foreach \x in {1,...,#2} { \DynkinWDot[\DynkinWt{#1}{\x-1}]{\x}{0}; }
  \ifthenelse{#3=0}{
    \DynkinLine{#2+2}{0}{#2+\E+1}{0};
    \DynkinDoubleLine{#2+\E+1}{0}{#2+\E+2}{0};
    \foreach \x in {1,...,#4} { \DynkinWDot[\DynkinWt{#1}{#2+\x-1}]{#2+1+\x}{0}; }
  }{
    \DynkinLine{#2   +2}{0}{#2+#3   +1}{0};
    \DynkinDots{#2+#3+1}{0}{#2+#3   +3}{0};
    \DynkinLine{#2+#3+3}{0}{#2+#3+\E+1}{0};
    \DynkinDoubleLine{#2+#3+\E+1}{0}{#2+#3+\E+2}{0};
    \foreach \x in {1,...,#3} { \DynkinWDot[\DynkinWt{#1}{#2   +\x-1}]{#2+1+\x}{0}; }
    \foreach \x in {1,...,\E} { \DynkinWDot[\DynkinWt{#1}{#2+#3+\x-1}]{#2+#3+2+\x}{0}; }
  }
}{#5} }
\newcommand{\dynkinBp}[5][]{ \dynkin{
  \pgfmathsetmacro{\E}{#4-1}
  \DynkinLine{1 }{0}{#2  }{0};
  \DynkinDots{#2}{0}{#2+2}{0};
  \DynkinXDot[\DynkinWt{#1}{0}]{1}{0};
  \foreach \x in {2,...,#2} { \DynkinWDot[\DynkinWt{#1}{\x-1}]{\x}{0}; }
  \ifthenelse{#3=0}{
    \DynkinLine{#2+2}{0}{#2+\E+1}{0};
    \DynkinDoubleLine{#2+\E+1}{0}{#2+\E+2}{0};
    \foreach \x in {1,...,#4} { \DynkinWDot[\DynkinWt{#1}{#2+\x-1}]{#2+1+\x}{0}; }
  }{
    \DynkinLine{#2   +2}{0}{#2+#3   +1}{0};
    \DynkinDots{#2+#3+1}{0}{#2+#3   +3}{0};
    \DynkinLine{#2+#3+3}{0}{#2+#3+\E+1}{0};
    \DynkinDoubleLine{#2+#3+\E+1}{0}{#2+#3+\E+2}{0};
    \foreach \x in {1,...,#3} { \DynkinWDot[\DynkinWt{#1}{#2   +\x-1}]{#2+1+\x}{0}; }
    \foreach \x in {1,...,\E} { \DynkinWDot[\DynkinWt{#1}{#2+#3+\x-1}]{#2+#3+2+\x}{0}; }
  }
}{#5} }

\newcommand{\dynkinC}[5][]{ \dynkin{
  \pgfmathsetmacro{\E}{#4-1}
  \DynkinLine{1 }{0}{#2  }{0};
  \DynkinDots{#2}{0}{#2+2}{0};
  \foreach \x in {1,...,#2} { \DynkinWDot[\DynkinWt{#1}{\x-1}]{\x}{0}; }
  \ifthenelse{#3=0}{
    \DynkinLine{#2+2}{0}{#2+\E+1}{0};
    \DynkinDoubleLine{#2+\E+2}{0}{#2+\E+1}{0};
    \foreach \x in {1,...,#4} { \DynkinWDot[\DynkinWt{#1}{#2+\x-1}]{#2+1+\x}{0}; }
  }{
    \DynkinLine{#2   +2}{0}{#2+#3   +1}{0};
    \DynkinDots{#2+#3+1}{0}{#2+#3   +3}{0};
    \DynkinLine{#2+#3+3}{0}{#2+#3+\E+1}{0};
    \DynkinDoubleLine{#2+#3+\E+2}{0}{#2+#3+\E+1}{0};
    \foreach \x in {1,...,#3} { \DynkinWDot[\DynkinWt{#1}{#2   +\x-1}]{#2   +1+\x}{0}; }
    \foreach \x in {1,...,\E} { \DynkinWDot[\DynkinWt{#1}{#2+#3+\x-1}]{#2+#3+2+\x}{0}; }
  }
}{#5} }
\newcommand{\dynkinCp}[5][]{ \dynkin{
  \pgfmathsetmacro{\E}{#4-1}
  \DynkinLine{1 }{0}{#2  }{0};
  \DynkinDots{#2}{0}{#2+2}{0};
  \foreach \x in {1,...,#2} { \DynkinWDot[\DynkinWt{#1}{\x-1}]{\x}{0}; }
  \ifthenelse{#3=0}{
    \DynkinLine{#2+2}{0}{#2+\E+1}{0};
    \DynkinDoubleLine{#2+\E+2}{0}{#2+\E+1}{0};
    \foreach \x in {1,...,\E} { \DynkinWDot[\DynkinWt{#1}{#2+\x-1}]{#2+1+\x}{0}; }
    \DynkinXDot[\DynkinWt{#1}{#2+#4-1}]{#2+#4+1}{0};
  }{
    \DynkinLine{#2   +2}{0}{#2+#3   +1}{0};
    \DynkinDots{#2+#3+1}{0}{#2+#3   +3}{0};
    \DynkinLine{#2+#3+3}{0}{#2+#3+#4+1}{0};
    \DynkinDoubleLine{#2+#3+\E+3}{0}{#2+#3+\E+2}{0};
    \foreach \x in {1,...,#3} { \DynkinWDot[\DynkinWt{#1}{#2   +\x-1}]{#2   +1+\x}{0}; }
    \foreach \x in {1,...,\E} { \DynkinWDot[\DynkinWt{#1}{#2+#3+\x-1}]{#2+#3+2+\x}{0}; }
    \DynkinXDot[\DynkinWt{#1}{#2+#3+#4-1}]{#2+#3+#4+2}{0};
  }
}{#5} }

\newcommand{\dynkinD}[5][]{ \dynkin{
  \DynkinLine{1 }{0}{#2  }{0};
  \DynkinDots{#2}{0}{#2+2}{0};
  \foreach \x in {1,...,#2} { \DynkinWDot[\DynkinWt{#1}{\x-1}]{\x}{0}; }
  \ifthenelse{#3=0}{
    \DynkinLine{#2   +2}{0}{#2+#4+1}{ 0   };
    \DynkinLine{#2+#4+1}{0}{#2+#4+2}{ 0.75};
    \DynkinLine{#2+#4+1}{0}{#2+#4+2}{-0.75};
    \foreach \x in {1,...,#4} { \DynkinWDot[\DynkinWt{#1}{#2+\x-1}]{#2+1+\x}{0}; }
    \DynkinWDot[right:\DynkinWt{#1}{#2+#4  }]{#2+#4+2}{ 0.75};
    \DynkinWDot[right:\DynkinWt{#1}{#2+#4+1}]{#2+#4+2}{-0.75};
  }{
    \DynkinLine{#2      +2}{0}{#2+#3   +1}{ 0   };
    \DynkinDots{#2+#3   +1}{0}{#2+#3   +3}{ 0   };
    \DynkinLine{#2+#3   +3}{0}{#2+#3+#4+2}{ 0   };
    \DynkinLine{#2+#3+#4+2}{0}{#2+#3+#4+3}{ 0.75};
    \DynkinLine{#2+#3+#4+2}{0}{#2+#3+#4+3}{-0.75};
    \foreach \x in {1,...,#3} { \DynkinWDot[\DynkinWt{#1}{#2   +\x-1}]{#2   +\x+1}{0}; }
    \foreach \x in {1,...,#4} { \DynkinWDot[\DynkinWt{#1}{#2+#3+\x-1}]{#2+#3+\x+2}{0}; }
    \DynkinWDot[right:\DynkinWt{#1}{#2+#3+#4  }]{#2+#3+#4+3}{ 0.75};
    \DynkinWDot[right:\DynkinWt{#1}{#2+#3+#4+1}]{#2+#3+#4+3}{-0.75};
  }
}{#5} }
\newcommand{\dynkinDp}[5][]{ \dynkin{
  \DynkinLine{1 }{0}{#2  }{0};
  \DynkinDots{#2}{0}{#2+2}{0};
  \foreach \x in {1,...,#2} { \DynkinWDot[\DynkinWt{#1}{\x-1}]{\x}{0}; }
  \ifthenelse{#3=0}{
    \DynkinLine{#2   +2}{0}{#2+#4+1}{ 0   };
    \DynkinLine{#2+#4+1}{0}{#2+#4+2}{ 0.75};
    \DynkinLine{#2+#4+1}{0}{#2+#4+2}{-0.75};
    \foreach \x in {1,...,#4} { \DynkinWDot[\DynkinWt{#1}{#2+\x-1}]{#2+1+\x}{0}; }
    \DynkinWDot[right:\DynkinWt{#1}{#2+#4  }]{#2+#4+2}{ 0.75};
    \DynkinXDot[right:\DynkinWt{#1}{#2+#4+1}]{#2+#4+2}{-0.75};
  }{
    \DynkinLine{#2      +2}{0}{#2+#3   +1}{ 0   };
    \DynkinDots{#2+#3   +1}{0}{#2+#3   +3}{ 0   };
    \DynkinLine{#2+#3   +3}{0}{#2+#3+#4+2}{ 0   };
    \DynkinLine{#2+#3+#4+2}{0}{#2+#3+#4+3}{ 0.75};
    \DynkinLine{#2+#3+#4+2}{0}{#2+#3+#4+3}{-0.75};
    \foreach \x in {1,...,#3} { \DynkinWDot[\DynkinWt{#1}{#2   +\x-1}]{#2   +\x+1}{0}; }
    \foreach \x in {1,...,#4} { \DynkinWDot[\DynkinWt{#1}{#2+#3+\x-1}]{#2+#3+\x+2}{0}; }
    \DynkinXDot[right:\DynkinWt{#1}{#2+#3+#4  }]{#2+#3+#4+3}{ 0.75};
    \DynkinWDot[right:\DynkinWt{#1}{#2+#3+#4+1}]{#2+#3+#4+3}{-0.75};
  }
}{#5} }

\newcommand{\dynkinBD}[5][]{ \dynkin{
  \pgfmathsetmacro{\E}{#4-1}
  \DynkinLine{1 }{0}{#2  }{0};
  \DynkinDots{#2}{0}{#2+2}{0};
  \foreach \x in {1,...,#2} { \DynkinWDot[\DynkinWt{#1}{\x-1}]{\x}{0}; }
  \ifthenelse{#3=0}{
    \DynkinLine{#2   +2  }{ 0  }{#2+#4+0.5}{ 0   };
    \DynkinLine{#2+#4+1.5}{ 0  }{#2+#4+2  }{ 0   };
    \DynkinLine{#2+#4+2  }{ 0  }{#2+#4+3  }{ 0.75};
    \DynkinLine{#2+#4+2  }{ 0  }{#2+#4+3  }{-0.75};
    \DynkinLine{#2+#4+1.5}{-1.5}{#2+#4+2  }{-1.5 };
    \DynkinDoubleLine{#2+#4+2}{-1.5}{#2+#4+3}{-1.5};
    \foreach \x in {1,...,\E} { \DynkinWDot[\DynkinWt{#1}{#2+\x-1}]{#2+1+\x}{0}; }
    \DynkinWDot[      \DynkinWt{#1}{#2+#4-1}]{#2+\E+3}{0};
    \DynkinWDot[right:\DynkinWt{#1}{#2+#4  }]{#2+#4+3}{ 0.75};
    \DynkinWDot[right:\DynkinWt{#1}{#2+#4+1}]{#2+#4+3}{-0.75};
    \DynkinWDot[below:\DynkinWt{#1}{#2+#4+2}]{#2+#4+2}{-1.5 };
    \DynkinWDot[below:\DynkinWt{#1}{#2+#4+3}]{#2+#4+3}{-1.5 };
    \DynkinCases{0.66}{#2+#4+1}{-1.5}{#2+#4+1}{0.7};
  }{
    \DynkinLine{#2      +2}{0}{#2+#3   +1}{ 0   };
    \DynkinDots{#2+#3   +1}{0}{#2+#3   +3}{ 0   };
    \DynkinLine{#2+#3   +3}{0}{#2+#3+#4+1.5}{ 0   };
    \DynkinLine{#2+#3+#4+2.5}{0}{#2+#3+#4+3}{0};
    \DynkinLine{#2+#3+#4+3}{0}{#2+#3+#4+4}{ 0.75};
    \DynkinLine{#2+#3+#4+3}{0}{#2+#3+#4+4}{-0.75};
    \DynkinLine{#2+#3+#4+2.5}{-1.5}{#2+#3+#4+3}{-1.5};
    \DynkinDoubleLine{#2+#3+#4+3}{-1.5}{#2+#3+#4+4}{-1.5};
    \foreach \x in {1,...,#3} { \DynkinWDot[\DynkinWt{#1}{#2   +\x-1}]{#2   +\x+1}{0}; }
    \foreach \x in {1,...,\E} { \DynkinWDot[\DynkinWt{#1}{#2+#3+\x-1}]{#2+#3+\x+2}{0}; }
    \DynkinWDot[      \DynkinWt{#1}{#2+#3+#4-1}]{#2+#3+#4+3}{0};
    \DynkinWDot[right:\DynkinWt{#1}{#2+#3+#4  }]{#2+#3+#4+4}{ 0.75};
    \DynkinWDot[right:\DynkinWt{#1}{#2+#3+#4+1}]{#2+#3+#4+4}{-0.75};
    \DynkinWDot[below:\DynkinWt{#1}{#2+#3+#4+2}]{#2+#3+#4+3}{-1.5};
    \DynkinWDot[below:\DynkinWt{#1}{#2+#3+#4+3}]{#2+#3+#4+4}{-1.5};
    \DynkinCases{0.66}{#2+#3+#4+2}{-1.5}{#2+#3+#4+2}{0.7};
  }
}{#5} }
\newcommand{\dynkinBDp}[5][]{ \dynkin{
  \pgfmathsetmacro{\E}{#4-1}
  \DynkinLine{1 }{0}{#2  }{0};
  \DynkinDots{#2}{0}{#2+2}{0};
  \DynkinXDot[\DynkinWt{#1}{0}]{1}{0};
  \foreach \x in {2,...,#2} { \DynkinWDot[\DynkinWt{#1}{\x-1}]{\x}{0}; }
  \ifthenelse{#3=0}{
    \DynkinLine{#2   +2  }{ 0  }{#2+#4+0.5}{ 0   };
    \DynkinLine{#2+#4+1.5}{ 0  }{#2+#4+2  }{ 0   };
    \DynkinLine{#2+#4+2  }{ 0  }{#2+#4+3  }{ 0.75};
    \DynkinLine{#2+#4+2  }{ 0  }{#2+#4+3  }{-0.75};
    \DynkinLine{#2+#4+1.5}{-1.5}{#2+#4+2  }{-1.5 };
    \DynkinDoubleLine{#2+#4+2}{-1.5}{#2+#4+3}{-1.5};
    \foreach \x in {1,...,\E} { \DynkinWDot[\DynkinWt{#1}{#2+\x-1}]{#2+1+\x}{0}; }
    \DynkinWDot[      \DynkinWt{#1}{#2+#4-1}]{#2+\E+3}{0};
    \DynkinWDot[right:\DynkinWt{#1}{#2+#4  }]{#2+#4+3}{ 0.75};
    \DynkinWDot[right:\DynkinWt{#1}{#2+#4+1}]{#2+#4+3}{-0.75};
    \DynkinWDot[below:\DynkinWt{#1}{#2+#4+2}]{#2+#4+2}{-1.5 };
    \DynkinWDot[below:\DynkinWt{#1}{#2+#4+3}]{#2+#4+3}{-1.5 };
    \DynkinCases{0.66}{#2+#4+1}{-1.5}{#2+#4+1}{0.7};
  }{
    \DynkinLine{#2      +2}{0}{#2+#3   +1}{ 0   };
    \DynkinDots{#2+#3   +1}{0}{#2+#3   +3}{ 0   };
    \DynkinLine{#2+#3   +3}{0}{#2+#3+#4+1.5}{ 0   };
    \DynkinLine{#2+#3+#4+2.5}{0}{#2+#3+#4+3}{0};
    \DynkinLine{#2+#3+#4+3}{0}{#2+#3+#4+4}{ 0.75};
    \DynkinLine{#2+#3+#4+3}{0}{#2+#3+#4+4}{-0.75};
    \DynkinLine{#2+#3+#4+2.5}{-1.5}{#2+#3+#4+3}{-1.5};
    \DynkinDoubleLine{#2+#3+#4+3}{-1.5}{#2+#3+#4+4}{-1.5};
    \foreach \x in {1,...,#3} { \DynkinWDot[\DynkinWt{#1}{#2   +\x-1}]{#2   +\x+1}{0}; }
    \foreach \x in {1,...,\E} { \DynkinWDot[\DynkinWt{#1}{#2+#3+\x-1}]{#2+#3+\x+2}{0}; }
    \DynkinWDot[      \DynkinWt{#1}{#2+#3+#4-1}]{#2+#3+#4+3}{0};
    \DynkinWDot[right:\DynkinWt{#1}{#2+#3+#4  }]{#2+#3+#4+4}{ 0.75};
    \DynkinWDot[right:\DynkinWt{#1}{#2+#3+#4+1}]{#2+#3+#4+4}{-0.75};
    \DynkinWDot[below:\DynkinWt{#1}{#2+#3+#4+2}]{#2+#3+#4+3}{-1.5};
    \DynkinWDot[below:\DynkinWt{#1}{#2+#3+#4+3}]{#2+#3+#4+4}{-1.5};
    \DynkinCases{0.66}{#2+#3+#4+2}{-1.5}{#2+#3+#4+2}{0.7};
  }
}{#5} }

\newcommand{\dynkinE}[3][]{ \dynkin[-0.5]{
  \pgfmathsetmacro{\E}{#2-1}
  \DynkinLine{1}{0}{#2-1}{ 0}
  \DynkinLine{3}{0}{3   }{-1};
  \foreach \x in {1,...,\E} { \DynkinWDot[\DynkinWt{#1}{\x-1}]{\x}{0}; }
  \DynkinWDot[right:\DynkinWt{#1}{#2-1}]{3}{-1};
}{#3} }
\newcommand{\dynkinEp}[3][]{ \dynkin[-0.5]{
  \pgfmathsetmacro{\E}{#2-2}
  \DynkinLine{1}{0}{#2-1}{ 0}
  \DynkinLine{3}{0}{3   }{-1};
  \foreach \x in {1,...,\E} { \DynkinWDot[\DynkinWt{#1}{\x-1}]{\x}{0}; }
  \DynkinXDot[      \DynkinWt{#1}{#2-2}]{#2-1}{ 0};
  \DynkinWDot[right:\DynkinWt{#1}{#2-1}]{   3}{-1};
}{#3} }

\newcommand{\dynkinF}[2][]{ \dynkin{
  \DynkinLine{1}{0}{2}{0};
  \DynkinDoubleLine{2}{0}{3}{0};
  \DynkinLine{3}{0}{4}{0};
  \foreach \x in {1,...,4} { \DynkinWDot[\DynkinWt{#1}{\x-1}]{\x}{0}; }
}{#2} }

\newcommand{\dynkinG}[2][]{ \dynkin{
  \DynkinTripleLine{1}{0}{2}{0};
  \foreach \x in {1,2} { \DynkinWDot[\DynkinWt{#1}{\x-1}]{\x}{0}; }
}{#2} }

\newcommand{\dynkinSLR}[5][]{ \dynkin{
  \DynkinLine{1 }{0}{#2  }{0};
  \DynkinDots{#2}{0}{#2+2}{0};
  \foreach \x in {1,...,#2} { \DynkinWDot[\DynkinWt{#1}{\x-1}]{\x}{0}; }
  \ifthenelse{#3=0}{
    \DynkinLine{#2+2}{0}{#2+#4+1}{0};
    \foreach \x in {1,...,#4} { \DynkinWDot[\DynkinWt{#1}{#2+\x-1}]{#2+1+\x}{0}; }
  }{
    \DynkinLine{#2   +2}{0}{#2+#3   +1}{0};
    \DynkinDots{#2+#3+1}{0}{#2+#3   +3}{0};
    \DynkinLine{#2+#3+3}{0}{#2+#3+#4+2}{0};
    \foreach \x in {1,...,#3} { \DynkinWDot[\DynkinWt{#1}{#2   +\x-1}]{#2+1+\x}{0}; }
    \foreach \x in {1,...,#4} { \DynkinWDot[\DynkinWt{#1}{#2+#3+\x-1}]{#2+#3+2+\x}{0}; }
  }
}{#5} }
\newcommand{\dynkinSLRp}[5][]{ \dynkin{
  \pgfmathsetmacro{\E}{#4-1}
  \DynkinLine{1 }{0}{#2  }{0};
  \DynkinDots{#2}{0}{#2+2}{0};
  \foreach \x in {1,...,#2} { \DynkinWDot[\DynkinWt{#1}{\x-1}]{\x}{0}; }
  \ifthenelse{#3=0}{
    \DynkinLine{#2+2}{0}{#2+#4+1}{0};
    \foreach \x in {1,...,\E} { \DynkinWDot[\DynkinWt{#1}{#2+\x-1}]{#2+1+\x}{0}; }
    \DynkinXDot[\DynkinWt{#1}{#2+#4-1}]{#2+#4+1}{0};
  }{
    \DynkinLine{#2   +2}{0}{#2+#3   +1}{0};
    \DynkinDots{#2+#3+1}{0}{#2+#3   +3}{0};
    \DynkinLine{#2+#3+3}{0}{#2+#3+#4+2}{0};
    \foreach \x in {1,...,#3} { \DynkinWDot[\DynkinWt{#1}{#2   +\x-1}]{#2   +\x+1}{0}; }
    \foreach \x in {1,...,\E} { \DynkinWDot[\DynkinWt{#1}{#2+#3+\x-1}]{#2+#3+\x+2}{0}; }
    \DynkinXDot[\DynkinWt{#1}{#2+#3+#4-1}]{#2+#3+#4+2}{0};
  }
}{#5} }

\newcommand{\dynkinSLC}[5][]{ \dynkin[-0.22]{
  \DynkinLine{1 }{ 0.75}{#2  }{ 0.75};
  \DynkinLine{1 }{-0.75}{#2  }{-0.75};
  \DynkinDots{#2}{ 0.75}{#2+2}{ 0.75};
  \DynkinDots{#2}{-0.75}{#2+2}{-0.75};
  \foreach \x in {1,...,#2} {
    \DynkinWDot[      \DynkinWt{#1}{         \x-1}]{\x}{ 0.75};
    \DynkinWDot[below:\DynkinWt{#1}{#2+#3+#4+\x-1}]{\x}{-0.75};
    \DynkinInvolution{\x}{-0.75}{\x}{0.75};
  }
  \ifthenelse{#3=0}{
    \DynkinLine{#2+2}{ 0.75}{#2+#4+1}{ 0.75};
    \DynkinLine{#2+2}{-0.75}{#2+#4+1}{-0.75};
    \foreach \x in {1,...,#4} {
      \DynkinWDot[      \DynkinWt{#1}{         #2+\x-1}]{#2+1+\x}{ 0.75};
      \DynkinWDot[below:\DynkinWt{#1}{#2+#3+#4+#2+\x-1}]{#2+1+\x}{-0.75};
      \DynkinInvolution{#2+1+\x}{-0.75}{#2+1+\x}{0.75};
    }
  }{
    \DynkinLine{#2   +2}{ 0.75}{#2+#3   +1}{ 0.75};
    \DynkinLine{#2   +2}{-0.75}{#2+#3   +1}{-0.75};
    \DynkinDots{#2+#3+1}{ 0.75}{#2+#3   +3}{ 0.75};
    \DynkinDots{#2+#3+1}{-0.75}{#2+#3   +3}{-0.75};
    \DynkinLine{#2+#3+3}{ 0.75}{#2+#3+#4+2}{ 0.75};
    \DynkinLine{#2+#3+3}{-0.75}{#2+#3+#4+2}{-0.75};
    \foreach \x in {1,...,#3} {
      \DynkinWDot[      \DynkinWt{#1}{         #2+\x-1}]{#2+1+\x}{ 0.75};
      \DynkinWDot[below:\DynkinWt{#1}{#2+#3+#4+#2+\x-1}]{#2+1+\x}{-0.75};
      \DynkinInvolution{#2+1+\x}{-0.75}{#2+1+\x}{0.75};
    }
    \foreach \x in {1,...,#4} {
      \DynkinWDot[      \DynkinWt{#1}{         #2+#3+\x-1}]{#2+#3+2+\x}{ 0.75};
      \DynkinWDot[below:\DynkinWt{#1}{#2+#3+#4+#2+#3+\x-1}]{#2+#3+2+\x}{-0.75};
      \DynkinInvolution{#2+#3+2+\x}{-0.75}{#2+#3+2+\x}{0.75};
    }
  }
}{#5} }
\newcommand{\dynkinSLCp}[5][]{ \dynkin[-0.22]{
  \pgfmathsetmacro{\E}{#4-1}
  \DynkinLine{1 }{ 0.75}{#2  }{ 0.75};
  \DynkinLine{1 }{-0.75}{#2  }{-0.75};
  \DynkinDots{#2}{ 0.75}{#2+2}{ 0.75};
  \DynkinDots{#2}{-0.75}{#2+2}{-0.75};
  \foreach \x in {1,...,#2} {
    \DynkinWDot[      \DynkinWt{#1}{         \x-1}]{\x}{ 0.75};
    \DynkinWDot[below:\DynkinWt{#1}{#2+#3+#4+\x-1}]{\x}{-0.75};
    \DynkinInvolution{\x}{-0.75}{\x}{0.75};
  }
  \ifthenelse{#3=0}{
    \DynkinLine{#2+2}{ 0.75}{#2+#4+1}{ 0.75};
    \DynkinLine{#2+2}{-0.75}{#2+#4+1}{-0.75};
    \foreach \x in {1,...,\E} {
      \DynkinWDot[      \DynkinWt{#1}{      #2+\x-1}]{#2+1+\x}{ 0.75};
      \DynkinWDot[below:\DynkinWt{#1}{#2+#4+#2+\x-1}]{#2+1+\x}{-0.75};
      \DynkinInvolution{#2+1+\x}{-0.75}{#2+1+\x}{0.75};
    }
    \DynkinXDot[      \DynkinWt{#1}{      #2+#4-1}]{#2+#4+1}{0.75};
    \DynkinXDot[below:\DynkinWt{#1}{#2+#4+#2+#4-1}]{#2+#4+1}{-0.75};
    \DynkinInvolution{#2+#3+#4+1}{-0.75}{#2+#3+#4+1}{0.75};
  }{
    \DynkinLine{#2   +2}{ 0.75}{#2+#3   +1}{ 0.75};
    \DynkinLine{#2   +2}{-0.75}{#2+#3   +1}{-0.75};
    \DynkinDots{#2+#3+1}{ 0.75}{#2+#3   +3}{ 0.75};
    \DynkinDots{#2+#3+1}{-0.75}{#2+#3   +3}{-0.75};
    \DynkinLine{#2+#3+3}{ 0.75}{#2+#3+#4+2}{ 0.75};
    \DynkinLine{#2+#3+3}{-0.75}{#2+#3+#4+2}{-0.75};
    \foreach \x in {1,...,#3} {
      \DynkinWDot[      \DynkinWt{#1}{         #2+\x-1}]{#2+1+\x}{ 0.75};
      \DynkinWDot[below:\DynkinWt{#1}{#2+#3+#4+#2+\x-1}]{#2+1+\x}{-0.75};
      \DynkinInvolution{#2+1+\x}{-0.75}{#2+1+\x}{0.75};
    }
    \foreach \x in {1,...,\E} {
      \DynkinWDot[      \DynkinWt{#1}{         #2+#3+\x-1}]{#2+#3+2+\x}{ 0.75};
      \DynkinWDot[below:\DynkinWt{#1}{#2+#3+#4+#2+#3+\x-1}]{#2+#3+2+\x}{-0.75};
      \DynkinInvolution{#2+#3+2+\x}{-0.75}{#2+#3+2+\x}{0.75};
    }
    \DynkinXDot[      \DynkinWt{#1}{         #2+#3+#4-1}]{#2+#3+#4+2}{ 0.75};
    \DynkinXDot[below:\DynkinWt{#1}{#2+#3+#4+#2+#3+#4-1}]{#2+#3+#4+2}{-0.75};
    \DynkinInvolution{#2+#3+#4+2}{-0.75}{#2+#3+#4+2}{0.75};
  }
}{#5} }

\newcommand{\dynkinSU}[5][]{ \dynkin[-0.25]{
  \DynkinLine{1 }{ 0.75}{#2  }{ 0.75};
  \DynkinLine{1 }{-0.75}{#2  }{-0.75};
  \DynkinDots{#2}{ 0.75}{#2+2}{ 0.75};
  \DynkinDots{#2}{-0.75}{#2+2}{-0.75};
  \foreach \x in {1,...,#2} {
    \DynkinWDot[      \DynkinWt{#1}{         \x-1}]{\x}{ 0.75};
    \DynkinWDot[below:\DynkinWt{#1}{#2+#3+#4+\x-1}]{\x}{-0.75};
    \DynkinInvolution{\x}{-0.75}{\x}{0.75};
  }
  \ifthenelse{#3=0}{
    \DynkinLine{#2   +2}{ 0.75}{#2+#4+1}{ 0.75};
    \DynkinLine{#2   +2}{-0.75}{#2+#4+1}{-0.75};
    \DynkinLine{#2+#4+1}{ 0.75}{#2+#4+2}{ 0   };
    \DynkinLine{#2+#4+1}{-0.75}{#2+#4+2}{ 0   };
    \foreach \x in {1,...,#4} {
      \DynkinWDot[      \DynkinWt{#1}{      #2+\x-1}]{#2+1+\x}{ 0.75};
      \DynkinWDot[below:\DynkinWt{#1}{#2+#4+#2+\x-1}]{#2+1+\x}{-0.75};
      \DynkinInvolution{#2+1+\x}{-0.75}{#2+1+\x}{0.75};
    }
    \DynkinWDot[right:\DynkinWt{#1}{#2+#4+#2+#4}]{#2+#4+2}{0};
  }{
    \DynkinLine{#2      +2}{ 0.75}{#2+#3   +1}{ 0.75};
    \DynkinLine{#2      +2}{-0.75}{#2+#3   +1}{-0.75};
    \DynkinDots{#2+#3   +1}{ 0.75}{#2+#3   +3}{ 0.75};
    \DynkinDots{#2+#3   +1}{-0.75}{#2+#3   +3}{-0.75};
    \DynkinLine{#2+#3   +3}{ 0.75}{#2+#3+#4+2}{ 0.75};
    \DynkinLine{#2+#3   +3}{-0.75}{#2+#3+#4+2}{-0.75};
    \DynkinLine{#2+#3+#4+2}{ 0.75}{#2+#3+#4+3}{ 0   };
    \DynkinLine{#2+#3+#4+2}{-0.75}{#2+#3+#4+3}{ 0   };
    \foreach \x in {1,...,#3} {
      \DynkinWDot[      \DynkinWt{#1}{         #2+\x-1}]{#2+1+\x}{ 0.75};
      \DynkinWDot[below:\DynkinWt{#1}{#2+#3+#4+#2+\x-1}]{#2+1+\x}{-0.75};
      \DynkinInvolution{#2+1+\x}{-0.75}{#2+1+\x}{0.75};
    }
    \foreach \x in {1,...,#4} {
      \DynkinWDot[      \DynkinWt{#1}{         #2+#3+\x-1}]{#2+#3+2+\x}{ 0.75};
      \DynkinWDot[below:\DynkinWt{#1}{#2+#3+#4+#2+#3+\x-1}]{#2+#3+2+\x}{-0.75};
      \DynkinInvolution{#2+#3+2+\x}{-0.75}{#2+#3+2+\x}{0.75};
    }
    \DynkinWDot[right:\DynkinWt{#1}{#2+#3+#4+#2+#3+#4}]{#2+#3+#4+3}{0};
  }
}{#5} }
\newcommand{\dynkinSUp}[5][]{ \dynkin[-0.25]{
  \DynkinLine{1 }{ 0.75}{#2  }{ 0.75};
  \DynkinLine{1 }{-0.75}{#2  }{-0.75};
  \DynkinDots{#2}{ 0.75}{#2+2}{ 0.75};
  \DynkinDots{#2}{-0.75}{#2+2}{-0.75};
  \foreach \x in {1,...,#2} {
    \DynkinWDot[      \DynkinWt{#1}{         \x-1}]{\x}{ 0.75};
    \DynkinWDot[below:\DynkinWt{#1}{#2+#3+#4+\x-1}]{\x}{-0.75};
    \DynkinInvolution{\x}{-0.75}{\x}{0.75};
  }
  \ifthenelse{#3=0}{
    \DynkinLine{#2   +2}{ 0.75}{#2+#4+1}{ 0.75};
    \DynkinLine{#2   +2}{-0.75}{#2+#4+1}{-0.75};
    \DynkinLine{#2+#4+1}{ 0.75}{#2+#4+2}{ 0   };
    \DynkinLine{#2+#4+1}{-0.75}{#2+#4+2}{ 0   };
    \foreach \x in {1,...,#4} {
      \DynkinWDot[      \DynkinWt{#1}{      #2+\x-1}]{#2+1+\x}{ 0.75};
      \DynkinWDot[below:\DynkinWt{#1}{#2+#4+#2+\x-1}]{#2+1+\x}{-0.75};
      \DynkinInvolution{#2+1+\x}{-0.75}{#2+1+\x}{0.75};
    }
    \DynkinXDot[right:\DynkinWt{#1}{#2+#4+#2+#4}]{#2+#4+2}{0};
  }{
    \DynkinLine{#2      +2}{ 0.75}{#2+#3   +1}{ 0.75};
    \DynkinLine{#2      +2}{-0.75}{#2+#3   +1}{-0.75};
    \DynkinDots{#2+#3   +1}{ 0.75}{#2+#3   +3}{ 0.75};
    \DynkinDots{#2+#3   +1}{-0.75}{#2+#3   +3}{-0.75};
    \DynkinLine{#2+#3   +3}{ 0.75}{#2+#3+#4+2}{ 0.75};
    \DynkinLine{#2+#3   +3}{-0.75}{#2+#3+#4+2}{-0.75};
    \DynkinLine{#2+#3+#4+2}{ 0.75}{#2+#3+#4+3}{ 0   };
    \DynkinLine{#2+#3+#4+2}{-0.75}{#2+#3+#4+3}{ 0   };
    \foreach \x in {1,...,#3} {
      \DynkinWDot[      \DynkinWt{#1}{         #2+\x-1}]{#2+1+\x}{ 0.75};
      \DynkinWDot[below:\DynkinWt{#1}{#2+#3+#4+#2+\x-1}]{#2+1+\x}{-0.75};
      \DynkinInvolution{#2+1+\x}{-0.75}{#2+1+\x}{0.75};
    }
    \foreach \x in {1,...,#4} {
      \DynkinWDot[      \DynkinWt{#1}{         #2+#3+\x-1}]{#2+#3+2+\x}{ 0.75};
      \DynkinWDot[below:\DynkinWt{#1}{#2+#3+#4+#2+#3+\x-1}]{#2+#3+2+\x}{-0.75};
      \DynkinInvolution{#2+#3+2+\x}{-0.75}{#2+#3+2+\x}{0.75};
    }
    \DynkinXDot[right:\DynkinWt{#1}{#2+#3+#4+#2+#3+#4}]{#2+#3+#4+3}{0};
  }
}{#5} }

\newcommand{\dynkinSLH}[5][]{ \dynkin{
  \DynkinLine{1 }{0}{#2  }{0};
  \DynkinDots{#2}{0}{#2+2}{0};
  \foreach \x in {1,...,#2} {
    \ifthenelse{\isodd{\x}}{
      \DynkinBDot[\DynkinWt{#1}{\x-1}]{\x}{0};
    }{
      \DynkinWDot[\DynkinWt{#1}{\x-1}]{\x}{0};
    }
  }
  \pgfmathparse{{#3}[0]}
  \ifthenelse{\pgfmathresult=0} {
    \DynkinLine{#2+2}{0}{#2+#4+1}{0};
    \foreach[count=\y] \x in {#4,...,1} {
      \ifthenelse{\isodd{\y}} {
        \DynkinBDot[\DynkinWt{#1}{#2+\x-1}]{#2+\x+1}{0};
      }{
        \DynkinWDot[\DynkinWt{#1}{#2+\x-1}]{#2+1+\x}{0};
      }
    }
  }{
    \pgfmathsetmacro{\M}{{#3}[0]}
    \DynkinLine{#2   +2}{0}{#2+\M   +1}{0};
    \DynkinDots{#2+\M+1}{0}{#2+\M   +3}{0};
    \DynkinLine{#2+\M+3}{0}{#2+\M+#4+2}{0};
    \foreach \x in {1,...,\M} {
      \pgfmathsetmacro{\MC}{\x+{#3}[1]}
      \ifthenelse{\isodd{\MC}}{
        \DynkinBDot[\DynkinWt{#1}{#2+\x-1}]{#2+\x+1}{0};
      }{
        \DynkinWDot[\DynkinWt{#1}{#2+\x-1}]{#2+\x+1}{0};
      }
    }
    \foreach[count=\y] \x in {#4,...,1} {
      \ifthenelse{\isodd{\y}} {
        \DynkinBDot[\DynkinWt{#1}{#2+\M+\x-1}]{#2+\M+2+\x}{0};
      }{
        \DynkinWDot[\DynkinWt{#1}{#2+\M+\x-1}]{#2+\M+2+\x}{0};
      }
    }
  }
}{#5} }
\newcommand{\dynkinSLHp}[5][]{ \dynkin{
  \DynkinLine{1 }{0}{#2  }{0};
  \DynkinDots{#2}{0}{#2+2}{0};
  \foreach \x in {1,...,#2} {
    \ifthenelse{\isodd{\x}}{
      \DynkinBDot[\DynkinWt{#1}{\x-1}]{\x}{0};
    }{
      \DynkinWDot[\DynkinWt{#1}{\x-1}]{\x}{0};
    }
  }
  \pgfmathparse{{#3}[0]}
  \ifthenelse{\pgfmathresult=0} {
    \DynkinLine{#2+2}{0}{#2+#4+1}{0};
    \pgfmathsetmacro{\E}{#4-2}
    \foreach[count=\y] \x in {\E,...,1} {
      \ifthenelse{\isodd{\y}} {
        \DynkinBDot[\DynkinWt{#1}{#2+\x-1}]{#2+\x+1}{0};
      }{
        \DynkinWDot[\DynkinWt{#1}{#2+\x-1}]{#2+1+\x}{0};
      }
    }
    \DynkinXDot[\DynkinWt{#1}{#2+#4-2}]{#2+#4  }{0};
    \DynkinBDot[\DynkinWt{#1}{#2+#4-1}]{#2+#4+1}{0};
  }{
    \pgfmathsetmacro{\M}{{#3}[0]}
    \DynkinLine{#2   +2}{0}{#2+\M   +1}{0};
    \DynkinDots{#2+\M+1}{0}{#2+\M   +3}{0};
    \DynkinLine{#2+\M+3}{0}{#2+\M+#4+2}{0};
    \foreach \x in {1,...,\M} {
      \pgfmathsetmacro{\MC}{\x+{#3}[1]}
      \ifthenelse{\isodd{\MC}}{
        \DynkinBDot[\DynkinWt{#1}{#2+\x-1}]{#2+\x+1}{0};
      }{
        \DynkinWDot[\DynkinWt{#1}{#2+\x-1}]{#2+\x+1}{0};
      }
    }
    \pgfmathsetmacro{\E}{#4-2}
    \foreach[count=\y] \x in {\E,...,1} {
      \ifthenelse{\isodd{\y}} {
        \DynkinBDot[\DynkinWt{#1}{#2+\M+\x-1}]{#2+\M+2+\x}{0};
      }{
        \DynkinWDot[\DynkinWt{#1}{#2+\M+\x-1}]{#2+\M+2+\x}{0};
      }
    }
    \DynkinXDot[\DynkinWt{#1}{#2+\M+#4-2}]{#2+\M+#4+1}{0};
    \DynkinBDot[\DynkinWt{#1}{#2+\M+#4-1}]{#2+\M+#4+2}{0};
  }
}{#5} }

\newcommand{\dynkinSOs}[5][]{ \dynkin{
  \DynkinLine{1 }{0}{#2  }{0};
  \DynkinDots{#2}{0}{#2+2}{0};
  \foreach \x in {1,...,#2} {
    \ifthenelse{\isodd{\x}}{
      \DynkinBDot[\DynkinWt{#1}{\x-1}]{\x}{0};
    }{
      \DynkinWDot[\DynkinWt{#1}{\x-1}]{\x}{0};
    }
  }
  \pgfmathparse{{#3}[0]}
  \ifthenelse{\pgfmathresult=0} {
    \DynkinLine{#2   +2}{0}{#2+#4+1}{ 0   };
    \DynkinLine{#2+#4+1}{0}{#2+#4+2}{ 0.75};
    \DynkinLine{#2+#4+1}{0}{#2+#4+2}{-0.75};
    \foreach[count=\y] \x in {#4,...,1} {
      \ifthenelse{\isodd{\y}} {
        \DynkinBDot[\DynkinWt{#1}{#2+\x-1}]{#2+\x+1}{0};
      }{
        \DynkinWDot[\DynkinWt{#1}{#2+\x-1}]{#2+1+\x}{0};
      }
    }
    \DynkinBDot[right:\DynkinWt{#1}{#2+#4  }]{#2+#4+2}{ 0.75};
    \DynkinWDot[right:\DynkinWt{#1}{#2+#4+1}]{#2+#4+2}{-0.75};	
  }{
    \pgfmathsetmacro{\M}{{#3}[0]}
    \DynkinLine{#2      +2}{0}{#2+\M   +1}{ 0   };
    \DynkinDots{#2+\M   +1}{0}{#2+\M   +3}{ 0   };
    \DynkinLine{#2+\M   +3}{0}{#2+\M+#4+2}{ 0   };
    \DynkinLine{#2+\M+#4+2}{0}{#2+\M+#4+3}{ 0.75};
    \DynkinLine{#2+\M+#4+2}{0}{#2+\M+#4+3}{-0.75};
    \foreach \x in {1,...,\M} {
      \pgfmathsetmacro{\MC}{\x+{#3}[1]}
      \ifthenelse{\isodd{\MC}}{
        \DynkinBDot[\DynkinWt{#1}{#2+\x-1}]{#2+\x+1}{0};
      }{
        \DynkinWDot[\DynkinWt{#1}{#2+\x-1}]{#2+\x+1}{0};
      }
    }
    \foreach[count=\y] \x in {#4,...,1} {
      \ifthenelse{\isodd{\y}} {
        \DynkinBDot[\DynkinWt{#1}{#2+\M+\x-1}]{#2+\M+2+\x}{0};
      }{
        \DynkinWDot[\DynkinWt{#1}{#2+\M+\x-1}]{#2+\M+2+\x}{0};
      }
    }
    \DynkinBDot[right:\DynkinWt{#1}{#2+\M+#4  }]{#2+\M+#4+3}{ 0.75};
    \DynkinWDot[right:\DynkinWt{#1}{#2+\M+#4+1}]{#2+\M+#4+3}{-0.75};
  }
}{#5} }
\newcommand{\dynkinSOsp}[5][]{ \dynkin{
  \DynkinLine{1 }{0}{#2  }{0};
  \DynkinDots{#2}{0}{#2+2}{0};
  \foreach \x in {1,...,#2} {
    \ifthenelse{\isodd{\x}}{
      \DynkinBDot[\DynkinWt{#1}{\x-1}]{\x}{0};
    }{
      \DynkinWDot[\DynkinWt{#1}{\x-1}]{\x}{0};
    }
  }
  \pgfmathparse{{#3}[0]}
  \ifthenelse{\pgfmathresult=0} {
    \DynkinLine{#2   +2}{0}{#2+#4+1}{ 0   };
    \DynkinLine{#2+#4+1}{0}{#2+#4+2}{ 0.75};
    \DynkinLine{#2+#4+1}{0}{#2+#4+2}{-0.75};
    \foreach[count=\y] \x in {#4,...,1} {
      \ifthenelse{\isodd{\y}} {
        \DynkinBDot[\DynkinWt{#1}{#2+\x-1}]{#2+\x+1}{0};
      }{
        \DynkinWDot[\DynkinWt{#1}{#2+\x-1}]{#2+1+\x}{0};
      }
    }
    \DynkinBDot[right:\DynkinWt{#1}{#2+#4  }]{#2+#4+2}{ 0.75};
    \DynkinXDot[right:\DynkinWt{#1}{#2+#4+1}]{#2+#4+2}{-0.75};	
  }{
    \pgfmathsetmacro{\M}{{#3}[0]}
    \DynkinLine{#2      +2}{0}{#2+\M   +1}{ 0   };
    \DynkinDots{#2+\M   +1}{0}{#2+\M   +3}{ 0   };
    \DynkinLine{#2+\M   +3}{0}{#2+\M+#4+2}{ 0   };
    \DynkinLine{#2+\M+#4+2}{0}{#2+\M+#4+3}{ 0.75};
    \DynkinLine{#2+\M+#4+2}{0}{#2+\M+#4+3}{-0.75};
    \foreach \x in {1,...,\M} {
      \pgfmathsetmacro{\MC}{\x+{#3}[1]}
      \ifthenelse{\isodd{\MC}}{
        \DynkinBDot[\DynkinWt{#1}{#2+\x-1}]{#2+\x+1}{0};
      }{
        \DynkinWDot[\DynkinWt{#1}{#2+\x-1}]{#2+\x+1}{0};
      }
    }
    \foreach[count=\y] \x in {#4,...,1} {
      \ifthenelse{\isodd{\y}} {
        \DynkinBDot[\DynkinWt{#1}{#2+\M+\x-1}]{#2+\M+2+\x}{0};
      }{
        \DynkinWDot[\DynkinWt{#1}{#2+\M+\x-1}]{#2+\M+2+\x}{0};
      }
    }
    \DynkinBDot[right:\DynkinWt{#1}{#2+\M+#4  }]{#2+\M+#4+3}{ 0.75};
    \DynkinXDot[right:\DynkinWt{#1}{#2+\M+#4+1}]{#2+\M+#4+3}{-0.75};
  }
}{#5} }

\newcommand{\dynkinSO}[5][]{ \dynkin{
  \pgfmathsetmacro{\E}{#4-1}
  \DynkinLine{1 }{0}{#2  }{0};
  \DynkinDots{#2}{0}{#2+2}{0};
  \foreach \x in {1,...,#2} { \DynkinWDot[\DynkinWt{#1}{\x-1}]{\x}{0}; }
  \ifthenelse{#3=0}{
    \DynkinLine{#2   +2  }{ 0  }{#2+#4+0.5}{ 0   };
    \DynkinLine{#2+#4+1.5}{ 0  }{#2+#4+2  }{ 0   };
    \DynkinLine{#2+#4+2  }{ 0  }{#2+#4+3  }{ 0.75};
    \DynkinLine{#2+#4+2  }{ 0  }{#2+#4+3  }{-0.75};
    \DynkinLine{#2+#4+1.5}{-1.5}{#2+#4+2  }{-1.5 };
    \DynkinDoubleLine{#2+#4+2}{-1.5}{#2+#4+3}{-1.5};
    \foreach \x in {1,...,\E} { \DynkinWDot[\DynkinWt{#1}{#2+\x-1}]{#2+1+\x}{0}; }
    \DynkinWDot[      \DynkinWt{#1}{#2+#4-1}]{#2+\E+3}{0};
    \DynkinWDot[right:\DynkinWt{#1}{#2+#4  }]{#2+#4+3}{ 0.75};
    \DynkinWDot[right:\DynkinWt{#1}{#2+#4+1}]{#2+#4+3}{-0.75};
    \DynkinWDot[below:\DynkinWt{#1}{#2+#4+2}]{#2+#4+2}{-1.5 };
    \DynkinWDot[below:\DynkinWt{#1}{#2+#4+3}]{#2+#4+3}{-1.5 };
    \DynkinCases{0.66}{#2+#4+1}{-1.5}{#2+#4+1}{0.7};
  }{
    \DynkinLine{#2      +2}{0}{#2+#3   +1}{ 0   };
    \DynkinDots{#2+#3   +1}{0}{#2+#3   +3}{ 0   };
    \DynkinLine{#2+#3   +3}{0}{#2+#3+#4+1.5}{ 0   };
    \DynkinLine{#2+#3+#4+2.5}{0}{#2+#3+#4+3}{0};
    \DynkinLine{#2+#3+#4+3}{0}{#2+#3+#4+4}{ 0.75};
    \DynkinLine{#2+#3+#4+3}{0}{#2+#3+#4+4}{-0.75};
    \DynkinLine{#2+#3+#4+2.5}{-1.5}{#2+#3+#4+3}{-1.5};
    \DynkinDoubleLine{#2+#3+#4+3}{-1.5}{#2+#3+#4+4}{-1.5};
    \foreach \x in {1,...,#3} { \DynkinWDot[\DynkinWt{#1}{#2   +\x-1}]{#2   +\x+1}{0}; }
    \foreach \x in {1,...,\E} { \DynkinWDot[\DynkinWt{#1}{#2+#3+\x-1}]{#2+#3+\x+2}{0}; }
    \DynkinWDot[      \DynkinWt{#1}{#2+#3+#4-1}]{#2+#3+#4+3}{0};
    \DynkinWDot[right:\DynkinWt{#1}{#2+#3+#4  }]{#2+#3+#4+4}{ 0.75};
    \DynkinWDot[right:\DynkinWt{#1}{#2+#3+#4+1}]{#2+#3+#4+4}{-0.75};
    \DynkinWDot[below:\DynkinWt{#1}{#2+#3+#4+2}]{#2+#3+#4+3}{-1.5};
    \DynkinWDot[below:\DynkinWt{#1}{#2+#3+#4+3}]{#2+#3+#4+4}{-1.5};
    \DynkinCases{0.66}{#2+#3+#4+2}{-1.5}{#2+#3+#4+2}{0.7};
  }
}{#5} }
\newcommand{\dynkinSOp}[5][]{ \dynkin{
  \pgfmathsetmacro{\E}{#4-1}
  \DynkinLine{1 }{0}{#2  }{0};
  \DynkinDots{#2}{0}{#2+2}{0};
  \DynkinXDot[\DynkinWt{#1}{0}]{1}{0};
  \foreach \x in {2,...,#2} { \DynkinWDot[\DynkinWt{#1}{\x-1}]{\x}{0}; }
  \ifthenelse{#3=0}{
    \DynkinLine{#2   +2  }{ 0  }{#2+#4+0.5}{ 0   };
    \DynkinLine{#2+#4+1.5}{ 0  }{#2+#4+2  }{ 0   };
    \DynkinLine{#2+#4+2  }{ 0  }{#2+#4+3  }{ 0.75};
    \DynkinLine{#2+#4+2  }{ 0  }{#2+#4+3  }{-0.75};
    \DynkinLine{#2+#4+1.5}{-1.5}{#2+#4+2  }{-1.5 };
    \DynkinDoubleLine{#2+#4+2}{-1.5}{#2+#4+3}{-1.5};
    \foreach \x in {1,...,\E} { \DynkinWDot[\DynkinWt{#1}{#2+\x-1}]{#2+1+\x}{0}; }
    \DynkinWDot[      \DynkinWt{#1}{#2+#4-1}]{#2+\E+3}{0};
    \DynkinWDot[right:\DynkinWt{#1}{#2+#4  }]{#2+#4+3}{ 0.75};
    \DynkinWDot[right:\DynkinWt{#1}{#2+#4+1}]{#2+#4+3}{-0.75};
    \DynkinWDot[below:\DynkinWt{#1}{#2+#4+2}]{#2+#4+2}{-1.5 };
    \DynkinWDot[below:\DynkinWt{#1}{#2+#4+3}]{#2+#4+3}{-1.5 };
    \DynkinCases{0.71}{#2+#4+1}{-1.5}{#2+#4+1}{0.7};
  }{
    \DynkinLine{#2      +2}{0}{#2+#3   +1}{ 0   };
    \DynkinDots{#2+#3   +1}{0}{#2+#3   +3}{ 0   };
    \DynkinLine{#2+#3   +3}{0}{#2+#3+#4+1.5}{ 0   };
    \DynkinLine{#2+#3+#4+2.5}{0}{#2+#3+#4+3}{0};
    \DynkinLine{#2+#3+#4+3}{0}{#2+#3+#4+4}{ 0.75};
    \DynkinLine{#2+#3+#4+3}{0}{#2+#3+#4+4}{-0.75};
    \DynkinLine{#2+#3+#4+2.5}{-1.5}{#2+#3+#4+3}{-1.5};
    \DynkinDoubleLine{#2+#3+#4+3}{-1.5}{#2+#3+#4+4}{-1.5};
    \foreach \x in {1,...,#3} { \DynkinWDot[\DynkinWt{#1}{#2   +\x-1}]{#2   +\x+1}{0}; }
    \foreach \x in {1,...,\E} { \DynkinWDot[\DynkinWt{#1}{#2+#3+\x-1}]{#2+#3+\x+2}{0}; }
    \DynkinWDot[      \DynkinWt{#1}{#2+#3+#4-1}]{#2+#3+#4+3}{0};
    \DynkinWDot[right:\DynkinWt{#1}{#2+#3+#4  }]{#2+#3+#4+4}{ 0.75};
    \DynkinWDot[right:\DynkinWt{#1}{#2+#3+#4+1}]{#2+#3+#4+4}{-0.75};
    \DynkinWDot[below:\DynkinWt{#1}{#2+#3+#4+2}]{#2+#3+#4+3}{-1.5};
    \DynkinWDot[below:\DynkinWt{#1}{#2+#3+#4+3}]{#2+#3+#4+4}{-1.5};
    \DynkinCases{0.71}{#2+#3+#4+2}{-1.5}{#2+#3+#4+2}{0.7};
  }
}{#5} }
\newcommand{\dynkinSOpq}[5][]{ \dynkin{
  \pgfmathsetmacro{\E}{#3-1}
  \DynkinLine{1 }{0}{#2  }{0};
  \DynkinDots{#2}{0}{#2+2}{0};
  \DynkinLine{#2     +2  }{ 0  }{#2+2   +1  }{ 0   };
  \DynkinDots{#2+2   +1  }{ 0  }{#2+2   +3  }{ 0   };
  \DynkinLine{#2+2   +3  }{ 0  }{#2+2+#3+1.5}{ 0   };
  \DynkinLine{#2+2+#3+2.5}{ 0  }{#2+2+#3+3  }{ 0   };
  \DynkinLine{#2+2+#3+3  }{ 0  }{#2+2+#3+4  }{ 0.75};
  \DynkinLine{#2+2+#3+3  }{ 0  }{#2+2+#3+4  }{-0.75};
  \DynkinLine{#2+2+#3+2.5}{-1.5}{#2+2+#3+3  }{-1.5 };
  \DynkinDoubleLine{#2+2+#3+3}{-1.5}{#2+2+#3+4}{-1.5};
  \foreach \x in {1,...,#2} { \DynkinWDot[\DynkinWt{#1}{     \x-1}]{     \x  }{0}; }
  \foreach \x in {1,...,\E} { \DynkinBDot[\DynkinWt{#1}{#2+2+\x-1}]{#2+2+\x+2}{0}; }
  \DynkinWDot[      \DynkinWt{#1}{#2     -1}]{#2     +2}{ 0   };
  \DynkinBDot[      \DynkinWt{#1}{#2       }]{#2     +3}{ 0   };
  \DynkinBDot[      \DynkinWt{#1}{#2+2+#3-1}]{#2+2+#3+3}{ 0   };
  \DynkinBDot[right:\DynkinWt{#1}{#2+2+#3  }]{#2+2+#3+4}{ 0.75};
  \DynkinBDot[right:\DynkinWt{#1}{#2+2+#3+1}]{#2+2+#3+4}{-0.75};
  \DynkinBDot[below:\DynkinWt{#1}{#2+2+#3+2}]{#2+2+#3+3}{-1.5 };
  \DynkinBDot[below:\DynkinWt{#1}{#2+2+#3+3}]{#2+2+#3+4}{-1.5 };
  \DynkinCases{0.66}{#2+2+#3+2}{-1.5}{#2+2+#3+2}{0.7};
  \DynkinUnderbrace{$#4$}{   1}{-0.6}{#2   +2}{-0.6};
}{#5} }

\newcommand{\dynkinSOpqshortp}[4][]{ \dynkin{
  \pgfmathsetmacro{\E}{#3-1}
  \DynkinLine{1 }{0}{#2  }{0};
  \DynkinDots{#2}{0}{#2+2}{0};
  \DynkinLine{#2   +2  }{ 0  }{#2+#3+0.5}{ 0   };
  \DynkinLine{#2+#3+1.5}{ 0  }{#2+#3+2  }{ 0   };
  \DynkinLine{#2+#3+2  }{ 0  }{#2+#3+3  }{ 0.75};
  \DynkinLine{#2+#3+2  }{ 0  }{#2+#3+3  }{-0.75};
  \DynkinLine{#2+#3+1.5}{-1.5}{#2+#3+2  }{-1.5 };
  \DynkinDoubleLine{#2+#3+2}{-1.5}{#2+#3+3}{-1.5};
  \DynkinXDot[\DynkinWt{#1}{0}]{1}{0};
  \foreach \x in {2,...,#2} { \DynkinWDot[\DynkinWt{#1}{   \x-1}]{     \x  }{0}; }
  \foreach \x in {1,...,\E} { \DynkinBDot[\DynkinWt{#1}{#2+\x-1}]{#2+\x+1}{0}; }
  \DynkinBDot[      \DynkinWt{#1}{#2+#3-1}]{#2+#3+2}{ 0   };
  \DynkinBDot[right:\DynkinWt{#1}{#2+#3  }]{#2+#3+3}{ 0.75};
  \DynkinBDot[right:\DynkinWt{#1}{#2+#3+1}]{#2+#3+3}{-0.75};
  \DynkinBDot[below:\DynkinWt{#1}{#2+#3+2}]{#2+#3+2}{-1.5 };
  \DynkinBDot[below:\DynkinWt{#1}{#2+#3+3}]{#2+#3+3}{-1.5 };
  \DynkinCases{0.66}{#2+#3+1}{-1.5}{#2+#3+1}{0.7};
}{#4} }

\newcommand{\dynkinEIV}[2][]{ \dynkin[-0.5]{
    \DynkinLine{0}{0}{4}{ 0};
    \DynkinLine{2}{0}{2}{-1};
    \DynkinWDot[      \DynkinWt{#1}{0}]{0}{ 0};
    \DynkinBDot[      \DynkinWt{#1}{1}]{1}{ 0};
    \DynkinBDot[      \DynkinWt{#1}{2}]{2}{ 0};
    \DynkinBDot[      \DynkinWt{#1}{3}]{3}{ 0};
    \DynkinWDot[      \DynkinWt{#1}{4}]{4}{ 0};
    \DynkinBDot[right:\DynkinWt{#1}{5}]{2}{-1};
}{#2} }
\newcommand{\dynkinEIVp}[2][]{ \dynkin[-0.5]{
    \DynkinLine{0}{0}{4}{ 0};
    \DynkinLine{2}{0}{2}{-1};
    \DynkinWDot[      \DynkinWt{#1}{0}]{0}{ 0};
    \DynkinBDot[      \DynkinWt{#1}{1}]{1}{ 0};
    \DynkinBDot[      \DynkinWt{#1}{2}]{2}{ 0};
    \DynkinBDot[      \DynkinWt{#1}{3}]{3}{ 0};
    \DynkinXDot[      \DynkinWt{#1}{4}]{4}{ 0};
    \DynkinBDot[right:\DynkinWt{#1}{5}]{2}{-1};
}{#2} }

\newcommand{\dynkinEVII}[2][]{ \dynkin[-0.5]{
    \DynkinLine{0}{0}{5}{ 0};
    \DynkinLine{2}{0}{2}{-1};
    \DynkinWDot[      \DynkinWt{#1}{0}]{0}{ 0};
    \DynkinBDot[      \DynkinWt{#1}{1}]{1}{ 0};
    \DynkinBDot[      \DynkinWt{#1}{2}]{2}{ 0};
    \DynkinBDot[      \DynkinWt{#1}{3}]{3}{ 0};
    \DynkinWDot[      \DynkinWt{#1}{4}]{4}{ 0};
    \DynkinWDot[      \DynkinWt{#1}{5}]{5}{ 0};
    \DynkinBDot[right:\DynkinWt{#1}{6}]{2}{-1};
}{#2} }
\newcommand{\dynkinEVIIp}[2][]{ \dynkin[-0.5]{
    \DynkinLine{0}{0}{5}{ 0};
    \DynkinLine{2}{0}{2}{-1};
    \DynkinWDot[      \DynkinWt{#1}{0}]{0}{ 0};
    \DynkinBDot[      \DynkinWt{#1}{1}]{1}{ 0};
    \DynkinBDot[      \DynkinWt{#1}{2}]{2}{ 0};
    \DynkinBDot[      \DynkinWt{#1}{3}]{3}{ 0};
    \DynkinWDot[      \DynkinWt{#1}{4}]{4}{ 0};
    \DynkinXDot[      \DynkinWt{#1}{5}]{5}{ 0};
    \DynkinBDot[right:\DynkinWt{#1}{6}]{2}{-1};
}{#2} }

\newcommand{\mr}[1]{\mathrm{#1}}
\newcommand{\bs}[1]{\boldsymbol{#1}}

\newcommand{\cA}{\mathcal{A}}
\newcommand{\cB}{\mathcal{B}}

\newcommand{\cD}{\mathcal{D}}
\newcommand{\cE}{\mathcal{E}}

\newcommand{\cH}{\mathcal{H}}

\newcommand{\cJ}{\mathcal{J}}

\newcommand{\cL}{\mathcal{L}}

\newcommand{\cN}{\mathcal{N}}
\newcommand{\cO}{\mathcal{O}}

\newcommand{\cQ}{\mathcal{Q}}

\newcommand{\cT}{\mathcal{T}}
\newcommand{\cU}{\mathcal{U}}
\newcommand{\cV}{\mathcal{V}}
\newcommand{\cW}{\mathcal{W}}

\newcommand{\cZ}{\mathcal{Z}}

      \newcommand{\fb}{\mathfrak{b}}

      \newcommand{\fg}{\mathfrak{g}}
      \newcommand{\fh}{\mathfrak{h}}

      \newcommand{\fk}{\mathfrak{k}}
      
      \newcommand{\fm}{\mathfrak{m}}

      \newcommand{\fp}{\mathfrak{p}}
      \newcommand{\fq}{\mathfrak{q}}

      \newcommand{\ft}{\mathfrak{t}}

      \newcommand{\fw}{\mathfrak{w}}

      \newcommand{\fz}{\mathfrak{z}}

\newcommand{\bC}{\mathbb{C}}

\newcommand{\bF}{\mathbb{F}}

\newcommand{\bH}{\mathbb{H}}

\newcommand{\bJ}{\mathbb{J}}
        \newcommand{\bk}{\Bbbk}

\newcommand{\bN}{\mathbb{N}}
\newcommand{\bO}{\mathbb{O}}
\newcommand{\bP}{\mathbb{P}}

\newcommand{\bR}{\mathbb{R}}

\newcommand{\bT}{\mathbb{T}}
\newcommand{\bU}{\mathbb{U}}
\newcommand{\bV}{\mathbb{V}}
\newcommand{\bW}{\mathbb{W}}

\newcommand{\bZ}{\mathbb{Z}}

\newcommand{\sW}{\mathscr{W}}

\makeatletter

\DeclareMathAlphabet{\mathrmsl}{OT1}{cmr}{m}{sl}            

\DeclareSymbolFont{script}{U}{eus}{m}{n}                    
  \DeclareSymbolFontAlphabet{\mathscr}{script}

\DeclareMathSymbol{\EuWedge}{0}{script}{"5E}                
\DeclareMathSymbol{\StdWedge}{2}{symbols}{"5E}              
  
\DeclareFontFamily{U}{mathx}{\hyphenchar\font45}            
\DeclareFontShape{U}{mathx}{m}{n}{<-> mathx10}{}
\DeclareSymbolFont{mathx}{U}{mathx}{m}{n}
\DeclareMathAccent{\widebar}{0}{mathx}{"73}

\newcommand{\mathopsl}[1]{                                  
  \operatorname{\mathrm{#1}} }

\newcommand{\ifbracs}[1]{                                   
  \ifx  \relax#1\relax                                        
  \else                 (#1)                                  
  \fi }
\newcommand{\ifbdot}[2][\bdot]{                             
  \ifx  \relax#2\relax  #1                                    
  \else                 #2                                    
  \fi }
\newcommand{\ifcdot}[1]{                                    
  \ifx  \relax#1\relax                                        
  \else                 \cdot #1                              
  \fi }
\newcommand{\ifscolon}[3][;]{                               
  \ifx  \relax#2\relax                                        
    \ifx  \relax#3\relax                                      
    \else                 #3
    \fi 
  \else
    \ifx  \relax#3\relax  #2
    \else                 #2 #1 #3
    \fi
  \fi }
\newcommand{\setbuild}[4]{                                  
  #1                                                          
  \ifx  \relax#4\relax  { #3 }                                
  \else                 { #3 \mid #4 }
  \fi
  #2 }
\newcommand{\Setbuild}[4]{                                  
  \left#1
  \ifx  \relax#4\relax  { #3 }
  \else                 { \left.#3\,\middle|\,#4\right. }
  \fi
  \right#2 }
\newcommand{\nullspacer}{\phantom{\hspace{0em}}}            
    
\newcommand{\bdot}[1][0.2]{                                 
  \hspace{#1em}\cdot\hspace{#1em} }
\newcommand{\acts}{\cdot}                                   
\renewcommand{\b}[1]{\bm\hat{#1}}                           
\newcommand{\defeq}{:=}                                     
\newcommand{\eqdef}{=:}                                     
\newcommand{\at}[1]{\vert_{#1}}                             

\newcommand{\colvect}[1]{                                   
  \arraycolsep=0.1em
  \begin{bmatrix} #1 \end{bmatrix} }
\newcommand{\colvectpunct}[2][-1.5em]{                      
  \raisebox{#1}{\!$#2$} }                                     
\newcommand{\inlinematrix}[1]{                              
  \left[ \begin{smallmatrix}
    #1
  \end{smallmatrix} \right] }
\newcommand{\by}{{\times}}                                  

\newcommand{\setof}[2]{\setbuild{\{}{\}}{#1}{#2}}           
\newcommand{\Setof}[2]{\Setbuild{\{}{\}}{#1}{#2}}           
\newcommand{\linspan}[3][]{                                 
  \setbuild{\langle}{\rangle}{#2}{#3}_{#1} }

\newcommand{\Hom}[3][]{\mathopsl{Hom}_{#1}(#2,#3)}          
\renewcommand{\dim}[1][]{\mathopsl{dim}_{#1}}               
\newcommand{\rk}[1]{\mathopsl{rk}#1}                        

\newcommand{\im}{\mathopsl{im}}                             

\newcommand{\isom}{\cong}                                   
\newcommand{\id}[1][]{\mathrm{id}_{#1}}                     
\renewcommand{\det}[1][]{\mathopsl{det}_{#1}}               
\newcommand{\tr}[1][]{\mathopsl{tr}_{#1}}                   
\newcommand{\pf}[1]{\mathopsl{pf}#1}                        
\newcommand{\adj}[1]{\mathopsl{adj}#1}                      
\newcommand{\norm}[1]{\|\ifbdot{#1}\|}                      
\newcommand{\intprod}{\lrcorner\hspace{0.6mm}}              
\newcommand{\alt}[1][]{\mathopsl{alt}_{#1}}                 
\newcommand{\sym}[1][]{\mathopsl{sym}_{#1}}                 
\newcommand{\trfree}{\circ}                                 
\newcommand{\conj}[1]{\overline{#1}}                        

\newcommand{\Union}[2]{{\textstyle\bigcup}_{#1}^{#2}}       
\newcommand{\union}[1][]{\cup_{#1}}                         
\newcommand{\Intsct}[2]{{\textstyle\bigcap}_{#1}^{#2}}      
\newcommand{\intsct}[1][]{\cap_{#1}}                        

\newcommand{\Dsum}[2]{\bigoplus_{#1}^{#2}}                  
\newcommand{\dsum}[1][]{\oplus_{#1}}                        
\newcommand{\Tens}[2][]{                                    
  \displaystyle\mathlarger{\otimes}
  _{#1}^{\ifbdot[\bullet]{#2}} }
\newcommand{\tens}[1][]{\otimes_{#1}}                       

\newcommand{\etens}[1][]{                                   
  \mathbin{\boxtimes_{#1}} }
\newcommand{\ltens}{\hspace{0.05em}}                        
\newcommand{\Symm}[2][]{                                    
  {\mathopsl{S}}
  _{#1}^{\ifbdot[\bullet]{#2}}}
\newcommand{\symm}[1][]{\odot_{#1}}                         
\renewcommand{\wedge}[1][]{\StdWedge_{#1}}                  
\newcommand{\Wedge}[2][]{                                   
  \displaystyle\EuWedge
  _{#1}^{\ifbdot[\bullet]{#2}} }
\newcommand{\Cartan}[2][]{                                  
  \displaystyle\mathlarger{\circledcirc}_{#1}                 
  \phantom{\hspace{0em}}^{\ifbdot[\bullet]{#2}} }             
\newcommand{\cartan}[1][]{\circledcirc_{#1}}                

\newcommand{\grp}[3][]{                                     
  \mathrm{#2}{#1} \ifbracs{#3} }
\newcommand{\alg}[3][]{                                     
  \mathfrak{#2}{#1} \ifbracs{#3} }
\newcommand{\Stab}[2][]{\mathopsl{Stab}_{#1}#2}             
\newcommand{\Norm}[2][]{N_{#1}(#2)}                         
\newcommand{\ad}{\mathopsl{ad}}                             
\newcommand{\Ad}{\mathopsl{Ad}}                             
\newcommand{\AD}{\underline{\mathopsl{Ad}}}                 
\newcommand{\gr}[2][]{\mathopsl{gr}_{#1}#2}                 
\newcommand{\grpcenter}[2][]{Z_{#1}(#2)}                    
\newcommand{\liecenter}[2][]{\fz_{#1}(#2)}                  
\newcommand{\sspart}[1]{#1_{\mathrm{ss}}}                   

\newcommand{\coroot}[1]{#1^{\vee}}                          

\newcommand{\cartanint}[2]{                                 
  \tfrac{ 2\killing{#1}{#2} }{ \killing{#2}{#2} } }
\newcommand{\height}[2][]{\mathopsl{ht}_{#1}(#2)}           
\newcommand{\lowest}[1][\fg]{\rho_{#1}}                     
\newcommand{\refl}[2]{\sigma_{#1}\ifbracs{#2}}              
\newcommand{\longest}[1][0]{w_{#1}}                         

\newcommand{\grpdot}[1][]{\cdot_{#1}}                       

\newcommand{\opp}[1]{\b{#1}}                                
\newcommand{\liechain}[3][\fp^{\perp}]{                     
  C_{\ifbdot[\bullet]{#2}} \ifbracs{\ifscolon{#1}{#3}} }
\newcommand{\liehom}[3][\fp^{\perp}]{                       
  H_{\ifbdot[\bullet]{#2}} \ifbracs{\ifscolon{#1}{#3}} }
\newcommand{\liecochain}[3][\fp^{\perp}]{                   
  C^{\ifbdot[\bullet]{#2}} \ifbracs{\ifscolon{#1}{#3}} }
\newcommand{\liecohom}[3][\fp^{\perp}]{                     
  H^{\ifbdot[\bullet]{#2}} \ifbracs{\ifscolon{#1}{#3}} }
\newcommand{\liebdy}[1][]{{\partial_{#1}}}                  
\newcommand{\liediff}[1][]{{\partial^*_{#1}}}               
\newcommand{\quab}[1][]{{\square_{#1}}}                     

\newcommand{\jmult}[1][]{\circ_{#1}}                        
\newcommand{\jespace}[3][\bJ]{#1_{#2} \ifbracs{#3}}         
\newcommand{\jtrip}[3]{                                     
  \{ \ifbdot{#1}, \ifbdot{#2}, \ifbdot{#3} \} }

\newcommand{\algbracadornment}{}                            
\newcommand{\bracbuild}[7][]{                               
  #2 \ifbdot{#5} #4 \ifbdot{#6} #3                            
    \ifx \relax#1\relax ^{\algbracadornment}                  
    \else               ^{#1}                                 
    \fi                                                       
  _{#7} \nullspacer }                                         

\newcommand{\liebrac}[2]{                                   
  \bracbuild[\hspace{0em}]{[}{]}{,}{#1}{#2}{} }

\newcommand{\liebracw}[2]{                                  
  \bracbuild[\hspace{0em}]{[}{]}{\wedge}{#1}{#2}{} }
\newcommand{\Killing}[2]{                                   
  \bracbuild[\hspace{0em}]{\left\langle}{\right\rangle}{,}{#1}{#2}{} }
\newcommand{\killing}[2]{                                   
  \bracbuild[\hspace{0em}]{\langle}{\rangle}{,}{#1}{#2}{} }

\newcommand{\poisson}[2]{                                   
  \bracbuild[\hspace{0em}]{\{}{\}}{,}{#1}{#2}{} }

\newcommand{\Algbrac}[3][]{                                 
  \bracbuild[#1]{\big\llbracket}{\big\rrbracket}{,}{#2}{#3}{} }
\newcommand{\algbrac}[3][]{                                 
  \bracbuild[#1]{\llbracket}{\rrbracket}{,}{#2}{#3}{} }

\newcommand{\algbracw}[4][]{                                
  \bracbuild[#1]{\llbracket}{\rrbracket}{\wedge}{#2}{#3}{#4} }

\newcommand{\algbraco}[4][]{                                
  \bracbuild[#1]{\llbracket}{\rrbracket}{\symm}{#2}{#3}{#4} }

\newcommand{\To}{\longrightarrow}                           
\newcommand{\From}{\longleftarrow}                          
\newcommand{\injto}{\hookrightarrow}                        
\newcommand{\surjto}{\twoheadrightarrow}                    
\newcommand{\Injto}{\lhook\joinrel\longrightarrow}          
\newcommand{\Surjto}{\:-\mathrel{\mkern-8mu}\surjto}        

\newcommand{\ses}[4][0]{                                    
  #1 \To #2 \Injto #3 \Surjto #4 \To #1 }                     

\newcommand{\s}[3][M]{                                      
  \Omega^{#2} \ifbracs{\ifscolon{#1}{#3}} }                   
\newcommand{\assocbdl}[3]{#1\times_{#2}#3}                  
\newcommand{\cpxbdl}[1]{\bC #1}                             
\newcommand{\jetbdl}[2]{\cJ^{#1}(#2)}                       
\renewcommand{\d}{\mathrm{d}}                               
\newcommand{\p}{\partial}                                   
\newcommand{\D}{\nabla}                                     

\newcommand{\Dspace}[1][]{                                  
  [\D]                                                        
  \ifx \relax#1\relax ^{\algbracadornment}                    
  \else               ^{#1}
  \fi }
\newcommand{\weyld}[2][]{\delta^{#1}_{#2}}                  
\newcommand{\bgg}[2][]{\cD^{\ifscolon[,]{#1}{#2}}}          
\newcommand{\bggrepr}[2][]{L^{\ifscolon[,]{#1}{#2}}}        
\newcommand{\bggproj}[2][]{\pi^{\ifscolon[,]{#1}{#2}}}      
\newcommand{\bggpi}[2][]{\Pi^{\ifscolon[,]{#1}{#2}}}        
\newcommand{\harm}{\circ}                                   

\newcommand{\vol}[1][]{\mathrm{vol}_{#1}}                   
\newcommand{\grad}[2][]{\mathopsl{grad}_{#1} #2}            
\newcommand{\lie}[1]{\mathcal{L}_{#1}}                      
\renewcommand{\div}[2][]{\mathopsl{div}_{#1}#2}             
\newcommand{\hodge}[1]{\operatorname{\ast}#1}               

\newcommand{\conf}[1][c]{\mathsf{#1}}                       

\newcommand{\Tor}[2][]{ T^{#1}_{#2} }                       
\newcommand{\Curv}[3][\D]{R^{#1}_{#2}\ifcdot{#3}}           
\newcommand{\Weyl}[3][]{W^{#1}_{#2}\ifcdot{#3}}             
\newcommand{\Nullset}[2][]{                                 
  \mathrm{Null}\ifbracs{#1}
  \ifx  \relax#2\relax 
  \else                \at{#2}
  \fi }
\newcommand{\sRic}[2][\D]{\mr{Ric}^{#1}\ifbracs{#2}}        
\newcommand{\Scal}[1][\D]{\mr{Scal}\ifbracs{#1}}            
\newcommand{\nRic}[3][\D]{r^{#1}_{#2}\ifbracs{#3}}          
\newcommand{\CY}[3][\D]{C^{#1}_{#2}\ifbracs{#3}}            
\newcommand{\Nijen}[2][J]{N^{#1}_{#2}}                      

\renewcommand{\Re}[1]{                                      
  \mathopsl{Re}\! \Setbuild{[}{]}{#1}{} }
\renewcommand{\Im}[1]{                                      
  \mathopsl{Im}\! \Setbuild{[}{]}{#1}{} }

\renewcommand{\max}[2][]{                                   
  \mathopsl{max}_{#1} \! \setof{#2}{} }
\newcommand{\ve}{\varepsilon}                               
\renewcommand{\tilde}[1]{\widetilde{#1}}                    
\newcommand{\nat}{\natural}                                 
\renewcommand{\setminus}{\smallsetminus}                    

\let\old@sum\sum                                            
\renewcommand{\sum}[2]{                                     
  {\textstyle\old@sum}_{#1}^{#2} }
\newcommand{\Sum}[2]{\old@sum_{#1}^{#2}}                    

\let\old@prod\prod                                          
\renewcommand{\prod}[2]{                                    
  {\textstyle\old@prod}_{#1}^{#2} }
\newcommand{\Prod}[2]{\old@prod_{#1}^{#2}}                  

\newcommand{\RP}[1][n]{\mathbb{RP}^{#1}}                    
\newcommand{\CP}[1][n]{\mathbb{CP}^{#1}}                    
\newcommand{\HP}[1][n]{\mathbb{HP}^{#1}}                    
\newcommand{\OP}[1][2]{\mathbb{OP}^{#1}}                    
\newcommand{\FP}[1][n]{\mathbb{FP}^{#1}}                    
\newcommand{\Grass}[2]{\grp[_{#1}]{Gr}{#2}}                 
\newcommand{\Sph}[1][n]{\mathbb{S}^{#1}}                    
\newcommand{\CSph}[1][n]{\cpxbdl{\Sph[#1]}}                 

\newcommand{\erepn}[1][\bC]{{#1}_{27}}                      
\newcommand{\eerepn}[1][\bC]{{#1}_{56}}                     
\newcommand{\spinrepn}[1][]{\$_{#1}}                        
\newcommand{\cpxrepn}[1]{#1_{\bC}}                          
\newcommand{\realrepn}[1]{#1_{\bR}}
\newcommand{\Realrepn}[2][1]{                               
  \left(#2\right)_{ \!\!\raisebox{#1ex}{$\scriptstyle\bR$}} }
\newcommand{\cpxconj}{\bC\mathtt{c}}                        
\newcommand{\pr}[2][]{\bP_{#1}(#2)}                         
\newcommand{\type}[2]{$\grp[_{#2}]{#1}{}$}                  

\newcommand{\erealform}[1]{\grp[_{\mr{#1}}]{E}{}}           

\newcommand{\cpx}[1]{{\mathbf #1}}                          
\newcommand{\qsum}[1][a]{\sum{#1=1}{3}}                     
\newcommand{\abs}[1]{\lvert #1 \rvert}                      

\makeatother

\makeatletter

\theoremstyle{plain}
  \newtheorem{thm}{Theorem}[chapter]      \newtheorem*{thm*}{Theorem}
  \newtheorem{prop}[thm]{Proposition}     \newtheorem*{prop*}{Proposition}
  \newtheorem{lem}[thm]{Lemma}            \newtheorem*{lem*}{Lemma}
  \newtheorem{cor}[thm]{Corollary}        \newtheorem*{cor*}{Corollary}

\theoremstyle{definition}
  \newtheorem{defn}[thm]{Definition}      \newtheorem*{defn*}{Definition}
  \newtheorem{conv}[thm]{Convention}      \newtheorem*{conv*}{Convention}
  \newtheorem{asmpt}[thm]{Assumption}     \newtheorem*{asmpt*}{Assumption}
  
\newtheoremstyle{italicremark}
  {\topsep}{\topsep}                      
  {\normalfont}{0pt}{\itshape}{.\ }{0pt}  
  {\thmname{#1} \thmnumber{#2}}           
\theoremstyle{italicremark}
  \newtheorem{rmk}[thm]{Remark}           \newtheorem*{rmk*}{Remark}
  \newtheorem{expl}[thm]{Example}         \newtheorem*{expl*}{Example}
  
\newenvironment{proofof}[1]
  { \begin{proof}[Proof of \pageref{#1} \ref{#1}] }
  { \end{proof} }
  
\newcommand{\noproof}{\qed}
  
\newenvironment{sketchproof}
  { \begin{proof}[Sketch proof] }
  { \end{proof} }
  
\newenvironment{enumeratepar}
  { \begin{enumerate*}[itemjoin=\smallskip\newline] }
  { \end{enumerate*} }
  
\newtheoremstyle{typeblock}
  {\topsep}{\topsep}
  {\normalfont}{0ex}{\bfseries}{.\ }{0ex}
  {#3}
\theoremstyle{typeblock}
  \newtheorem*{type@block}{}
  \newtheorem*{block}{}
\newcommand{\typeheader}[3][Type ]{#1$\bm{\mathrm{#2}_{#3}}$}
\newenvironment{typeblock}[3][Type ]  
  { \begin{type@block}[{\typeheader[#1]{#2}{#3}}] }
  { \end{type@block} }
\makeatother

\newcommand{\zbox}[2][c]{ \makebox[0em][#1]{$#2$} }

\newcommand{\itemref}[2]{\thref{#1}\ref{#1-#2}}

\newcommand{\proofref}[2]{\noindent \ref{#1-#2}}
\newcommand{\equivref}[4][\Rightarrow]{
  \noindent \ref{#2-#3}$#1$\ref{#2-#4}:}
  
\newcommand{\threfit}[1]{\pageref{#1} \ref{#1}}

\newcommand{\rmcite}[2][]{\emph{\cite[#1]{#2}}}
  
\newif\ifBufferDynkinGlobal
\newif\ifBufferDynkinLocal

\makeatletter
\newcommand*{\test@emph}{it}
\newcommand{\forcenoemph}[1]{
  \ifx  \f@shape \test@emph  \expandafter \emph{#1}%
  \else                      \expandafter #1%
  \fi }

\makeatother

\setlist[enumerate]{
  label=(\arabic*),
  font=\normalfont,
  itemsep=-0.1em,
  topsep=0.4em }
\setlist[itemize]{
  itemsep=-0.1em,
  topsep=0.4em }
  
\newenvironment{slimitemize}{ \begin{itemize}[leftmargin=1.4em] }
                            { \end{itemize} }

\renewcommand{\thetable}{\arabic{chapter}.\arabic{table}}

\newcommand{\listofcaptioned}{
  \listoffigures
  \vspace{2em}
  \bgroup
  \let\clearpage\relax
  \listoftables
  \egroup
}

\makeatletter
\let\old@thebibliography\thebibliography
\renewcommand{\thebibliography}[1]{
  \old@thebibliography{#1}
  \setlength{\itemsep}{0.0em}
  \small
}
\makeatother

\usepackage{fix-cm}

\renewcommand{\chaptername}{chapter}
\titleformat{\chapter}[display]
            {\flushright}
            {\raisebox{2.5em}{\textsc{\Large \chaptername}}
             {\fontsize{60}{0}\usefont{OT1}{cmr}{m}{n}\selectfont \thechapter}}
            {-0.5em}
            {\Huge\bfseries #1}
            
\titlespacing*{\chapter}{0pt}{0in}{40pt}

\newcommand*{\missingreference}{\color{red}{\bfseries??}}
\newcommand*{\missingcitation}{\color{green}{\bfseries??}}
\makeatletter
\def\@setref#1#2#3{%
  \ifx#1\relax
   \protect\G@refundefinedtrue
   \nfss@text{\reset@font\missingreference}%
   \@latex@warning{Reference `#3' on page \thepage \space
             undefined}%
  \else
   \expandafter#2#1\null
  \fi}
\def\@citex[#1]#2{\leavevmode
  \let\@citea\@empty
  \@cite{\@for\@citeb:=#2\do
    {\@citea\def\@citea{,\penalty\@m\ }%
     \edef\@citeb{\expandafter\@firstofone\@citeb\@empty}%
     \if@filesw\immediate\write\@auxout{\string\citation{\@citeb}}\fi
     \@ifundefined{b@\@citeb}{\hbox{\reset@font\missingcitation}%
       \G@refundefinedtrue
       \@latex@warning
         {Citation `\@citeb' on page \thepage \space undefined}}%
       {\@cite@ofmt{\csname b@\@citeb\endcsname}}}}{#1}}
\makeatother

\newcommand{\riem}{riemannian}               \newcommand{\Riem}{Riemannian}
\newcommand{\kahler}{K\"{a}hler}             \newcommand{\Kahler}{K\"{a}hler}
\newcommand{\qk}{quaternion-\kahler}         \newcommand{\Qk}{Quaternion-\kahler}
  \newcommand{\QK}{Quaternion-\Kahler}
\newcommand{\hk}{hyperk{\"a}hler}            
\newcommand{\LC}{Levi-Civita}
\newcommand{\Rspace}{R-space}                \newcommand{\Rspaces}{\Rspace s}

\newcommand{\ppg}{projective parabolic geometry}
  
  \newcommand{\PPG}{Projective Parabolic Geometry}
\newcommand{\ppgs}{projective parabolic geometries}
  \newcommand{\Ppgs}{Projective parabolic geometries}

\newcommand{\proj}{projective}               \newcommand{\Proj}{Projective}
\newcommand{\cproj}{c-projective}            \newcommand{\Cproj}{C-projective}
\newcommand{\qtn}{quaternion}                \newcommand{\Qtn}{Quaternion}
\newcommand{\qproj}{q-projective}            
\newcommand{\oct}{octonion}                  \newcommand{\Oct}{Octonion}

\newcommand{\formulae}{formulae}            
\newcommand{\torsionfree}{torsion-free}      
\newcommand{\tracefree}{trace-free}          
\newcommand{\einstein}{Einstein}             

\newcommand{\hamiltonian}{hamiltonian}       \newcommand{\Hamiltonian}{Hamiltonian}
\newcommand{\grassmannian}{grassmannian}     \newcommand{\Grassmannian}{Grassmannian}

\newcommand{\non}{non-}                      \newcommand{\Non}{Non-}
\newcommand{\wrt}{with respect to}           \newcommand{\Wrt}{With respect to}
\newcommand{\self}{self-}                    \newcommand{\Self}{Self-}
\newcommand{\pseudo}{pseudo-}                \newcommand{\Pseudo}{Pseudo-}
\newcommand{\co}{co-}                        
\newcommand{\lhs}{left-hand side}
\newcommand{\rhs}{right-hand side}

\newcommand{\sltriple}{$\alg[_2]{sl}{}$-triple}

\newcommand{\latin}[1]{\emph{#1}}
\newcommand{\ie}{\latin{i.e.}}
\newcommand{\cf}{\latin{cf}.}

\newcommand{\etc}[1][.]{\latin{etc}#1}
\newcommand{\etal}[1][.]{\latin{et al}#1}

\newcommand{\apriori}{\latin{a priori}}

\newcommand{\smalltitleword}[1]{\textit{\Large #1}}
\degree{Doctor of Philosophy}
\author{George Edward Frost}
\title{\smalltitleword{The} \\[-0.2em]
       \PPG\ \\[-0.4em]
       \smalltitleword{of} \\[-0.2em]
       \Riem, \Kahler\ \smalltitleword{and} \\
       \QK\ Metrics}
\department{Department of Mathematical Sciences}
\degreemonthyear{January 2016}
\norestrictions

\begin{document}

\maketitle
\label{c:abstract}

\BufferDynkinLocaltrue
\renewcommand{\dynkinnameoffset}{-0.75}

\begin{abstract}

We present a uniform framework generalising and extending the classical theories of \proj\ differential geometry, \cproj\ geometry, and almost \qtn ic\ geometry.  Such geometries, which we call \emph{\ppgs}, are abelian parabolic geometries whose flat model is an \Rspace\ $G\acts\fp$ in the infinitesimal isotropy representation $\bW$ of a larger \self dual\ symmetric \Rspace\ $H\acts\fq$.  We also give a classification of \ppgs\ with $H\acts\fq$ irreducible which, in addition to the aforementioned classical geometries, includes a geometry modelled on the Cayley plane $\OP[2]$ and conformal geometries of various signatures.

The larger \Rspace\ $H\acts\fq$ severely restricts the Lie-algebraic structure of a \ppg.  In particular, by exploiting a Jordan algebra structure on $\bW$, we obtain a $\bZ^2$-grading on the Lie algebra of $H$ in which we have tight control over Lie brackets between various summands.  This allows us to generalise known results from the classical theories.  For example, which \riem\ metrics are compatible with the underlying geometry is controlled by the first BGG operator associated to $\bW$.

In the final chapter, we describe \ppgs\ admitting a $2$-dimensional family of compatible metrics.  This is the usual setting for the classical \proj\ structures; we find that many results which hold in these settings carry over with little to no changes in the general case.

\end{abstract}

\tableofcontents
\listofcaptioned

\chapter{Introduction} 
\label{c:intro}

\BufferDynkinLocaltrue
\renewcommand{\dynkinnameoffset}{-0.75}

Given a \riem\ manifold $(M,g)$, it is natural to ask whether $M$ admits any other metrics $\b{g}$ with the same geodesics as $g$, viewed as unparameterised curves.  Such $\b{g}$ are said to be \emph{\proj ly\ equivalent} to $g$, leading to the notion of a \proj\ equivalence class of metrics.  This is the classical formulation of \proj\ differential geometry, as studied by authors such as Beltrami \cite{b1865-projequiv}, Dini \cite{d1869-2dlocalform}, Painlev{\'e} \cite{p1897-quadintegrals}, \LC\ \cite{l2009-dyneqns} and {\'E}.\ Cartan \cite{c1924-projconns, c1926-confproj}, to name just a few.

More properly, geodesics are a feature of linear connections rather than metrics.  Discarding the metrics $g, \b{g}$ then allows us to talk about a \emph{\proj\ equivalence class} $\Dspace[r]$ of linear connections $\D$ on $TM$, leading to the notion of a \emph{\proj\ manifold} $(M, \Dspace[r])$.  By results of Cartan \cite{c1924-projaffine, c1924-projconns, c1926-confproj} and Thomas \cite{t1925-projtractor}, there is an equivalence of categories between \proj\ structures on $M$ and \emph{Cartan connections} on the frame bundle of $M$.  Thus, in modern language, \proj\ differential geometry is an example of an \emph{abelian parabolic geometry} modelled on the \proj\ space $\RP[n]$, which is naturally a \proj\ manifold when equipped with the \proj\ equivalence class of its spherical metric.  This gives leverage to the theory of parabolic subalgebras and their representation theory: for example, $\Dspace[r]$ is identified with the space of \emph{Weyl connections}, which are induced by splitting the parabolic filtration of the Lie algebra $\alg{sl}{n+1,\bR}$ of \tracefree\ $(n+1) \by (n+1)$ real matrices.

There is now no guarantee that $\Dspace[r]$ contains the \LC\ connection of a metric.  In the classical picture with a background metric $g$, Sinjukov \cite{s1979-geodmappings} found that a second metric $\b{g}$ is \proj ly\ equivalent to $g$ if and only if an endomorphism $A(g,\b{g})$, constructed solely from $g, \b{g}$, satisfies a certain first-order differential equation.  The $2$-dimensional version of this equation was essentially known to Liouville \cite{l1889-2dmaineqn}; see also \cite{m2012-geodequivgr}.  An invariant version of Sinjukov's equation was later discovered by Eastwood and Matveev \cite{em2008-projmetrics}, which controls whether the \LC\ connection of a given metric lies in $\Dspace[r]$.  In parabolic language, this invariant equation coincides with the first BGG operator \cite{cd2001-curvedbgg, css2001-bgg} associated to the natural representation $\bW \defeq \Symm{2}\bR^{n+1}$ of $\alg{sl}{n+1,\bR}$.  The important point is that metrisability of $\Dspace[r]$ is controlled by a \proj ly\ invariant first-order linear differential equation with a representation-theoretic origin.

Moving to the holomorphic category, {\=O}tsuki and Tashiro \cite{ot1954-hgeodesics} found that \kahler\ metrics $g,\b{g}$ are \proj ly\ equivalent if and only if they are affinely equivalent, rendering \proj\ equivalence uninteresting in this context.  \Cproj\ geometry arises as the natural adaptation of \proj\ differential geometry to an almost complex manifold $(M,J)$: a curve $\gamma$ is a \emph{c-geodesic} of a connection $\D$ if and only if $\D_X X$ lies in the linear span $\linspan{X,JX}{}$ for all vectors $X$ tangent to $\gamma$, leading to a notion of \emph{\cproj\ equivalence}.  An \emph{(almost) \cproj\ structure} is then the choice of a \cproj\ equivalence class $\Dspace[c]$ of linear connections on $TM$.  The classical theory proceeds in much the same way as for \proj\ structures, and many results known for \proj\ structures were adapted to the \cproj\ setting; see \cite{dm1978-hproj, m1998-holomorphic, t1957-hproj}.  In particular, Domashev and Mike{\v s} \cite{dm1978-hproj} found that two \kahler\ metrics have the same c-geodesics if and only if a particular endomorphism $A(g,\b{g})$ satisfies a first-order linear differential equation similar to Sinjukov's equation.  There is also an interpretation in terms of \emph{hamiltonian $2$-forms}, as described by Apostolov \etal\ \cite{acg2003-ham2forms0, acg2006-ham2forms1, acgt2004-ham2forms2, acgt2008-ham2forms4, acgt2008-ham2forms3}.

Complex \proj\ space $\CP[n]$ is naturally a \cproj\ manifold when equipped with the \cproj\ equivalence class of its Fubini--Study metric.  Via the general theory of parabolic geometries \cite{cs2009-parabolic1}, we obtain an equivalence of categories between almost \cproj\ structures on $M$ and parabolic geometries modelled on $\CP[n]$ \cite{cemn2015-cproj, h2009-cproj}.  This again opens the door to methods from parabolic geometries and the BGG machinery, and one finds that metrisability of $\Dspace[c]$ is controlled by the first BGG operator associated to the real representation $\bW \defeq \realrepn{( \bC^{n+1} \etens \conj{\bC^{n+1}} )}$ of $\alg{sl}{n+1,\bC}$.  Using this language, the recent survey \cite{cemn2015-cproj} has obtained many results which mirror known results in \proj\ differential geometry.

The classical theory of almost \qtn ic\ manifolds can also be made to fit into this \proj\ picture.  An \emph{almost \qtn ic\ structure} on a manifold $M$ is a rank three subbundle $\cQ \leq \alg{gl}{TM}$ which is pointwise isomorphic to the unit \qtn s\ $\alg{sp}{1}$; see \cite{am1993-qtnlikestrs, s1982-qkmfds, s1986-qkdg}.  A connection is \emph{(almost) \qtn ic} if it preserves $\cQ$; it turns out \cite{am1996-qtnsubord} that the \qtn ic\ connections form an affine space modelled on $T^*M$, leading to a notion of \emph{\qtn ic\ equivalence}.  On the other hand, a curve $\gamma$ is called a \emph{q-geodesic} of $\D$ if $\D_X X$ lies in the \qtn ic\ span of $X$.  Fujimura \cite{f1976-Qconns} proved that connections have the same q-geodesics if and only if are \qtn ically\ equivalent, thus fitting \qtn ic\ geometry into the same framework as \proj\ and \cproj\ geometries.  We will see later that compatible metrics (\ie\ \qk\ metrics) are controlled by a first-order linear differential equation resembling Sinjukov's equation.

Salamon \cite{s1986-qkdg} originally described \qtn ic\ manifolds as manifolds modelled locally on the \qtn ic\ \proj\ space $\HP[n]$, which lends itself to a parabolic description.  The general theory gives an equivalence of categories between almost \qtn ic\ structures on $M$ and parabolic geometries modelled on $\HP[n]$, with \qtn ic\ connections corresponding to Weyl connections.  Metrisability is controlled by a first BGG operator, now associated to the representation $\bW \defeq \realrepn{( \Wedge[\bC]{2} \bC^{2n+2} )}$ of $\alg{sl}{n+1,\bH}$.

These three classical theories evidently have similar descriptions: all are abelian parabolic geometries modelled on a \proj\ space $\FP[n]$, with a well-defined metrisability problem controlled the first BGG operator associated to a representation $\bW$.  Moreover, many results in the three theories have proofs which differ only in places where the base field $\bF$ has influence.  The objectives of this thesis are as follows:
\begin{enumerate}
  \item \label{item:intro-goals-ppg}
  Construct a general framework in which the classical \proj\ structures may be described as special cases, using the language of abelian parabolic geometries;
    
  \item \label{item:intro-goals-W}
  Give a general interpretation for the representation $\bW$ to which the metrisability problem is associated, and develop algebraic tools for other key representations;

  \item \label{item:intro-goals-metric}
  Interpret solutions of the first BGG operator associated to the representation $\bW$ as (\pseudo\riem) metrics compatible with the underlying geometric structure;
  
  \item \label{item:intro-goals-adapt}
  Generalise results which have similar statements in the three classical cases and adapt them to the general framework;
  
  \item \label{item:intro-goals-pencil}
  Describe $2$-dimensional families of compatible metrics, in order to make contact with the classical approaches to the \proj\ structures.
\end{enumerate}

The first key observation to achieve \ref{item:intro-goals-ppg} is as follows.  For \proj\ differential geometry we have $\fg \defeq \alg{sl}{n+1,\bR}$ and $\bW \defeq \Symm{2}\bR^{n+1}$, with $\fp$ given by crossing the final node of the Satake diagram.  We may identify $\RP[n]$ with a generalised flag manifold $G/P$, where $P$ is the adjoint stabiliser of a lowest weight vector in $\bV^*$ for any irreducible $\fg$-representation $\bV$ whose highest weight is supported on the final node of the Satake diagram of $\fg$ \cite[Prop.\ 3.2.5]{cs2009-parabolic1}.  We thus obtain a \proj\ embedding $\RP[n] \injto \pr{\bV^*}$ for any such representation; since these representations are of the form $\bV_k \defeq \Symm{k}\bR^{n+1*}$ for $k>0$, we recover the Veronese embeddings $\RP[n] \injto \pr{\Symm{k}\bR^{n+1}}$.  The case $k=2$ is our previous representation $\bW$, thus giving an embedding into the projectivisation of the representation which sets up the metrisability problem.

Secondly, we notice that $\bR^{n+1} \dsum \bR^{n+1*}$ is naturally a symplectic vector space when equipped with the symplectic structure $\omega((u,\alpha), (v,\beta)) \defeq \beta(u) - \alpha(v)$.  The adjoint representation of $\fh \defeq \alg{sp}{\bR^{n+1} \dsum \bR^{n+1*}, \omega} \isom \alg{sp}{2n+2,\bR}$ then decomposes as
\begin{equation*}
  \fh = \Symm{2}(\bR^{n+1} \dsum \bR^{n+1*})
    \isom \Symm{2}\bR^{n+1} \dsum \alg{gl}{n+1,\bR} \dsum \Symm{2}\bR^{n+1*},
\end{equation*}
where $\alg{gl}{n+1,\bR}$ is the reductive Lie algebra of endomorphisms of $\fh$ which preserve the block decomposition of $\fh$.  Moreover $\alg{gl}{n+1,\bR} = \fg \dsum \bR$, so we may recover our original algebra $\fg = \alg{sl}{n+1,\bR}$ from $\fh$ as the semisimple part of $\alg{gl}{n+1,\bR}$.  Thus $\bW$ may be regarded as the infinitesimal isotropy representation $\fh/\fq$ of the abelian parabolic subalgebra $\fq \defeq (\fg \dsum \bR) \dsum \bW^*$.  This supplies a larger symmetric \Rspace\ $H\acts\fq$ (see \thref{defn:lie-para-rspace}), isomorphic to the space of lagrangian subspaces of $\bR^{2n+2}$, which contains $G\acts\fp \isom \RP[n]$.

These observations equally apply to the flat models $\CP[n]$ and $\HP[n]$ of \cproj\ and \qtn ic\ geometry.  For $\CP[n]$ we have $\fg = \alg{sl}{n+1,\bC}$ and $\bW = \realrepn{( \bC^{n+1} \etens \conj{\bC^{n+1}} )}$, and we find a \proj\ embedding $\CP[n] \injto \pr{\bW}$ as before.  The $\fg$-representation $\fh \defeq \bW \dsum (\fg\dsum\bR) \dsum \bW^*$ also has a graded Lie algebra structure, now isomorphic to $\alg{su}{n+1,n+1}$, yielding a larger symmetric \Rspace\ $H\acts\fq$ of  maximal isotropic subspaces of $\bC^{2n+2}$ of a hermitian inner product of signature $(n+1,n+1)$.  For $\HP[n]$ we have $\fg = \alg{sl}{n+1,\bH}$ and $\bW = \realrepn{( \Wedge[\bC]{2}\bC^{2n+2} )}$, with a Pl{\"u}cker-type embedding $\HP[n] \injto \pr{\bW}$.  There is again a graded Lie algebra structure on $\fh \defeq \bW \dsum (\fg\dsum\bR) \dsum \bW^*$, this time isomorphic to the real form $\alg[^*]{so}{4n+4}$ of $\alg{sl}{2n+2,\bC}$, yielding a larger symmetric \Rspace\ $H\acts\fq$ of isotropic \qtn ic\ subspaces of $\bH^{n+1}$.

The final necessary observation is that all three \Rspaces\ $H\acts\fq$ are \emph{\self dual} in the sense of \cite{bdpp2011-rspaces}; for abelian parabolic subalgebras $\fq, \opp{\fq} \leq \fh$ which satisfy $\fq \dsum \opp{\fq}^{\perp} = \fh$, this amounts to asking that $\fq, \opp{\fq}$ are conjugate by an element of $H$.  \Self duality\ provides the final ingredient in our general definition.

\begin{defn*} Let $H\acts\fq$ be a \self dual\ symmetric \Rspace\ with infinitesimal isotropy representation $\bW \defeq \fh/\fq$.  A \emph{\ppg} is a parabolic geometry modelled on the \Rspace\ $G\acts\fp$ given by the stabiliser of a lowest weight orbit in $\bW$. \end{defn*}

We will describe how to recover $G\acts\fp$ from $H\acts\fq$ later.  This ``top-down'' approach is beneficial because $H\acts\fq$ contains the Lie-algebraic structure of both $\fg$ and its representation $\bW$.  Note however that the \Rspace\ $G\acts\fp$ is not \apriori\ symmetric, as we find for the classical \proj\ structures, and some work involving the root data of $\fh$ and $\fg$ is required to show this.  The key step is relating \self duality\ of $H\acts\fq$ to a Jordan algebra structure on $\bW$, which is a commutative but \non associative\ algebra satisfying a power associativity relation.  An idempotent decomposition of $\bW$ then gives detailed information about the graded components of both $\fg$ and $\bW$, in particular allowing us to construct a $\bZ^2$-grading on $\fh$ \wrt\ algebraic Weyl structures $\fq$ and $\fp$.

\begin{thm*} With notation as above, $G\acts\fp$ is a symmetric \Rspace.  Moreover $\bW$ admits the structure of a Jordan algebra and decomposes into three graded components as a $\fp$-representation.  Thus $\fh \isom \bW \dsum (\fg \dsum \bR) \dsum \bW^*$ admits a $\bZ^2$-grading, with Lie brackets between the various summands given by Table \ref{tbl:ppg-alg-Z2gr}. \end{thm*}

It is not such a surprise that $\bW$ admits a Jordan algebra structure: the relation between \self dual\ symmetric \Rspaces\ and Jordan algebras has been extensively studied by Tits \cite{t1962-kktalg}, Koecher \cite{k1967-imbedding1, k1968-imbedding2}, Meyberg \cite{m1970-jordantriple} and Bertram \cite{b2002-projgeom, b2011-jordannonassoc, b2013-jordanlie}.  In outline, any Jordan algebra $(\bW,\jmult)$ can be embedded into the \emph{Kantor--Koecher--Tits algebra} $\fh \defeq \bW \dsum \alg{der}{\bW} \dsum \bW^*$, where $\alg{der}{\bW}$ is the Lie algebra generated by the multiplication maps $L_{x \jmult y} : z \mapsto (x \jmult y) \jmult z$.  Then $\fq \defeq \alg{der}{\bW} \dsum \bW^*$ becomes an abelian parabolic subalgebra of $\fh$, and the corresponding \Rspace\ can be shown to be \self dual \cite{m1970-jordantriple}.  Thus we obtain a $1$-to-$1$ correspondence between Jordan algebras and \self dual\ symmetric \Rspaces; Loos \cite{l1971-jrspace} extends this to a $1$-to-$1$ correspondence between so-called \emph{Jordan triple systems} and symmetric \Rspaces.

A Jordan algebra is called \emph{formally real} if the trace form $\tau(x,y) \defeq \tr(L_{x \jmult y})$ is positive definite.  The classification of formally real Jordan algebras was obtained by Jordan, von Neumann and Wigner \cite{jnw1934-jordanqm}, who showed that they comprise four infinite families and a single exceptional algebra.  The four families consist of the symmetric real matrices, the hermitian matrices, the \qtn-hermitian matrices, and the \emph{spin factors}, which may be described as a Clifford algebra.  The exceptional \emph{Albert algebra} consists of $3 \by 3$ \oct-hermitian matrices and was described by Albert \cite{a1934-albertalg}.

Of interest to us is the intimate relationship between formally real Jordan algebras and projective geometry \cite{b2002-octonions, b2002-projgeom}: the space $\mathopsl{idem}(\bW)$ of primitive idempotents may be stratified by their trace, and the idempotents with trace one may be viewed as points in a \proj\ space.  The Jordan algebras of symmetric, hermitian and \qtn-hermitian\ matrices yield the classical \proj\ spaces $\RP$, $\CP$ and $\HP$ \cite{m1978-jordan}, while the spin factors allow us to view the conformal sphere $\Sph[n]$ as ``one dimensional \proj\ geometry over $\bR^n$'' \cite{l2002-riemnormeddiv, ll2010-specialholonomy}.  For the Albert algebra, one obtains Moufang's \oct ic\ \proj\ plane $\OP[2]$ \cite{b2002-octonions, m1933-cayleyplane}.
Moreover by work of Hirzebruch \cite{h1965-riemsymm}, the trace form induces a \riem\ metric on $\mathopsl{idem}(\bW)$, which may then is a $2$-point homogeneous \riem\ symmetric space of rank one.  Thus we recover Cartan's classification \cite{c1926-riemsymm1, c1927-riemsymm2} of rank one \riem\ symmetric spaces.  There is also a complex analytic description in terms of bounded symmetric domains of tube type; see \cite{fkklr2012-symmdomains, k1970-jordandg, l1971-jrspace}.

Returning to the case of a \ppg, it turns out that the Jordan algebra structure of $\bW$ strictly confines the Lie brackets between the various summands of the $\bZ^2$-grading of $\fh$, allowing us to calculate many brackets independently of the \ppg\ in question.  This often means that, once one has reduced a result to the verification of an algebraic identity, the proof may be completed via a series of formal manipulations using the Killing form and Jacobi identity.  Although the reader may express dissatisfaction at these (often long) algebraic manipulations, the author would argue that the strictness of the algebraic framework deftly explains why many results in the classical theories have only subtly different proofs.  In particular, our fourth goal above may be achieved by framing a desired result in purely algebraic terms, where it can be solved by algebraic manipulations.

\begin{thm*} Consider a \ppg\ on $M$ with flat model $G\acts\fp$ and infinitesimal isotropy representation $\bW$.  Then solutions of the first BGG operator associated to $\bW$ induce metrics compatible with the underlying geometric structure.  In particular, \einstein\ metrics correspond to so-called normal BGG solutions. \end{thm*}

Using the classification of \self dual\ symmetric spaces, it is straightforward to classify the \ppgs\ with $H\acts\fq$ irreducible.  This classification may be phrased in terms of a pair of integers $(r,n)$, which arise from the idempotent decomposition of the Jordan algebra $\bW$ and play an important role in the algebraic theory.  In addition to the classical \proj\ structures coming from the formally real Jordan algebras $\bW$, the classification includes geometries modelled on the \grassmannian\ of $2$-planes, a symmetric \Rspace\ associated to the split real form of $\alg[_6]{e}{\bC}$, and conformal geometries of various signatures.

\vfill
\begin{block}[Overview] The first two chapters provide relevant background material.  Chapter \ref{c:lie} introduces parabolic subalgebras and \Rspaces, focusing primarily on their structure theory and Lie algebra homology.  Chapter \ref{c:para} reviews the theory of parabolic geometries, including their tractor calculus and the curved BGG machinery.

Next we study the three classical \proj\ structures.  Chapter \ref{c:proj} gives a detailed introduction to \proj\ differential geometry, both from the classical and parabolic perspectives.  Chapters \ref{c:cproj} and \ref{c:qtn} describe \cproj\ geometry and almost \qtn ic\ geometry in a similar way.  Hopefully the reader will excuse some repetition: these chapters serve primarily as a literature review and as general motivation, although there are some (apparently) original results in the \qtn ic\ case.

The general framework uniting these geometries is defined and studied in Chapter \ref{c:ppg}.  In particular, we undertake a detailed investigation of the algebraic structure of a \ppg.  Afterwards we examine the metrisability problem, obtaining results similar to the classical cases.

In light of goal \ref{item:intro-goals-pencil} above, Chapter \ref{c:mob2} is devoted to the study of \ppgs\ admitting a $2$-dimensional family of compatible metrics.  Notably, we obtain results on the geodesic flow of a metric, and a family of commuting vector fields.

Finally, the two appendices contain supplementary material.  Appendix \ref{c:app-alg} contains some algebraic identities, whose proofs are long and would have disrupted the flow of the text.  Appendix \ref{c:app-tbl} summaries some representation-theoretic data relating to \ppgs\ for the reader's convenience. \end{block}

\begin{block}[Notation] While notation should not pose a significant problem, a few points are worth mentioning.  For valence one tensors $\alpha, \beta$, our conventions for the wedge and symmetric product are $\alpha \wedge \beta \defeq \alpha \tens \beta - \beta \tens \alpha$ and $\alpha \symm \beta \defeq \tfrac{1}{2}( \alpha \tens \beta + \beta \tens \alpha )$.  We denote the external tensor product by $\etens$, and the Cartan product by $\cartan$.  Hamiltonian's \qtn s\ and Cayley's \oct s\ are denoted by $\bH$ and $\bO$ respectively.  Finally, unless stated otherwise, all differential geometric objects are smooth.  An index of notation is provided on page \pageref{c:app-symb}. \end{block}

\begin{block}[Acknowledgements] There are several people I must thank, without whom this project would not have been completed: my supervisor, Prof.\ David Calderbank, for his patience and constant suggestion of useful ideas; my wife, Loui, for her unwavering support and understanding of my strange working hours; and my parents, for their encouragement and reinforcing the importance of education.  

I am also grateful for financial support from the EPSRC, and to my examiners for many helpful suggestions. \end{block}

\chapter[Background from Lie theory]{Background from \\ Lie theory} 
\label{c:lie}

\BufferDynkinLocaltrue
\renewcommand{\dynkinnameoffset}{-0.75}

We begin with a review of some necessary ingredients from Lie theory.  We shall assume that the reader has a working knowledge of the structure theory and representation theory of semisimple Lie algebras; see \cite{cs2009-parabolic1,fh1991-repntheory,h1972-repntheory} for a thorough introduction.

The structure theory and representation theory of parabolic subalgebras of semisimple Lie algebras shall feature heavily throughout this thesis, so we spend some time describing the pertinent results in Section \ref{s:lie-para}.  In particular, the theory of \Rspaces\ and their \proj\ embeddings shall be important.

As we shall see in Section \ref{s:para-bgg}, the invariant differential operators associated to a parabolic geometry are related to the Lie algebra homology of its flat model.  We introduce Lie algebra homology and cohomology in Section \ref{s:lie-hom}, as well as describing the algorithm for its computation provided by Kostant's version of the Bott--Borel--Weil theorem.  The standard reference for this material is \cite{k1961-liehom}, but readable accounts may also be found in \cite{be1989-penrose,cs2009-parabolic1,h2008-bggcat}.

\vspace{0.2em}
\section{Parabolic subalgebras} 
\label{s:lie-para}

Let $\fg$ be a complex semisimple Lie algebra.  A maximal solvable subalgebra $\fb\leq\fg$ is called a \emph{Borel subalgebra}, while a subalgebra $\fp\leq\fg$ is called \emph{parabolic} if $\fp$ contains a Borel subalgebra.  However, the following equivalent definition is available; see \cite{bdpp2011-rspaces,cds2005-ricci}.

\begin{defn} \thlabel{defn:lie-para-para} A subalgebra $\fp$ of a semisimple Lie algebra $\fg$ is \emph{parabolic} if the Killing polar $\fp^{\perp}$ is a nilpotent subalgebra of $\fg$.  We say $\fp$ is an \emph{abelian parabolic} if $\fp^{\perp}\leq\fg$ is an abelian subalgebra. \end{defn}

\begin{lem} \emph{\cite[Thm.\ 1]{b1975-liealg}}  \thlabel{lem:lie-para-pperp} Let $\fp\leq\fg$ be parabolic.  Then $\fp^{\perp}$ is a nilpotent ideal of $\fp$. \end{lem}

\begin{proof} By invariance of the Killing form, we have $\killing{ \liebrac{\fp^{\perp}}{\fp} }{ \fp } = \killing{ \fp^{\perp} }{ \liebrac{\fp}{\fp} } \subseteq \killing{ \fp^{\perp} }{ \fp } = 0$ and hence $\liebrac{\fp}{\fp^{\perp}} \subseteq \fp^{\perp}$. \end{proof}

In fact, one can show \cite[Lem.\ 1.1]{d1976-lienilrad} that $\fp^{\perp}$ coincides with the nilpotent radical of $\fp$.  Then the quotient $\fp^0 \defeq \fp/\fp^{\perp}$ is a reductive Lie algebra, called the \emph{reductive Levi factor} of $\fp$.  It is always possible to choose a splitting of the projection $\fp\surjto\fp^0$, so we may identify $\fp^0$ with a subalgebra complementary to $\fp^{\perp}$ in $\fp$ such that $\fp \isom \fp^0\ltimes\fp^{\perp}$ is a semi-direct sum \cite{cds2005-ricci}.

The benefit of \thref{defn:lie-para-para} is that it works over any field of characteristic zero, whereas the definition via Borel subalgebras only works over $\bC$.  One may also use this idea to define parabolic subalgebras $\fp$ of a reductive Lie algebra $\fg$ \cite{n2014-coxconfig}, by asking that $\fp^{\perp}$ is a nilpotent subalgebra of $\fp \intsct \liebrac{\fg}{\fg}$.  However, we shall restrict attention to semisimple Lie algebras $\fg$, where $\fg = \liebrac{\fg}{\fg}$.

In Subsection \ref{ss:lie-para-std} we shall develop the basic structural theory of parabolic subalgebras from the root data of $\fg$.  This suggests a filtration associated to any parabolic subalgebra, which we discuss in Subsection \ref{ss:lie-para-gr}.  The representation theory of a parabolic subalgebra is discussed in Subsection \ref{ss:lie-para-repn}.  Finally, we discuss conjugacy classes of parabolics in Subsection \ref{ss:lie-para-rspace}, which forms the basis of central definitions and results in this thesis.

\subsection{Standard parabolics} 
\label{ss:lie-para-std}

Let $\fg$ be a complex semisimple Lie algebra, which we fix henceforth.  By appealing to the structure theory of $\fg$, we may identify a family of so-called \emph{standard parabolics}.  Choose a Cartan subalgebra $\ft\leq\fg$ with roots $\Delta \subset \ft^*$, and choose a positive subsystem $\Delta^{+}\subset\Delta$ with simple roots $\Delta^0$.  From the root space decomposition $\fg = \ft \dsum \Dsum{\alpha\in\Delta}{} \fg_{\alpha}$ of $\fg$, it is easy to see that $\fb \defeq \ft \dsum \Dsum{\alpha\in\Delta^{+}}{} \fg_{-\alpha}$ is a Borel subalgebra of $\fg$, called the \emph{standard Borel} \wrt\ $\ft$ and $\Delta^{+}$.

\begin{defn} \thlabel{defn:lie-para-std} A parabolic subalgebra $\fp\leq\fg$ is called a \emph{standard parabolic} \wrt\ $\ft$ and $\Delta^{+}$ if it contains the standard Borel.%
\footnote{Many authors, notably \cite{be1989-penrose,cs2009-parabolic1,fh1991-repntheory}, define the standard Borel to contain all positive root spaces.  We choose the opposite convention to better facilitate the treatment of homology and subsequently BGG operators; see Section \ref{s:lie-hom}.}
\end{defn}

Thus a standard parabolic is the direct sum of $\ft$, all negative root spaces, and some positive root spaces.  Since the Weyl group of $\fg$ acts transitively on the set of positive subsystems of $\Delta$, any Borel subalgebra of $\fg$ is conjugate to the standard one via the adjoint action of $G$.  Consequently any parabolic is conjugate to a (unique) standard parabolic, yielding the following \cite[Thm.\ 3.2.1]{cs2009-parabolic1}.

\begin{lem} \thlabel{lem:lie-para-csa} Let $\fp\leq\fg$ be a parabolic subalgebra.  Then $\ft$ and $\Delta^{+}$ can be chosen such that $\fp$ is a standard parabolic.  \noproof \end{lem}

The set of standard parabolics \wrt\ $\ft$ and $\Delta^{+}$ may also be enumerated using the set $\Delta^0 \subset \Delta$ of simple roots \cite[p.\ 384]{fh1991-repntheory}.

\vspace{-0.1em}
\begin{prop} \thlabel{prop:lie-para-std} There is a bijection between standard parabolics $\fp \leq \fg$ and subsets $\Sigma\subseteq\Delta^0$, given by mapping $\fp\leq\fg$ to $\Sigma_{\fp} \defeq \setof{ \alpha\in\Delta^0 }{ \fg_{\alpha} \not\leq \fp }$ and conversely by mapping a subset $\Sigma \subseteq \Delta^0$ to the standard parabolic with positive root spaces $\linspan{ \Delta^0\setminus \Sigma }{} \intsct \Delta^{+}$.  \noproof \end{prop}
\vspace{-0.1em}

Since elements of $\Delta^0$ are the nodes of the Dynkin diagram of $\fg$, this suggests an obvious notation for standard parabolics $\fp\leq\fg$: we represent $\fp$ by crossing the nodes corresponding to elements of the associated subset $\Sigma_{\fp} \subseteq \Delta^0$.  Thus the Dynkin diagram with no nodes crossed is the \non proper\ parabolic $\fg$, while crossing all nodes yields the standard Borel; other examples may be found in \cite[\S2.2]{be1989-penrose}.

\vspace{-0.1em}
\begin{lem} \thlabel{lem:lie-para-levi} The subspace $\linspan{ \coroot{\alpha} }{ \alpha\in \Delta^0\setminus\Sigma } \leq \ft$ forms a Cartan subalgebra for the semisimple part of $\fp^0$, while $\liecenter{\fp^0} = \linspan{ H\in\ft }{ \alpha(H)=0 ~\,\forall \alpha \in \Delta^0\setminus\Sigma }$.  \noproof \end{lem}
\vspace{-0.1em}

These two subspaces of $\ft$ are complementary and orthogonal \wrt\ the Killing form \cite[Thm.\ 3.2.1]{cs2009-parabolic1}.  In particular, the dimension of $\liecenter{\fp^0}$ is equal to the number of elements of the corresponding subset $\Sigma_{\fp}$, \ie\ the number of crossed nodes.

Let us also mention briefly how to deal with parabolic subalgebras of a real semisimple Lie algebra.  In this case $\fp\leq\fg$ is called a \emph{standard parabolic} if its complexification is a standard parabolic in the sense of \thref{defn:lie-para-std}, and any parabolic is conjugate to a standard parabolic by the adjoint action of a maximal compact subgroup $K\leq \grp{Int}{\fg}$.  It turns out that there is a bijection between standard parabolics of $\fg$ and subsets $\Sigma \subseteq \Delta^0$ which are disjoint from the set $\Delta^0_c$ of compact simple roots and stable under the \non compact\ root involution, which are themselves in bijection with subsets of restricted simple roots.  Therefore real standard parabolics are classified by crossing white nodes of the Satake diagram%
\footnote{Our convention is that black nodes of the Satake diagram represent compact simple roots, while white nodes represent \non compact\ simple roots.}
of $\fg$, with the caveat that we must also cross all nodes connected by an arrow.  Details and examples may be found in \cite[\S3.2.9]{cs2009-parabolic1}.

\subsection{Filtrations and gradings} 
\label{ss:lie-para-gr}

Let $V$ be a vector space over a field $\bk$, which for us will be $\bR$ or $\bC$.

\begin{defn} \thlabel{defn:lie-para-filtgr} A \emph{filtration} of $V$ is a family $\{V_i\}_{i\in\bZ}$ of subspaces satisfying $V_{i+1} \supset V_i$ for all $i\in\bZ$, and $\Union{i\in\bZ}{} V_i = V$ and $\Intsct{i\in\bZ}{} V_i = \{0\}$.  A \emph{grading} of $V$ is a vector space decomposition $V = \Dsum{i\in\bZ}{} V_{(i)}$. \end{defn}

Typically we are interested in filtrations for which $V_i\neq 0, V$ for only finitely many $i\in\bZ$.  Any filtration $\{V_i\}_{i\in\bZ}$ of $V$ gives rise to a graded vector space $\gr{V}$ with components $V_{(i)} \defeq V_i / V_{i-1}$ called the \emph{associated graded} of $V$.  Although there are no natural linear maps between $V$ and $\gr{V}$ in either direction, there is an isomorphism $V\isom \gr{V}$ given by choosing a splitting of each projection $V_i \surjto V_{(i)}$, thus identifying $V_{(i)}$ with a complement to $V_{i-1}$ in $V_i$.  Given such splittings of filtered $V,W$, a linear map $f:V \to W$ has \emph{homogeneity $k$} if $f(V_{(i)}) \subseteq W_{(i+k)}$ for all $i\in\bZ$, giving a grading of $\Hom{V}{W}$ by homogeneous degree.

Now consider a Lie algebra $\fg$ over $\bk$.  A \emph{filtration} of $\fg$ is a filtration $\{\fg_i\}_{i\in\bZ}$ of the underlying vector space such that $\liebrac{\fg_i}{\fg_j} \subseteq \fg_{i+j}$ for all $i,j\in\bZ$.  Then $\fg_i$ is a subalgebra of $\fg$ for all $i\leq 0$ and an ideal for $i<0$ \cite[Cor.\ 3.2.1]{cs2009-parabolic1}.  Given an unfiltered $\fg$ with a representation on a filtered vector space $V$, we obtain a filtration of $\fg$ by defining $\fg_i \defeq \setof{ X \in \fg }{ X\acts v \in V_{i+j} ~\,\forall v\in V_j }$.

Let $\fp$ be a parabolic subalgebra of $\fg$.  Then by definition $\fp^{\perp}$ is nilpotent, so that the lower central series $(\fp^{\perp})^1 \supset (\fp^{\perp})^2 \supset \cdots \supset 0$ terminates after a finite number of steps, where $(\fp^{\perp})^1 \defeq \fp^{\perp}$ and $(\fp^{\perp})^{k+1} \defeq \liebrac{\fp^{\perp}}{(\fp^{\perp})^k}$.  We shall say that $\fp$ has \emph{height $n$} if $(\fp^{\perp})^k = 0$ for all $k>n$.  It is straightforward to check that
\begin{equation*}
  \fg_0 \defeq \fp, \quad
  \fg_{-1} \defeq \fp^{\perp}
  \quad\text{and}\quad
  \fg_k \defeq \begin{cases} \liebrac{\fp^{\perp}}{\fg_{-k+1}} & k<0 \\
                             (\fg_{-1-k})^{\perp}              & k>0
               \end{cases}
\end{equation*}
defines a filtration of $\fg$, which we refer to as the \emph{$\fp^{\perp}$-filtration} of $\fg$.  If $\fp$ has height $n$ then clearly $\fg_{-(n+1)}=0$ and $\fg_{n}=\fg$, so that there are $2n$ proper filtration components.

The associated graded Lie algebra $\gr{\fg}$ has components $\fg_{(i)} \defeq \fg_i/\fg_{i-1}$ satisfying $\liebrac{\fg_{(i)}}{\fg_{(j)}} \subseteq \fg_{(i+j)}$ for all $i,j\in\bZ$.  Clearly if $\fp$ has height $n$ then $\fg_{(i)}=0$ for $\abs{i}>n$, so that there are $2n+1$ \non zero\ graded components.  Notice also that $\fg_{(0)} \defeq \fp/\fp^{\perp}$ is precisely the reductive Levi factor $\fp^0$ of $\fp$.

\begin{lem} \thlabel{lem:lie-para-weyl} There is a unique $\xi_0 \in \liecenter{\fp^0}$ such that $\liebrac{\xi_0}{X} = iX$ for all $X\in\fp_{(i)}$.  \noproof \end{lem}

The element $\xi_0\in\liecenter{\fp^0}$ is called the \emph{grading element} of $\fp$.  From the definition of the filtration, a choice of splitting of the projection $\fp \surjto \fp^0 \defeq \fp/\fp^{\perp}$ evidently induces splittings of all projections $\fg_i \surjto \fg_{(i)}$; such splittings always exist \cite[Lem.\ 2.2]{cds2005-ricci}.

\begin{defn} \thlabel{defn:lie-para-weyl} An \emph{algebraic Weyl structure} for $\fp$ is a choice of lift of the grading element $\xi_0\in\liecenter{\fp^0}$ to $\fp$ \wrt\ the projection $\fp \surjto \fp^0$. \end{defn}

Thus an algebraic Weyl structure induces an isomorphism $\fg \isom \gr{\fg}$; for abelian parabolics this amounts to an isomorphism $\fg \isom \fg/\fp \dsum \fp^0 \dsum \fp^{\perp}$.  Since the space of such lifts is an affine space modelled on $\fp^{\perp}$, we obtain the following \cite[Lem.\ 2.5]{cds2005-ricci}.

\begin{lem} \thlabel{lem:lie-para-exp} The subgroup $\exp{\fp^{\perp}}\leq P$ acts simply transitively on the affine space of algebraic Weyl structures for $\fp$.  \noproof \end{lem}

It follows that the stabiliser of an algebraic Weyl structure $\xi\in\fp$ is a subgroup of $P$ projecting isomorphically onto $P^0 \defeq P/\exp{\fp^{\perp}}$, so that $\xi$ also splits the quotient group homomorphism $P \surjto P^0$.  This is the basis of {\v C}ap and Slov{\'a}k's treatment of Weyl structures \cite{cs2003-weylstr}; see \cite[App.\ A]{cds2005-ricci} for a detailed comparison.

\begin{defn} \thlabel{defn:lie-para-opp} Parabolics $\fp, \opp{\fp} \leq \fg$ of the same height and with associated filtrations $\{\fg_i\}_{i\in\bZ}, \{\opp{\fg}\}_{i\in\bZ}$ are said to be \emph{opposite} if $\fg_i \intsct \opp{\fg}_i$ is complementary to $\fg_{i-1}$ in $\fg_i$. \end{defn}

For abelian parabolics $\fp, \opp{\fp}$, this amounts to asking that $\fp^{\perp} \intsct \opp{\fp} = 0$.  Generally, the complement $\fg_i \intsct \opp{\fg}_i$ to $\fg_{i-1}$ splits the projection $\fg_i \surjto \fg_{(i)}$; thus the choice of an opposite parabolic is equivalent to the choice of an algebraic Weyl structure \cite[Lem.\ 2.5]{cds2005-ricci}.

\begin{lem} \thlabel{lem:lie-para-expopp} $\exp{\fp^{\perp}}$ acts simply transitively on the set of parabolics opposite to $\fp$.  \noproof \end{lem}

Choose a Cartan subalgebra $\ft\leq\fg$ and a simple subsystem $\Delta^0$ \wrt\ which $\fp$ is a standard parabolic corresponding to a subset $\Sigma \subseteq \Delta^0 = \setof{ \alpha_1, \ldots, \alpha_k }{}$.  Each root $\alpha = \sum{i=1}{k} a_i\alpha_i$ has an associated \emph{$\Sigma$-height} $\height[\Sigma]{\alpha} \defeq \sum{i \, : \, \alpha_i\in\Sigma}{} \, a_i$, and \thref{lem:lie-para-csa} can be adapted to show that each graded component $\fg_{(i)}$ consists of those root spaces of $\Sigma$-height $i$.  In particular $\fg_{(0)} = \fp^0$ consists of the root spaces of height zero, while $\fp^{\perp}$ consists of root spaces with negative height.  It also follows that all parabolics conjugate to $\fp$ have the same height, being given by the $\Sigma$-height of the highest root of $\fg$.  Moreover \thref{lem:lie-para-expopp} implies that the data $(\ft,\Delta^{+},\Sigma)$ may be chosen in such a way that $\opp{\fp}$ is the standard parabolic corresponding to the data $(\ft,-\Delta^{+},-\Sigma)$.%
\footnote{That is, $\opp{\fp}$ is defined as in \thref{prop:lie-para-std} but with the roles of positive and negative roots exchanged.}

\begin{lem} \thlabel{lem:lie-para-abl} Suppose that $\fg$ is simple with abelian parabolic $\fp$.  Then $\liebrac{ \fp^{\perp} }{ \fg } = \fp$. \end{lem}

\begin{proof} Invariance of the Killing form gives $\killing{ \liebrac{\fp^{\perp}}{\fg} }{ \fp^{\perp} } = \killing{ \fg }{ \liebrac{\fp^{\perp}}{\fp^{\perp}} } = 0$, so that $\liebrac{ \fp^{\perp} }{ \fg } \subseteq (\fp^{\perp})^{\perp} = \fp$ by \non degeneracy.  Conversely, choose an algebraic Weyl structure for $\fp$ and hence an isomorphism $\fg \isom \fg/\fp \dsum \fp^0 \dsum \fp^{\perp}$.  Since $\fg$ is simple, $\liebrac{\fp^{\perp}}{\fp^0} = \fp^{\perp}$ and $\liebrac{\fp^{\perp}}{\fg/\fp} = \fp^0$ by \cite[Prop.\ 3.1.2(4)]{cs2009-parabolic1}; thus $\fp \isom \fp^0 \dsum \fp^{\perp} = \liebrac{\fp^{\perp}}{\fg/\fp \dsum \fp^0} \subseteq \liebrac{\fp^{\perp}}{\fg}$. \end{proof}

\vspace{0.1em}
\subsection{Representations of a parabolic} 
\label{ss:lie-para-repn}

Generally speaking the representation theory of a parabolic subalgebra $\fp\leq\fg$ is quite complicated, but there are significant simplifications for completely reducible representations.  Let $\fp^0 \defeq \fp/\fp^{\perp}$ be the reductive Levi factor of $\fp$.

\begin{lem}
\thlabel{lem:lie-para-prepn} \rmcite[Prop.\ 3.2.12(1)]{cs2009-parabolic1}
Every completely reducible $\fp$-representation $V$ is the trivial lift of a completely reducible $\fp^0$-representation.  Moreover, the grading element $\xi_0 \in \liecenter{\fp^0}$ acts by a scalar on each $\fp^0$-irreducible component of $V$.
\noproof \end{lem}

Indeed, in this case the corresponding linear map $\fp \to \alg{gl}{V}$ factors through the projection $\fp \surjto \fp^0$, giving a representation of $\fp^0$ on $V$ which pulls back to the given $\fp$-representation.  Since $\fp^0$ is not semisimple, its representations are not automatically completely reducible.  In fact, a $\fp^0$-representation is completely reducible if and only if its centre $\liecenter{\fp^0}$ acts diagonalisably \cite[\S3.2.12]{cs2009-parabolic1}.

The representation theory of $\fp$ can also be described in terms of highest weights.  For this, choose a Cartan subalgebra $\ft\leq\fg$ and a simple subsystem $\Delta^0$ \wrt\ which $\fp$ is the standard parabolic corresponding to a subset $\Sigma\subseteq\Delta^0$.  Recall that (isomorphism classes of) irreducible $\fg$-representations are in bijection with dominant integral weights \cite[Thm.\ 14.18]{fh1991-repntheory}, \ie\ those which can be written as a \non negative\ integral linear combination of the fundamental weights.  By \thref{lem:lie-para-prepn}, irreducible $\fp$-representations are given by a representation of the semisimple part $\sspart{\fp^0}$ and a linear functional on the centre $\liecenter{\fp^0}$.  We say that a weight $\lambda\in \ft^*$ is \emph{$\fp$-dominant} (respectively, \emph{$\fp$-integral}) if the Cartan number $\cartanint{\lambda}{\alpha}$ is \non negative (respectively, integral) for all $\alpha \in \Delta^0 \setminus \Sigma$.

\begin{prop} \emph{\cite[Cor.\ 3.2.12]{cs2009-parabolic1}} There is a bijection between isomorphism classes of irreducible $\fp$-representations and $\fp$-dominant and $\fp$-integral weights $\lambda \in \ft^*$. \noproof \end{prop}

In Dynkin diagram notation, the $\fp$-dominant and $\fp$-integral weights are precisely those with \non negative\ integer coefficients over the uncrossed nodes.  If we are only interested in $\fp$-representations there is no restriction on the coefficients over crossed nodes.  However if we ask that a representation integrates to a representation of a parabolic subgroup $P$ with Lie algebra $\fp$, it turns out that the coefficients over crossed nodes must be integers \cite[\S3.2.12]{cs2009-parabolic1}.

In the sequel we shall restrict to the following subclass of $\fp$-representations.

\begin{defn} \thlabel{defn:lie-para-prepn} A $\fp$-representation $V$ is \emph{filtered} if there is a finite $\fp$-invariant filtration
\begin{equation*}
  V = V_N \supset V_{N-1} \supset \cdots \supset V_0 \supset 0
\end{equation*}
such that each graded component $V_{(i)}$ is a completely reducible $\fp$-representation.  We henceforth redefine a \emph{$\fp$-representation} to mean a filtered $\fp$-representation. \end{defn}

Since each $V_{(i)}$ is completely reducible, the grading element $\xi_0 \in \liecenter{\fp^0}$ acts by a scalar on each irreducible component of the induced $\fp^0$-representation, called its \emph{(geometric) weight}. An algebraic Weyl structure $\xi$ for $\fp$ then splits the filtration on $V$ into the eigenspaces of $\xi$.

The restriction to $\fp$ of a $\fg$-representation $\bV$ is a filtered $\fp$-representation.  To see this, note that since $\fp^{\perp} \leq \fg$ is a nilpotent subalgebra, it acts nilpotently on $\bV$ by Engel's theorem.  Consequently we obtain a finite filtration
\vspace{-0.25em}
\begin{equation} \label{eq:lie-para-pTfilt}
  \bV = \bV_N \supset \fp^{\perp}\acts\bV \supset \cdots
              \supset (\fp^{\perp})^N \acts \bV = \bV_0 \supset 0
\vspace{-0.25em}
\end{equation}
of $\bV$, which we call the \emph{$\fp^{\perp}$-filtration} of $\bV$.  Clearly $\fp^{\perp}$ acts trivially on each graded component $\bV_{(i)}$, while the identity $\liebrac{\fp}{\fp^{\perp}} = \fp^{\perp}$ implies that \eqref{eq:lie-para-pTfilt} is $\fp$-invariant.  Moreover \thref{lem:lie-para-csa} allows us to choose a Cartan subalgebra $\ft \leq \fg$ and a simple subsystem $\Delta^0$ \wrt\ which $\fp$ is a standard parabolic, implying that $\liecenter{\fp^0} \leq \ft$ acts diagonalisably on $\bV$ and hence that the $\bV_{(i)}$ are completely reducible.

The lowest filtration component $\bV_0$ in \eqref{eq:lie-para-pTfilt} is sometimes called the \emph{socle} of $\bV$, with $N$ the \emph{height} of $\bV$. Dually, the first graded component $\liehom{0}{\bV} \defeq \bV_{(N)} = \bV / (\fp^{\perp} \acts \bV)$ is called the \emph{top} of $\bV$.  The homological notation will be explained in Section \ref{s:lie-hom}.

Let $\liecohom{0}{\bV} \defeq \setof{ v \in \bV }{ \alpha \acts v = 0 ~\, \forall \alpha \in \fp^{\perp} }$ denote the kernel of the $\fp^{\perp}$-action on $\bV$, which is a $\fp$-subrepresentation of $\bV$ since $\liebrac{\fp}{\fp^{\perp}} = \fp^{\perp}$.  We have the following relation between $\fg$-subrepresentations of $\bV$ and $\fp$-subrepresentations of $\liecohom{0}{\bV}$.

\begin{prop}
\thlabel{prop:lie-para-grepn}  \rmcite[Prop.\ 3.2.13]{cs2009-parabolic1}
There is a bijection between irreducible $\fg$-subrep\-resentations of $\bV$ and irreducible $\fp$-subrepresentations of $\liecohom{0}{\bV}$.  In particular if $\bV$ is the irreducible $\fg$-representation with lowest weight $\lambda$, then $\liecohom{0}{\bV}$ is the irreducible $\fp$-representation with the same lowest weight.
\noproof \end{prop}

\begin{cor} \thlabel{cor:lie-para-socle}
$\liecohom{0}{\bV}$ coincides with the socle $\bV_0$ of the $\fp^{\perp}$-filtration \eqref{eq:lie-para-pTfilt}.
\end{cor}

\begin{proof} Since the action of $\fp^{\perp}$ preserves the $\fg$-irreducible components of $\bV$, it suffices to consider the case that $\bV$ is an irreducible $\fg$-representation.  Then $\liecohom{0}{\bV}$ is an irreducible $\fp$-representation by \thref{prop:lie-para-grepn}.  However $\bV_0$ is a $\fp$-subrepresentation of $\liecohom{0}{\bV}$ by construction, giving equality as claimed. \end{proof}

Then since $f \in \bV^*_0 = \liecohom{0}{\bV^*}$ if and only if $(\alpha \acts f)(v) = -f(\alpha \acts v) = 0$ for all $\alpha \in \fp^{\perp}$ and $v \in \bV$, the socle $\bV^*_0$ of $\bV^*$ equals the annihilator of $\fp^{\perp} \acts \bV$.  Therefore $(\bV^*_0)^* \isom \bV / (\fp^{\perp} \acts \bV) = \liehom{0}{\bV}$.  Then if $\bV$ has highest weight $\lambda$ as a $\fg$-representation, \thref{prop:lie-para-grepn} implies immediately that $\liehom{0}{\bV}$ has highest weight $\lambda$ as a $\fp$-representation, so that $\liehom{0}{\bV}$ is obtained by ``putting the crosses in'' to the Dynkin diagram of $\bV$; we shall interpret this homologically in Section \ref{s:lie-hom}.

\thref{prop:lie-para-grepn} also allows us to compute the weight of an irreducible $\fp$-representa\-tion $V$ as follows.  Choosing data so as to identify $\fp$ with the standard parabolic associated to a subset $\Sigma$ of simple roots, the weight of $V$ is given by the $\Sigma$-height $\height[\Sigma]{\lambda}$ of the highest weight $\lambda$ of $V$.  Writing $\lambda$ in terms of the fundamental weights using the inverse Cartan matrix $C^{-1}$ of $\fg$, it follows easily that $V$ has weight $\lowest[\fp]^{\top} C^{-1} \lambda$, where $\lowest[\fp]$ is the $\fp$-dominant weight with a one over each crossed node.  Tables of the inverse Cartan matrices for complex simple Lie algebras may be found in \cite[Tbl.\ B.4]{cs2009-parabolic1}.

\begin{expl} \thlabel{expl:lie-para-e6wt} Consider the standard parabolic subalgebra $\fp$ of $\alg[_6]{e}{\bC}$ and its irreducible representation $V$ defined by
\begin{equation*}
  \fp = \dynkinEp{6}{lie-para-repn-e6p}
  \quad\text{acting on}\quad
   \dynkinEp[1,0,3,0,-5,1]{6}{lie-para-repn-e6V} = V.
\end{equation*}
Numbering the nodes of the Dynkin diagram ``clockwise'' starting with the left-most node, the weight of $V$ is
\begin{equation*}
  \lowest[\fp]^{\top} C^{-1} \lambda =
  \tfrac{1}{3} \big(\begin{smallmatrix} 0 & 0 & 0 & 0 & 1 & 0 \end{smallmatrix}\big) \!\!
               \left(\begin{smallmatrix} 4 & 5  & 6  & 4  & 2 & 3 \\
                                         5 & 10 & 12 & 8  & 4 & 6 \\
                                         6 & 12 & 18 & 12 & 6 & 9 \\
                                         4 & 8  & 12 & 10 & 5 & 6 \\
                                         2 & 4  & 6  & 5  & 4 & 3 \\
                                         3 & 6  & 9  & 6  & 3 & 6
                     \end{smallmatrix}\right) \!\!\!
               \left(\!\begin{smallmatrix} 1 \\ 0 \\ 3 \\ 0 \\ -5 \\ 1 \end{smallmatrix}\!\right)
  = 1.
\end{equation*}
\end{expl}

\subsection{\Rspaces\ and \proj\ embeddings} 
\label{ss:lie-para-rspace}

Let $G$ be a semisimple Lie group with Lie algebra $\fg$ and recall that the adjoint action of $G$ takes parabolic subalgebras to parabolic subalgebras.

\begin{defn} \thlabel{defn:lie-para-rspace} An \emph{\Rspace} for $G$ is an adjoint orbit of parabolic subalgebras.  An \Rspace\ is \emph{symmetric} if its parabolic subalgebras are abelian. \end{defn}

The \Rspace\ $G\acts\fp$ is also known as the \emph{generalised flag manifold} $G/P$, where $P\leq G$ is a parabolic subgroup with Lie algebra $\fp$; see \cite{be1989-penrose}.  The link is formalised as follows \cite[Lem.\ 1.6]{bdpp2011-rspaces}.

\begin{prop} \thlabel{prop:lie-para-flagmfd} Given an \Rspace\ $G\acts\fp$, the stabiliser $P' \defeq \Stab[G]{\fp'} \leq G$ of a parabolic $\fp' \in G\acts\fp$ is a parabolic subgroup such that $G/P'$ is diffeomorphic to $G\acts\fp$.  \noproof \end{prop}

In the complex setting, this result can be stretched considerably further: then $G/P$ is a compact \proj\ \kahler\ manifold with a holomorphic $G$-action \cite[\S3.2.6]{cs2009-parabolic1}.  In particular, the complex symmetric \Rspaces\ are precisely the hermitian symmetric spaces of compact type for the maximal compact subgroup \cite{bdpp2011-rspaces}.

Suppose now that $G$ is complex and connected, choose a Cartan subalgebra $\ft\leq\fg$ and a positive subsystem $\Delta^{+}$, and let $V$ be an irreducible $\fg$-representation with highest weight $\lambda\in\ft^*$.  Then $-\lambda$ is the lowest weight of the dual representation $V^*$, and since the weight space $V_{-\lambda}^*$ is $1$-dimensional we obtain a well-defined point in $\pr{V^*}$.

\begin{prop} \thlabel{prop:lie-para-proj} Let $V$ be the irreducible $\fg$-representation of highest weight $\lambda\in\ft^*$.  Then the stabiliser $\fp$ of the lowest weight space in $V^*$ is the standard parabolic defined by $\Sigma \defeq \setof{ \alpha\in\Delta^0 }{ \killing{\lambda}{\alpha} \neq 0 }$, hence giving a \proj\ embedding $G\acts\fp \injto \pr{V^*}$. \end{prop}

\begin{sketchproof} We outline the proof from \cite[Prop.\ 3.2.5]{cs2009-parabolic1}.  The dual representation $V^*$ has lowest weight $-\lambda$, so let $v_0 \in V^*_{-\lambda}$ be a lowest weight vector and define $\fp$ to be the stabiliser $\setof{X \in \fg}{X \acts V^*_{-\lambda} \subseteq V^*_{-\lambda}}$ of $V^*_{-\lambda}$.  It follows easily that $\fp$ is the standard parabolic with corresponding subset $\Sigma \defeq \setof{\alpha \in \Delta^0}{\fg_{\alpha} \not\leq \fp}$.

Choose a simple root $\alpha \in \Delta^0$, and an \sltriple\ $e \in \fg_{\alpha}$, $f \in \fg_{-\alpha}$ and $h \defeq \liebrac{e}{f} \in \ft$.  By considering the root reflection through $\alpha$, one then shows that $e \acts v_0 = 0$ if and only if $\killing{\alpha}{\lambda} = 0$; since $\fg_{\alpha} \leq \fp$ if and only if $\fg_{\alpha} \acts v_0 = 0$, the claimed form of $\Sigma$ follows.

To obtain the \proj\ embedding note that the adjoint action of $P\defeq \Stab[G]{\fp}$ preserves $V^*_{-\lambda}$; we easily conclude that $P \leq \Stab[G]{[v_0]}$, which is actually an equality by \thref{prop:lie-para-flagmfd}.  Therefore the holomorphic submersion $G \to G\acts[v_0]$ given by $g \mapsto g\acts[v_0] = [g\acts v_0]$ factors to a holomorphic bijection $G\acts\fp \isom G/P \to G\acts[v_0]$, which is a diffeomorphism and hence a biholomorphism by compactness of $G/P$. \end{sketchproof}

It is straightforward to determine the Dynkin notation for $\fp$ from the highest weight $\lambda \in \ft^*$ of $V$.  Given a simple root $\alpha\in\Delta^0$, the coefficient of the corresponding fundamental weight is $\cartanint{\lambda}{\alpha}$, which is the number we write over the node of the Dynkin diagram.  Thus $\killing{\lambda}{\alpha} \neq 0$ if and only if there is a \non zero\ coefficient over that node, meaning $\fp$ is given by crossing nodes in the support of $\lambda$.

\begin{cor} \thlabel{cor:lie-para-proj} There is a \proj\ embedding $G\acts\fp \injto \pr{V^*}$ for any irreducible $\fg$-representation $V$ whose highest weight is supported on the crossed nodes of $\fp$.  \noproof \end{cor}

\begin{expl} \thlabel{expl:lie-para-proj} Consider the \Rspace\ $G\acts\fp = \smash{ \dynkinApk{1}{1}{lie-para-rspace-sl} }$ of $\fg = \alg{sl}{n+1,\bC}$, where we have crossed the $k$th node.  The corresponding fundamental representation is the exterior power $\Wedge[\bC]{k}\bC^{n+1}$ of the standard representation of $\fg$ on $\bC^{n+1}$.  If $\{e_i\}_i$ and $\{\ve^i\}_i$ are the standard dual bases of $\bC^{n+1}$ and $\bC^{n+1*}$ then $v^k \defeq \ve^1 \wedge\cdots\wedge \ve^k \in \Wedge[\bC]{k} \bC^{n+1*}$ is a lowest weight vector for $\fg$, and the stabiliser of the line through $v^k$ coincides with the stabiliser of the $k$-dimensional subspace $\linspan{ \ve^1,\ldots,\ve^k }{} \leq \bC^{n+1*}$.  Thus $G\acts\fp$ is the complex \grassmannian\ $\Grass{k}{\bC^{n+1*}} \isom \Grass{n+1-k}{\bC^{n+1}}$.  In particular $k=n$ gives the \grassmannian\ of hyperplanes in $\bC^{n+1*}$, which is just $\CP[n]$. \end{expl}

Notice that there is a minimal \proj\ embedding $G\acts\fp \injto \pr{V^*}$, defined by taking $V$ to be the irreducible $\fg$-representation whose highest weight has a one over each crossed node of $\fp$.  Generally speaking, choosing $V$ to have a strictly larger weight results in $G\acts\fp$ having larger codimension in $\pr{V^*}$.  Thanks to a result of Kostant, it is possible to describe the image of $G\acts\fp$ inside $\pr{V^*}$ as an intersection of quadrics.%
\footnote{Kostant never published this result, and it appears (with attribution) in \cite{lt1979-flagalg,l1982-hwquadric}; also see \cite[p.\ 368]{p2007-lieinv}.  The author is grateful to Fran Burstall for his patient explanation of Kostant's results.}

For this, recall that the Cartan square $\Cartan{2}V^*$ is the highest weight subrepresentation of $\Symm{2}V^*$ and appears with multiplicity one.  Viewing $\Symm{2}V$ as the space of homogeneous quadratic polynomials on $V^*$, the projection $\Symm{2}V^* \surjto U^* \defeq \Symm{2}V^* \!/  \Cartan{2}V^*$ is dual to the inclusion $U \injto \Symm{2}V$, which identifies $U \leq \Symm{2}V$ with the annihilator of $\Cartan{2}V^*$.

Now if $v_0 \in V^*_{-\lambda}$ is a lowest weight vector for $\fg$ and $[v]\in G\acts[v_0]$, we have $v = g\acts v_0$ for some $g\in G$ and hence $v\tens v = (g\acts v_0) \tens (g\acts v_0) = g \acts (v_0 \tens v_0)$.  Since $v_0\tens v_0$ is a lowest weight vector for $\Cartan{2}V^*$, this gives an inclusion
\begin{equation} \label{eq:lie-para-kostant}
  G\acts[v_0] \subseteq \setof{ v \in V^* }{ f(v\tens v)=0 ~\,\forall f\in U }
\end{equation}
of $G\acts\fp \isom G\acts[v_0]$ into the intersection of quadrics cut out by $U$.  By computing the action of Casimir elements it is possible to deduce that \eqref{eq:lie-para-kostant} is an equality, hence describing exactly which quadratic equations cut out $G\acts\fp$ inside $\pr{V^*}$.

\begin{thm}[Kostant] \thlabel{thm:lie-para-kostant} Suppose that $G$ is a complex semisimple Lie group with Lie algebra $\fg$, and let $V$ be an irreducible $\fg$-representation of highest weight $\lambda$.  Let $G\acts[v_0]$ denote the lowest weight orbit in $V^*$ and let $U\leq\Symm{2}V$ be the annihilator of $\Cartan{2}V^*$.  Then:
\begin{enumerate}
  \item $G\acts[v_0]$ is the intersection of quadrics cut out by $U$, so that \eqref{eq:lie-para-kostant} is an equality.
  
  \item The ideal of $G\acts[v_0]$ is generated by $U$.
  
  \item The homogeneous coordinate ring of $G\acts[v_0]$ is $\Cartan{}V \defeq {\textstyle \Dsum{i=0}{\infty}}\, \Cartan{i}V$.  \noproof
\end{enumerate}
\end{thm}

\begin{expl} Continuing notation from \thref{expl:lie-para-proj}, consider the resulting embedding $G\acts\fp \injto \pr{\Wedge[\bC]{k} \bC^{n+1*}}$ which identifies $\fp$ with the stabiliser $\linspan{\ve^1, \ldots, \ve^k}{} \leq \bC^{n+1*}$.  In the case $k=n$ we have $G\acts\fp = \CP[n]$, and the Cartan square of $\Wedge[\bC]{n} \bC^{n+1*} \isom \bC^{n+1}$ coincides with its symmetric square, representing the fact that $\CP[n] \injto \pr{\bC^{n+1}}$ is the minimal \proj\ embedding.  In the case $k=n-1$ we have $G\acts\fp = \Grass{2}{\bC^{n+1}}$, with
\begin{equation*}
\vspace{-0.3em}
  \Symm{2} \Wedge{n-1}\bC^{n+1}
    = \Symm{2}\bigg( \dynkinA[0,0,0,1,0]{2}{0}{3}{} \bigg)
    = \dynkinname{ \dynkinA[0,0,0,2,0]{2}{0}{3}{} }
                 { \Cartan{2} \Wedge{n-1}\bC^{n+1} }
      \dsum
      \dynkinname{ \dynkinA[0,0,0,1,0,0,0]{2}{0}{5}{} }
                 { U = \Wedge{n-3}\bC^{n+1} }.
\vspace{-0.3em}
\end{equation*}
Identifying $U \isom \Wedge{4}\bC^{n+1*}$, we see that $G\acts[v_0]$ consists of those elements $[v] \in \pr{\Wedge{2}\bC^{n+1}}$ for which $v \wedge v = 0$, which is just the Pl{\"u}cker embedding $\Grass{2}{\bC^{n+1}} \injto \pr{\Wedge{2}\bC^{n+1}}$. \end{expl}

\section{Lie algebra homology and cohomology} 
\label{s:lie-hom}

Lie algebra homology and cohomology were introduced by Kostant \cite{k1961-liehom} to provide an algebraic backdrop for the Bott--Borel--Weil theorem, which computes the sheaf cohomology of generalised flag manifolds \cite{be1989-penrose,d1976-simplebott}.  Parabolic geometries are modelled on such manifolds, so we shall be interested in a ``curved'' analogue of Kostant's results provided by \cite{cd2001-curvedbgg,css2001-bgg}.  We introduce the necessary definitions and basic results in Subsection \ref{ss:lie-hom-defn}, before outlining the description of the Hasse diagram and Kostant's version of the Bott--Borel--Weil theorem in Subsection \ref{ss:lie-hom-bbw}.

\subsection{Basic definitions} 
\label{ss:lie-hom-defn}

Let $\fg$ be any Lie algebra and $V$ a $\fg$-representation.  We define the space $\liechain[\fg]{k}{V} \defeq \Wedge{k}\fg \tens V$ of \emph{$k$-chains} on $\fg$ with values in $V$ and a \emph{boundary map}
\begin{equation*}
\begin{aligned}
  \liebdy : \liechain[\fg]{k}{V} &\to \liechain[\fg]{k-1}{V} \\
  \beta \tens v
  &\mapsto \sum{i}{}\, (\beta \intprod \ve^i) \tens (e_i \acts v)
    + \sum{i<j}{}\, \liebrac{e_i}{e_j} \wedge ((\beta \intprod \ve^i) \intprod \ve^j) \tens v,
\end{aligned}
\end{equation*}
where $\{e_i\}_i$ is a basis of $\fg$ with dual basis $\{\ve^i\}_i$.  It can be checked directly \cite[Lem.\ 3.2]{cd2001-curvedbgg} that $\liebdy$ is independent of the choice of basis and satisfies $\liebdy^2 = 0$, so that $(\liechain[\fg]{}{V},\liebdy)$ forms a chain complex.

\begin{defn} \thlabel{defn:lie-hom-hom} The homology $\liehom[\fg]{}{V}$ of $(\liechain[\fg]{}{V}, \liebdy)$ is called the \emph{Lie algebra homology} of $\fg$ with values in $V$. \end{defn}

In particular, the zeroth homology $\liehom[\fg]{0}{V} = V/(\fg\acts V)$ is the space of \emph{\co invariants} of $V$.  Evidently $\liechain[\fg]{k}{V}$ carries a natural action of $\fg$ given by extending the adjoint action of $\fg$ on itself by the representation $V$.  Moreover it is clear that $\liebdy$ is $\fg$-equivariant, so there is an induced representation on $\liehom[\fg]{k}{V}$.

Now suppose that $\fg$ is semisimple, $\fp\leq\fg$ is parabolic and $V$ is a $\fg$-representation.  Since $\liebrac{\fp}{\fp^{\perp}} = \fp^{\perp}$, the chain space $\liechain{k}{V}$ is naturally a $\fp$-representation.

\begin{lem} \emph{\cite[Lem.\ 3.3]{cd2001-curvedbgg}} \thlabel{lem:lie-hom-homeq} $\liebdy : \liechain{k}{V} \to \liechain{k-1}{V}$ is $\fp$-equivariant, so that there is a natural representation of $\fp$ on $\liehom{k}{V}$.  \noproof \end{lem}

For $\fg$ simple and $\fp$ abelian, \thref{lem:lie-para-abl} immediately yields the following.

\begin{cor} \thlabel{cor:lie-hom-h0g} Let $\fg$ be simple with abelian parabolic $\fp$.  Then $\liehom{0}{\fg} = \fg/\fp$.  \noproof \end{cor}

Lie algebra cohomology may be defined by a dual approach.  For any Lie algebra $\fg$, the space $\liecochain[\fg]{k}{V} \defeq \Wedge{k}\fg^* \tens V$ of \emph{$k$-cochains} on $\fg$ with values in $V$ may be identified with $\liechain[\fg]{k}{V^*}^*$, with \emph{differential}%
\footnote{Note that our use of $\liebdy,\liediff$ is reversed compared to some authors' conventions: one often sees $\p:\liecochain[\fg]{k-1}{V} \to \liecochain[\fg]{k}{V}$ as the differential and $\p^* : \liechain[\fg]{k}{V} \to \liechain[\fg]{k-1}{V}$ as the boundary map.}
\begin{equation*} \begin{aligned}
  \liediff : \liecochain[\fg]{k-1}{V} &\to \liecochain[\fg]{k}{V} \\
  \beta\tens v &\mapsto \sum{i}{}\, \ve^i\wedge\beta \tens (e_i \acts v)
    + \sum{i<j}{}\, \ve^i \wedge \ve^j \wedge (\liebrac{e_j}{e_i} \intprod \beta) \tens v
\end{aligned}
\end{equation*}
given by (minus) the transpose of $\liebdy : \liechain[\fg]{k}{V^*} \to \liechain[\fg]{k-1}{V^*}$.  Then $(\liediff)^2 = 0$ again, so that $(\liecochain[\fg]{k}{V},\liediff)$ is a cochain complex whose cohomology $\liecohom[\fg]{}{V}$ is called the \emph{Lie algebra cohomology} of $\fg$ with values in $V$.  In particular, the zeroth cohomology $\liecohom[\fg]{0}{V}$ is the kernel of the $\fg$-action on $V$.

Returning to the case that $\fg$ is semisimple and $\fp\leq\fg$ is parabolic, the duality $(\fg/\fp)^* \isom \fp^{\perp}$ induced by the Killing form means that $\liecochain{k}{V} \isom \liechain[\fg/\fp]{k}{V}$.  Choosing a parabolic $\opp{\fp}$ opposite to $\fp$, it follows that $\liediff : \liecochain{k-1}{V} \to \liecochain{k}{V}$ is $\opp{\fp}$-equivariant (but notably \emph{not} $\fp$-equivariant).  Therefore $\liecohom{k}{V}$ is not naturally a $\fp$-representation, but only a $\fp^0$-representation \wrt\ the chosen algebraic Weyl structure.  For this reason, we prefer to work with homology over cohomology.

Continuing to work with an algebraic Weyl structure, it is possible to find positive definite inner products on $\fg$ and $V$ \wrt\ which $\liebdy : \liechain{k}{V} \to \liechain{k-1}{V}$ is (minus) the adjoint of $\liediff : \liechain{k-1}{V} \to \liechain{k}{V}$, where here and below we use that $\liechain{k}{V} = \liecochain[\fg/\fp]{k}{V}$.  This allows us to define a laplacian-like operator \cite{k1961-liehom}.

\begin{defn} \thlabel{defn:lie-hom-quabla} The $\fp^0$-homomorphism $\quab \defeq \liebdy\liediff + \liediff\liebdy : \liechain{k}{V} \to \liechain{k}{V}$ is called the \emph{algebraic laplacian}. \end{defn}

Kostant proves \cite[Prop.\ 2.1]{k1961-liehom} that $\quab$ also induces a Hodge decomposition
\vspace{-0.4em}
\begin{equation} \label{eq:lie-hom-hodge}
  \liechain{k}{V} = \im\liebdy \dsum \ker\quab \dsum \im\liediff,
\vspace{-0.4em}
\end{equation}
of chain spaces, where $\ker\liebdy = \im\liebdy \dsum \ker\quab$ and $\ker\liediff = \ker\quab \dsum \im\liediff$; also see \cite[\S3.3.1]{cs2009-parabolic1} Consequently $\liehom{k}{V} \isom \ker\quab$, identifying homology classes with harmonic representatives in $\ker\quab = \ker\liebdy \intsct \ker\liediff$.  In particular we have isomorphisms
\vspace{-0.3em}
\begin{equation} \label{eq:lie-hom-isom}
  \liehom{k}{V} \isom \liecohom[\fg/\fp]{k}{V}
    \isom \liehom[\fg/\fp]{k}{V^*}^* \isom \liecohom{k}{V^*}^*
\vspace{-0.3em}
\end{equation}
of $\fp^0$-representations for all $k \in \bZ$; see \cite[p.\ 12]{cd2001-curvedbgg}.

\subsection{Computation of homology components} 
\label{ss:lie-hom-bbw}

Given a $\fg$-representation $V$, it will often be useful to compute the $\fp^0$-irreducible components of $\liehom{k}{V}$ explicitly.  We saw how to do this for $\liehom{0}{V}$ in Subsection \ref{ss:lie-para-repn}: it is the top of the $\fp^{\perp}$-filtration of $V$, so is the $\fp$-representation with highest weight given by ``putting the crosses in''.  For arbitrary degree, the Hodge decomposition \eqref{eq:lie-hom-hodge} implies that the $\fp^0$-irreducible components of $\liehom{k}{V}$ may be identified with the $\fp^0$-irreducible components of $\liechain{k}{V}$ which lie in the kernel of $\quab$.
	
If $V$ has highest weight $\lambda\in\ft^*$ for $\fg$, we may restrict $V$ to a $\fp^0$-representation and decompose into irreducible $\fp^0$-subrepresentations.  The subspace $V^{\mu}$ generated by all highest weight vectors of weight $\mu\in\ft^*$ for $\fp^0$ is called the \emph{$\mu$-isotypical component}; evidently only finitely many isotypical components are \non zero.  One can show that $V$ is the direct sum of its isotypical components, with $V^{\mu}$ isomorphic to a direct sum of a number of copies of the irreducible representation with highest weight $\mu$ \cite[\S2.2.14]{cs2009-parabolic1}.  The action of $\quab$ on $V^{\mu}$ is given by Kostant's spectral theorem \cite{k1961-liehom}; see also \cite{cs2007-curvedcasimirs}.

\begin{thm} \thlabel{thm:lie-hom-spectral} Let $V$ be an irreducible $\fg$-representation of highest weight $\lambda$.  Then $\quab : \liechain{k}{V} \to \liechain{k}{V}$ acts on the $\mu$-isotypical component of $\liechain{k}{V}$ by multiplication by $-\tfrac{1}{2}( \norm{\lambda + \lowest}^2 - \norm{\mu + \lowest}^2)$, where $\lowest$ is the lowest form of $\fg$.%
\footnote{The \emph{lowest form} $\lowest[\fg]$ of $\fg$ is the sum of the fundamental weights of $\fg$, so is represented by a one over each node of the Dynkin diagram.}
\noproof \end{thm}

To analyse the resulting weight condition $\norm{\lambda + \lowest} = \norm{\mu + \lowest}$, it is necessary to introduce some machinery called the \emph{Hasse diagram} of $\fp$.  We shall not need the details, so only give a synopsis; see \cite{be1989-penrose,cs2009-parabolic1} for details.

First observe that the Weyl group $\sW_{\fg}$ of $\fg$ may be given the structure of a directed graph, with vertices the elements of $\sW_{\fg}$ and an edge $w \overset{\alpha}{\to} w'$ if and only if $\ell(w') = \ell(w)+1$ and $w' = \sigma_{\alpha}w$ for some $\alpha \in \Delta^{+}$.  Moreover since $\sW_{\fg}$ acts transitively on the set of Weyl chambers, the vertex set of $\sW_{\fg}$ is in bijection with the Weyl orbit of $\lowest\in\ft^*$.  Thus, considering the form of the Cartan matrix, all elements of $\sW_{\fg}$ are obtained by repeatedly applying the simple root reflection rules from Figure \ref{fig:lie-hom-refl} to $\lowest$; see \cite[\S4.1]{be1989-penrose}.

Now suppose that $\fp\leq\fg$ is a standard parabolic corresponding to a subset $\Sigma \subseteq \Delta^0$ and let $\sW_{\fp}$ be the Weyl group of the semisimple part of $\fp^0$, which may be naturally viewed as a subgraph of $\sW_{\fg}$.  The \emph{Hasse diagram} $\sW^{\fp}$ of $\fp$ is the subgraph of $\sW_{\fg}$ with vertices the elements whose action sends any $\fg$-dominant weight to a $\fp$-dominant weight.  By \cite[Prop.\ 5.13]{k1961-liehom}, every $w\in\sW_{\fg}$ may be written as $w = w_{\fp}w^{\fp}$ for elements $w_{\fp} \in \sW_{\fp}$ and $w^{\fp} \in \sW^{\fp}$ of minimal length.  Moreover the stabiliser of the lowest form $\lowest[\fp]$ of $\fp$ in $\sW_{\fg}$ is precisely $\sW_{\fp}$, so that the Weyl orbit of $\lowest[\fp]$ is in bijection with the vertex set of $\sW^{\fp}$.  Thus the vertices of the Hasse diagram may be computed by repeatedly applying simple reflections to $\lowest[\fp]$.

\begin{figure}[h]
\newcounter{DynkinReflCount}  
\setcounter{DynkinReflCount}{1}
  \newcommand{\ReflShell}[5]{
    \sigma_{\alpha_i} \Big(~
        \dynkin{ \DynkinEllipsis{-2}{0}{0}{0}\DynkinLine{-0.5}{0}{0}{0}
                 \DynkinLine{0}{0}{1}{0} #1
                 \DynkinLine{2}{0}{2.5}{0} \DynkinEllipsis{2}{0}{4}{0}
                 \DynkinWDot[$a$]{0}{0} \DynkinWDot[$b$]{1}{0} \DynkinWDot[$c$]{2}{0} }
               {lie-hom-hasse-refl\theDynkinReflCount}
      ~\Big)
    \stepcounter{DynkinReflCount}
    &= ~\dynkin{ \DynkinEllipsis{-2}{0}{0}{0} \DynkinLine{-0.5}{0}{0}{0}
                 \DynkinLine{0}{0}{2}{0} #2
                 \DynkinEllipsis{3}{0}{5}{0} \DynkinLine{3}{0}{3.5}{0}
                 \DynkinWDot[$#3$]{0}{0} \DynkinWDot[$#4$]{1.5}{0} \DynkinWDot[$#5$]{3}{0} }
               {lie-hom-hasse-refl\theDynkinReflCount}
    \stepcounter{DynkinReflCount}
  }
  \begin{align*}
    \ReflShell{\DynkinLine{1}{0}{2}{0}}{\DynkinLine{1.5}{0}{3}{0}}
              {a{+}b}{-b}{b{+}c} \\
    \ReflShell{\DynkinDoubleLine{1}{0}{2}{0}}{\DynkinDoubleLine{1.5}{0}{3}{0}}
              {a{+}b}{-b}{2b{+}c} \\
    \ReflShell{\DynkinDoubleLine{2}{0}{1}{0}}{\DynkinDoubleLine{3}{0}{1.5}{0}}
              {a{+}b}{-b}{b{+}c}
  \end{align*}
  \caption[Recipes for performing simple root reflections]
          {Recipes for performing a simple root reflection at the central node.  If $b$ is the coefficient that node, we add $b$ to adjacent nodes, with multiplicity if there are multiple edges, and then replace $b$ by $-b$.}
  \label{fig:lie-hom-refl}
\end{figure}

\begin{expl} \thlabel{expl:lie-hom-e6hasse} The first few columns of the Hasse diagram of the parabolic $\fp \leq \fg = \alg[_6]{e}{\bC}$ from \thref{expl:lie-para-e6wt} are given in the first diagram of Figure \ref{fig:lie-hom-e6expl}, where the edges are labelled by the node number corresponding to the simple reflection. \end{expl}

A careful analysis of the weight condition $\norm{\lambda + \lowest} = \norm{\mu + \lowest}$ now leads to Kostant's version of the Bott--Borel--Weil theorem \cite[Thm.\ 5.14]{k1961-liehom}, which calculates the irreducible components of $\liehom{}{V}$ using the Hasse diagram of $\fp$.  For notational convenience, we define the \emph{affine action} of $\sW_{\fg}$ on weights by $w\acts\mu \defeq w(\mu+\lowest) - \lowest$.

\begin{thm} \thlabel{thm:lie-hom-bbw} Let $V$ a complex irreducible $\fg$-representation with highest weight $\lambda$ and consider the $\mu$-isotypical component $\liehom{}{V}^{\mu}$ of $\liehom{}{V}$.  Then:
\begin{enumerate}
  \item $\liehom{}{V}^{\mu} \neq 0$ if and only if $\mu = w\acts\lambda$ for some $w \in \sW^{\fp}$.
  
  \item Each isotypical component $\liehom{}{V}^{w\acts\lambda}$ is irreducible, and the multiplicity of $w\acts\lambda$ as a weight of $\liechain{}{V}$ is one.
  
  \item The component $\liehom{}{V}^{w\acts\lambda}$ is contained in $\liehom{\ell(w)}{V}$.  \noproof
\end{enumerate}
\end{thm}

Therefore the set of $\fp^0$-irreducible components of $\liehom{}{V}$ is in bijection with the vertex set of $\sW^{\fp}$.  The calculation of the $\fp^0$-irreducible components of $\liehom{k}{V}$ is now completely algorithmic: first, determine the Hasse diagram of $\fp$ up to the $(k+1)$st column; for each weight in this column, take a sequence leading back to $\lowest[\fp]$ labelled from left-to-right by simple roots $\alpha_{i_1}, \ldots, \alpha_{i_k}$; then the corresponding component of $\liehom{k}{V}$ has highest weight given by the affine action $\alpha_{i_1} \acts (\cdots \alpha_{i_k} \acts \lambda) = (\alpha_{i_1} \cdots \alpha_{i_k}) \acts \lambda$ applied to the $\fp$-dominant weight induced by $\lambda$.

\begin{figure}[t]
  \newcommand{\ehassediagramexpl}[9]{  
    \arraycolsep=0.2em
    \hspace{1em}
    \begin{array}{c} \dynkinEp[#2]{6}{lie-hom-#1-hasse1} \end{array}
    \dynkin{ \DynkinConnector{0.2}{0.2}{1.3}{0.2}
             \DynkinLabel{$5$}{0.75}{-0.05} }{lie-hom-#1-hasse2}
    \begin{array}{c} \dynkinEp[#3]{6}{lie-hom-#1-hasse3} \end{array}
    \dynkin{ \DynkinConnector{0.2}{0.2}{1.3}{0.2}
             \DynkinLabel{$4$}{0.75}{-0.05} }{lie-hom-#1-hasse4}
    \begin{array}{c} \dynkinEp[#4]{6}{lie-hom-#1-hasse5} \end{array}
    \\[-2.25em]
    \zbox{ \dynkin{ \draw[thin]
      ( 7*\DynkinStep, 2*\DynkinStep) --                           ( 8*\DynkinStep, 2*\DynkinStep)
      ( 8*\DynkinStep, 2*\DynkinStep) arc(90:-90:1*\DynkinStep)    ( 8*\DynkinStep, 0*\DynkinStep)
      ( 8*\DynkinStep, 0*\DynkinStep) --                           (-8*\DynkinStep, 0*\DynkinStep)
      (-8*\DynkinStep, 0*\DynkinStep) arc(-90:90:-1.5*\DynkinStep) (-8*\DynkinStep,-3*\DynkinStep)
      (-8*\DynkinStep,-3*\DynkinStep) edge[->]                     (-7*\DynkinStep,-3*\DynkinStep);
      \DynkinLabel{$3$}{-7.5}{-3.25}
    }{lie-hom-#1-hasse6} }
    \hspace{5.8em}
    \\[-4.05em]
    \begin{array}{c} \dynkinEp[#5]{6}{lie-hom-#1-hasse7} \end{array}
    \hspace{-0.3ex}
    \dynkin{ \DynkinConnector{0.2}{0.3}{1.3}{1.3}
             \DynkinLabel{$2$}{0.55}{0.5}
             \DynkinConnector{0.2}{0.1}{1.3}{-0.9}
             \DynkinLabel{$6$}{0.55}{-1.5} }{lie-hom-#1-hasse8}
    \begin{array}{c} \dynkinEp[#6]{6}{lie-hom-#1-hasse9} \\[1em]
                     \dynkinEp[#7]{6}{lie-hom-#1-hasse10} \end{array}
    \hspace{-0.55ex}
    \dynkin[-0.18]{ \DynkinConnector{0.2}{ 1.52}{1.3}{1.52}
                    \DynkinLabel{$1$}{0.75}{1.25}
                    \DynkinConnector{0.2}{ 1.1}{1.3}{-0.7}
                    \DynkinLabel{$6$}{0.4}{-0.8}
                    \DynkinConnector{0.2}{-1.12}{1.3}{-1.12}
                    \DynkinLabel{$2$}{0.75}{-2.25} }{lie-hom-#1-hasse11}
    \begin{array}{c} \dynkinEp[#8]{6}{lie-hom-#1-hasse12} \\[1em]
                     \dynkinEp[#9]{6}{lie-hom-#1-hasse13} \end{array}
    \hspace{-0.8ex}
    \dynkin[-0.18]{ \DynkinConnector{0.2}{ 1.52}{1.3}{1.52}
                    \DynkinLabel{$6$}{0.75}{1.25}
                    \DynkinConnector{0.2}{-0.7}{1.3}{ 1.1}
                    \DynkinLabel{$1$}{0.35}{-0.35}
                    \DynkinConnector{0.2}{-1.12}{1.3}{-1.12}
                    \DynkinLabel{$3$}{0.75}{-2.25}
                    \DynkinLabel{$\cdots$}{2.2}{-1.8}
                    \DynkinLabel{$\cdots$}{2.2}{0.8} }{lie-hom-#1-hasse14}
    \hspace{1em} }
  \begin{multline*}
    \ehassediagramexpl{hasse}
                      {0,0,0,0,1,0}{0,0,0,1,-1,0}{0,0,1,-1,0,0}
                      {0,1,-1,0,0,1}{1,-1,0,0,0,1}{0,1,0,0,0,-1}
                      {-1,0,0,0,0,1}{1,-1,1,0,0,-1}
    \hspace{0.4em}
  \end{multline*}
  \vspace{-1em}
  \begin{multline*}
    \ehassediagramexpl{quabla}
                    {0,0,0,0, 0,1}{0,0,0,1,-2,1}{0,0,1,0,-3,1}
                    {0,1,0,0,-4,2}{1,0,0,0,-5,3}{0,3,0,0,-6,0}
                    {0,0,0,0,-6,4}{1,2,0,0,-7,1}
    \hspace{0.18em}
  \end{multline*}
  \vspace{-0.7em}
  \caption[Hasse diagram and homology of a parabolic subalgebra of $\mathfrak{e}_6(\bC)$]
          {\textsc{Top:} the Hasse diagram of the parabolic subalgebra of $\alg[_6]{e}{\bC}$ from \thref{expl:lie-hom-e6hasse}.  \textsc{Bottom:} The components of the homology of the adjoint representation of $\alg[_6]{e}{\bC}$ from \thref{expl:lie-hom-e6hom}, up to degree five.}
  \label{fig:lie-hom-e6expl}
\end{figure}

\begin{expl} \thlabel{expl:lie-hom-e6hom} Continuing notation from \thref{expl:lie-hom-e6hasse}, we can compute the homology $\liehom{}{\fg}$ valued in the adjoint representation of $\fg$.  The components of $\liehom{}{\fg}$ up to degree five are given in the corresponding column of the second diagram of Figure \ref{fig:lie-hom-e6expl}, where we have retained the arrow labelling from \thref{expl:lie-hom-e6hasse} purely for clarity. \end{expl}

Finally, it is straightforward to extend \thref{thm:lie-hom-bbw} to other cases of interest; we summarise the results from \cite[Prop.\ 3.3.6]{cs2009-parabolic1}.  Firstly for a family $\{V_i\}_{i\in I}$ of complex irreducible $\fg$-representations, there is a natural $\fp$-representation isomorphism $\liehom{k}{\Dsum{i\in I}{} V_i} \isom \Dsum{i\in I}{} \liehom{k}{V_i}$ for each degree $k$.  If on the other hand $\fg_1,\fg_2$ are complex semisimple Lie algebras with parabolic subalgebras $\fp_i \leq \fg_i$ and irreducible representations $V_i$, the external tensor product $V_1 \etens V_2$ is an irreducible representation of $\fg_1 \dsum \fg_2$ for which
\vspace{-0.2em}
\begin{equation*}
  \liehom[\fp_1^{\perp}\dsum\fp_2^{\perp}]{k}{V_1\etens V_2}
    \isom \Dsum{i+j=k}{} \,\big( \liehom[\fp_1^{\perp}]{i}{V_1} \etens
                              \liehom[\fp_2^{\perp}]{j}{V_2} \big)
\vspace{-0.2em}
\end{equation*}
as a representation of $\fp_1 \dsum \fp_2$.

If now $\fg$ is a real semisimple Lie algebra and $V$ a complex $\fg$-representation, it is easy to see that the complexification of the real homology $\liehom{k}{V}$ is naturally isomorphic to the complex homology $\liehom[\cpxrepn{\fp}^{\perp}]{k}{V}$ as a representation of $\fp \leq \cpxrepn{\fp}$.  If on the other hand $V$ is a real representation, there is a natural isomorphism $\cpxrepn{ \liehom{k}{V} } \isom \liehom[\cpxrepn{\fp}^{\perp}]{k}{\cpxrepn{V}}$ of $\cpxrepn{\fp}$-representations and one of two cases may occur.  Indeed, for simplicity assume that $V$ is irreducible and admits no $\fg$-invariant complex structure, and let $W\leq \liehom[ \cpxrepn{\fp}^{\perp} ]{k}{ \cpxrepn{V} }$ be a $\cpxrepn{\fp}$-irreducible subrepresentation.  Then either:
\begin{itemize}
  \item $W = \conj{W}$ is the complexification of a real irreducible component of $\liehom{k}{V}$; or
  
  \item $\conj{W}$ is an irreducible subrepresentation of $\liehom[ \cpxrepn{\fp}^{\perp} ]{k}{ \cpxrepn{V} }$, with $W\dsum\conj{W}$ the complexification of a single complex irreducible component in $\liehom{k}{V}$.
\end{itemize}
In either case, no irreducible component of the real homology $\liehom{k}{V}$ admits a \qtn ic\ structure.  We shall see applications of these results in later chapters.

\chapter[Background from parabolic geometry]{Background from \\ parabolic geometry} 
\label{c:para}

\BufferDynkinLocaltrue
\renewcommand{\dynkinnameoffset}{-0.75}

\vspace{0.6em}
Having described the structure theory of parabolic subalgebras, we turn now to Cartan geometries and parabolic geometries.  We shall assume that the reader is familiar with the basic concepts of differential geometry, such as vector bundles, principal bundles and principal connections; see
\cite{b1987-einstein,e1997-riemannian,kn1996-founddg1,kn1996-founddg2,l1997-riemcurv,l2003-smoothmfds,s1997-cartan}
for comprehensive accounts.  Unless stated otherwise, all objects will be assumed to be smooth and all principal bundles carry right actions.  Given a manifold $M$ and a vector bundle $E \surjto M$, the space of sections of $\Wedge{k}T^*M \tens E$ will be denoted by $\s{k}{E}$.

Intuitively, a Cartan geometry is a curved analogue of a homogeneous space, while a parabolic geometry is a Cartan geometry modelled on a generalised flag manifold.  Then the theory of parabolic subalgebras developed in Chapter \ref{c:lie} imbues a parabolic geometry with a rich algebraic structure, which (in most circumstances) can be exploited to obtain an equivalence of categories between parabolic geometries of a certain type and simpler underlying geometric structures; we describe this in Section \ref{s:para-para}.  There is also a well-developed theory of invariant differential operators on parabolic geometries, which will be important for us in later chapters; we describe this theory in Section \ref{s:para-bgg}.

\section{Cartan geometries and parabolic geometries} 
\label{s:para-para}

We begin by reviewing the basic theory of Cartan geometries and parabolic geometries in Subsections \ref{ss:para-para-cartan} and \ref{ss:para-para-para}, from the modern perspective of principal bundles and principal connections.  This differs from Cartan's original approach \cite{c1924-projaffine,c1924-projconns,c1926-confproj} which was phrased in terms of gauge transformations \cite[\S 5.1]{s1997-cartan}.

\subsection{Cartan geometries} 
\label{ss:para-para-cartan}

Let $G$ be a (real) Lie group with Lie subgroup $P\leq G$, and let $\fp\leq\fg$ be their Lie algebras.  The left-invariant vector fields on $G$ induce a naturally defined trivialisation $TG \isom G \times \fg$, whose inverse can be conveniently encoded in a canonical $\fg$-valued $1$-form $\omega^G \in \s[G]{1}{\fg}$ called the \emph{Maurer--Cartan form} of $G$, defined by $\omega^G_g(\xi) \defeq L_{g^{-1}*}(\xi)$ for all $g\in G$ and $\xi \in T_g G$.  Clearly $\omega^G$ reproduces the generators of left-invariant vector fields and defines an isomorphism $T_g G \isom \fg$ for each $g \in G$, and moreover we have the \emph{Maurer--Cartan equation} $\d\omega^G + \tfrac{1}{2} \liebracw{\omega^G}{\omega^G} = 0$.  Viewing this as a ``zero curvature'' condition on the canonical principal $P$-bundle $G \surjto G/P$,%
\footnote{But note that $\omega^G$ is not a principal $P$-connection, since it is $\fg$-valued rather than $\fp$-valued.}
a Cartan geometry is a curved geometry modelled locally on the homogeneous space $G/P$.

\begin{defn} \thlabel{defn:para-para-cartan} A \emph{Cartan geometry} $(F^P \surjto M,\, \omega)$ of type $G/P$ on a manifold $M$ of dimension $\dim M = \dim(G/P)$ is a principal $P$-bundle $p : F^P \surjto M$, called the \emph{Cartan bundle}, equipped with a $\fg$-valued $1$-form $\omega\in\s[F^P]{1}{\fg}$, called the \emph{Cartan connection}, such that:
\begin{enumerate}
  \item \label{defn:para-para-cartan-Pinv}
  $\omega$ is $P$-invariant, \ie\ $R_g^*\omega = \Ad(g^{-1})\omega$ for all $g\in P$;
  
  \item \label{defn:para-para-cartan-fund}
  $\omega(X^{\xi}) = \xi$ for all fundamental vector fields $X^{\xi} \in \s[F^P]{0}{TF^P}$ with $\xi \in \fp$; and
  
  \item \label{defn:para-para-cartan-cc}
  $\omega_u$ defines a linear isomorphism $T_u F^P \to \fg$ at each $u \in F^P$.
\end{enumerate}
The homogeneous space $G/P$ with its Maurer--Cartan form $\omega^G \in \s[G]{1}{\fg}$ is called the \emph{flat model} of the Cartan geometry.  We refer to item \ref{defn:para-para-cartan-cc} as the \emph{Cartan condition}. \end{defn}

A principal $P$-bundle $F^P \surjto M$ does not determine a unique Cartan condition.  Indeed, the possible Cartan connections on $F^P$ form an open subset of an affine space modelled on the $P$-invariant horizontal subspace of $\s[F^P]{1}{\fg}$; see \cite[\S 1.5.2]{cs2009-parabolic1}.

\begin{expl} \thlabel{expl:para-para-cexpl}
\begin{enumeratepar}
  \item \label{expl:para-para-cexpl-euc}
  Let $G$ be the euclidean group $\grp{O}{n} \ltimes \bR^n$, and let $P=\grp{O}{n}$.  Since a $P$-structure is a \riem\ metric $g$, a Cartan geometry of type $G/P$ is equivalent \cite[\S 6.3]{s1997-cartan} to a principal $P$-connection on the orthonormal frame bundle of $g$, and hence a metric connection on $TM$.  If \torsionfree, this connection coincides with the \LC\ connection of $g$.
  
  \item \label{expl:para-para-cexpl-conf}
  Let $G=\grp{SO}{n+1,1}$ be the Lorentz group.  Inside $\bR^{n+1,1}$ is the light-cone of \non zero\ null vectors, whose projectivisation is the conformal $n$-sphere $\Sph[n]$.  The action of $G$ on $\bR^{n+1,1}$ preserves the light-cone, so descends to an action on $\Sph[n]$ which identifies $G$ with the M{\"o}bius group of conformal transformations of $\Sph[n]$.  The pointwise stabiliser of this action is a subgroup $P$ isomorphic to $\grp{CO}{n}\ltimes \bR^{n*}$, which identifies $\Sph[n]$ with the homogeneous space $G/P$.  Using the lorentzian metric on $\bR^{n+1,1}$, a Cartan geometry of type $G/P$ induces a conformal connection on $TM$.  A careful treatment of conformal geometry as a Cartan geometry may be found in \cite{bc2010-conf,cs2009-parabolic1,s1997-cartan}.
\end{enumeratepar}
\end{expl}

Suppose that $F\surjto M$ is any principal $P$-bundle.  The following definition is standard but of vital importance to later developments.

\begin{defn} \thlabel{defn:para-para-assoc} Let $V$ be a $P$-representation.  The \emph{associated bundle} $V_M \defeq \assocbdl{F}{P}{V}$ is the quotient of $F\times V$ by the right $P$-action defined by $(u,v) \acts g \defeq (u \acts g, g^{-1} \acts v)$. \end{defn}

Then $V_M$ is a vector bundle over $M$ with standard fibre $V$.  Moreover the map $\s[F]{0}{V}^P \to \s{0}{V_M}$ given by mapping a $P$-equivariant function $f:F\to V$ to the section $s$ of $V_M$ defined by $s(p(u)) \defeq [u,f(u)]$, the class of $(u,f(u))$ in $V_M$, is an isomorphism.  We may then think of sections of $V_M$ as $P$-equivariant functions $F \to V$.

For a Cartan geometry $(F^P \surjto M,\, \omega)$, the Cartan condition allows us to identify many geometric bundles with bundles associated to $F^P$, hence linking the geometry of $F^P$ with the geometry of $M$.  As a fundamental example, item \ref{defn:para-para-cartan-cc} above determines a trivialisation $TF^P \isom F^P \times \fp$, while \ref{defn:para-para-cartan-fund} implies that the vertical bundle of $F^P \surjto M$ may be identified with $\fp_M = \assocbdl{F^P}{P}{\fg}$.  In this picture the fundamental vector fields generating the $P$-action are $X^{\xi} = \omega^{-1}(\xi)$ for $\xi\in\fp$, so the natural trivialisation of the vertical bundle is provided by the \emph{constant vector fields} $\omega^{-1}(\xi)$ for $\xi\in\fp$; on the flat model these fields are just the left-invariant vector fields on $F^P = G$.  Moreover, differential forms on $F^P$ are determined uniquely by their values on the $\omega^{-1}(\xi)$.

Items \ref{defn:para-para-cartan-Pinv} and \ref{defn:para-para-cartan-cc} of \thref{defn:para-para-cartan} also imply that $\omega_u : T_u F^P \to \fg$ descends to an isomorphism $\omega_u \bmod \fp : T_{p(u)}M \to \fg/\fp$ for each $u \in F^P$, thus identifying $TM$ with the associated bundle $(\fg/\fp)_M = \assocbdl{F^P}{P}{\fg/\fp}$.  This means that $M$ inherits the ``first order'' geometry of $G/P$.  By functoriality of the associated bundle construction, this also identifies all tensor bundles with associated bundles.


\begin{defn} \thlabel{defn:para-para-curv} The \emph{curvature} of a Cartan geometry $(F^P \surjto M,\, \omega)$ on $M$ is the $\fg$-valued $2$-form
\begin{equation*}
  K \defeq \d\omega + \tfrac{1}{2} \liebracw{\omega}{\omega} \in \s[F^P]{2}{\fg}.
\end{equation*}
A Cartan geometry is \emph{flat} if its curvature vanishes identically. \end{defn}

Via the isomorphism between sections of associated bundles and $P$-equivariant functions, $K$ induces a \emph{curvature function} $\kappa : F^P \to \Wedge{2}\fg^* \tens \fg$ defined by
\begin{align*}
  \kappa_u(\xi,\eta) &= K(\omega_u^{-1}(\xi), \omega_u^{-1}(\eta)) \\
    &= \liebrac{\xi}{\eta} - \omega_u\left( \liebrac{\omega_u^{-1}(\xi)}{\omega_u^{-1}(\eta)} \right).
\end{align*}
Thus the curvature $K$ is the obstruction to $\omega_u$ defining a Lie algebra homomorphism $T_u F^P \to \fg$.  The $P$-invariance of $\omega$ implies that $\xi \mapsto \omega^{-1}(\xi)$ is $P$-equivariant, so that differentiating gives $\liebrac{\omega^{-1}(\xi)}{\omega^{-1}(\eta)} = \omega^{-1}(\liebrac{\xi}{\eta})$.  It follows that $K$ is $P$-invariant and horizontal (that is, $X \intprod K = 0$ for any vertical vector field on $F^P$), so may be viewed as a $2$-form $K_M \in \s{2}{\fg_M}$ on $M$; equivalently $\kappa$ takes values in $\Wedge{2}(\fg/\fp)^* \tens \fg$.  In the case of \itemref{expl:para-para-cexpl}{euc}, $K_M$ coincides with the curvature of the metric connection.

\begin{prop} \thlabel{prop:para-para-flat} \emph{\cite[Thm.\ 5.5.1]{s1997-cartan}} A Cartan geometry $(F^P \surjto M,\, \omega)$ is flat if and only if it is locally isomorphic%
\footnote{That is, there is a principal $P$-bundle isomorphism $\Psi : G/P \to F^P$ such that $\Psi^*\omega = \omega^G$.}
to its flat model $(G \surjto G/P,\, \omega^G)$. \noproof \end{prop}

Clearly the restriction $(F^P\at{U} \surjto U,\, \omega\at{U})$ to an open set $U\subseteq M$ is also a Cartan geometry of type $G/P$.  Then \thref{prop:para-para-flat} means that every point in $M$ has a neighbourhood $U$ such that $(F^P\at{U} \surjto U,\, \omega\at{U})$ is isomorphic to the restriction of $(G \surjto G/P,\, \omega^G)$ to a neighbourhood of $0 \in G/P$.

\begin{defn} \thlabel{defn:para-para-tor} The \emph{torsion} of $\omega$ is the $\fg/\fp$-valued $2$-form $T \in \s[F^P]{2}{\fg/\fp}$ defined by projecting values of $K$ to $\fg/\fp$.  A Cartan geometry is \emph{\torsionfree} if its torsion vanishes. \end{defn}

Since $K$ is $P$-invariant and horizontal, the torsion descends to a $2$-form $T_M \in \s{2}{TM}$.  In \thref{expl:para-para-cexpl}, \torsionfree ness\ amounts to \torsionfree ness\ of the metric connection or conformal connection induced on $TM$.

\subsection{Parabolic geometries} 
\label{ss:para-para-para}

Due to the algebraic properties of parabolic subalgebras discussed in Section \ref{s:lie-para}, Cartan geometries modelled on generalised flag manifolds have a rich algebraic structure.

\begin{defn} \label{defn:para-para-para} A \emph{parabolic geometry} is a Cartan geometry modelled on a generalised flag manifold $G/P$, with $G$ a semisimple Lie group and $P\leq G$ a parabolic subgroup. \end{defn}

We retain any adjectives applied to either the Cartan geometry or the parabolic subalgebra $\fp$ of $P$.  In particular, a parabolic geometry is \emph{abelian} if $\fp$ is an abelian parabolic; for example, conformal geometry from \itemref{expl:para-para-cexpl}{conf} is abelian.

A natural question for Cartan geometries is whether they are equivalent to simpler underlying first-order structures.  For conformal geometry, the \emph{conformal equivalence problem} states that a conformal manifold is equivalent to a certain kind of conformal Cartan geometry \cite[\S 7.3]{s1997-cartan}.  We shall discuss this further in Subsection \ref{ss:para-calc-equiv}, but as a first step we can obtain the underlying geometric structure.

Let $(F^P \surjto M,\,\omega)$ be a parabolic geometry of type $G/P$ on $M$, and suppose that $\fp$ has height $n$.  Via the trivialisation $TF^P \isom F^P\times\fg$ induced by $\omega$, the $\fp^{\perp}$-filtration of $\fg$ induces a filtration
\begin{equation} \label{eq:para-para-TFPfilt}
  TF^P = TF^P_n \supset \cdots \supset TF^P_{-n} \supset 0
\end{equation}
by smooth subbundles, where $TF^P_i \defeq \omega^{-1}(\fg_i)$ has rank equal to $\dim(\fg_i)$.
Since $\exp\fp^{\perp} \leq P$ acts freely on $F^P$, the orbit space $F^0 \defeq F^P/\exp\fp^{\perp}$ is a principal $P^0$-bundle.  Moreover $\fp$ preserves the filtration \eqref{eq:para-para-TFPfilt}, hence inducing filtrations $TF^0 = TF^0_n \supset \cdots \supset TF^0_0 \supset 0$ and $TM = TM_n \supset \cdots \supset TM_1 \supset 0$ of $TF^0$ and $TM$ by smooth subbundles.  For convenience we choose an algebraic Weyl structure for $\fp$, thus identifying $\fg \isom \gr{\fg}$ and splitting the natural projection $P \surjto P^0$.  For each $u_0 \in F^0$ and $\xi \in T_{u_0}F^0_i$, choose lifts $u \in F^P$ of $u_0$ and $\tilde{\xi} \in T_uF^P$ of $\xi$.  It turns out \cite[Prop.\ 3.1.5]{cs2009-parabolic1} that
\begin{equation*}
  \omega^0_{(i)}(\xi) \defeq \omega( \tilde{\xi} ) \bmod \fg_{i-1}
\end{equation*}
is a well-defined smooth section of $\Hom{TF^0_i}{\fg_{(i)}}$, which is $P^0$-equivariant and has kernel $TF^0_{i-1}$.  The graded vector bundle $\gr{(TM)}$ is associated to $\fg/\fp \isom \fg_{(1)} \dsum \cdots \dsum \fg_{(n)}$ and the adjoint action of $P^0$ on $\fg/\fp$ preserves this grading, giving a homomorphism $\Ad : P^0 \to \grp{GL}{\fg/\fp}_{\mr{gr}}$ into the group of grading-preserving automorphisms of $\fg/\fp$.  This allows us to talk about reductions of $\gr{(TM)}$ to structure group $P^0$.

The data $(TM_i,\, \omega^0_{(i)})$ constitutes an \emph{infinitesimal flag structure} of type $G/P$ on $M$.  In fact, these data are equivalent to a filtration $TM = TM_n \supset \cdots \supset TM_1 \supset 0$ in which the rank of $TM_i$ equals $\dim(\fg_i/\fp)$, together with a reduction of $\gr{(TM)}$ to structure group $P^0$ \wrt\ $\Ad : P^0 \to \grp{GL}{\fg/\fp}_{\mr{gr}}$; see \cite[Prop.\ 3.1.6]{cs2009-parabolic1}.  For abelian parabolic geometries this amounts to a reduction of $TM$ to structure group $P^0$.

\begin{expl} \thlabel{expl:para-para-flag} In the case of conformal geometry, $G = \grp{SO}{n+1,1}$ and $P$ is the stabiliser of a given isotropic line.  The choice of an algebraic Weyl structure and the adjoint action induces an isomorphism $P^0 \isom \grp{CO}{n}$, so that an infinitesimal flag structure of type $G/P$ is a first-order $\grp{CO}{n}$-structure, \ie\ a conformal structure. \end{expl}

\smallskip
\section{Tractor calculus and the equivalence of categories} 
\label{s:para-calc}

It is often beneficial to ``linearise'' the Cartan connection, by inducing linear connections on associated \emph{tractor bundles}; we discuss this in Subsection \ref{ss:para-calc-tractor}.  Geometric analogues of algebraic Weyl structures are introduced in Subsection \ref{ss:para-calc-weyl}, after which we conclude the discussion of the parabolic equivalence problem in Subsection \ref{ss:para-calc-equiv}.

\subsection{Tractor bundles and tractor connections} 
\label{ss:para-calc-tractor}

Let $F \surjto M$ be a principal $P$-bundle with principal $P$-connection $\theta \in \s[F]{1}{\fp}$, and let $V$ be a $P$-representation.  If $\cH \leq TF$ is the horizontal distribution of $\theta$, then $\cH \times \{0\} \leq TF \times V$ is $P$-invariant and therefore descends to a horizontal distribution on the associated bundle $V_M \defeq \assocbdl{F}{P}{V}$, so that $\theta$ induces a linear connection on each associated bundle.  In the case of a Cartan geometry $(F^P \surjto M ,\, \omega)$, the Cartan connection $\omega$ is not a principal $P$-connection so cannot be used to induce connections on associated bundles.  We can remedy this by forming the \emph{extended Cartan bundle} $F^G \defeq \assocbdl{F^P}{P}{G}$ associated to the restriction of the adjoint action of $G$ on itself.

\begin{prop} \rmcite[Thm.\ 1.5.6]{cs2009-parabolic1} The extended Cartan bundle $F^G$ is a principal $G$-bundle, with a unique principal $G$-connection induced by the Cartan connection $\omega$. \noproof \end{prop}

We shall also denote this principal $G$-connection by $\omega \in \s[F^G]{1}{\fg}$.  Then for any $G$-representation $\bV$, the associated vector bundle $\cV \defeq \assocbdl{F^G}{G}{\bV}$ inherits a linear connection from $\omega$.  Clearly restriction to $P$ give an isomorphism $\cV \isom \assocbdl{F^P}{P}{\bV}$, so that $\cV$ may be viewed as a bundle associated to $F^P$.

\begin{defn} Let $\bV$ be the restriction to $P$ of a $G$-representation.  Then $\cV \defeq \assocbdl{F^P}{P}{\bV}$ is the \emph{tractor bundle} associated to $\bV$, while the linear connection $\D^{\bV}$ induced by $\omega$ is the \emph{tractor connection}. \end{defn}

In order to obtain a formula for the tractor connection on each tractor bundle, we define the \emph{invariant derivative} by
\begin{equation} \label{eq:para-calc-invd} \begin{aligned}
  \D^{\omega} : \s[F^P]{0}{\bV} &\to \s[F^P]{0}{\fg^*\tens\bV} \\
  \D^{\omega}_{\xi} s &= \d s(\omega^{-1}(\xi))
\end{aligned} \end{equation}
for all $\xi\in\fg$, which evidently depends only on the Cartan connection $\omega$ and the representation $\bV$.  It is straightforward to see that if $f$ is $P$-equivariant then so is $\D^{\omega}f$, so that we may view $\D^{\omega} : \s{0}{\cV} \to \s{0}{\fg_M^*\tens\cV}$ as a map of sections.  Given $\xi,\eta \in \fg$, one can verify by direct calculation that $\D^{\omega}$ satisfies a Leibniz rule
\begin{equation*}
  \D^{\omega}_{\xi} (s_1 \tens s_2)
    = \D^{\omega}_{\xi} s_1 \tens s_2 + s_1 \tens \D^{\omega}_{\xi} s_2
\end{equation*}
on tensor products $\bV_1 \tens \bV_2$, and a Ricci identity
\begin{equation} \label{eq:para-calc-ricci}
  \D^{\omega}_{\xi} (\D^{\omega}_{\eta} s) - \D^{\omega}_{\eta} (\D^{\omega}_{\xi} s)
    = \D^{\omega}_{\liebrac{\xi}{\eta}} s - \D^{\omega}_{\kappa(\xi,\eta)} s,
\end{equation}
where $\kappa : F^P \to \Wedge{2} \fg_M^* \tens \fg_M$ is the curvature function of $\omega$.  The final term in \eqref{eq:para-calc-ricci} is first order in general due to torsion, and if $\omega$ is \torsionfree\ then $\D^{\omega}_{\liebrac{X}{Y}} s = -\kappa(X,Y) \acts s$.  Moreover the map $\s{0}{\cV} \ni s \mapsto (s,\D^{\omega}s) \in \s{0}{\cV \dsum (\fg^*_M\tens\cV)}$ identifies the first jet bundle $\jetbdl{1}{\cV}$ of $\cV$ with an associated bundle \cite[Prop.\ 1.3]{cd2001-curvedbgg}.

\begin{lem} \thlabel{lem:para-calc-verttriv} $\D^{\omega}$ is vertically trivial in the sense that $\D^{\omega}_{\xi} f + \xi \acts f = 0$ for all $P$-equivariant functions $f:F^P\to\bV$ and all $\xi\in\fp$. \noproof \end{lem}

The modified map $(s,\xi) \mapsto \D^{\omega}_{\xi} s + \xi\acts s$ then vanishes for $\xi\in\fp$, so takes values in $(\fg/\fp)^*\tens\bV$ and therefore defines a linear connection on $\bV$ by the Cartan condition.

\begin{prop} \thlabel{prop:para-calc-tractor} The linear connection induced by $(s,\xi) \mapsto \D^{\omega}_{\xi} s + \xi\acts s$ coincides with the tractor connection $\D^{\bV}$ on $\cV$. \noproof \end{prop}

Direct calculation using \thref{prop:para-calc-tractor} and equation \eqref{eq:para-calc-invd} yields the following.

\begin{cor} \thlabel{cor:para-calc-tractorcurv} The curvature of $\D^{\bV}$ is given by $\Curv[\bV]{X,Y}{s} = K_M(X,Y) \acts s$ for all $X,Y \in \s{0}{TM}$ and $s \in \s{0}{\cV}$. \noproof \end{cor}

An important example of a tractor bundle is the \emph{adjoint tractor bundle} $\fg_M \defeq \assocbdl{F^P}{P}{\,\fg}$ given by the restriction of the adjoint action of $G$ on $\fg$.  The projection $\fg \surjto \fg/\fp$ exhibits $TM$ as a natural quotient of $\fg_M$, while dually the inclusion $(\fg/\fp)^* \injto \fg^*$ exhibits $T^*M$ as a natural subbundle of $\fg^*_M$.  For parabolic geometries, $TM$ and $T^*M$ are bundles of nilpotent Lie algebras modelled on $\fg/\fp$ and $\fp^{\perp}$ respectively.

Viewed as a $P$-equivariant bilinear map $\liebrac{}{} : \fg\times\fg \to \fg$, the Lie bracket on $\fg$ induces a bilinear map $\algbrac{}{} : \s{0}{\fg_M} \times \s{0}{\fg_M} \to \s{0}{\fg_M}$ of sections.

\begin{defn} \thlabel{defn:para-calc-bracs} The map $\algbrac{}{} : \s{0}{\fg_M} \times \s{0}{\fg_M} \to \s{0}{\fg_M}$ of sections of $\fg_M$ is called the \emph{algebraic bracket}. \end{defn}

Note also that the pointwise action of $\fg$ on $\bV$ extends to a map $\bdot : \s{0}{\fg_M} \times \s{0}{\cV} \to \s{0}{\cV}$ of sections of the associated bundles, called the \emph{algebraic action}; this identifies $\cV$ with a bundle of representations for the Lie algebra bundle $\fg_M$.

\begin{prop} Let $\bV$ be a $\fg$-representation with tractor bundle $\cV$.  Then:
\begin{enumerate}
  \item $\algbrac{s_1}{s_2} \acts t = s_1 \acts (s_2 \acts t) - s_2 \acts (s_1 \acts t)$; and
  
  \item The algebraic bracket and algebraic action are tractor-parallel, \ie\
    \begin{align*}
    \D^{\fg}_X \algbrac{s_1}{s_2}
      &= \algbrac{ \D^{\fg}_X s_1 }{ s_2 } + \algbrac{ s_1 }{ \D^{\fg}_X s_2 } \\
    \text{and}\quad
    \D^{\bV}_X (s \acts t)
      &= (\D^{\fg}_X s) \acts t + s \acts (\D^{\bV}_X t)
  \end{align*}
\end{enumerate}
for all $s,s_1,s_2 \in \s{0}{\fg_M}$ and $t \in \s{0}{\cV}$. \noproof \end{prop}

\subsection{Weyl structures} 
\label{ss:para-calc-weyl}

Recall that an algebraic Weyl structure is a choice of lift of the unique grading element $\xi^0 \in \liecenter{\fp^0}$ to $\fp$, and this choice splits the $\fp^{\perp}$-filtration of $\fg$.  Since $\xi^0$ is unique and $P$-invariant, we obtain a unique \emph{grading section} $\xi^0_M$, which is a section of $\fp^0_M \defeq \assocbdl{F^P}{P}{\fp^0}$.

\begin{defn} \thlabel{defn:para-calc-weyl} A \emph{Weyl structure} is a smooth lift of the grading section $\xi^0_M \in \s{0}{\fp^0}$ to a section of $\fp_M \defeq \assocbdl{F^P}{P}{\fp}$. \end{defn}

Equivalently, a Weyl structure is a smooth choice of an algebraic Weyl structure at each point of $M$.  By \thref{lem:lie-para-exp}, the space of Weyl structures is an affine bundle modelled on $\fp^{\perp}_M \isom T^*M$.

Applying the results of Subsection \ref{ss:lie-para-repn} pointwise, a Weyl structure induces isomorphisms $\xi_{\cV} : \gr\cV \to \cV$ of each associated bundle with its associated graded bundle; in particular we obtain an isomorphism $\xi_{\fg} : TM \dsum \fp^0_M \dsum T^*M \to \fg_M$.  This allows us to use the invariant derivative \eqref{eq:para-calc-invd} to induce linear connections on tractor bundles.

\begin{defn} Let $\bV$ be a $\fg$-representation with associated bundle $\cV \defeq \assocbdl{F^P}{P}{\bV}$.
\begin{enumerate}
  \item The \emph{Ricci-corrected Weyl connection} is the linear connection $\D^{(1)}$ on $\cV$ defined by restricting the invariant derivative $\D^{\omega}$ to $TM$ using $\xi_{\fg}$, \ie\ $\D^{(1)}_X s \defeq \D^{\omega}_{\xi_{\fg}X} s$.
  
  \item The \emph{Weyl connection} is the linear connection $\D$ on $\cV$ induced by $\D^{(1)}$ on $\gr\cV$, \ie\ $\D s \defeq (\xi_{\cV}^{-1} \circ \D^{(1)} \circ \xi_{\cV})s$ for all $s \in \s{0}{\cV}$.
\end{enumerate}
\end{defn}

If $\bV$ is a completely reducible $P$-representation then $\exp\fp^{\perp}$ acts trivially, so that $\cV \isom \gr\cV$ canonically.  In this case the connections $\D^{(1)}$ and $\D$ agree.

We may equivalently view $\xi_M$ as an isomorphism $F^P \times (\fg/\fp \dsum \fp^0 \dsum \fp^{\perp}) \to F^P \times \fg$, also denoted by $\xi_M$.  Using this we may decompose the Cartan connection $\omega$ as
\begin{equation} \label{eq:para-calc-cconn}
  \omega = \xi_M\omega_{\fg/\fp} \dsum \omega_{\fp},
  \quad\text{where}\quad
  \omega_{\fp} = \xi_M\omega_{\fp^0} \dsum \omega_{\fp^{\perp}},
\end{equation}
$\omega_{\fg/\fp} \defeq \omega \bmod \fp \in \s{1}{\fg/\fp}$ is the \emph{solder form}, and $\omega_{\fp^0} \defeq \omega_{\fp} \bmod \fp^{\perp} \in \s{1}{\fp^0}$.  Since the space $\fw$ of algebraic Weyl structures is a homogeneous space for $\exp\fp^{\perp}$, the fundamental vector fields generated by elements of $\fp^{\perp}$ induce a Maurer--Cartan form $\eta : T\fw \to \fp^{\perp}$, trivialising $T\fw \isom \fw\times\fp^{\perp}$ via the affine structure.  We may then view $\xi_M$ as a $P$-equivariant function $\xi_{\fw} : F^P \to \fw$.  It turns out to be fruitful to write
\begin{equation} \label{eq:para-calc-cconnp}
  \omega_{\fp} = (\xi_M\omega_{\fp^0} - \xi_{\fw}^*\eta)
    + (\omega_{\fp^{\perp}} + \xi_{\fw}^*\eta)
\end{equation}
for the following reasons \cite[Prop.\ 4.2]{cds2005-ricci}.

\begin{prop} \thlabel{prop:para-calc-weylconn} Let $\xi_M$ be a Weyl structure.  Then:
\begin{enumerate}
  \item \label{prop:para-calc-weylconn-D1}
  $\omega_{\fp}$ is a principal $P$-connection on $F^P$, inducing $\D^{(1)}$ on associated bundles.
  
  \item \label{prop:para-calc-weylconn-D}
  $\xi_M\omega_{\fp^0} - \xi_{\fw}^*\eta$ is a principal $P$-connection on $F^P$, inducing $\D$ on associated bundles.
  
  \item \label{prop:para-calc-weylconn-r}
  $\omega_{\fp^{\perp}} + \xi_{\fw}^*\eta$ is a horizontal $P$-invariant $\fp^{\perp}$-valued $1$-form on $F^P$; if $\nRic{}{}$ is the $T^*M$-valued $1$-form induced on $M$ then $\D_X^{(1)} s = \D_X s + \nRic{X}{}\acts s$. \noproof
\end{enumerate}
\end{prop}

\begin{defn} \thlabel{defn:para-calc-nric} The $T^*M$-valued $1$-form $\nRic{}{} \in \s{1}{T^*M}$ induced by $\omega_{\fp^{\perp}} + \xi_{\fw}^*\eta$ on $M$ is called the \emph{normalised Ricci tensor} of the Weyl structure $\xi_M$. \end{defn}

Since $\D^{(1)}$ is induced by the invariant derivative on $\cV$, \thref{lem:para-calc-verttriv} and \itemref{prop:para-calc-weylconn}{r} allow us to write the tractor connection on $\cV$ as
\begin{equation} \label{eq:para-calc-DV}
  \D^{\bV}_X s = X\acts s + \D_X s + \nRic{X}{}\acts s.
\end{equation}
Using the Weyl structure, the Lie algebra differential $\liediff : \Wedge{k}\fp^{\perp} \tens \bV \to \Wedge{k+1}\fp^{\perp} \tens \bV$ induces a bundle map $\liediff : \Wedge{k}T^*M\tens\cV \to \Wedge{k+1}T^*M\tens\cV$.  Equation \eqref{eq:para-calc-DV} then becomes $\D^{\bV} s = \liediff s + \D s + \nRic{}{} \acts s$, as in \cite{cg2002-tractorcalc,hsss2010-prolconns}.

\begin{rmk} If we choose a point $x\in M$ and hence an algebraic Weyl structure, a Weyl structure is equivalent to a $P^0$-equivariant section of the projection $F^P \surjto F^0 \defeq F^P / \exp(\fp^{\perp})$; see \cite[App.\ A]{cds2005-ricci}.  This is the original approach of {\v C}ap and Slov{\'a}k \cite{cs2003-weylstr}. \end{rmk}

It is natural to ask how the Weyl connection $\D$ and normalised Ricci tensor $\nRic{}{}$ change when we change the Weyl structure $\xi_M$.  Viewing $\xi_M$ as the $P$-equivariant function $\xi_{\fw} : F^P \to \fw$ from above, any other Weyl structure may be written $\xi_M' = (\Ad q) \xi_M$ for some $P$-equivariant function $q : F^P \to \exp{\fp^{\perp}}$.  Let $q(t) : F^P \to \exp{\fp^{\perp}}$ be a $P$-equivariant curve with $q(0) = \id$ and $q'(0) = \gamma \in \s{1}{}$.  Then for any object $F(\xi_M)$ depending on the Weyl structure, we define the \emph{first-order variation of $F$} by
\begin{equation*} 
  (\weyld{\gamma} F)(\xi_M) \defeq
    \frac{\d}{\d t}\, F\!\left( (\Ad(q(t)^{-1}))\,\xi_M \right) \Big|_{t=0}
    \colvectpunct[-0.5em]{~.}
\end{equation*}
By the fundamental theorem of calculus, $F$ is independent of $\xi_M$ if and only if $\weyld{\gamma} F = 0$ for all $\gamma \in \s{1}{}$.  Moreover for any $\gamma \in \s{1}{}$, Taylor's theorem allows us to write
\begin{equation} \label{eq:para-calc-taylor}
  F \big( (\Ad \gamma)\xi_M \big) = F(\xi_M) + (\weyld{\gamma} F)(\xi_M)
    + \tfrac{1}{2}(\weyld{\gamma}\weyld{\gamma} F)(\xi_M) + \cdots,
\end{equation}
where we view $\gamma$ as a $P$-equivariant function $F^P \to \fp^{\perp}$.

\begin{prop} \thlabel{prop:para-calc-change} \emph{\cite[App.\ B]{cds2005-ricci}} Let $\xi_M$ be a Weyl structure and let $\gamma \in \s{1}{}$.  Then:
\begin{enumerate}
  \item \label{prop:para-calc-change-D}
  $\weyld{\gamma} \D_X = \algbrac{X}{\gamma}_{\fp^0} + \D_X \gamma$; and
  
  \item \label{prop:para-calc-change-r}
  $\weyld{\gamma} \nRic{X}{} = -\D_X\gamma + \algbrac{X}{\gamma}_{\fp^{\perp}}$,
\end{enumerate}
where $\algbrac{X}{\gamma}_{\fp^0}$ and $\algbrac{X}{\gamma}_{\fp^{\perp}}$ are the projections of the algebraic bracket $\algbrac{X}{\gamma}$ to the appropriate summands using $\xi_M$. \noproof \end{prop}

For completely reducible $P$-representations, $\exp\fp^{\perp}$ acts trivially and we do not see the $\D_X\gamma$ terms.  By projecting to $F^0$ and writing $\gamma = (\gamma_{-k}, \ldots, \gamma_{-1})$ \etc, we can also recover the componentwise \formulae\ given by {\v C}ap and Slov{\'a}k \cite{cs2003-weylstr}; see also \cite[\S 5.1]{cs2009-parabolic1}.

\subsection{The equivalence of categories} 
\label{ss:para-calc-equiv}

A \emph{morphism} between two Cartan geometries $(F_1^P \surjto M, \, \omega_1)$ and $(F_2^P \surjto M, \, \omega_2)$ of type $G/P$ over $M$ is a principal $P$-bundle morphism $\Psi:F_1^P \to F_2^P$ for which $\Psi^* \omega_2 = \omega_1$.  This makes Cartan geometries of type $G/P$ over $M$ into a category, and the conformal equivalence problem may be restated as an equivalence of categories between conformal manifolds and ``normal'' conformal Cartan geometries.  Modulo a minor restriction, a similar statement holds for all parabolic geometries, where the conformal structure is replaced with the underlying infinitesimal flag structure from Subsection \ref{ss:para-para-para}.

\begin{defn} \thlabel{defn:para-calc-regcurv} A parabolic geometry $(F^P \surjto M, \, \omega)$ with curvature form $K_M \in \s{2}{\fg_M}$ is \emph{regular} if $K_M(TM_i, TM_j) \subseteq (\fg_M)_{i+j-1}$ for all $i,j > 0$. \end{defn}

In particular, abelian parabolic geometries are automatically regular.  More generally, a Weyl structure yields an isomorphism $TM \isom \gr{(TM)}$.  The algebraic bracket satisfies $\algbrac{TM_{(i)}}{TM_{(j)}} \subseteq TM_{(i+j-1)}$ for all $i>0$, making $(\gr{(TM)}, \algbrac{}{})$ into a bundle of nilpotent Lie algebras modelled on $\fg/\fp$.
Regularity is equivalent to local triviality of $(\gr{(TM)}, \, \algbrac{}{})$, together with a reduction of the frame bundle of $\gr{(TM)}$ \wrt\ $\Ad : P^0 \to \grp{Aut}{\fg/\fp}_{\mr{gr}}$; see \cite[p.\ 252]{cs2009-parabolic1}.  Thus regularity ensures a close relationship between the filtration of $TM$ and the reduction of $\gr{(TM)}$.

The next natural question to ask is which Cartan connections induce the same underlying infinitesimal flag structure.  The difference $\Phi \defeq \b{\omega} - \omega$ of Cartan connections is a horizontal $\fg$-valued $1$-form on $F^P$, so may be viewed as a section of $T^*M \tens \fg_M$.

\begin{lem} \rmcite[Prop.\ 3.1.10]{cs2009-parabolic1} Cartan connections $\omega$ and $\b{\omega}$ on $F^P\surjto M$ induce the same infinitesimal flag structure if and only if $\Phi(TM_i) \subseteq (\fg_M)_{i-\ell}$ for some $\ell\geq 1$. \noproof \end{lem}

In this case the difference $\b{K} - K$ of curvatures maps $TM_i \times TM_j$ to $(\fg_M)_{i+j- \ell}$, thus yielding a map $\gr{(\b{K} - K)} : \gr{(TM)} \times \gr{(TM)} \to \gr{(\fg_M)}$.  On the other hand, regularity implies that $\Phi \defeq \b{\omega} - \omega$ descends to a map $\gr{\Phi} : \gr{(TM)} \to \gr{(\fg_M)}$.  If $\liediff : \gr{(T^*M \tens \fg_M)} \to \gr(\Wedge{2}T^*M \tens \fg_M)$ is the graded bundle map induced by the Lie algebra differential, it is straightforward to check that $\gr[(\ell)](\b{K} - K) = \liediff(\gr[(\ell)] \Phi)$, where $\gr[(\ell)]$ denotes the $\ell$th graded component.  This suggests a normalisation condition where $\gr(K)$ takes values in a natural subbundle complementary to $\im\liediff$.  By Kostant's Hodge decomposition \eqref{eq:lie-hom-hodge}, such a subbundle is provided by $\ker \liebdy$.

\begin{defn} \thlabel{defn:para-calc-normal} A Cartan connection $\omega$ is \emph{normal} if its curvature $K$ satisfies $\liebdy K = 0$.  In this case, the \emph{harmonic curvature} $K_{\circ}$ is the image of $K$ in $\liehom{2}{\fg_M}$. \end{defn}

Normality is particularly simple for abelian parabolic geometries; we will give a characterisation in terms of torsion before \thref{thm:para-bgg-ablcurv}.

\smallskip
\begin{thm} If $K(TM_i, TM_j) \subseteq (\fg_M)_{i+j-\ell}$ for some $\ell\geq 1$ then $\liediff(\gr[(\ell)] K) = 0$.  Moreover if $\omega$ is normal then $\gr[(\ell)] K$ is a section of $\ker\quab$. \noproof \end{thm}

Choosing an algebraic Weyl structure, it follows that $\gr[(\ell)] K$ coincides with the graded component of harmonic curvature $K_{\circ}$ of degree $\ell$.

\begin{cor} For regular normal parabolic geometries, $K = 0$ if and only if $K_{\circ} = 0$; thus the harmonic curvature is a complete obstruction to local flatness. \noproof \end{cor}

A regular Cartan connection can always be modified to produce a normal Cartan connection inducing the same underlying infinitesimal flag structure \cite[Thm.\ 3.1.13]{cs2009-parabolic1}.  To obtain a uniqueness result, note that $\Wedge{k}\fp^{\perp} \tens \fg \isom \Hom{\Wedge{k}(\fg/\fp)}{\fg}$ inherits a natural filtration by ``homogeneous degree'' from the filtrations of $\Wedge{k}(\fg/\fp)$ and $\fg$.  Since $\liebdy$ is filtration preserving, this descends to a filtration of $\liehom{k}{\fg}$.

\begin{thm} \thlabel{thm:para-calc-equiv} \rmcite[Thm.\ 3.22]{cs2000-prolongation} Suppose that $\liehom{1}{\fg}_{1} = 0$.  Then associating the underlying infinitesimal flag structure to any parabolic geometry induces an equivalence of categories between regular normal parabolic geometries of type $G/P$ over $M$ and regular infinitesimal flag structures of type $G/P$ over $M$. \noproof \end{thm}

For abelian parabolic geometries with $\liehom{1}{\fg}_1 = 0$, \thref{thm:para-calc-equiv} yields an equivalence of categories between normal parabolic geometries and first-order $P^0$-structures.

\begin{rmk} \thlabel{rmk:para-calc-hom1} The condition $\liehom{1}{\fg}_1 = 0$ ensures that the reduction of the frame bundle of $\gr{(TM)}$ to structure group $P^0$ contains geometric information, other than an orientation.  It turns out that $\liehom{1}{\fg}_k = 0$ for all $k > 0$ unless $\fp$ or its complexification contain a simple ideal isomorphic to either
\vspace{-0.8em}
\begin{equation} \label{eq:para-calc-problem}
  \dynkin{ \DynkinLine{0}{0}{1}{0};
           \DynkinDots{1}{0}{3}{0};
           \DynkinLine{3}{0}{4}{0};
           \DynkinXDot{0}{0};
           \DynkinWDot{1}{0};
           \DynkinWDot{3}{0};
           \DynkinWDot{4}{0}; }{para-calc-equiv-problem1}
  \quad\text{or}\quad
  \dynkin{ \DynkinLine{0}{0}{1}{0};
           \DynkinDots{1}{0}{3}{0};
           \DynkinLine{3}{0}{4}{0};
           \DynkinDoubleLine{5}{0}{4}{0};
           \DynkinXDot{0}{0};
           \DynkinWDot{1}{0};
           \DynkinWDot{3}{0};
           \DynkinWDot{4}{0};
           \DynkinWDot{5}{0}; }{para-calc-equiv-problem2};
\vspace{-0.4em}
\end{equation}
see \cite[Prop.\ 5.1]{y1993-diffsimplegr}.  It is possible \cite[\S3.1.16]{cs2009-parabolic1} to obtain a similar equivalence of categories in the case that $\liehom{1}{\fg}_1 \neq 0$, by associating to each parabolic geometry a stronger underlying structure called a \emph{$P$-frame bundle of degree one}.  Only the first case of \eqref{eq:para-calc-problem} is problematic for us, so we postpone this discussion until Subsection \ref{ss:proj-para-cartan}. \end{rmk}

\section{BGG operators and curvature decomposition} 
\label{s:para-bgg}

Parabolic geometries also have a well-developed theory of invariant differential operators thanks to the \emph{curved Bernstein--Gelfand--Gelfand (BGG) sequence}, which we introduce in Subsection \ref{ss:para-bgg-bgg}.  This allows us to describe the decomposition of the Cartan curvature \wrt\ a Weyl structure in Subsection \ref{ss:para-bgg-curv}.  Finally, it turns out that the solution space of the first BGG operator can be prolonged to a closed system of PDEs on an auxiliary tractor bundle; we discuss this is Subsection \ref{ss:para-bgg-prol}.

\subsection{The curved BGG sequence} 
\label{ss:para-bgg-bgg}

The invariant differential operators between tractor bundles over the flat model $G/P$ have a complete description in terms of \emph{generalised Verma modules}, which are $\fg$-representations $M_{\fp}(\lambda) \defeq \cU(\fg) \tens[\cU(\fp)] \bV_{\lambda}$ classified by the $\fp$-highest weights modules $\bV_\lambda$.  As in Subsection \ref{ss:lie-hom-bbw}, the affine action of the Hasse diagram of $\fp$ determines a graph of weights of $\fp$, and it turns out there is a \non trivial\ homomorphism $M_{\fp}(\lambda) \to M_{\fp}(\mu)$ whenever $\lambda,\mu$ are joined by a single arrow in this diagram.  Thus the Hasse diagram determines a sequence of generalised Verma module homomorphisms for each $\bV$, which may be viewed as homomorphisms between irreducible components of $\liehom{}{\bV}$ via \thref{thm:lie-hom-bbw}.  Taken together, these yield a resolution
\begin{equation*}
  0 \From \liehom{\dim\fp^{\perp}}{\bV} \From \,\cdots\,
    \From \liehom{1}{\bV} \From \liehom{0}{\bV} \From \bV \From 0
\end{equation*}
of $\bV$ called the \emph{Bernstein--Gelfand--Gelfand complex}.  On the other hand, a differential operator between two tractor bundles over $G/P$ is invariant if and only if it is dual to a homomorphism of generalised Verma modules.  In this picture the Hasse diagram of $\fp$ determines a complex of invariant differential operators between irreducible components of homology modules, which can be computed algorithmically using Kostant's version of the Bott--Borel--Weil theorem; see \thref{thm:lie-hom-bbw}.

Now let $(F^P \surjto M, \, \omega)$ be a parabolic geometry.  Following work of Baston in conformal geometry \cite{b1990-verma,b1991-ahs1,b1991-ahs2}, a \emph{curved Bernstein--Gelfand--Gelfand sequence} for general parabolic geometries was constructed by {\v C}ap, Slov{\'a}k and Sou{\v c}ek \cite{css2001-bgg}.  The construction was later simplified by Calderbank and Diemer \cite{cd2001-curvedbgg}; we follow this account here.

Fix a $\fg$-representation $\bV$ with tractor bundle $\cV \defeq \assocbdl{F^P}{P}{\bV}$.  By \thref{lem:lie-hom-homeq}, the Lie algebra boundary map $\liebdy : \Wedge{k}\fp^{\perp} \tens \bV \to \Wedge{k-1}\fp^{\perp} \tens \bV$ induces a complex of bundle maps $\liebdy : \Wedge{k}T^*M \tens \cV \to \Wedge{k-1}T^*M \tens \cV$, with homology $\liehom{k}{\cV} \defeq \assocbdl{F^P}{P}{\liehom{k}{\bV}}$.

\begin{thm} \emph{\cite[Thm.\ 3.6]{cd2001-curvedbgg}} There is a natural sequence
\begin{equation} \label{eq:para-bgg-bgg}
  \s{0}{\liehom{0}{\cV}} \xrightarrow{ \bgg[\bV]{0} }
  \s{0}{\liehom{1}{\cV}} \xrightarrow{ \bgg[\bV]{1} }
  \s{0}{\liehom{2}{\cV}} \xrightarrow{ \bgg[\bV]{2} } \cdots
\end{equation}
of linear differential operators whose symbols depend on $(G/P,\bV)$ but not on $(M,\omega)$.  Moreover if $\omega$ is flat then \eqref{eq:para-bgg-bgg} is locally exact. \noproof \end{thm}

The BGG operators $\bgg[\bV]{k}$ are defined by constructing linear differential operators $\Wedge{k}T^*M \tens \cV \to \Wedge{k+1}T^*M \tens \cV$ which vanish on $\im \liebdy$ and take values in $\ker \liebdy$, thus yielding differential operators $\liehom{k}{\cV} \to \liehom{k+1}{\cV}$ on homology.  On the level of $\fp$-representations, the Lie algebra differential $\liediff : \Wedge{k}\fp^{\perp}\tens\bV \to \Wedge{k+1}\fp^{\perp}\tens\bV$ is a likely candidate, but unfortunately does not determine a map of associated bundles due to its lack of $P$-invariance.%
\footnote{At least, not without choosing a Weyl structure.}
Ignoring this for now, by choosing an algebraic Weyl structure for $\fp$ we can identify $\liehom{k}{\bV} \isom \ker\quab$ using Kostant's Hodge decomposition \eqref{eq:lie-hom-hodge}.  Since $\quab$ is invertible on its image and commutes with $\liediff$, the projection onto $\ker\quab$ is
\begin{equation} \label{eq:para-bgg-liediffPi}
  \id - \quab^{-1}\quab
    = \id - \quab^{-1}\liebdy \circ \liediff - \liediff \circ \quab^{-1}\liebdy.
\end{equation}
The lack of $P$-equivariance can be remedied by replacing $\liediff$ with the exterior covariant derivative $\d^{\bV} : \s{k}{\cV} \to \s{k+1}{\cV}$ induced by the tractor connection $\D^{\bV}$, where
\begin{align*}
  (\d^{\bV} s)(X_0,\ldots,X_k)
  &= \sum{i=0}{k} (-1)^i \, \D^{\bV}_{X_i} s(X_0, \ldots, \hat{X_i}, \ldots, X_k) \\
  &\qquad + \sum{i<j}{} (-1)^{i+j} s(\liebrac{X_i}{X_j}, X_0, \ldots, \hat{X_i}, \ldots, \hat{X_j}, \ldots, X_k)
\end{align*}
and hat denotes omission.  Extending $\d^{\bV}$ to $\s{0}{\Wedge{k}\fg_M^*\tens\cV}$ in the obvious way, one can show that $\d^{\bV}$ is a first-order $P$-invariant modification of the bundle map $\Wedge{k}\fg_M^*\tens\cV \to \Wedge{k+1}\fg_M^*\tens\cV$ induced by the Lie algebra differential $\liediff : \Wedge{k}\fg^* \tens \bV \to \Wedge{k+1}\fg^* \tens \bV$, prompting the definition of the following first-order operator \cite[Eqn.\ (4.1)]{cd2001-curvedbgg}.

\begin{defn} The operator $\quab_M \defeq \d^{\bV} \circ \liebdy + \liebdy \circ \d^{\bV} : \s{k}{\cV} \to \s{k}{\cV}$ is called the \emph{first-order laplacian}. \end{defn}

By choosing a Weyl structure, the algebraic laplacian on $\Wedge{k}\fp^{\perp} \tens \bV$ induces a bundle map $\quab : \Wedge{k}T^*M \tens \cV \to \Wedge{k}T^*M \tens \cV$.  One can show that
\begin{equation} \label{eq:para-bgg-quabladiff}
  (\quab_M - \quab)s = \sum{i}{} \, \ve^i \acts (\D_{e_i} s + \nRic{e_i}{} \acts s)
\end{equation}
for any local frame $\{e_i\}_i$ of $M$ with dual coframe $\{\ve^i\}_i$.  Thus $(\quab_M - \quab) s$ has strictly lower weight than $s$; in particular the restriction of $\quab_M$ to each graded component of $\cV \isom \gr{\cV}$ coincides with $\quab$, since there $\fp^{\perp}_M \isom T^*M$ acts trivially.  By writing
\begin{equation} \label{eq:para-bgg-neumann}
  \quab_M = \quab(\id - \cN)
  \quad\text{where}\quad
  \cN \defeq -\quab^{-1}(\quab_M - \quab),
\end{equation}
we arrive at the following.

\begin{prop} \thlabel{prop:para-bgg-quablainv} $\quab_M$ is invertible on $\im \liebdy$, with finite-order differential inverse given by the Neumann series $\quab_M^{-1} = \left( \sum{k\geq 0}{} \, \cN^k \right) \quab^{-1}$. \noproof \end{prop}

\thref{prop:para-bgg-quablainv} suggests that we consider the differential operator%
\footnote{We will often suppress mention of $\bV$ in this and later \formulae.}
\begin{equation*}
  \bggpi[\bV]{k} \defeq
    \id - \quab_M^{-1} \liebdy \circ \d^{\bV} - \d^{\bV} \circ \quab_M^{-1} \liebdy
\end{equation*}
as an analogue of \eqref{eq:para-bgg-liediffPi}.  From the definition and algebraic properties of $\liebdy$, it straightforward to see that $\bggpi[\bV]{k}$ maps $\Wedge{k}T^*M \tens \cV$ to itself; vanishes on $\im \liebdy$; takes values in $\ker \liebdy$; and induces the identity map on homology.  Therefore we have natural differential projections to and representations of homology classes given by
\vspace{-0.2em}
\begin{equation} \label{eq:para-bgg-projrepr} \begin{alignedat}{4}
  \bggproj[\bV]{k} &\defeq \mr{proj} \circ \bggpi[\bV]{k}
    &&: \s{k}{\cV} \to \s{0}{\liehom{k}{\cV}} \\
  \bggrepr[\bV]{k} &\defeq \bggpi[\bV]{k} \circ \mr{repr}
    &&: \s{0}{\liehom{k}{\cV}} \to \s{k}{\cV},
\end{alignedat}
\vspace{-0.2em}
\end{equation}
where $\mr{proj}$ is the projection to homology and $\mr{repr}$ is the choice of a representative.  Moreover $\bggproj{k} \circ \quab_M = 0$, while $\bggrepr{k}$ provides the unique representative in $\ker \quab_M$.

\begin{defn} The BGG operator is $\bgg[\bV]{k} \defeq \bggproj[\bV]{k+1} \circ \d^{\bV} \circ \bggrepr[\bV]{k}$. \end{defn}

Since $\bggproj{k+1}$ and $\bggrepr{k}$ may be differential operators themselves, the BGG operators are generally higher than first order.  The first BGG operator $\bgg{0}$ is always finite order, with the order of $\bgg{0}$ equal to the difference in (geometric) weights of the $\fp$-representations $\liehom{0}{\bV}$ and $\liehom{1}{\bV}$; see \cite[\S 3]{n2010-weightedjet}.

In general the composition $\bgg{k+1} \circ \bgg{k}$ does not vanish due to the curvature of $\omega$; it is straightforward to compute that $\bgg{k+1} \circ \bgg{k} = \bggproj{k+2} \circ \Curv[\bV]{}{} \circ \bggrepr{k}$, so that we recover the BGG complex on the flat model $G/P$.

There is also a bilinear pairing $\liehom{k}{\cV_1} \times \liehom{\ell}{\cV_2} \to \liehom{k+\ell}{\cV_3}$ associated to any triple of $\fg$-representations $\bV_1,\bV_2,\bV_3$ with a bilinear pairing $\bV_1 \times \bV_2 \to \bV_3$.  Indeed, the associated wedge product $\wedge : \bV_1 \times \bV_2 \to \bV_3$ induces
\vspace{-0.2em}
\begin{equation} \label{eq:para-calc-bggcup}
  \sqcup \defeq \bggproj[\bV_3]{k+\ell} \circ \wedge \circ (\bggrepr[\bV_1]{k}, \bggrepr[\bV_2]{\ell}),
\vspace{-0.2em}
\end{equation}
a bilinear differential pairing on homology.  Up to curvature corrections, the BGG operators $\bgg[\bV]{k}$ satisfy a Leibniz rule over $\sqcup$, with $\sqcup$ a cup product on homology \cite[Prop.\ 5.7]{cd2001-curvedbgg}.  One can also show that the symbol of $\sqcup$ depends only on $(G/P,\bV_1,\bV_2,\bV_3)$ but not on $(M,\omega)$.  We will see an explicit example (for $k=\ell=0$) in Chapter \ref{c:mob2}.

\begin{rmk} One can define other multilinear differential pairings associated to the BGG complex, and these pairings have a rich algebraic structure which is encapsulated in a \emph{curved $A_{\infty}$-algebra}; see \cite[\S6]{cd2001-curvedbgg} and \cite{k1999-Ainfinity}. \end{rmk}

\subsection{Curvature decomposition} 
\label{ss:para-bgg-curv}

The $P$-invariant decomposition of the Cartan connection $\omega$ \wrt\ a Weyl structure $\xi_M$ provided by \eqref{eq:para-calc-cconn} and \eqref{eq:para-calc-cconnp} also induces a $P$-invariant decomposition of the Cartan curvature $K \defeq \d\omega + \tfrac{1}{2}\algbracw{\omega}{\omega}{} \in \s[F^P]{2}{\fg}$.  Indeed, writing
\begin{equation*}
  \omega = \xi_M\omega_{\fg/\fp} + \underbrace{(\xi_M\omega_{\fp^0}
    - \xi_{\fw}^*\eta)}_{\eqdef\, \omega_{\fw}} + \underbrace{(\omega_{\fp^{\perp}}
    + \xi_{\fw}^*\eta)}_{\eqdef\, \rho}
\end{equation*}
yields a $P$-invariant decomposition $K = \xi_M K_{\fg/\fp} + \xi_M K_{\fp^0} + K_{\fp^{\perp}}$, where
\vspace{-0.3em}
\begin{equation} \label{eq:para-bgg-curvdecomp} \begin{aligned}
  K_{\fg/\fp}
    &\defeq \xi_M( \d\omega_{\fg/\fp}
      + \tfrac{1}{2} \algbracw{ \omega_{\fg/\fp} }{ \omega_{\fg/\fp} }{ } ) \\
    & \qquad + \algbracw{ \xi_M\omega_{\fg/\fp} }{ \omega_{\fw} }{}
             + \algbracw{ \xi_M\omega_{\fg/\fp} }{ \rho }{}_{\fg/\fp} \\
  K_{\fp^0}
    &\defeq ( \d\omega_{\fw} + \tfrac{1}{2} \algbracw{ \omega_{\fw} }{ \omega_{\fw} }{ } ) \\
    & \qquad + \algbracw{ \xi_M\omega_{\fg/\fp} }{ \rho }{}_{\fp^0} \\
  \text{and}\quad 
  K_{\fp^{\perp}}
    &\defeq ( \d\rho + \tfrac{1}{2} \algbracw{ \rho }{ \rho }{ } )
      + \algbracw{ \omega_{\fw} }{ \rho }{} \\
    & \qquad + \algbracw{ \xi_M\omega_{\fg/\fp} }{ \rho }{}_{\fp^{\perp}}.
\end{aligned}
\vspace{-0.3em}
\end{equation}
Moreover since $K$ is horizontal, \eqref{eq:para-bgg-curvdecomp} descends to a decomposition $K_M = \Tor[\D]{} + \Weyl[\D]{}{} + \CY{}{}$ of the curvature form $K_M \in \s{2}{\fg_M}$ of $\omega$.

\begin{defn} \thlabel{defn:para-bgg-curvcmpts} The components $\Tor[\D]{} \in \s{2}{TM}$, $\Weyl[\D]{}{} \in \s{2}{\fp^0_M}$ and $\CY[\D]{}{} \in \s{2}{T^*M}$ are known respectively as the \emph{torsion}, the \emph{Weyl curvature} and the \emph{Cotton--York tensor} of the Weyl structure $\xi_M$. \end{defn}

It follows immediately from \thref{cor:para-calc-tractorcurv} that the curvature $\Curv[\bV]{}{}$ of the tractor connection $\D^{\bV}$ on $\cV \defeq \assocbdl{F^P}{P}{\bV}$ acts via
\begin{equation*}
  \Curv[\bV]{X,Y}{s} = (\Tor[\D]{X,Y} + \Weyl[\D]{X,Y}{} + \CY{X,Y}{}) \acts s
\end{equation*}
for all $s \in \s{0}{\cV}$.  If the Cartan connection $\omega$ is normal then $\p K = 0$ by definition, and since $\p$ preserves the $\fp^{\perp}$-filtration of $\fg$ we must also have $\p\Tor[\D]{} = 0$, $\p\Weyl[\D]{}{} = 0$ and $\p\CY{}{} = 0$ as elements of $\s{1}{TM}$, $\s{1}{\fp^0_M}$ and $\s{1}{T^*M}$ respectively.

In general it is quite difficult to describe the components $K_{\fg/\fp}$, $K_{\fp^0}$ and $K_{\fp^{\perp}}$ explicitly.  As a first step, \itemref{prop:para-calc-weylconn}{D} implies the expressions on the first lines of these components in \eqref{eq:para-bgg-curvdecomp} descend to the Cartan torsion $T_M$, the curvature $\Curv{}{}$ of the Weyl connection $\D$, and the exterior covariant derivative $\d^{\D} \nRic{}{}$ of the normalised Ricci tensor.  The remaining terms are described as follows.  \cite[Thm.\ 5.2.9]{cs2009-parabolic1}.

\begin{thm} \thlabel{thm:para-bgg-curv} Choose a Weyl structure and let $\liediff : T^*M\tens\fg_M \to \Wedge{2}T^*M\tens\fg_M$ be the resulting bundle map induced by the Lie algebra differential $\liediff$.  Then
\begin{equation*}
  K_M = (\Tor[\D]{}, \Curv{}{}, \CY{}{}) + \liediff \nRic{}{}
\end{equation*}
as an element of $\s{2}{TM\dsum\fp^0_M\dsum T^*M}$.  In particular, the harmonic curvature $K_{\harm}$ of $\omega$ coincides with the components of $K_M$ lying in $\ker\quab_M$. \noproof \end{thm}

In the sequel we shall be exclusively interested in abelian parabolic geometries, for which the curvature decomposition is straightforward to describe.  Then since $\fg/\fp$ and $\fp^{\perp}$ are irreducible $\fp$-representations, \eqref{eq:para-bgg-curvdecomp} becomes
\begin{equation*}
\vspace{-0.075em}
\begin{aligned}
  K_{\fg/\fp}
    &= \xi_M \d\omega_{\fg/\fp}
      + \algbracw{ \xi_M\omega_{\fg/\fp} }{ \omega_{\fw} }{} \\
  K_{\fp^0}
    &= ( \d\omega_{\fw}
      + \tfrac{1}{2} \algbracw{ \omega_{\fw} }{ \omega_{\fw} }{ } )
      + \algbracw{ \xi_M\omega_{\fg/\fp} }{ \rho }{} \\
  \text{and}\quad 
  K_{\fp^{\perp}}
    &= \d\rho
      + \algbracw{ \omega_{\fw} }{ \rho }{}.
\end{aligned}
\vspace{-0.075em}
\end{equation*}
The solder form $\omega_{\fg/\fp}$ defines an isomorphism $TM \isom \gr{(TM)}$, so that the section of $T^*M\tens TM$ induced by $\xi_M \omega_{\fg/\fp}$ is just the identity map.  We conclude that
\vspace{-0.075em}
\begin{equation} \label{eq:para-bgg-ablcurv}
  \Tor[\D]{} = T_M, \quad
  \Weyl[\D]{}{} = \Curv[\D]{}{} + \algbracw{\id}{\nRic{}{}}{}
  \quad\text{and}\quad
  \CY{}{} = \d^{\D} \nRic{}{},
\vspace{-0.075em}
\end{equation}
where $\algbracw{\id}{\nRic{}{}}{X,Y} \defeq \algbrac{X}{\nRic{Y}{}} - \algbrac{Y}{\nRic{X}{}}$.  In particular $\Tor{} \defeq \Tor[\D]{}$ is just the torsion of $\D$, which is independent of $\D$ by \thref{prop:para-calc-change} and the fact that $\algbrac{ \algbrac{X}{\gamma} }{ Y } = \algbrac{ \algbrac{Y}{\gamma} }{ X }$ for all $X,Y \in \s{0}{TM}$ and $\gamma \in \s{1}{}$.

Normality of the Cartan connection is equivalent to having $\liebdy\Tor{} = 0$; in this case, since $\liebdy \Weyl[\D]{}{} = 0$ also, $\nRic{}{}$ is the unique solution of $\liebdy\Curv{}{} + \liebdy \algbracw{\id}{\nRic{}{}}{} = 0$.  Using that $\liebdy \nRic{}{} = 0$ automatically, $\algbracw{\id}{\nRic{}{}}{} = \liediff \nRic{}{}$ and that $\quab_M$ agrees with the algebraic laplacian $\quab = \liebdy \liediff + \liediff \liebdy$ on $T^*M$, it follows that $\nRic{}{} = -\quab_M^{-1} \liebdy \Curv{}{}$.

\begin{thm} \thlabel{thm:para-bgg-ablcurv} Let $(F^P\surjto M, \, \omega)$ be a normal abelian parabolic geometry of type $G/P$, and let $\xi_M$ be a Weyl structure.  Then:
\begin{enumerate}
  \item \label{thm:para-bgg-ablcurv-weyl}
  The Weyl connections are precisely those with $\liebdy$-closed torsion, which is invariant and coincides with the degree one component of the harmonic curvature.
  
  \item \label{thm:para-bgg-ablcurv-cmpts}
  $\nRic{}{} = -\quab_M^{-1} \liebdy \Curv{}{}$, and the remaining components of the harmonic curvature are given by the components of $\Weyl[\D]{}{} \defeq \Curv{}{} + \algbracw{\id}{\nRic{}{}}{}$ and $\CY{}{} \defeq \d^{\D} \nRic{}{}$ in $\ker \quab_M$.
  
  \item \label{thm:para-bgg-ablcurv-weyld}
  Under infinitesimal change of Weyl structure, we have $\weyld{\gamma}\Tor{} = 0$, $\weyld{\gamma} \Weyl[\D]{}{} = \algbrac{\Tor{}}{\gamma}$, $\weyld{\gamma}\nRic{}{} = \D\gamma$, $\weyld{\gamma} \Curv{}{} = -\algbracw{\id}{\D\gamma}{} + \algbrac{\Tor{}}{\gamma}$ and $\weyld{\gamma} \CY{}{} = \Weyl[\D]{}{\gamma}$. \noproof
\end{enumerate}
\end{thm}

For proof of parts \ref{thm:para-bgg-ablcurv-weyl} and \ref{thm:para-bgg-ablcurv-cmpts} of \thref{thm:para-bgg-ablcurv}, see \cite[\S 4.7]{cs2003-weylstr} and \cite[\S 5.2.3]{cs2009-parabolic1}.  Part \ref{thm:para-bgg-ablcurv-weyld} follows easily from \thref{prop:para-calc-change} and the relations \eqref{eq:para-bgg-ablcurv}.

\subsection{Prolongation of BGG operators} 
\label{ss:para-bgg-prol}

The first BGG operator $\bgg{\bV}$ on a tractor bundle $\cV \defeq \assocbdl{F^P}{P}{\bV}$ is of finite type \cite{cd2001-curvedbgg,css2001-bgg}, so that its kernel is finite dimensional.  Such differential operators are typically studied by \emph{prolongation} (\ie, by further differentiation) to obtain a closed system of PDEs describing the solutions.  Branson \etal\ described a (\non invariant) method of prolongation for abelian parabolic geometries \cite{bceg2006-geometricprol} which generalised examples from conformal geometry \cite{h2008-confbgg} and \proj\ differential geometry \cite{em2008-projmetrics,eg2011-projectivebgg}.  An invariant prolongation of the whole BGG sequence was later obtained by Hammerl \etal\ \cite{hsss2010-prolconns,hsss2012-prolconns2}.

\begin{thm} \thlabel{thm:para-bgg-prol} \rmcite[Thm.\ 4.2]{hsss2010-prolconns}
There is a unique $(\im\liebdy)$-valued differential operator \mbox{$\Phi : \s{k}{\cV} \to \s{k+1}{\cV}$} such that the splitting operator $\bggrepr{\bV} : \liehom{k}{\cV} \to \Wedge{k} T^*M \tens \cV$ induces an isomorphism between $\ker \bgg{k}$ and $\ker(\d^{\bV} + \Phi) \intsct (\ker\liebdy)$.
\end{thm}

To describe this prolongation, note that the operators $\d^{\bV}$, $\bggrepr{k}$ and $\bgg{k}$ from the BGG sequence may be arranged as in Figure \ref{fig:para-bgg-prolsq}.  The key result is as follows \cite[Thm.\ 3.2]{hsss2010-prolconns}.

\begin{prop} \thlabel{prop:para-bgg-prolcomm} Figure \ref{fig:para-bgg-prolsq} commutes if and only if $\liebdy \circ \d^{\bV} \circ \d^{\bV}$ vanishes on $\im\bggrepr{k}$.  In this case, $\bggrepr{k}$ restricts to an isomorphism between $\ker\bgg{k}$ and $(\ker\d^{\bV}) \intsct (\ker\liebdy)$. \noproof \end{prop}

It will not generally be the case that $\liebdy \circ \d^{\bV} \circ \d^{\bV} = 0$.  To arrange this we will modify each $\d^{\bV}$ by a finite-order differential operator $\Phi$ so as not to change the resulting BGG operator $\bgg{k}$, and such that $\d^{\cV} \defeq \d^{\bV} + \Phi$ satisfies $\liebdy \circ \d^{\cV} \circ \d^{\cV} = 0$.  We thus obtain a bijection between $\ker\bgg{k}$ and $(\ker \d^{\cV}) \intsct (\ker\liebdy)$ as above.

\begin{figure}[h]
  \begin{equation*} \begin{CD}
    \Wedge{k}T^*M \tens \cV  @>\d^{\bV}>>          \Wedge{k+1}T^*M \tens \cV \\
    @A\bggrepr{k}AA     @AA\bggrepr{k+1}A  \\
    \liehom{k}{\cV}          @>>\bgg{k}>        \liehom{k+1}{\cV}.
  \end{CD} \end{equation*}
  \vspace{-0.3em}
  \caption[A commuting square in the curved BGG sequence]
          {The $k$th square of the curved BGG sequence on $\cV$.}
  \label{fig:para-bgg-prolsq}
\end{figure}

To construct $\Phi$, note that the $\fp^{\perp}$-filtration $\bV = \bV_N \supset \cdots \supset \bV_0 \supset 0$ of $\bV$ induces a $P$-invariant filtration of the Lie algebra chain space $\liechain{k}{\bV} \defeq \Wedge{k}\fp^{\perp} \tens \bV$ in the obvious way.  This in turn induces a filtration%
\footnote{This is the ``diagonal'' filtration from \cite[p.\ 12]{hsss2010-prolconns}.}
$A = A_k \supset \cdots \supset A_0 \supset 0$ of
\begin{equation*}
  A \defeq \Hom{\Wedge{k}\fp^{\perp} \tens \bV}{\Wedge{k+1}\fp^{\perp} \tens \bV},
\end{equation*}
where $\Phi \in A_i$ if and only if $\Phi(\Wedge{k}\fp^{\perp}) \tens \bV_j \subseteq \Wedge{k+1}\fp^{\perp} \tens \bV_{i+j}$ for all $j$.  This is clearly $P$-invariant, so induces a filtration of $\cA \defeq \assocbdl{F^P}{P}{A}$ by smooth subbundles $\cA_i$.

\begin{lem} \thlabel{lem:para-bgg-bggchange} Let $\Phi \in \s{0}{\cA_1}$ with $\im\Phi \subseteq \im\liebdy$.  Then the BGG operators $\bgg{k}$ and $\bggrepr{k}$ are unchanged by the replacement of $\d^{\bV}$ with $\d^{\bV} + \Phi$. \noproof \end{lem}

Such a replacement clearly does not affect the conclusion of \thref{prop:para-bgg-prolcomm}, so it remains to construct a suitable $\Phi$.  Recalling that $\d^{\bV} \circ \d^{\bV} = \Curv[\d^{\bV}]{}{}$ is just the curvature of $\d^{\bV}$, we define
\begin{equation*} 
  \Phi(s) \defeq - \quab_M^{-1} \liebdy ( \Curv[\hspace{0.2ex} \d^{\bV}]{}{s} )
\end{equation*}
for each $s \in \s{k}{\cV}$, which is well-defined by \thref{prop:para-bgg-quablainv}.  Then $\Phi$ takes values in $\im\liebdy$, since $\quab_M$ (and hence $\quab_M^{-1}$) commutes with $\liebdy$, and moreover $\Phi \in \s{0}{\cA_1}$ by properties of $\liebdy$ and $\d^{\bV}$.  Therefore by \thref{lem:para-bgg-bggchange} we may replace $\d^{\bV}$ with $\d^{\cV} \defeq \d^{\bV} + \Phi$ at each stage without affecting the BGG operators, for which
\begin{align*}
  \liebdy( \Curv[\hspace{0.2ex} \d^{\cV}]{}{s} )
  &=  \liebdy ( \Curv[\hspace{0.2ex} \d^{\bV}]{}{s} )
    - \liebdy \d^{\bV} \quab_M^{-1} \liebdy ( \Curv[\d^{\bV}]{}{s} ) \\
  &=  \liebdy ( \Curv[\hspace{0.2ex} \d^{\bV}]{}{s} )
    - \quab_M \quab_M^{-1} \liebdy ( \Curv[\d^{\bV}]{}{s} )
   = 0 \\[-2.1em]
\end{align*}
since $\quab_M = \liebdy \circ \d^{\bV}$ on $\im\liebdy \subseteq \ker\liebdy$.  By \thref{prop:para-bgg-prolcomm}, the splitting operator constructed from $\d^{\cV}$ then provides the prolongation of the BGG operator $\bgg{k}$.  It is easy to see that $\Phi$ is the unique differential correction with the desired properties.

\begin{defn} \thlabel{defn:para-bgg-prolconn} $\d^{\cV} \defeq \d^{\bV} - \quab_M^{-1} \liebdy (\Curv[\d^{\bV}]{}{})$ is the \emph{prolongation operator} of $\bgg{k}$. \end{defn}

Note that $\Phi \defeq -\quab_M^{-1} \liebdy (\Curv[\d^{\bV}]{}{})$ is not algebraic in general.  Indeed, each term $\quab_M - \quab$ appearing in the Neumann series \eqref{eq:para-bgg-neumann} for $\quab_M^{-1}$ is typically first-order, so that $\Phi$ has order bounded above by the height $N$ of the $\fp^{\perp}$-filtration of $\bV$.

When the representation $\bV$ is given, the correction $\Phi$ can be computed using a Weyl structure; independence from the choice of Weyl structure is ensured by the uniqueness of $\Phi$.  Then the Lie algebra differential and the algebraic laplacian induce commuting bundle maps $\liediff$ and $\quab = \liebdy \liediff + \liediff \liebdy$.  By \thref{thm:lie-hom-spectral}, $\quab$ acts by a scalar on each irreducible graded component of $\cV$, while \eqref{eq:para-bgg-quabladiff} and the Neumann series of \thref{prop:para-bgg-quablainv} allow us to compute the action of $\quab_M^{-1}$ on graded components.

For first BGG operators, where $\Curv[\d^{\bV}]{}{} = \d^{\bV} \circ \d^{\bV}$ is just the tractor curvature $\Curv[\bV]{}{}$, we obtain a modification $\Phi \defeq -\quab_M^{-1}\liebdy(\Curv[\bV]{}{})$.  In general $\Phi$ is not algebraic, so that $\d^{\cV}$ is not necessarily a connection on $\cV$.  In cases where $\Phi$ is algebraic, we refer to $\D^{\cV} \defeq \d^{\cV}$ as the \emph{prolongation connection} on $\cV$.  Due to a differential Bianchi identity, it will turn out that $\Phi$ is indeed algebraic in all cases that we are interested in; see Section \ref{s:ppg-bgg}.

\begin{rmk} Hammerl \etal\ \cite[\S 1.4]{hsss2012-prolconns2} also provide an iterative method for constructing a prolongation \emph{connection} $\D^{\cV}$ on $\cV$.  The process starts with $\D^0 \defeq \D^{\bV}$.  Then if $\phi_i$ is the \non zero\ graded component of $\liebdy (\Curv[\D^i]{}{})$ of highest weight, one defines $\D^{i+1} \defeq \D^i - \quab^{-1}\phi_i$, eventually reaching an invariant prolongation connection $\D^{\cV}$ whose curvature satisfies $\liebdy (\Curv[\D^{\cV}]{}{}) = 0$.  However, in the sequel we shall calculate using $\quab_M^{-1}$ due to the aforementioned differential Bianchi identity. \end{rmk}
\chapter{\Proj\ differential geometry} 
\label{c:proj}

\renewcommand{\algbracadornment}{r}
\BufferDynkinLocaltrue
\renewcommand{\dynkinnameoffset}{-0.75}

\Proj\ differential geometry is a classical subject which studies the behaviour of unparametrised geodesics.  This leads to an equivalence relation among \riem\ metrics, where two metrics are \emph{\proj ly\ equivalent} if they have the same geodesics (as unparametrised curves).  We begin in Section \ref{s:proj-class} with a review of the classical approach to the theory, including the so-called \emph{main equation} describing the set of \proj ly\ equivalent metrics.  We describe \proj\ differential geometry as an abelian parabolic geometry in Section \ref{s:proj-para}, where the pertinent Lie algebras are $\fg = \alg{sl}{n+1,\bR}$ with parabolic $\fp$ given by crossing the last node.  The flat model is then $G\acts\fp \isom \RP[n]$ which, as we saw in \thref{rmk:para-calc-hom1}, is problematic for the general equivalence of categories.  We indicate how to obtain an equivalence using $P$-frame bundles in Subsection \ref{ss:proj-para-cartan}.

In the parabolic picture, the main equation may be interpreted as the first BGG operator associated to the representation $\bW = \Symm{2}\bR^{n+1}$ of $\fg$.  We describe the resulting differential equation in Section \ref{s:proj-bgg}, as well as obtaining the prolongation given by Eastwood and Matveev \cite{em2008-projmetrics} in index-free notation.  Finally, we make an observation regarding the representation $\bW$ which is crucial for later generalisations: the direct sum $\fh \defeq \bW \dsum (\fg \dsum \bR) \dsum \bW^*$ is a graded Lie algebra isomorphic to $\alg{sp}{2n+2,\bR}$.

\section{Classical definition and results} 
\label{s:proj-class}

We begin by reviewing the classical formulation of \proj\ differential geometry.  For greater generality we allow \non degenerate\ metrics of arbitrary signature, reserving the adjective \emph{\riem} for positive definite metrics.  In what follows, suppose that $(M,g)$ is an $n$-dimensional (\pseudo)\riem\ manifold with \LC\ connection $\D$.

\begin{defn} \thlabel{defn:proj-class-equiv} A smooth curve $\gamma\subset M$ is a \emph{geodesic} of $g$ if $\D_X X\in \linspan{X}{}$ for every vector field $X$ tangent to $\gamma$.  Two metrics $g,\b{g}$ are called \emph{\proj ly\ equivalent} if they have the same geodesics (as unparameterised curves). \end{defn}

Note that we do not require geodesics to be affinely parameterised.  The term \emph{geodesically equivalent} is also used in the literature; however, we will use the former to emphasise the underlying \proj\ geometry.  The following characterisation of \proj ly\ equivalent metrics is well-known \cite{bm2003-geombenenti, em2008-projmetrics, tm2003-geodintegrability}.

\begin{lem} \thlabel{lem:proj-class-algbrac} Metrics $g,\b{g}$ are \proj ly\ equivalent if and only if their \LC\ connections $\D,\b{\D}$ are related by
\vspace{-0.2em}
\begin{equation} \label{eq:proj-class-algbrac} \begin{aligned}
  \b{\D}_X Y &\phantom{:}= \D_X Y + \algbrac{X}{\alpha} \acts Y \\
  \text{where}\quad
  \algbrac{X}{\alpha} \acts Y &\defeq \tfrac{1}{2}\left( \alpha(X)Y + \alpha(Y)X \right)
\end{aligned}
\vspace{-0.2em}
\end{equation}
for some $\alpha\in\s{1}{}$ and all $X,Y\in\s{0}{TM}$. \end{lem}

\begin{proof} Since $\D, \b{\D}$ have the same unparametrised geodesics, it follows from \thref{defn:proj-class-equiv} that $\b{\D}_X X - \D_X X = \alpha(X)X$ for some smooth function $\alpha:TM\to\bR$.  By basic properties of linear connections we see that in fact $\alpha\in\s{1}{}$, so that $X\mapsto \alpha(X)X$ defines a quadratic form $\s{0}{TM}\times\s{0}{TM} \to \s{0}{TM}$.  Polarisation yields
\vspace{-0.2em}
\begin{equation} \label{eq:proj-class-algbrac-1}
  (\b{\D}_X Y - \D_X Y) + (\b{\D}_Y X - \D_Y X)
    = \alpha(X)Y + \alpha(Y)X
\vspace{-0.2em}
\end{equation}
which, upon noting that the \lhs\ in \eqref{eq:proj-class-algbrac-1} equals $2(\b{\D}_X Y - \D_X Y)$ since both $\D,\b{\D}$ are \torsionfree, gives \eqref{eq:proj-class-algbrac}. \end{proof}

\begin{rmk} \thlabel{rmk:proj-class-half} We call the endomorphism $\algbrac{X}{\alpha} \in \alg{gl}{TM}$ the \emph{algebraic bracket} of $X$ and $\alpha$.  It is clear that $\algbrac{X}{\alpha}\acts Y$ is symmetric in $X,Y$, and we may write $\algbrac{}{\alpha} = \id\symm\alpha \in \s{1}{\alg{gl}{TM}}$.  Note that the factor $\tfrac{1}{2}$ in the definition \eqref{eq:proj-class-algbrac} is a \non standard\ normalisation convention; the reason for this choice will become apparent later. \end{rmk}

Fix a local frame $\{e_i\}_i$ of $TM$ with dual coframe $\{\ve^i\}_i$.  Using the standard formula $\div[g]{X} = \tfrac{1}{2} \p_{e_i}( \log(\det g) \ve^i(X) )$ for the divergence of $X$ \wrt\ $g$ (see for example \cite[p.\ 436]{l2003-smoothmfds}), taking a trace in \eqref{eq:proj-class-algbrac} yields
\vspace{-0.2em}
\begin{equation} \label{eq:proj-class-exact1form}
  \alpha = \tfrac{1}{n+1} \,\d\! \left( \log \frac{\det \b{g}}{\det g} \right)
  \colvectpunct[-0.6em]{.}
\vspace{-0.2em}
\end{equation}
In particular, $\alpha$ is an exact $1$-form which depends only on $g$ and $\b{g}$.

A na{\"i}ve first approach to studying \proj ly\ equivalent pairs $g,\b{g}$ might be via the endomorphism $G$ satisfying $\b{g} = g(G\bdot,\bdot)$; equivalently $G = \sharp \circ \b{\flat}$, where $\flat = \sharp^{-1}$ and $\b{\flat} = \b{\sharp}^{-1}$ are the usual musical isomorphisms of $g$ and $\b{g}$ respectively.  However, it turns out to be more fruitful to instead study the endomorphism
\vspace{-0.2em}
\begin{equation} \label{eq:proj-class-endoA}
  A(g,\b{g}) \defeq \left(\frac{\det\b{g}}{\det g}\right)^{\!\!1/(n+1)} \b{\sharp}\circ\flat.
\vspace{-0.2em}
\end{equation}
Clearly $A(g,\b{g})$ is invertible with inverse $A(\b{g},g)$, and is \self adjoint\ \wrt\ both $g, \b{g}$.  The endomorphisms $G$ and $A \defeq A(g,\b{g})$ are related by $A = (\det G)^{1/(n+1)} G^{-1}$, or equivalently $G = (\det A)^{-1} A^{-1}$.  Moreover $\b{g}$ can be recovered from the pair $(g,A)$ as
\vspace{-0.2em}
\begin{equation} \label{eq:proj-class-gbar}
  \b{g} = (\det A)^{-1}g(A^{-1}\bdot,\bdot).
\vspace{-0.2em}
\end{equation}
The following result explains the key benefit of using $A$ over $G$: it satisfies a first-order linear differential equation.  We reproduce the proof from \cite[Thm.\ 2]{bm2003-geombenenti} here for completeness, which the authors attribute to \cite{e1997-riemannian} (although see \thref{rmk:proj-class-sinjukov} below).

\begin{prop} \thlabel{prop:proj-class-maineqn} Let $g,\b{g}$ be metrics with \LC\ connections $\D,\b{\D}$ respectively.  Then $g,\b{g}$ are \proj ly\ equivalent if and only if $A=A(g,\b{g})$ defined by \eqref{eq:proj-class-endoA} satisfies the first-order linear differential equation
\vspace{-0.2em}
\begin{equation} \label{eq:proj-class-maineqn}
  g( (\D_X A)\bdot, \bdot) = X^{\flat}\symm\mu
\vspace{-0.2em}
\end{equation}
for some $\mu\in\s{1}{}$ and all $X\in\s{0}{TM}$.  In this case $\b{\D} = \D+\algbrac{}{\alpha}$, where $\alpha\in\s{1}{}$ satisfies $\mu = -\alpha(A\bdot) = \d(\tr A)$. \end{prop}

\begin{proof} By \thref{lem:proj-class-algbrac}, $g,\b{g}$ are \proj ly\ equivalent if and only if $\b{\D}=\D+\algbrac{}{\alpha}$ for some $\alpha\in\s{1}{}$.  In this case, differentiation gives
\vspace{-0.2em}
\begin{equation*}
\begin{aligned}
  (\D_X\b{g})(Y,Z)
  &= (\b{\D}_X \b{g})(Y,Z) - (\algbrac{X}{\alpha}\acts \b{g})(Y,Z) \\
  &= \b{g}(\algbrac{X}{\alpha} \acts Y,Z) + \b{g}(Y,\algbrac{X}{\alpha} \acts Z) \\
  &= \alpha(X)\b{g}(Y,Z) + \tfrac{1}{2}\alpha(Y)\b{g}(X,Z)
     + \tfrac{1}{2}\alpha(Z)\b{g}(X,Y) \\
  &= (\det A)^{-1}\left[ \alpha(X)g(A^{-1}Y,Z) + \tfrac{1}{2}\alpha(Y)g(A^{-1}X,Z)
     + \tfrac{1}{2}\alpha(Z)g(A^{-1}X,Y) \right].
\end{aligned}
\vspace{-0.2em}
\end{equation*}
On the other hand, \eqref{eq:proj-class-gbar} and the identity $\D_X A^{-1} = -A^{-1}\circ \D_X A \circ A^{-1}$ give
\vspace{-0.2em}
\begin{equation*}
\begin{aligned}
  \D_X\b{g}
  &= \D_X\left((\det A)^{-1}g(A^{-1}\bdot,\bdot)\right) \\
  &= -(\det A)^{-2} \d(\det A)(X) g(A^{-1}\bdot,\bdot)
     -(\det A)^{-1}g(\D_X A \circ A^{-1} \bdot, A^{-1}\bdot) \\
  &= (\det A)^{-1}\left[ -\d(\log\det A)(X)g(A^{-1}\bdot,\bdot)
     - g(\D_XA\circ A^{-1}\bdot, A^{-1}\bdot) \right] \\
  &= (\det A)^{-1}\left[ \alpha(X) g(A^{-1}\bdot,\bdot)
     - g(\D_XA\circ A^{-1}\bdot, A^{-1}\bdot) \right],
\end{aligned}
\vspace{-0.2em}
\end{equation*}
where the last line follows by \eqref{eq:proj-class-exact1form}.  Comparing these two expressions and precomposing with $A^{-1}$ in both slots yields \eqref{eq:proj-class-maineqn}, with $\mu\defeq -\alpha(A\bdot)$.  Raising an index using $g$ yields $(\D_X A)Y = \tfrac{1}{2} \left( g(X,Y) \mu^{\sharp} + \mu(Y)X \right)$, so that taking a trace over $Y$ in the last display equation gives $\mu = \d(\tr A)$ as required. \end{proof}

\begin{rmk} \thlabel{rmk:proj-class-sinjukov} Equation \eqref{eq:proj-class-maineqn} was known to Sinjukov \cite{s1979-geodmappings}, and is referred to as the \emph{Sinjukov equation} in a sizeable amount of the literature (see for example \cite{bm2011-splittingglueing, bm2003-geombenenti, km2010-projlichnerowicz}).  Equivalent equations may also be found in \cite{hmv2009-geodmappings}.  To avoid any historical misattribution, we shall call it the \emph{main equation} of \proj\ differential geometry.

Mike{\v s} later prolonged the main equation to obtain a closed differential system controlling the metrics \proj ly\ equivalent to $g$ \cite{m1996-geodprol}.  We shall study the invariant version of this prolongation, described by Eastwood and Matveev \cite{em2008-projmetrics}, in Section \ref{s:proj-bgg}. \end{rmk}

It is clear from \eqref{eq:proj-class-maineqn} that \proj ly equivalent metrics $g$, $\b{g}$ are affinely equivalent if and only if $A(g,\b{g})$ is a multiple of the identity, if and only if $\mu = \d(\tr A)$ vanishes.

\thref{prop:proj-class-maineqn} states that a solution $(g,A)$ of the main equation is equivalent to a \proj ly\ equivalent pair of metrics $g, \b{g}\defeq(\det A)^{-1} g(A^{-1}\bdot,\bdot)$.  Linearity of \eqref{eq:proj-class-maineqn} also implies that $(g, A_t \defeq A - t\id)$ is a solution for all $t\in\bR$ with $\det A_t \neq 0$.  We then have a $1$-parameter family of metrics $g_t \defeq (\det A_t)^{-1} g(A_t^{-1}\bdot, \bdot)$ \proj ly\ equivalent to $g$; we will study this \emph{metrisability pencil} in more detail in Chapter \ref{c:mob2}.

Finally, let us indicate the link between the classical theory above and the description as a parabolic geometry outlined in the following section.  Since geodesics are really a feature of connections rather than metrics, call two linear connections $\D,\b{\D}$ \emph{\proj ly equivalent} if $\b{\D}_X = \D_X+\algbrac{X}{\alpha}$ for some $\alpha \in \s{1}{}$, with $\algbrac{X}{\alpha}$ the algebraic bracket from \eqref{eq:proj-class-algbrac}.  A \emph{\proj\ structure} on $M$ is an equivalence class $\Dspace$ of \proj ly\ equivalent connections; $\algbrac{}{\alpha}$ defines an embedding $\s{1}{} \injto \s{1}{\alg{gl}{TM}}$, exhibiting the projective class $\Dspace$ as an affine space modelled on $\s{1}{}$.

\begin{lem} \thlabel{lem:proj-class-connsonL} There is a bijection between connections $\D\in \Dspace$ and connections on the line bundle $\cL \defeq (\Wedge{n}TM)^{2/(n+1)}$, where $\D,\b{\D}\in \Dspace$ are related by $\alpha\in\s{1}{}$ if and only if $\alpha$ is the change of induced connection on $\cL$. \end{lem}

\begin{proof} Firstly, note that since $\Wedge{n}TM$ is an oriented line bundle it admits oriented roots; thus $\cL$ is well-defined.  Now recall that an element of $\alg{gl}{TM}$ acts on $\Wedge{n}TM$ by its trace.  Using \eqref{eq:proj-class-algbrac}, this means that $\algbrac{X}{\alpha}$ acts on $\Wedge{n}TM$ by multiplication by $\tfrac{1}{2}(n+1)\alpha(X)$, and hence on $\cL$ by multiplication by $\alpha(X)$.  It follows that the difference of the connections induced on $\cL$ by $\D,\b{\D}$ is precisely $\alpha$. \end{proof}

We may also consider \proj\ structures $\Dspace$ whose connections have torsion.  The symmetry of $\algbrac{X}{\alpha} \acts Y$ in $X,Y$ implies that all connections in $\Dspace$ have the same torsion $\Tor{}$, which is then an invariant of the \proj\ structure.

\vspace{-0.25em}
\section{Description as a parabolic geometry} 
\label{s:proj-para}
\vspace{-0.15em}

Real \proj\ space $\RP$ may be described equivalently as either the set of lines through the origin in $\bR^{n+1}$, or as the quotient of the $n$-sphere $\Sph[n]$ by the $\bZ_2$-action induced by its antipodal map.  The sphere has a canonical \emph{round metric} for each signature $(p,q)$, given by pulling back the euclidean metric of signature $(p+1,q)$ on $\bR^{n+1}$ to $\Sph[n]$, \wrt\ which the antipodal map is a local isometry.  Thus $\RP$ inherits a canonical metric $g_{\mr{FS}}$ of constant curvature and signature $(p,q)$ which we call the \emph{Fubini--Study metric}.  Let $\D^{\mr{FS}}$ be the \LC\ connection of $g_{\mr{FS}}$.  It is straightforward to describe the geodesics of the \proj\ structure $[\D^{\mr{FS}}]$: the embedded \proj\ lines $\RP[1] \injto \RP$ are totally geodesic, so the geodesics of $[\D^{\mr{FS}}]$ are the smooth curves contained in such a line.  In an affine chart these curves lie within a line in $\bR^n$ so that, from the point of view of geodesics, $\RP$ is the natural compactification of $\bR^n$.  For this reason, $\RP$ is a good candidate for the flat model of \proj\ differential geometry.

Of course, \proj\ differential geometry was first described as a Cartan geometry by Cartan \cite{c1924-projaffine, c1924-projconns}, and the \proj\ equivalence problem was solved by Thomas \cite{t1925-projtractor} shortly afterwards.  A modern version of Cartan and Thomas' work may be found in \cite[\S 8]{s1997-cartan} and \cite{beg1994-thomasbundle} respectively.  Readable accounts of \proj\ differential geometry from the point of view of parabolic geometry may be found in \cite[\S 4.1.5]{cs2009-parabolic1} and \cite{br2007-W=0, e2008-projdiffgeom, em2008-projmetrics}.

\subsection{The flat model $\RP$} 
\label{ss:proj-para-flat}

As a generalised flag manifold, $\RP$ may be identified with $G/P$ for
\vspace{-0.2em}
\begin{equation*}
  G \defeq \grp{PGL}{n+1,\bR} \defeq \grp{GL}{n+1,\bR} / \linspan{\id}{}
\vspace{-0.2em}
\end{equation*}
and $P\leq G$ the (projection to $G$ of) the $\grp{GL}{n+1,\bR}$-stabiliser of a chosen line $\linspan{v_0}{}$ in $\bR^{n+1}$.  The Lie algebra of $G$ is the split real form $\fg \defeq \alg{sl}{n+1,\bR}$ of $\alg{sl}{n+1,\bC}$, while the parabolic subalgebra $\fp$ corresponds to crossing the last node:
\vspace{-0.1em}
\begin{equation} \label{eq:proj-para-satake}
\smash{
  \fp \defeq \dynkinSLRp{2}{0}{2}{proj-para-flat-p}
  \leq
  \dynkinSLR{2}{0}{2}{proj-para-flat-g} = \fg.
}
\vspace{-0.1em}
\end{equation}
Thus \proj\ differential geometry is an abelian parabolic geometry, with Killing polar $\fp^{\perp} \isom \bR^{n*}$ and reductive Levi factor $\fp^0 \defeq \fp/\fp^{\perp}$.  A choice of algebraic Weyl structure splits the $\fp^{\perp}$-filtration $\fg\supset\fp\supset\fp^{\perp}\supset 0$ of $\fg$, determining an isomorphism $\fg \isom \fg/\fp \oplus \fp^0 \oplus \fp^{\perp}$.  Evidently such a choice is equivalent to a choice of a subspace of $\bR^{n+1}$ complementary to $\linspan{v_0}{}$, thus yielding a decomposition
\begin{equation} \label{eq:proj-para-sldecomp}
  \alg{sl}{n+1,\bR} = \Setof{ \begin{pmatrix}
                                -\tr A & \alpha \\
                                X     & A
                              \end{pmatrix} }
                            {\, X\in\bR^n, \, \alpha\in\bR^{n*}, \, A\in\alg{gl}{n,\bR} }
\end{equation}
of $\fg$.  There is a corresponding decomposition of $\grp{PGL}{n+1,\bR}$ in which the Levi subgroup $P^0\defeq P/\exp{\fp^{\perp}}$ consists of (equivalence classes of) block-diagonal matrices in $G$; clearly the adjoint action of $P^0$ on $\fg/\fp\isom\bR^n$ induces an isomorphism $\fp^0\isom \alg{gl}{n,\bR}$.  Then $\fp^\perp\isom\bR^{n*}$ consists of matrices with only the $\alpha$-block, $\fp^0\isom \alg{gl}{n,\bR}$ of block-diagonal matrices, and $\fp\isom \alg{gl}{n,\bR}\ltimes \bR^{n*}$ of block upper-triangular matrices.

Suppose we have chosen an algebraic Weyl structure, so that $\fg\isom \bR^n\dsum \alg{gl}{n,\bR}\dsum \bR^{n*}$.  Using the description \eqref{eq:proj-para-sldecomp} we can write down the Lie bracket between elements of the three summands: if $X, Y\in\bR^n$, $A, B\in\alg{gl}{n,\bR}$ and $\alpha, \beta\in\bR^{n*}$ then
\vspace{-0.1em}
\begin{equation*}
  \liebrac{X}{Y} = 0 = \liebrac{\alpha}{\beta}, \quad
  \liebrac{A}{B} = AB-BA, \quad
  \liebrac{A}{X} = AX
  \quad\text{and}\quad
  \liebrac{A}{\alpha} = -\alpha\circ A,
\vspace{-0.1em}
\end{equation*}
so that $\bR^n,\bR^{n*}$ form abelian subalgebras on which $\alg{gl}{n,\bR}$ acts in the natural way; meanwhile $\bR^n$ and $\bR^{n*}$ bracket into $\fp^0\isom\alg{gl}{n,\bR}$ according to
\vspace{-0.1em}
\begin{equation} \label{eq:proj-para-gliebrac}
  [[X,\alpha],Y] = \alpha(X)Y + \alpha(Y)X.
\vspace{-0.1em}
\end{equation}
After appropriately normalising the inclusion $\liebrac{\bR^n}{\bR^{n*}} \injto \alg{gl}{n,\bR}$, \thref{prop:para-calc-change} implies that a change of Weyl structure is precisely a change of connection within the \proj\ class according to \eqref{eq:proj-class-algbrac}.  We will justify this normalisation in Chapter \ref{c:ppg}.

Finally, the Satake diagram \eqref{eq:proj-para-satake} and \thref{prop:lie-para-proj} tells us that there is a \proj\ embedding $G\acts\fp \injto \pr{\bV^*}$ for any irreducible $\fg$-representation $\bV$ whose highest weight is supported on the right-most node.  These representations are precisely the symmetric powers $\bV_k \defeq \Symm{k}\bR^{n+1*}$ of the \co standard\ representation, so that the embedding $\RP \injto \pr{\Symm{k}\bR^{n+1}}$ is the (degree $k$) Veronese embeddings defined by $[v] \mapsto [v \tens\cdots\tens v]$; see \cite{h2013-alggeom}.  Of course, the Veronese embeddings are not minimal for $k>1$, and Kostant's \thref{thm:lie-para-kostant} identifies $\RP$ with the intersection of quadrics given by projection away from the Cartan square in $\Symm{2}\bV_k$.  For later developments we are mostly interested in the case $k=2$, corresponding to the representation $\bW^* \defeq \bV_2 = \Symm{2}\bR^{n+1*}$.  By the description of the Cartan product in $\fg = \alg{sl}{n+1,\bR}$ from \cite{e2005-cartanprod}, the $\fg$-representation $\Symm{2}\bW^*$ decomposes into irreducible pieces as
\vspace{-0.1em}
\begin{equation*}
  \Symm{2}\bW
    = \Symm{2} \left( \dynkinSLR[2,0,0,0]{2}{0}{2}{proj-para-flat-W} \right)
    = \dynkinname{ \dynkinSLR[4,0,0,0]{2}{0}{2}{proj-para-flat-W1} }
                 { \Symm{4}\bR^{n+1} }
        \oplus
      \dynkinname{ \dynkinSLR[0,2,0,0,0]{3}{0}{2}{proj-para-flat-W2} }
                 { \Cartan{2}\Wedge{2}\bR^{n+1} },
\vspace{-0.1em}
\end{equation*}
where the first summand is the Cartan square.  Thus
\vspace{-0.1em}
\begin{equation*}
  \bU^* \defeq \dynkinSLR[0,0,0,2,0]{2}{0}{3}{proj-para-flat-U}
    = \Cartan{2}\Wedge{2}\bR^{n+1*}
\vspace{-0.1em}
\end{equation*}
is the space of homogeneous quadratic equations which cut out $\RP$ as an intersection of quadrics.  In particular, $\bU^*$ is an irreducible $\fg$-representation.

\subsection{Recovering the Cartan connection} 
\label{ss:proj-para-cartan}

Recall from Subsection \ref{ss:para-calc-equiv} that the general equivalence of categories result for a parabolic geometry requires that the first homology $\liehom{1}{\fg}$ has no component in homogeneity one.  We shall see why this condition is problematic for $\fg=\alg{sl}{n+1,\bR}$ shortly, but it is already possible to see that the underlying infinitesimal flag structure carries no information.  Indeed, the choice of an algebraic Weyl structure determines an isomorphism of the Levi subgroup $P^0 \defeq P/\exp \fp^{\perp}$ with $\grp{GL}{n,\bR}$.  Then the principal $P^0$-bundle $F^0\defeq F^P/\exp \fp^{\perp}$, from which one obtains the Cartan bundle $F^P$ by prolongation, is simply the full $\grp{GL}{n,\bR}$-frame bundle of $M$.  Of course, a principal $\grp{GL}{n,\bR}$-bundle with its solder form contains no geometric information.

To see why the homogeneity condition is problematic for $\alg{sl}{n+1,\bR}$, suppose that $\fg$ is a simple Lie algebra.  Choose a Cartan subalgebra $\ft$ and a positive subsystem $\Delta^{+}$ \wrt\ which $\beta \in \ft^*$ is the highest root of $\fg$, and suppose that $\fp \leq \fg$ is a standard parabolic corresponding to a subset $\Sigma \subseteq \Delta^0$ of simple roots.  It is proved in \cite[Prop.\ 3.3.7]{cs2009-parabolic1} that $\liehom{1}{\fg}$ can have no irreducible components in homogeneity one unless $\Sigma=\{\alpha\}$ consists of a single simple root; in this case the irreducible components are in bijection with the simple roots $\alpha_i$ for which the Cartan integers $\cartanint{\beta}{\alpha_i}$ and $\cartanint{\beta}{\alpha}$ are equal and \non zero.  The highest root of $\fg = \alg{sl}{n+1,\bR}$ is
\vspace{-0.1em}
\begin{equation} \label{eq:proj-para-gadjoint}
  \fg = \dynkinSLR[1,0,0,1]{2}{0}{2}{proj-para-cartan-adj},
\vspace{-0.1em}
\end{equation}
so that only the parabolics
\vspace{-0.1em}
\begin{equation*}
\smash{
  \dynkin{ \DynkinLine{0}{0}{1}{0};
           \DynkinDots{1}{0}{3}{0};
           \DynkinLine{3}{0}{4}{0};
           \DynkinXDot{0}{0};
           \DynkinWDot{1}{0};
           \DynkinWDot{3}{0};
           \DynkinWDot{4}{0};
         }{proj-para-cartan-homog1}
  \quad\text{and}\quad
  \dynkinSLRp{2}{0}{2}{proj-para-cartan-homog2}
}
\vspace{-0.1em}
\end{equation*}
have a \non trivial\ component of $\liehom{1}{\fg}$ of homogeneity one.  Of course, the resulting \Rspaces\ are dual \proj\ spaces and correspond to our choice from \eqref{eq:proj-para-satake}.

Given the problems above, we must work a little harder to achieve an equivalence of categories between \proj\ structures on $M$ and normal Cartan connections.  This entails constructing a so-called \emph{$P$-frame bundle of degree one} over $M$, which is a principal $P$-bundle $F^P$ over $M$ together with a $(\fg/\fp^{\perp})$-valued $1$-form $\theta$ such that:
\begin{itemize}
  \item $\theta$ is invariant under the $P$-action induced on $\fg/\fp^{\perp}$ by the adjoint action;
  
  \item $\ker \theta$ is the space of vector fields generated by the action of $\fp^{\perp}\subset \fp$; and
  
  \item $\theta$ maps the vertical bundle of $F^P$ to $\fp/\fp^{\perp}$, via $\eta(X^{\xi}) = \xi + \fp^{\perp}$ for all $\xi\in\fp$.
\end{itemize}
It turns out that a $P$-frame bundle of degree one contains just enough geometric information to obtain an equivalence of categories for \proj\ differential geometry.  We outline the proof of this equivalence, whose full proof may be found in \cite[Prop.\ 4.1.5]{cs2009-parabolic1}.

\begin{thm} \thlabel{thm:proj-para-equiv} There is an equivalence of categories between \proj\ structures on $M$ and $P$-frame bundles of degree one of type $G/P$.  A \proj\ structure is \torsionfree\ if and only if the corresponding $P$-frame bundle is normal. \end{thm}

\begin{sketchproof} Fix an algebraic Weyl structure for $\fg$, so that $\fg/\fp \isom \bR^n$, $\fp^0\isom \alg{gl}{n,\bR}$ and $\fp^{\perp}\isom\bR^{n*}$.  Suppose first that $\Dspace$ is a \proj\ structure on $M$ and consider the frame bundle $F^0 \defeq \grp{GL}{TM}$ of $M$ with its canonical solder form $\eta \in\s[F^0]{1}{\fg/\fp}$.  Then each $\D\in\Dspace$ induces a principal $\fp^0$-connection $1$-form $\gamma^{\D}$ on the full frame bundle $F^0$ of $M$.  We construct a principal $P$-bundle $F^P$ over $M$ with fibre $F^P_u \defeq \setof{\gamma^{\D}_u}{\D \in\Dspace}$ over $F^0$, with action of $g=g_0\exp\alpha \in P \isom P^0 \ltimes \exp{\fp^{\perp}}$ given by
\vspace{-0.2em}
\begin{equation} \label{eq:proj-para-cartanPaction}
  (\gamma^{\D}_u \acts g)(X) = \gamma^{\D}_{u\acts g_0}(X) + \liebrac{\alpha}{\eta(X)}.
\vspace{-0.2em}
\end{equation}
One then checks that $F^P$ is indeed a principal $P$-bundle over $M$, and for each $\gamma^{\D}_u \in F^P_u$ the $1$-form $\eta+\gamma^{\D}_u$ pulls back to a $P$-frame form $\theta$ on $F^P$.  Moreover, as we vary $\alpha \in \fp^{\perp}$ in \eqref{eq:proj-para-cartanPaction} the $\gamma^{\D}\acts g$ run over all connections in the \proj\ class.

Conversely, suppose that $(F^P,\theta)$ is a $P$-frame bundle of degree one.  The quotient bundle $F^0 \defeq F^P/\exp\fp^{\perp}$ is the frame bundle of $M$, and projecting the values of $\theta$ onto $\fg/\fp$ yields the solder form $\eta$ on $F^0$.  Choosing a local section $\sigma : U \to F^P$, it turns out that $\sigma$ pulls back the $\fp^0$-component of $\theta$ to a principal $\fp^0$-connection $\gamma^{\sigma}$ on $F^0\at{U}$.  Replacing the local section $\sigma$ by $\b{\sigma} = \sigma\acts g_0\exp \alpha$ for local functions $g_0:U\to P^0$ and $\alpha : U \to \fp^{\perp}$, the principal connections on $F^0$ change according to $\gamma^{\b{\sigma}} = \gamma^{\sigma} + \liebrac{}{\alpha}$.  Equation \eqref{eq:proj-para-gliebrac} then implies that the space of linear connections determined by all such $\sigma$ constitute a \proj\ structure $\Dspace$ over $U\subset M$.  Patching these local \proj\ structures together yields a \proj\ structure on $M$.

If $\gamma^{\D}$ is the principal connection $1$-form corresponding to $\D\in\Dspace$, viewed as the $\fp^0$-component of the $P$-frame form $\theta$, its torsion coincides with the torsion of the $P$-frame bundle.  A $P$-frame bundle is normal if and only if its torsion lies in the kernel of the principal bundle map induced by the Lie algebra homology map $\p : \Wedge{2}\fp^{\perp}\tens\fg \to \fp^{\perp}\tens\fg$, which turns out to be injective for this choice of $\fg$.  Thus the $P$-frame bundle is normal if and only if the torsion of the \proj\ structure vanishes. \end{sketchproof}

\subsection{Representations of $\alg{sl}{n+1,\bR}$} 
\label{ss:proj-para-repns}

For later use let us describe some important $\fg$- and $\fp$-representations, as well as their associated bundles.  The adjoint representation of $\fg$ has highest weight given by \eqref{eq:proj-para-gadjoint}, so that the isotropy representation $\fg/\fp$ and its dual $(\fg/\fp)^*\isom\fp^{\perp}$ have highest weights
\vspace{-0.05em}
\begin{equation*}
  \fg/\fp = \dynkinSLRp[1,0,0,1]{2}{0}{2}{proj-para-repns-gp}
  \quad\text{and}\quad
  \fp^{\perp} = \dynkinSLRp[0,0,0,1,-2]{2}{0}{3}{proj-para-repns-pperp}
\vspace{-0.05em}
\end{equation*}
respectively as $\fp$-representations.  By the Cartan condition, the corresponding associated bundles are $TM$ and $T^*M$, and a Weyl structure gives an isomorphism $\fg_M\isom TM\dsum \alg{gl}{TM}\dsum T^*M$.  The top exterior power of $\fg/\fp$ is the derivative of the group character $\det:\grp{GL}{n,\bR} \to \bR$, so its highest weight has coefficient $n+1$ over the crossed node; it follows that the line bundle $\cL\defeq(\Wedge{n}TM)^{2/(n+1)}$ of \thref{lem:proj-class-connsonL} is associated to the $1$-dimensional representation
\vspace{0.225em}
\begin{equation} \label{eq:proj-para-repnL}
  L = \dynkinSLRp[0,0,0,2]{2}{0}{2}{proj-para-repns-L}.
\vspace{0.225em}
\end{equation}
Note that $L$ is the zeroth homology $\liehom{0}{\bW^*}$, where $\bW \defeq \Symm{2}\bR^{n+1}$ is the $\fg$-represent-ation from Subsection \ref{ss:proj-para-flat}.  Since the other fundamental representations of $\fg$ are exterior powers of the standard representation, every tensor bundle on $M$ can be associated to a particular $\fp$-representation.

We are yet to describe another obvious representation of $\fg$: the standard representation on $\bT \defeq \bR^{n+1}$, which has highest weight
\vspace{0.225em}
\begin{equation} \label{eq:proj-para-stdrepn}
  \bT = \dynkinSLR[1,0,0,0]{2}{0}{2}{proj-para-repns-T}.
\vspace{0.225em}
\end{equation}
The reason for this omission is that due to our choice of group $G = \grp{PGL}{n+1,\bR}$ with Lie algebra $\fg=\alg{sl}{n+1,\bR}$, not all $\fg$-representation integrate to globally defined representations of $G$.  Indeed, if $n+1$ is even then the connected component of the identity in $\tilde{G} \defeq \grp{SL}{n+1,\bR}$ is a double cover \cite{s1997-cartan} of $G = \grp{PGL}{n+1,\bR}$, and a $\fg$-representation integrates to $G$ if and only if it integrates to $\tilde{G}$ with the element $-\id\in\tilde{G}$ acting trivially.  If the representation in question has highest weight $\lambda = \sum{i}{} \lambda_i\omega_i$, written here in terms of the fundamental weights, $-\id$ acts trivially if and only if $\sum{i}{} \lambda_i$ is even \cite{fh1991-repntheory}, thus giving a simple integrability criterion which of course fails for $\bT$.  Note that there is no problem when $n+1$ is odd, since then $G=\tilde{G}$.

We can still form bundles associated to $\fg$-representations which do not integrate to $G$ by locally extending the Cartan bundle to $\tilde{G}$.  More precisely, we let $\tilde{P} \leq \tilde{G}$ denote the stabiliser of a chosen line in $\bR^{n+1}$ (previously denoted $v_0$), which corresponds to crossing the last node of the Satake diagram.  We then form the extended Cartan bundle $F^{\tilde{P}} = \assocbdl{F^P}{\tilde{P}}{\tilde{P}}$ with structure group $\tilde{P}$, and there is a canonical extension of the unique normal Cartan connection on $F^P$ to a normal Cartan connection on $F^{\tilde{P}}$.  We can then integrate all $\fg$-representations to $\tilde{G}$, form the associated bundles, and quotient by the $\bZ_2$-action coming from the double cover $F^{\tilde{P}}\surjto F^P$ if necessary.  Note that $\tilde{G}/\tilde{P} \isom \Sph[n]$ and, by modifying the proof of \thref{thm:proj-para-equiv} slightly \cite[Prop.\ 4.1.5]{cs2009-parabolic1}, we obtain an equivalence of categories between \emph{oriented} \proj\ structures on $M$ and $\tilde{P}$-frame bundles of degree one of type $\tilde{G}/\tilde{P}$.

For the standard representation \eqref{eq:proj-para-stdrepn}, the associated bundle $\cT \defeq \assocbdl{F^{\tilde{P}}}{\tilde{P}}{\bT}$ is a rank $n+1$ vector bundle called the \emph{standard tractor bundle}.  Since an element $\inlinematrix{ (\det C)^{-1} & \alpha \\ 0 & C } \in \tilde{P}$ acts on the line $\bT_0 \defeq \linspan{v_0}{}$ in $\bT$ stabilised by $\tilde{P}$ by multiplication by $(\det C)^{-1}$,
\vspace{-0.15em}
\begin{equation*}
  \assocbdl{F^{\tilde{P}}}{\tilde{P}}{\bT_0}
    \isom \dynkinSLRp[0,0,0,-1]{2}{0}{2}{proj-para-repns-T0}
    = \cL^{-1/2}
\vspace{-0.15em}
\end{equation*}
is a square root of the line bundle $\cL^*$ from \thref{lem:proj-class-connsonL}.  The relations $\bT_0 = \fp^{\perp} \acts \bT$ and $\fp^{\perp} \acts \bT_0 = 0$ exhibit $\bT_0$ as the socle of the height one $\fp^{\perp}$-filtration $\bT \supset \bT_0 \supset 0$ of $\bT$; the decomposition \eqref{eq:proj-para-sldecomp} of $\fg$ implies that the top $\bT / \bT_0 = \liehom{0}{\bT}$ has associated bundle
\vspace{-0.15em}
\begin{equation*}
  \assocbdl{F^{\tilde{P}}}{\tilde{P}}{(\bT/\bT_0)}
  \isom \dynkinSLRp[1,0,0,0]{2}{0}{2}{proj-para-repns-TT0}
  = \cL^{-1/2}\tens TM.
\vspace{-0.15em}
\end{equation*}
In particular, a Weyl structure yields a decomposition $\cT \isom (\cL^{-1/2}\tens TM)\dsum\cL^{-1/2}$.

\begin{rmk} \thlabel{rmk:proj-para-tractorjet} In fact, the extension of the Cartan bundle from $F^P$ to $F^{\tilde{P}}$ is equivalent to a choice of square root of $\cL^*$ \cite{beg1994-thomasbundle, cg2000-irrtractor}.  Indeed, the $1$-jet bundle $\jetbdl{1}{\cL^{-1/2}}$ of $\cL^{-1/2}$ fits into the canonical short exact sequence
\vspace{-0.15em}
\begin{equation*}
  \ses{ \cL^{-1/2} \tens T^*M }
      { \jetbdl{1}{\cL^{-1/2}} }
      { \cL^{-1/2} },
\vspace{-0.15em}
\end{equation*}
a splitting of which is equivalent to a linear connection on $\cL^{-1/2}$.  Choosing such a splitting yields an isomorphism $\jetbdl{1}{\cL^{-1/2}} \isom (\cL^{-1/2}\tens T^*M) \dsum \cL^{-1/2}$, so that we recover the standard tractor bundle by defining $\cT \defeq \jetbdl{1}{\cL^{-1/2}}$.  It remains to recover the normal Cartan connection from the tractor connection on $\cT$, which is the subject of Thomas' work \cite{t1925-projtractor}. \end{rmk}

\subsection{Harmonic curvature} 
\label{ss:proj-para-harm}

Having identified the tensor bundles, we may analyse the harmonic curvature of the canonical Cartan connection, which lies in the Lie algebra homology $\liehom{2}{\fg}$.  The Hasse diagram computing this homology is given in Figure \ref{fig:proj-para-hasse}, so that the harmonic curvature has a single irreducible component lying in
\vspace{-0.05em}
\begin{equation*}
  \dynkinSLRp[1,0,0,1,1,-4]{2}{0}{4}{proj-para-harm-weyl}
    = \Wedge{2}T^*M \cartan \alg{sl}{TM},
\vspace{-0.05em}
\end{equation*}
where $\cartan$ is the Cartan product and $\alg{sl}{TM} = T^*M\cartan TM$.  This single piece may be identified with the Weyl curvature $\Weyl{}{}$, which is a totally \tracefree\ and \proj ly\ invariant $\alg{sl}{TM}$-valued $2$-form.  This is as we expect: normality of the Cartan connection implies that the Weyl connections are \torsionfree, so there is no torsion component of the harmonic curvature.

\begin{figure}[t]
  \begin{equation*}
  \arraycolsep=0.1em
  \begin{array}{ccccccc}
    \dynkinSLRp[1,0,0,1]{2}{0}{2}{proj-para-harm-hasse1}
    &
    \dynkin{ \DynkinConnector{0.2}{0}{1.8}{0}; }{proj-para-harm-hasse2}
    &
    \dynkinSLRp[1,0,0,2,-3]{2}{0}{3}{proj-para-harm-hasse3}
    &
    \dynkin{ \DynkinConnector{0.2}{0}{1.8}{0}; }{proj-para-harm-hasse4}
    &
    \dynkinSLRp[1,0,0,1,1,-4]{2}{0}{4}{proj-para-harm-hasse5}
    &
    \dynkin{ \DynkinConnector{0.2}{0}{1.8}{0};
             \DynkinLabel{$\cdots$}{3}{-0.75}; }{proj-para-harm-hasse6}
  \end{array}
  \end{equation*}
  \caption[The Hasse diagram computing $\liehom{}{\alg{sl}{n+1,\bR}}$]
          {The Hasse diagram of the adjoint representation of $\fg = \alg{sl}{n+1,\bR}$, which computes the homology $\liehom{}{\fg}$.}
  \label{fig:proj-para-hasse}
\end{figure}

Since \proj\ differential geometry is abelian, \thref{thm:para-bgg-curv} states that the curvature any Weyl connection $\D$ decomposes as $\Curv{}{} = \Weyl{}{} - \algbracw{\id}{\nRic{}{}}{}$, where $\nRic{}{}$ is the normalised Ricci tensor defined by $\nRic{}{} \defeq -\quab_M^{-1}\liebdy \Curv{}{}$.  Recalling that the Cartan curvature vanishes if and only if the harmonic curvature vanishes, we recover the well-known classical result which states that a \proj\ manifold is locally diffeomorphic to $\RP$ (\ie\ is \proj ly\ flat) if and only if its Weyl curvature vanishes \cite{br2007-W=0, w1921-projflat}.  Of course, $(\p\Curv{}{})_X(Y) = \ve^i(\Curv{e_i,X}{Y})$ is simply the Ricci curvature $\sRic{}{}$ of $\D$.  The following result is well-known (see for example \cite{e2008-projdiffgeom}) and included for completeness.

\begin{prop} \thlabel{prop:proj-para-nric} The normalised Ricci tensor $\nRic{}{}$ of $\D\in\Dspace$ is given by
\begin{equation} \label{eq:proj-para-nric}
  \nRic{}{}
    = -\tfrac{2}{n-1} \sym (\p\Curv{}{}) - \tfrac{2}{n+1} \alt (\p\Curv{}{}).
\end{equation}
In particular if $\D$ has a nowhere-vanishing parallel section of $\cL$, then $\nRic{}{} = -\tfrac{2}{n-1} \liebdy \Curv{}{}$. \end{prop}

\begin{proof} Applying the Lie algebra boundary map to $\Curv{}{}=\Weyl{}{} - \algbracw{\id}{\nRic{}{}}{}$ yields
\begin{equation*}
  \liebdy\Curv{}{}
    = -\tfrac{n-1}{2}\sym{\nRic{}{}} - \tfrac{n+1}{2}\alt{\nRic{}{}},
\end{equation*}
so that $\sym{\nRic{}{}} = -\tfrac{2}{n-1} \sym(\liebdy \Curv{}{})$ and $\alt{\nRic{}{}} = -\tfrac{2}{n+1} \alt(\liebdy \Curv{}{})$.  Then \eqref{eq:proj-para-nric} follows by writing $\nRic{}{} = \sym{\nRic{}{}} + \alt{\nRic{}{}}$.  If $\ell \in \s{0}{\cL}$ is nowhere-vanishing and $\D$-parallel, we have $\Curv{}{\ell} = 0$ and, since $\Weyl{}{}$ is totally \tracefree, $\Curv{X,Y}{\ell} = 0 = \algbracw{\id}{\nRic{}{}}{X,Y} \acts\ell = \nRic{Y}{X}\ell - \nRic{X}{Y}\ell$.  Since $\ell$ is nowhere-vanishing we conclude that $\nRic{}{}$ is symmetric; the last claim now follows easily. \end{proof}

Recall the Cotton--York tensor $\CY{}{} \defeq \d^{\D}\nRic{}{}$ of $\D$, which is a $T^*M$-valued $1$-form.  The following identities involving $\Weyl{}{}$ and $\CY{}{}$ are well-known in \proj\ differential geometry (see for example \cite{em2008-projmetrics}) and will be useful in Section \ref{s:proj-bgg}.

\begin{prop} \thlabel{prop:proj-para-calc} Let $\D\in\Dspace$ be a Weyl connection.  Then:
\begin{enumerate}
  \item \label{prop:proj-para-calc-1}
  There are Bianchi identities
  \begin{equation*} \begin{aligned}
    \Weyl{X,Y}{Z} + \Weyl{Y,Z}{X} + \Weyl{Z,X}{Y} &= 0 \\
    \text{and}\quad
    \CY{X,Y}{Z} + \CY{Y,Z}{X} + \CY{Z,X}{Y} &= 0
  \end{aligned}
  \end{equation*}
  
  \item \label{prop:proj-para-calc-2}
  $\ve^i(\D_{e_i}\Weyl{X,Y}{}) = -\tfrac{1}{2}(n-2) \CY{X,Y}{}$ \wrt\ any local frame $\{e_i\}_i$ of $TM$ with dual coframe $\{\ve^i\}_i$.
\end{enumerate}
\end{prop}

\begin{proof} \proofref{prop:proj-para-calc}{1} Since $\algbrac{X}{\alpha} \acts Y$ is symmetric in $X,Y$, the algebraic Bianchi identity for $\Curv{}{}$ gives the desired Bianchi identity for $\Weyl{}{}$.  The differential Bianchi identity $\d^{\D} \Curv{}{} = 0$ yields $\d^{\D}\Weyl{}{} + \algbracw{\id}{\CY{}{}}{} = 0$ and hence
\vspace{-0.2em}
\begin{equation} \label{eq:proj-para-diffbianchi} \begin{aligned}
  &(\D_X\Weyl{}{})_{Y,Z} + (\D_Y\Weyl{}{})_{Z,X} + (\D_Z\Weyl{}{})_{X,Y} \\
    & \hspace{4em}
    = -\algbrac{X}{\CY{Y,Z}{}} - \algbrac{Y}{\CY{Z,X}{}} - \algbrac{Z}{\CY{X,Y}{}}
\end{aligned}
\vspace{-0.2em}
\end{equation}
for all $X,Y,Z\in\s{0}{TM}$.  Taking the trace of both sides, terms of the form $\ve^i\left( (\D_X\Weyl{}{})_{Y,Z}\acts e_i \right)$ vanish since $\Weyl{}{}$ is totally \tracefree.  Since $\ve^i\big( \algbrac{X}{\CY{Y,Z}{}} \acts e_i \big) = \tfrac{1}{2}(n+1) \CY{Y,Z}{X}$, the Bianchi identity for $\CY{}{}$ follows.  

\smallskip

\proofref{prop:proj-para-calc}{2} Taking a trace over $X$ in \eqref{eq:proj-para-diffbianchi} and using that $\liebdy \Weyl{}{}=0$, we obtain
\vspace{-0.2em}
\begin{equation*}
  \ve^i\left( \D_{e_i}\Weyl{X,Y}{} \right)
    = -\ve^i \circ\left( \algbrac{e_i}{\CY{X,Y}{}} + \algbrac{X}{\CY{Y,e_i}{}}
      + \algbrac{Y}{\CY{e_i,X}{}} \right).
\vspace{-0.2em}
\end{equation*}
Evaluating on a vector field $Z$ and expanding the algebraic brackets, this simplifies to
\vspace{-0.2em}
\begin{equation*}
  \ve^i\left( \D_{e_i}\Weyl{X,Y}{} \right) \acts Z
    = -\tfrac{1}{2}\left( (n-1)\CY{X,Y}{Z} + \CY{Y,Z}{X} + \CY{Z,X}{Y} \right),
\vspace{-0.2em}
\end{equation*}
which yields the desired result by the Bianchi identity for $\CY{}{}$. \end{proof}

\vspace{-0.15em}
\section{Metrisability of \proj\ structures} 
\label{s:proj-bgg}

We say that the \proj\ class $\Dspace$ is \emph{metrisable} if it contains a metric connection, in which case its geodesics are those of a (\pseudo\riem) metric.

Suppose that $g\in\s{0}{\Symm{2} T^*M}$ is a metric.  We would like to reduce the question of whether the \LC\ connection of $g$ lies in the \proj\ class $\Dspace$ to a differential equation on $\D$, as we did for the endomorphism $A$ defined by \eqref{eq:proj-class-endoA} in the classical description.  For this, observe that there is a natural decomposition
\vspace{-0.2em}
\begin{equation} \label{eq:proj-bgg-tracesummands}
  T^*M \tens \Symm{2}TM = (\id\symm TM) \dsum (T^*M\tens[\trfree] \Symm{2}TM)
\vspace{-0.2em}
\end{equation}
into the image and kernel of the natural trace $T^*M\tens\Symm{2}TM \to TM$.  This suggests that a differential system involving successive differentials and traces of the inverse metric $g^{-1} \in \s{0}{\Symm{2}TM}$ may close.  This is indeed the case \cite{em2008-projmetrics}, and it is a matter of tensoring with an appropriate weight to achieve \proj\ invariance.

\begin{prop} \thlabel{prop:proj-bgg-metriceqn} The first-order linear differential equation $(\D h)_{\trfree}=0$ is projectively invariant on sections of $\cL^*\tens \Symm{2}TM$, where the subscript \emph{``$\trfree$''} denotes the trace-free part. \end{prop}

\begin{proof} We calculate the variation $\weyld{\gamma}\D h = \algbrac{\bdot}{\gamma} \acts h$ \wrt\ $\gamma\in\s{1}{}$.  Viewing $h$ as a bundle map $T^*M \to \cL^*\tens TM$, the Leibniz rule gives
\begin{equation*} \begin{aligned}
  (\algbrac{X}{\alpha} \acts h)(\beta)
  &= \algbrac{X}{\alpha} \acts h(\beta,\bdot)
    - h(\algbrac{X}{\alpha} \acts \beta,\bdot) \\
  &= -\alpha(X)h(\beta,\bdot) + \tfrac{1}{2}( \alpha(X)h(\beta,\bdot)
    + \alpha(h(\beta,\bdot))X ) \\
  &\quad\qquad + \tfrac{1}{2}h(\alpha(X)\beta+\beta(X)\alpha, \bdot) \\
  &= (X\symm h(\alpha,\bdot))(\beta).
\end{aligned}
\end{equation*}
Thus $\weyld{\gamma} \D h$ lies in the summand $\id\symm TM$ of \eqref{eq:proj-bgg-tracesummands}, giving \proj\ invariance. \end{proof}

The equation $(\D h)_{\trfree}=0$ will be called the \emph{linear metric equation}, which may equivalently be written as
\begin{equation} \label{eq:proj-bgg-metriceqnZ}
  \D h = \id\symm Z^{\D}
\end{equation}
for some section $Z^{\D}$ of $\cL^*\tens TM$ depending on $\D$.  Taking a trace in \eqref{eq:proj-bgg-metriceqnZ} easily yields $Z^{\D} = \tfrac{2}{n+1}\p(\D h)$.  Moreover, the proof of \thref{prop:proj-bgg-metriceqn} implies that $\weyld{\gamma}Z^{\D} = h(\gamma,\bdot)$, \ie\ $Z^{\D}\mapsto Z^{\D}+h(\gamma,\bdot)$ under change of Weyl connection $\D\mapsto \D + \algbrac{}{\gamma}$.  \thref{prop:proj-bgg-metriceqn} should be compared with \thref{prop:proj-class-maineqn}; this comparison will be made more precise in \thref{rmk:proj-bgg-linearmain} below.

Since $\Wedge{n}TM \isom \cL^{(n+1)/2}$, a metric $g$ induces a section of $\cL^*\tens\Symm{2}TM$ defined by
\begin{equation*}
  h \defeq (\det g)^{1/(n+1)} \ltens g^{-1},
\end{equation*}
which we call the \emph{linear metric} associated to $g$.  Then $\det h = (\det g)^{-1/(n+1)} \in \s{0}{\cL^*}$, so that we may recover $g=(\det h)^{-1} \ltens h^{-1}$ from $h$; \cf\ equation \eqref{eq:proj-class-gbar}.

\begin{cor} \thlabel{cor:proj-bgg-metriceqn} There is a bijection between \non degenerate\ solutions of the linear metric equation and metric connections in $\Dspace$. \end{cor}

\begin{proof} If $h$ is a \non degenerate\ solution of the linear metric equation \eqref{eq:proj-bgg-metriceqnZ} with $\D h = \id\symm Z^{\D}$, then $h^{-1}(Z^{\D},\bdot) \in\s{1}{}$ and hence $\D^g \defeq \D - \algbrac{}{h^{-1}(Z^{\D},\bdot)} \in \Dspace$ is independent of $\D$ and satisfies $\D^g h=0$.  It follows that $\D^g$ is the \LC\ connection of $g \defeq (\det h)^{-1} \ltens h^{-1}$.  Conversely if $\D^g\in \Dspace$ is the \LC\ connection of some metric $g$, we have $\D^g h=0$ for its linear metric $h \defeq (\det g)^{1/(n+1)} \ltens g^{-1}$.  Since the linear metric equation is projectively invariant, this $h$ is a solution for any $\D \in \Dspace$. \end{proof}

\thref{cor:proj-bgg-metriceqn} reduces the metrisability problem for a \proj\ structure to the study of a \proj ly\ invariant first-order linear differential equation.  We shall henceforth refer to all solutions of the linear metric equation as linear metrics.

\begin{rmk} \thlabel{rmk:proj-bgg-linearmain} Given a ``background'' metric $g$ and a \proj ly\ equivalent metric $\b{g}$, the associated linear metrics $h,\b{h}$ define an endomorphism $A$ of $T^*M$ by $\b{h}=h(A\bdot,\bdot)$.  Evidently (the transpose of) $A$ is precisely the endomorphism \eqref{eq:proj-class-endoA} featured in the main equation \eqref{eq:proj-class-maineqn}, clarifying the relation between the main equation and the linear metric equation.  We shall return to this topic in Chapter \ref{c:mob2}. \end{rmk}

The linear metric equation is over-determined, so its solution space may be prolonged to a closed differential system on an auxiliary bundle.  In the presence of a background metric, a (\non invariant) prolongation was obtained by Mike{\v s} \cite{m1996-geodprol} (see also \cite[p.\ 151]{hmv2009-geodmappings}).  An invariant prolongation was found by Eastwood and Matveev \cite{em2008-projmetrics}, which we translate into index-free notation below.

\begin{thm} \thlabel{thm:proj-bgg-metricprol} There is a linear isomorphism between the space of solutions of the metric equation and the parallel sections of the projectively invariant connection 
\vspace{-0.5em}
\begin{equation} \label{eq:proj-bgg-metricprol}
  \D^{\cW}_X \! \colvect{ h \\ Z \\ \lambda } =
  \colvect{ \D_X h - X\symm Z \\
            \D_X Z - h(\nRic{X}{}, \bdot) - \lambda X \\
            \D_X \lambda - \nRic{X}{Z} }
  - \tfrac{2}{n} \!
  \colvect{ 0 \\
            -\Weyl{e_i,X}{h(\ve^i,\bdot)} \\
            h(\CY{e_i,X}{},\ve^i) }
\vspace{-0.5em}
\end{equation}
on sections $(h,Z,\lambda)$ of $\cW \defeq (\cL^*\tens \Symm{2}TM)\dsum (\cL^*\tens TM)\dsum \cL^*$. \end{thm}

\begin{proof} Suppose that $h\in\s{0}{\cL^*\tens \Symm{2}TM}$ is a solution of the linear metric equation \eqref{eq:proj-bgg-metriceqnZ}, so that $\D h = \id \symm Z^{\D}$ for some $Z^{\D}\in \s{0}{\cL^*\tens TM}$.  Recalling that normality implies zero torsion, differentiating again and using the Ricci identity gives
\vspace{-0.5em}
\begin{equation} \label{eq:proj-bgg-metricprol-1}
  \Weyl{X,Y}{h} = \Curv{X,Y}{h} + \algbracw{\id}{\nRic{}{}}{X,Y}\acts h
                = - X\symm Q^{\D}_Y + Y\symm Q^{\D}_X,
\vspace{-0.5em}
\end{equation}
where $Q^{\D}_X \defeq \D_X Z^{\D} - h(\nRic{X}{},\bdot)$.  Taking a trace over $X$ in \eqref{eq:proj-bgg-metricprol-1} \wrt\ a local frame $\{e_i\}_i$ with dual frame $\{\ve^i\}_i$ gives
\vspace{-0.5em}
\begin{equation*}
  \Weyl{e_i,Y}{h(\ve^i,\bdot)}
    = \tfrac{1}{2}\left( -nQ^{\D}_Y - Q^{\D}_Y + Q^{\D}_Y + (\tr Q^{\D})Y \right)
    = \tfrac{n}{2}\left( -Q^{\D}_Y + \lambda^{\D}Y \right),
\vspace{-0.5em}
\end{equation*}
where $\lambda^{\D} \defeq \tfrac{1}{n}(\tr Q^{\D}) \in\s{0}{\cL^*}$.  Rearranging and substituting for $Q^{\D}$ then gives the second slot of $\D^{\cW}$.  Differentiating again and skew-symmetrising, we obtain
\vspace{-0.5em}
\begin{equation*}
\begin{split}
  &\Weyl{X,Y}{Z^{\D}} - \algbracw{\id}{\nRic{}{}}{X,Y}{Z^{\D}} - h(\CY{X,Y}{},\bdot) \\
  &\qquad
  \begin{aligned}
    &= (X\symm Z^{\D})(\nRic{Y}{},\bdot) - (Y\symm Z^{\D})(\nRic{X}{},\bdot)
        + (\D_X\lambda^{\D})Y - (\D_Y\lambda^{\D})X \\
    &\qquad - \tfrac{2}{n} \big( (\D_X \Weyl{e_i,Y}{}) \acts h(\ve^i,\bdot)
      + \Weyl{e_i,Y}{ \D_X h(\ve^i,\bdot) } \big) \\
    &\qquad + \tfrac{2}{n} \big( (\D_Y \Weyl{e_i,X}{}) \acts h(\ve^i,\bdot)
      + \Weyl{e_i,X}{ \D_X h(\ve^i,\bdot) } \big)
  \end{aligned}
\end{split}
\vspace{-0.5em}
\end{equation*}
Expanding the algebraic bracket on the \lhs\ and tracing over $X$ yields
\vspace{-0.5em}
\begin{equation} \label{eq:proj-bgg-metricprol-2} \begin{split}
  &(n-1)\left( \D_Y\lambda^{\D} - \nRic{Y}{Z^{\D}} \right) - h(\CY{e_i,Y}{},\ve^i) \\
  &\qquad
  \begin{aligned}
    =& -\tfrac{2}{n}\ve^j\big(
      (\D_{e_j}\Weyl{e_i,Y}{})\acts h(\ve^i,\bdot)
      - (\D_Y\Weyl{e_i,e_j}{})\acts h(\ve^i,\bdot) \\
    &\phantom{-\tfrac{2}{n}\ve^j \big((}
      \vphantom{a} + \Weyl{e_i,Y}{\D_{e_j} h(\ve^i,\bdot)}
      - \Weyl{e_i,e_j}{\D_Y h(\ve^i,\bdot)} \big).
  \end{aligned}
\end{split}
\end{equation}
On the \rhs\ the second and fourth terms vanish, because $\Weyl{e_i,e_j}{}$ is skew in $i,j$ while $h(\ve^i,\ve^j)$ is symmetric in $i,j$; the third term also vanishes by straightforward calculation, using that $\Weyl{}{}$ acts trivially on $\cL^*$.  The first term equals $\tfrac{n-2}{n} h(\CY{e_i,Y}{},\ve^i)$ by \itemref{prop:proj-para-calc}{2}, so that \eqref{eq:proj-bgg-metricprol-2} becomes
\begin{equation*}
  (n-1)\left( \D_Y\lambda^{\D} - \nRic{Y}{Z^{\D}} \right)
    = \tfrac{2n-2}{n} h(\CY{e_i,Y}{},\ve^i).
\end{equation*}
Dividing by $n-1$ and rearranging now gives the final slot of $\D^{\cW}$.

The \proj\ invariance of $\D^{\cW}$ can be checked directly: the first piece in \eqref{eq:proj-bgg-metricprol} is \proj ly\ invariant because $\weyld{\gamma}Z^{\D} = h(\gamma,\bdot)$ and $\weyld{\gamma}\lambda^{\D} = \gamma(Z^{\D})$; \proj\ invariance of the curvature correction follows from \itemref{thm:para-bgg-ablcurv}{weyld}. \end{proof}

\smallskip

A few remarks and observations regarding \thref{thm:proj-bgg-metricprol} are in order.  We may write the bundle $\cW \defeq (\cL^*\tens \Symm{2}TM)\dsum (\cL^*\tens TM) \dsum \cL^*$ in the form
\begin{equation*}
  \cW \isom \Symm{2}\big( (\cL^{-1/2}\tens TM)\dsum \cL^{-1/2}\, \big)
    \isom \Symm{2}\cT
\end{equation*}
where $\cT$ is the standard tractor bundle from Subsection \ref{ss:proj-para-repns}.  Therefore $\cW$ is associated to the representation $\bW \defeq \Symm{2}\bR^{n+1}$ whose projectivisation appears as the codomain of the Veronese embedding $\RP[n] \injto \pr{\Symm{2} \bR^{n+1}}$.  Moreover, the first piece of the prolongation connection $\D^{\cW}$ in \eqref{eq:proj-bgg-metricprol} may be reconciled with the connection induced on $\Symm{2}\cT$ by the standard tractor connection; see \cite[Eqn.\ (2.8)]{gm2015-weylnullity}.  The point is that the prolongation connection $\D^{\cW}$ is, up to a curvature correction, the tractor connection on the bundle associated to the projective embedding of the flat model.

The observations above may be understood in terms of BGG operators.  The first BGG operator associated to the representation $\bW$ is a differential operator
\begin{equation*}
  \bgg{\bW} :
    \dynkinname{ \dynkinSLRp[2,0,0,0,0]{2}{0}{3}{proj-bgg-metricdom} }
               { L^*\tens\Symm{2}TM }
  \To
    \dynkinname{ \dynkinSLRp[2,0,0,1,-2]{2}{0}{3}{proj-bgg-metricim} }
               {\zbox{ (L^*\tens\Symm{2}TM) \cartan T^*M }},
\end{equation*}
which is easily seen to be first-order by using the last row of the inverse Cartan matrix.  In this picture, the prolongation connection $\D^{\cW}$ is described by the general scheme of Subsection \ref{ss:para-bgg-prol}.  In particular, the BGG operator is finite type and has a linear solution space whose dimension, called the \emph{mobility} of $\Dspace$, is bounded above by $\dim(\Symm{2}\bR^{n+1}) = \tfrac{1}{2}(n+1)(n+2)$; this was known to Sinjukov \cite{s1979-geodmappings}.

However, something special happens here which does not happen for generic first BGG operators:  the graded $\fg$-representation $\fh \defeq \bW \dsum (\fg\dsum\bR) \dsum \bW^*$ may be written
\begin{equation*}
  \fh = \Symm{2}\bR^{n+1} \dsum (\bR^{n+1}\tens\bR^{n+1*}) \dsum \Symm{2}\bR^{n+1*}
    \isom \Symm{2}(\bR^{n+1}\dsum\bR^{n+1*}).
\end{equation*}
This is precisely the adjoint representation of the Lie algebra $\alg{sp}{\bR^{n+1} \dsum \bR^{n+1*},\omega}$, where $\omega\in\Wedge{2}(\bR^{n+1}\dsum\bR^{n+1*})$ is the canonical symplectic form defined by $\omega(u+\alpha,v+\beta) \defeq \alpha(v) - \beta(u)$, thus giving an isomorphism of $\fh$ with the graded vector space underlying $\alg{sp}{2n+2,\bR}$.  The graded Lie algebra structure on $\alg{sp}{2n+2,\bR}$ induces a graded Lie algebra structure on $\fh$ in which the summands $\bW$ and $\bW^*$ are abelian subalgebras, while the grading implies that $\fq \defeq (\fg\dsum\bR) \ltimes \bW^*$ and $\opp{\fq} \defeq \bW \rtimes (\fg\dsum\bR)$ are opposite abelian parabolics.  The corresponding Satake diagrams are
\begin{equation*}
  \fq = \dynkinCp{2}{0}{3}{proj-bgg-q}
   \leq \dynkinC{2}{0}{3}{proj-bgg-h} = \alg{sp}{2n+2,\bR},
\end{equation*}
so that $H\acts\fq$ is the \grassmannian\ of lagrangian subspaces of $\bR^{2n+2}$.

Conversely, choose an algebraic Weyl structure for $\fh$.  The infinitesimal isotropy representation $\fh/\fq$ of $\fh$ descends to a representation $\bW$ of the semisimple part $\fg \defeq \liebrac{\fq^0}{\fq^0} \isom \alg{sl}{n+1,\bR}$ of the Levi factor $\fq^0 \defeq \fq/\fq^{\perp} \isom \alg{gl}{n+1,\bR}$ of $\fq$.  Since $\bW$ is induced by the adjoint representation of $\fh=\alg{sp}{2n+2,\bR}$, it has highest weight
\begin{equation*}
  \bW = \dynkinSLR[2,0,0,0]{2}{0}{2}{proj-bgg-W}
\end{equation*}
as a $\fg$-representation.  The dual representation $\bW^*$ has highest weight supported on the final node, so determines a symmetric \Rspace\ $G\acts\fp$ by crossing this node.  Of course, $G\acts\fp$ is just $\RP[n]$.  We shall see later that similar isomorphisms exist in \cproj\ and \qtn ic\ geometries, and are key for generalising the classical \proj\ structures.

\chapter{\Cproj\ geometry} 
\label{c:cproj}

\renewcommand{\algbracadornment}{c}
\BufferDynkinLocaltrue
\renewcommand{\dynkinnameoffset}{-0.35}

{\=O}tsuki and Tashiro \cite{ot1954-hgeodesics} observed that two \kahler\ metrics are \proj ly\ equivalent if and only if they are affinely equivalent, implying that the na{\"i}ve application of \proj\ differential geometry to complex manifolds is essentially uninteresting.  \Cproj\ geometry arises as a natural generalisation, where we replace geodesics with a complex analogue.  We begin by reviewing the rudiments of almost complex geometry in Section \ref{s:cproj-cpx}, before describing the classical theory of \cproj\ geometry in Section \ref{s:cproj-class}.  This develops in two parallel threads: a similar formulation to the theory of Section \ref{s:proj-class} using \kahler\ metrics, or via so-called \emph{hamiltonian $2$-forms}.

We describe \cproj\ geometry as an abelian parabolic geometry in Section \ref{s:cproj-para}, with $\fg$ the real Lie algebra underlying $\alg{sl}{n+1,\bC}$ and $\fp$ the stabiliser of a given complex line.   Thus $G\acts\fp \isom \CP$, and the general theory gives an equivalence of categories between (almost) \cproj\ structures and normal Cartan geometries of type $\CP[n]$.

The metrisability of a \cproj\ structure proceeds in much the same way as in Section \ref{s:proj-bgg}, with compatible \kahler\ metrics given by solutions of the first BGG operator associated to $\bW \defeq \realrepn{(\bC^{n+1} \etens \conj{\bC^{n+1}})}$.  The first BGG operator associated to $\bW^*$, called the \emph{\cproj\ hessian}, is also important to the theory; we discuss these two BGG operators in Section \ref{s:cproj-bgg}.

\section{Background on almost complex geometry} 
\label{s:cproj-cpx}

We begin by reviewing the basic theory of almost complex manifolds, primarily to fix conventions and notation; comprehensive introductions may be found in \cite{s2001-symplectic, j2003-calibrated, kn1996-founddg2, v2002-hodge}.  An \emph{almost complex structure} on a manifold $M$ is an endomorphism $J\in\s{0}{\alg{gl}{TM}}$ satisfying $J^2=-\id$.  For each $x\in M$, the choice of such a $J$ equips the tangent space $T_xM$ with the structure of a complex vector space in which we identify multiplication by $\cpx{i}\defeq \sqrt{-1}\in \bC$ with the application of $J\at{x}\in\alg{gl}{T_xM}$.  In particular, it follows that $\dim M = 2n$ is even.  The pair $(M,J)$ is called an \emph{almost complex manifold}; equivalently, an almost complex structure is a reduction of the frame bundle of $M$ to structure group $\grp{GL}{n,\bC} \leq \grp{GL}{2n,\bR}$.

By the Newlander--Nirenberg theorem \cite{nn1957-cpxcoords}, an almost complex manifold admits holomorphic coordinates if and only if the \emph{Nijenhuis torsion}
\begin{equation} \label{eq:cproj-cpx-nijenhuis} \begin{split}
  \Nijen{X,Y} &\defeq (J\cL_X J-\cL_{JX}J)Y \\
    &\phantom{:}= [X,Y] + J[JX,Y] + J[X,JY] - [JX,JY]
\end{split}
\end{equation}
of $J$ vanishes.  In this case, we drop the prefix ``almost'' in ``almost complex structure'' and ``almost complex manifold'', and say that $J$ is \emph{integrable}.  Equivalently, a complex structure is a \torsionfree\ $\grp{GL}{n,\bC}$-structure.

Since $J^2=-\id$, it has eigenvalues $\pm\cpx{i}$ on each tangent space.  Thus $J$ is diagonalisable on the complexification $\cpxbdl{TM} \defeq TM\tens[\bR]\bC$ of the tangent bundle; if $T^{1,0}M$ and $T^{0,1}M$ are respectively the $(\pm\cpx{i})$-eigensubbundles of $J$, we have
\begin{equation*}
  \bC TM = T^{1,0}M\dsum T^{0,1}M = \setof{X}{JX =  \cpx{i}X}
                              \dsum \setof{X}{JX = -\cpx{i}X}
\end{equation*}
and clearly $\conj{T^{1,0}M}=T^{0,1}M$.  There is a corresponding decomposition of the complexified cotangent bundle: we define $\Wedge{1,0}M \defeq (T^{1,0}M)^*$, which identifies $\Wedge{1,0}M$ with the annihilator of $T^{0,1}M$ in $\cpxbdl{TM}$, and similarly for $\Wedge{0,1}M$.  It follows that
\begin{equation} \label{eq:cproj-cpx-CT*M}
  \bC T^*M = \Wedge{1,0}M\dsum \Wedge{0,1}M
         = \setof{\omega}{J\omega = -\cpx{i}\omega}
     \dsum \setof{\omega}{J\omega =  \cpx{i}\omega},
\end{equation}
where $J\omega \defeq -\omega\circ J$.  Note that some authors (for example \cite{j2003-calibrated}) define $\Wedge{1,0}M$ to be the $(+\cpx{i})$-eigensubbundle of $J$, in which case $\Wedge{1,0}M=(T^{0,1}M)^*$ instead.
Given arbitrary $X\in\s{0}{\cpxbdl{TM}}$ and $\omega\in\s{0}{\cpxbdl{T^*M}}$, the $(1,0)$- and $(0,1)$-parts are given by
\begin{equation*} \begin{aligned}
  X^{1,0} &\defeq \tfrac{1}{2}\left( X - \cpx{i}JX \right), &\quad
  X^{0,1} &\defeq \tfrac{1}{2}\left( X + \cpx{i}JX \right), \\
  \text{and}\quad
  \omega^{1,0} &\defeq \tfrac{1}{2}\left( \omega + \cpx{i}J\omega \right), &\quad
  \omega^{0,1} &\defeq \tfrac{1}{2}\left( \omega - \cpx{i}J\omega \right).
\end{aligned} \end{equation*}
If $X\in\s{0}{T^{1,0}M}$ then taking $X_{\bR} \defeq \tfrac{1}{2}(X+\widebar{X})$ gives $X = \tfrac{1}{2}(X_{\bR} - \cpx{i}JX_{\bR})$.  It follows that the map $TM \to T^{1,0}M$ given by $X\mapsto \tfrac{1}{2}(X-\cpx{i}JX)$ is a vector bundle isomorphism; similar statements hold for $(0,1)$-vectors and complex forms.

The decomposition \eqref{eq:cproj-cpx-CT*M} induces a decomposition
\vspace{-0.2em}
\begin{equation} \label{eq:cproj-cpx-pqforms}
  \Wedge{k}\bC T^*M
    = \Wedge{k} \left( \Wedge{1,0}M\dsum \Wedge{0,1}M \right)
    = \Dsum{p+q=k}{} \left( \Wedge{p}\Wedge{1,0}M
                      \tens \Wedge{q}\Wedge{0,1}M \right)
\vspace{-0.2em}
\end{equation}
of $\Wedge{k} \cpxbdl{TM}$ for each $k \in \bN$, where the $(p,q)$th summand $\Wedge{p,q}M := \Wedge{p}\Wedge{1,0}M \dsum \Wedge{q}\Wedge{0,1}M$ in \eqref{eq:cproj-cpx-pqforms} is called the space of \emph{$(p,q)$-forms} on $M$; clearly $\omega \in\s{0}{\Wedge{p,q}M}$ has $p+q$ covariant indices, $p$ of type $(1,0)$ and $q$ of type $(0,1)$.  The decomposition \eqref{eq:cproj-cpx-pqforms} may be described in representation-theoretic terms: if $F^0$ is the reduction of the frame bundle of $M$ to structure group $P^0\defeq \grp{GL}{n,\bC}$, then $\Wedge[\bC]{k}\cpxbdl{T^*M}$ is the bundle associated to $F^0$ and $\Wedge{k} \bR^{2n*}\tens\bC$.  The $\grp{GL}{n,\bC}$-representation $\Wedge{k}\bR^{2n*}\tens\bC$ decomposes as
\vspace{-0.2em}
\begin{equation*}
  \Wedge{k}\bR^{2n*}\tens\bC
    = \Dsum{p+q=k}{}( \Wedge[\bC]{p}\bC^{n*} \tens[\bC] \Wedge[\bC]{q}\conj{\bC^{n*}} ),
\vspace{-0.2em}
\end{equation*}
where overline denotes the complex conjugate representation.  Of course, $\Wedge{p,q}M$ is the bundle associated to $\Wedge[\bC]{p}\bC^{n*} \tens \Wedge[\bC]{q}\conj{\bC^{n*}}$.  We shall often need the decomposition
\vspace{-0.2em}
\begin{equation} \label{eq:cproj-cpx-2forms}
  \Wedge[\bC]{2}\cpxbdl{T^*M}
    = ( \Wedge{2,0}M\dsum\Wedge{0,2}M ) \dsum \Wedge{1,1}M
\vspace{-0.2em}
\end{equation}
of complex $2$-forms, where the two summands are respectively the complexifications of the bundles $\Wedge[-]{2}T^*M$ and $\Wedge[+]{2}T^*M$ of $J$-anti-invariant and $J$-invariant $2$-forms.  Decompositions analogous to \eqref{eq:cproj-cpx-pqforms} and \eqref{eq:cproj-cpx-2forms} exist for $k$-vectors and symmetric $k$-tensors, for which we write $\Wedge[\bC]{k} \cpxbdl{TM} = {\textstyle \Dsum{p+q=k}{}} T^{p,q}M$ and $\Symm[\bC]{k} \cpxbdl{T^*M} = {\textstyle \Dsum{p+q=k}{}} \Symm{p,q}T^*M$.

The natural connections on an almost complex manifold $(M,J)$ are those which preserve $J$, \ie\ $\D J=0$.  The space of such connections is affine, modelled on $1$-forms with values in the complex-linear endomorphisms $\alg{gl}{TM,J}$ of $TM$.  The (complexification of the) torsion $\Tor[\D]{}\in\s{2}{TM}$ of $\D$ splits into $(2,0)$-, $(1,1)$- and $(0,2)$-parts as in \eqref{eq:cproj-cpx-2forms}, and it is easy to show that the $(0,2)$-part is proportional to the Nijenhuis torsion \eqref{eq:cproj-cpx-nijenhuis} of $J$.  Moreover $\D$ may be deformed in such a way as to remove the $(2,0)$- and $(1,1)$-parts \cite{cemn2015-cproj}, so that $(M,J)$ admits a \torsionfree\ complex connection if and only if $J$ is integrable, if and only if its Nijenhuis torsion \eqref{eq:cproj-cpx-nijenhuis} vanishes.

Finally, the class of (\pseudo)\riem\ metrics $g$ of interest are the \emph{hermitian} metrics, \ie\ those satisfying $g(JX,JY) = g(X,Y)$.  If $\b{g}$ is an arbitrary metric, $g(X,Y) \defeq \tfrac{1}{2}( \b{g}(X,Y) + \b{g}(JX,JY) )$ is hermitian, so that almost complex manifolds always admit hermitian metrics.  The complexification of a hermitian metric $g$ defines a hermitian inner product on each $\cpxbdl{T_xM}$, thus giving a section of $\Symm{1,1} \cpxbdl{T^*M}$; equivalently $g$ is a section of the $J$-invariant subbundle $\Symm[+]{2}T^*M$.  It follows that $\omega\defeq g(J\bdot,\bdot)$ is a $J$-invariant $2$-form, called the \emph{hermitian form} of $g$.  Evidently $g$ and $\omega$ may be viewed respectively as the real and imaginary parts of a hermitian metric on $\cpxbdl{T^*M}$.  If $\omega$ is closed, we call $g$ a \emph{(\pseudo)\kahler\ metric} and $\omega$ its \emph{\kahler\ form}.

Suppose that $J$ is integrable, and that $\D$ is the \LC\ connection of a hermitian metric $g$.  Then it is straightforward to check that $\d\omega=0$ if and only if $\D\omega=0$.  Since $\omega=g(J\bdot,\bdot)$, it follows that $g$ is \kahler\ if and only if $\D J=0$, \ie\ if and only if $\D$ is a complex connection.%
\footnote{Evidently any two of the conditions $\D g=0$, $\D J=0$ and $\D\omega=0$ imply the third.}
The condition $\D J=0$ and the exchange identity imply that the \riem\ curvature tensor $\Curv{}{}$ is a $J$-invariant $\alg{gl}{TM,J}$-valued $2$-form, which also implies that the Ricci curvature $\sRic[g]{}$ of $g$ is symmetric and $J$-invariant.  Moreover since $\D$ preserves more than just a \riem\ metric, $\Curv{}{}$ admits a finer decomposition than a generic \riem\ curvature tensor; we shall not need this decomposition and refer the interested reader to \cite[\S 2D]{b1987-einstein} and \cite[Eqn.\ (3)]{acg2006-ham2forms1}.

\medskip
\section{Classical definition and results} 
\label{s:cproj-class}

\Cproj\ geometry%
\footnote{The classical literature uses the term ``H-\proj'' rather than ``\cproj''.  The latter was recently adopted to reflect the fact that a generic \cproj\ class does \emph{not} contain a holomorphic connection; see \cite[Rmk.\ 1]{mr2013-conification}.}
was introduced by {\=O}tsuki and Tashiro in \cite{ot1954-hgeodesics} as a generalisation of \proj\ differential geometry to the holomorphic category.  They observed that two \kahler\ metrics on a complex manifold $(M,J)$ that are \proj ly\ equivalent in the sense of \thref{defn:proj-class-equiv} are affinely equivalent, meaning that \proj\ equivalence is uninteresting for \kahler\ manifolds.  Many results analogous to those in \proj\ differential geometry have natural generalisations \cite{m1998-holomorphic}, as well as generalisations to almost complex structures \cite{t1957-hproj}.  We study this formulation in Subsection \ref{ss:cproj-class-maineqn}.

Later, \cproj\ differential geometry was unknowingly rediscovered under the guise of \emph{\hamiltonian\ $2$-forms} by Apostolov \etal\ while studying \kahler\ metrics with special curvature properties \cite{acg2003-ham2forms0, acg2006-ham2forms1, acgt2004-ham2forms2, acgt2008-ham2forms4, acgt2008-ham2forms3}.  The authors observed that the existence of a \hamiltonian\ $2$-form imbued the \kahler\ metric with an isometric torus action, which allowed them to give both local and global classifications of \kahler\ metrics admitting \hamiltonian\ $2$-forms.  An overview of \hamiltonian\ $2$-forms is given in Subsection \ref{ss:cproj-class-ham2forms}.

For greater generality we allow \non degenerate\ metrics of arbitrary signature, reserving the term \emph{\riem\ \kahler} for positive definite \kahler\ metrics.

\medskip
\subsection{Classical approach} 
\label{ss:cproj-class-maineqn}

Let $(M,g,J)$ be a $2n$-dimensional \kahler\ manifold with \LC\ connection $\D$.

\begin{defn} \thlabel{defn:cproj-class-cgeodesic} A smooth curve $\gamma\subset M$ is a \emph{c-geodesic} of $g$ if $\D_X X \in\linspan{X,JX}{}$ for every vector field $X$ tangent to $\gamma$.  If $\b{g}$ is a second \kahler\ metric, $g$ and $\b{g}$ are called \emph{\cproj ly\ equivalent} if they have the same c-geodesics (as unparameterised curves). \end{defn}

Evidently any two metrics on a Riemann surface (\ie\ $n=1$) are \cproj ly\ equivalent, so that a different approach is required to avoid triviality.  We shall discuss this further in Subsection \ref{ss:cproj-bgg-hess}; for now we assume that $n>1$.  The following characterisation of \cproj ly\ equivalent \kahler\ metrics was obtained by {\=O}tsuki and Tashiro \cite{ot1954-hgeodesics} and Ishihara and Tachibana \cite{it1960-infholo}; its proof follows a similar argument to that of \thref{lem:proj-class-algbrac}.

\begin{lem} \thlabel{lem:cproj-class-algbrac} \Kahler\ metrics $g,\b{g}$ are \cproj ly\ equivalent if and only if their \LC\ connections $\D,\b{\D}$ are related by
\begin{equation} \label{eq:cproj-class-algbrac} \begin{aligned}
  \b{\D}_X Y &\phantom{:}= \D_X Y + \algbrac{X}{\alpha} \acts Y \\
  \text{where}\quad
  \algbrac{X}{\alpha} \acts Y &\defeq \tfrac{1}{2}(\alpha(X)Y + \alpha(Y)X - \alpha(JX)JY - \alpha(JY)JX)
\end{aligned}
\end{equation}
for some $\alpha\in\s{1}{}$ and all $X,Y\in\s{0}{TM}$. \noproof \end{lem}

\begin{rmk} The endomorphism $\algbrac{X}{\alpha}$ will be called the \emph{algebraic bracket} of $X$ and $\alpha$, which is easily seen to be complex-linear and symmetric in $X,Y$. \end{rmk}

As in \proj\ differential geometry, taking a trace in \eqref{eq:cproj-class-algbrac} yields
\vspace{0.05em}
\begin{equation} \label{eq:cproj-class-exact1form}
  \alpha = \tfrac{1}{2n+2} \,\d\! \left( \log \frac{\det\b{g}}{\det g} \right)
  \colvectpunct[-0.6em]{,}
\vspace{0.05em}
\end{equation}
so that $\alpha$ is an exact $1$-form; \cf\ equation \eqref{eq:proj-class-exact1form}.  Given \cproj ly\ equivalent \kahler\ metrics $g,\b{g}$, we follow the lead of \eqref{eq:proj-class-endoA} by considering the endomorphism
\vspace{0.05em}
\begin{equation} \label{eq:cproj-class-endoA}
  A(g,\b{g}) \defeq \left( \frac{\det\b{g}}{\det g} \right)^{\!\!1/(2n+2)} \b{\sharp}\circ\flat.
\vspace{0.05em}
\end{equation}
It is straightforward to see that $A(g,\b{g})$ is invertible with inverse $A(\b{g},g)$ and is \self adjoint\ \wrt\ both $g$ and $\b{g}$; moreover since $g,\b{g}$ are hermitian, $A \defeq A(g,\b{g})$ is complex-linear.  The second \kahler\ metric $\b{g}$ may be recovered from the pair $(g,A)$ as
\vspace{0.05em}
\begin{equation} \label{eq:cproj-class-gbar}
  \b{g} = (\det A)^{-1/2}\, g(A^{-1}\bdot,\bdot);
\vspace{0.05em}
\end{equation}
\cf\ equation \eqref{eq:proj-class-gbar}.  Calculating as in \thref{prop:proj-class-maineqn}, Domashev and Mike{\v s} \cite{dm1978-hproj} found that $A$ satisfies a first-order linear differential equation.

\smallskip
\begin{prop} \thlabel{prop:cproj-class-maineqn} Let $g,\b{g}$ be \kahler\ metrics with \LC\ connections $\D,\b{\D}$ respectively.  Then $g,\b{g}$ are \cproj ly\ equivalent if and only if $A \defeq A(g,\b{g})$ defined by \eqref{eq:cproj-class-endoA} satisfies the first-order linear differential equation
\begin{equation} \label{eq:cproj-class-maineqn}
  g((\D_X A)\bdot,\bdot) = X^{\flat}\symm\mu + JX^{\flat}\symm J\mu
\end{equation}
for some $\mu \in\s{1}{}$ and all $X\in\s{0}{TM}$.  In this case $\b{\D} = \D + \algbrac{}{\alpha}$, where $\alpha\in\s{1}{}$ satisfies $\mu = -\alpha(A\bdot) = \tfrac{1}{2}\d(\tr A)$. \noproof \end{prop}

\begin{rmk} \thlabel{rmk:cproj-class-maineqn} Although the lineage of \eqref{eq:cproj-class-maineqn} is less contentious than that of \eqref{eq:proj-class-maineqn}, we refer to it as the \emph{main equation} of \cproj\ geometry to develop the analogy between the two theories.  Domashev and Mike{\v s} obtained a prolongation of the main equation \cite[Thm.\ 2]{dm1978-hproj}; also see the survey \cite{m1998-holomorphic}.
\end{rmk}

Thus a solution $(g,A)$ of \eqref{eq:cproj-class-maineqn} is equivalent to a pair of \cproj ly\ equivalent \kahler\ metrics $g, \b{g}\defeq (\det A)^{-1/2} g(A^{-1}\bdot,\bdot)$.  Linearity then yields a \emph{metrisability pencil} of \cproj ly\ equivalent \kahler\ metrics induced by $A_t \defeq A - t\id$ and \eqref{eq:cproj-class-gbar}.

Finally, let us indicate the link between the classical theory above and the description as a parabolic geometry outlined in Section \ref{s:cproj-para}.  \thref{lem:cproj-class-algbrac} shows that c-geodesics are really a feature of linear connections, so let us call $\D,\b{\D}$ \cproj ly\ equivalent if $\b{\D}_X = \D_X + \algbrac{X}{\alpha}$ for some $\alpha \in \s{1}{}$, where $\algbrac{X}{\alpha}$ is the algebraic bracket from \eqref{eq:cproj-class-algbrac}.  Here we make no assumption on the torsions of $\D,\b{\D}$.  A \cproj\ structure on $M$ is then an equivalence class $\Dspace$ of \cproj ly\ equivalent connections, with $\alpha \mapsto \algbrac{}{\alpha}$ defining an embedding of $\s{1}{} \injto \s{1}{\alg{gl}{TM,J}}$.  Thus $\Dspace$ is an affine space modelled on $\s{1}{}$, leading to the following analogue of \thref{lem:proj-class-connsonL}.

\begin{lem} \thlabel{lem:cproj-class-connsonL} There is a bijection between connections $\D \in \Dspace$ and connections on the line bundle $\cL \defeq (\Wedge{2n}TM)^{1/(n+1)}$, where $\D,\b{\D}$ are related by $\alpha\in\s{1}{}$ if and only if $\alpha$ is the change of induced connection on $\cL$. \noproof \end{lem}

A quantity is \emph{\cproj ly\ invariant} if it is independent of $\D\in \Dspace$.  By the symmetry of $\algbrac{X}{\alpha}\acts Y$ in $X,Y$, the relation $\b{\D}_X Y = \D_X Y + \algbrac{X}{\alpha}\acts Y$ implies that $\D,\b{\D}$ have the same torsion, which is then a \cproj\ invariant of $\Dspace$.

\subsection{\Hamiltonian\ $2$-forms} 
\label{ss:cproj-class-ham2forms}

An equivalent approach to \cproj\ geometry was developed independently by Apostolov \etal\ \cite{acg2003-ham2forms0, acg2006-ham2forms1} while studying \kahler\ metrics with special curvature properties.

\begin{defn} \thlabel{defn:cproj-class-ham2form} Let $(M,g,J)$ be a $2n$-dimensional \kahler\ manifold.  A \emph{\hamiltonian\ $2$-form} on $M$ is a $J$-invariant $2$-form $\phi$ satisfying
\begin{equation} \label{eq:cproj-class-ham2form}
  \D_X\phi = \tfrac{1}{2}( \d\sigma\wedge JX^{\flat} - J\d\sigma\wedge X^{\flat} )
\end{equation}
for some smooth function $\sigma\in\s{0}{}$ and all $X\in\s{0}{TM}$.%
\footnote{Equation \eqref{eq:cproj-class-ham2form} is trivial when $n=1$, so we additionally assume that $\sigma$ is a Killing potential, \ie\ that $J\grad[g]{\sigma}$ is a Killing vector field.}
\end{defn}

It is straightforward to see that the main equation \eqref{eq:cproj-class-maineqn} and the \hamiltonian\ $2$-form equation \eqref{eq:cproj-class-ham2form} are equivalent.  Indeed, an endomorphism $A$ satisfies the main equation if and only if $\phi\defeq g(JA\bdot,\bdot)$ is a \hamiltonian\ $2$-form.  Taking a trace in \eqref{eq:cproj-class-ham2form} \wrt\ the \kahler\ form $\omega$ yields $\sigma = \tr[\omega]\phi$, so that $\d\sigma$ coincides with $\mu = \tfrac{1}{2} \d(\tr A)$ from \thref{prop:cproj-class-maineqn}.  The \cproj ly\ equivalent metric determined by $\phi$ is given by \eqref{eq:cproj-class-gbar}.  A more detailed account of the relationship between may be found in \cite{cmr2015-cprojmob}.

Since the \hamiltonian\ $2$-form equation is over-determined, it may be prolonged to a closed system.  For this the authors of \cite{acg2006-ham2forms1} inverted \eqref{eq:cproj-class-ham2form} using the metric $g$, thus obtaining a closed differential system on $\Wedge[+]{2}TM \dsum TM \dsum (M\times \bR)$.  This system is equivalent to the prolongation given by Domashev and Mike{\v s} discussed in \thref{rmk:cproj-class-maineqn}.

Apostolov \etal\ were able to obtain a local classification of \hamiltonian\ $2$-forms \cite[Thm.\ 1]{acg2006-ham2forms1} by utilising an intrinsic torus action afforded by a pair of \cproj ly\ equivalent metrics.  Firstly, if $\phi$ is a \hamiltonian\ $2$-form then so is $\phi_t \defeq \phi-t\omega$ for all $t\in\bR$; \cf\ the paragraph following \thref{rmk:cproj-class-maineqn}.  The \emph{pfaffian} $\pf{\phi_t} \defeq \hodge(\tfrac{1}{n!}\phi_t^{\wedge n}) \in\s{0}{}$ of $\phi_t$, induced by the Hodge star $\hodge$ of $g$, equips $M$ with a $1$-parameter family of vector fields $K_t \defeq J\grad[g]{(\pf{\phi_t})}$.  It turns out that the $K_t$ are commuting \hamiltonian\ Killing vector fields \cite[Prop.\ 3]{acg2006-ham2forms1}, and the rank of the family $t\mapsto K_t$ is constant on a dense open subset of $M$.  The $K_t$ therefore generate an isometric torus action, from which a local classification may be obtained by specialising Pedersen and Poon's description \cite{pp1991-pedpoon} of \kahler\ manifolds admitting such a torus action.  This classification may be viewed as a complex version of \LC's classification of \riem\ metrics admitting a \proj\ equivalent metric \cite{l2009-dyneqns}.  A global classification was obtained for compact manifolds in \cite{acgt2004-ham2forms2}, by appealing to Guillemin's construction \cite{cdg2002-guillemin, g1994-torickahler, g1994-hamiltonianTn} of toric \kahler\ manifolds via the Delzant polytope \cite{f1993-toric} of the acting torus.  Applications of \hamiltonian\ $2$-forms in \kahler\ geometry were developed in \cite{acgt2008-ham2forms4, acgt2008-ham2forms3}.  A local classification of \pseudo\kahler\ metrics admitting a \cproj ly\ equivalent metric was recently obtained by Bolsinov \etal\ \cite{bmr2015-localcproj} by a combination of the aforementioned toric methods and appealing to the local classification of \proj ly\ equivalent metrics \cite{bm2011-splittingglueing, bm2013-priemlocalform, bmp2009-priemlocalform2d}.

\vspace{-0.4em}
\section{Description as a parabolic geometry} 
\label{s:cproj-para}
\vspace{-0.2em}

Complex \proj\ space $\CP$ may be defined as either the set of complex lines through the origin in $\bC^{n+1}$, or as the quotient of the $(2n+1)$-sphere $\Sph[2n+1]$ by its $\Sph[1]$-action; the latter is the famous Hopf fibration $\Sph[2n+1] \surjto \CP[n]$ \cite{b1987-einstein}.  The standard metric on $\bC^{n+1}$ of a given signature induces a \kahler\ metric $g_{\mr{FS}}$ of constant holomorphic sectional curvature on $\CP$, called the \emph{Fubini--Study metric} of that signature \cite{kn1996-founddg2}.  Much like in the case of $\RP$, the embedded complex \proj\ lines $\CP[1] \injto \CP$ are totally geodesic with respect to the \LC\ connection $\D^{\mr{FS}}$ of $g_{\text{FS}}$.  It follows that the c-geodesics of $[\D^{g_{\text{FS}}}]$ are the smooth curves contained within such a complex line, so that \wrt\ an affine chart $\bC^n \injto \CP$ the c-geodesics lie within a complex line in $\bC^n$.  Thus, from the point of view of c-geodesics, $\CP$ is the natural compactification of $\bC^n$ and therefore is suitable as the flat model of \cproj\ geometry.

We shall assume henceforth that $n>1$, so that we consider almost complex manifolds $(M,J)$ of dimension $2n \geq 4$.  Evidently any two connections on a riemann surface are \cproj ly\ equivalent, so that the case $n=1$ requires a different approach to be \non trivial; we shall discuss this in \thref{rmk:cproj-bgg-dim2}.

\subsection{The flat model $\CP$} 
\label{ss:cproj-para-flat}

As a generalised flag manifold $\CP$ is isomorphic to $G/P$ for the real Lie groups
\vspace{-0.3em}
\begin{equation*}
  G \defeq \grp{PGL}{n+1,\bC} \defeq \grp{GL}{n+1,\bC} / \linspan[\bC]{\id}{}
\vspace{-0.3em}
\end{equation*}
and $P\leq G$ the (projection to $G$ of) the $\grp{GL}{n+1,\bC}$-stabiliser of a chosen complex line $\linspan[\bC]{v_0}{}$ in $\bC^{n+1}$.  It is important to note that we view both $G$ and $P$ as \emph{real} Lie groups.  The Lie algebra of $G$ is the underlying real Lie algebra $\fg\defeq \alg{sl}{n+1,\bC}$, with complexification $\cpxrepn{\fg} \defeq \alg{sl}{n+1,\bC} \dsum \alg{sl}{n+1,\bC}$, while the parabolic subalgebra $\fp$ corresponds to crossing the last node of each factor:
\vspace{0.25em}
\begin{equation} \label{eq:cproj-para-satake}
\smash{
  \fp = \dynkinSLCp{2}{0}{2}{cproj-para-flat-p} \leq
        \dynkinSLC{2}{0}{2}{cproj-para-flat-g} = \fg.
}
\vspace{0.25em}
\end{equation}
Therefore \cproj\ geometry is an abelian parabolic geometry, with Killing polar $\fp^{\perp}\isom \bC^{n*} \defeq \Hom[\bR]{\bC^n}{\bR}$ and reductive Levi factor $\fp^0 \defeq \fp/\fp^{\perp}$.  An algebraic Weyl structure splits the $\fp^{\perp}$-filtration $\fg\supset \fp\supset \fp^{\perp}\supset 0$ of $\fg$, determining an isomorphism $\fg \isom \fg/\fp \dsum \fp^0 \dsum \fp^{\perp}$.  Such an isomorphism is equivalent to a complex subspace of $\bC^{n+1}$ complementary to $\linspan[\bC]{v_0}{}$, thus yielding a decomposition
\begin{equation} \label{eq:cproj-para-sldecomp}
  \alg{sl}{n+1,\bC} = \Setof{
  \begin{pmatrix}
    -\tr A & \alpha \\
     ~X    & A
  \end{pmatrix}}{\, X\in\bC^n,~ \alpha\in\bC^{n*},~ A\in\alg{gl}{n,\bC} }
\end{equation}
akin to \eqref{eq:proj-para-sldecomp}.  In the corresponding decomposition of $\grp{PGL}{n+1,\bC}$, the Levi subgroup $P^0 \defeq P/\exp{\fp^{\perp}}$ consists of (equivalence classes of) block-diagonal matrices in $G$; clearly the adjoint action of $P^0$ on $\fg/\fp\isom\bC^n$ induces an isomorphism $\fp^0 \isom \alg{gl}{n,\bC}$.   Then $\fp^{\perp}\isom\bC^{n*}$ consists of matrices with only the $\alpha$-block, $\fp^0\isom\alg{gl}{n,\bC}$ of block-diagonal matrices, and $\fp \isom \alg{gl}{n,\bC}\ltimes\bC^{n*}$ of block upper-triangular matrices.

Choose an algebraic Weyl structure, so that $\fg \isom \bC^n\dsum\alg{gl}{n,\bC}\dsum\bC^{n*}$.  Using the decomposition \eqref{eq:cproj-para-sldecomp}, we can write down the Lie brackets between the three summands: if $X,Y\in\bC^n$, $A,B\in\alg{gl}{n,\bC}$ and $\alpha,\beta\in\bC^{n*}$ then
\begin{equation*}
  \liebrac{X}{Y} = 0 = \liebrac{\alpha}{\beta}, \quad
  \liebrac{A}{B} = AB - BA, \quad
  \liebrac{A}{X} = AX, \quad
  \liebrac{A}{\alpha} = -\alpha\circ A,
\end{equation*}
so that $\bC^n,\bC^{n*}$ form abelian subalgebras on which $\alg{gl}{n,\bC}$ acts naturally.  The bracket $\bC^n \times \bC^{n*} \to \alg{gl}{n,\bC}$ requires a little more work to compute: if $J$ is the complex structure used to identify $\bC^n\isom \bR^{2n}$, Hrdina \cite{h2009-cproj} calculates that
\vspace{-0.1em}
\begin{equation*}
  \liebrac{ \liebrac{X}{\alpha} }{ Y } = \alpha(X)Y + \alpha(Y)X - \alpha(JX)JY - \alpha(JY)JX.
\vspace{-0.1em}
\end{equation*}
Up to normalisation conventions, we see in particular that change of algebraic Weyl structure corresponds to change of connection in the \cproj\ class via \eqref{eq:cproj-class-algbrac}.

Observe that $\CP$ has a \proj\ embedding defined as follows.  The complex representation $\bC^{n+1}$ of $\cpxrepn{\fg} = \alg{sl}{n+1,\bC} \dsum \alg{sl}{n+1,\bC}$ has a complex structure given by multiplication by $\cpx{i}$, which we view as a linear map on the underlying real vector space.  There is a conjugate-linear isomorphism $\sigma : \bC^{n+1} \to \conj{\bC^{n+1}}$ given by complex conjugation, so that $\sigma^2=\id$ and $\sigma\cpx{i}\sigma^{-1}$ is the complex structure on $\conj{\bC^{n+1}}$.  Noting that representations of $\cpxrepn{\fg}$ are external tensor products of representations of each factor,%
\footnote{We will relate $\fg$- and $\cpxrepn{\fg}$-representations more carefully in the next subsection.}
it is clear that $\cpx{i} \etens \sigma\cpx{i}$ is a real structure on $\bC^{n+1}\etens\conj{\bC^{n+1}}$; the fixed-point set of $\sigma \etens \sigma$ is the underlying real representation $\bW \defeq \realrepn{( \bC^{n+1}\etens\conj{\bC^{n+1}} )}$ of $\fg=\alg{sl}{n+1,\bC}$.  Identifying $[v]\in\CP$ with the real $2$-plane $\linspan[\bC]{v,\cpx{i}v}{}$ and $\conj{\bC^{n+1}}$ with $\sigma(\bC^{n+1})$, the map $[v] \mapsto [v\etens\sigma v]$ defines an embedding $\CP \injto \pr{\bC^{n+1} \etens \conj{\bC^{n+1}}} = \pr{\cpxrepn{\bW}}$.  Since $(\sigma\etens\sigma)(v\etens\sigma v) = \sigma v \etens \sigma^2 v = v \etens \sigma v$, this map takes values in the underlying real representation $\bW$, thus yielding a \proj\ embedding $\CP \injto \pr{\bW}$.

Evidently this embedding is not minimal, so that Kostant's \thref{thm:lie-para-kostant} tells us that $\CP$ is the intersection of quadrics in $\pr{\bW}$ given by projection away from the Cartan square in $\Symm{2}\bW^*$.  In notation explained properly in the next subsection,
\vspace{-0.05em}
\begin{equation*}
  \Symm{2}\bW^*
    = \Symm{2}\hspace{-0.6ex}
      \left[ \dynkinSLC[0,0,0,1, 0,0,0,1]{2}{0}{2}{cproj-para-flat-Wd} \right]
    = \dynkinSLC[0,0,0,2, 0,0,0,2]{2}{0}{2}{cproj-para-flat-O2Wd}
    \dsum
      \dynkinSLC[0,0,0,1,0, 0,0,0,1,0]{2}{0}{3}{cproj-para-flat-Ud}
    \colvectpunct[-0.7em]{\, ,}
\vspace{-0.05em}
\end{equation*}
where the first summand is the Cartan square $\Cartan{2}\bW^* = \realrepn{( \Symm[\bC]{2}\bC^{n+1*} \etens \Symm[\bC]{2}\conj{\bC^{n+1*}} )}$ and
\vspace{-0.2em}
\begin{equation*}
  \bU^* \defeq \dynkinSLC[0,0,0,1,0, 0,0,0,1,0]{2}{0}{3}{cproj-para-flat-U}
    = \realrepn{( \Wedge[\bC]{2}\bC^{n+1} \etens \Wedge[\bC]{2}\conj{\bC^{n+1}} )}
\vspace{-0.2em}
\end{equation*}
is the space of homogeneous quadratic polynomials which cuts out $\CP$ as an intersection of quadrics.  As in Subsection \ref{ss:proj-para-flat}, $\bU$ is an irreducible $\fg$-representation.

Finally, since $\fg=\alg{sl}{n+1,\bC}$ is again one of the problematic Lie algebras for the equivalence of categories guaranteed by \thref{thm:para-calc-equiv}.  A prolongation can be constructed by methods similar to \thref{thm:proj-para-equiv}; see \cite[Thm.\ 2.8]{cemn2015-cproj}.

\begin{thm} \thlabel{thm:cproj-para-equiv} On any almost complex manifold $(M,J)$ of dimension $2n \geq 4$, there is an equivalence of categories between almost \cproj\ structures and regular normal parabolic geometries of type $\CP[n]$.  The flat model is $\CP$ with its canonical complex structure and \cproj\ structure determined by the Fubini--Study metric $g_{\mr{FS}}$. \noproof \end{thm}

Normality implies that the Weyl connections have torsion of type $(0,2)$ only \cite{cemn2015-cproj}.  The discussion following \thref{lem:cproj-class-connsonL} then implies that this torsion $\Tor{}$ is a \cproj ly\ invariant $TM$-valued $(0,2)$-form, being proportional to the Nijenhuis torsion $\Nijen{}$ of $J$.

\subsection{Representations of $\alg{sl}{n+1,\bC}$} 
\label{ss:cproj-para-repns}

The description of $\fg$- and $\fp$-representations proceeds more smoothly if we consider the complexification $\cpxrepn{\fg} \defeq \alg{sl}{n+1,\bC} \dsum \alg{sl}{n+1,\bC}$ of $\fg$.  Representations of $\cpxrepn{\fg}$ are external tensor products of representations of the two factors, with one factor notated on each line of the Satake diagram \eqref{eq:cproj-para-satake}.  However, some care must be taken.  Firstly, a representation of $\cpxrepn{\fg}$ is the complexification of a representation of $\fg$ if and only if it is \self conjugate.  Secondly, a complex representation of $\cpxrepn{\fg}$ may be reducible over $\bC$ but irreducible over $\bR$: for example, consider the $\cpxrepn{\fg}$-representations%
\footnote{Note that we denote representations of the complexifications $\cpxrepn{\fg}, \cpxrepn{\fp}$ on their respective Dynkin diagrams, rather than on the Satake diagrams \eqref{eq:cproj-para-satake} of $\fg,\fp$.}
\begin{equation*}
  \dynkinAA["$a$","$b$","$c$","$d$", "$a$","$b$","$c$","$d$"]{2}{0}{2}{cproj-para-repns-conj1}
  \quad\text{and}\quad
        \dynkinAA["$a$","$b$","$c$","$d$", 0,0,0,0]{2}{0}{2}{cproj-para-repns-conj2}
  \dsum \dynkinAA[0,0,0,0, "$a$","$b$","$c$","$d$"]{2}{0}{2}{cproj-para-repns-conj3}
  \colvectpunct[-0.7em]{\, .}
\end{equation*}
The first is irreducible over both $\bR$ and $\bC$, being the complexification of an irreducible $\fg$-representation%
\footnote{With highest weight given by ``putting the arrows in''.};
while the second is reducible over $\bC$ but irreducible over $\bR$, since the second summand is the complex conjugate of the first.

With these technicalities in mind let us describe some important representations of $\fg$ and $\fp$, as well as their associated bundles.  By the Cartan condition, the complexified isotropy representation $\cpxrepn{(\fg/\fp)}$ and its dual $\cpxrepn{(\fg/\fp)^*} \isom \cpxrepn{\fp^{\perp}}$ have associated bundles
\begin{align} \label{eq:cproj-para-TM}
  \cpxbdl{TM} &=
    \dynkinAAp[1,0,0,0,1, 0,0,0,0,0]{2}{0}{3}{cproj-para-repns-gp1} \hspace{0.25ex}\dsum
    \dynkinAAp[0,0,0,0,0, 1,0,0,0,1]{2}{0}{3}{cproj-para-repns-gp2}
    \hspace{0.3ex} = T^{1,0}M \dsum T^{0,1}M \\[0.8em]
  \text{and}\quad \label{eq:cproj-para-T*M}
  \cpxbdl{T^*M} &=
    \dynkinAAp[0,0,0,1,-2, 0,0,0,0,0]{2}{0}{3}{cproj-para-repns-pperp1} \dsum
    \dynkinAAp[0,0,0,0,0, 0,0,0,1,-2]{2}{0}{3}{cproj-para-repns-pperp2}
    = \Wedge{1,0}M \dsum \Wedge{0,1}M.
  %
\end{align}
The decomposition \eqref{eq:cproj-cpx-pqforms} of $\Wedge[\bC]{k} \cpxbdl{T^*M}$ into $(p,q)$-forms complex irreducible subrepresentations is given by taking exterior powers of \eqref{eq:cproj-para-T*M}; similar statements hold for the decompositions of $\Wedge[\bC]{k}\cpxbdl{TM}$, $\Symm[\bC]{k}\cpxbdl{T^*M}$ and $\Symm[\bC]{k}\cpxbdl{TM}$, with $\Wedge{p,p}M$ and $\Wedge{p,q}M \dsum \Wedge{q,p}M$ complexifications of real representations for all $p+q \in \setof{1,\ldots,2n}{}$.

The line bundle $\cL \defeq (\Wedge{2n}TM)^{1/(n+1)}$ from \thref{lem:cproj-class-connsonL} is associated to the $(n+1)$st root $L$ of $\Wedge{2n}(\fg/\fp)$.  Since the top exterior power of $\fg/\fp$ is the derivative of the real character $\det:\grp{GL}{n,\bC} \to \bR$, equation \eqref{eq:cproj-para-TM} implies that $L$ has highest weight
\begin{equation} \label{eq:cproj-para-repnL}
  L \defeq \dynkinSLCp[0,0,0,1, 0,0,0,1]{2}{0}{2}{cproj-para-repns-L}
    \colvectpunct[-0.6em]{\, .}
\end{equation}
Note that $L = \liehom{0}{\bW^*}$, where $\bW = \realrepn{( \bC^{n+1} \etens \conj{\bC^{n+1}} )}$ is the $\fg$-representation from Subsection \ref{ss:cproj-para-flat} admitting a \proj\ embedding $\CP \injto \pr{\bW}$.

As Subsection \ref{ss:proj-para-repns}, our choice of group $G = \grp{PGL}{n+1,\bC}$ with Lie algebra $\fg = \alg{sl}{n+1,\bC}$ means that not all $\fg$-representations integrate to global representations of $G$.  Indeed, a representation integrates if and only if the central circle in $\tilde{G} = \grp{SL}{n+1,\bC}$ acts trivially.  As before we can locally avoid this problem by forming the extended Cartan bundle $F^{\tilde{P}} \defeq \assocbdl{F^P}{P}{\tilde{P}}$, where $\tilde{P} \leq \tilde{G}$ is the parabolic stabiliser of a given complex line in $\bC^{n+1}$ (previously denoted $\linspan[\bC]{v_0}{}$).  Then all $\fg$-representations integrate to $\tilde{G}$-representations \wrt\ which we may form associated bundles, before taking the (local) quotient to $F^P$.  We obtain a $G$-integrability criterion by complexifying: since $\cpxrepn{G} \cong \grp{PGL}{n+1,\bC} \times \grp{PGL}{n+1,\bC}$, a $\cpxrepn{\fg}$-representation integrates to $\cpxrepn{G}$ if and only if the sum of the weight coefficients on each branch is even \cite{fh1991-repntheory}.

The standard representation of $\fg = \alg{sl}{n+1,\bC}$ on $\bT \defeq \bC^{n+1}$ has highest weight
\begin{equation*}
  \bT \defeq \dynkinAA[1,0,0,0, 0,0,0,0]{2}{0}{2}{cproj-para-tractor-T}
\end{equation*}
and therefore does not integrate to $G$.  The associated bundle $\cT \defeq \assocbdl{F^{\tilde{P}}}{\tilde{P}}{\bT}$ is called the \emph{standard tractor bundle} of \cproj\ geometry, a complex bundle of (real) rank $2n+2$ over $M$.  The decomposition \eqref{eq:cproj-para-sldecomp} implies that an element $\inlinematrix{ (\det[\bC] C)^{-1} & \alpha \\ 0 & C } \in\tilde{P}$ acts on $\bT_0 \defeq \linspan[\bC]{v_0}{} \leq \bT$ by multiplication by $(\det[\bC] C)^{-1}$, so that
\begin{equation} \label{eq:cproj-para-tractordecomp}
  \bT_0 = \dynkinAAp[0,0,0,-1, 0,0,0,0]{2}{0}{2}{cproj-para-tractor-T0}
  \quad\text{and}\quad
  \bT/\bT_0 = \dynkinAAp[1,0,0,0, 0,0,0,0]{2}{0}{2}{cproj-para-tractorTT0}
\end{equation}
are the socle and top of the $\fp^{\perp}$-filtration $\bT \supset \bT_0 \supset 0$ of $\bT$.  In particular, $L^{-1,0} \defeq \bT_0$ is an $(n+1)$st root of the bundle $\Wedge{n,0}M$ of holomorphic $n$-forms.  Writing $L^{0,-1} \defeq \conj{L^{-1,0}}$ and $L^{1,0} \defeq (L^{-1,0})^*$, the representation $L^{p,q} \defeq (L^{1,0})^{\tens p} \etens (L^{0,1})^{\tens q}$ has highest weight
\begin{equation*}
  L^{p,q} \defeq \dynkinAAp[0,0,0,"$p$", 0,0,0,"$q$"]{2}{0}{2}{cproj-para-tractor-Lpq}
\end{equation*}
for all $(p,q) \in \bZ^2$; thus $L^{1,1} \isom \cpxrepn{L}$, the complexification of \eqref{eq:cproj-para-repnL}.  A choice of Weyl structure then gives a grading $\cT \isom (\cL^{-1,0} \tens T^{1,0}M) \dsum \cL^{-1,0}$, where $\cL^{p,q}$ is the bundle associated to $L^{p,q}$.  As in \thref{rmk:proj-para-tractorjet}, the choice of a representation $L^{-1,0}$ with $L^{-1,0} \etens \conj{L^{-1,0}} \isom \cpxrepn{L^*}$ is equivalent to the extension of the Cartan bundle to structure group $\tilde{P}$; see \cite[\S 3.1]{cemn2015-cproj}.

\subsection{Harmonic curvature} 
\label{ss:cproj-para-harm}

We next compute the harmonic curvature of the canonical Cartan connection, which lies in the Lie algebra homology $\liehom{2}{\fg}$.  The adjoint representation of $\fg$ is given by
\vspace{0.1em}
\begin{equation*}
  \alg{sl}{n+1,\bC} = \Realrepn{ \dynkinAA[1,0,0,1, 0,0,0,0]{2}{0}{2}{cproj-para-harm-adj1}
                           \dsum \dynkinAA[0,0,0,0, 1,0,0,1]{2}{0}{2}{cproj-para-harm-adj2} }
    \colvectpunct[-1.4em]{\, ,}
\vspace{0.1em}
\end{equation*}
so that the beginning of the Hasse diagram computing this homology is given by taking the real $\fp$-representation underlying the direct sum of each complex $\cpxrepn{\fp}$-representation in Figure \ref{fig:cproj-para-hasse} and its complex conjugate.  Each complex representation in the third column of Figure \ref{fig:cproj-para-hasse} contributes a component to the harmonic curvature, which is a section of the real bundle associated to the direct sum of that representation and its complex conjugate.  From top to bottom these bundles are
\vspace{0.1em}
\begin{equation*} \begin{aligned}
  \realrepn{\big( \Wedge{2,0}M \cartan[\bC] \alg{sl}{T^{1,0}M}
       &\dsum \Wedge{0,2}M \cartan[\bC] \alg{sl}{T^{0,1}M} \big)}, \\
  \realrepn{\big( \Wedge{1,1}M \cartan[\bC] \alg{sl}{T^{1,0}M}
       &\dsum \Wedge{1,1}M \cartan[\bC] \alg{sl}{T^{0,1}M} \big)} \\
  \text{and}\quad
  \realrepn{\big( \Wedge{0,2}M \cartan[\bC] T^{1,0}M
       &\dsum \Wedge{2,0}M \cartan[\bC] T^{0,1}M \big)},
\end{aligned}
\vspace{0.1em}
\end{equation*}
which may be identified with subbundles of $\Wedge[-]{2}T^*M \cartan \alg{sl}{TM,J}$, $\Wedge[+]{2}T^*M \cartan \alg{sl}{TM,J}$ and $\Wedge[-]{2}T^*M \cartan TM$ respectively.  The components of the harmonic curvature lying in $\Wedge[+]{2}T^*M \cartan \alg{sl}{TM,J}$ and $\Wedge[-]{2}T^*M \cartan \alg{sl}{TM,J}$ are called the $(1,1)$- and $(2,0)$-parts $\Weyl[\D+]{}{}$ and $\Weyl[\D-]{}{}$ of the Weyl curvature $\Weyl[\D]{}{}$ of $\D$.  The Weyl curvature is totally \tracefree, complex-linear and satisfies the algebraic symmetries of a curvature tensor.  However, unlike in \proj\ differential geometry, $\Weyl[\D]{}{}$ is not \cproj ly\ invariant in general; this is discussed further below.  The $(2,0)$-part $\Weyl[\D-]{}{}$ vanishes if and only if the \cproj\ class contains a holomorphic connection \cite{cemn2015-cproj}.

The component in $\Wedge[-]{2}T^*M \cartan TM$ may be identified with the torsion $\Tor{}$ of a connection in $\Dspace$, which we identify with the Nijenhuis torsion $\Nijen{}$ of $J$.  This component vanishes if and only if $J$ is integrable.  As usual, the harmonic curvature vanishes entirely if and only if $(M,J,\Dspace)$ is locally isomorphic to the flat model $\CP[n]$.

By \thref{thm:para-bgg-ablcurv}, the curvature tensor $\Curv{}{}$ of a Weyl connection $\D \in \Dspace$ decomposes as $\Curv{}{} = \Weyl[\D]{}{} - \algbracw{\id}{\nRic{}{}}{}$, where $\nRic{}{} \defeq -\quab_M^{-1}\p\Curv{}{}$ is the normalised Ricci tensor of $\D$.  Since $\weyld{\gamma} \Weyl[\D]{}{} = \algbrac{\Tor{}}{}$ by \itemref{thm:para-bgg-ablcurv}{weyld} and $\Tor{}$ has type $(0,2)$, we conclude that $\Weyl[+]{}{} \defeq \Weyl[\D+]{}{}$ is \cproj ly\ invariant and $\weyld{\gamma} \Weyl[\D-]{}{} = \algbrac{\Tor{X,Y}}{\gamma}$; in particular, $\Weyl[\D-]{}{}$ is \cproj ly\ invariant if and only if $J$ is integrable.

\vspace{1.5em}
\begin{figure}[h]
  \begin{equation*}
  \arraycolsep=0.1em
  \begin{array}{cccccc}
    \dynkinAAp[1,0,0,1, 0,0,0,0]{2}{0}{2}{cproj-para-harm-hasse1}     
    &
    \begin{array}{c}                                
    \dynkin{ \DynkinConnector{0}{ 0.2}{1.5}{ 2.1};
             \DynkinConnector{0}{-0.2}{1.5}{-2.1};
           }{cproj-para-harm-hasse2}
    \end{array}
    &
    \begin{array}{c}                                
    \dynkinAAp[1,0,0,2,-3, 0,0,0,0, 0]{2}{0}{3}{cproj-para-harm-hasse3} \\[2.5em]
    \dynkinAAp[1,0,0,0, 1, 0,0,0,1,-2]{2}{0}{3}{cproj-para-harm-hasse4}
    \end{array}
    &
    \begin{array}{c}                                
    \dynkin{ \DynkinConnector{0}{ 2.7}{1.5}{ 4.6};
             \DynkinConnector{0}{ 2.3}{1.5}{ 0.4};
             \DynkinConnector{0}{-1.9}{1.5}{ 0  };
             \DynkinConnector{0}{-2.3}{1.5}{-4.2};
           }{cproj-para-harm-hasse5}
    \end{array}
    &
    \begin{array}{c}                                
      \dynkinAAp[1,0,0,1,1,-4, 0,0,0,0,0, 0]{2}{0}{4}{cproj-para-harm-hasse6} \\[2.5em]
      \dynkinAAp[1,0,0,0,2,-3, 0,0,0,0,1,-2]{2}{0}{4}{cproj-para-harm-hasse7} \\[2.5em]
      \dynkinAAp[1,0,0,0,0, 1, 0,0,0,1,0,-3]{2}{0}{4}{cproj-para-harm-hasse8}
    \end{array}
    &
    \begin{array}{c}                                
    \dynkin{ \DynkinConnector{0}{ 4.6}{1.5}{ 2.7};
             \DynkinConnector{0}{ 0.4}{1.5}{ 2.3};
             \DynkinConnector{0}{-0  }{1.5}{-1.9};
             \DynkinConnector{0}{-4.2}{1.5}{-2.3};
             \DynkinLabel{$\cdots$}{2.5}{ 1.7};
             \DynkinLabel{$\cdots$}{2.5}{-2.9};
           }{cproj-para-harm-hasse9}
    \end{array}
  \end{array}
  \end{equation*}
  \vspace{0.5em}
  \caption[The Hasse diagram computing $\liehom{}{\alg{sl}{n+1,\bC}}$]
          {The Hasse diagram of $\alg{sl}{n+1,\bC}\etens\bC$ (drawn here for $n\geq 5$) which, together with its complex-conjugate, computes the homology $\cpxrepn{\liehom{}{\fg}}$.}
  \vspace{1em}
  \label{fig:cproj-para-hasse}
\end{figure}

The Cotton--York tensor $\CY{}{} \defeq \d^{\D} \nRic{}{}$ of $\D$ splits into components $\CY[\D\pm]{}{}$ according to the decomposition \eqref{eq:cproj-cpx-2forms} of $\Wedge{2}T^*M$.  For later use we collect together some curvature identities, the proof of which are similar, modulo handling torsion, to the corresponding results of \thref{prop:proj-para-nric,prop:proj-para-calc}.  Detailed proofs may be found in \cite{cemn2015-cproj}.

\medskip
\begin{prop} \thlabel{prop:cproj-para-calc} Let $\D\in \Dspace$ be a Weyl connection.  Then:
\begin{enumerate}
  \item \label{prop:cproj-para-calc-Wpm}
  $\Weyl[+]{}{}$ is \cproj ly\ invariant and valued in $\alg{sl}{TM,J}$, while $\Weyl[\D-]{}{}$ is c-projec-tively invariant if and only if $J$ is integrable.  Both $\Weyl[\D\pm]{}{}$ satisfy $\liebdy \Weyl[\D\pm]{}{} = 0$, with $\Weyl[+]{}{}$ totally \tracefree; if $J$ is integrable then $\Weyl[\D-]{}{}$ is also totally \tracefree.
  
  \item \label{prop:cproj-para-calc-Wb}
  $\Weyl[+]{}{}$ and $\Weyl[\D-]{}{}$ satisfy the Bianchi identities
    \vspace{-0.4em}
    \begin{equation} \label{eq:cproj-para-Wbianchi} \begin{aligned}
    \Weyl[+]{X,Y}{Z} + \Weyl[+]{Y,Z}{X} + \Weyl[+]{Z,X}{Y} &= 0 \\
    \text{and}\quad
    \Weyl[\D-]{X,Y}{Z} + \Weyl[\D-]{Y,Z}{X} + \Weyl[\D-]{Z,X}{Y} &= (\d^{\D}\Tor{})_{X,Y,Z}.
    \end{aligned}
    \vspace{-0.4em}
    \end{equation}
  
  \item \label{prop:cproj-para-calc-Cb}
  $\CY{}{}$ and $\CY{}{}\circ J$ also satisfy Bianchi identities
    \vspace{-0.4em}
    \begin{equation} \label{eq:cproj-para-CYbianchi} \begin{aligned}
    \CY{X,Y}{Z} + \CY{Y,Z}{X} + \CY{Z,X}{Y} &= 0 \\
    \text{and}\quad
    \CY{X,Y}{JZ} + \CY{Y,Z}{JX} + \CY{Z,X}{JY} &= 0.
    \end{aligned}
    \vspace{-0.4em}
    \end{equation}
  
  \item \label{prop:cproj-para-calc-db}
  We have $\ve^i(\D_{e_i}\Weyl[+]{X,Y}{}) = (n-1)\CY[\D+]{X,Y}{}$ and $\ve^i(\D_{e_i}\Weyl[\D-]{X,Y}{}) = (n-2)\CY[\D-]{X,Y}{}$ \wrt\ any local frame $\{e_i\}_i$ with dual coframe $\{\ve^i\}_i$.
  
  \item \label{prop:cproj-para-calc-nric}
  $\nRic{}{}$ is related to the Ricci curvature $\liebdy \Curv{}{}$ of $\D$ by
    \vspace{-0.4em}
    \begin{equation*}
    \nRic{}{} = -\tfrac{1}{n-1}(\sym \liebdy \Curv{}{}) - \tfrac{1}{n+1}(\alt \liebdy \Curv{}{})
      + \tfrac{2}{n^2-1}(\sym \liebdy \Curv{}{})^{+} \\
    \vspace{-0.4em}
    \end{equation*}
  In particular, if $\nRic{}{}$ is symmetric and $J$-invariant then $\nRic{}{} = -\tfrac{1}{n+1} \liebdy \Curv{}{}$. \noproof
\end{enumerate}
\end{prop}

\section{Associated BGG operators} 
\label{s:cproj-bgg}

Metrisability of a \cproj\ structure may be handled in a similar way to metrisability of real \proj\ structures.  As in the real \proj\ case, the flat model $\CP$ embeds into the projectivisation of an irreducible $\fg$-representation $\bW$, and the first BGG operator associated to $\bW$ has kernel isomorphic to the space of compatible \kahler\ metrics.  We will make this correspondence explicit in Subsection \ref{ss:cproj-bgg-metric} below.

The dual representation $\bW^*$ also has an associated first BGG operator, which, for reasons discussed in Subsection \ref{ss:cproj-bgg-hess}, we did not consider for \proj\ differential geometry.  In \cproj\ geometry this operator, called the \emph{\cproj\ hessian}, controls which Weyl connections have symmetric $J$-invariant normalised Ricci tensor.

For the remainder of this chapter we assume that $J$ is integrable, so that $\Weyl[-]{}{} \defeq \Weyl[\D-]{}{}$ and hence $\Weyl{}{} \defeq \Weyl[\D]{}{}$ are \cproj ly\ invariant.

\subsection{Metrisability of \cproj\ structures} 
\label{ss:cproj-bgg-metric}

Let $(M,J,\Dspace)$ be a \cproj\ manifold of dimension $2n$.  Describing the metrisability of $\Dspace$ proceeds in much the same way as in the real \proj\ case, except that now we are concerned with (\pseudo)\kahler\ metrics.  Such metrics are smooth sections of $\Symm[+]{2} T^*M$, suggesting we consider the natural decomposition
\begin{equation} \label{eq:cproj-bgg-tracesummands}
  T^*M \tens \Symm[+]{2}TM = (\id\symm TM)_{+} \dsum (T^*M \tens[\trfree] \Symm[+]{2}TM)
\end{equation}
in place of \eqref{eq:proj-bgg-tracesummands}, where the first summand is the image of $Z\mapsto \id\symm Z+J\symm JZ$ and the second summand is the kernel of the natural trace $T^*M\tens \Symm[+]{2}TM \surjto TM$.  Projection onto $T^*M\tens[\trfree]\Symm[+]{2}TM$ in \eqref{eq:cproj-bgg-tracesummands} shall be denoted by the subscript ``$\trfree$''.  The proof of the following proposition is analogous to the proof of \thref{prop:proj-bgg-metriceqn}.

\begin{prop} \thlabel{prop:cproj-bgg-metriceqn} The first-order linear differential equation $(\D h)_{\trfree}=0$ is c-projec-tively invariant on sections of $\cL^*\tens \Symm[+]{2} TM$. \noproof \end{prop}

We refer to the equation $(\D h)_{\trfree}=0$ as the \emph{linear metric equation} of \cproj\ geometry, and its solutions as \emph{linear metrics}.  We may equivalently write
\begin{equation} \label{eq:cproj-bgg-metriceqnZ}
  \D_X h = X\symm Z^{\D} + JX\symm JZ^{\D}
\end{equation}
for some $Z^{\D} \in \s{0}{\cL^*\tens TM}$ and all $X\in\s{0}{TM}$, which should be compared with \eqref{eq:proj-bgg-metriceqnZ} and \eqref{eq:cproj-class-maineqn}. Taking a trace in \eqref{eq:cproj-bgg-metriceqnZ} easily yields $Z^{\D} = \tfrac{2}{n}\p(\D h)$.

Identifying $\Wedge{2n}TM \isom \cL^{n+1}$, a \non degenerate\ \kahler\ metric $g$ induces a section of $\cL^*\tens\Symm[+]{2}TM$ defined by
\begin{equation*}
  h \defeq (\det g)^{1/(2n+2)}g^{-1},
\end{equation*}
which we call the \emph{linear metric} associated to $g$.  Then $\det h = (\det g)^{-1/(n+1)}$ is a section of $\cL^2$, so that we may recover $g = (\det h)^{-1/2} \ltens h^{-1}$ from $h$; \cf\ equation \eqref{eq:cproj-class-gbar} and the discussion prior to \thref{cor:proj-bgg-metriceqn}.

\begin{cor} \thlabel{cor:cproj-bgg-metriceqn} There is a linear isomorphism between solutions of the linear metric equation and almost complex metric connections in $\Dspace$. \noproof \end{cor}

Thus we have reduced the metrisability problem for \cproj\ structures to the study of a \cproj ly\ invariant first-order linear differential equation.  Since the linear metric equation is over-determined, we next seek its prolongation.

\begin{thm} \thlabel{thm:cproj-bgg-metricprol} There is a linear isomorphism between the space of solutions of the linear metric equation and parallel sections of the \cproj\ invariant connection
\vspace{-0.2em}
\begin{equation} \label{eq:cproj-bgg-metricprol}
  \D^{\cW}_X \! \colvect{ h \\ Z \\ \lambda }
  = \colvect{ \D_X h - X\symm Z - JX\symm JZ \\
              \D_X Z - h(\nRic{}{},\bdot) - \lambda \ltens X \\
              \D_X \lambda - \nRic{X}{Z}
            }
  - \tfrac{1}{n} \!
    \colvect{ 0 \\
                -\Weyl[+]{e_i,X}{h(\ve^i,\bdot)} \\
                h(\CY[\D+]{e_i,X}{},\ve^i) }
\vspace{-0.2em}
\end{equation}
on sections $(h,Z,\lambda)$ of $\cW \defeq (\cL^*\tens \Symm[+]{2}TM)\dsum (\cL^* \tens TM)\dsum \cL^*$. \end{thm}

\begin{proof} The only necessary modification of the proof of \thref{thm:proj-bgg-metricprol} is the following observation.  Since every Weyl connection $\D\in \Dspace$ is complex, we may assume that every local frame $\{e_i\}_i$ satisfies $e_{i+n}=Je_i$.  Then since $\Weyl[-]{X,JY}{} = J\circ\Weyl[-]{X,Y}{}$ by \thref{prop:cproj-para-calc}, we have $\Weyl[-]{e_i,X}{h(\ve^i,\bdot)} = \Weyl[-]{Je_i,X}{h(J\ve^i,\bdot)} = -\Weyl[-]{e_i,X}{h(\ve^i,\bdot)}$ so only the $J$-invariant piece $\Weyl[+]{}{}$ of the Weyl curvature contributes to the curvature correction.  A similar observation applies to the Cotton--York tensor $\CY{}{}$. \end{proof}

In the presence of a background \kahler\ metric, writing the prolongation connection \eqref{eq:cproj-bgg-metricprol} \wrt\ its \LC\ connection yields the prolongation obtained by Domashev and Mike{\v s}.  The invariant version was obtained by Calderbank in the context of \hamiltonian\ $2$-forms \cite{c2012-ham2vecs}, and in the language of \thref{thm:cproj-bgg-metricprol} in \cite{cemn2015-cproj}.

Unsurprisingly, the linear metric may be interpreted as a first BGG operator.  As described in Subsection \ref{ss:cproj-para-flat}, the flat model $\CP$ admits a \proj\ embedding $\CP \injto \pr{\bW}$ for $\bW \defeq \realrepn{( \bC^{n+1}\etens\conj{\bC^{n+1}} )}$.  Identifying $\bC^{n+1}$ with the standard representation $\bT \defeq \bC^{n+1}$ and using \eqref{eq:cproj-para-tractordecomp}, an algebraic Weyl structure yields
\vspace{0.2em}
\begin{equation} \label{eq:cproj-bgg-tractorW} \begin{aligned}
  \cpxrepn{\bW}
  &= \big( (L^{-1,0}\tens T^{1,0}M) \dsum L^{-1,0} \big)
    \etens \big( (L^{0,-1}\tens T^{0,1}M) \dsum L^{0,-1} \big) \\
  &\isom \big( L^{-1,-1}\tens\Symm{1,1}\cpxrepn{(\fg/\fp)} \big) \dsum
    \big( L^{-1,-1}\tens\cpxrepn{(\fg/\fp)} \big) \dsum L^{-1,-1},
\end{aligned}
\vspace{0.2em}
\end{equation}
so that the algebraic Weyl structure gives an isomorphism of the real bundle associated to $\bW$ and the bundle $\cW = (\cL^*\tens\Symm[+]{2}TM) \dsum (\cL^*\tens TM) \dsum \cL^*$ from \thref{thm:cproj-bgg-metricprol}.  The first BGG operator associated to $\bW$ is a differential operator
\vspace{0.25em}
\begin{equation*}
  \bgg{\bW} :
  \dynkinname{ \dynkinSLCp[1,0,0,0, 1,0,0,0]{2}{0}{2}{cproj-bgg-metric-dom} }
             {\zbox{ \cL^*\tens\Symm[+]{2}TM }}
  \to
  \dynkinname{ \Realrepn[0.8]{ \dynkinAAp[1,0,0,1,-2, 1,0,0,0, 0]{2}{0}{3}{cproj-metric-im1}
                         \dsum \dynkinAAp[1,0,0,0, 0, 1,0,0,1,-2]{2}{0}{3}{cproj-metric-im2} } }
             { (\cL^*\tens\Symm[+]{2}TM) \cartan T^*M }
    \colvectpunct[-1.5em]{~,}
\vspace{0.15em}
\end{equation*}
which is easily seen to be first order.  The prolongation connection $\D^{\cW}$ calculated in \thref{thm:cproj-bgg-metricprol} is then precisely the prolongation connection of this BGG operator.  The dimension of the space of parallel sections of $\D^{\cW}$ is called the \emph{mobility} of the \cproj\ structure, which is bounded above by $\dim\bW = (n+1)^2$.

\begin{rmk} \thlabel{rmk:cproj-ppg-altW} We could equally identify the summand $L^{-1,-1} \tens (T^{1,0}M \etens T^{0,1}M)$ in \eqref{eq:cproj-bgg-tractorW} with $L^{-1,-1}\tens\Wedge{1,1} \cpxbdl{TM}$; indeed the $\fp$-representations $\Symm[+]{2}(\fg/\fp)$ and $\Wedge[+]{2}(\fg/\fp)$ are isomorphic via $J$.  This explains the existence of the theory of hamiltonian $2$-forms. \end{rmk}

Recall from Section \ref{s:proj-bgg} that in \proj\ differential geometry, the $\alg{sl}{n+1,\bR}$-representation $\Symm{2}\bR^{n+1} \dsum \alg{gl}{n+1,\bR} \dsum \Symm{2}\bR^{n+1*}$ has a graded Lie algebra structure isomorphic to $\alg{sp}{2n+2,\bR}$.  It turns out that a similar phenomenon occurs in \cproj\ geometry: since $\cpxrepn{\fg} = \alg{sl}{n+1,\bC}\dsum\alg{sl}{n+1,\bC}$, we have
\vspace{-0.4em}
\begin{equation*}
  \cpxrepn{\fg} \dsum\bC = \alg{s}{\alg{gl}{n+1,\bC}\dsum\alg{gl}{n+1,\bC}}
\end{equation*}
and the complexification of $\fh \defeq \bW \dsum (\fg\dsum\bR) \dsum \bW^*$ may be written as
\begin{equation} \label{eq:cproj-bgg-hdecomp} \begin{split}
  \cpxrepn{\fh}
  &\isom (\bC^{n+1}\etens \conj{\bC^{n+1*}}) \dsum
    \big( (\bC^{n+1}\cartan\bC^{n+1*}) \dsum (\conj{\bC^{n+1}}\cartan\conj{\bC^{n+1*}}) \big) \\
  &\qquad \dsum (\bC^{n+1*}\etens\conj{\bC^{n+1*}}) \\
  &\isom (\bC^{n+1}\dsum\conj{\bC^{n+1}}) \cartan (\bC^{n+1*}\dsum\conj{\bC^{n+1*}}),
\end{split}
\end{equation}
where we have used that the Cartan product $\bC^{n+1}\cartan\bC^{n+1*}$ is the space of \tracefree\ complex-linear maps of $\bC^{n+1}$ and that $\bC^{n+1}\cartan\conj{\bC^{n+1}} = \bC^{n+1}\etens \conj{\bC^{n+1}}$.  This is precisely the adjoint representation of the complex Lie algebra $\alg{sl}{\bC^{n+1}\dsum \conj{\bC^{n+1*}}}$, so that $\cpxrepn{\fh}$ has a graded Lie algebra structure isomorphic to $\alg{sl}{2n+2,\bC}$.  The grading implies that $\cpxrepn{\bW}, \cpxrepn{\bW}^*$ are abelian subalgebras of $\cpxrepn{\fh}$, while $\cpxrepn{\fq} \defeq (\cpxrepn{\fg}\dsum\bC) \ltimes \bW^*$ and $\cpxrepn{\opp{\fq}} \defeq \bW \rtimes (\cpxrepn{\fg}\dsum\bC)$ are opposite abelian parabolics.  Given the evident real structure present in \eqref{eq:cproj-bgg-hdecomp}, we see that $\fh$ is isomorphic to the real form
\vspace{0.05em}
\begin{equation*}
\smash{
  \fq = \dynkinSUp{2}{0}{2}{cproj-bgg-metric-q} \hspace{-0.5em}
    \leq~ \dynkinSU{2}{0}{2}{cproj-bgg-metric-h}
    \hspace{-1.5ex}= \alg{su}{n+1,n+1} \isom \fh
}
\vspace{0.05em}
\end{equation*}
of $\alg{sl}{2n+2,\bC}$, so that $H\acts\fq$ is \grassmannian\ of maximal isotropic subspaces of $\bC^{2n+2}$ for a hermitian inner product of signature $(n+1,n+1)$; see \cite[Ex.\ 2.2.2]{bdpp2011-rspaces} and \cite[Ex.\ 2.3.4(3)]{cs2009-parabolic1} for details regarding $\alg{su}{n+1,n+1}$.

\subsection{The \cproj\ hessian} 
\label{ss:cproj-bgg-hess}

In the previous subsection, we showed that the metrisability of a \cproj\ structure is controlled by the first BGG operator associated to the $\fg$-representation $\bW \defeq \realrepn{( \bC^{n+1} \etens \conj{\bC^{n+1}} )}$.  There is also a first BGG operator associated to $\bW^*$, which we did not consider in the real \proj\ case.  By way of justification for this omission recall that in \proj\ differential geometry $\bW \defeq \Symm{2}\bR^{n+1}$, so that the first BGG operator associated to $\bW^*$ is a linear differential operator
\vspace{-0.1em}
\begin{equation*}
  \bgg{\bW^*} : \dynkinSLRp[0,0,0,0,2]{2}{0}{3}{} \to \dynkinSLRp[0,0,0,3,-4]{2}{0}{3}{}.
\vspace{-0.2em}
\end{equation*}
Using the inverse Cartan matrix of $\alg{sl}{n+1,\bR}$, we easily see that $\bgg{\bW^*}$ is third-order.  Choosing an algebraic Weyl structure for $\bW$ we see that $\cW^* \isom \cL \dsum (\cL\tens T^*M) \dsum (\cL\tens\Symm{2}T^*M)$ is the full $2$-jet bundle $\jetbdl{2}{\cL}$ of $\cL$, so that the first BGG operator associated to $\bW^*$ is simply the third-order Ricci-corrected derivative $\bgg{\bW^*}(\ell) = \D^3 \ell + 2(\D\nRic{}{})\ell + 4\D\ell \ltens \nRic{}{}$; see \cite[p.\ 169]{cds2005-ricci}.  For $\ell\in\s{0}{\cL}$ \non vanishing, writing this equation \wrt\ the special connection $\D^{\ell} \in \Dspace$ defined by $\D^{\ell}\ell = 0$ yields $\bgg{\bW^*}(\ell) = 2(\D \nRic{}{})\ell$, so that the kernel of $\bgg{\bW^*}$ may be interpreted as a space of ``\einstein\ scales'' \cite{cgh2012-projbgg,gm2015-weylnullity}; there is also a relation to Tanno equations \cite{t1978-someeqns}.

The \cproj\ hessian has a significantly different character, which is arguably more important to the underlying \cproj\ geometry.  Here $\bW^* = \realrepn{( \bC^{n+1*}\etens\conj{\bC^{n+1*}} )}$ so that the associated first BGG operator is a linear differential operator
\vspace{0.25em}
\begin{equation*}
  \bgg{\bW^*} :
    \dynkinname{ \dynkinSLCp[0,0,0,0,1, 0,0,0,0,1]{2}{0}{3}{cproj-bgg-hess-dom} }
               { \cL }
  \To
    \dynkinname{ \Realrepn[0.8]{ \dynkinAAp[0,0,0,1,-3, 0,0,0,0, 1]{2}{0}{3}
                                           {cproj-bgg-hess-im1}
                           \dsum \dynkinAAp[0,0,0,0, 1, 0,0,0,1,-3]{2}{0}{3}
                                           {cproj-bgg-hess-im2} } }
               { \cL\tens\Symm[-]{2}T^*M }
\vspace{0.15em}
\end{equation*}
called the \emph{\cproj\ hessian}.  Using the inverse Cartan matrix of $\cpxrepn{\fg}$, it is straightforward to check that $\bgg{\bW^*}$ must be second order; it follows that $\bgg{\bW^*}$ is projection onto $\cL\tens \Symm[-]{2}T^*M$ of the Ricci-corrected second derivative \cite{cds2005-ricci}, \ie\
\vspace{0.2em}
\begin{equation} \label{eq:cproj-bgg-hesseqn}
  \bgg[\bW^*]{}_{X,Y}(\ell) = (\D^2_{X,Y}\ell + \ell\ltens\nRic{X}{Y})_{-}.
\vspace{0.2em}
\end{equation}
It is straightforward to check that $\bgg{\bW^*}$ is \cproj ly\ invariant.  By \thref{lem:cproj-class-connsonL}, a section $\ell\in\s{0}{\cL}$ uniquely determines a connection $\D^{\ell}$ in the \cproj\ class by $\D^{\ell}\ell = 0$.  The kernel of $\bgg{\bW^*}$ has the following geometric interpretation \cite[Prop.\ 4.9]{cemn2015-cproj}.

\smallskip
\begin{prop} \thlabel{prop:cproj-bgg-hesseqn} A nowhere-vanishing section $\ell\in\s{0}{\cL}$ satisfies $\bgg{\bW^*}(\ell) = 0$ if and only if the normalised Ricci tensor $\nRic[\D^{\ell}]{}{}$ of $\D^{\ell}$ is symmetric and $J$-invariant. \end{prop}

\begin{proof} Since $\ell$ is nowhere-vanishing, $\Curv[\D^{\ell}]{}{\ell}=0$ if and only if $\nRic[\D^{\ell}]{}{}$ is symmetric.  Calculating \wrt\ $\D^{\ell}$, \eqref{eq:cproj-bgg-hesseqn} gives $\bgg{\bW^*}(\ell) = 0$ if and only if $\nRic[\D^{\ell}]{}{}$ has vanishing $J$-anti-invariant part. \end{proof}

\begin{cor} \thlabel{cor:cproj-bgg-hessmetric} Let $h \in \s{0}{\cL^*\tens\cB}$ be a linear metric.  Then $(\det h)^{1/2}$ lies in the kernel of the \cproj\ hessian. \end{cor}

\begin{proof} By \thref{cor:cproj-bgg-metriceqn}, the \LC\ connection $\D^g$ of the \kahler\ metric $g \defeq (\det h)^{-1/2} \ltens h^{-1/2}$ lies within the \cproj\ class and satisfies $\D^g(\det h)^{-1/2}=0$.  Since the normalised Ricci tensor $\nRic[g]{}{}$ of $g$ is symmetric and $J$-invariant, the result follows by \thref{prop:cproj-bgg-hesseqn}. \end{proof}

\begin{rmk} \thlabel{rmk:cproj-bgg-dim2} As remarked in the introduction to Section \ref{s:cproj-para}, all connections on a Riemann surface are \cproj\ and hence the case $n=1$ must be dealt with differently.  Indeed there is an isomorphism $\CP[1] \isom \Sph[2]$, so that $1$-dimensional \cproj\ geometry coincides with $2$-dimensional conformal geometry until we specify an additional piece of structure.  This structure is called a \emph{M{\"o}bius structure}, and is essentially the choice of a \cproj\ hessian; see \cite[\S 6.3]{bc2010-conf} and \cite{c1998-2dewgeom} for details. \end{rmk}

Since $\bgg{\bW^*}$ is a first BGG operator, its solution space is isomorphic to the space of parallel sections of a connection $\D^{\cW^*}$ on $\cW^* \isom \cL \dsum (\cL\tens T^*M) \dsum (\cL\tens \Symm[+]{2}T^*M)$.

\begin{thm} \thlabel{thm:cproj-bgg-hessprol} There is a linear isomorphism between the space of solutions of the \cproj\ hessian $\bgg{\bW^*}$ and parallel sections of the \cproj ly\ invariant connection
\begin{equation} \label{eq:cproj-bgg-hessprol}
  \D^{\cW^*}_X \! \colvect{ \ell \\ \eta \\ \theta }
  = \colvect{ \D_X\ell - \eta(X) \\
              \D_X\eta + \ell\ltens\nRic{X}{} - \theta(X,\bdot) \\
              \D_X\theta + \nRic{X}{}\symm\eta + J\nRic{X}{}\symm J\eta
            }
  + \colvect{ 0 \\
              0 \\
              (\Weyl[+]{\cdot,J\cdot}{\eta})(JX)
                + \ell\ltens\CY[\D+]{\cdot,J\cdot}{JX}
            }
\end{equation}
on sections $(\ell,\eta,\theta)$ of $\cW^*\isom \cL\dsum (\cL\tens T^*M)\dsum (\cL\tens \Symm[+]{2}T^*M)$. \end{thm}

\begin{proof} Define $\eta^{\D} \defeq \D\ell$ and $\theta^{\D} \defeq (\D^2\ell + \ell\ltens\nRic{}{})_{+}$, so that $\theta^{\D}= \D\eta^{\D} + \ell\ltens\nRic{}{}$ if and only if $\ell$ is in the kernel of $\bgg{\bW^*}$.  It remains to establish the third slot of $\D^{\cW^*}$.

Note first that $\D_X\theta^{\D} = \D_X(\D\eta^{\D}) + \ell\ltens(\D_X\nRic{}{}) + \eta^{\D}(X)\nRic{}{}$, which contracting with $Y\in\s{0}{TM}$ and alternating in $X,Y$ yields
\begin{equation} \label{eq:cproj-bgg-hessprol-1} \begin{split}
  (\d^{\D}\theta^{\D})_{X,Y}
  &= \Weyl{X,Y}{\eta^{\D}} + \ell\ltens\CY{X,Y}{}
    + \algbracw{\id}{\nRic{}{}}{X,Y}\acts\eta^{\D} \\
  &\qquad + \eta^{\D}(X)\nRic{Y}{} - \eta^{\D}(Y)\nRic{X}{}.
\end{split} \end{equation}
Expanding the algebraic bracket and simplifying, the last three terms on the \rhs\ above evaluate to
\begin{align*}
  &\algbracw{\id}{\nRic{}{}}{X,Y} \acts \eta^{\D} + \eta^{\D}(X)\nRic{Y}{}
    - \eta^{\D}(Y)\nRic{X}{} \\
  &\quad = (\nRic{Y}{}\symm \eta^{\D} + J\nRic{Y}{}\symm J\eta^{\D})(X,Z)
    - (\nRic{X}{}\symm \eta^{\D} + J\nRic{X}{}\symm J\eta^{\D})(Y,Z).
\end{align*}
Since $\theta^{\D}$ is symmetric and $J$-invariant, $(\D_X\theta^{\D})(JY,Z)$ is skew in $Y,Z$ and hence
\begin{align*}
  2(\D_X\theta^{\D})(Y,Z) &= (\D_X\theta^{\D})(Y,Z) + (\D_X\theta^{\D})(JZ,JY) \\
  &= (\d^{\D}\theta^{\D})_{X,Y}(Z) + (\D_Y\theta^{\D})(X,Z) \\
  &\qquad + (\d^{\D}\theta^{\D})_{X,JZ}(JY) + (\D_{JZ}\theta^{\D})(X,JY) \\
  &= (\d^{\D}\theta^{\D})_{X,Y}(Z) + (\d^{\D}\theta^{\D})_{X,JY}(JZ) + (\d^{\D}\theta^{\D})_{Y,JZ}(JX).
\end{align*}
Substituting \eqref{eq:cproj-bgg-hessprol-1} in the last display, we find that
\begin{equation} \label{eq:cproj-bgg-hessprol-2} \begin{split}
  2(\D_X\theta^{\D})(Y,Z)
    &= (\Weyl[+]{X,Y}{\eta^{\D}})(Z) +
      (\Weyl[+]{X,JZ}{\eta^{\D}})(JY) + (\Weyl[+]{Y,JZ}{\eta^{\D}})(JX)
      \\
    &\qquad + \ell\ltens\CY[\D+]{X,Y}{Z} + \ell\ltens\CY[\D+]{X,JZ}{JY} +
      \ell\ltens\CY[\D+]{Y,JZ}{JX}  \\
    &\qquad - 2(\nRic{X}{}\symm\eta^{\D}
      + J\nRic{X}{}\symm J\eta^{\D})(Y,Z).
  \end{split}
\end{equation}
where the $J$-anti-invariant parts $\Weyl[-]{}{}$ and $\CY[\D-]{}{}$ vanishes automatically because the \lhs\ is $J$-invariant.  Using the Bianchi identity \eqref{eq:cproj-para-Wbianchi} for $\Weyl[+]{}{}$, the first line on the \rhs\ of \eqref{eq:cproj-bgg-hessprol-2} simplifies to
\begin{equation} \label{eq:cproj-bgg-hessprol-3} \begin{aligned}
  &(\Weyl[+]{X,Y}{\eta^{\D}})(Z) + (\Weyl[+]{X,JZ}{\eta^{\D}})(JY)
    + (\Weyl[+]{Y,JZ}{\eta^{\D}})(JX) \\
  &\quad = -\eta^{\D}( -\Weyl[+]{Y,Z}{X} - \Weyl[+]{Z,X}{Y} - \Weyl[+]{JZ,JY}{X} \\
  &\quad \hspace{4em} - \Weyl[+]{JY,X}{JZ} + \Weyl[+]{Y,JZ}{JX} ) \\
  &\quad = 2(\Weyl[+]{Y,Z}{\eta^{\D}})(X)
\end{aligned} \end{equation}
while the second line simplifies to
\begin{equation} \label{eq:cproj-bgg-hessprol-4} \begin{aligned}
  &\ell\ltens\CY[\D+]{X,Y}{Z} + \ell\ltens\CY[\D+]{X,JZ}{JY} + \ell\ltens\CY[\D+]{Y,JZ}{JX} \\
  &\quad = -\ell\ltens\CY[\D+]{Y,Z}{X} - \ell\ltens\CY[\D+]{Z,X}{Y}
    - \ell\ltens\CY[\D+]{JZ,JY}{X} \\
  &\quad\qquad - \ell\ltens\CY[\D+]{JY,X}{JZ} + \ell\ltens\CY[\D+]{Y,JZ}{JX} \\
  &\quad = 2\ell\ltens\CY[\D+]{Y,Z}{X}
\end{aligned} \end{equation}
by \eqref{eq:cproj-para-CYbianchi}.  Substituting \eqref{eq:cproj-bgg-hessprol-3} and \eqref{eq:cproj-bgg-hessprol-4} into \eqref{eq:cproj-bgg-hessprol-2}, we arrive at the third slot of \eqref{eq:cproj-bgg-hessprol}.  \Cproj\ invariance may be checked directly, as in the proof of \thref{thm:proj-bgg-metricprol}, or by observing that $\D^{\cW^*}$ is the prolongation connection of a first BGG operator. \end{proof}

It is straightforward to check that the tractor parts of $\D^{\cW^*}$ and $\D^{\cW}$ are dual, as expected.  In particular, the prolongation connections $\D^{\cW}$ and $\D^{\cW^*}$ are dual on the flat model.

\chapter{\Qtn ic\ geometry} 
\label{c:qtn}

\renewcommand{\algbracadornment}{q}
\BufferDynkinLocaltrue
\renewcommand{\dynkinnameoffset}{-0.65}

\addtolength{\abovedisplayskip}{-0.2em}
\addtolength{\belowdisplayskip}{-0.2em}

Following Berger's classification of \riem\ holonomy groups \cite{b1955-holonomy} and developments in almost complex geometry, it was natural to consider manifolds with a ``\qtn ic'' structure.  According to definitions given by Salamon \cite{s1986-qkdg} and Alekseevsky and Marchiafava \cite{am1993-qtnlikenote1,am1993-qtnlikenote2}, an almost \qtn ic\ manifold is a first order $G$-structure with structure group $\grp{GL}{n,\bH} \grpdot \grp{Sp}{1}$.  We review the classical theory of almost \qtn ic\ manifolds in Section \ref{s:qtn-class}, as well as describing a \proj\ interpretation which develops in parallel with Sections \ref{s:proj-class} and \ref{s:cproj-class}.

In Section \ref{s:qtn-para} we describe almost \qtn ic\ geometry as an abelian parabolic geometry modelled on $G\acts \fp = \HP$, so that $\fg = \alg{sl}{n+1,\bH}$ and $\fp$ is given by crossing the penultimate node.  The flat model has a \proj\ embedding $\HP \injto \pr{\bW}$ for $\bW \defeq \realrepn{(\Wedge[\bC]{2}\bC^{2n+2})}$, and, as for \proj\ and \cproj\ geometries, the first BGG operator associated to $\bW$ controls metrisability of the \qtn ic\ structure.  The first BGG operator associated to $\bW^*$ is again a second-order hessian equation; we investigate both of these operators in Section \ref{s:qtn-bgg}.

\section{Background on almost \qtn ic\ geometry} 
\label{s:qtn-class}

After briefly reviewing the \qtn\ algebra in Subsection \ref{ss:qtn-class-qtns}, we study the basic theory of almost \qtn ic\ manifolds in Subsection \ref{ss:qtn-class-qtnic}.  This theory is ostensibly better suited to description as a parabolic geometry than the classical theories of \proj\ differential geometry or \cproj\ geometry, and we develop this parabolic viewpoint in Section \ref{s:qtn-para}.  Before that however we give a ``\proj'' interpretation of almost \qtn ic\ geometry, in the style of Sections \ref{s:proj-class} and \ref{s:cproj-class}.

\subsection{Background on \qtn s} 
\label{ss:qtn-class-qtns}

The \emph{\qtn s} are the elements of the (unique up to isomorphism) normed associative division algebra $\bH$ of dimension four over $\bR$.%
\footnote{Hopefully the reader will allow us to forgo the obligatory anecdote about Hamilton and the bridge---but see \cite[\S 1.1.1]{w2000-qtnalggeom} for a nice historical review.}
As a real algebra, $\bH$ is spanned by the real unit $\cpx{1}$ and three imaginary units $\{\cpx{i_a}\}_{a=1}^3$ satisfying Hamilton's \qtn ic\ relations
\vspace{-0.1em}
\begin{equation} \label{eq:qtn-class-qtnrels}
  \cpx{i}_a^2 = -\id, \quad
  \cpx{i}_a\cpx{i}_b = \cpx{i}_c = -\cpx{i}_b\cpx{i}_a
\vspace{-0.1em}
\end{equation}
for all cyclic permutations $(a,b,c)$ of $(1,2,3)$.  In particular, the \qtn s\ are \non commutative.  A generic \qtn\ is of the form $q = q_0\cpx{1} + \qsum q_a\cpx{i}_a$ for real numbers $q_0,q_a$.  Then the euclidean norm gives a notion of a norm in $\bH$, while reversing the signs of imaginary components gives a notion of conjugation.  A \qtn\ is \emph{imaginary} if $\conj{q} = -q$, and the imaginary \qtn s\ form a $3$-dimensional subspace $\alg{sp}{1} \leq \bH$.  In particular if $q \in \alg{sp}{1}$ then $q^2 = -\norm{q}^2 \cpx{1}$, so that the square roots of $-1$ form a $2$-sphere in $\bH$.  Thus the choice of $\{\cpx{i_a}\}_{a=1}^3$ satisfying \eqref{eq:qtn-class-qtnrels} is not canonical, being equivalent to an oriented orthonormal basis of $\alg{sp}{1}$; such a basis determines the isomorphism $\bH \isom \bR^4$.

Due to \non commutativity, the automorphism group of $\bH$ is quite large.  Firstly $q \mapsto \conj{q}$ defines an automorphism of $\bH$ to its opposite algebra; moreover each $q \in \grp{Sp}{1} \defeq \setof{ q\in\bH }{ \norm{q}=1 }$ determines an automorphism $p \mapsto qp\conj{q}$ by conjugation.  The latter yields a homomorphism $\grp{Sp}{1} \to \grp{Aut}{\bH}$ with kernel $\pm\cpx{1}$ which is orthogonal \wrt\ the euclidean norm on $\bH$.  Restricting to $\alg{sp}{1}$ gives a surjective homomorphism $\grp{Sp}{1} \surjto \grp{SO}{\alg{sp}{1}}$ which realises the exceptional isomorphism $\grp{Sp}{1} \isom \grp{Spin}{3}$.  On the other hand if $\phi \in \grp{Aut}{\bH}$ then $\phi$ preserves the $2$-sphere of imaginary units, so is an orthogonal map.  For orthogonal $p,q \in \alg{sp}{1}$ we have $0 \neq pq \in \alg{sp}{1}$, so that $\{p,q,pq\}$ is an oriented orthogonal basis of $\alg{sp}{1}$; in fact all such bases arise in this way.  Due to the large automorphism group, it is interesting to consider a larger class of maps than just the \qtn-linear ones.  Namely, we consider all real linear maps $f:V\to W$ such that $f(vq) = f(v)\phi(q)$ for some $\phi \in\grp{Aut}{\bH}$ and all $v\in V$ and $q\in\bH$.  Here $\bH$ acts on $V,W$ on the right to ensure the usual matrix multiplication conventions.

Finally, a \emph{\qtn ic\ structure} on a real vector space $V$ is a $3$-dimensional subspace $Q \leq \alg{gl}{V}$ admitting a basis $\{J_a\}_{a=1}^3$ satisfying analogues of the \qtn ic\ relations \eqref{eq:qtn-class-qtnrels}.  In particular, the choice of basis $\{J_a\}_{a=1}^3$ for $Q$ is an additional choice.  Evidently a \qtn ic\ vector space is a right $\bH$-module, so is isomorphic to $\bH^n$ for some $n$; in particular, $\dim V = 4n$ is divisible by four.  Choosing a particular $J_a$ allows us to form the complex vector space $\cpxrepn{V}$ in which $\cpx{i} \defeq \sqrt{-1} \in\bC$ acts via $J_a$.  Then $(V,Q)$ may be identified with $\cpxrepn{V}$ together with the conjugate-linear map $\cpx{j}$ determined by the action of $J_b$; this alternative viewpoint will frequently be useful.

\vspace{0.1em}
\subsection{Classical theory} 
\label{ss:qtn-class-qtnic}

We now review the basic theory of almost \qtn ic\ manifolds as described by Alekseevsky and Marchiafava \cite{am1993-qtnlikenote1,am1993-qtnlikenote2,am1996-qtnsubord}.  Let $M$ be a smooth manifold of dimension $m > 4$.  An \emph{almost \qtn ic\ structure} on $M$ is a rank three subbundle $\cQ \leq \alg{gl}{TM}$ with fibres isomorphic to $\alg{sp}{1}$.  A linear connection $\D$ on $TM$ is then called \emph{almost \qtn ic} if $\D$ preserves $\cQ$, \ie\ $\D_X \cQ \subseteq \cQ$ for all $X \in \s{0}{TM}$.  We define an \emph{almost \qtn ic\ manifold} to be such a pair $(M,\cQ)$ which admits an almost \qtn ic\ connection.
Equivalently an almost \qtn ic\ structure is a reduction of the frame bundle of $M$ to structure group $P^0 \defeq \grp{GL}{n,\bH}\grpdot\grp{Sp}{1}$, the quotient of $\grp{GL}{n,\bH}\times\grp{Sp}{1}$ by its natural $\bZ_2$-action, which admits a principal $P^0$-connection satisfying a certain torsion condition; see \cite[\S4]{am1993-qtnlikenote1}.

\begin{rmk} \thlabel{rmk:qtn-class-dim4} In dimension four the structure group $\grp{GL}{1,\bH} \grpdot \grp{Sp}{1}$ is isomorphic to the conformal group $\grp{CO}{4,\bR}$, so that a $4$-dimensional \qtn ic\ manifold is just a $4$-dimensional conformal manifold \cite{g1969-holSpSp}.  To retain the features of almost \qtn ic\ geometry one must stipulate that the anti-\self dual\ part of the conformally invariant Weyl curvature vanishes, so that we define a $4$-dimensional almost \qtn ic\ manifold as a \self dual\ conformal $4$-manifold; see \cite{b1987-einstein} and \cite[\S 4.1.9]{cs2009-parabolic1}, as well as the extensive literature on \self dual\ conformal manifolds. \end{rmk}

Let $(M,\cQ)$ be an almost \qtn ic\ manifold.  Since the fibres of $\cQ$ are isomorphic to $\alg{sp}{1}$, in a neighbourhood of any point in $M$ there is a \emph{local \qtn ic\ frame} $\{J_a\}_{a=1}^3$ of $\cQ$ satisfying the usual multiplicative properties \eqref{eq:qtn-class-qtnrels} of the \qtn s:
\begin{equation} \label{eq:qtn-class-Ja}
  J_aJ_b = -\delta_{ab}\id + \qsum[c] \epsilon_{abc}J_c
\end{equation}
for all $a,b,c \in \{1,2,3\}$, where $\epsilon_{abc}$ is \LC's alternating symbol.  Then the tangent space $T_xM$ is a \qtn ic\ vector space in which application of $J_a\at{x}$ corresponds to multiplication by $\cpx{i_a}$, so that $\dim{M} = 4n$ is divisible by four.  Note that the $J_a$ need not extend to global almost complex structures on $M$: it is well known that \qtn ic\ \proj\ space $\HP[n]$ is an almost \qtn ic\ manifold, but admits no almost complex structures \cite{m1962-qtnprojspace}.

Identifying $\alg{sp}{1} \isom \alg{so}{3}$, the natural action of $\grp{SO}{3}$ on $\alg{gl}{TM}$ preserves $\cQ$.  We therefore have an $\grp{SO}{3}$-freedom in choosing a local \qtn ic\ frame $\{J_a\}_{a=1}^3$; let us fix such a frame once and for all, noting that we must take care to ensure that any expression involving the $J_a$ is invariant under this $\grp{SO}{3}$ action.  For notational convenience we will assume from now on that $(a,b,c)$ is a cyclic permutation of $(1,2,3)$, and occasionally employ the summation convention for \qtn ic\ indices.

It is easy to show that a linear connection $\D$ on $TM$ is almost \qtn ic\ if and only if there are $1$-forms $\{\alpha_a\}_{a=1}^3$ such that
\begin{equation} \label{eq:qtn-class-qtnconn}
  \D J_a = \alpha_b \tens J_c - \alpha_c \tens J_b.
\end{equation}
Then $\D$ is almost \qtn ic\ if and only if the \emph{fundamental form} $\Omega^{\cQ} \defeq \qsum J_a\tens J_a$ is $\D$-parallel.  This form was discovered by Kraines \cite{k1966-qtntopo} and is often written as a $4$-form using a compatible metric.  Note in particular that $\Omega^{\cQ}$ is invariant under the natural action of $\grp{SO}{3}$ on $\cQ$, so that the class $\Dspace$ of almost \qtn ic\ connections is determined by $\cQ$.  Since an almost \qtn ic\ structure has structure group $P^0 \defeq \grp{GL}{n,\bH} \grpdot \grp{Sp}{1}$, the difference of two almost \qtn ic\ connections $\D,\b{\D}$ compatible with $\cQ$ is a $\fp^0$-valued $1$-form, where $\fp^0 \defeq \alg{gl}{n,\bH} \dsum \alg{sp}{1}$.  Fujimura gives the following characterisation of \qtn ically\ equivalent connections \cite{f1976-Qconns}, which is a \qtn ic\ analogue of \thref{lem:proj-class-algbrac,lem:cproj-class-algbrac}; also see \cite{am1994-gradqtnvfs, am1994-qtnlaplacian, am1996-qtnsubord}.

\begin{lem} \thlabel{lem:qtn-class-algbrac} Let $\D$ be an almost \qtn ic\ connection on $(M,\cQ)$ and suppose that $\b{\D}$ is some other linear connection.  Then $\b{\D}$ is almost \qtn ic\ if and only if
\begin{equation} \label{eq:qtn-class-algbrac} \begin{split}
  \b{\D}_X Y &\phantom{:}= \D_X Y + \algbrac{X}{\alpha} \acts Y,
  \quad\text{where} \\
  \algbrac{X}{\alpha} \acts Y &\defeq \tfrac{1}{2}( \alpha(X)Y + \alpha(Y)X
    - \qsum\left[ \alpha(J_aX)J_aY + \alpha(J_aY)J_aX \right])
\end{split} \end{equation}
for some $\alpha\in\s{1}{}$ and all $X,Y\in\s{0}{TM}$. \noproof \end{lem}

Using \eqref{eq:qtn-class-algbrac} and the \qtn ic\ relations \eqref{eq:qtn-class-Ja}, it is easy to see that $\liebrac{ \algbrac{X}{\alpha} }{ J_a } = \alpha(J_bX)J_c - \alpha(J_cX)J_b$,
so that $\algbrac{X}{\alpha}$ takes values in the normaliser of $\cQ$.

\begin{rmk} \thlabel{rmk:qtn-class-algbrac} In keeping with \thref{lem:proj-class-algbrac,lem:cproj-class-algbrac}, we shall call the endomorphism $\algbrac{X}{\alpha}$ the \emph{algebraic bracket} of $X$ and $\alpha$.  Clearly $\algbrac{X}{\alpha}$ is symmetric in $X,Y$.  Note also that since we may write $\algbrac{}{\alpha} = \id\symm\alpha + \qsum J_a\odot J_a\alpha$, the algebraic bracket is independent of the choice of local \qtn ic\ frame $\{J_a\}_{a=1}^3$. \end{rmk}

\thref{lem:qtn-class-algbrac} exhibits $\Dspace$ as an affine space modelled on $1$-forms, with an embedding $\s{1}{} \injto \s{1}{\fp^0}$ given by the algebraic bracket.  This may be understood as follows.  Since an almost \qtn ic\ structure is a $P^0$-structure,  the difference of two $P^0$-connections is a $1$-form valued in the first prolongation $(\fp^0 \tens \bH^{n+1*}) \cap (\bH^{n+1} \tens \Symm{2}\bH^{n+1*})$ of $\fp^0$.  In our case, the first prolongation is isomorphic to $\bH^{n*}$, with isomorphism given explicitly by $\alpha \mapsto \algbrac{}{\alpha}$.  Alekseevsky and Marchiafava \cite{am1996-qtnsubord} give an algebraic proof of this fact, while Salamon gives a representation-theoretic proof \cite{s1986-qkdg}.

As in \proj\ differential geometry and \cproj\ geometry, the following characterisation of almost \qtn ic\ connections is available; \cf\ \thref{lem:proj-class-connsonL,lem:cproj-class-connsonL}.

\begin{lem} \thlabel{lem:qtn-class-connsonL} There is a bijection between linear connections on the line bundle $\cL\defeq (\Wedge{4n}TM)^{1/(2n+2)}$ and connections in $\Dspace$. \noproof \end{lem}

The $1$-forms $\{\alpha_a\}_{a=1}^3$ from \eqref{eq:qtn-class-qtnconn} vanish identically if and only if the induced connection on $\cQ$ is flat, in which case its holonomy reduces to $\grp{GL}{n,\bH} \leq \grp{GL}{n,\bH} \grpdot \grp{Sp}{1}$.  Such manifolds are called \emph{hypercomplex}, since $(M,J)$ is an almost complex manifold \wrt\ each unit norm $J\in\s{0}{\cQ}$; see \cite{j1991-qtnquotient,p2000-hypercpx}.

A quantity on $(M,\cQ)$ is \emph{\qtn ically invariant} if it is independent of the choice of connection in $\Dspace$.  Since $\algbrac{X}{\alpha} \acts Y$ is symmetric in $X,Y$, two almost \qtn ic\ connections $\D,\b{\D}$ have the same torsion; thus the torsion is a \qtn ic\ invariant.  In terms of first-order $P^0$-structures, the torsion coincides with the \emph{intrinsic torsion} of a $P^0$-structure, as determined by the Spencer complex.  We drop the prefix ``almost'' in ``almost \qtn ic'' if $(M,\cQ)$ admits a \torsionfree\ almost \qtn ic\ connection, in which case $\D$ is simply called a \emph{\qtn ic\ connection}.  In this case the \qtn ic\ structure constitutes a ``$1$-integrable'' first-order $P^0$-structure.  Note that, unlike in almost complex geometry, the vanishing of the intrinsic torsion does not guarantee that any member of a \qtn ic\ local frame is integrable, nor that we may choose \qtn ic\ coordinates.  Indeed, Kulkarni \cite{k1978-uniformization} proves that a simply connected almost \qtn ic\ manifold admits \qtn ic\ coordinates if and only if it is diffeomorphic to $\HP$ with its canonical \qtn ic\ structure.

Finally, the natural class of (\pseudo)\riem\ metrics on an almost \qtn ic\ manifold are those $g$ with $g(JX,JY) = g(X,Y)$ for all unit norm $J\in\s{0}{\cQ}$.  Such metrics are called \emph{Q-hermitian}; evidently it suffices that $g(J_aX,J_aY) = g(X,Y)$ for all $a=1,2,3$ and $X,Y\in\s{0}{TM}$.  A Q-hermitian metric $g$ is \emph{(\pseudo)\qk} if its \LC\ connection $\D$ is almost \qtn ic; since a \LC\ connection is \torsionfree, it is immediate that such an almost \qtn ic\ structure is \qtn ic.  Equivalently, $\D$ preserves the subbundle $\cQ$, or the fundamental form $\Omega^{\cQ}$ is $\D$-parallel.  The following is well-known; see \cite[Thm.\ 14.39]{b1987-einstein} for two different proofs.

\begin{prop} \thlabel{prop:qtn-class-einstein} Every \qk\ metric is \einstein. \noproof \end{prop}

In particular a \qk\ manifold of dimension $4n \geq 8$ has constant scalar curvature, and \qk\ manifolds are often called \emph{positive} or \emph{negative} depending on the sign of their scalar curvature.  A \qk\ metric is equivalent to a reduction of structure group to the maximal compact subgroup $\grp{Sp}{n} \grpdot \grp{Sp}{1}$ of $P^0$, which is one of the groups on Berger's holonomy list \cite{b1955-holonomy}.  Surveys of the state of \qk\ geometry may be found in \cite{am1993-qtnlikestrs,s1982-qkmfds,s1999-qkgeom} and \cite[Chpt.\ 14]{b1987-einstein}.

If the \qtn ic\ structure underlying a \qk\ manifold $(M,g,\cQ)$ is hypercomplex, $(M,g,\cQ)$ is called \emph{(\pseudo)\hk} \cite{h1992-hkmfds,hklr1987-hksupersymm}; in this case $g$ has vanishing Ricci tensor \cite{b1955-affineconn}, with \riem\ holonomy contained in $\grp{Sp}{n}$.  Note that some authors require \qk\ manifolds to have \non zero\ scalar curvature.

\smallskip
\subsection{\Proj\ interpretation} 
\label{ss:qtn-class-proj}

We can give a \proj\ interpretation of \qk\ metrics whose \LC\ connections lie in the same \qtn ic\ class by emulating the theory of Section \ref{s:proj-class} and Subsection \ref{ss:cproj-class-maineqn}.  Let $(M,g,\cQ)$ be a $4n$-dimensional (\pseudo)\qk\ manifold with \LC\ connection $\D$.

\begin{defn} \thlabel{defn:qtn-class-qgeodesic} A smooth curve $\gamma\subset M$ is called a \emph{q-geodesic} of $\D$ if $\D_X X \in \linspan{X,J_aX}{a=1,2,3}$ for every vector field $X$ tangent to $\gamma$.  \Qk\ metrics $g,\b{g}$ are \emph{\qproj ly\ equivalent} if they have the same (unparameterised) q-geodesics. \end{defn}

If $n=1$ this definition is vacuous, so we shall assume that $n > 1$.  Following Fujimura \cite{f1976-Qconns}, the following characterisation of q-projectivity is available.

\begin{lem} \thlabel{lem:qtn-class-qproj} Two \qk\ metrics are \qproj ly\ equivalent if and only if their \LC\ connections are \qtn ically\ equivalent. \noproof \end{lem}

Thus \qproj\ equivalence is determined by the underlying \qtn ic\ geometry.  Denote the \LC\ connections of \qproj ly\ equivalent \qk\ metrics $g,\b{g}$ by $\D,\b{\D}$.  Then $\b{\D} = \D + \algbrac{}{\alpha}$ for some $\alpha\in\s{1}{}$, and taking a trace in \eqref{eq:qtn-class-algbrac} yields
\vspace{0.25em}
\begin{equation*} 
  \alpha = \tfrac{1}{4n+4} \,\d\! \left( \log\frac{\det\b{g}}{\det g} \right)
  \colvectpunct[-0.6em]{,}
\vspace{0.25em}
\end{equation*}
so that $\alpha$ is an exact $1$-form as in \eqref{eq:proj-class-exact1form} and \eqref{eq:cproj-class-exact1form}.  We consider the endomorphism
\vspace{0.2em}
\begin{equation} \label{eq:qtn-class-endoA}
  A(g,\b{g}) \defeq \left( \frac{\det\b{g}}{\det g} \right)^{\!\!1/(4n+4)} \b{\sharp}\circ\flat
\vspace{0.2em}
\end{equation}
analogous to \eqref{eq:proj-class-endoA} and \eqref{eq:cproj-class-endoA}.  When $g,\b{g}$ are clear from the context, we will write $A \defeq A(g,\b{g})$.  Evidently $A(g,\b{g})$ is invertible with inverse $A(\b{g},g)$, and is \self adjoint\ \wrt\ both $g$ and $\b{g}$; moreover since $g,\b{g}$ are Q-hermitian, $A$ is \qtn-linear.  The second \qk\ metric $\b{g}$ may be recovered from $(g,A)$ as
\begin{equation} \label{eq:qtn-class-gbar}
  \b{g} = (\det A)^{1/4}\, g(A^{-1}\bdot, \bdot).
\end{equation}
Calculating as in \thref{prop:proj-class-maineqn,prop:cproj-class-maineqn}, we find that $A$ satisfies the anticipated first-order linear differential equation.

\begin{prop} \thlabel{prop:qtn-class-maineqn} Let $g,\b{g}$ be \qk\ metrics with \LC\ connections $\D,\b{\D}$ respectively.  Then $g,\b{g}$ are \qproj ly\ equivalent if and only if $A \defeq A(g,\b{g})$ defined by \eqref{eq:qtn-class-endoA} satisfies the first-order linear differential equation
\begin{equation} \label{eq:qtn-class-maineqn}
  g((\D_X A)\bdot, \bdot) = X^{\flat}\symm\mu + \qsum J_aX^{\flat}\symm J_a\mu
\end{equation}
for some $\mu\in\s{1}{}$ and all $X\in\s{0}{TM}$.  In this case $\b{\D} = \D + \algbrac{}{\alpha}$, where $\alpha\in\s{1}{}$ satisfies $\mu = -\alpha(A\bdot) = \tfrac{1}{4}\d(\tr A)$. \noproof \end{prop}

\section{Description as a parabolic geometry} 
\label{s:qtn-para}

\Qtn ic\ \proj\ space $\HP$ is defined as the set of \qtn ic\ lines through the origin in $\bH^{n+1}$, or equivalently as the base of the \qtn ic\ Hopf fibration of the $(4n+3)$-sphere $\Sph[4n+3]$ by its natural $\Sph[3]$-action \cite{n2014-hopf}; other characterisations are provided in \cite{am1994-gradqtnvfs,m1976-qtntanno}.   The standard Q-hermitian metric on $\bH^{n+1}$ descends to a \qk\ metric $g_{\mr{FS}}$ on $\HP$, called the \emph{Fubini--Study} metric.  The embedded \qtn ic\ \proj\ lines $\HP[1] \injto \HP$ are the q-geodesics \wrt\ the class of \qtn ic\ connections determined by the \LC\ $\D^{\mr{FS}}$ of $g_{\mr{FS}}$, so are totally geodesic by a result of Gray \cite{g1965-minvars,g1969-holSpSp}.  \Wrt\ an affine chart $\bH^n \injto \HP$ these q-geodesics lie within a \qtn ic\ line in $\bH^n$, so that from the point of view of q-geodesics $\HP$ appears as the natural compactification of $\HP$.  Therefore $\HP$ is a good candidate for the flat model of almost \qtn ic\ geometry; indeed, Salamon defines \cite{s1986-qkdg} an almost \qtn ic\ manifold as one locally modelled on $\HP$.

As we shall explain in Subsection \ref{ss:qtn-para-flat}, the complexification of the flat model $\HP$ may be viewed as the \grassmannian\ of complex $2$-planes.  A careful treatment of parabolic geometries modelled on \grassmannian s\ may be found in \cite{be1991-paraconf,gs1999-qtntwistor}; see also \cite{hsss2012-prolconns2}.  The first of these references in particularly amenable to our discussion of representations in Subsection \ref{ss:qtn-para-repns}.  Throughout this section we shall assume that $n>1$, so that we consider almost \qtn ic\ manifolds $(M,\cQ)$ of dimension $4n \geq 8$.

\vspace{0.1em}
\subsection{The flat model $\HP$} 
\label{ss:qtn-para-flat}

Recalling Salamon's description \cite{s1986-qkdg}, \qtn ic\ \proj\ space is a good candidate for the flat model of almost \qtn ic\ geometry.  As a generalised flag manifold $\HP$ may be identified with $G/P$ for the real Lie groups%
\footnote{\Non commutativity\ of $\bH$ means that the centre of $\grp{GL}{n+1,\bH}$ consists of all \emph{real} multiples of the identity matrix.}
\begin{equation*}
  G \defeq \grp{PGL}{n+1,\bH} \defeq \grp{GL}{n+1,\bH} / \linspan{\id}{},
\end{equation*}
and $P\leq G$ the (projection to $G$ of) the $\grp{GL}{n+1,\bH}$-stabiliser of a given quaternionic line $\linspan[\bH]{v_0}{}$ in $\bH^{n+1}$.

The Lie algebra of $G$ is the real form $\fg \defeq \alg{sl}{n+1,\bH}$ of $\alg{sl}{2n+2,\bC}$, which consists of all $(n+1)\by(n+1)$ \qtn ic\ matrices with vanishing real trace.  To understand the parabolic subalgebra $\fp \leq \fg$, we identify $\bH^{n+1}$ with $\bC^{2n+2}$ together with a conjugate-linear map $\cpx{j} : \bC^{2n+2} \to \bC^{2n+2}$ which satisfies $\cpx{j}^2 = -\id$ and anti-commutes with multiplication by $\cpx{i} \in \bC$.  In this picture the \qtn ic\ line $\linspan[\bH]{v}{} \leq \bH^{n+1}$ is identified with the complex $2$-plane $\linspan[\bC]{v,\cpx{j}v}{}$, which may in turn be identified with the span of the $2$-vector $v\wedge\cpx{j}v \in \Wedge[\bC]{2}\bC^{2n+2}$.  It follows that the Satake diagram of $\fp$ is given by crossing the penultimate node:
\begin{equation} \label{eq:qtn-para-satake}
  \fp = \dynkinSLHp{2}{0}{3}{qtn-para-flat-p}
  \leq \dynkinSLH{2}{0}{3}{qtn-para-flat-g} = \fg.
\end{equation}
Therefore almost \qtn ic\ geometry is an abelian parabolic geometry, with Killing polar $\fp^{\perp} \isom \bH^{n*} \defeq \Hom[\bR]{\bH^n}{\bR}$ and reductive Levi factor $\fp^0 \defeq \fp/\fp^{\perp}$.  An algebraic Weyl structure for the $\fp^{\perp}$-filtration $\fg \supset \fp \supset \fp^{\perp} \supset 0$ of $\fg$ is evidently equivalent to a choice of \qtn ic\ subspace of $\bH^{n+1}$ complementary to $\linspan[\bH]{v_0}{}$, thus yielding a decomposition
\begin{equation} \label{eq:qtn-para-sldecomp}
  \alg{sl}{n+1,\bH} = \Setof{
  \begin{pmatrix}
    q & \alpha \\
    X      & A
  \end{pmatrix}}{ \begin{gathered} X\in\bH^n,~ \alpha\in\bH^{n*},~ A\in\alg{gl}{n,\bH}, \\
                                   ~ q\in\bH \text{\,~s.t.\,} \Re{q} = -\tr{A} \end{gathered} }
  \colvectpunct[-0.9em]{.}
\end{equation}
In the induced decomposition of $\grp{PGL}{n+1,\bH}$, the Levi subgroup $P^0 \defeq P /\exp{\fp}^{\perp}$ consists of (equivalence classes of) block-diagonal matrices in $G$; clearly the adjoint action of $P^0$ on $\fg/\fp \isom \bH^n$ induces an isomorphism $\fp^0 \isom \alg{gl}{n,\bH} \dsum \alg{sp}{1}$.  Then $\fp^{\perp} \isom \bH^{n*}$ consists of matrices with only the $\alpha$-block, $\fp^0$ of block-diagonal matrices, and $\fp \isom (\alg{gl}{n,\bH}\dsum\alg{sp}{1}) \ltimes \bH^{n*}$ of block upper-triangular matrices.  There is a corresponding decomposition on the group level.

Choose an algebraic Weyl structure, so that $\fg \isom \bH^n \dsum (\alg{gl}{n,\bH} \dsum \alg{sp}{1}) \dsum \bH^{n*}$.  The decomposition \eqref{eq:qtn-para-sldecomp} allows us to calculate the Lie bracket between elements in the three summands: if $X,Y\in\bH^n$, $A,B\in\alg{gl}{n,\bH}\dsum \alg{sp}{1}$ and $\alpha,\beta\in\bH^{n*}$ then
\begin{equation*} \begin{gathered}
  \liebrac{X}{Y} = 0 = \liebrac{\alpha}{\beta}, \quad
  \liebrac{A}{B} = AB - BA, \\
  \liebrac{A}{X} = AX
  \quad\text{and}\quad
  \liebrac{A}{\alpha} = -\alpha \circ A,
\end{gathered} \end{equation*}
so that $\bH^n, \bH^{n*}$ form abelian subalgebras on which $\fp^0$ acts naturally.  The bracket $\bH^n \times \bH^{n*} \to \fp^0$ may be computed by adapting Hrinda's calculation in \cproj\ geometry \cite{h2009-cproj}; this calculation is surely known in the literature, but the author could not find a reference.

\begin{lem} \thlabel{lem:qtn-class-liebrac} The Lie bracket $\liebrac{ \liebrac{X}{\alpha} }{ Y }$ between $X,Y\in\bH^n$ and $\alpha\in\bH^{n*}$ satisfies
\begin{equation*}
  \liebrac{ \liebrac{X}{\alpha} }{ Y }
    = \alpha(X)Y + \alpha(Y)X - \qsum\left[ \alpha(J_aX)J_aY + \alpha(J_aY)J_aX \right],
\end{equation*}
where $\{J_a\}_{a=1}^3$ is a \qtn ic\ structure on $\bH^n$. \end{lem}

\begin{proof} We must first identify $\fg=\alg{sl}{n+1,\bH}$ with a subalgebra of $\alg{gl}{4n+4,\bR}$ as follows.  Given $A=(a_{ij})\in\fg$, write $\bH\ni a_{ij} = a_{ij}^0\cpx{1} + a_{ij}^1\cpx{i}_1 + a_{ij}^2\cpx{i}_2 + a_{ij}^3\cpx{i}_3$ for $a_{ij}^k\in\bR$.  Then $A$ may be identified with the $(4n\by 4n)$-matrix given by replacing each entry $a_{ij}\in\bH$ with the $(4\by 4)$ block
\begin{equation*}
\newcommand{\ph}{\phantom{-}}  
\newcommand{\svdots}{               \raisebox{3pt}{$\scalebox{0.65}{$\vdots$}$}}
\newcommand{\sddots}{\hspace{-0.2ex}\raisebox{3pt}{$\scalebox{0.58}{$\ddots$}$}}
  \left( \begin{smallmatrix}
    \ph a_{ij}^0 & \ph a_{ij}^1 & \ph a_{ij}^2 & \ph a_{ij}^3 \\
       -a_{ij}^1 & \ph a_{ij}^0 &    -a_{ij}^3 & \ph a_{ij}^2 \\
       -a_{ij}^2 & \ph a_{ij}^3 & \ph a_{ij}^0 &    -a_{ij}^1 \\
       -a_{ij}^3 &    -a_{ij}^2 & \ph a_{ij}^1 & \ph a_{ij}^0
  \end{smallmatrix} \right)
  \colvectpunct[-1.5em]{.}
\end{equation*}
In particular, the action of the real unit $\cpx{1}$ is given by the $(4n\by 4n)$ identity matrix $I$, while the \qtn ic\ units $\{\cpx{i_a}\}_{a=1}^3$ act via the real endomorphisms $\{J_a\}_{a=1}^3$ given by a direct sum of $n$ copies of the $(4\by 4)$ blocks
\begin{equation} \label{eq:qtn-para-brac-1}
\newcommand{\ph}{\phantom{-}}
  J_1 =
  \left(\!\! \begin{smallmatrix}
    \ph 0 & \ph 1 & \ph 0 & \ph 0 \\
       -1 & \ph 0 & \ph 0 & \ph 0 \\
    \ph 0 & \ph 0 & \ph 0 &    -1 \\
    \ph 0 & \ph 0 & \ph 1 & \ph 0
  \end{smallmatrix} \,\right)
  \colvectpunct[-0.7em]{,} \quad
  J_2 =
  \left(\!\! \begin{smallmatrix}
    \ph 0 & \ph 0 & \ph 1 & \ph 0 \\
    \ph 0 & \ph 0 & \ph 0 & \ph 1 \\
       -1 & \ph 0 & \ph 0 & \ph 0 \\
    \ph 0 &    -1 & \ph 0 & \ph 0
  \end{smallmatrix} \,\right)
  \quad\text{and}\quad
  J_3 =
  \left(\! \begin{smallmatrix}
    \ph 0 & \ph 0 & \ph 0 & \ph 1 \\
    \ph 0 & \ph 0 &    -1 & \ph 0 \\
    \ph 0 & \ph 1 & \ph 0 & \ph 0 \\
       -1 & \ph 0 & \ph 0 & \ph 0
  \end{smallmatrix} \,\right)
  \colvectpunct[-0.7em]{.}
\end{equation}
Choose $X,Y \in \bH^n$ and $\alpha \in \bH^{n*}$, and write $X_i,Y_i \in \bH$ and $\alpha_i \in \bH^*$ for their components.  Using the decomposition \eqref{eq:qtn-para-sldecomp}, the bracket $\liebrac{ \liebrac{X}{\alpha} }{ Y }$ is given by
\bgroup  \newcommand{\alx}[1]{\alpha_1{#1}_1 \,+\; \cdots \;+\, \alpha_n{#1}_n}  
         \newcommand{\svdots}{               \raisebox{3pt}{$\scalebox{0.65}{$\vdots$}$}}
         \newcommand{\sddots}{\hspace{-0.2ex}\raisebox{3pt}{$\scalebox{0.58}{$\ddots$}$}}
\begin{align} \label{eq:qtn-para-brac-2}
  \liebrac{ \liebrac{X}{\alpha} }{ Y }
  &= \left[ \left[
    \left( \begin{smallmatrix}
      0      & 0      & \cdots & 0      \\
      X_1    & 0      & \cdots & 0      \\
      \svdots & \svdots & \sddots & \svdots \\
      X_n    & 0      & \cdots & 0
    \end{smallmatrix} \right)
    \colvectpunct[-0.9em]{,}
    \left( \begin{smallmatrix}
      0      & \alpha_1 & \cdots & \alpha_n \\
      0      & 0        & \cdots & 0        \\
      \svdots & \svdots   & \sddots & \svdots   \\
      0      & 0        & \cdots & 0
    \end{smallmatrix}\right)
    \right] \colvectpunct[-0.9em]{,}
    \left( \begin{smallmatrix}
    0      & 0      & \cdots & 0      \\
    Y_1    & 0      & \cdots & 0      \\
    \svdots & \svdots & \sddots & \svdots \\
    Y_n    & 0      & \cdots & 0
    \end{smallmatrix} \right)
    \right] \notag \\
  &= \left( \begin{smallmatrix}
                &0     &            & 0      & \cdots & 0 \\
     \alx{X}Y_1 &+     & \alx{Y}X_1 & 0      & \cdots & 0 \\
                &\svdots&            & \svdots & \sddots & \svdots \\
     \alx{X}Y_n &+     & \alx{Y}X_n & 0      & \cdots & 0
  \end{smallmatrix} \right)
  \colvectpunct[-1.0em]{.}
  \notag \\[-2.2em]
\end{align}
\egroup
Identifying $\alpha \in \bH^*$ and $X_i \in \bH$ with the $(4\by 4)$-matrices
\begin{equation*}
\newcommand{\al}{\alpha} \newcommand{\ph}{\phantom{-}}  
  \al_i = 
  \left( \begin{smallmatrix}
    \ph \al_i^0 & \ph \al_i^1 & \ph \al_i^2 & \ph \al_i^3 \\
       -\al_i^1 & \ph \al_i^0 &    -\al_i^3 & \ph \al_i^2 \\
       -\al_i^2 & \ph \al_i^3 & \ph \al_i^0 &    -\al_i^1 \\
       -\al_i^3 &    -\al_i^2 & \ph \al_i^1 & \ph \al_i^0
  \end{smallmatrix} \right)
  \quad\text{and}\quad
  X_i =
  \left( \begin{smallmatrix}
    \ph X_i^0   & \ph X_i^1   & \ph X_i^2   & \ph X_i^3   \\
       -X_i^1   & \ph X_i^0   &    -X_i^3   & \ph X_i^2   \\
       -X_i^2   & \ph X_i^3   & \ph X_i^0   &    -X_i^1   \\
       -X_i^3   &    -X_i^2   & \ph X_i^1   & \ph X_i^0
  \end{smallmatrix} \right)
\end{equation*}
as above, multiplying out the matrices in one of the constant factors $\alpha_iX_i \in \bH$ shows that the factor $(\alpha_1 X_1 + \cdots + \alpha_n X_n)$ is given by summing over terms of the form
\bgroup  \newcommand{\alx}[2]{\alpha_i^{#1}X_i^{#2}}  
\vspace{-0.15em}
\begin{equation} \label{eq:qtn-para-brac-3} \begin{aligned}
  \alpha_i X_i
  &= (\alx{0}{0} - \alx{1}{1} - \alx{2}{2} - \alx{3}{3})I \\
  &\qquad + (\alx{0}{1} + \alx{1}{0} + \alx{2}{3} - \alx{3}{2})J_1 \\
  &\qquad + (\alx{0}{2} - \alx{1}{3} + \alx{2}{0} + \alx{3}{1})J_2 \\
  &\qquad + (\alx{0}{3} + \alx{1}{2} - \alx{2}{1} + \alx{3}{0})J_3,
\end{aligned}
\vspace{-0.15em}
\end{equation}
\egroup
where we identify $I$ with the $(4\by 4)$ identity matrix and $\{J_a\}_{a=1}^3$ with the matrices from \eqref{eq:qtn-para-brac-1}.  Equating the entries on the leading diagonal with the real part, we note that
\bgroup  \newcommand{\al}{\alpha} \newcommand{\ph}{\phantom{-}}  
         \newcommand{\vp}{\vphantom{^0_0}}
\begin{align*}
  &\alpha_i^0 X_i^1 + \alpha_i^1 X_i^0 + \alpha_i^2 X_i^3 - \alpha_i^3 X_i^2 \\
  &\quad \begin{aligned}
    &= -\Re{ \!
      \left( \begin{smallmatrix} 
          \ph \al_i^0 & \ph \al_i^1 & \ph \al_i^2 & \ph \al_i^3 \\
         -\al_i^1 & \ph \al_i^0 &    -\al_i^3 & \ph \al_i^2 \\
         -\al_i^2 & \ph \al_i^3 & \ph \al_i^0 &    -\al_i^1 \\
         -\al_i^3 &    -\al_i^2 & \ph \al_i^1 & \ph \al_i^0
      \end{smallmatrix} \right) \hspace{-1ex}
      \left(\!\! \begin{smallmatrix} 
          \ph 0\vp & \ph 1\vp & \ph 0\vp & \ph 0\vp \\
             -1\vp & \ph 0\vp & \ph 0\vp & \ph 0\vp \\
          \ph 0\vp & \ph 0\vp & \ph 0\vp &    -1\vp \\
          \ph 0\vp & \ph 0\vp & \ph 1\vp & \ph 0\vp
      \end{smallmatrix} \right) \hspace{-1ex}
      \left( \begin{smallmatrix} 
          \ph X_i^0 & \ph X_i^1 & \ph X_i^2 & \ph X_i^3 \\
         -X_i^1 & \ph X_i^0 &    -X_i^3 & \ph X_i^2 \\
         -X_i^2 & \ph X_i^3 & \ph X_i^0 &    -X_i^1 \\
         -X_i^3 &    -X_i^2 & \ph X_i^1 & \ph X_i^0
      \end{smallmatrix} \right) \! } \\
    &= -\Re{ \alpha_i J_1 X_i }
  \end{aligned}
\end{align*}
\egroup  
and similarly for the other terms in \eqref{eq:qtn-para-brac-3}.  Therefore
\vspace{-0.1em}
\begin{equation*}
  \alpha_i X_i = \Re{ (\alpha_i X_i)I - (\alpha_i J_1 X_i)J_1 - (\alpha_i J_2 X_i)J_2 - (\alpha_i J_3 X_i)J_3 };
\vspace{-0.1em}
\end{equation*}
upon recalling that $\alpha(X) \in \bR$ by the definition $\bH^{n*} \defeq \Hom{\bH^n}{\bR}$, this yields
\vspace{-0.1em}
\begin{equation*}
  (\alpha_1 X_1 + \cdots + \alpha_n X_n)Y_i
    = \alpha(X)Y_i - \alpha(J_1X)J_1Y_i - \alpha(J_2X)J_2Y_i - \alpha(J_3X)J_3Y_i
\vspace{-0.1em}
\end{equation*}
for all $i=1,\ldots,n$.  The desired formula now follows from \eqref{eq:qtn-para-brac-2}. \end{proof}

In particular \thref{lem:qtn-class-liebrac} implies that, up to normalisation conventions, change of Weyl structure corresponds precisely to change of connection in the \qtn ic\ class.

Note that $\HP$ has a \proj\ embedding defined as follows.  As above we view $\bH^{n+1}$ as $\bC^{2n+2}$ with a conjugate-linear map $\cpx{j}$, so that $\HP$ is identified with the set of $\cpx{j}$-invariant $2$-planes in $\bC^{2n+2}$, \ie\ with the subspace
\vspace{-0.1em}
\begin{equation*}
  \HP = \Setof{ W\in\Grass{2}{\bC^{2n+2}} }{ \cpx{j}W = W }
\vspace{-0.1em}
\end{equation*}
of the \grassmannian\ $\Grass{2}{\bC^{2n+2}}$.  Note that since $\cpx{j}^2 = -\id$, the even wedge powers $\Wedge[\bC]{2k}\bC^{2n+2}$ admit real structures given by $\Wedge{2k}\cpx{j} : \Wedge[\bC]{2k}\bC^{2n+2} \to \Wedge[\bC]{2k}\bC^{2n+2}$, so we may form the underlying real representations $\realrepn{( \Wedge[\bC]{2k} \bC^{2n+2} )}$.  Identifying $[v]\in\HP$ with the $2$-plane $\linspan[\bC]{v,\cpx{j}v}{}$, the map $[v]\mapsto[v\wedge\cpx{j}v]$ defines a \proj\ embedding $\HP \injto \pr{\Wedge[\bC]{2}\bC^{2n+2}}$.  Since $(\Wedge{2}\cpx{j})(v\wedge\cpx{j}v) = -\cpx{j}v\wedge v = v\wedge\cpx{j}v$, this embedding takes values in the underlying real representation $\bW \defeq \realrepn{( \Wedge[\bC]{2}\bC^{2n+2} )}$, thus yielding an embedding $\HP \injto \pr{\bW}$ which generalises the Pl{\"u}cker embedding.

A dimension count shows that $\HP[n]$ has \non zero\ codimension in $\pr{\bW}$.  Indeed, Kostant's \thref{thm:lie-para-kostant} tells us that $\HP$ is the intersection of quadrics in $\pr{\bW}$ with defining equations given by projecting away from the Cartan square in $\Symm{2}\bW$.  Explicitly,
\begin{align*}
  \Symm{2}\bW^*
     = \Symm{2} \bigg( \dynkinSLH[0,0,0,1,0]{2}{0}{3}{qtn-para-flat-Wd} \bigg)
  &\isom
    \dynkinname{ \dynkinSLH[0,0,0,2,0]{2}{0}{3}{qtn-para-flat-OWd} }
               { \Cartan{2}\bW^* }
    \dsum
    \dynkinname{ \dynkinSLH[0,0,0,1,0,0,0]{2}{0}{5}{qtn-para-flat-W4} }
               { \realrepn{( \Wedge[\bC]{4}\bC^{2n+2*} )} },
\end{align*}
so that
$
  \bU^* \defeq 
     \realrepn{( \Wedge[\bC]{4}\bC^{2n+2*} )}
$
is the space of homogeneous quadratic polynomials which cuts out $\HP$ as an intersection of quadrics.  As in the real and complex cases, $\bU^*$ is an irreducible $\fg$-representation.

Finally, since $\fg=\alg{sl}{n+1,\bH}$ is not one of the problematic Lie algebras for the equivalence of categories provided by \thref{thm:para-calc-equiv}, the following is immediate.

\begin{thm} \thlabel{thm:qtn-para-equiv} On any manifold $M$ of dimension $4n \geq 8$, there is an equivalence of categories between almost \qtn ic\ structures and regular normal parabolic geometries of type $\HP$.  The flat model is $\HP$ with its canonical \qtn ic\ structure determined by the Fubini--Study metric $g_{\mr{FS}}$. \noproof \end{thm}

Normality implies that $\Tor{}$ is the intrinsic torsion of the almost \qtn ic\ structure generated by the Weyl connections, which we describe more carefully in Subsection \ref{ss:qtn-para-harm}.  In dimension four we obtain a similar equivalence of categories between regular normal parabolic geometries of type $\HP[1] \isom \Sph[4]$ and \self dual\ conformal structures.

\medskip
\subsection{Representations of $\alg{sl}{n+1,\bH}$} 
\label{ss:qtn-para-repns}

In order to describe $\fg$- and $\fp$-representations, it will be convenient to consider the complexification $\cpxrepn{\fg} \defeq \alg{sl}{2n+2,\bC}$ of $\fg$.  Then completely reducible $\cpxrepn{\fp}$-representations are trivial extensions of representations of the complexified reductive Levi factor $\cpxrepn{\fp^0} \isom \alg{gl}{2n,\bC} \dsum \alg{sl}{2,\bC}$, so are external tensor products of representations of the two factors.

The fundamental representations of $\cpxrepn{\fg}$ are the exterior powers of $\bC^{2n+2}$.  Identifying $\bH^{n+1}$ with $\bC^{2n+2}$ as above, the odd exterior and symmetric powers of $\bC^{2n+2}$ evidently admit \qtn ic\ structures, while the even powers admit real structures.  Thus if a $\cpxrepn{\fg}$-representation has even coefficients over the odd nodes of the Dynkin diagram, it is the complexification of a $\fg$-representation.

Let us now describe some important $\fg$- and $\fp$-representations, as well as their associated bundles.  The complexified isotropy representation $\cpxrepn{(\fg/\fp)}$ and its dual $\cpxrepn{(\fg/\fp)}^* \isom \cpxrepn{\fp^{\perp}}$ decompose as external tensor products according to%
\footnote{Note that we denote representations of the complexifications $\cpxrepn{\fg}, \cpxrepn{\fp}$ on their respective Dynkin diagrams, rather than on the Satake diagrams \eqref{eq:qtn-para-satake} of $\fg, \fp$.}
\begin{align*}
  \cpxrepn{(\fg/\fp)} &= \dynkinApp[1,0,0,0,0,1]{2}{0}{4}{qtn-para-repn-fg} = E \etens H \\
  \quad\text{and}\quad
  \cpxrepn{\fp^{\perp}} &= \dynkinApp[0,0,0,1,-2,1]{2}{0}{4}{qtn-para-repn-pperp} = E^* \etens H^*,
\end{align*}
where
\vspace{0.1em}
\begin{equation} \label{eq:qtn-para-EH}
\arraycolsep=0.1em
\begin{array}{rlrl}  
  E &\defeq \dynkinApp[1,0,0,0,0,0]{2}{0}{4}{qtn-para-repn-E}, &
  \qquad H &\defeq \dynkinApp[0,0,0,0,0,1]{2}{0}{4}{qtn-para-repn-H} \\
  \text{and}\quad
  E^* \!\!&\hspace{0.65ex}= \dynkinApp[0,0,0,1,-1,0]{2}{0}{4}{qtn-para-repn-Ed}, &
  \qquad H^* \!\!&\hspace{0.65ex}= \dynkinApp[0,0,0,0,-1,1]{2}{0}{4}{qtn-para-repn-Hd}
\end{array}
\end{equation}
are the natural representations of $\cpxrepn{\alg{gl}{n,\bH}} = \alg{gl}{2n,\bC}$ on $E \defeq \bC^{2n} \isom \bH^n$ and of $\cpxrepn{\alg{sp}{1}} = \alg{sl}{2,\bC}$ on $H \defeq \bC^2 \isom \bH$,%
\footnote{The \non commutativity\ of the \qtn s\ makes it necessary to view $E$ as a \emph{left} $\alg{gl}{2n,\bC}$-module, whereas we view $H$ as a \emph{right} $\alg{sl}{2,\bC}$-module.}
each of which carry a \qtn ic\ structure.  By the Cartan condition, the complexified tangent and cotangent bundles decompose as external tensor products $\cpxbdl{TM} \isom \cE \etens \cH$ and $\cpxbdl{T^*M} \isom \cE^* \etens \cH^*$, where $\cE,\cH$ are the bundles associated to $E,H$.%
  \footnote{Thus almost \qtn ic\ geometry is an \emph{almost \grassmannian\ geometry}; see \cite{gs1999-qtntwistor} and \cite[\S 4.1.3]{cs2009-parabolic1}.}
Note that $\cE,\cH$ need not be defined globally, for reasons discussed below.  Decompositions of $\Wedge[\bC]{k} \cpxbdl{TM}$, $\Symm[\bC]{k} \cpxbdl{TM}$ and their duals may then be given by representation-theoretic means, for example using Schur functors \cite{fh1991-repntheory}; decompositions for small values of $k$ may be obtained by more direct methods.  The case $k=2$ shall be important for us in the sequel, for which
\vspace{0.2em}
\begin{equation} \label{eq:qtn-para-wedge2} \begin{split}
  \Wedge[\bC]{2}( E\etens H ) &\isom
    \dynkinname{ \dynkinApp[0,1,0,0,0,0,2]{3}{0}{4}{qtn-para-repn-wedge1} }
               { \Wedge[\bC]{2}E \etens \Symm[\bC]{2}H } \dsum
    \dynkinname{ \dynkinApp[2,0,0,0,0,1,0]{3}{0}{4}{qtn-para-repn-wedge2} }
               { \Symm[\bC]{2}E \etens \Wedge[\bC]{2}H } \\[0.5em]
  \quad\text{and}\quad
  \Symm[\bC]{2}( E \etens H ) &\isom
    \dynkinname{ \dynkinApp[0,1,0,0,0,1,0]{3}{0}{4}{qtn-para-repn-symm1} }
               { \Wedge[\bC]{2}E \etens \Wedge[\bC]{2}H } \dsum
    \dynkinname{ \dynkinApp[2,0,0,0,0,0,2]{3}{0}{4}{qtn-para-repn-symm2} }
               { \Symm[\bC]{2}E \etens \Symm[\bC]{2}H }.
\end{split}
\vspace{0.2em}
\end{equation}
Each summand in \eqref{eq:qtn-para-wedge2} admits a real structure, so is the complexification of an underlying real representation.  In particular $\alg{sl}{2,\bC}$ acts trivially on $\Wedge[\bC]{2}E \etens \Wedge[\bC]{2}H$, so its associated bundle may be identified with the subbundle $\Symm[+]{2}TM$ of Q-invariant bilinear forms on $T^*M$.  The complementary subbundle is associated to $\realrepn{( \Symm[\bC]{2}E \etens \Symm[\bC]{2}H )}$ and shall be denoted by $\Symm[-]{2}TM$, although note that $\Symm[-]{2}TM$ does \emph{not} consist of Q-anti-invariant bilinear forms.  Salamon \cite[Prop.\ 9.2]{s1989-riemholo} provides similar decompositions of the bundles of $3$- and $4$-forms, while Swann \cite{s1991-hkqk,s1991-hkqkthesis} refines these decompositions and gives a decomposition for the bundle of $5$-forms; summaries may also be found in \cite{ms2005-QKkilling,s1999-qkgeom}.

The line bundle $\cL \defeq (\Wedge{4n}TM)^{1/(2n+2)}$ from \thref{lem:qtn-class-connsonL} is associated to the $(2n+2)$nd root $L$ of $\Wedge{4n}(\fg/\fp)$, which has highest weight
\vspace{0.15em}
\begin{equation} \label{eq:qtn-para-repnL}
  L \defeq \dynkinSLHp[0,0,0,0,1,0]{2}{0}{4}{qtn-para-repn-L}.
\vspace{0.15em}
\end{equation}
Since $E,H$ have complex dimension $2n,2$ respectively, both $\Wedge[\bC]{2}H$ and $\Wedge[\bC]{2n}E$ are complex line bundles isomorphic to the complexification of $L$; \cf\ the spinor bundles from \cite{be1991-paraconf,gs1999-qtntwistor}.  As in real \proj\ and \cproj\ cases, $L$ is the zeroth homology $\liehom{0}{\bW^*}$, where $\bW = \realrepn{( \Wedge[\bC]{2}\bC^{2n+2} )}$ is the $\fg$-representation from Subsection \ref{ss:qtn-para-flat} admitting a \proj\ embedding $\HP \injto \pr{\bW}$.

Following the discussion of Subsection \ref{ss:qtn-class-qtnic}, the bundle $\cQ$ of almost complex structures is pointwise isomorphic to $\alg{sp}{1}$ and has a natural action of $\grp{SO}{3}$.  Extending the adjoint action of $\alg{so}{3} \isom \alg{sp}{1} \isom \realrepn{(\Symm[\bC]{2}H)}$ trivially to $\fp^0 \isom \alg{gl}{n,\bH} \dsum \alg{sp}{1}$, we may view $\cQ$ as the bundle associated to the adjoint representation
\vspace{0.15em}
\begin{equation*} 
  Q = \dynkinSLHp[0,0,0,0,0,2]{2}{0}{4}{qtn-para-repn-Q} = \realrepn{( \Symm[\bC]{2}H )}
\vspace{0.15em}
\end{equation*}
of the simple factor $\alg{sp}{1}$.

As before, our choice of group $G = \grp{PGL}{n+1,\bH}$ with Lie algebra $\fg = {\alg{sl}{n+1,\bH}}$ means that not all $\fg$-representations integrate to global $G$-representations, which causes problems for forming some tractor bundles.  Since the centre of $\grp{GL}{n+1,\bH}$ consists of all real multiples of the identity, the topology of $\grp{PGL}{n+1,\bH}$ is much like the topology of $\grp{PGL}{n+1,\bR}$ discussed in Subsection \ref{ss:proj-para-repns}.  Thus an irreducible $\fg$-representation with highest weight $\lambda$ integrates to $G$ if and only if the central $\bZ_2$ of $\tilde{G} \defeq \grp{SL}{n+1,\bH}$ acts trivially, which is the case if and only if the coefficient sum of $\lambda$ is even \cite{fh1991-repntheory}.  We then form the extended Cartan bundle $F^{\tilde{P}} \defeq \assocbdl{F^P}{P}{\tilde{P}}$, where $\tilde{P}\leq\tilde{G}$ is the parabolic stabiliser of a given \qtn ic\ line in $\bH^{n+1}$ (previously denoted $\linspan[\bH]{v_0}{}$), form tractor bundles \wrt\ $F^{\tilde{P}}$, and take the local quotient by $\bZ_2$ if necessary.

Obvious examples of representations which do not integrate to $G$ are $E,H$ defined by \eqref{eq:qtn-para-EH}, as well as the standard representation of $\fg$ on
\vspace{-0.1em}
\begin{equation*}
  \bT \defeq \dynkinA[1,0,0,0,0]{2}{0}{3}{} = \bH^{n+1}.
\vspace{-0.1em}
\end{equation*}
It follows from the decomposition \eqref{eq:qtn-para-sldecomp} that an element $\inlinematrix{ q & \alpha \\ 0 & A } \in \tilde{P}$ acts on the \qtn ic\ line $\bT_0 \defeq \linspan[\bH]{v_0}{} \leq \bT$ by multiplication with $q$.  Therefore
\vspace{-0.1em}
\begin{equation*}
  \bT_0 = \dynkinApp[0,0,0,-1,1]{2}{0}{3}{} = H^*
  \quad\text{and}\quad
  \bT/\bT_0 = \dynkinApp[1,0,0,0,0]{2}{0}{3}{} = E
\vspace{-0.1em}
\end{equation*}
as $\tilde{P}$-representations, which we identify with the socle and top of the $\fp^{\perp}$-filtration $\bT \supset \bT_0 \supset 0$.  Thus $\bT \isom E \dsum H^*$ \wrt\ an algebraic Weyl structure; the associated bundle $\cT$ is called the \emph{standard tractor bundle} of \qtn ic\ geometry.

In terms of structure groups, the extended bundle $F^{\tilde{P}}$ is equivalent to the double cover of the \qtn ic\ frame bundle $F^0$ with structure group $P^0$, giving a bundle $\tilde{F}^0$ with structure group $\tilde{P}^0 \defeq \grp{GL}{n,\bH}\times\grp{Sp}{1}$.  In fact, all $\fg$-representations integrate to $G$ provided that a certain cohomology class vanishes \cite{mr1976-qtncohom,s1986-qkdg,s1991-hkqk,s2010-twisting}.  In particular this class vanishes when $n+1$ is odd, since then $G = \tilde{G}$, which is the origin of Salamon's result stating that a \qk\ manifold of dimension $8n$ is spin \cite{s1982-qkmfds}.

\subsection{Harmonic curvature} 
\label{ss:qtn-para-harm}

We turn now to computing the harmonic curvature of the canonical Cartan connection, which lies in the second Lie algebra homology $\liehom{2}{\fg}$.  The Hasse diagram computing this homology is given in Figure \ref{fig:qtn-para-hasse}, from which we see that the harmonic curvature has two components, one for each summand in the third column.  From top to bottom the complexifications of these $\fp$-representations are
\begin{equation} \label{eq:qtn-para-harm} \begin{gathered}
  (\Wedge[\bC]{2}E^* \tens \Wedge[\bC]{2}H^*) \cartan[\bC] (E \tens H) \cartan[\bC] \Symm[\bC]{2} H \\
  \text{and}\quad
  (\Symm[\bC]{2}E^* \tens \Wedge[\bC]{2}H^*) \cartan[\bC] (E \cartan[\bC] E^*),
\end{gathered} \end{equation}
whose underlying real representations have associated bundles $\Symm[+]{2}T^*M \cartan TM \cartan \cQ$ and $\Wedge[+]{2}T^*M \cartan \realrepn{( \cE\cartan[\bC]\cE^* )}$.

The component in $\Symm[+]{2}T^*M \cartan TM \cartan \cQ$ coincides with the intrinsic torsion $\Tor{}$ of any Weyl connection.  To see this, observe that there is a natural map $\Symm[+]{2}T^*M \tens \cQ \injto \Wedge{2} T^*M$ defined by $g \tens J \mapsto g(J\bdot,\bdot)$, which is clearly injective.  Moreover for each $J\in\s{0}{\cQ}$ there is a unique $J'\in\s{0}{\cQ}$ such that $\{J,J',J\circ J'\}$ is a local \qtn ic\ frame of $\cQ$.  \Wrt\ this frame it is easy to check that $g(J\bdot,\bdot)$ is $J$-anti-invariant, so that this map takes values in the irreducible summand $\Wedge[-]{2}T^*M$.  A dimension count then implies that $\Symm[+]{2}T^*M \tens \cQ \isom \Wedge[-]{2}T^*M$, so that $\Tor{}$ may be identified with a section of $\Wedge[-]{2}T^*M \cartan TM$.  {\v Cap} and Slov{\'a}k prove \cite[Prop.\ 4.1.8]{cs2009-parabolic1} that $\Tor{}$ takes values in the irreducible subbundle given by the intersection of the kernels of the natural contractions $\Symm[+]{2}T^*M \tens TM \tens \cQ \surjto T^*M\tens \cQ$ and $\Symm[+]{2}T^*M \tens TM \tens \cQ \surjto \Symm[+]{2}T^*M \tens TM$.

\vspace{1.3em}
\begin{figure}[!h]
  \begin{equation*}
  \arraycolsep=0.07em
  \raisebox{1.2ex}{$\begin{array}{ccccccc}                 
    \begin{array}{c}
      \dynkinSLHp[1,0,0,0,1]{2}{0}{3}{qtn-para-harm-hasse1}
    \end{array}
    &
    \dynkin{ \DynkinConnector{0.2}{-0.26}{1.3}{-0.26}; }{}
    &
    \begin{array}{c}
      \dynkinSLHp[1,0,0,1,-2,2]{2}{0}{4}{qtn-para-harm-hasse3}
    \end{array}
    &
    \dynkin{ \DynkinConnector{0.2}{-0.1}{1.3}{ 0.7};
             \DynkinConnector{0.2}{-0.5}{1.3}{-1.3};
           }{qtn-para-harm-hasse4}
    &
    \begin{array}{c}
      \dynkinSLHp[1,0,0,1,0,-3,3]{2}{0}{5}{qtn-para-harm-hasse5} \\[0.5em]
      \dynkinSLHp[1,0,0,0,3,-4,0]{2}{0}{5}{qtn-para-harm-hasse6}
    \end{array}
    &
    \dynkin{ \DynkinConnector{0.2}{ 0.7}{1.3}{-0.1};
             \DynkinConnector{0.2}{-1.3}{1.3}{-0.5};
             \DynkinLabel{$\cdots$}{2.2}{-1.1}
           }{qtn-para-harm-hasse7}
  \end{array}$}
  \end{equation*}
  \vspace{0.5em}
  \caption[The Hasse diagram computing $\liehom{}{\alg{sl}{n+1,\bH}}$]
          {The Hasse diagram of the adjoint representation of $\fg = \alg{sl}{n+1,\bH}$, which computes the homology $\liehom{}{\fg}$, drawn here for $n\geq 4$.}
  \vspace{1.3em}
  \label{fig:qtn-para-hasse}
\end{figure}

The component of \eqref{eq:qtn-para-harm} in $\Wedge[+]{2}T^*M \cartan \realrepn{( \cE\cartan[\bC] \cE^* )}$ is the Weyl curvature $\Weyl[\D]{}{}$ of $\D$.  Since the complexification of $\alg{gl}{TM}$ is associated to $(E \etens H) \tens (E^* \etens H^*)$, the description of the Cartan product in \cite{e2005-cartanprod} gives
\vspace{0.3em}
\begin{equation} \label{eq:qtn-para-glEH} \begin{aligned}
  (E \etens H) \tens (E^*\etens H^*)
  &\isom ((E\cartan[\bC]E^*)\etens(H\cartan[\bC]H^*)) \dsum (\bC \etens (H\cartan[\bC]H^*)) \\
  &\hspace{2.5em} \dsum ((E\cartan[\bC]E^*) \etens \bC) \dsum (\bC \etens \bC).
\end{aligned}
\vspace{0.3em}
\end{equation}
The second summand may be identified with $\cpxrepn{Q}$ via a dimension count, while the final two summands respectively produce the complexifications of the \tracefree\ part $\alg{sl}{TM,\cQ}$ and trace part of $\alg{gl}{TM,\cQ}$.  It follows that $\Weyl[\D]{}{}$ is Q-invariant, $\alg{sl}{TM}$-valued, and totally \tracefree.

By \thref{thm:para-bgg-ablcurv}, the curvature tensor $\Curv{}{}$ of a Weyl connection $\D \in\Dspace$ decomposes as $\Curv{}{} = \Weyl[\D]{}{} - \algbracw{\id}{\nRic{}{}}{}$, where $\nRic{}{} \defeq -\quab_M^{-1} \liebdy \Curv{}{}$ is the normalised Ricci tensor of $\D$.  Then $\weyld{\gamma} \Weyl[\D]{}{} = \algbrac{\Tor{}}{\gamma}$ by \itemref{thm:para-bgg-ablcurv}{weyld}, so that $\Weyl[\D]{}{}$ is \qtn ically\ invariant if and only if $\Tor{}$ vanishes, if and only if each Weyl connection is \qtn ic.

The Cotton--York tensor $\CY{}{} \defeq \d^{\D} \nRic{}{}$ of $\D$ splits into components $\CY[\D\pm]{}{}$ according to the decomposition \eqref{eq:qtn-para-wedge2} of $\Wedge{2}T^*M$.  For later use we collect some curvature identities; their proofs are similar to the corresponding identities in Subsections \ref{ss:proj-para-harm} and \ref{ss:cproj-para-harm}.

\begin{prop} \thlabel{prop:qtn-para-calc} Let $\D\in\Dspace$ be a Weyl connection.  Then:
\begin{enumerate}
  \item \label{prop:qtn-para-calc-Wid}
  $\Weyl[\D]{}{}$ satisfies the Bianchi identity
  \begin{equation} \label{eq:qtn-para-Wbianchi}
    \Weyl[\D]{X,Y}{Z} + \Weyl[\D]{Y,Z}{X} + \Weyl[\D]{Z,X}{Y} = (\d^{\D}\Tor{})_{X,Y,Z}.
  \end{equation}
  
  \item \label{prop:qtn-para-calc-Wprop}
  $\Weyl[\D]{}{}$ is a Q-invariant and valued in $\alg{sl}{TM,\cQ}$; thus $\Weyl[\D]{}{}$ is totally \tracefree, and $\Weyl[\D]{}{J} = 0$ and $\Weyl[\D]{JX,JY}{} = \Weyl[\D]{X,Y}{}$ for all unit norm $J\in\s{0}{\cQ}$ and all $X,Y\in\s{0}{TM}$.  Moreover $\Weyl[\D]{}{}$ is \qtn ically\ invariant if and only if the intrinsic torsion $\Tor{}$ of $\Dspace$ vanishes.
  
  \item \label{prop:qtn-para-calc-CYid}
  $\CY{}{}$ satisfies the Bianchi identities
    \begin{equation} \label{eq:qtn-para-CYbianchi} \begin{aligned}
    \CY{X,Y}{Z} + \CY{Y,Z}{X} + \CY{Z,X}{Y} &= 0 \\
    \text{and}\quad
    \CY{X,Y}{JZ} + \CY{Y,Z}{JX} + \CY{Z,X}{JY} &= 0
  \end{aligned}
  \end{equation}
  for each $J\in\s{0}{\cQ}$.
  
  \item \label{prop:qtn-para-calc-WCY}
  We have $\ve^i(\D_{e_i}\Weyl[\D]{X,Y}{}) = (2n-1)\CY[\D+]{X,Y}{}$ \wrt\ any local frame $\{e_i\}_i$ of $TM$ with dual coframe $\{\ve^i\}_i$.
  
  \item \label{prop:qtn-para-calc-nric}
  $\nRic{}{}$ is related to the Ricci curvature $\liebdy \Curv{}{}$ of $\D$ by
  \begin{equation*}
    \nRic{}{} = -\tfrac{1}{2n}(\sym{\liebdy \Curv{}{}}) - \tfrac{1}{2n+2}(\alt{\liebdy \Curv{}{}})
      + \tfrac{1}{n^2+2n}(\sym{\liebdy \Curv{}{}})^{+}.
  \end{equation*}
  In particular if $\nRic{}{}$ is symmetric and Q-invariant, then $\nRic{}{} = -\tfrac{1}{2n+4} \liebdy \Curv{}{}$. \noproof
\end{enumerate}
\end{prop}

We have already proved the claims in \itemref{prop:qtn-para-calc}{Wprop}; more direct proofs can be found in \cite[Prop.\ 1.3(2)]{am1996-qtnsubord} and \cite[Prop.\ 4]{b2012-asympEW}.  The Bianchi identities \eqref{eq:qtn-para-CYbianchi} follow upon precomposing the differential Bianchi identity $\d^{\D}\Curv{}{} = 0$ with $J$, before evaluating the algebraic bracket \wrt\ the local \qtn ic\ frame determined by $J$ and taking a trace.  Note that Q-invariance of $\Weyl[\D]{}{}$ means that we do not get an identity for $\CY[\D-]{}{}$ in \itemref{prop:qtn-para-calc}{CYid}, as we did in \itemref{prop:cproj-para-calc}{db}.

\vspace{0.05em}
\section{Associated BGG operators} 
\label{s:qtn-bgg}

Unsurprisingly, the metrisability of a \qtn ic\ structure may be handled in a similar way to the metrisability of \proj\ and \cproj\ structures.  The flat model $\HP$ embeds into the projectivisation of an irreducible $\fg$-representation $\bW$, and the first BGG operator associated to $\bW$ has kernel isomorphic to the space of compatible \qk\ metrics.  We study this BGG operator in Subsection \ref{ss:qtn-bgg-metric}.

As in \cproj\ geometry, the first BGG operator associated to the dual representation $\bW^*$ is a hessian operator whose theory proceeds in much the same way as that of Subsection \ref{ss:cproj-bgg-hess}.  We study this \emph{\qtn ic\ hessian} in Subsection \ref{ss:qtn-bgg-hess}, which controls which Weyl connections have symmetric Q-invariant normalised Ricci tensor.

For the remainder of this chapter we assume that the intrinsic torsion $\Tor{}$ of $\cQ$ vanishes, in which case Weyl curvature $\Weyl{}{} \defeq \Weyl[\D]{}{}$ is \qtn ically\ invariant.

\subsection{Metrisability of \qtn ic\ structures} 
\label{ss:qtn-bgg-metric}

Let $(M,\cQ,\Dspace)$ be a \qtn ic\ manifold of dimension $4n$.  By definition, compatible metrics are \qk\ or locally \hk, depending on their scalar curvature, so are smooth sections of the Q-hermitian subbundle $\Symm[+]{2}T^*M$.  Following the programmes of Sections \ref{s:proj-bgg} and \ref{s:cproj-bgg}, we are interested in the natural decomposition
\begin{equation} \label{eq:qtn-bgg-tracesummands}
  T^*M \tens \Symm[+]{2}TM = (\id\symm TM)_{+} \dsum (T^*M \tens[\trfree] \Symm[+]{2}TM),
\end{equation}
where the first summand is the image of $Z\mapsto \id\symm Z + \qsum J_a\symm J_aZ$ and the second summand is the kernel of the natural trace $T^*M \tens \Symm[+]{2}TM \surjto TM$.  We shall denote projection onto $T^*M \tens[\trfree] \Symm[+]{2}TM$ in \eqref{eq:qtn-bgg-tracesummands} by the subscript ``$\trfree$''.  The proof of the following is similar to \thref{prop:proj-bgg-metriceqn,prop:cproj-bgg-metriceqn}.

\begin{prop} \thlabel{prop:qtn-bgg-metriceqn} The first-order linear differential equation $(\D h)_{\trfree} = 0$ is \qtn ically\ invariant on sections of $\cL^*\tens\Symm[+]{2}TM$. \noproof \end{prop}

We refer to the equation $(\D h)_{\trfree} = 0$ as the \emph{linear metric equation} of \qtn ic\ geometry, and its solutions as \emph{linear metrics}.  We may equivalently write
\begin{equation} \label{eq:qtn-bgg-metriceqnZ}
  \D_X h = X\symm Z^{\D} + \qsum J_aX \symm J_aZ^{\D}
\end{equation}
for some $Z^{\D} \in\s{0}{\cL^*\tens TM}$ and all $X\in\s{0}{TM}$, which should be compared to \eqref{eq:qtn-class-maineqn} and the corresponding equations \eqref{eq:proj-class-maineqn} and \eqref{eq:cproj-class-maineqn} in \proj\ differential geometry and \cproj\ geometry.  Taking a trace in \eqref{eq:qtn-bgg-metriceqnZ} yields $Z^{\D} = \tfrac{1}{2n-1} \p(\D h)$.

Identifying $\cL^{2n+2} \isom \Wedge{4n}TM$, a \non degenerate\ \qk\ metric $g$ induces a section of $\cL^*\tens\Symm[+]{2}TM$ defined by
\begin{equation*}
  h \defeq (\det g)^{1/(4n+4)}\, g^{-1},
\end{equation*}
which we call the \emph{linear metric} associated to $g$.  Then $\det h = (\det g)^{-1/(n+1)}$ is a section of $\cL^4$, from which we may recover $g = (\det h)^{-1/4} \ltens h^{-1}$; \cf\ equation \eqref{eq:qtn-class-gbar}.

\begin{cor} \thlabel{cor:qtn-bgg-metriceqn} There is a linear isomorphism between solutions of the linear metric equation and \qtn ic\ metric connections in $\Dspace$. \noproof \end{cor}

Thus we have reduced the metrisability problem for \qtn ic\ structures to the study of a \qtn ically\ invariant first-order linear differential equation.  Its prolongation proceeds as for \thref{thm:proj-bgg-metricprol,thm:cproj-bgg-metricprol}.

\begin{thm} \thlabel{thm:qtn-bgg-metricprol} There is a linear isomorphism between solutions $h$ of the metric equation and parallel sections of the \qtn ically\ invariant connection
\vspace{-0.2em}
\begin{equation} \label{eq:qtn-bgg-metricprol}
  \D^{\cW}_X \! \colvect{ h \\ Z \\ \lambda } =
  \colvect{ \D_X h - X\symm Z -\qsum J_aX\symm J_aZ \\
            \D_X Z - h(\nRic{X}{},\bdot) - \lambda\ltens X \\
            \D_X \lambda - \nRic{X}{Z} }
  - \tfrac{1}{2n} \!
  \colvect{ 0 \\
            -\Weyl{e_i,X}{h(\ve^i,\bdot)} \\
            h(\CY[\D+]{e_i,X}{},\ve^i) }
\vspace{-0.2em}
\end{equation}
defined on sections $(h,Z,\lambda)$ of $\cW \defeq (\cL^*\tens \Symm[+]{2}TM)\dsum (\cL^*\tens TM)\dsum \cL^*$. \end{thm}

\begin{proof} We may choose a local \qtn ic\ frame $\{J_a\}_{a=1}^3$ and apply the observations in the proof of \thref{thm:cproj-bgg-metricprol} for each $J_a$.  Then only the $J_a$-invariant pieces of the Weyl and Cotton--York tensors contribute to the curvature correction for each $a$, so that only the Q-invariant pieces contribute. \end{proof}

The linear metric equation may be interpreted as a first BGG operator as follows.  As described in Subsection \ref{ss:qtn-para-flat}, the flat model $\HP$ enjoys a \proj\ embedding $\HP \injto \pr{\bW}$ for $\bW \defeq \realrepn{( \Wedge[\bC]{2}\bC^{2n+2} )}$.  Identifying $\bC^{2n+2}$ with the standard representation $\bT = \bH^{n+1}$, an algebraic Weyl structure gives a decomposition
\begin{align*}
  \cpxrepn{\bW}
  &= \Wedge[\bC]{2}( E \tens H^* )
    \isom \Wedge[\bC]{2}E \dsum (E \tens H^*) \dsum \Wedge[\bC]{2}H^* \\
  &\isom \left( \Wedge[\bC]{2}H^* \tens (\Wedge[\bC]{2}E \tens \Wedge[\bC]{2}H) \right)
    \dsum \left( \Wedge[\bC]{2}H^* \tens \cpxbdl{TM} \right) \dsum \Wedge[\bC]{2}H^*
  \\[-2.2em]
\end{align*}
of the complexification of $\bW$, hence giving an isomorphism of the bundle associated to the underlying real representation $\bW$ with the bundle $\cW$ form \thref{thm:qtn-bgg-metricprol}.  The first BGG operator associated to $\bW$ is a differential operator
\begin{equation*}
  \bgg{\bW} :
    \dynkinname{ \dynkinSLHp[0,1,0,0,0,0,0]{3}{0}{4}{qtn-bgg-metric-dom} }
               { L^*\tens\Symm[+]{2}TM }
  \To
    \dynkinname{ \dynkinSLHp[0,1,0,0,1,-2,1]{3}{0}{4}{qtn-bgg-metric-im} }
               {\zbox{ (L^*\tens\Symm[+]{2}TM) \cartan T^*M }},
\end{equation*}
which is clearly first order.  The prolongation connection $\D^{\cW}$ from \thref{thm:qtn-bgg-metricprol} is precisely the prolongation connection of this BGG operator.  We shall call the dimension of the space of parallel sections of $\D^{\cW}$ the \emph{mobility} of the \qtn ic\ structure.  On the flat model, the space of solutions is pointwise isomorphic to the $\fg$-representation $\bW = \realrepn{( \Wedge[\bC]{2}\bC^{2n+2} )}$, so that the mobility is bounded above by $\dim\bW = (n+1)(2n+1)$.

A BGG solution is called \emph{normal} if it is parallel for the tractor connection, \ie\ if the tractor and curvature correction parts of the prolongation connection independently act trivially.  We have the following characterisation of normal solutions of the linear metric equation; a similar characterisation for \proj\ differential geometry is given in \cite{a2008-projholonomy1,cgm2014-projeinstein,gm2014-einsteingeods}, while the \cproj\ case is handled in \cite[Prop.\ 4.8]{cemn2015-cproj}.

\begin{prop} \thlabel{prop:qtn-bgg-normal} A \non degenerate\ $\D^{\cW}$-parallel section $(h,Z,\lambda)$ of $\cW$ is normal if and only if the corresponding Q-hermitian metric $g \defeq (\det h)^{-1/4} \ltens h^{-1}$ is \einstein. \end{prop}

\begin{proof} Let $\D$ be the \LC\ connection of $g$.  Calculating \wrt\ $\D$, if $(h,Z,\lambda)$ is a normal solution then \eqref{eq:qtn-bgg-metricprol} reads
\begin{equation*}
  \D^{\cW}_X \! \colvect{ h \\ Z^{\D} \\ \lambda^{\D} }
    = \colvect{ 0 \\
                -h(\nRic{X}{},\bdot) - \lambda^{\D} \ltens X \\
                \D_X \lambda^{\D} }
    = 0,
\end{equation*}
since $Z^{\D} = \tfrac{1}{2n-1}\p(\D h) = 0$.  We deduce that $\lambda^{\D}$ is a constant multiple of the global trivialisation $(\det h)^{1/4}$ of $\cL$ determined by $h$, say $\lambda = c (\det h)^{1/4}$, so that $h(\nRic{X}{},\bdot) = -c(\det h)^{1/4}\ltens X$.  Applying $h^{-1}$ to both sides gives $\nRic{X}{} = -c g(X,\bdot)$ for all $X\in\s{0}{TM}$, so that $g$ is \einstein\ by \itemref{prop:qtn-para-calc}{nric}.

Conversely suppose that $g$ is an \einstein\ metric with $\nRic{}{} = -cg$.  We must show that the curvature correction in \eqref{eq:qtn-bgg-metricprol} vanishes.  Calculating \wrt\ the \LC\ connection $\D$ of $g$, evidently $\CY{}{} \defeq \d^{\D}\nRic{}{} = -c(\d^{\D} g)$ vanishes.  For the Weyl term $\Weyl{e_i,X}{h(\ve^i,\bdot)}$, since $\Weyl{}{}$ acts trivially on $\cL$ it suffices to show that $\Weyl{e_i,X}{g^{-1}(\ve^i,\bdot)} = 0$.  By the curvature decomposition, we have
\begin{equation} \label{eq:qtn-bgg-normal-2}
  \Weyl{e_i,X}{g^{-1}(\ve^i,\bdot)}
    = \Curv{e_i,X}{g^{-1}(\ve^i,\bdot)} + c\, \algbracw{\id}{g}{e_i,X} \acts g^{-1}(\ve^i,\bdot).
\end{equation}
For the first term on the \rhs\ of \eqref{eq:qtn-bgg-normal-2},
\begin{equation*}
  g( \Curv{e_i,X}{g^{-1}(\ve^i,\bdot)}, Y )
    = ( \Curv{e_i,X}{\ve^i} )(Y)
    = -(\p\Curv{}{})_X(Y)
    = (2n+4)\nRic{X}{Y}
\end{equation*}
for all $X,Y\in\s{0}{TM}$, so that $\Curv{e_i,X}{g^{-1}(\ve^i,\bdot)} = -(2n+4)cX$.  For the second term, direct calculation yields
\begin{align*}
  c\, \algbracw{\id}{g}{e_i,X} \acts g^{-1}(\ve^i,\bdot)
  &= \tfrac{1}{2}c\, \big(
    \ve^i(X)e_i - \qsum[a] [ g(X,J_ae_i)J_a\ve^{i\sharp} - \ve^i(J_aX)J_ae_i ] \\
  &\hspace{3.1em}
    - \ve^i(e_i)X + \qsum[a] [ g(e_i,J_aX)J_a\ve^{i\sharp} - \ve^i(J_ae_i)J_aX ] \big) \\
  &= -(2n+4)cX,
\end{align*}
so that substitution in \eqref{eq:qtn-bgg-normal-2} shows that $\Weyl{e_i,X}{h(\ve^i,\bdot)} = 0$ as well. \end{proof}

Since by \thref{prop:qtn-class-einstein} all \non degenerate\ \qk\ metrics are \einstein, the following is immediate from \thref{prop:qtn-bgg-normal}.

\begin{cor} \thlabel{cor:qtn-bgg-normal} All \non degenerate\ solutions of the metric equation are normal. \noproof \end{cor}

Given the developments of Sections \ref{s:proj-bgg} and \ref{s:cproj-bgg}, we should expect that \qtn ic\ geometry admits a ``big'' Lie algebra of the form $\fh \defeq \bW \dsum (\fg\dsum\bR) \dsum \bW^*$.  This is indeed the case: complexifying gives $\cpxrepn{\fg} = \alg{sl}{2n+2,\bC}$ and $\cpxrepn{\fg}\dsum\bC = \alg{gl}{2n+2,\bC}$, so that
\begin{equation*} \begin{split}
  \cpxrepn{\fh}
  &\isom \Wedge[\bC]{2}\bC^{2n+2} \dsum (\bC^{2n+2} \tens \bC^{2n+2*})
    \dsum \Wedge[\bC]{2}\bC^{2n+2*} \\
  &= \Wedge[\bC]{2}( \bC^{2n+2} \dsum \bC^{2n+2*} ),
\end{split} \end{equation*}
which is the adjoint representation of the complex Lie algebra $\alg{so}{\bC^{2n+2} \dsum \bC^{2n+2*}}$.  Thus $\cpxrepn{\fh}$ has a graded Lie algebra structure isomorphic to $\alg{so}{4n+4,\bC}$, with abelian subalgebras $\cpxrepn{\bW}, \cpxrepn{\bW^*}$ and opposite abelian parabolics $\cpxrepn{\fq} \defeq (\cpxrepn{\fg}\dsum\bC) \ltimes \bW^*$ and $\cpxrepn{\opp{\fq}} \defeq \bW \rtimes (\cpxrepn{\fg}\dsum\bC)$.  To introduce the appropriate real form $\fh$, we note that, up to isomorphism, there is a unique Q-hermitian bilinear form on $\bH^{n+1}$ with skew-symmetric real part and symmetric imaginary part.  The algebra preserving this form is
\vspace{0.3em}
\begin{equation} \label{eq:qtn-bgg-so*}
  \alg[^*]{so}{4n+4} \defeq
    \Setof{ \!\begin{pmatrix} A & B \\ -B^{\top} & \,\conj{\!A} \end{pmatrix} }
          {\, A^{\top} = -A, ~ \,\conj{\!B}{}^{\top} = B }
    \colvectpunct[-0.9em]{,}
\vspace{0.3em}
\end{equation}
which is evidently a (real) subalgebra of $\alg{so}{4n+4,\bC}$.  The \qtn ic\ structure on $\bW$ implies that $\alg[^*]{so}{4n+4}$ is the appropriate real form of $\alg{so}{4n+4,\bC}$, so that the Satake diagrams of $\fq\leq\fh$ are
\begin{equation*}
  \fq = \dynkinSOsp{2}{0}{2}{qtn-bgg-metric-q}
    \hspace{-0.3em} \leq \hspace{0.2em} \dynkinSOs{2}{0}{2}{qtn-bgg-metric-h}
    \hspace{-0.4em} = \fh.
\end{equation*}
Then $H\acts\fq$ is the \grassmannian\ of maximal isotropic \qtn ic\ subspaces of $\bH^{n+1}$.  In terms of associated bundles, $\fh_M \isom \cW \dsum \alg{gl}{\cT} \dsum \cW^*$, where $\alg{gl}{\cT}$ is the Lie algebra bundle of \qtn-linear automorphisms of $\cT$.

\bigskip
\subsection{The \qtn ic\ hessian} 
\label{ss:qtn-bgg-hess}

As in \cproj\ geometry, the first BGG operator associated to the $\fg$-representation $\bW^* \defeq \realrepn{( \Wedge[\bC]{2}\bC^{2n+2*} )}$ is a second order hessian operator: the first BGG operator is
\vspace{0.2em}
\begin{equation*}
  \bgg{\bW^*} :
    \dynkinname{ \dynkinSLHp[0,0,0,0,1,0]{2}{0}{4}{qtn-bgg-hess-dom} }
               { \cL }
  \To
    \dynkinname{ \dynkinSLHp[0,0,0,2,-3,2]{2}{0}{4}{qtn-bgg-hess-im} }
               { \cL \tens \Symm[-]{2}T^*M },
\vspace{0.2em}
\end{equation*}
where $\Symm[-]{2}T^*M \defeq \realrepn{( \Symm[\bC]{2}\cE^* \tens \Symm[\bC]{2}\cH^* )}$ is the complement to the Q-invariant subbundle $\Symm[+]{2}T^*M$ in $\Symm{2}T^*M$.  We call $\bgg{\bW^*}$ the \emph{\qtn ic\ hessian}.  Using the inverse Cartan matrix of $\cpxrepn{\fg}$ it is straightforward to see that $\bgg{\bW^*}$ is second order, so that $\bgg{\bW^*}$ is given by projection of the Ricci-corrected second derivative onto $\cL\tens\Symm[-]{2}T^*M$, \ie\
\begin{equation} \label{eq:qtn-bgg-hesseqn}
  \bgg{\bW^*}_{X,Y}(\ell) = (\D^2_{X,Y}\ell + \ell\ltens\nRic{X}{Y})_{-}.
\end{equation}
A simple calculation shows that $\bgg{\bW^*}$ is \qtn ically\ invariant.  By \thref{lem:qtn-class-connsonL}, a section $\ell \in\s{0}{\cL}$ uniquely determines a \qtn ic\ connection $\D^{\ell}$ defined by $\D^{\ell} \ell = 0$.  We may then characterise solutions of the \qtn ic\ hessian in a similar way to the \cproj\ hessian; \cf\ \thref{prop:cproj-bgg-hesseqn,cor:cproj-bgg-hessmetric}.

\begin{prop} \thlabel{prop:qtn-bgg-hesseqn} A nowhere-vanishing section $\ell\in\s{0}{\cL}$ satisfies $\bgg{\bW^*}(\ell) = 0$ if and only if the normalised Ricci tensor $\nRic[\D^{\ell}]{}{}$ of $\D^{\ell}$ is symmetric and Q-invariant. \noproof \end{prop}

\begin{cor} \thlabel{cor:qtn-bgg-hessmetric} Let $h\in\s{0}{\cL^*\tens\Symm[+]{2}TM}$ be a linear metric.  Then $(\det h)^{1/4}$ lies in the kernel of the \qtn ic\ hessian. \noproof \end{cor}

\begin{rmk} \thlabel{rmk:qtn-bgg-dim4} As remarked above, $4$-dimensional almost \qtn ic\ structure is just a $4$-dimensional conformal structure, corresponding to the identification of $\Sph[4]$ with $\HP[1]$.  The aforementioned \self duality\ condition is equivalent to the choice of a \qtn ic\ hessian or \emph{M{\"o}bius structure}; see \cite[\S 6.4]{bc2010-conf}. \end{rmk}

Since $\bgg{\bW^*}$ is a first BGG operator, its solution space is linearly isomorphic to the space of parallel sections of a connection $\D^{\cW^*}$ on $\cW^* \isom \cL \dsum (\cL\tens T^*M) \dsum (\cL\tens\Symm[+]{2}T^*M)$.

\begin{thm} \thlabel{thm:qtn-bgg-hessprol} There is a linear isomorphism between the space of solutions of the \qtn ic\ hessian $\bgg{\bW^*}$ and the parallel sections of the \qtn ically\ invariant connection
\vspace{0.75em}
\begin{equation} \label{eq:qtn-bgg-hessprol}
  \D^{\cW^*}_X \! \colvect{ \ell \\ \eta \\ \theta }
    = \colvect{ \D_X\ell - \eta(X) \\
                \D_X\eta - \theta(X,\bdot) + \ell \ltens \nRic{X}{} \\
                \D_X\theta + \nRic{X}{}\symm\eta
                  + \qsum J_a\nRic{X}{}\symm J_a\eta
              }
\vspace{0.3em}
\end{equation}
on sections $(\ell,\eta,\theta)$ of $\cW^* \isom \cL \dsum (\cL\tens TM) \dsum (\cL\tens\Symm[+]{2} T^*M)$. \end{thm}

\begin{proof} Choose a unit norm section $J\in\s{0}{\cQ}$ and form the local \qtn ic\ frame $\{J_a\}_{a=1}^3$ with $J_1=J$.  We define $\eta^{\D} \defeq \D\ell \in\s{1}{\cL}$ and $\theta^{\D} \defeq (\D^2\ell + \ell\ltens\nRic{}{})_{-}$, so that $\bgg{\bW^*}(\ell) = 0$ if and only if $\theta^{\D} = \D\eta^{\D} + \ell\ltens \nRic{}{}$.  In this case $\theta^{\D}$ is Q-invariant, and by applying \eqref{eq:qtn-class-qtnconn} we find that the same is true for $\D_X\theta^{\D}$.  Imitating the proof of \thref{thm:cproj-bgg-hessprol} and using \thref{prop:qtn-para-calc}, we conclude that $\bgg{\bW^*}(\ell) = 0$ if and only if $(\ell,\eta^{\D},\theta^{\D})$ is parallel for the \qtn ically\ invariant connection
\vspace{0.3em}
\begin{equation*}
  \tilde{\D}^{\cW^*}_X \!
    \colvect{ \ell \\ \eta \\ \theta }
    = \colvect{ \D_X\ell - \eta(X) \\
                \D_X\eta - \theta(X,\bdot) + \ell \ltens \nRic{X}{} \\
                \D_X\theta + \nRic{X}{}\symm\eta
                  + \qsum J_a\nRic{X}{}\symm J_a\eta
              }
    + \colvect{ 0 \\
                0 \\
                (\Weyl{\cdot,J\cdot}{\eta})(JX)
                  + \ell\ltens\CY[\D+]{\cdot,J\cdot}{JX}
              }
\vspace{0.3em}
\end{equation*}
for each such $J$.  Since the tractor part is independent of $J$, it follows that the curvature correction $\Psi(Y,Z) \defeq \eta^{\D}(\Weyl{JY,Z}{JX}) - \ell\ltens\CY[\D+]{JY,Z}{JX}$ is also independent of the choice of $J$.  Exploiting this independence and the fact that $\Psi$ is Q-invariant, taking $J=J_a$ in the definition of $\Psi$ gives
\begin{align*}
  \Psi(J_bY,Z)
  &= \eta^{\D}(\Weyl{J_cY,Z}{J_aX}) - \ell \ltens \CY[\D+]{J_cY,Z}{J_aX} \\
  &= \eta^{\D}(\Weyl{J_aJ_cY,J_aZ}{J_aX}) - \ell \ltens \CY[\D+]{J_aJ_cY,J_aZ}{J_aX} \\
  &= \eta^{\D}(\Weyl{J_bJ_cY,J_aZ}{J_bX}) - \ell \ltens \CY[\D+]{J_bJ_cY,J_aZ}{J_bX} \\
  &= \eta^{\D}(\Weyl{J_aY,J_aZ}{J_bX}) - \ell \ltens \CY[\D+]{J_aY,J_aZ}{J_bX} \\
  &= \eta^{\D}(\Weyl{J_bY,J_bZ}{J_bX}) - \ell \ltens \CY[\D+]{J_bY,J_bZ}{J_bX} \\
  &= \Psi(Y,J_bZ)
\end{align*}
for all cyclic permutations $(a,b,c)$ of $(1,2,3)$, which contradicts the Q-invariance of $\Psi$.  Therefore $\Psi = 0$ and hence $\tilde{\D}^{\cW^*} = \D^{\cW^*}$ as required. \end{proof}

\begin{cor} \thlabel{cor:qtn-bgg-hessprol} All solutions of the \qtn ic\ hessian are normal. \noproof \end{cor}

Note that the authors of \cite{be1991-paraconf} refer to solutions of $\bgg{\bW^*}$ as \emph{\einstein\ scales}.  The metric prolongation connection $\D^{\cW}$ defined by \eqref{eq:qtn-bgg-metricprol} is dual to $\D^{\cW^*}$ modulo a curvature correction.  In particular, \thref{cor:qtn-bgg-normal} implies that these connections are dual if $(M,\cQ)$ admits a \non degenerate\ linear metric.  

\chapter[\Ppgs]{Projective parabolic \\ geometries} 
\label{c:ppg}

\renewcommand{\algbracadornment}{}
\BufferDynkinLocaltrue
\renewcommand{\dynkinnameoffset}{-0.75}

The primary goal of this thesis is to describe the classical \proj\ structures using a common framework.  The general definition is presented in \thref{defn:ppg-defn-defn} and exploits the existence of a ``big'' \Rspace\ $H\acts\fq$ for each classical \proj\ structure, induced by the graded Lie algebra $\fh \defeq \bW \dsum (\fg\dsum\bR) \dsum \bW^*$ in each case.  Section \ref{s:ppg-defn} is dedicated to setting up this definition, which requires some preliminary work regarding duality for \Rspaces.  We will see in Section \ref{s:ppg-alg} that the algebraic structure of the classical \proj\ structures generalises in almost every detail to the \ppg\ framework; in particular, there is a $\bZ^2$-graded Lie algebra in which we have tight control over the possible Lie brackets.  This is made possible by a Jordan algebra structure on the infinitesimal isotropy representation $\bW \defeq \fh/\fq$ of $\fq$.

The constraints on a \ppg\ are quite strict, resulting in a short classification which we describe in Section \ref{s:ppg-class}.  Via the algebraic work of Section \ref{s:ppg-alg}, the classification over $\bC$ can be phrased entirely in terms of a pair of integers $(r,n)$ which have a purely Lie-theoretic origin; in particular $\dim M = rn$.

Finally, we describe the BGG operators associated to the representations $\bW$ and $\bW^*$ in Section \ref{s:ppg-bgg}.  As for the classical structures, a \ppg\ has a well-defined metrisability problem controlled by the first BGG operator associated to $\bW$, while the first BGG operator associated to $\bW^*$ is a hessian-type equation.  The algebraic work of Section \ref{s:ppg-alg} gives considerable information about these operators and their solutions, allowing us to obtain relatively explicit prolongations.

\section{Definition and \self duality} 
\label{s:ppg-defn}

In the previous three chapters we studied the theories of \proj\ differential geometry, \cproj\ geometry and \qtn ic\ geometry, both from their classical perspectives and as abelian parabolic geometries.  Even the most inattentive reader will have noticed a great deal of similarity within these developments.  We collect the appropriate observations in Subsection \ref{ss:ppg-defn-key}, before introducing a notion of duality for \Rspaces\ in Subsection \ref{ss:ppg-defn-sd}.  This allows us to present a general framework in Subsection \ref{ss:ppg-defn-defn}, which generalises all three classical structures.

\subsection{Key features of classical \proj\ geometries} 
\label{ss:ppg-defn-key}

As explained in the corresponding chapters, \proj\ differential geometry, c-projec-tive geometry and almost \qtn ic\ geometry are all abelian parabolic geometries, respectively modelled on the symmetric \Rspaces\ $\RP[n]$, $\CP[n]$ and $\HP[n]$.  Each of these models enjoys a \proj\ embedding $G\acts\fp \injto \pr{\bW}$ for an appropriate irreducible $\fg$-representation $\bW$, with associated graded representation
\begin{equation*}
  \gr\bW \isom (L^*\tens B) \dsum (L^*\tens\fg/\fp) \dsum L^*.
\end{equation*}
Here $L^*$ is the $1$-dimensional socle of the $\fp^{\perp}$-filtration of $\bW$ and $B$ is an irreducible $\fp^0$-subrepresentation of $\Symm{2}(\fg/\fp)$, such that $L^*\tens B \isom \liehom{0}{\bW}$ may be identified with the zeroth Lie algebra homology.  For \proj, \cproj\ and \qtn ic\ geometries, we had $\bW = \Symm{2}\bR^{n+1}$, $\bW = \realrepn{(\bC^{n+1} \etens \conj{\bC^{n+1}})}$ and $\bW = \realrepn{(\Wedge[\bC]{2} \bC^{2n+2})}$ respectively, where $\etens$ denotes the external tensor product.

We also observed in each case that $\fh \defeq \bW \dsum (\fg\dsum\bR) \dsum \bW^*$ admits the structure of a graded Lie algebra: for \proj, \cproj\ and almost \qtn ic\ geometries, we had $\fh \isom \alg{sp}{2n+2,\bR}$, $\fh \isom \alg{su}{n+1,n+1}$ and $\fh \isom \alg[^*]{so}{4n+4}$ respectively.  Moreover it follows that the trivial central extension $\fq^0 \defeq \fg\dsum\bR$ is a reductive Lie algebra, with $\fq \defeq \fq^0 \ltimes \bW^*$ and $\opp{\fq} \defeq \bW \rtimes \fq^0$ opposite abelian parabolic subalgebras of $\fh$; in particular, $H\acts\fq$ is a symmetric \Rspace.  The crucial observation about the \Rspaces\ $H\acts\fq$ is as follow.  For the classical \proj\ structures, the Satake diagrams of $\fq \leq \fh$ are
\smallskip
\begin{equation*} \begin{gathered}
  \dynkinCp{2}{0}{3}{ppg-defn-key-C}
    \leq \alg{sp}{2n+2,\bR}, \quad
  \dynkinSUp{2}{0}{2}{ppg-defn-key-A} \hspace{-0.7em}
    \leq \alg{su}{n+1,n+1} \\
  \text{and}\quad
  \dynkinSOsp{3}{0}{3}{ppg-defn-key-D} \hspace{-0.6em}
    \leq \alg[^*]{so}{4n+4}.
\end{gathered}
\medskip
\end{equation*}
If $\longest[\fh]$ is the longest element of the Weyl group of $\fh$, $-\longest[\fh]$ induces an automorphism of the Satake diagram, which is the identity for $\alg{sp}{2n+2,\bR}$ and $\alg[^*]{so}{4n+4}$, and the involution indicated by the arrows for $\alg{su}{n+1,n+1}$.  In each case the single crossed node of $\fq$ is preserved by this automorphism.  To interpret this in terms of \proj\ embeddings as in Subsection \ref{ss:lie-para-rspace}, choose a Cartan subalgebra $\ft$ and simple subsystem \wrt\ which $\fq$ is a standard parabolic.  Then by \thref{cor:lie-para-proj} there is a \proj\ embedding $H\acts\fq \injto \pr{\bV^*}$ for any irreducible $\fh$-representation $\bV$ whose highest weight $\lambda \in \ft^*$ is supported on the single crossed node of $\fq$.  The highest weight of $\bV^*$ is $-\longest[\fh](\lambda)$, which in each case is also supported on the single crossed node of $\fq$, and thus the \Rspaces\ induced by the \proj\ embeddings into $\pr{\bV^*}$ and $\pr{\bV}$ are isomorphic; we will interpret this as a \self duality\ condition in the next subsection.

\subsection{Duality for \Rspaces} 
\label{ss:ppg-defn-sd}

We continue now to describe a form of duality for \Rspaces, following \cite{bdpp2011-rspaces}.  Let $H\acts\fq$ be an \Rspace\ and consider the set $(H\acts\fq)^*$ of parabolic subalgebras of $\fh$ which are opposite to some $\fq' \in H\acts\fq$.  By the following \cite[Prop.\ 2.3]{bdpp2011-rspaces}, $(H\acts\fq)^*$ is an \Rspace.

\begin{lem} $(H\acts\fq)^*$ is a single conjugacy class of parabolic subalgebras of $\fh$. \end{lem}

\begin{proof} Let $\fq_1, \fq_2 \in (H\acts\fq)^*$, with $\fq_1, \fq_2$ opposite to $\fq, h\acts\fq$ for $h\in H$.  Then $h^{-1}\acts\fq_2$ is opposite to $\fq$, giving $qh^{-1} \acts \fq_2 = \fq_1$ for a unique element $q \in \exp\fq^{\perp}$ by \thref{lem:lie-para-expopp}. \end{proof}

Clearly $(H\acts\fq)^*$ has the same height as $H\acts\fq$.  Moreover since we obviously have $(H\acts\fq)^{**} = H\acts\fq$, we make the following definition \cite[Defn.\ 2.4]{bdpp2011-rspaces}.

\begin{defn} $(H\acts\fq)^*$ is the \emph{dual} of $H\acts\fq$.  An \Rspace\ is called \emph{\self dual} if it coincides with its dual. \end{defn}

There is always a \non canonical\ diffeomorphism between $H\acts\fq$ and its dual.  Indeed, the Cartan involution $\theta$ corresponding to a maximal compact subgroup of $H$ induces a diffeomorphism $H\acts\fq \ni \fq' \mapsto \theta(\fq') \in (H\acts\fq)^*$; see \cite[p.\ 8]{bdpp2011-rspaces}.

\begin{cor} \thlabel{cor:ppg-defn-sdqq} $H\acts\fq$ is \self dual\ if and only if $\fq$ is conjugate to every opposite parabolic. \end{cor}

\begin{proof} Let $\opp{\fq}$ be opposite to $\fq$.  Then $H\acts\fq$ is \self dual\ if and only if for all $h_1 \acts \opp{\fq} \in (H\acts\fq)^*$ and $h_2 \acts \fq \in H\acts\fq$, there exists an $h \in H$ such that $hh_1 \acts \opp{\fq} = h_2\acts\fq$. \end{proof}

It is straightforward to detect \self duality\ from the Satake diagram of $\fq$.  For this recall that if $\longest[\fh] \in \sW_{\fh}$ is the longest element of the Weyl group of $\fh$, then $-\longest[\fh]$ permutes the simple roots, hence inducing an involution of the Satake diagram of $\fh$.

\begin{lem} \thlabel{lem:ppg-defn-longest} \emph{\cite[p.\ 8]{bdpp2011-rspaces}} $H\acts\fq$ is \self dual\ if and only if the set of crossed nodes of the Satake diagram of $H\acts\fq$ is preserved by $-\longest[\fh]$. \end{lem}

\begin{proof} First note that parabolics $\fq, \fq' \leq \fh$ are opposite if and only if their complexifications are opposite, so we may assume that $H\acts\fq$ is complex.

Now choose a Cartan subalgebra $\ft\leq\fh$ and a simple subsystem $\Delta^0$ \wrt\ which $\fq$ is a standard parabolic, corresponding to a subset $\Sigma \subseteq \Delta^0$.  Then since $H\acts\fq$ and $(H\acts\fq)^*$ each contain a unique standard parabolic, it suffices to see when these coincide.  Using that $\longest[\fh]$ is an involution, the standard parabolic determined by $-\longest[\fh](\Sigma) \subseteq \Delta^0$ is conjugate, via $\longest[\fh]$, to the standard parabolic opposite to $\fq$, \ie\ the parabolic subalgebra of $\fh$ consisting of $\ft$ and root spaces $\fh_{\alpha}$ for which $\fh_{-\alpha} \leq \fq$.  In particular the parabolic determined by $-\longest[\fh](\Sigma)$ lies in $(H\acts\fq)^*$, with equality if and only if $-\longest[\fh](\Sigma) = \Sigma$. \end{proof}

Equivalently, for any $\fh$-representation $\bV$ whose highest weight is supported on the crossed nodes of $H\acts\fq$, the \Rspaces\ determined by \proj\ embeddings into $\pr{\bV^*}$ and $\pr{\bV}$ are isomorphic.  While this gives a convenient characterisation of \self duality\ when we know the Satake diagram, it will be useful to have a more theoretical condition.  

\begin{lem} \thlabel{lem:ppg-defn-opp} $(\exp(x)\acts\fq)^{\perp} = \exp(x) \acts \fq^{\perp}$ for all $x\in\fh$. \end{lem}

\begin{proof} Since $\exp(x)$ is an automorphism of $\fh$, we have $\killing{ \exp(x)\acts\fq^{\perp} }{ \exp(x)\acts\fq } = \killing{ \fq^{\perp} }{ \fq } = 0$ by invariance of the Killing form. \end{proof}

By combining results from \cite[\S 4]{bdpp2011-rspaces}, we arrive at the following characterisation of \self duality for symmetric \Rspaces.

\begin{prop} \thlabel{prop:ppg-defn-sd} Let $H\acts\fq$ be a symmetric \Rspace.  Then the following are equivalent:
\begin{enumerate}
  \item \label{prop:ppg-defn-sd-sd}
  $H\acts\fq$ is \self dual;
  
  \item \label{prop:ppg-defn-sd-ker}
  There is an $f \in \fq^{\perp}$ such that $\ker(\ad f)^2 = \fq$;
  
  \item \label{prop:ppg-defn-sd-lift}
  For any parabolic $\opp{\fq}$ opposite to $\fq$, there exists $e \in \opp{\fq}^{\perp}$ and $f \in \fq^{\perp}$ such that $\liebrac{e}{f}$ equals twice the algebraic Weyl structure $\xi$ induced by $\opp{\fq}$.
\end{enumerate}
\end{prop}

\begin{proof}
Since $\fq^{\perp}$ is abelian we have $\liebrac{\fh}{\fq^{\perp}} \leq \fq$ and $\liebrac{\fq}{\fq^{\perp}} \leq \fq^{\perp}$, and hence $(\ad x)^2$ vanishes on $\fq$ for all $x \in \opp{\fq}^{\perp}$.  Since $\Ad{\exp(x)} = \exp(\ad x)$, we have
\vspace{0.05em}
\begin{equation} \label{eq:ppg-defn-bch}
  \exp(x)\acts y = y + \liebrac{x}{y} + \tfrac{1}{2}(\ad x)^2(y)
\vspace{0.05em}
\end{equation}
for all $y\in\fh$ by the series expansion of the exponential map.

\smallskip

\equivref{prop:ppg-defn-sd}{sd}{ker}
This is mostly \cite[Lem.\ 4.1]{bdpp2011-rspaces}.  Denote by $\Omega_{\fq}$ the set of parabolic subalgebras of $\fh$ which are opposite to $\fq$, which is a dense open subset of $H\acts\fq$ by \cite{t1965-celldecomp}.  Therefore if $\opp{\fq} \in \Omega_{\fq}$, the set $\Omega_{\fq} \intsct \Omega_{\opp{\fq}}$ is also open and dense; \thref{lem:lie-para-expopp} then provides an $f \in \fq^{\perp}$ such that $\exp(f)\acts\opp{\fq} \in \Omega_{\fq} \intsct \Omega_{\opp{\fq}}$.

For arbitrary $x\in\fq^{\perp}$, we have $\exp(x)\acts\opp{\fq} \in \Omega_{\fq} \intsct \Omega_{\opp{\fq}}$ if and only if $\exp(x)\acts\opp{\fq}$ is opposite to $\opp{\fq}$, if and only if $\exp(x)\acts\opp{\fq}^{\perp} \intsct \opp{\fq} = 0$.  Taking $y \in \opp{\fq}^{\perp}$, the first two terms in \eqref{eq:ppg-defn-bch} live in $\opp{\fq}$, while the third lives in $\fq^{\perp}$.  Since $\fq, \opp{\fq}$ are opposite, we have $\fh = \opp{\fq} \dsum \fq^{\perp}$ and hence $\exp(x) \acts y \in \exp(x)\acts\opp{\fq}^{\perp} \intsct \opp{\fq}$ if and only if $(\ad x)^2(y) = 0$.  Therefore $\exp(x)\acts\opp{\fq} \in \Omega_{\fq} \intsct \Omega_{\opp{\fq}}$ if and only if $(\ad x)^2$ is injective on $\opp{\fq}^{\perp}$; since also $\fh = \opp{\fq}^{\perp} \dsum \fq$ and $\ker(\ad x)^2 \subseteq \fq$, this is equivalent to $\ker(\ad x)^2 = \fq$.

\smallskip

\equivref{prop:ppg-defn-sd}{ker}{lift}
This is a special case of \cite[Prop.\ 4.3]{bdpp2011-rspaces}.  If $\ker(\ad f)^2 = \fq$ for some $f \in \fq^{\perp}$ then $\fq$, $\opp{\fq}$ and $\exp(f)\acts\opp{\fq}$ are mutually opposite.  Then by \thref{lem:lie-para-expopp} there are unique elements $e \in \opp{\fq}^{\perp}$ and $x \in (\exp(f)\acts\opp{\fq})^{\perp} = \exp(f)\acts\opp{\fq}^{\perp}$ such that $\exp(f) \acts \opp{\fq} = \exp(e) \acts \fq$ and $\opp{\fq} = \exp(x) \acts \fq$.  If $\xi$ is the algebraic Weyl structure of the pair $(\fq,\opp{\fq})$, it is straightforward to see that the algebraic Weyl structure of $(\fq,\exp(f)\acts\opp{\fq})$ is $\exp(f)\acts\xi = \xi + f$, where the last equality follows by \eqref{eq:ppg-defn-bch}.  Since $\exp(f)\acts\opp{\fq} = \exp(e)\acts\fq$, the algebraic Weyl structure of $(\opp{\fq}, \exp(f)\acts\opp{\fq})$ is similarly $\exp(e)\acts(-\xi) = -\xi + e$.  Writing the pair $(\fq, \exp(f)\acts\opp{\fq})$ as $(\exp(-x)\acts\opp{\fq}, \exp(e)\acts\fq)$, it follows that $\exp(f)\acts\xi = -\exp(-x) \exp(e) \acts \xi$.  Using \eqref{eq:ppg-defn-bch} and that $x \in \exp(e)\acts\fq^{\perp}$, we have $\exp(-x) \exp(e) \acts \xi = \exp(e)\acts\xi - x$. Therefore
\begin{equation*}
  x = \exp(f)\acts\xi + \exp(e)\acts\xi = 2\xi - e + f
\end{equation*}
and hence we obtain
\begin{equation} \label{eq:ppg-defn-sd-1}
  \liebrac{ \exp(f)\acts\xi \, }{ \, \exp(e)\acts\xi }
    =  \liebrac{ x }{ \exp(e)\acts\xi }
    =  x = 2\xi - e + f.
\end{equation}
On the other hand, we have
\begin{equation} \label{eq:ppg-defn-sd-2} \begin{aligned}
  \liebrac{ \exp(f)\acts\xi }{ \exp(e)\acts\xi }
    &= \liebrac{ \xi + f }{ \xi - e } \\
    &= -\liebrac{\xi}{e} + \liebrac{f}{\xi} - \liebrac{f}{e}
     = -e + f + \liebrac{e}{f}.
\end{aligned} \end{equation}
Comparing \eqref{eq:ppg-defn-sd-1} and \eqref{eq:ppg-defn-sd-2} now gives $\liebrac{e}{f} = 2\xi$.

\smallskip

\equivref{prop:ppg-defn-sd}{lift}{sd}
Let $\opp{\fq}$ be opposite to $\fq$ with corresponding algebraic Weyl structure $\xi$, and suppose that $e \in \opp{\fq}^{\perp}$ and $f \in \fq^{\perp}$ satisfy $\liebrac{e}{f} = 2\xi$.  Then applying \eqref{eq:ppg-defn-bch} gives
\begin{align*}
  \exp(e)\exp(-f)\exp(e) \acts \xi
  &= \exp(e)\exp(-f) \acts (\xi - e) \\
  &= \exp(x) \acts (\xi -f -e -2\xi + f) \\
  &= \exp(x) \acts (-\xi - e) \\
  &= -\xi.
\end{align*}
It follows that $\exp(e)\exp(-f)\exp(e)$ is a graded automorphism of $\gr\fh$, sending the graded component $\fh_{(i)}$ to $\fh_{(-i)}$ for all $i \in \setof{+1,0,-1}{}$.  In particular since $\xi$ induces isomorphisms $\fq \isom \fh_{(0)} \dsum \fh_{(-1)}$ and $\opp{\fq} \isom \fh_{(1)} \dsum \fh_{(0)}$, we have $\exp(e)\exp(-f)\exp(e) \acts \opp{\fq} = \fq$.  Thus $\opp{\fq}$ is conjugate to $\fq$ and hence $H\acts\fq$ is \self dual by \thref{cor:ppg-defn-sdqq}. \end{proof}

Elements $f \in \fq^{\perp}$ satisfying $\ker(\ad f)^2 = \fq$ are called \emph{regular}.  Equivalently, $(\ad f)^2$ factors to an isomorphism $F_f : \fh/\fq \to \fq^{\perp}$.  With notation as in \itemref{prop:ppg-defn-sd}{lift}, we immediately see that $F_f(e) = \liebrac{ f }{ \liebrac{f}{e} } = -\liebrac{f}{2\xi} = -2f$ and hence $e = -2F_f^{-1}(f)$.

\subsection{General definition} 
\label{ss:ppg-defn-defn}

A symmetric \Rspace\ $H\acts\fq$ induces an \Rspace\ of smaller dimension using the \proj\ embeddings of Subsection \ref{ss:lie-para-rspace}.  For this, let $\fq^0 \defeq \fq/\fq^{\perp}$ be the reductive Levi factor of $\fq$, and let $\bW \defeq \fh/\fq$.  We may decompose $\fq^0$ into its semisimple part $\fg \defeq \liebrac{\fq^0}{\fq^0}$ and centre $\liecenter{\fq^0}$, whose dimension equals the number of simple factors of $\fh$ by \thref{lem:lie-para-levi}.  Since $\fq$ is a subalgebra of $\fh$, the adjoint representation of $\fh$ induces a representation of $\fq$ on $\bW$ via $y \acts (x+\fq) = \liebrac{y}{x} + \fq$.  Since $\liebrac{\fq^{\perp}}{\fh} = \fq$ by \thref{lem:lie-para-abl}, the action of $\fq$ on $\bW = \fh/\fq$ descends to an action $\fq^0 \defeq \fq/\fq^{\perp}$, thus defining a representation of $\fg$ on $\bW$.  There is an \Rspace\ $G\acts\fp$ associated to both of the \proj\ embeddings $G\acts\fp \injto \pr{\bW}$ and $G\acts\fp \injto \pr{\bW^*}$, given respectively by crossing the nodes on which the highest weights of $\bW^*$ and $\bW$ are supported.  We take this ``top down'' view of $H\acts\fq$ and $G\acts\fp$ for our general definition.

\begin{defn} \thlabel{defn:ppg-defn-defn} Let $H\acts\fq$ be a symmetric \Rspace\ with infinitesimal isotropy representation $\bW \defeq \fh/\fq$.
\begin{enumerate}
  \item \label{defn:ppg-defn-defn-iso}
  The \Rspace\ $G\acts\fp \injto \pr{\bW}$ will be called the \emph{isotropy \Rspace} of $H\acts\fq$.
  
  \item \label{defn:ppg-defn-defn-ppg}
  A \emph{\ppg} is a parabolic geometry modelled on the isotropy \Rspace\ of a \self dual\ symmetric \Rspace.
\end{enumerate}
\end{defn}

Recall that for the classical \proj\ structures the \Rspace\ $G\acts\fp$ is abelian, with the $\fp^{\perp}$-filtration of $\bW$ of height two.  These are not \apriori\ true for a general \ppg, and establishing these properties will be the focus of Section \ref{s:ppg-alg}.

For later use, we record the following method \cite[Prop.\ 3.2.2]{cs2009-parabolic1} for determining the Dynkin (or Satake) diagram of $\fg$ from that of $H\acts\fq$.

\begin{prop} \thlabel{prop:ppg-defn-dynkin} Let $H\acts\fq$ be a symmetric \Rspace.  Then the Satake diagram of $\fg$ is given by removing all crossed nodes from $\fq$ and their associated edges. \end{prop}

\begin{proof} It suffices to consider the case that $\fh$ is complex, so choose a Cartan subalgebra $\ft\leq\fh$ and a simple subsystem $\Delta^0$ \wrt\ which $\fq$ is the standard parabolic corresponding to a subset $\Sigma \subseteq \Delta^0$.  These choices induce an algebraic Weyl structure of $\fq$, identifying $\fg$ with a subalgebra of $\fh$.  By \thref{lem:lie-para-levi} we know that $\Delta^0\setminus\Sigma$ forms a simple subsystem for $\fg$, so it remains to describe the Cartan matrix.

Suppose first that $\alpha,\beta \in \Delta^0\setminus\Sigma$ lie in different simple factors of $\fg$.  Then $\alpha,\beta$ are orthogonal \wrt\ the Killing form of $\fg$.  On the other hand, since $\liebrac{\fh_{\alpha}}{\fh_{\beta}} = 0$ we cannot have $\alpha+\beta$ as a root of $\fh$.  But $\alpha-\beta$ is also not a root, so by considering the $\alpha$-root string through $\beta$ we see that $\killing{\alpha}{\beta} = 0$ in $\fh$.

Suppose now that $\alpha, \beta$ lie in the same simple factor of $\fg$.  Since the Killing form of $\fg$ is determined uniquely up to scale on each simple factor, it must coincide with the restriction of the Killing form of $\fh$.  Thus the Cartan integers \wrt\ $\fg$ coincide with the Cartan integers \wrt\ $\fh$. \end{proof}

Since \thref{prop:ppg-defn-dynkin} will form the basis of our classification, we give an example.

\begin{expl} \thlabel{expl:ppg-defn-dynkin}
(1) Let $\fh = \alg[_6]{e}{\bC}$ and consider its \Rspace\
\vspace{-0.4em}
\begin{equation*}
    H\acts\fq \defeq \dynkinEp{6}{ppg-defn-defn-e6q}
      \Injto \pr{\erepn},
\end{equation*}
where $\erepn$ is the $27$-dimensional representation of $\alg[_6]{e}{\bC}$.  The longest element of $\fh$ induces the usual automorphism of the Dynkin diagram, so that $H\acts\fq$ is not \self dual\ by \thref{lem:ppg-defn-longest}.  Nevertheless, \thref{prop:ppg-defn-dynkin} tells us that the Dynkin type of the semisimple part $\fg$ of $\fq^0$ is
\begin{equation*}
  \fg = \dynkin{ \DynkinLine{0}{0}{2}{ 0   };
                 \DynkinLine{2}{0}{3}{ 0.75};
                 \DynkinLine{2}{0}{3}{-0.75};
                 \DynkinWDot{0}{ 0   };
                 \DynkinWDot{1}{ 0   };
                 \DynkinWDot{2}{ 0   };
                 \DynkinWDot{3}{ 0.75};
                 \DynkinWDot{3}{-0.75}; }{ppg-defn-defn-e6g}
    = \alg{so}{10,\bC}.
\end{equation*}

\smallskip

\noindent
(2) By way of consistency with the classical \proj\ structures, the reader can check that for $\fh = \alg{sp}{2n+2,\bR}$ we obtain $\fg = \alg{sl}{n+1,\bR}$; for $\fh = \alg{su}{n+1,n+1}$ we obtain $\fg = \alg{sl}{n+1,\bC}$; and for $\fh = \alg[^*]{so}{4n+4}$ we obtain $\fg = \alg{sl}{n+1,\bH}$.
\end{expl}

\section{Algebraic structure} 
\label{s:ppg-alg}

Fix a \self dual\ symmetric \Rspace\ $H\acts\fq$ with infinitesimal isotropy representation $\bW \defeq \fh/\fq$ and isotropy \Rspace\ $G\acts\fp$.  Our goal in this section is to investigate the algebraic structure of a \ppg\ in more detail.  Indeed, we do not yet know whether $G\acts\fp$ consists of abelian parabolics as in the classical cases, or the structure of $\bW$.  It suffices to consider the case that $\fh$ is simple by the following \cite[Lem.\ 3.2.3]{cs2009-parabolic1}.

\begin{lem} \thlabel{lem:ppg-alg-simple} Let $\fh$ be semisimple with abelian parabolic $\fq$, such that \wrt\ any algebraic Weyl structure no simple ideal of $\fh$ is contained in $\fq^0$.  Then each simple factor of $\fh$ has an abelian parabolic, such that $\fq$ is the direct sum of these parabolics. \noproof \end{lem}

The key step in the determination of the algebraic structures of $\fg$ and $\bW$ is the relationship between \self dual\ symmetric \Rspaces\ and Jordan algebras, which we describe in Subsection \ref{ss:ppg-alg-jordan}.  In Subsection \ref{ss:ppg-alg-str} we show that the structure found for the classical cases carries over in almost every detail; in particular, we see in Subsection \ref{ss:ppg-alg-Z2} that this gives $\fh$ the structure of a $\bZ^2$-graded algebra, which allows us to determine Lie brackets between many of its summands.  Finally, in Subsection \ref{ss:ppg-alg-L} we study a $1$-dimensional representation $L$ parametrising the space of Weyl connections, which also provides key information regarding traces of certain Lie brackets.

\subsection{Relation to Jordan algebras} 
\label{ss:ppg-alg-jordan}

We shall see that the algebraic structure of $\bW$ is neatly described in terms of Jordan algebras.  Fix a field $\bk$, which for us will be $\bR$ or $\bC$.

\begin{defn} \thlabel{defn:ppg-alg-jordan} A \emph{Jordan algebra} over $\bk$ is a commutative (but \non associative) $\bk$-algebra $(\bJ,\jmult)$ satisfying the \emph{Jordan identity}
\begin{equation} \label{eq:ppg-alg-jordanid}
  (x \jmult y) \jmult (x \jmult x) = x \jmult (y \jmult (x \jmult x))
\end{equation}
for all $x,y \in \bJ$.  The Jordan identity implies that $\bJ$ is ``power associative'', meaning that $(x^m \jmult y) \jmult x^n = x^m \jmult (y \jmult x^n)$ for all $x,y \in \bJ$ and $m,n \in \bZ_{>0}$. \end{defn}

\begin{expl} Given any associative $\bk$-algebra $\bJ$, it can be given the structure of a Jordan algebra via $x \jmult y \defeq \tfrac{1}{2}(xy + yx)$, where juxtaposition denotes the original multiplication in $\bJ$.  Jordan algebras obtained in this way are called \emph{special}. \end{expl}

Jordan algebras were introduced by Jordan \cite{j1933-jordanalg} as an algebraic framework for the system of observable quantities in quantum mechanics.  A detailed history of Jordan algebras and their role in physics may be found in McCrimmon's excellent book \cite{m2007-jordan}.  Our interaction with Jordan algebras shall be limited to the close relationship between \self dual\ symmetric \Rspaces\ and Jordan algebras, where we will almost exclusively favour Lie-theoretic language.  By \thref{prop:ppg-defn-sd}, \self duality\ is equivalent to the existence of a regular element $f \in \fq^{\perp}$ with $F_f \defeq (\ad f)^2 : \fh/\fq \to \fq^{\perp}$ an isomorphism.

\begin{thm}[Meyberg] \thlabel{thm:ppg-alg-Wjordan} Suppose that $H\acts\fq$ is a \self dual\ symmetric \Rspace\ and choose an algebraic Weyl structure $\xi^{\fq}$.  Let $f \in \bW \defeq \fh/\fq$ be a regular element with corresponding isomorphism $F_f \defeq (\ad f)^2 : \bW \to \bW^*$.  Then
\begin{equation*}
  x \jmult y \defeq \tfrac{1}{2} \liebrac{ \liebrac{x}{f} }{ y }
\end{equation*}
makes $(\bW,\jmult)$ into a semisimple Jordan algebra with identity element $e \defeq -2 F_f^{-1}(f)$. \end{thm}

\begin{sketchproof} Since $\bW \leq \fh$ is an abelian subalgebra \wrt\ the algebraic Weyl structure, commutativity of $\jmult$ follows immediately from the Jacobi identity.  The Jordan identity \eqref{eq:ppg-alg-jordanid} follows by a direct calculation via the Jacobi identity which we omit; see \cite{m1970-jordantriple}.  To see that $e$ is a unit for $\jmult$, note that $e = -2F_f^{-1}(f)$ is the element from \itemref{prop:ppg-defn-sd}{lift} which satisfies $\liebrac{e}{f} = 2\xi^{\fq}$.  Therefore $e \jmult x = \tfrac{1}{2} \liebrac{ \liebrac{e}{f} }{ x } = \liebrac{\xi^{\fq}}{x} = x$ for all $x \in \bW$ as required. \end{sketchproof}

Further details may be found in \cite{k1970-jordandg, l1971-jrspace, m1970-jordantriple}; see also \cite[\S 1.4]{b2011-jordannonassoc} and \cite[\S 7]{b2013-jordanlie}.  Dually, $\bW^*$ becomes a unital semisimple Jordan algebra with product $\alpha \jmult \beta \defeq \tfrac{1}{2} \liebrac{ \liebrac{\alpha}{e} }{ \beta }$ for all $\alpha, \beta \in \bW^*$.  The \sltriple\ $(e,\liebrac{e}{f},f)$ integrates to a group homomorphism $\varphi : \grp{SL}{2,\bC} \to H$ such that the adjoint action of $\varphi \left( \inlinematrix{ 0 & 1 \\ -1 & 0 } \right)$ induces an isomorphism of Jordan algebras between $\bW$ and $\bW^*$; see \cite[\S 2]{ks2015-smallrepns}.
Of course, the regular element $f \in \bW^*$ need not be unique.  Jordan products induced by different regular elements correspond to different \emph{isotopes} of $\bW$; see \cite{m2007-jordan}.

\begin{rmk} \thlabel{rmk:ppg-alg-tkk} Following work of Tits \cite{t1962-kktalg}, Koecher \cite{k1967-imbedding1, k1968-imbedding2} and Meyberg \cite{m1970-jordantriple}, a Jordan algebra $(\bW,\jmult)$ can conversely be embedded into a semisimple Lie algebra $\fh$ as follows.  By defining
\begin{equation*} 
  \jtrip{x}{y}{z} \defeq (x \jmult y) \jmult z - y \jmult (x \jmult z) + x \jmult (y \jmult z),
\end{equation*}
we give $(\bW, \jtrip{}{}{})$ the structure of a \emph{Jordan triple system}, meaning that $\jtrip{}{}{}$ is symmetric in the last two entries and $L_{x,y} : \bW \to \bW$ defined by $L_{x,y}(z) \defeq \jtrip{x}{y}{z}$ satisfies $\liebrac{ L_{x,y} }{ L_{z,w} } = L_{z,\jtrip{x}{y}{w}} - L_{\jtrip{y}{x}{z},w}$.  The space $\alg{der}{\bW} \defeq \setof{ L_{x,y} }{ x,y \in \bW }$ is then a Lie algebra under commutator, acting naturally on $\bW$ and $\bW^*$ in the obvious way.  If the trace form $\killing{x}{y} \defeq \tr(L_{x,y})$ is \non degenerate, the direct sum
\vspace{0.1em}
\begin{equation*}
  \fh \defeq \bW \dsum \alg{der}{\bW} \dsum \bW^*
\vspace{0.1em}
\end{equation*}
is a graded semisimple Lie algebra with the brackets
\vspace{0.1em}
\begin{equation*}
  \liebrac{x}{\alpha} = L_{x,\alpha}, \quad
  \liebrac{x}{y} = 0 = \liebrac{\alpha}{\beta}, \quad
  \liebrac{A}{x} = Ax
  \quad\text{and}\quad
  \liebrac{A}{\alpha} = -\alpha \circ A,
\vspace{0.1em}
\end{equation*}
called the \emph{Kantor--Koecher--Tits algebra} of $\bW$; see \cite{b2013-jordanlie, c2008-jtriples} and the references therein. \end{rmk}

An element $e$ in a Jordan algebra $(\bJ,\jmult)$ is called an \emph{idempotent} if $e \jmult e = e$.  Given multiple idempotents $e_1, \ldots, e_n$, they are \emph{orthogonal} if $e_i \jmult e_j = \delta_{ij} e_i$ for all $i,j$.  As for commutative rings, there is a \emph{Peirce decomposition} of $\bJ$ defined as follows.  Denote by $m_e : \bJ \to \bJ$ the multiplication map by $e \in \bJ$.  Then if $e$ is idempotent, Albert proves \cite[Eqn.\ (13)]{a1947-jordanstr} that $m_e$ has minimal polynomial
\begin{equation*}
  \tfrac{1}{2} m_e (2m_e-\id) (m_e-\id),
\end{equation*}
\ie\ $m_e$ has eigenvalues $0$, $1/2$ and $1$.  Denoting the eigenspaces by $\jespace{e}{\lambda}$, it follows that $\jespace{e}{0}$ and $\jespace{e}{1}$ are mutually orthogonal Jordan subalgebras of $\bJ$, while \cite[Thm.\ 6]{a1947-jordanstr}
\begin{equation} \label{eq:ppg-alg-idemid} \begin{gathered}
  \jespace{e}{1/2} \jmult \jespace{e}{1/2}
    \leq \jespace{e}{0} \dsum \jespace{e}{1}, \quad
  \jespace{e}{1} \jmult \jespace{e}{1/2} \leq \jespace{e}{1/2} \\
  \text{and}\quad
  \jespace{e}{0} \jmult \jespace{e}{1/2} \leq \jespace{e}{1/2}.
\end{gathered}
\end{equation}
On the other hand, if $e_1, e_2$ are orthogonal idempotents of $\bJ$ then one finds that $(a \jmult e_1) \jmult e_2 = (a \jmult e_2) \jmult e_1$ for all $a \in \bJ$, which lies in the intersection $\jespace{e_1}{1/2} \intsct \jespace{e_2}{1/2}$.  Suppose now that $\{ e_1, \ldots, e_n \}$ is a set of of idempotents for which $e_1 + \cdots + e_n = e$ is the identity in $\bJ$, and write $\bJ_i \defeq \jespace{e_i}{1}$ and $\bJ_{ij} \defeq \jespace{e_i}{1/2} \intsct \jespace{e_j}{1/2}$ for all $i \neq j$.  Then $\bJ_i$ is a $1$-dimensional vector space, spanned by $e_i$, while the $\bJ_{ij}$ are \non empty\ by the above.

\begin{thm} \thlabel{thm:ppg-alg-peirce} Suppose that $(\bJ,\jmult)$ is a Jordan algebra with unit $e \in \bJ$, and let $e_1, \ldots, e_n \in \bJ$ be pairwise orthogonal idempotents satisfying $e_1 + \cdots + e_n = e$.  Then there is a Peirce decomposition
\begin{equation*}
  \bJ = \Dsum{i=1}{n} \left( \bJ_i \dsum {\textstyle \Dsum{j>i}{}} \, \bJ_{ij} \right),
\end{equation*}
where $\bJ_i \defeq \jespace{e_i}{1}$ and $\bJ_{ij} \defeq \jespace{e_i}{1/2} \intsct \jespace{e_j}{1/2}$ for $i \neq j$. \noproof \end{thm}

The Peirce decomposition was first considered by Jordan, von Neumann and Wigner in the seminal paper \cite{jnw1934-jordanqm} for so-called ``totally real'' Jordan algebras.  The general case was developed by Albert \cite[p.\ 559]{a1947-jordanstr}.  For later use, \eqref{eq:ppg-alg-idemid} and the Peirce decomposition allow us to prove the following multiplication properties for the eigenspaces \cite[Thm.\ 12]{a1947-jordanstr}.

\begin{cor} \thlabel{cor:ppg-alg-idemmult} The idempotent eigenspaces satisfy
\begin{equation*} \begin{gathered}
  \bJ_i \jmult \bJ_i \leq \bJ_i, \quad
  \bJ_i \jmult \bJ_{ij} \leq \bJ_{ij}, \quad
  \bJ_i \jmult \bJ_j = \bJ_{ij} \jmult \bJ_{k\ell} = 0, \quad \\
  \bJ_{ij} \jmult \bJ_{jk} \leq \bJ_{ik}
  \quad \text{and} \quad
  \bJ_{ij} \jmult \bJ_{ij} \leq \bJ_i \dsum \bJ_j
\end{gathered}
\end{equation*}
for all distinct indices $i,j,k,\ell$. \noproof \end{cor}

The Peirce decomposition provides the crucial structural step in analysing the algebraic structure of the isotropy \Rspace\ $G\acts\fp$; this is the subject of the next subsection.

\vspace{0.3em}
\subsection{Structure of $\fg$ and $\bW$} 
\label{ss:ppg-alg-str}

The aim of this subsection is to describe the structure of $G\acts\fp$ and the $\fg$-representation $\bW$.  This is done in the following.

\begin{thm} \thlabel{thm:ppg-alg-str} Let $H\acts\fq$ be a \self dual\ symmetric \Rspace\ with infinitesimal isotropy representation $\bW \defeq \fh/\fq$ and isotropy \Rspace\ $G\acts\fp$.  Then:
\begin{enumerate}
  \item \label{thm:ppg-alg-str-g}
  $G\acts\fp$ is a symmetric \Rspace; and
  
  \item \label{thm:ppg-alg-str-W}
  The $\fp^{\perp}$-filtration of the $\fg$-representation $\bW$ induced by $\fh/\fq$ has height two.
\end{enumerate}
\end{thm}

The proof of \thref{thm:ppg-alg-str} will require some preliminary work involving the structure theory of $\fh$ and $\fg$, as well as their parabolic subalgebras $\fq$ and $\fp$.  It suffices to see that both $\fg$ and $\bW$ decompose into a direct sum of three eigenspaces \wrt\ an algebraic Weyl structure for $\fp$.  The graded properties of $\fg$ and $\bW$ are clearly unchanged by complexification, so we may assume that $H\acts\fq$ is complex.  Moreover by \thref{lem:ppg-alg-simple}, it suffices to consider the case when $\fh$ is simple.

\begin{rmk} \thlabel{rmk:ppg-alg-str} Using the classification of symmetric \Rspaces\ \cite{be1989-penrose}, one can easily check that \thref{thm:ppg-alg-str} does not require that $H\acts\fq$ is \self dual.  F.\ Burstall has recently provided a much simpler proof than the one presented here, that benefits from not requiring \self duality. \end{rmk}

Choose algebraic Weyl structures $\xi^{\fq}\in\liecenter{\fq^0}$ for $\fq\leq\fh$ and $\xi^{\fp} \in \liecenter{\fp^0}$ for $\fp\leq\fg$.  Then $\xi^{\fq}$ and $\xi^{\fp}$ are linearly independent, since otherwise $\xi^{\fp}$ would act trivially on the whole of $\fg$, implying that $\xi^{\fp} = 0$.  As in \thref{lem:lie-para-levi}, we can extend $\setof{\xi^{\fq},\xi^{\fp}}{}$ to a basis for a Cartan subalgebra $\ft$ with root system $\Delta \subset \ft^*$ and positive subsystem $\Delta^{+} \subset \Delta$ of $\fh$.  This identifies $\fq$ and $\fp$ with standard parabolic subalgebras of $\fh$ and $\fg$.  Evidently any two such Cartan subalgebras are conjugate by the action of $P^0$; since $P^0$ preserves the graded structure of $\fg$ and $\bW$, the choice of $\ft$ is inconsequential for our purposes.

\begin{conv} \thlabel{conv:ppg-alg-killing} We normalise the Killing form of $\fh$ so that $\killing{\alpha}{\alpha} = 2$ for all long roots $\alpha \in \Delta \subset \ft^*$. \end{conv}

Let $\Delta^0 \subset \Delta^{+}$ be the set of simple roots and let $\beta \in \ft^*$ be the corresponding highest root of $\fh$.  It is well-known that $\beta$ is a long root, whence $\killing{\beta}{\beta} = 2$ by \thref{conv:ppg-alg-killing}.  Since $\fh$ is simple, $\fq$ is the standard parabolic given by crossing a single simple root $\alpha_r \in \Delta^0$.  Writing $\beta = \sum{i}{} n_i \alpha_i$ as a sum of simple roots, by assumption we have $\height[\fq]{\beta} = n_r = 1$; equivalently, $\beta(\xi_r)=1$ for $\xi_1,\ldots,\xi_{\ell} \in \ft$ dual to $\alpha_1, \ldots, \alpha_{\ell} \in \Delta^0$.

As in the discussion following \thref{lem:lie-para-expopp}, the algebraic Weyl structure $\xi^{\fq}$ provides an isomorphism $\fh \isom \bW \dsum \fq^0 \dsum \bW^*$, with the grading given by $\fq$-height.  In particular $\height[\fq]{\alpha} \in \setof{+1,0,-1}{}$ for all $\alpha \in \Delta$, with $\height[\fq]{\beta} = 1$ and therefore $\fh_{\beta} \leq \bW$.

\begin{lem} \thlabel{lem:ppg-alg-arlong} \rmcite[Lem.\ 3.7]{br1990-twistorsymm}.
$\alpha_r$ is a long root.
\end{lem}

\begin{proof} Consider the dual root system $\Delta^*$, consisting of roots $\alpha^* \defeq \tfrac{2}{\killing{\alpha}{\alpha}} \alpha$ for $\alpha \in \Delta$.  A simple subsystem for $\Delta^*$ is then given by the duals $\alpha_i^*$ of the simple roots $\alpha_i \in \Delta^0$.  Write the highest root $\beta = \sum{i}{} n_i\alpha_i \in \Delta$. Therefore
\vspace{0.1em}
\begin{equation*}
  \beta^*
    = \frac{2}{\killing{\beta}{\beta}} \beta
    = \Sum{i}{} \frac{2}{\killing{\alpha_i}{\alpha_i}} \frac{n_i \killing{\alpha_i}{\alpha_i}}{\killing{\beta}{\beta}} \alpha_i
    = \Sum{i}{} \frac{n_i \killing{\alpha_i}{\alpha_i}}{\killing{\beta}{\beta}} \alpha_i^*,
\vspace{0.1em}
\end{equation*}
so that the $n_i \tfrac{\killing{\alpha_i}{\alpha_i}}{\killing{\beta}{\beta}}$ are integers. However by assumption we have $\height[\fq]{\beta} = n_r = 1$, so that we necessarily have $\alpha_r$ long. \end{proof}

Consider the fundamental weight $\omega_r$ corresponding to $\alpha_r$. Then if $\alpha = \sum{i}{} n_i \alpha_i \in \Delta$ is any root, $\killing{\omega_r}{\alpha} = \sum{i}{} n_i \killing{\omega_r}{\alpha_i} = \tfrac{1}{2} n_r \killing{\alpha_r}{\alpha_r} = n_r = \height[\fq]{\alpha}$. Therefore the $\fq$-height is given by the inner product with $\omega_r$.

Recall that roots $\alpha,\beta \in \Delta$ are \emph{strongly orthogonal} if neither of $\alpha\pm\beta$ are roots.  In particular, strongly orthogonal roots $\alpha,\beta$ generate commuting \sltriple s\ \cite{ak2002-strongorth}.

\begin{lem} \thlabel{lem:ppg-alg-strorth} Roots $\alpha, \beta \in \Delta$ of $\fq$-height one are strongly orthogonal if and only if they are orthogonal. \end{lem}

\begin{proof} It is well-known that strong orthogonality implies orthogonality: indeed, if neither of $\alpha \pm \beta$ are roots then the $\beta$-root string through $\alpha$ consists solely of $\alpha$, so we must have $\killing{\alpha}{\beta} = 0$.  Conversely, suppose that $\alpha, \beta$ are orthogonal.  Then $\alpha + \beta$ cannot be a root since it would have $\fq$-height two; the $\beta$-root string through $\alpha$ then consists solely of $\alpha$ again, from which $\killing{\alpha}{\beta} = 0$ implies that $\alpha - \beta$ is not a root either. \end{proof}

We are going to construct a maximal sequence $(\beta_n,\ldots,\beta_0)$ of orthogonal long roots of $\fq$-height one.  We will fix a particular sequence shortly; for now, take $(\beta_n, \ldots, \beta_0)$ to be any maximal sequence of such roots.  In particular, the $\beta_i$ are mutually strongly orthogonal by \thref{lem:ppg-alg-strorth} and satisfy $\killing{\beta_i}{\beta_j} = 2 \delta_{ij}$.

\begin{asmpt} \thlabel{asmpt:ppg-alg-n>0} We assume henceforth that $n>0$. \end{asmpt}

We will interpret this later as excluding \ppgs\ over zero-dimensional manifolds; in the mean time, the reader may consult the tables in \cite[\S3--5]{ak2002-strongorth} for reassurance.
To construct our particular choice of $(\beta_n, \ldots, \beta_0)$, we need the following lemma from \cite[\S 2]{rrs1992-abelianrad}; we shall need details from the proof, so include it for the reader's convenience.  Let $\sW_{\fg}$ be the Weyl group of $\fg$, viewed as a subgroup of the Weyl group $\sW_{\fh}$ of $\fh$ via the algebraic Weyl structure $\xi^{\fq}$.

\begin{lem} \thlabel{lem:ppg-alg-trans} Let $\sW_{\fg}$ be the Weyl group of $\fg$.  Then:
\begin{enumerate}
  \item \label{lem:ppg-alg-trans-one}
  $\sW_{\fg}$ acts transitively on the set of $\fq$-height one long roots.
    
  \item \label{lem:ppg-alg-trans-seq}
  For each $n\geq 0$, $\sW_{\fg}$ acts transitively on all orthogonal sequences $(\beta_n, \ldots, \beta_0)$ of $\fq$-height one roots of a fixed length.
  
  \item \label{lem:ppg-alg-trans-len}
  Every maximal orthogonal sequence of $\fq$-height one long roots has the same length.
\end{enumerate}
\end{lem}

\begin{proof} \proofref{lem:ppg-alg-trans}{one}
This is a special case of \cite[Lem.\ 2.6]{rrs1992-abelianrad}.  Let $\alpha \in \Delta^0\setminus\{\alpha_r\}$ be a $\fq$-height zero simple root.  Then for any $\beta \in \Delta$ with $\height[\fq]{\beta} = 1$, we have
\vspace{-0.15em}
\begin{equation*}
  \killing{ \refl{\alpha}{\beta} }{ \omega_r }
    = \killing{ \beta - \cartanint{\beta}{\alpha}\alpha }{ \omega_r }
    = \killing{\beta}{\omega_r} - \cartanint{\beta}{\alpha} \killing{\alpha}{\omega_r}
    = 1,
\vspace{-0.15em}
\end{equation*}
so that $\refl{\alpha}{\beta}$ is also of $\fq$-height one.  Since these simple reflections generate the Weyl group $\sW_{\fg}$ of $\fg$ (viewed as a subgroup of $\sW_{\fh}$), it follows that $\sW_{\fg}$ preserves the set of $\fq$-height one roots.  Now consider two orbits, represented by roots $\beta,\beta'$ which we suppose are chosen to be maximal \wrt\ the root partial ordering.  Then for any $\alpha \in \Delta^0\setminus\{\alpha_r\}$ we have that $\refl{\alpha}{\beta} = \beta - \cartanint{\beta}{\alpha}\alpha$ is a root of $\fq$-height one, so that $\killing{\beta}{\alpha} \geq 0$ by maximality of $\beta$ within its orbit.  Moreover $\killing{\beta}{\alpha_r} \geq 0$ since otherwise $\refl{\alpha_r}{\beta} = \beta + \alpha_r$ would be a root of $\fq$-height two.  Therefore $\killing{\beta}{\alpha} \geq 0$ for all $\alpha \in \Delta^0$, so that $\beta$ is dominant.  Similarly $\beta'$ is dominant; but there is a unique dominant long root, so that the two orbits must coincide.

\smallskip

\proofref{lem:ppg-alg-trans}{seq}
This is \cite[Prop.\ 2.8]{rrs1992-abelianrad}.  By \ref{lem:ppg-alg-trans-one} we may conjugate the sequence by an element of $\sW_{\fg}$ and assume that $\beta_n = \beta$, the highest root of $\fh$.
If $\killing{\beta_n}{\alpha_r} \neq 0$ then since $\beta_n + \alpha_r$ is not a root, we must have $\killing{\beta_n}{\alpha_r} > 0$.  Since $\beta_n$ is a dominant weight for $\fh$, it is a \non negative\ integral linear combination of the fundamental weights.  Then $\killing{\beta_n}{\alpha_r} > 0$ implies that the coefficient of $\omega_r$ is positive, so that $\killing{\beta_n}{\alpha} > 0$ for all positive roots $\alpha \in \Delta^{+}$.  However any root of $\fq$-height one is necessary positive, having a coefficient one of the simple root $\alpha_r$, so since $\killing{\beta_n}{\beta_{n-1}} = 0$ we must have $\killing{\beta_n}{\alpha_r} = 0$.

\vspace{0.15em}
Now $\alpha_r$ and the remaining $\beta_i$ lie in the root subsystem $\linspan{\beta_n}{}^{\perp} \subset \Delta$.  Although this may not be irreducible, its irreducible component containing $\alpha_r$ also contains the remaining $\beta_i$ since they have $\fq$-height one.  Writing $\Delta_0$ for the roots of $\fg$, the Weyl group of $\Delta_0 \intsct \linspan{\beta_n}{}^{\perp}$, viewed as a subgroup of $\sW_{\fg}$, acts transitively on the $\fq$-height one roots in $\linspan{\beta_n}{}^{\perp}$.  Therefore we may assume that $\beta_{n-1}$ is the highest root of the irreducible component of $\linspan{\beta_n}{}^{\perp}$ containing $\alpha_r$.  Proceeding inductively, we see that $\sW_{\fg}$ acts transitively as required.

\smallskip\vspace{0.2em}

\proofref{lem:ppg-alg-trans}{len}
This follows immediately from \ref{lem:ppg-alg-trans-seq}. \end{proof}

\vspace{0.2em}

We may therefore assume that our maximal sequence $(\beta_n, \ldots, \beta_0)$ is constructed as follows: $\beta_n = \beta$, the highest root of $\fh$, and we inductively define $\beta_i$ to be the highest root of the irreducible component of $\linspan{\beta_n, \ldots, \beta_{i+1}}{}^{\perp}$ containing $\alpha_r$.  This construction is familiar in the context of hermitian symmetric spaces \cite[Eqn.\ (3.2)]{w1972-hermitiansymm}, and apparently is originally due to Harish-Chandra.

\vspace{0.2em}

\begin{lem} \thlabel{lem:ppg-alg-b0ar} $\beta_0 = \alpha_r$.
\vspace{0.2em}
\end{lem}

\begin{proof} Since $H\acts\fq$ is \self dual, the standard parabolic $\fq$ is conjugate to the standard opposite parabolic $\opp{\fq}$, which occurs if and only if the corresponding parabolic normalisers $Q \defeq \Norm[H]{\fq}$ and $\opp{Q} \defeq \Norm[H]{\opp{\fq}}$ of $H$ are conjugate.  There is a natural isomorphism of the Weyl group $\sW_{\fh}$ of $\fh$ with the quotient $\Norm[H]{\ft} / \grpcenter[H]{\ft}$ of the adjoint normaliser of the Cartan subalgebra $\ft$ by the centraliser of $\ft$ \cite[Thm.\ 3.2.19(1)]{cs2009-parabolic1}, and it follows from the Bruhat decomposition \cite[\S 6.4]{p2007-lieinv} that $Q, \opp{Q}$ are in fact conjugate by an element $w \in \sW_{\fh}$; see \itemref{rmk:ppg-alg-bruhat}{weyl} below.  It is proved in \cite[Prop.\ 3.12]{rrs1992-abelianrad} that we may take $w$ to be the element $w_n \defeq \refl{\beta_0}{} \circ \cdots \circ \refl{\beta_n}{}$.  In particular, $w_n$ takes roots of $\fq$-height $k \in \setof{+1, 0, -1}{}$ to roots of $\fq$-height $-k$.

\vspace{0.15em}
Since $\alpha_r \leq \beta_0$ in the root partial ordering, if $\beta_0 \neq \alpha_r$ then maximality of the sequence $(\beta_n, \ldots, \beta_0)$ implies that $\killing{\alpha_r}{\beta_0} \neq 0$.  We must then have $\killing{\alpha_r}{\beta_0} > 0$, since $\beta_0 + \alpha_r$ is not a root.  Therefore the Cartan integer $\killing{\alpha_r}{\beta_0} = 1$; indeed, $\killing{\alpha_r}{\beta_0} \in \setof{1,2}{}$ and $\killing{\alpha_r}{\beta_0} = 2$ would imply that $\beta_0 = \alpha_r$ by the Cauchy--Schwarz inequality.  In particular, $\refl{\alpha_r}{\beta_0} = \beta_0 - \alpha_r$ is a positive root of $\fq$-height zero.  Since $\killing{\alpha_r}{\beta_i} = 0$ for $i>0$ by construction, it follows that
\vspace{-0.3em}
\begin{equation*}
  w_n(\alpha_r) \defeq \alpha_r - \Sum{i=0}{n} \killing{\beta_i}{\alpha_r} \beta_i
    = \alpha_r - \beta_0,
\vspace{-0.3em}
\end{equation*}
so that $w_n$ takes a root of $\fq$-height one to a root of $\fq$-height zero.  But this contradicts \self duality\ of $H\acts\fq$ by the first paragraph, implying that $\beta_0 = \alpha_r$ as required. \end{proof}

\begin{rmk} \thlabel{rmk:ppg-alg-bruhat}
\begin{enumeratepar}
  \item \label{rmk:ppg-alg-bruhat-nosd}
  One can check that $\beta_0 \neq \alpha_r$ when $H\acts\fq$ is not \self dual.  As an example, crossing node $k\leq n+1$ in $\alg{sl}{2n+2,\bC}$ yields a maximal orthogonal sequence given by $\beta_i = \alpha_{k-i} + \cdots + \alpha_{2n+2-k+i}$ for $i = k-1, \ldots, 0$.  In particular, $\beta_0 = \alpha_k + \cdots + \alpha_{2n+2-k} \neq \alpha_k$ unless $k=n+1$, \ie\ unless we cross the central node.
  
  \item \label{rmk:ppg-alg-bruhat-weyl}
  More carefully, it suffices to see that $w$ and $w_n$ lie in the same coset of the Cartan subgroup $\grpcenter[H]{\ft}$ in $\sW_{\fh} \isom \Norm[H]{\ft} / \grpcenter[H]{\ft}$.  This is because elements of $\grpcenter[H]{\ft}$ act on $H$ preserving the grading, thus having no effect on the conjugacy between $Q$ and $\opp{Q}$. \end{enumeratepar}
\end{rmk}

Choose an \sltriple\ $\setof{e_i,h_i,f_i}{}$ corresponding to each $\beta_i$, where $e_i \in \fh_{\beta_i}\leq\bW$, $f_i \in \fh_{-\beta_i}\leq\bW^*$ and $h_i \defeq \liebrac{e_i}{f_i} \in \ft$.  By \thref{conv:ppg-alg-killing}, the corresponding \co root\ is $\coroot{\beta_i} \defeq \tfrac{2}{\killing{\beta_i}{\beta_i}} h_i = h_i$.  Moreover since the $\beta_i$ are strongly orthogonal, we have $\liebrac{e_i}{f_j} = \delta_{ij}h_i$ for all $i,j$.  It follows that
\vspace{-0.3em}
\begin{equation*}
  e \defeq \Sum{i=0}{n} e_i \in \bW, \quad
  f \defeq \Sum{i=0}{n} f_i \in \bW^*,
  \quad\text{and}\quad
  h \defeq \Sum{i=0}{n} h_i \in \ft
\vspace{-0.3em}
\end{equation*}
satisfies $h=\liebrac{e}{f}$ and hence $\{e,h,f\}$ is also an \sltriple\ for $\fh$.  Thus the $\{e_i, h_i, f_i\}$ generate commuting subalgebras isomorphic to $\alg{sl}{2,\bC}$, with the subalgebra generated by $\{e,h,f\}$ the diagonal subalgebra in their direct sum.

\begin{lem} \thlabel{lem:ppg-alg-isom} $F_f \defeq (\ad f)^2 : \bW \to \bW^*$ is a linear isomorphism. \end{lem}

\begin{proof} We will show that the linear map $\tfrac{1}{4}(\ad e)^2 : \bW^* \to \bW$ is inverse to $F_f$.  For $\alpha\in\Delta$ a root of $\fq$-height one and all $x\in\fh_{\alpha} \leq \bW$, the Jacobi identity yields
\begin{align} \label{eq:ppg-alg-isom-1}
  (\ad e)^2(\ad f)^2(x) &= \liebrac{ e }{ \liebrac{ e }{ \liebrac{ f }{ \liebrac{f}{x} }}} \notag \\
  &= \liebrac{e}{\liebrac{h}{\liebrac{f}{x}}} + \liebrac{e}{\liebrac{f}{\liebrac{h}{x}}} \notag \displaybreak \\
  &= -2\liebrac{e}{\liebrac{f}{x}} + 2\liebrac{h}{\liebrac{h}{x}} \notag \\
  &= - 2\liebrac{h}{x} + 2\liebrac{h}{\liebrac{h}{x}} \notag \\
  &= 2\alpha(h) \left( \alpha(h)-1 \right)x.
\end{align}
Since $\killing{\beta_i}{\beta_j} = 2\delta_{ij}$, the root reflections $\refl{\beta_i}{}$ through $\beta_i$ mutually commute and hence
\begin{equation*}
  w_n(\alpha) \defeq
    (\refl{\beta_0}{} \circ\cdots\circ \refl{\beta_n}{})(\alpha)
    = \alpha - \Sum{i=0}{n} \killing{\alpha}{\beta_i} \beta_i
\end{equation*}
is a root of $\fq$-height $1-\Sum{i=0}{n} \killing{\alpha}{\beta_i} \in \setof{+1,0,-1}{}$.  Since the $\beta_i$ are all long roots, $\killing{\alpha}{\beta_i} = \cartanint{\alpha}{\beta_i}$ is integral for all $i = n, \ldots, 0$.  Therefore either $\killing{\alpha}{\beta_i} = 2$ for a unique $i$, with $\killing{\alpha}{\beta_j} = 0$ for $j\neq i$; or $\killing{\alpha}{\beta_i} = 1$ for at most two $i$, with $\killing{\alpha}{\beta_j} = 0$ otherwise; or $\killing{\alpha}{\beta_i} = 0$ for all $i$.  The last case cannot occur, since by \thref{lem:ppg-alg-b0ar} the root subsystem $\linspan{\beta_n, \ldots, \beta_0}{}^{\perp} \subseteq \linspan{\alpha_r}{}^{\perp}$ contains no roots of $\fq$-height one.  The case where $\killing{\alpha}{\beta_i} = 1$ for a unique $i$ also cannot occur, since this contradicts \self duality\ of $H\acts\fq$ as in the proof of \thref{lem:ppg-alg-b0ar}.  The remaining two cases both give $\alpha(h) = \sum{i=0}{n} \killing{\alpha}{\beta_i} = 2$, yielding $(\ad e)^2 \circ (\ad f)^2 = 4\, \id$ by \eqref{eq:ppg-alg-isom-1}. \end{proof}

Note that \thref{lem:ppg-alg-isom} gives an alternative proof of the implication \ref{prop:ppg-defn-sd-sd}$\Rightarrow$\ref{prop:ppg-defn-sd-ker} from \thref{prop:ppg-defn-sd}.  Indeed, we have explicitly constructed a regular element $f \in \fq^{\perp}$.

\begin{cor} \thlabel{cor:ppg-alg-xiq} $h = \liebrac{e}{f}$ is equal to twice the algebraic Weyl structure $\xi^{\fq}$ of $\fq$. \end{cor}

\begin{proof} \thref{lem:ppg-alg-isom} equivalently states that $f \in \bW^*$ is a regular element.  Since $F_f(e) = \liebrac{ \liebrac{e}{f} }{ f } = \liebrac{h}{f} = -2f$, we obtain that $e = -2F_f^{-1}(f)$ is the corresponding element from \itemref{prop:ppg-defn-sd}{lift} with $\liebrac{e}{f} = 2\xi^{\fq}$. \end{proof}

By Meyberg's \thref{thm:ppg-alg-Wjordan}, the \sltriple\ $\setof{e,h,f}{}$ induces dual Jordan algebra structures on $\bW$ and $\bW^*$, defined by
\begin{equation*}
  x \jmult y \defeq \tfrac{1}{2} \liebrac{ \liebrac{x}{f} }{ y }
  \quad\text{and}\quad
  \alpha \jmult \beta \defeq \tfrac{1}{2} \liebrac{ \liebrac{\alpha}{e} }{ \beta }
\end{equation*}
for all $x,y \in \bW$ and $\alpha,\beta \in \bW^*$.  Since we have chosen the $\beta_i$ with root spaces contained in $\bW$, we shall work predominantly with the Jordan multiplication on $\bW$.

\begin{lem} The $e_i \in \bW$ are mutually orthogonal Jordan idempotents. \end{lem}

\begin{proof} \thref{lem:ppg-alg-strorth} gives $e_i \jmult e_j = \tfrac{1}{2} \liebrac{ \liebrac{e_i}{f} }{ e_j } = \tfrac{1}{2} \liebrac{h_i}{e_j} = \delta_{ij} e_i$ for all $i,j$. \end{proof}

By \thref{thm:ppg-alg-peirce} we consequently have a Peirce decomposition%
\footnote{The idea to utilise the Peirce decomposition and strongly orthogonal roots comes from \cite{ks2015-smallrepns}.}
\begin{equation} \label{eq:ppg-alg-peirce}
  \bW = \Dsum{i=0}{n} \Big( \bW_i \dsum \mathsmaller{\Dsum{j>i}{}} \bW_{ij} \Big),
\end{equation}
where the Jordan eigenspaces are defined by
\begin{equation} \label{eq:ppg-alg-peircecmpt} \begin{aligned}
  \bW_i &\defeq \setof{ x\in\bW }{ e_i \jmult x = x } \\
  \text{and}\quad
  \bW_{ij} &\defeq \setof{ x\in\bW }{ e_i \jmult x = \tfrac{1}{2}x = e_j \jmult x }.
\end{aligned} \end{equation}
In terms of root data, if $x \in \fh_{\alpha} \leq \bW_i$ then $x = e_i \jmult x = \tfrac{1}{2} \liebrac{h_i}{x} = \tfrac{1}{2} \killing{\alpha}{\beta_i} x$, giving $\killing{\alpha}{\beta_i} = 2$ and hence $\alpha=\beta_i$ by the Cauchy--Schwarz inequality.  On the other hand for $x \in \fh_{\alpha} \leq \bW_{ij}$, a similar line of argument gives $\killing{\alpha}{\beta_i} = \killing{\alpha}{\beta_j} = 1$.  Combined with the multiplication properties of \thref{cor:ppg-alg-idemmult}, it follows that
\begin{align*}
  \bW_i &= \fh_{\beta_i} = \linspan{ e_i }{} \\
  \text{and}\quad
  \bW_{ij} &= \linspan{ \fh_{\alpha} }{ \killing{\beta_i}{\alpha} = \killing{\beta_j}{\alpha} = 1
                                        ~\text{and}~
                                        \killing{\beta_k}{\alpha} = 0 ~\, \forall k\neq i,j }.
\end{align*}
In particular, $\bW_i$ is $1$-dimensional spanned by $e_i$.  There is a dual Peirce decomposition $\bW^* = \Dsum{i=0}{n} \big( \bW^*_i \dsum \mathsmaller{\Dsum{j>i}{}} \bW^*_{ij} \big)$ defined in much the same way.

Now we turn to describing the root data of the semisimple part $\fg$ of $\fq^0$.  The roots of $\fg$ are precisely the roots of $\fh$ of $\fq$-height zero, with simple subsystem $\Delta^0 \setminus \{\alpha_r\}$.  By \thref{lem:lie-para-levi}, the corresponding Cartan subalgebra is
\begin{equation} \label{eq:ppg-alg-gcsa}
  \ft_0 \defeq \linspan{\, \coroot{\alpha} }{ \alpha \in \Delta^0\setminus\{\alpha_r\} \,}
    = \Setof{ H\in\ft }{ \omega_r(H) = 0 }.
\end{equation}

The isotropy \Rspace\ $G\acts\fp$ is, by definition, induced by the stabiliser of the lowest weight orbit in $\bW$.  Identifying $\bW = \fh/\fq = \liehom[\fq^{\perp}]{0}{\bW}$ using \thref{cor:lie-hom-h0g}, the highest weight of $\fg$ on $\bW$ is the restriction $\beta_n\at{\ft_0}$ of the highest root of $\fh$ to $\ft_0$.  Note that $\beta_n\at{\ft_0} = 0$ if and only if $\killing{\beta_n}{\alpha} = 0$ for all $\alpha \in \Delta^0 \setminus \{\alpha_r\}$; since $n>0$, we have $\killing{\beta_n}{\alpha_r} = 0$ automatically, and thus $\beta_n\at{\ft_0} \neq 0$.  Moreover $\alpha_r\at{\ft_0} \neq 0$: the expansion of $\alpha_r$ in fundamental weights is $\alpha_r = 2\omega_r + \sum{i\neq r}{} \cartanint{\alpha_r}{\alpha_i} \omega_i$, which has a \non zero\ coefficient for at least one other $\omega_i$ by connectedness of the Dynkin diagram of $\fh$.

\begin{lem} \thlabel{lem:ppg-alg-wgar} Let $\longest[\fg]$ be the longest element of $\sW_{\fg}$.  Then $\longest[\fg](\beta_n) = \alpha_r$. \end{lem}

\begin{proof}  Since $\alpha_r$ and $\beta_n$ are both long roots of $\fq$-height one, \itemref{lem:ppg-alg-trans}{one} provides an element $w \in \sW_{\fg}$ such that $w(\beta_n) = \alpha_r$.  Then $w \leq \longest[\fg]$ in the Bruhat order on $\sW_{\fg}$, implying that $\alpha_r = w(\beta_n) \geq \longest[\fg](\beta_n)$ since $\beta_n$ is dominant.  But by definition $\alpha_r$ is the least root of $\fq$-height one, so that $\longest[\fg](\beta_n) = \alpha_r$ as required. \end{proof}

Since the highest weight of the $\fg$-representation $\bW^*$ is $-\longest[\fg](\beta_n\at{\ft_0})$, it follows from \thref{prop:lie-para-proj} that the standard parabolic $\fp \leq \fg$ is associated to the subset
\begin{equation} \label{eq:ppg-alg-psigma}
  \Sigma \defeq \Setof{ \alpha \in \Delta^0\setminus\{\alpha_r\} }
                      { \killing{\alpha}{\alpha_r} \neq 0 }
    \subseteq \Delta^0
\end{equation}
of $\fq$-height zero simple roots.

As a final ingredient, consider the Lie algebra automorphism $\sigma : \fh \to \fh$ induced by the \sltriple\ $(e_0, h_0, f_0)$, defined by
\vspace{0.1em}
\begin{equation} \label{eq:ppg-alg-sigma}
  \sigma(y) \defeq \exp(e_0)\exp(-f_0)\exp(e_0) \acts y.
\vspace{0.1em}
\end{equation}
The inverse is given explicitly by $\sigma^{-1}(y) = \exp(-e_0)\exp(f_0)\exp(-e_0) \acts y$.

\begin{lem} \thlabel{lem:ppg-alg-sigma} \emph{\cite[Prop.\ 2.33]{cn2015-paraproj}} Consider the automorphism $\sigma$ defined by \eqref{eq:ppg-alg-sigma}.  Then:
\begin{enumerate}
  \item \label{lem:ppg-alg-sigma-refl}
  $\sigma$ preserves $\ft$, where it restricts to the root reflection through $\beta_0 = \alpha_r$.
  
  \item \label{lem:ppg-alg-sigma-root}
  For every $\alpha\in\Delta$, we have $\sigma(\fh_{\alpha}) = \fh_{\sigma(\alpha)}$.
\end{enumerate}
\end{lem}

\begin{proof} \proofref{lem:ppg-alg-sigma}{refl}
The automorphism $\sigma$ restricts to the identity on $\ker\beta_0$ and maps $\fh_{\beta_0}$ to $\fh_{-\beta_0}$, so that $\sigma(H) = H - \beta_0(H) h_0 = \refl{\beta_0}{H}$ for all $H\in\ft$.

\smallskip

\proofref{lem:ppg-alg-sigma}{root} For all $x\in\fh_{\alpha}$ and $H\in\ft$ we have
\begin{align*}
  \liebrac{H}{\sigma(x)}
    &= \sigma( \liebrac{\sigma^{-1}(H)}{x} )
     = \sigma( \liebrac{\sigma(H)}{x} ) \\
    &\hspace{4em}
     = \sigma( \alpha(\sigma(H)) x )
     = (\sigma(\alpha))(H) \sigma(x),
\end{align*}
giving $\sigma(\fh_{\alpha}) \subseteq \fh_{\sigma(\alpha)}$.  Since both root spaces are $1$-dimensional and $\sigma$ is an automorphism, this is an equality. \end{proof}

We are finally in a position to prove \thref{thm:ppg-alg-str}.

\begin{proofof}{thm:ppg-alg-str} It suffices to describe the eigenspaces of the algebraic Weyl structure $\xi^{\fp}$ of $\fp$, so our first task is to identify $\xi^{\fp}$.  We first calculate how $\sigma(\xi^{\fq}) = 2\sigma(h)$ acts on root spaces $\fh_{\alpha}$ contained in $\fg$.  Then $\height[\fq]{\alpha} = 0$, so we may write $\alpha = \Sum{i\neq r}{} a_i\alpha_i$ for some $a_i \in \bZ$ with the same signs.  For all $x\in\fh_{\alpha}$ we have
\vspace{0.15em}
\begin{equation*} \begin{aligned}
  \liebrac{ \sigma(h) }{ x }
  &= \sigma( \liebrac{h}{\sigma^{-1}(x)} ) \\
  &= \sigma( (\sigma^{-1}(\alpha))(h) \sigma^{-1}(x) ) \\
  &= \alpha( h - 2h_0 ) x \\
  &= \alpha(h)x - 2\sum{i\neq r}{} \, a_i \killing{\alpha_r}{\alpha_i} x
\end{aligned}
\vspace{0.15em}
\end{equation*}
by \thref{lem:ppg-alg-b0ar}.  The first term here vanishes, since it equals $\alpha(h)x = \liebrac{h}{x} = 2(\height[\fq]{\alpha})x = 0$ by \thref{cor:ppg-alg-xiq}.  Since $\alpha_i + k\alpha_r$ is not a root for $k>1$, it follows from \eqref{eq:ppg-alg-psigma} that $\killing{\alpha_r}{\alpha_i} = -1$ for all $\alpha_i \in \Sigma$.  On the other hand, $\killing{\alpha_r}{\alpha_i} = 0$ for $\alpha_i \in \Delta^0 \setminus (\{\alpha_r\} \union \Sigma)$, and the above gives $\liebrac{ \sigma(h) }{ x } = 2\left( \sum{i \,:\, \alpha_i \in \Sigma}{} \, a_i \right) x = 2(\height[\fp]{\alpha}) x$.  Therefore $\sigma(h)$ acts by twice the $\fp$-height on root spaces contained in $\fg$.  Since $h \in \ft \leq \fg\dsum\liecenter{\fq^0}$, we conclude that $\sigma(h) \in \liecenter{\fp^0} \dsum \liecenter{\fq^0}$.  In particular since $\liecenter{\fq^0}$ is $1$-dimensional and spanned by $\xi^{\fq}$, we have $\xi^{\fp} = \sigma(\xi^{\fq}) - k\xi^{\fq}$ for some $k\in\bR$ for which $\xi^{\fp} \in \ft_0$.  By the description \eqref{eq:ppg-alg-gcsa} of $\ft_0$, we require that $k$ satisfies
\begin{align*}
  0 = \omega_r(\xi^{\fp})
  &= \omega_r(\sigma(\xi^{\fq}) - k\xi^{\fq}) \\
  &= \omega_r( (1-k)\xi^{\fq} - \beta_0(\xi^{\fq})h_0 ) \\
  &= \tfrac{1}{2}\omega_r\left( (1-k)(h_n + \cdots + h_0)
    - \beta_0(h_n + \cdots + h_0)h_0 \right) \\
  &= \tfrac{1}{2}(1-k)(n+1) - 1.
\end{align*}
Therefore $k = \tfrac{n-1}{n+1}$ and consequently
\begin{equation} \label{eq:ppg-alg-xip}
  \xi^{\fp} = \sigma(\xi^{\fq}) - \tfrac{n-1}{n+1} \xi^{\fq}
    = \tfrac{2}{n+1}\xi^{\fq} - h_0
\end{equation}
We can now complete the proof.

\smallskip

\proofref{thm:ppg-alg-str}{g} Since $\xi^{\fq}$ acts trivially on $\fg \leq \fq^0$, we have $\liebrac{\xi^{\fp}}{x} = \sigma(\liebrac{\xi^{\fq}}{\sigma^{-1}(x)}) = \height[\fq]{\sigma(\alpha)} x$ and thus $\height[\fp]{\alpha} = \height[\fq]{\sigma(\alpha)}$ for all $\alpha\in\Delta$ with $\fh_{\alpha} \leq \fg$.  Since $\fq$ is an abelian parabolic, it remains to prove that each eigenspace occurs.  For this, note that $\fp^0 \leq \fq^0$ has $\fp$-height zero, while $\sigma(\bW_{0i})$ and $\sigma(\bW^*_{0i})$ consist respectively of root spaces of $\fp$-height $\pm 1$.

\smallskip

\proofref{thm:ppg-alg-str}{W} It remains to identify the eigenspaces of $\xi^{\fp}$ on $\bW$.  Given a root $\alpha$ of $\fq$-height one, the descriptions \eqref{eq:ppg-alg-peircecmpt} give $\fh_{\alpha} \leq \bW_i$ if and only if $\alpha = \beta_i$, so that $\sigma(\alpha) = \alpha$ for $i \neq 0$ and $\sigma(\alpha) = -\alpha$ for $i = 0$.  On the other hand, $\fh_{\alpha} \leq \bW_{ij}$ if and only if $\killing{\alpha}{\beta_i} = 1 = \killing{\alpha}{\beta_j}$ and $\killing{\alpha}{\beta_k} = 0$ for all $k \neq i,j$, so that $\sigma(\alpha) = \alpha$ for $j \neq 0$ and $\sigma(\alpha) = \alpha - \beta_0$ for $j = 0$.  It follows that $\sigma(\xi^{\fq})$ acts on $\bW_i$ by the identity for $i>0$, and by $-1$ for $i=0$; while it acts on $\bW_{ij}$ by $+1$ for $j>i>0$, and trivially for $j>i=0$.  Since $\xi^{\fq}$ acts on all root spaces contained in $\bW$ by the identity, we conclude from above that
\begin{equation} \label{eq:ppg-alg-eigenvals} \begin{alignedat}{6}
  \liebrac{\xi^{\fp}}{x} &= \phantom{-}\tfrac{2}{n+1}x &\qquad&
    \forall x \in \Dsum{i=1}{n} \Big( \bW_i \dsum \mathsmaller{\Dsum{j>i}{}} \bW_{ij} \Big) \\
  \liebrac{\xi^{\fp}}{x} &= -\tfrac{n-1}{n+1}x &\qquad&
    \forall x \in \Dsum{i=1}{n} \bW_{0i} \\
  \text{and}\qquad
  \liebrac{\xi^{\fp}}{x} &= -\tfrac{2n}{n+1}x &\qquad&
    \forall x \in \bW_0,
\end{alignedat} \end{equation}
which is the eigenspace decomposition of $\bW$ as a $\fp^0$-representation.  Therefore $\bW$ has three graded components, hence height two as a $\fg$-representation. \end{proofof}

\begin{rmk} We observed above that a simple root $\alpha \in \Delta^0\setminus\{\alpha_r\}$ lies in $\Sigma$ if and only if $\killing{\alpha_i}{\alpha_r} = -1$.  Since both $\alpha$ and $\alpha_r$ are simple roots, $\alpha-\alpha_r$ is not a root and hence $\killing{\alpha}{\alpha_r} \leq 0$.  Moreover since $\fq$ is an abelian parabolic $\alpha+2\alpha_r$ is not a root.  Consequently $\killing{\alpha}{\alpha_r}=-1$ is the only \non zero\ option.  In particular, to form $\fp$ we cross the nodes of the Dynkin diagram of $\fg$ that were connected to crossed nodes in the Dynkin diagram of $\fq$.  This makes the classification in Section \ref{s:ppg-class} very straightforward. \end{rmk}

\smallskip
\subsection{The $\bZ^2$-grading} 
\label{ss:ppg-alg-Z2}

\thref{thm:ppg-alg-str} describes the structure of $\fg$ and $\bW$ as $\fp^0$-representations, \wrt\ algebraic Weyl structures for $\fq$ and $\fp$.  We can use this to describe $\fh$ as a $\fp^0$-representation, and the Lie brackets between various irreducible components.  For this, we first identify the summands of $\bW$ with certain $\fp^0$-representations.

Consider an irreducible $\fg$-representation $\bV$ whose $\fp^{\perp}$-filtration has height $k$ and socle $\bV_0$.  Since $\liehom{0}{\bV} \defeq \bV / (\fp^{\perp}\acts\bV)$, there is a natural linear map
\begin{equation} \label{eq:ppg-alg-hommap} \begin{aligned}
  \Psi : \liehom{0}{\bV} &\to \Hom{\Tens{k}\fp^{\perp}}{\bV_0} \\
  \Psi[v] (\alpha_1 \tens\cdots\tens \alpha_k) &\defeq \alpha_1\acts (\alpha_2 \,\cdots\, (\alpha_k \acts v)).
\end{aligned} \end{equation}
Then $\Psi$ is well-defined since $\bV_0$ is, by definition, the kernel of the $\fp^{\perp}$-action; moreover since $\fp^{\perp}$ is abelian, $\Psi$ takes values in $\Hom{\Symm{k}\fp^{\perp}}{\bV_0}$.  In the case that $\bV$ is an irreducible $\fg$-representation, $\liehom{0}{\bV}$ is an irreducible $\fp$-representation and therefore by Schur's lemma $\Psi$ is an isomorphism onto its image, which is an irreducible $\fp$-subrepresentation of $\Hom{\Symm{k}\fp^{\perp}}{\bV_0}$.  Applying this to the case $\bV = \bW$, we deduce the following.

\begin{prop} \thlabel{prop:ppg-alg-gr} Let $L^* \defeq \bW_0$ be the $\fp^{\perp}$-socle of $\bW$.  Then \wrt\ any algebraic Weyl structure for $\fp$, we have
\vspace{0.2em}
\begin{equation*}
  \bW \isom (L^*\tens B) \dsum (L^* \tens \fg/\fp) \dsum L^*
\vspace{0.2em}
\end{equation*}
for a $\fp$-subrepresentation $B \leq \Symm{2}(\fg/\fp)$. \end{prop}

\begin{proof} By \thref{lem:ppg-alg-simple}, it suffices to suppose that $\fh$ is simple, in which case $\bW$ is an irreducible $\fg$-representation and $L^*$ is $1$-dimensional by the theory of Subsection \ref{ss:ppg-alg-str}.  By \itemref{thm:ppg-alg-str}{W} there are three graded components of $\bW$, and hence $k=2$ in \eqref{eq:ppg-alg-hommap}.  Then by definition, the lowest and highest weight summands are $L^*$ and $\liehom{0}{\bW} \defeq \bW/(\fp^{\perp}\acts\bW)$ respectively, and using the map $\Psi$ from \eqref{eq:ppg-alg-hommap} we identify $\liehom{0}{\bW}$ with $L^*\tens B \leq \Hom{\Symm{2}\fp^{\perp}}{L^*}$ for some irreducible $\fp$-subrepresentation $B \leq \Symm{2}(\fg/\fp)$.  Thus it remains to identify the graded component $\bW_{(1)} \defeq (\fp^{\perp}\acts\bW) / L^*$.

Since $\fp$ is an abelian parabolic by \itemref{thm:ppg-alg-str}{g}, the action of $\fp^{\perp}$ in the filtration $\bW \supset \fp^{\perp}\acts\bW \supset L^* \supset 0$ of $\bW$ can only lower the weight by at most one.  Then the action of $\fg$ on $L^*$ induces a linear map $L^*\tens\fg \to \fp^{\perp}\acts\bW$.  By composing with the quotient by $L^*$ we obtain a linear map $L^*\tens\fg \surjto (\fp^{\perp}\acts \bW) / L^*$ which vanishes identically on $L^*\tens\fp$ since the action of $\fp$ is filtration preserving, thus giving a linear map $q : L^*\tens(\fg/\fp) \to (\fp^{\perp}\acts\bW)/L^*$ which surjects by construction.  Finally, observe that an element $\lambda\tens (X+\fp) \in L^*\tens \fg/\fp$ lies in $\ker q$ if and only if $X\acts\lambda \in L^*$, if and only if $X\in\fp$, so that $q$ is also injective and hence an isomorphism. \end{proof}

We will frequently omit the tensor product symbol when tensoring with elements of $L$ or $L^*$.  The following is immediate from the isomorphism $\bW_{(1)} \isom L^* \tens \fg/\fp$.

\begin{cor} \thlabel{cor:ppg-alg-gr} Under the identifications of \threfit{prop:ppg-alg-gr}, the Lie bracket between $\lambda \in L^*$ and $X \in \fg/\fp$ is given by $\liebrac{\lambda}{X} = \lambda \ltens X \in L^*\tens \fg/\fp$. \noproof \end{cor}

In terms of the Peirce decomposition \eqref{eq:ppg-alg-peirce}, we may identify%
\footnote{This is why we labelled the highest root of $\fh$ by $\beta_n$, rather than $\beta_0$ as is arguably more logical.}
\begin{align*}
  L^* \tens B       &= \Dsum{i=1}{n} \left( \bW_i \dsum {\textstyle \Dsum{j>i}{}} \, \bW_{ij} \right) \\
  L^* \tens \fg/\fp &= \Dsum{i=1}{n} \bW_{0i} \\
  L^*               &= \bW_0.
\end{align*}
Moreover, by \eqref{eq:ppg-alg-eigenvals} any algebraic Weyl structure for $\fp$ acts on these summands by multiplication by $\tfrac{2}{n+1}$, $-\tfrac{n-1}{n+1}$ and $-\tfrac{2n}{n+1}$ respectively.  This is precisely the idempotent decomposition of $\bW$ \wrt\ the single idempotent $e_0$.

It follows from \itemref{lem:ppg-alg-trans}{seq} that the Jordan eigenspaces $\bW_{ij}$ all have the same dimension, say $r \defeq \dim\bW_{ij}$.%
\footnote{Alternatively, Albert proves this directly using Jordan-theoretic methods \cite[Thm.\ 13]{a1947-jordanstr}.}
This leads to the following dimension \formulae.

\begin{cor} \thlabel{cor:ppg-alg-dims} We have $\dim\bW = \tfrac{1}{2}(n+1)(rn+2)$ and $\dim(\fg/\fp) = rn$, and consequently $\dim B = \tfrac{1}{2}n(rn-r+2)$. \noproof \end{cor}

For \ppgs, the Cartan condition then implies that $M$ has dimension $rn$.  This justifies why we ignored the case $n=0$ in \thref{asmpt:ppg-alg-n>0}, since it corresponds to a \ppg\ over a zero-dimensional manifold.

\begin{defn} \thlabel{defn:ppg-sd-params} We shall refer to $r$ and $n$ respectively as the \emph{scalar parameter} and \emph{\proj\ dimension} of a \ppg. \end{defn}

Recall that a Lie algebra $\fh$ is said be \emph{$\bZ^2$-graded} if $\fh = \Dsum{(i,j)\in\bZ^2}{} \fh_{(i,j)}$ such that $\liebrac{\fh_{(i_1,j_1)}}{\fh_{(i_2,j_2)}} \subseteq \fh_{(i_1+i_2,j_1+j_2)}$.  That is, each $\fh_{(i)} \defeq \Dsum{j\in\bZ}{} \fh_{(j,i-j)}$ is a graded component of $\fh$.  We have the following as a direct consequence of \thref{thm:ppg-alg-str,prop:ppg-alg-gr}, which is the key result that allows us to compute Lie brackets in $\fh$.

\begin{thm} \thlabel{thm:ppg-alg-Z2gr} The choice of algebraic Weyl structures for $\fq\leq\fh$ and $\fp\leq\fg$ induces the $\bZ^2$-grading of the Lie algebra $\fh$ pictured in Figure \ref{fig:ppg-alg-Z2gr}. \noproof \end{thm}

\begin{figure}[h]
\begin{equation} \label{eq:ppg-alg-Z2gr}
\arraycolsep=0.3em
\begin{array}{ccccccc}
  \fh \\
  \rotatebox{270}{\hspace{-0.7em}$\isom$} \\
  \bW   &\isom& (L^*\tens B) &\dsum& (L^*\tens \fg/\fp)             &\dsum& L^*          \\
  \dsum &     & \dsum        &     & \dsum                          &     & \dsum        \\
  \fq^0 &\isom& \fg/\fp      &\dsum& (\fp^0\oplus\liecenter{\fq^0}) &\dsum& \fp^{\perp}  \\
  \dsum &     & \dsum        &     & \dsum                          &     & \dsum        \\
  \bW^* &\isom& L            &\dsum& (L\tens\fp^{\perp})            &\dsum& (L\tens B^*).
\end{array}
\end{equation}
\caption[The $\bZ^2$-grading on $\fh$]
        {The $\bZ^2$-grading on $\fh$ described by \thref{thm:ppg-alg-Z2gr}.}
\label{fig:ppg-alg-Z2gr}
\end{figure}

Using the $\bZ^2$-grading, each $(a,b) \in \bZ$ induces a $\bZ$-grading of $\fh$: this is the grading with algebraic Weyl structure $b\xi^{\fq} + a\xi^{\fp}$, which corresponds to the grading given by stepping a line of gradient $b/a$ through \eqref{eq:ppg-alg-Z2gr}.  For example, choosing $(1,0)$ yields the grading $\fh \isom \bW\dsum\fq^0\dsum\bW^*$, while the slope $(1,1)$ yields a $|2|$-grading
\vspace{0.1em}
\begin{equation*} \begin{aligned}
  \fh
  &\isom
    (L^*\tens B) \dsum \big( \fg/\fp \dsum (L^*\tens\fg/\fp) \big) \dsum
    \big( L \dsum (\fp^0\dsum\liecenter{\fq^0}) \dsum L^* \big) \\
  &\qquad \dsum
    \big( (L\tens\fp^{\perp}) \dsum \fp^{\perp} \big) \dsum (L\tens B^*).
\end{aligned}
\vspace{0.1em}
\end{equation*}
Note that this map from $\bZ^2$ to $\bZ$-gradings of $\fh$ is not injective: for example, the pairs $(ka,kb)$ yield the same ``diagonal'' $\bZ$-grading for all $k\in\bN$.

The automorphism $\sigma$ defined by \eqref{eq:ppg-alg-sigma} may be viewed as a ``reflection'' of the diagram \eqref{eq:ppg-alg-Z2gr} as follows.%
\footnote{But note that perhaps $\sigma^2 \neq \id$ in general, since \itemref{lem:ppg-alg-sigma}{root} only guarantees that $\sigma^2(\fh_{\alpha}) = \fh_{\alpha}$.}
From the relation \eqref{eq:ppg-alg-xip}, for any $x \in L^*\tens\fg/\fp$ we have $\liebrac{\xi^{\fq}}{x} = x$ and $\liebrac{\xi^{\fp}}{x} = -\tfrac{n-1}{n+1} x$ by \eqref{eq:ppg-alg-eigenvals}.  Consequently \itemref{lem:ppg-alg-sigma}{root} gives
\vspace{0.1em}
\begin{equation*} \begin{gathered}
  \liebrac{ \xi^{\fq} }{ \sigma(x) }
    = \sigma( \liebrac{\sigma(\xi^{\fq})}{x} )
    = \sigma( \liebrac{ \xi^{\fp} + \tfrac{n-1}{n+1} \xi^{\fq} }{ x } )
    = 0 \\
  \text{and} \quad
  \liebrac{ \xi^{\fp} }{ \sigma(x) }
    = \sigma( \liebrac{ \xi^{\fq} - \tfrac{n-1}{n+1} \sigma(\xi^{\fq}) }{ x } )
    = \sigma( \liebrac{ \xi^{\fq} - \tfrac{n-1}{n+1}(\xi^{\fp} + \tfrac{n-1}{n+1} \xi^{\fq}) }
                      { x } )
    = x,
\end{gathered}
\vspace{0.1em}
\end{equation*}
implying that $\sigma(x) \in \fg/\fp$.  Therefore $\sigma(L^* \tens \fg/\fp) = \fg/\fp$; similarly we have $\sigma(L\tens\fp^{\perp}) = \fp^{\perp}$, $\sigma(L) = L^*$, while $\sigma$ preserves $L^*\tens B$, $L\tens B^*$ and $\fp^0 \dsum \liecenter{\fq^0}$.  In particular, $\sigma$ exchanges the horizontal and vertical gradings of $\fh$.  Consequently $\sigma$ induces dual Jordan algebra structures on $\tilde{\bW} \defeq (L^*\tens B) \dsum \fg/\fp \dsum L$ and $\tilde{\bW}^* \defeq L^* \dsum \fp^{\perp} \dsum (L\tens B^*)$, with Jordan products defined by $x \mathbin{\tilde{\jmult}} y \defeq \sigma( \sigma^{-1}(x) \jmult \sigma^{-1}(y) )$.

Up to the overall normalisation fixed by \thref{conv:ppg-alg-killing}, the Killing form on $\fh$ is the orthogonal direct sum of the Killing form on $\fg$, the standard inner product on $\liecenter{\fq^0}$ and the (symmetrised) pairing $\bW \times \bW^* \to \bC$.  We chose the normalisation \thref{conv:ppg-alg-killing} so that the Killing form between dual summands is simply the natural contraction.  In particular,
\vspace{0.1em}
\begin{equation*}
  \killing{h}{\theta} = \theta\intprod h, \quad
  \killing{Z}{\eta} = \eta(Z), \quad
  \killing{\ell}{\lambda} = \ell\lambda \quad\text{and}\quad
  \killing{X}{\alpha} = \alpha(X)
\vspace{0.1em}
\end{equation*}
for $(h,Z,\lambda) \in \bW$, $(X,\alpha) \in \fg/\fp \dsum \fp^{\perp}$ and $(\ell,\eta,\theta) \in\bW^*$.  We can now calculate a large number of algebraic brackets between the summands of \eqref{eq:ppg-alg-Z2gr}.

\begin{thm} \thlabel{thm:ppg-alg-bracs} Consider elements
\begin{equation} \label{eq:ppg-alg-bracselems} \begin{alignedat}{4}
  (h,Z,\lambda) &\in \bW &&\isom (L^*\tens B) \dsum (L^* \tens \fg/\fp) \dsum L^*, \\
  (X,A,\alpha) &\in \fq^0 &&\isom \fg/\fp \dsum (\fp^0\dsum\liecenter{\fq^0}) \dsum \fp^{\perp} \\
  \text{and}\quad
  (\ell,\eta,\theta) &\in \bW^* &&\isom L \dsum (L\tens\fp^{\perp}) \dsum (L\tens B^*),
\end{alignedat} \end{equation}
as well as their primed counterparts.  Then with normalisation conventions as above, Lie brackets between the various summands in $\fh$ are given by Table \ref{tbl:ppg-alg-Z2gr}. \end{thm}

\vspace{-0.7em}
\bgroup \begin{table}[!h]
\begin{equation*}
\def\arraystretch{1.15}  
\arraycolsep=0.5ex         
\newcommand{\e}[2]{      
  \ifx  \relax#1\relax  \cellcolor{gray!10}
  \else
    \ifx  \relax#2\relax  #1
    \else                 \begin{array}{c} #1 \\[-0.45em] \mathsmaller{\in\, #2} \end{array}
    \fi
  \fi }
\newcommand{\E}[1]{\e{#1}{}}
\newcommand{\g}[2]{\e{#1}{#2} \cellcolor{gray!20}}
\newcommand{\w}[2]{\e{#1}{#2} \cellcolor{gray!40}}
\newcommand{\vpad}{\vphantom{\e{\liebrac{A}{\ell}}{L}}}
\newcommand{\cent}{\fp^0\dsum\liecenter{\fq^0}}
\begin{array}{|c|*{3}{ccc|}}
  \cline{2-10} \multicolumn{1}{c|}{}
              & \E{~h'~}                                 & \E{Z'}
              & \E{~\lambda'~}                           & \E{X'}
              & \E{A'}                                   & \E{\alpha'}
              & \E{\ell'}                                & \E{\eta'}
              & \E{\theta'}
  \\ \hline 
  \E{h}       & \e{0}{}                                  & \e{0}{}
              & \e{0}{}                                  & \e{0}{}
              & \w{\liebrac{h}{A'}}{L^*\tens B}          & \e{h(\alpha',\bdot)}{L^*\tens\fg/\fp}
              & \e{0}{}                                  & \e{-h(\eta',\bdot)}{\fg/\fp}
              & \g{\liebrac{h}{\theta'}}{\cent}
  \\
  \E{Z}       & \e{}{}                                   & \e{0}{}
              & \e{0}{}                                  & \g{\liebrac{Z}{X'}}{L^*\tens B}
              & \w{\liebrac{Z}{A'}}{L^*\tens\fg/\fp}     & \e{\alpha'(Z)}{L^*}
              & \e{-\ell'\ltens Z}{\fg/\fp}              & \g{\liebrac{Z}{\eta'}}{\cent}
              & \e{-\theta'(Z,\bdot)}{\fp^{\perp}}
  \\
  \E{\lambda} & \e{}{}                                   & \e{}{}
              & \e{0}{}                                  & \e{\lambda\ltens X'}{L^*\tens\fg/\fp}
              & \w{\tfrac{2(\tr A')}{r(n+1)}\lambda}{L^*}& \e{0}{}
              & \g{\liebrac{\lambda}{\ell'}}{\cent}      & \e{-\lambda\ltens\eta'}{\fp^{\perp}}
              & \e{0}{}
  \\ \hline
  \E{X}       & \e{}{}                                   & \e{}{}
              & \e{}{}                                   & \e{0}{}
              & \e{-A'X}{\fg/\fp}                        & \g{\liebrac{X}{\alpha'}}{\fp^0}
              & \e{0}{}                                  & \e{\eta'(X)}{L}
              & \e{\theta'(X,\bdot)}{L\tens\fp^{\perp}}
  \\
  \E{A}       & \e{}{}                                   & \e{}{}
              & \e{}{}                                   & \e{}{}
              & \e{AA'{-}A'A}{\fp^0}                     & \e{-\alpha'\circ A}{\fp^{\perp}}
              & \raisebox{-0.3ex}{$
                  \w{\tfrac{2(\tr A)}{r(n+1)}\ell'}{L}$} & \w{\liebrac A{\eta'}}{L\tens\fp^{\perp}}
              & \w{\liebrac{A}{\theta'}}{L\tens B^*}
  \\
  \E{\alpha}  & \e{}{}                                   & \e{}{}
              & \e{}{}                                   & \e{}{}
              & \e{}{}                                   & \e{0}{}
              & \e{\ell'\ltens\alpha}{L\tens\fp^{\perp}} & \g{\liebrac{\alpha}{\eta'}}{L\tens B^*}
              & \e{0}{}
  \\ \hline
  \E{\ell}    & \e{}{} \vpad                             & \e{}{}
              & \e{}{}                                   & \e{}{}
              & \e{}{}                                   & \e{}{}
              & \e{0}{}                                  & \e{0}{}
              & \e{0}{}
  \\
  \E{\eta}    & \e{}{} \vpad                             & \e{}{}
              & \e{}{}                                   & \e{}{}
              & \e{}{}                                   & \e{}{}
              & \e{}{}                                   & \e{0}{}
              & \e{0}{}
  \\
  \E{\theta}  & \e{}{} \vpad                             & \e{}{}
              & \e{}{}                                   & \e{}{}
              & \e{}{}                                   & \e{}{}
              & \e{}{}                                   & \e{}{}
              & \e{0}{}
  \\ \hline
\end{array}
\end{equation*}
\vspace{-0.5em}
\caption[Lie brackets between the $\bZ^2$-graded summands of $\fh$]
        {Lie brackets between the various summands of $\fh$ according to the $\bZ^2$-grading \eqref{eq:ppg-alg-Z2gr}, where the entry labelled by row $x$ and column $y$ is the bracket $\liebrac{x}{y}$.  Elements are defined in \eqref{eq:ppg-alg-bracselems}, the colouring pertains to the proof of \thref{thm:ppg-alg-Z2gr}, and the empty part of the table may be determined by skew-symmetry.}
\label{tbl:ppg-alg-Z2gr}
\end{table}  \egroup
\vspace{-0.3em}

\begin{proof} Since $\bW$ and $\bW^*$ are abelian Lie algebras, all brackets in the top-left and bottom-right squares of Table \ref{tbl:ppg-alg-Z2gr} are zero. Moreover $\liebrac{X}{X'}=0$ and $\liebrac{\alpha}{\alpha'}=0$ by \itemref{thm:ppg-alg-str}{g}, while the form of the $\bZ^2$-grading implies that the brackets $\liebrac{h}{X'}$, $\liebrac{\lambda}{\alpha'}$, $\liebrac{h}{\ell'}$, $\liebrac{\lambda}{\theta'}$, $\liebrac{X}{\ell'}$ and $\liebrac{\alpha}{\theta'}$ all vanish. This accounts for all the zeroes in Table \ref{tbl:ppg-alg-Z2gr}.

Next, we have $\liebrac{\lambda}{X} = \lambda\ltens X \in L^*\tens\fg/\fp$ by \thref{cor:ppg-alg-gr}.  From the $\bZ^2$-grading we know that $\liebrac{\eta}{\lambda} \in \fp^{\perp}$.  Since the Killing form between $L\tens\fp^{\perp}$ and $L^*\tens\fg/\fp$ is just the contraction, it follows that $\killing{\liebrac{\eta}{\lambda}}{X} = \killing{\eta}{\liebrac{\lambda}{X}} = \eta(\lambda\ltens X)$ for all $X\in\fg/\fp$, so that $\liebrac{\eta}{\lambda} = \lambda\ltens\eta \in \fp^{\perp}$ by \non degeneracy.  Similarly $\liebrac{X}{\eta} \in L$, for which $\killing{\liebrac{X}{\eta}}{\lambda} = \killing{X}{\liebrac{\eta}{\lambda}} = \eta(\lambda \ltens X)$ and hence $\liebrac{X}{\eta} = \eta(X) \in L$.

The reflection $\sigma$ allows us to ascertain the Lie bracket $\liebrac{\ell}{\alpha}$: since $L = \sigma(L^*)$ and $\fp^{\perp} = \sigma(L \tens \fp^{\perp})$, we may write $\ell = \sigma(\lambda)$ and $\alpha = \sigma(\eta)$, for which $\liebrac{\ell}{\alpha} = \sigma(\liebrac{\lambda}{\eta}) = -\sigma(\lambda \ltens \eta)$.  Since $\sigma$ is defined in terms of the adjoint action, the Leibniz rule gives $\sigma(\lambda \ltens \eta) = \sigma(\lambda) \tens \sigma(\eta) = \ell \ltens \alpha$ and hence $\liebrac{\ell}{\alpha} = -\ell \ltens \alpha$.  Arguing as in the previous paragraph, it follows also that $\liebrac{Z}{\alpha} = \alpha(Z) \in L^*$ and $\liebrac{Z}{\ell} = -\ell \ltens Z \in \fg/\fp$.

The bracket $\liebrac{h}{\alpha}$ may be calculated using the map \eqref{eq:ppg-alg-hommap}.  By construction, $h \in L^*\tens B$ is identified with the linear map $\alpha \tens \beta \mapsto h(\alpha,\beta) \defeq \liebrac{ \alpha }{ \liebrac{\beta}{h} }$.  Since $\liebrac{\beta}{h} \in L^* \tens \fg/\fp$, the above gives $h(\alpha, \beta) = -\alpha( \liebrac{\beta}{h} ) \in L^*$, which implies that $\liebrac{h}{\alpha} = h(\alpha,\bdot)$.  For $\liebrac{\theta}{X} \in L \tens \fp^{\perp}$, we identify $L \tens B^* = \liecohom{0}{\bW^*} \isom \liehom[\fg/\fp]{0}{\bW^*}$ by \eqref{eq:lie-hom-isom}.  Thus $\theta$ may be identified with the linear map $X\tens Y \mapsto \theta(X,Y) = \liebrac{ X }{ \liebrac{Y}{\theta} }$, giving $\liebrac{\theta}{X} = -\theta(X,\bdot)$.  Writing $\eta = \ell \ltens \alpha = \liebrac{\alpha}{\ell}$, we then have $\liebrac{h}{\eta} = \liebrac{ \liebrac{h}{\alpha} }{ \ell } = -\ell \ltens h(\alpha,\bdot) = -h(\eta,\bdot)$ by the previous computations; $\liebrac{\theta}{Z} = \theta(Z,\bdot)$ similarly.

By the grading on $\fg$, we view $\fp^0$ as a subalgebra of $\alg{gl}{\fg/\fp} \isom \alg{gl}{rn,\bC}$ and the brackets of $\fp^0$ with $\fg/\fp$ and $\fp^{\perp}$ are the natural actions; whence $\liebrac{B}{X} = BX$ and $\liebrac{B}{\alpha} = -\alpha\circ B$ for all $B\in\fp^0$.  Since $\liecenter{\fq^0}$ is abelian and acts trivially on $\fg$, the brackets involving $X$, $A$ and $\alpha$ follow.  Finally, the entries highlighted in light grey depend on the Lie algebra $\fh$ and may be determined on a case-by-case basis as in \thref{rmk:ppg-alg-remaining}; entries highlighted in dark grey will follow from the forthcoming \thref{prop:ppg-alg-Lwedge}. \end{proof}

Note that elements $h \in L^*\tens B$ may be viewed as $L^*$-valued symmetric bilinear forms on $\fp^{\perp}$; indeed \eqref{eq:ppg-alg-hommap} gives $h(\alpha,\beta) = \liebrac{ \liebrac{h}{\alpha} }{ \beta }$ for all $\alpha,\beta \in \fp^{\perp}$, for which $h(\alpha,\beta) = h(\beta,\alpha)$ since $\fp^{\perp}$ is abelian by \itemref{thm:ppg-alg-str}{g}.  Similarly, elements $\theta \in L\tens B^*$ may be viewed as $L$-valued symmetric bilinear forms on $\fg/\fp$.

\subsection{Characterisation of $L$} 
\label{ss:ppg-alg-L}

Here we describe the role of $L$, as well as completing the entries of Table \ref{tbl:ppg-alg-Z2gr}.

\begin{prop} \thlabel{prop:ppg-alg-Lwedge} $\Wedge{rn}(\fg/\fp)$ and $L^{r(n+1)/2}$ are isomorphic as $\fp^0$-representations. \end{prop}

\begin{proof} Since $\dim(\fg/\fp) = rn$, both $\Wedge{rn}(\fg/\fp)$ and $L$ are $1$-dimensional and hence irreducible $\fp^0$-representations.  By Schur's lemma the action is by a scalar, so it suffices to compare the actions of the grading element of $\fp$.  By \eqref{eq:ppg-alg-eigenvals}, the algebraic Weyl structure of $\fp$ acts on $L \isom \bW^*_0$ by multiplication by $2n/(n+1)$.  Since $\fg/\fp$ has weight one, $\Wedge{rn}(\fg/\fp)$ has weight
\begin{equation*}
  rn = \frac{r(n+1)}{2} \cdot \frac{2n}{n+1}
\end{equation*}
and hence the isomorphism follows. \end{proof}

\vspace{-0.15em}
The isomorphism afforded by \thref{prop:ppg-alg-Lwedge} may be used to gain information about the Lie bracket $\liebrac{X}{\alpha} \in \fp^0$ between elements $X\in\fg/\fp$ and $\alpha\in\fp^{\perp}$.
\vspace{-0.15em}

\begin{lem} \thlabel{lem:ppg-alg-Lbrac} $\liebrac{ \liebrac{X}{\alpha} }{ \ell } = \alpha(X)\ell$ and $\liebrac{ \liebrac{X}{\alpha} }{ \lambda } = -\alpha(X)\lambda$ for all $\ell \in L$ and $\lambda \in L^*$. \end{lem}

\begin{proof} We have $\liebrac{ \liebrac{X}{\alpha} }{ \ell } = \liebrac{ \liebrac{X}{\ell} }{ \alpha } + \liebrac{ X }{ \liebrac{\alpha}{\ell} }$ by the Jacobi identity, where the first summand vanishes since the bracket between $\fg/\fp$ and $L$ is trivial.  The second summand equals $\liebrac{ X }{ \ell\ltens\alpha } = \alpha(X)\ell$ by Table \ref{tbl:ppg-alg-Z2gr}. \end{proof}

\begin{cor} \thlabel{cor:ppg-alg-trace} $\liebrac{A}{\ell} = \tfrac{2}{r(n+1)}(\tr A)\ell \in L$ for all $A \in \fp^0$ and $\ell \in L$.  In particular, $\liebrac{X}{\alpha}$ has trace $\tfrac{1}{2}r(n+1)$ as an endomorphism of $\fg/\fp$ for all $X \in \fg/\fp$ and $\alpha \in \fp^{\perp}$. \end{cor}

\begin{proof} The action of $A \in P^0$ on $\Wedge{rn}(\fg/\fp)$ is given by multiplication by the determinant of $A$.  Since the derivative of the determinant is the trace, the action of $A \in \fp^0$ is given by multiplication by the trace of $A$.  The first claim now follows from \thref{prop:ppg-alg-Lwedge}, while the second claim follows from \thref{lem:ppg-alg-Lbrac}. \end{proof}

\begin{rmk} \thlabel{rmk:ppg-alg-remaining} Note that \thref{cor:ppg-alg-trace} proves the dark grey entries in Table \ref{tbl:ppg-alg-Z2gr}.  Indeed, the brackets $\liebrac{A}{\ell}$ and $\liebrac{A}{\lambda}$ follow immediately, while we have
\begin{align*}
  \liebrac{A}{h} (\alpha,\beta)
  &= \liebrac{A}{h(\alpha,\beta)} - h(\liebrac{A}{\alpha},\beta) - h(\alpha,\liebrac{A}{\alpha}) \\
  &= -\tfrac{r(n+1)}{2} h(\alpha,\beta) + h(\alpha\circ A, \beta) + h(\alpha, \beta\circ A)
\end{align*}
by the Leibniz rule.  It also allows us to calculate the light grey entries in terms of the Lie bracket $\fg/\fp \times \fp^{\perp} \to \fp^0$, which depend on the Lie algebra $\fh$ in question.  For example, for $Z \in L^*\tens B$ and $X \in \fg/\fp$ we have $\liebrac{Z}{X} \in L^*\tens B$ by the $\bZ^2$-grading.  Therefore by the Jacobi identity
\begin{align*}
  \liebrac{Z}{X}(\alpha,\beta)
  &\defeq \liebrac{ \liebrac{ \liebrac{Z}{X} }
                            { \alpha } }
                  { \beta } \\
  &\phantom{:}= \liebrac{ \liebrac{\alpha(Z)}{X} }
                        { \beta }
    + \liebrac{ \liebrac{ Z }
                        { \liebrac{X}{\alpha} } }
              { \beta } \\
  &\phantom{:}= \liebrac{ \alpha{Z}X }{ \beta } + \liebrac{ \beta(Z) }{ \liebrac{X}{\alpha} }
    + \liebrac{ Z }{ \liebrac{ \liebrac{X}{\alpha} }{ \beta } } \\
  &\phantom{:}= \alpha(Z)\beta(X) + \alpha(X)\beta(Z) + \liebrac{\liebrac{X}{\alpha}}{\beta}(Z)
\end{align*}
for all $\alpha,\beta \in \fp^{\perp}$.  The brackets $\liebrac{\alpha}{\eta}$, $\liebrac{\lambda}{\ell}$, $\liebrac{Z}{\eta}$ and $\liebrac{h}{\theta}$ may be determined similarly.  For example, in the case $\fh = \alg{sp}{2n+2,\bC}$ one finds that $\fp^0 = \alg{gl}{n,\bC}$ and $\liebrac{X}{\alpha} = \tfrac{1}{2}(\alpha(X)\id + X\tens\alpha)$, yielding $\liebrac{Z}{X} = X \symm Z$ and $\liebrac{\alpha}{\eta} = \alpha \symm \eta$. \end{rmk}

For the classical \proj\ structures, there was a bijection between the affine space of Weyl structures and the induced connections on the bundle $\cL$ associated to $L$.

\begin{cor} \thlabel{cor:ppg-alg-connsonL} For any \ppg, there is a bijection between Weyl structures and their induced connections on $\cL$. \end{cor}

\begin{proof} Let $\D, \b{\D}$ be the Weyl connections associated to two Weyl structures, and suppose that $\b{\D} = \D + \algbrac{}{\gamma}$ for some $\gamma \in \s{1}{}$.  Then on sections of $\cL$ we have $\b{\D}\ell = \D\ell + \ell \ltens \gamma$, so that the induced connections coincide if and only if $\gamma = 0$. \end{proof}

In the terminology of \cite[\S 3.5]{cs2003-weylstr}, this means that $\cL$ is a \emph{bundle of scales} for the underlying parabolic geometry.  \thref{cor:ppg-alg-trace,cor:ppg-alg-connsonL} may then be viewed as special cases of Proposition 3.2(2) and Theorem 3.8(1) from \cite{cs2003-weylstr}.  We shall not need this level of generality here.

\begin{rmk} \thlabel{rmk:ppg-alg-Z2nonsd} The $\bZ^2$-grading from \thref{thm:ppg-alg-Z2gr} and the Lie brackets calculated in \thref{thm:ppg-alg-bracs} only require that $\fp$ is an abelian parabolic and $\bW$ has height two.  Therefore they go through in the \non \self dual\ case alluded  to in \thref{rmk:ppg-alg-str}.  The key difference is that we do not have a Peirce decomposition for $\bW$, and cannot identify an integer $r \defeq \dim \bW_{ij}$.  In this case, we find that $L \isom (\Wedge{m}(\fg/\fp))^{(2w-1)/mw}$, where $m \defeq \dim(\fg/\fp)$ and $w \defeq \killing{\omega_r}{\omega_r}$. \end{rmk}

\section{Classification} 
\label{s:ppg-class}

Having described the algebraic structure associated to a \ppg, we are overdue writing down a classification.  This is a straightforward application of the classification of symmetric \Rspaces, which is well-known; see for example \cite[Prop.\ 3.2.3]{cs2009-parabolic1}.  We first give a classification over $\bC$ in Subsection \ref{ss:ppg-class-cpx}, before selecting appropriate real forms in Subsection \ref{ss:ppg-class-real}.

\subsection{Classification over $\bC$} 
\label{ss:ppg-class-cpx}

Suppose that $\fh$ is a complex semisimple Lie algebra with parabolic $\fq$, and choose a Cartan subalgebra $\ft$ and simple subsystem $\Delta^0$ \wrt\ which $\fq$ is the standard parabolic associated to a subset $\Sigma \subseteq \Delta^0$.  These choices carry an algebraic Weyl structure for $\fq$ with them, hence splitting the $\fq^{\perp}$-filtration of $\fh$, with grading given by the $\Sigma$-height of roots of $\fh$.  In particular, the height of $\fq$ is therefore equal to the $\Sigma$-height of the highest root.  It follows that $\fq$ is an abelian parabolic if and only if the crossed node corresponds to a one in the highest root, leading to a short and simple classification of irreducible symmetric \Rspaces.  Tables of highest roots may be found in \cite[Tbl.\ B.2]{cs2009-parabolic1} and \cite{h1972-repntheory}.

To classify the flat models of \ppgs, it suffices then to decide which symmetric \Rspaces\ are \self dual, and consider their isotropy \Rspaces.  This is straightforward thanks to \thref{prop:ppg-defn-dynkin}, so it remains to go through each Dynkin type, determine whether $H\acts\fq$ is \self dual, and if so find the algebra $\fg$, the $\fg$-representation $\bW$ and the parabolic $\fp$.  The integer $n$, which by definition equals the maximal number of orthogonal roots of $\fq$-height one, may be read off from the tables of \cite[\S 3-5]{ak2002-strongorth}; \thref{cor:ppg-alg-dims} then gives $r = \dim(\fg/\fp) / n$.  Equivalently $r \defeq \dim\bW_{ij}$ and $n = \dim(\fg/\fp) / r$, so that one may calculate $(r,n)$ from the Peirce decomposition.

\begin{typeblock}{A}{k} The highest root of $\fh = \alg{sl}{k+1,\bC}$ is the sum of all simple roots, so we may cross any node to obtain an abelian parabolic $\fq$.  Crossing the $\ell$th node determines a symmetric \Rspace\ $H\acts\fq$ with minimal embedding $H\acts\fq \injto \pr{\Wedge[\bC]{\ell}\bC^{k+1}}$ which, since the usual involution flips the Dynkin diagram, is \self dual\ if and only if $k=2m-1$ is odd and we cross the $m$th node.  Thus there is a single \self dual\ symmetric \Rspace\ of this type, corresponding to $\fh = \alg{sl}{2m,\bC}$ and $H\acts\fq \injto \pr{\Wedge[\bC]{m}\bC^{2m}}$, given by the \grassmannian\ $\Grass{m}{\bC^{2m}}$ of complex $m$-planes in $\bC^{2m}$.

The Dynkin diagram of the semisimple part $\fg$ of $\fq^0$ is given by
\vspace{0.3em}
\begin{equation*}
\smash{
  \fg = \dynkinAA{2}{0}{2}{ppg-class-cpx-Ag}
    = \alg{sl}{m,\bC} \dsum \alg{sl}{m,\bC}.
}
\end{equation*}
In particular, $\fg$ is not simple.  The adjoint representation of $\fh$ and the $\fg$-representation $\bW$ induced by $\fh/\fq$ have highest weights
\begin{equation*}
  \fh = \dynkinASU[1,0,0,0, 1,0,0,0, 0]{2}{0}{2}{ppg-class-cpx-Ah}
  \quad\text{and}\quad
  \bW = \dynkinAA[1,0,0,0, 1,0,0,0]{2}{0}{2}{ppg-class-cpx-AW}
    = \bC^m \etens \conj{\bC^m}.
\end{equation*}
Then the highest weight of $\bW^*$ is supported on the final node of each branch, so that
\vspace{-0.5\baselineskip} \begin{equation*} \vspace{-0.5\baselineskip}
  G\acts\fp = \dynkinAAp{2}{0}{2}{ppg-class-cpx-Ap}
    = \CP[m-1]\times\CP[m-1]
\end{equation*}
is a product of \proj\ spaces.  The Jordan algebra $\bW$ is the (complexification of) the space of $m \by m$ hermitian matrices under multiplication, with $r = 2$ and $n = m-1 = \tfrac{1}{2}(k-1)$.  In future, it will be convenient to talk about type \type{A}{2n+1}. \end{typeblock}

\begin{typeblock}{B}{k} The highest root of $\fh = \alg{so}{2k+1,\bC}$ has a coefficient of one for only the first simple root, so there is a symmetric \Rspace\ $H\acts\fq$ determined by crossing this node.  The corresponding minimal \proj\ embedding is $H\acts\fq \injto \pr{\bC^{2k+1}}$, identifying $H\acts\fq$ with the complex conformal sphere $\cpxbdl{\Sph[2k-1]}$.

The Dynkin diagram of $\fg$ is given by removing the first node of $\fh$, so that $\fg = \alg{so}{2k-1,\bC}$.  The adjoint representation of $\fh$ has highest weight
\begin{equation*}
  \fh = \dynkinB[0,1,0,0,0,0]{3}{0}{3}{ppg-class-cpx-Bh}
  \quad\text{and hence}\quad
  \bW = \dynkinB[1,0,0,0,0]{2}{0}{3}{ppg-class-cpx-BW}
    = \bC^{2k-1}
\end{equation*}
as a $\fg$-representation.  Then, since the fundamental representations in type \type{B}{} are isomorphic to their duals, we have $G\acts\fp = \smash{ \dynkinBp{2}{0}{3}{ppg-class-cpx-Bp} }$, which is the complexified conformal sphere $\CSph[2k-3]$.  The Jordan algebra $\bW$ should be viewed as $\bC^{2k-2} \dsum \bC$, the algebra of \emph{spin factors} \cite{m2007-jordan}, with product given by Clifford multiplication; this may equally be viewed as $1$-dimensional \proj\ geometry over $\bC^{2k-3}$ as in \cite{ll2010-specialholonomy}.  The integers $(r,n)$ are given by $r=2k-3$ and $n=1$. \end{typeblock}

\begin{typeblock}{C}{k} The long simple root alone has coefficient one in the highest root of $\fh = \alg{sp}{2k,\bC}$, giving a single symmetric \Rspace\ $H\acts\fq$ given by crossing this node.  It has minimal \proj\ embedding $H\acts\fq \injto \pr{\Wedge[0]{k}\bC^{2k}}$ into the space of primitive $k$-forms, which identifies $H\acts\fq$ with the complex \grassmannian\ of langrangian subspaces of $\bC^{2k}$.

The Dynkin diagram of $\fg$ is given by removing the long node of $\fh$, so that $\fg = \alg{sl}{k,\bC}$.  The adjoint representation of $\fh$ has highest weight
\vspace{0.1em}
\begin{equation*}
  \fh = \dynkinC[2,0,0,0,0]{2}{0}{3}{ppg-class-cpx-Ch}
  \quad\text{and hence}\quad
  \bW = \dynkinA[2,0,0,0]{2}{0}{2}{ppg-class-cpx-CW}
    = \Symm[\bC]{2}\bC^{k}
\vspace{0.1em}
\end{equation*}
as a $\fg$-representation.  It follows that $G\acts\fp = \smash{ \dynkinAp{2}{0}{2}{ppg-class-cpx-Cp} }$, which is the complex \proj\ space $\CP[k-1]$.  The Jordan algebra $\bW$ is the space of $k \by k$ complex symmetric matrices, with $r=1$ and $n=k-1$; thus we shall talk about type \type{C}{n+1} in future. \end{typeblock}

\begin{typeblock}{D}{k} The highest root of $\fh = \alg{so}{2k,\bC}$ has a one over either the first node, or either half-spin node, so we may cross any of these.  Crossing the first node yields the \self dual\ \Rspace\ $H\acts\fq = \cpxbdl{\Sph[2k-2]}$, with $\fg = \alg{so}{2k-2,\bC}$ and $G\acts\fp = \CSph[2k-4]$.  The Jordan algebra is $\bW = \bC^{2k-2}$, with $r=2k-4$ and $n=1$.  Types \type{B}{k+2} and \type{D}{k+2} may therefore be treated concurrently, which we refer to as type \type{BD}{n+4}.

On the other hand, crossing either spin node yields two isomorphic symmetric \Rspaces\ $H\acts\fq$ with minimal \proj\ embeddings $H\acts\fq \injto \pr{\spinrepn^{\pm}}$.  The usual involution preserves this crossed node if and only if $k$ is even, so we obtain a \self dual\ \Rspace\ only when $k = 2m$.  The Dynkin diagram of $\fg$ is then given by removing one spin node, yielding $\fg = \alg{sl}{2m,\bC}$.  The adjoint representation of $\fh$ is
\vspace{0.3em}
\begin{equation*}
  \fh = \smash{ \dynkinD[0,1,0,0,0,0,0]{3}{0}{2}{ppg-class-cpx-Dh} }
  \quad\text{and hence}\quad
  \bW = \dynkinA[0,1,0,0,0,0]{3}{0}{3}{ppg-class-cpx-DW}
    = \Wedge[\bC]{2}\bC^{2m}
\vspace{0.2em}
\end{equation*}
as a $\fg$-representation.  Therefore $G\acts\fp = \dynkinApp{3}{0}{3}{ppg-class-cpx-Dp}$, which is the \grassmannian\ of complex $2$-planes in $\bC^{2m-2}$.  The Jordan algebra $\bW$ is the complexification of the space of $m \by m$ \qtn-hermitian matrices under multiplication, with $r=4$ and $n = m-1 = \tfrac{1}{2}(k-2)$; it will be convenient to talk about type \type{D}{2n+2} in future. \end{typeblock}

\begin{typeblock}{E}{6} Temporarily writing the coefficients of simple roots in the Dynkin diagram, the highest root of $\fh = \alg[_6]{e}{\bC}$ is given by
\begin{equation*}
  \dynkinE[1,2,3,2,1,2]{6}{ppg-class-cpx-E6h}
    \colvectpunct[-0.7em]{\, ,}
\end{equation*}
so that we get isomorphic symmetric \Rspaces\ by crossing either the left-most or right-most node.  However the usual involution exchanges these two nodes, so there are no \self dual\ symmetric \Rspaces\ in this type. \end{typeblock}

\begin{typeblock}{E}{7} Returning to usual fundamental weight notation, the highest root of $\fh = \alg[_7]{e}{\bC}$ is given by
\begin{equation*}
  \dynkinE[0,0,0,0,0,1,0]{7}{ppg-class-cpx-E7root}
    \colvectpunct[-0.7em]{\, ,}
\end{equation*}
so we get a single symmetric \Rspace\ $H\acts\fq$ by crossing the right-most node, which is \self dual\ since the usual involution is just the identity.  Its minimal \proj\ embedding is $H\acts\fq \injto \pr{\eerepn^*}$, where $\eerepn$ is the $56$-dimensional representation of $\alg[_7]{e}{\bC}$, so that $H\acts\fq$ is an exceptional manifold of dimension $27$.

The Dynkin diagram of $\fg$ is given by removing the right-most node, so that $\fg = \alg[_6]{e}{\bC}$.  The adjoint representation of $\fh$ has highest weight
\begin{equation*}
  \fh = \dynkinE[1,0,0,0,0,0,0]{7}{ppg-class-cpx-E7h}
  \quad\text{and hence}\quad
  \bW = \dynkinE[1,0,0,0,0,0]{6}{ppg-class-cpx-E7W} = \erepn.
\end{equation*}
The symmetric \Rspace\ $G\acts\fp$ is given by
\vspace{-0.35em}
\begin{equation*}
  G\acts\fp = \dynkinEp{6}{ppg-class-cpx-E7p}
    \colvectpunct[-0.7em]{\, ,}
\vspace{0.35em}
\end{equation*}
which is Rosenfeld's \proj\ plane $\pr{\bC\tens\bO}$ over the bi-\oct s\ $\bC\tens\bO$; see \cite[\S 1.5]{r1997-liegeom}, \cite[\S 4.3]{b2002-octonions}, and \cite[pp.\ 313--316]{b1987-einstein} for further details. \end{typeblock}

\begin{typeblock}[{\typeheader[Types ]{E}{8}, \typeheader[]{F}{4} and \typeheader[]{G}{2}}]{}{} Temporarily writing the coefficients of simple roots over the nodes of the Dynkin diagram, the highest roots of $\alg[_8]{e}{\bC}$, $\alg[_4]{f}{\bC}$ and $\alg[_2]{g}{\bC}$ are
\begin{equation*}
  \dynkinE[2,3,4,5,6,4,2,3]{8}{ppg-class-cpx-E8root}, \quad
  \dynkinF[2,3,4,2]{ppg-class-cpx-F4root}
  \quad\text{and}\quad
  \dynkinG[3,2]{ppg-class-cpx-G2root}
\end{equation*}
respectively, so there are no symmetric \Rspaces\ in these types. \end{typeblock}

In summary, we have the following classification of irreducible \ppgs\ over the complex numbers.  Here we choose the rank of $\fh$ to ensure a simple expression for $n$ and consequently the dimension of $G\acts\fp$.

\begin{thm} \thlabel{thm:ppg-class-cpxclass} The data $\setof{\fh,\fg,\bW}{}$ describing the flat model of an irreducible \ppg\ is contained in Table \ref{tbl:ppg-class-cpx}. \end{thm}

\begin{table}[!h]
  \centering
  \begin{tabular}{| c | *{5}{>{$}c<{$}|} } \hline
    Type & \fh & \fg & \bW & r & n
    \\ \hline
    \type{C}{n+1} & \alg{sp}{2n+2,\bC} & \alg{sl}{n+1,\bC} & \Symm[\bC]{2}\bC^{n+1} & 1 & n
    \\
    \type{A}{2n+1} & \alg{sl}{2n+2,\bC} & \alg{sl}{n+1,\bC}\dsum\alg{sl}{n+1,\bC} & \bC^{n+1}\etens\conj{\bC^{n+1}} & 2 & n
    \\
    \type{D}{2n+2} & \alg{so}{4n+4,\bC} & \alg{sl}{2n+2,\bC} & \Wedge[\bC]{2}\bC^{2n+2} & 4 & n
    \\
    \type{E}{7} & \alg[_7]{e}{\bC} & \alg[_6]{e}{\bC} & \erepn & 8 & 2
    \\
    \type{BD}{n+4} & \alg{so}{n+4,\bC} & \alg{so}{n+2,\bC} & \bC^{n+2} & n & 1
    \\ \hline
  \end{tabular}
  \caption[Classification of complex \ppgs]
          {The classification of the complexified flat models of \ppgs.}
  \label{tbl:ppg-class-cpx}
\end{table} \vspace{-0.5\baselineskip}

\medskip
Note in particular that there is a unique $r$ for each admissible simple type, so that the classification may be phrased entirely in terms of $r$.

\begin{rmk} Note that there are some special isomorphisms between the symmetric \Rspaces\ $G\acts\fp$.  Namely, the complexified conformal spheres $\CSph[1]$, $\CSph[2]$ and $\CSph[4]$ are isomorphic to $\CP[1]$, $\CP[1]\times\CP[1]$ and $\Grass{2}{\bC^4}$; moreover $\CSph[8]$ may be viewed as Rosenfeld's $\pr{\bC\tens\bO}$.  This is the origin of \thref{rmk:cproj-bgg-dim2,rmk:qtn-bgg-dim4} which identified $\CP$ and $\HP$ with \self dual\ conformal structures of dimensions two and four. \end{rmk}

We collect some useful numerical data in Table \ref{tbl:ppg-class-dims}.  By \thref{cor:ppg-alg-trace}, the bracket $\liebrac{X}{\alpha}$ has trace $\tfrac{1}{2} r(n+1) \alpha(X)$; the value $\tr\liebrac{X}{\alpha} = \tfrac{1}{2}(n+1) \alpha(X)$ in type \type{C}{n+1} justifies our unusual normalisation convention for the algebraic bracket \eqref{eq:proj-class-algbrac} in Chapter \ref{c:proj}.

Finally, the highest weights of the graded components of $\bW$ and $\bW^*$ are collected in Tables \ref{tbl:app-tbl-W} and \ref{tbl:app-tbl-Wd} respectively.  Note in particular that we recover the $1$-dimensional representations $L$ defined by \eqref{eq:proj-para-repnL}, \eqref{eq:cproj-para-repnL} and \eqref{eq:qtn-para-repnL} for the classical cases.

\vspace{0.4\baselineskip}
\begin{table}[h]
  \tabcolsep=1.5ex                                
  \centering
  \begin{tabular}{| c | *{6}{>{$}c<{$}|} }
    \hline
    Type           & r & n & \dim\bW                & \dim(\fg/\fp)
                   & \dim B                         & \tfrac{1}{2}r(n+1)
    \\ \hline
    \type{C}{n+1}  & 1 & n & \tfrac{1}{2}(n+1)(n+2) & n
                   & \tfrac{1}{2}n(n+1)             & \tfrac{1}{2}(n+1)  \\
    \type{A}{2n+1} & 2 & n & (n+1)^2                & 2n
                   & n^2                            & n+1                \\
    \type{D}{2n+2} & 4 & n & (n+1)(2n+1)            & 4n
                   & n(2n-1)                        & 2n+2               \\
    \type{E}{7}    & 8 & 2 & 27                     & 16
                   & 10                             & 12                 \\
    \type{BD}{n+4} & n & 1 & n+2                    & n
                   & 1                              & n                  \\ \hline
  \end{tabular}
  \caption[Numerical data for the complex \ppgs]
          {Some useful numerical data for each complex \ppg, written in terms of the integers $(r,n)$.}
  \label{tbl:ppg-class-dims}
\end{table} \vspace{-0.4\baselineskip}

\subsection{Real forms} 
\label{ss:ppg-class-real}

Using Table \ref{tbl:ppg-class-cpx}, it is a simple task to classify irreducible \ppgs\ over $\bR$: we look for real forms where $\fq$ is formed by crossing a single white nodes in $\fh$, which is not joined to any other nodes by an arrow.  Tables of simple real Lie algebras may be found in \cite[Tbl.\ B.4]{cs2009-parabolic1}.  This classification is also summarised in Table \ref{tbl:app-tbl-real}.

\begin{typeblock}{C}{n+1} There is a single permitted real form, the split real form $\fh = \alg{sp}{2n+2,\bR}$ of $\alg{sp}{2n+2,\bC}$.  Then $\fg = \alg{sl}{n+1,\bR}$, with Satake diagrams
\begin{equation*}
\smash{
  \fh = \dynkinC{2}{0}{3}{ppg-class-real-Ch}
  \quad\text{and}\quad
  \fg = \dynkinSLR{2}{0}{2}{ppg-class-real-Cg}.
}
\end{equation*}
The \Rspace\ $G\acts\fp$ is the real \proj\ space $\RP$ while $\bW = \Symm{2}\bR^{n+1}$ is the space of symmetric real $(n+1) \by (n+1)$ matrices, so that this \ppg\ is \proj\ differential geometry as studied in Chapter \ref{c:proj}. \end{typeblock}

\begin{typeblock}{A}{2n+1} There are two permitted real forms:
\begin{slimitemize}
  \item The split real form $\fh = \alg{sl}{2n+2,\bR}$ of $\alg{sl}{2n+2,\bC}$, so that $\fg = \alg{sl}{n+1,\bR} \dsum \alg{sl}{n+1,\bR}$ and the corresponding Satake diagrams are
  \vspace{0.4em}
  \begin{equation*}
  \smash{
    \fh = \dynkinASU{2}{0}{2}{ppg-class-real-A1h}
    \,\text{and}\quad
    \fg = \dynkinAA{2}{0}{2}{ppg-class-real-A1g}
    \colvectpunct[-0.75em]{~.}
  }
  \end{equation*}
  The \Rspace\ $G\acts\fp$ is $\RP[n] \times \RP[n]$, a product of \proj\ spaces; the resulting geometry should perhaps be called ``para-\cproj'' geometry.  The Jordan algebra $\bW$ is the external tensor product $\bR^{n+1} \etens \conj{\bR^{n+1}}$.
  
  \item The real form $\fh = \alg{su}{n+1,n+1}$ of $\alg{sl}{2n+2,\bC}$, so that $\fg = \alg{sl}{n+1,\bC}$ is the underlying real Lie algebra and the corresponding Satake diagrams are
  \begin{equation*}
    \fh = \dynkinSU{2}{0}{2}{ppg-class-real-A2h}
    \,\text{and}\quad
    \fg = \dynkinSLC{2}{0}{2}{ppg-class-real-A2g}
    \colvectpunct[-0.7em]{\, .}
  \end{equation*}
  The \Rspace\ $G\acts\fp$ is complex \proj\ space $\CP[n]$, viewed as a real manifold with complex structure, corresponding to \cproj\ geometry as studied in Chapter \ref{c:cproj}.  The Jordan algebra $\bW$ is the real representation underlying $\bC^{n+1} \etens \conj{\bC^{n+1}}$, which may be viewed as the space of $(n+1) \by (n+1)$ hermitian matrices.
\end{slimitemize}
\end{typeblock}
  
\begin{typeblock}{D}{2n+2} There are two suitable real forms:
\begin{slimitemize}
  \item The split real form $\fh = \alg{so}{4n+4,\bR}$ of $\alg{so}{4n+4,\bC}$, so that $\fg = \alg{sl}{2n+2,\bR}$ is the split real form of $\alg{sl}{2n+2,\bC}$ and the Satake diagrams of interest are
  \begin{equation*}
  \smash{
    \fh = \dynkinD{2}{0}{2}{ppg-class-real-D1h}
    \,\text{and}\quad
    \fg = \dynkinA{3}{0}{3}{ppg-class-real-D1g}.
  }
  \end{equation*}
  The \Rspace\ $G\acts\fp$ is the grassmannian of real $2$-planes in $\bR^{2n+2}$, with Jordan algebra $\bW = \Wedge{2} \bR^{n+1}$.
  
  \item The real form $\fh = \alg[^*]{so}{4n+4}$ of $\alg{so}{4n+4,\bC}$, with $\fg = \alg{sl}{n+1,\bH}$ and Satake diagrams
  \begin{equation*}
    \fh = \dynkinSOs{2}{0}{2}{ppg-class-real-D2h}
    \,\text{and}\quad
    \fg = \dynkinSLH{3}{0}{3}{ppg-class-real-D2g}.
  \end{equation*}
  The Lie algebra $\alg[^*]{so}{4n+4}$ is described in \eqref{eq:qtn-bgg-so*}.  The \Rspace\ $G\acts\fp$ is the \qtn ic\ \proj\ space $\HP$, so that this \ppg\ is almost \qtn ic\ geometry as studied in Chapter \ref{c:qtn}.  The Jordan algebra $\bW$ is the real representation underlying $\Wedge[\bC]{2} \bC^{2n+2}$, which may be viewed as the space of $(n+1) \by (n+1)$ \qtn-hermitian matrices.
\end{slimitemize}
\end{typeblock}

\begin{typeblock}{E}{7} There are again two permitted real forms:
\begin{slimitemize}
  \item The split real form $\fh = \alg[_{7(7)}]{e}{}$ of $\alg[_7]{e}{\bC}$ (sometimes called $\erealform{V}$), so that $\fg = \alg[_{6(6)}]{e}{}$ (sometimes called $\erealform{I}$) and the Satake diagrams are
  \begin{equation*}
  \smash{
    \fh = \dynkinE{7}{ppg-class-real-EVh}
    \quad\text{and}\quad
    \fg = \dynkinE{6}{ppg-class-real-EVg}
    \colvectpunct[-0.7em]{\, .}
  }
  \end{equation*}
  The \Rspace\ $G\acts\fp$ appears to not have a name, but could perhaps be understood in terms of Freudenthal's magic square \cite{f1955-magicsquare}.   
  
  \item The real form $\fh = \alg[_{7(-25)}]{e}{}$ of $\alg[_7]{e}{\bC}$ (sometimes called $\erealform{VII}$), so that $\fg = \alg[_{6(-26)}]{e}{}$ (sometimes called $\erealform{IV}$) with Satake diagrams
  \begin{equation*}
  \smash{
    \fh = \dynkinEVII{ppg-class-real-EVIIh}
    \quad\text{and}\quad
    \fg = \dynkinEIV{ppg-class-real-EVIIg}
    \colvectpunct[-0.7em]{\, .}
  }
  \end{equation*}
  The \Rspace\ $G\acts\fp$ is the \oct ic\ (Cayley) plane $\OP$, a \non Desarguesian\ \proj\ plane discovered by Moufang \cite{m1933-cayleyplane}; see also Baez's treatise \cite{b2002-octonions} on the \oct s.  The Jordan algebra $\bW$ is the exceptional Albert algebra
  \begin{equation*}
    \alg{alb}{3} \defeq \Setof{ \left( \begin{smallmatrix}
      a             & \alpha        & \beta  \\
      \conj{\alpha} & b             & \gamma \\
      \conj{\beta}  & \conj{\gamma} & c
    \end{smallmatrix} \right)
      }{ a,b,c \in \bR, ~
         \alpha,\beta,\gamma \in \bO }
  \end{equation*}
  of $3 \by 3$ \oct-hermitian matrices.
\end{slimitemize}
\end{typeblock}

\begin{typeblock}{BD}{n+4} There are a number of real forms, parameterised by the number of white nodes of the Satake diagram.  If there are $p$ white nodes, we have the indefinite real form $\fh = \alg{so}{p+2,n-p+2}$ of $\alg{so}{n+4,\bC}$.  Then $\fg$ is formed by removing the left-most white node, giving $\fg = \alg{so}{p+1,n-p+1}$ and Satake diagrams
\begin{equation*}
  \fh = \dynkinSOpq{3}{2}{p+2}{}
  \,\text{and}\quad
  \fg = \dynkinSOpq{2}{2}{p+1}{}
  \colvectpunct[-1.7em]{\! .}
\vspace{-0.6em}
\end{equation*}
The \Rspaces\ $H\acts\fq$ and $G\acts\fp$ are conformal spheres of signatures $(p+1,n-p+1)$ and $(p,n-p)$ respectively.  Note that we allow $p=0$ and $p=n+2$, leading to positive-definite and negative-definite conformal geometries.  The Jordan algebra $\bW$ is the indefinite real inner product space $\bR^{p+1, \, n+1-p}$, equipped with Clifford multiplication. \end{typeblock}

This classification is summarised in Table \ref{tbl:app-tbl-real}.  Except in Subsection \ref{ss:mob2-evals-riem}, we will not need to distinguish between different real forms of each type.

\begin{rmk} \thlabel{rmk:ppg-class-formal} Notice that our classification includes the compact rank one \riem\ symmetric spaces $\RP[n]$, $\CP[n]$, $\HP[n]$, $\OP[2]$ and $\Sph[n]$, which are precisely the flat models $G\acts\fp$ admitting a positive definite metric.  By a result of Hirzebruch \cite{h1965-riemsymm}, we obtain a one-to-one correspondence between rank one \riem\ symmetric spaces and formally real Jordan algebras.  Roughly, the \riem\ metric on $G\acts\fp$ has isometry group a maximal compact subgroup $K\leq G$, giving $G\acts\fp \isom G/P \isom K/(K \intsct P)$ by the second isomorphism theorem.  Then the resulting Cartan decomposition $\fg = \fk \dsum \fm$ and homogeneity of $G/P$ allow us to equip $\fm$ with the structure of a formally real Jordan algebra.  Koecher's theorem establishes a one-to-one correspondence between formally real Jordan algebras and so-called symmetric cones; $G\acts\fp$ may be viewed as the space of primitive idempotents in the corresponding cone.

In particular, we could have classified the flat models which admit positive definite metrics either via Cartan's classification of rank one \riem\ symmetric spaces \cite{c1926-riemsymm1, c1927-riemsymm2} or by Jordan, von Neumann and Wigner's classification of formally real Jordan algebras \cite{jnw1934-jordanqm}.  The data describing $G\acts\fp$ as a \riem\ symmetric space is listed in Table \ref{tbl:ppg-class-riem}. \end{rmk}

\begin{table}[!h]
  \begin{tabular}{| c | *{5}{>{$}c<{$}|} } \hline
    Type & \fh & G & \bW & G\acts\fp & K \\ \hline
    \type{C}{n+1} & \alg{sp}{2n+2,\bR} & \grp{PGL}{n+1,\bR} & \Symm{2}\bR^{n+1} & \RP[n] & \grp{PSO}{n+1} \\
    \type{A}{2n+1} & \alg{su}{n+1,n+1} & \grp{PGL}{n+1,\bC} & \realrepn{(\bC^{n+1} \etens \conj{\bC^{n+1}})} & \CP[n] & \grp{PSU}{n+1} \\
    \type{D}{2n+2} & \alg[^*]{so}{4n+4} & \grp{PGL}{n+1,\bH} & \realrepn{(\Wedge[\bC]{2} \bC^{2n+2})} & \HP[n] & \grp{PSp}{n+1} \\
    \type{E}{7} & \alg[_{7(-25)}]{e}{} & \grp[_{6(-26)}]{E}{} & \erepn[\bR] & \OP[2] & \grp[_{4(-52)}]{F}{} \\
    \type{BD}{n+4} & \alg{so}{n+2,2} & \grp{SO}{n+1,1} & \bR^{n+1, \, 1} & \Sph[n] & \grp{SO}{n+1} \\ \hline
  \end{tabular}
  \medskip
  \caption[The flat models admitting a positive definite metric]
          {The flat models $G\acts\fp$ admitting a positive definite metric, together with the corresponding Jordan algebra $\bW$ and maximal compact subgroup $K \leq G$.}
  \label{tbl:ppg-class-riem}
\end{table}
\vspace{-0.5\baselineskip}

\section{Calculus and associated BGG operators} 
\label{s:ppg-bgg}

In this subsection, we develop the theory of the important BGG operators for a \ppg.  For the classical \proj\ structures, the first BGG operator associated to $\bW$ controls the family of compatible metrics, while the first BGG operator associated to $\bW^*$ is a hessian-type equation.

We begin by deriving some curvature identities in Subsection \ref{ss:ppg-bgg-calc}, some of which hold for general abelian parabolic geometries and others which result from the $\bZ^2$-grading of Table \ref{tbl:ppg-alg-Z2gr}.  This allows us to describe the first BGG operator associated to $\bW$, which again controls the space of compatible metrics.  We do this in Subsection \ref{ss:ppg-bgg-metric}, in particular obtaining a relatively explicit prolongation.

In Subsection \ref{ss:ppg-bgg-nric} we study the normalised Ricci tensor of a compatible metric, in particular showing that it defines a section of $\cB^*$.  In \cproj\ geometry, this corresponds to the fact that the normalised Ricci tensor of a \kahler\ metric is symmetric and $J$-invariant.  We also characterise so-called \emph{normal solutions}, which correspond to compatible \einstein\ metrics.
Finally, we describe the first BGG operator associated to $\bW^*$ in Subsection \ref{ss:ppg-bgg-hess}, which is again a hessian-type equation.  Notably, every compatible metric determines a solution of the hessian.

We fix a \ppg\ over $M$ with parameters $(r,n)$, where we continue to assume that $n>0$.

\subsection{\Proj\ parabolic calculus} 
\label{ss:ppg-bgg-calc}

\thref{thm:para-bgg-ablcurv} says that the curvature tensor $\Curv{}{}$ of a Weyl structure $\D$ decomposes into Weyl and normalised Ricci parts according to
\vspace{0.1em}
\begin{equation} \label{eq:ppg-bgg-curvdecomp}
  \Curv[\D]{}{} = \Weyl[\D]{}{} - \algbracw{\id}{\nRic{}{}}{},
\vspace{0.1em}
\end{equation}
where by construction $\liebdy \Weyl[\D]{}{} = 0$ and $\nRic{}{} \defeq -\quab_M^{-1} \liebdy \Curv{}{}$.  The algebraic work of Section \ref{s:ppg-alg} allows us to extract other useful curvature identities.  For this we fix a local frame $\{e_i\}_i$ of $M$ with dual coframe $\{\ve^i\}_i$, often omitting the summation symbol when summing over $i = 1, \ldots, rn$.

\begin{prop} \thlabel{prop:ppg-bgg-Wbianchi} There is the Bianchi identity
\vspace{0.1em}
\begin{equation} \label{eq:ppg-bgg-Wbianchi}
  \Weyl[\D]{X,Y}{Z} + \Weyl[\D]{Y,Z}{X} + \Weyl[\D]{Z,X}{Y} = (\d^{\D} \Tor{})_{X,Y,Z}
\vspace{0.1em}
\end{equation}
\wrt\ any Weyl structure. \end{prop}

\begin{proof} Since $TM$ is abelian, \eqref{eq:ppg-bgg-curvdecomp} and the Bianchi identity for $\Curv{}{}$ give the result. \end{proof}

\begin{cor} \thlabel{cor:ppg-bgg-Wtrace} The Weyl curvature satisfies $\tr(\Weyl[\D]{X,Y}{}) = \ve^i( (\d^{\D} \Tor{})_{e_i,X,Y} )$.  In particular, if $\Tor{} = 0$ then $\Weyl[\D]{}{}$ is totally \tracefree. \end{cor}

\begin{proof} The result follows immediately from $\liebdy \Weyl[\D]{}{} = 0$ and \thref{prop:ppg-bgg-Wbianchi}. \end{proof}

By \thref{cor:ppg-alg-connsonL}, every nowhere-vanishing section $\ell$ of $\cL$ uniquely determines a Weyl structure $\D^{\ell}$ by decreeing that $\D^{\ell} \ell = 0$.  If $\ell_1, \ell_2$ are two such sections related by $\ell_1 = e^f \ell_2$ for some $f \in \s{0}{}$, a general section is of the form $\ell = h\ell_1 = he^f \ell_2$ for $h \in \s{0}{}$.  Then $\D^{\ell_1} \ell = \d h \tens \ell_1$ and
\vspace{0.1em}
\begin{equation*}
  \D^{\ell_2} \ell = \d(he^f) \tens \ell_2 = \d h \tens \ell_1 + \d f \tens \ell
    = \D^{\ell_1} \ell + \d f \tens \ell,
\vspace{0.1em}
\end{equation*}
so that $\D^{\ell_2} = \D^{\ell_1} + \algbrac{}{\d f}$ as Weyl connections.  In particular $\D^{\ell_1}$ and $\D^{\ell_2}$ are related by an exact $1$-form, so the class of Weyl connections obtained in this manner form an affine space modelled on the exact $1$-forms; for this reason they are referred to as \emph{exact Weyl structures} in the literature \cite{cs2003-weylstr, cs2009-parabolic1}.

\begin{cor} \thlabel{cor:ppg-bgg-exactweyl} $\nRic[\D^{\ell}]{X}{Y} - \nRic[\D^{\ell}]{Y}{X} = -\tr(\Weyl[\D^{\ell}]{X,Y}{})$ for every exact Weyl structure $\D^{\ell} \in \Dspace$.  In particular, if $\Tor{} = 0$ then $\nRic[\D^{\ell}]{}{}$ is symmetric. \end{cor}

\begin{proof} Since $\D^{\ell} \ell = 0$ by construction, $\Curv[\D^{\ell}]{X,Y}{\ell} = 0$ also and hence $\algbracw{\id}{\nRic[\D^{\ell}]{}{}}{X,Y} \acts \ell = \Weyl[\D^{\ell}]{X,Y}{\ell}$ by the decomposition \eqref{eq:ppg-bgg-curvdecomp}.  Applying \thref{lem:ppg-alg-Lbrac,cor:ppg-alg-trace} then proves the first claim, while the second follows from \thref{cor:ppg-bgg-Wtrace}. \end{proof}

The differential Bianchi identity $\d^{\D}\Curv[\D]{}{} = 0$ also yields some useful curvature identities, generalising \thref{prop:proj-para-calc,prop:cproj-para-calc,prop:qtn-para-calc}.  Since we will soon impose the condition $\Tor{} = 0$ anyway, we make this assumption now for convenience.  In this case \itemref{thm:para-bgg-ablcurv}{weyld} says that the Weyl curvature $\Weyl[\D]{}{}$ is invariant, so we write $\Weyl{}{} \defeq \Weyl[\D]{}{}$.

\begin{prop} \thlabel{prop:ppg-bgg-calc} Suppose that $\Tor{} = 0$.  Then \wrt\ any Weyl structure:
\begin{enumerate}
  \item \label{prop:ppg-bgg-calc-db}
  $\d^{\D}\Weyl{}{} = -\algbracw{\id}{\CY{}{}}{} = -\liediff \CY{}{}$;
  
  \item \label{prop:ppg-bgg-calc-cyb}
  There is a Bianchi identity $\CY{X,Y}{Z} + \CY{Y,Z}{X} + \CY{Z,X}{Y} = 0$; and
  
  \item \label{prop:ppg-bgg-calc-dbt}
  $\ve^i \left( (\D_{e_i} \Weyl{}{})_{X,Y} \right) = -(\quab \CY{}{})_{X,Y}$.
\end{enumerate}
\end{prop}

\begin{proof} \proofref{prop:ppg-bgg-calc}{db} Applying $\d^{\D}$ to the decomposition \eqref{eq:ppg-bgg-curvdecomp}, the differential Bianchi identity $\d^{\D}\Curv{}{} = 0$ yields $\d^{\D}\Weyl{}{} = \d^{\D} \algbracw{\id}{\nRic{}{}}{} = -\algbracw{\id}{\CY{}{}}{}$ as required.

\smallskip

\proofref{prop:ppg-bgg-calc}{cyb} Taking a trace of both sides in \ref{prop:ppg-bgg-calc-db}, the \lhs\ vanishes since $\Weyl{}{}$ is totally \tracefree\ by \thref{cor:ppg-bgg-Wtrace}.  On the \rhs\ we obtain
\begin{align*}
  0 &= \tr\left[ \algbrac{X}{\CY{Y,Z}{}} + \algbrac{Y}{\CY{Z,X}{}} + \algbrac{Z}{\CY{X,Y}{}} \right] \\
  &= \tfrac{1}{2}r(n+1) \left[ \CY{X,Y}{Z} + \CY{Y,Z}{X} + \CY{Z,X}{Y} \right]
\end{align*}
by \thref{cor:ppg-alg-trace}.  Since $r(n+1) \neq 0$, the result follows.

\smallskip

\proofref{prop:ppg-bgg-calc}{dbt} We apply $\liebdy$ to the differential Bianchi identity.  On the \lhs\ we obtain
\begin{align*}
  \liebdy( \d^{\D}\Weyl[\D]{}{} )_{X,Y}
    &= \sum{i}{} \, \ve^i \acts (\d^{\D}\Weyl[\D]{}{})_{e_i,X,Y} \\
    &= \sum{i}{} \, \ve^i \acts \big( \D_{e_i}\Weyl[\D]{X,Y}{} + \D_X\Weyl[\D]{Y,e_i}{}
       + \D_Y\Weyl[\D]{e_i,X}{} \\
    &\hspace{5em}
       - \Weyl[\D]{\liebrac{e_i}{X},Y}{} - \Weyl[\D]{\liebrac{X}{Y},e_i}{}
       - \Weyl[\D]{\liebrac{Y}{e_i},X}{} \big) \displaybreak \\
    &= \sum{i}{} \, \ve^i \acts \big( (\D_{e_i}\Weyl[\D]{}{})_{X,Y} + (\D_X\Weyl[\D]{}{})_{Y,e_i}
       + (\D_Y\Weyl[\D]{}{})_{e_i,X} \big) \\
    &= \sum{i}{} \, \ve^i \acts (\D_{e_i}\Weyl{}{})_{X,Y}
\end{align*}
since $\liebdy \Weyl{}{} = 0$.  On the \rhs, since $T^*M$ acts trivially on itself, we have
\begin{equation*}
  \liebdy \algbracw{\id}{\CY{}{}}{} = \liebdy \liediff \CY{}{} = \quab \CY{}{}.
\end{equation*}
Therefore $\ve^i\left( (\D_{e_i}\Weyl{}{})_{X,Y} \right) = \liebdy(\d^{\D}\Weyl{}{})_{X,Y} = -(\quab \CY{}{})_{X,Y}$ as required. \end{proof}

Here $\CY{}{}$ is a section of $\Wedge{2}T^*M \tens T^*M$, which in general is not an irreducible $P$-bundle.  Therefore $\quab$ may scale each component independently, as we see in \itemref{prop:cproj-para-calc}{db} for \cproj\ geometry; in particular, there is not a general formula for $\quab \CY{}{}$ in terms of the integers $(r,n)$.

We previously described the components of the harmonic curvature for the classical \proj\ structures, and completeness demands that we do the same for types \type{E}{7} and \type{BD}{n+4}.  By Subsection \ref{ss:ppg-class-cpx}, the Hasse diagrams computing the homologies $\liehom{}{\fg}$ are given by Figures \ref{fig:ppg-bgg-hasseE} and \ref{fig:ppg-bgg-hasseBD} respectively.

\vspace{0.4em}
\begin{figure}[h]
  \begin{equation*}
  \arraycolsep=0.2em
  \begin{array}{*{6}{c}}
    \dynkinEp[0,0,0,0,0,1]{6}{ppg-bgg-metric-hasseE-1}
    &
    \dynkin{ \DynkinConnector{0.2}{0.3}{1.8}{0.3}; }{}
    &
    \dynkinEp[0,0,0,1,-2,1]{6}{ppg-bgg-metric-hasseE-3}
    &
    \dynkin{ \DynkinConnector{0.2}{0.3}{1.8}{0.3}; }{}
    &
    \dynkinEp[0,0,1,0,-3,1]{6}{ppg-bgg-metric-hasseE-5}
    &
    \dynkin{ \DynkinConnector{0.2}{0.3}{1.8}{0.3};
             \DynkinLabel{$\cdots$}{2.8}{-0.4}; }{ppg-bgg-metric-hasseE-6}
  \end{array}  
  \end{equation*}
  \vspace{-0.5em}
  \caption[The Hasse diagram computing $\liehom{}{\alg[_6]{e}{\bC}}$]
          {The Hasse diagram computing the homology $\liehom{}{\alg[_6]{e}{\bC}}$.}
  \label{fig:ppg-bgg-hasseE}
\end{figure}

\vspace{-1em}

\begin{figure}[h]
  \begin{multline*}
  \arraycolsep=0.2em
    \hspace{2.7em}
    \begin{array}{c} \dynkinBDp[0,1,0,0,0,0,0,0,0]{3}{0}{2}
                               {ppg-bgg-metric-hasseBD-1} \end{array}
    \dynkin{ \DynkinConnector{0.2}{0.68}{1.3}{0.68}
             \DynkinLabel{}{0.75}{-0.05} }{ppg-bgg-metric-hasseBD-2}
    \begin{array}{c} \dynkinBDp[-2,2,0,0,0,0,0,0,0]{3}{0}{2}
                               {ppg-bgg-metric-hasseBD-3} \end{array}
    \\[-3.68em]
    \hspace{3em}
    \zbox{ \dynkin{ \draw[thin]
      ( 4*\DynkinStep, 2*\DynkinStep)
        -- ( 5*\DynkinStep, 2*\DynkinStep)
      ( 5*\DynkinStep, 2*\DynkinStep) arc(90:-90:1.45*\DynkinStep)
           ( 4*\DynkinStep, 0*\DynkinStep)
      ( 5*\DynkinStep,-0.9*\DynkinStep)
        -- (-9*\DynkinStep,-0.9*\DynkinStep)
      (-9*\DynkinStep,-0.9*\DynkinStep) arc(-90:90:-0.9*\DynkinStep)
           (-9*\DynkinStep,-4*\DynkinStep)
      (-9*\DynkinStep,-2.7*\DynkinStep) edge[->]
           (-8*\DynkinStep,-2.7*\DynkinStep);
      \DynkinLabel{}{-7.5}{-3.25}
    }{ppg-bgg-metric-hasseBD-4} }
    \\[-3.6em]
    \begin{array}{c} \dynkinBDp[-4,0,2,0,0,0,0,0,0,0]{4}{0}{2}
                               {ppg-bgg-metric-hasseBD-5} \end{array}
    \dynkin{ \DynkinConnector{0.2}{0.68}{1.3}{0.68}
             \DynkinLabel{$\cdots$}{2.2}{-0.02} }
           {ppg-bgg-metric-hasseBD-6}
    \hspace{8em}
  \end{multline*}
  \vspace{-0.9em}
  \caption[The Hasse diagram computing $\liehom{}{\alg{so}{n+2,\bC}}$]
          {The Hasse diagram computing the homology $\liehom{}{\alg{so}{n+2,\bC}}$.}
  \label{fig:ppg-bgg-hasseBD}
\end{figure}
\vspace{0.2em}

We find a single component of harmonic curvature in type \type{E}{7}, which lives in the bundle associated to
\begin{align*}
  \dynkinEp[0,0,1,0,-3,1]{6}{ppg-bgg-calc-harmE-1}
    &= \dynkinEp[0,0,0,0,0,1]{6}{ppg-bgg-calc-harmE-2} \cartan[\bC]
      \Wedge[\bC]{2} \!\left( \dynkinEp[0,0,1,0,-2,0]{6}{ppg-bgg-calc-harmE-3} \right) \\
    &= \cpxbdl{( TM \cartan \Wedge{2}T^*M )}
\end{align*}
and may therefore be identified with the Cartan torsion $\Tor{}$.  In particular, the $\Tor{}$ is the only obstruction to local flatness for type \type{E}{7}, giving flat $M$ locally isomorphic to a real form of $\pr{\bO\tens\bC}$.  This makes the assumption $\Tor{} = 0$ rather unfortunate, but seemingly necessary to obtain any interesting results in the following subsections.

For type \type{BD}{n+4} we again find a single component of harmonic curvature, which of course coincides with (the complexification of) the conformal Weyl curvature $\Weyl{}{}$.  In particular, the torsion component $\Tor{}$ vanishes automatically.  It is well-known that $\Weyl{}{}$ gives a complete obstruction to local conformal flatness \cite{b1987-einstein, bc2010-conf, s1997-cartan}.

\vspace{0.3em}
\subsection{Metrisability of \ppgs} 
\label{ss:ppg-bgg-metric}
\vspace{0.2em}

Recall that metrisability of each classical \proj\ structure is controlled by an invariant first-order differential equation, which we called the \emph{linear metric equation}.  In each case, this equation was the first BGG operator associated to the irreducible $\fg$-representation $\bW$.

\begin{defn} \thlabel{defn:ppg-bgg-metriceqn} The differential equation $\bgg{\bW}(h) = 0$ determined by the first BGG operator $\bgg{\bW}$ associated to $\bW$ is called the \emph{linear metric equation}, and its solutions will be called \emph{linear metrics}. \end{defn}

We will see shortly that the (\non degenerate) solutions of $\bgg[\bW]{}$ correspond to metrics compatible with the underlying geometric structure.  For this, fix a Weyl structure and hence an isomorphism $\cW \isom (\cL^*\tens\cB) \dsum (\cL^*\tens TM) \dsum \cL^*$ by \thref{prop:ppg-alg-gr}.  Here $\cL^*\tens\cB \isom \liehom{0}{\cW}$ and $\cL \isom \liehom{0}{\cW^*}$, for $\cB$ associated to a $\fp^0$-subrepresentation $B \leq \Symm{2}(\fg/\fp)$.  By applying the results of Subsection \ref{ss:ppg-alg-Z2} to the associated bundles, sections $h$ of $\cL^*\tens\cB$ may be viewed as $\cL^*$-valued symmetric bilinear forms on $T^*M$, defined by $h(\alpha,\beta) \defeq \algbrac{ \algbrac{h}{\alpha} }{ \beta }$.  Here $\algbrac{}{}$ is the algebraic bracket on $\fh_M$.

\begin{defn} \thlabel{defn:ppg-bgg-nondeg} A section $h$ of $\cL^*\tens\cB$ is \emph{\non degenerate} if it is \non degenerate\ as an $\cL^*$-valued symmetric bilinear-form on $T^*M$. \end{defn}

Recall that \thref{prop:ppg-alg-Lwedge} provides an isomorphism $\cL \isom (\Wedge{rn}TM)^{2/r(n+1)}$.  Then $(\det h)^{1/r}$ is a section of $\Wedge{rn}(\cL^*\tens\cB)^{1/r} \isom \cL$, and evidently $h$ is \non degenerate\ if and only if $(\det h)^{1/r}$ is nowhere-vanishing.  In this case it follows that $g \defeq (\det h)^{1/r} \ltens h^{-1}$ is a \non degenerate\ section of $\Symm{2}T^*M$, \ie\ a metric on $M$.  Recall that an element $f \in \bW^*$ is called regular if $(\ad f)^2$ factors to an isomorphism $\bW \to \bW^*$.  \Non degenerate\ sections of $\cL^* \tens \cB$ are a source of regular elements.

\begin{prop} \thlabel{prop:ppg-bgg-nondeg} Let $h$ be a \non degenerate\ section of $\cL^* \tens \cB$.  Then at every point, $f \defeq h^{-1} + (\det h)^{1/r}$ is a regular element of $\bW^*$. \end{prop}

\begin{proof} Since $h$ is \non degenerate, $\lambda \defeq (\det h)^{-1/r}$ is a nowhere-vanishing section of $\cL^*$, thus yielding a section $\lambda^{-1}$ of $\cL$.  Applying \thref{prop:ppg-defn-sd} pointwise, if $\xi$ is the Weyl structure of $\fq$ inducing the splitting $\fh_M \isom \cW \dsum \fq^0_M \dsum \cW^*$ then it suffices to find a section $e$ of $\cW$ for which $\algbrac{e}{f} = 2\xi$; we show that $e \defeq h + \lambda$ is such a section.

We first compute the two brackets $\algbrac{h}{h^{-1}}$ and $\algbrac{\lambda}{\lambda^{-1}}$.  Since both are sections of $\fp^0_M \dsum \liecenter{\fq^0}_M$ by the $\bZ^2$-grading, we can write $\algbrac{h}{h^{-1}} = A + a\xi$ and $\algbrac{\lambda}{\lambda^{-1}} = B + b\xi$ for some $A,B \in \s{0}{\fp^0_M}$ and $a,b \in \s{0}{}$.  By Table \ref{tbl:ppg-alg-Z2gr} we have
\vspace{0.1em}
\begin{equation*} \begin{gathered}
  \algbrac{ \algbrac{h}{h^{-1}} }{ X }
    = -\algbrac{ h }{ h^{-1}(X,\bdot) }
    = -X \\
  \text{and} \quad
  \algbrac{ \algbrac{\lambda}{\lambda^{-1}} }{ X }
    = \algbrac{ \lambda \ltens X }{ \lambda^{-1} }
    = X
\end{gathered}
\vspace{0.1em}
\end{equation*}
for all $X \in \s{0}{TM}$, giving $A = -\id$ and $B = \id$ since $\xi$ acts trivially on $TM \leq \fq^0_M$.

The Jacobi identity implies that $\algbrac{h}{h^{-1}}$ acts trivially on sections $\ell$ of $\cL$.  Then
\vspace{0.2em}
\begin{equation*}
  0 = \algbrac{ -\id + a\xi }{ \ell }
    = \left( \tfrac{2}{r(n+1)} \tr(-\id) - a \right) \ell
    = \left( -\tfrac{2n}{n+1} - a \right) \ell
\vspace{0.1em}
\end{equation*}
by \thref{cor:ppg-alg-dims,cor:ppg-alg-trace}, whence $a = -\tfrac{2n}{n+1}$.  On the other hand, $\algbrac{\lambda}{\lambda^{-1}}$ acts trivially on sections $\theta$ of $\cL \tens \cB^*$.  Then for all $X,Y \in \s{0}{TM}$, we have
\vspace{0.1em}
\begin{equation*} \begin{aligned}
  \algbrac{ \algbrac{\lambda}{\lambda^{-1}} }{ \theta }(X,Y)
    &= \algbrac{ \algbrac{\lambda}{\lambda^{-1}} }{ \theta(X,Y) } \\
    &\qquad
      - \theta( \algbrac{ \algbrac{\lambda}{\lambda^{-1}} }{ X }, Y )
      - \theta( X, \algbrac{ \algbrac{\lambda}{\lambda^{-1}} }{ Y } ) \\
    &= \algbrac{ \id + b\xi }{ \theta(X,Y) } \\
    &\qquad
      - \theta( \algbrac{ \lambda \ltens X }{ \lambda^{-1} }, Y )
      - \theta( X, \algbrac{ \lambda \ltens Y }{ \lambda^{-1} } ) \\
    &= \left( \tfrac{2}{r(n+1)} (\tr \id) - b + 2 \right) \theta(X,Y) \\
    &= \left( \tfrac{2n}{n+1} - b + 2 \right) \theta(X,Y) = 0,
\end{aligned}
\vspace{0.1em}
\end{equation*}
whence $b = \tfrac{2n}{n+1} + 2$.
By Table \ref{tbl:ppg-alg-Z2gr} and the previous calculations, we then obtain
\vspace{0.1em}
\begin{equation*} \begin{aligned}
  \algbrac{ h+\lambda }{ h^{-1}+\lambda^{-1} }
    &= \algbrac{h}{h^{-1}} + \algbrac{\lambda}{\lambda^{-1}} \\
    &= \left( -\id - \tfrac{2n}{n+1}\xi \right)
      + \left( \id + \big(\tfrac{2n}{n+1} + 2\big)\xi \right)
    = 2\xi
\end{aligned}
\vspace{0.1em}
\end{equation*}
as required.  The result now follows by \thref{prop:ppg-defn-sd}. \end{proof}

\begin{rmk} \thlabel{rmk:ppg-bgg-nondeg} In the \non \self dual\ theory mentioned in \thref{rmk:ppg-alg-str,rmk:ppg-alg-Z2nonsd}, \thref{prop:ppg-bgg-nondeg} implies that $\cL^* \tens \cB$ admits \non degenerate\ sections only when $H\acts\fq$ is \self dual.  This is the first place that \self duality\ becomes a necessary assumption. \end{rmk}

Since we are looking for compatible metric connections in the class $\Dspace$ of Weyl connections, we shall assume henceforth that the torsion $\Tor{}$ of any Weyl connection vanishes.  As mentioned previously, this assumption has the regrettable side-effect of local flatness for \ppgs\ modelled on real forms of $\alg[_6]{e}{\bC}$.  Since $\weyld{\gamma} \Weyl[\D]{}{} = \algbrac{\Tor{}}{\gamma} = 0$, the Weyl curvature $\Weyl{}{} \defeq \Weyl[\D]{}{}$ is an invariant.

Our first task is to give an explicit formula for the linear metric equation.  For this, we shall need the tractor connection $\D^{\bW}$ on $\cW$ and its curvature $\Curv[\bW]{}{}$.

\begin{lem} \thlabel{lem:ppg-bgg-metrictractor} The tractor connection $\D^{\bW}$ on $\cW$ and its curvature $\Curv[\bW]{}{}$ may be written
\begin{align}
	\notag \\[-0.95\baselineskip]
  \label{eq:ppg-bgg-metrictractor}
  \D^{\bW}_X \colvect{ h \\ Z \\ \lambda }
  &= \colvect{ \D_X h - \algbrac{Z}{X} \\
               \D_X Z - \lambda\ltens X - h(\nRic{X}{},\bdot) \\
               \D_X \lambda - \nRic{X}{Z} } \\[0.5em]
  \label{eq:ppg-bgg-metriccurv}
  \Curv[\bW]{X,Y}{ \colvect{ h \\ Z \\ \lambda } }
  &= \colvect{ \Weyl[\D]{X,Y}{h} \\
               \Weyl[\D]{X,Y}{Z} - h(\CY{X,Y}{},\bdot) \\
               \CY{X,Y}{Z} }
	\\[-0.95\baselineskip] \notag
\end{align}
\wrt\ any Weyl structure. \end{lem}

\begin{proof} This is immediate from the general \formulae\ $\D^{\bW}_X s = X\acts s + \D_X s + \nRic{X}{}\acts s$ and, since the torsion is assumed to vanish, $\Curv[\bW]{X,Y}{s} = \Weyl[\D]{X,Y}{s} + \CY{X,Y}{}\acts s$. \end{proof}

Via the Weyl structure, the algebraic laplacian induces a bundle map $\quab : \cW \to \cW$ which acts on $\fp$-irreducible components according to Kostant's Spectral \thref{thm:lie-hom-spectral}.

\begin{lem} \thlabel{lem:ppg-bgg-Wquabla} $\quab$ acts trivially on $\cL^*\tens\cB$, by multiplication with $\tfrac{1}{2}(rn-r+2)$ on $\cL^*\tens TM$, and by multiplication with $rn$ on $\cL^*$. \end{lem}

\begin{proof} The first slot is clear, since $\cL^*\tens\cB$ is the zeroth homology of $\cW$.  Since $\liebdy$ acts trivially on $\cW$, the algebraic laplacian is given by $\quab s = \liebdy \liediff s$ for all $s \in \s{0}{\cW}$.  Then for the second slot,%
\footnote{We continue to suppress summation signs when contracting \wrt\ a local (co)frame.}
\vspace{0.1em}
\begin{equation*}
  \quab Z = \liebdy\liediff Z
    = \liebdy( \ve^i \tens \algbrac{e_i}{Z} )
    = \algbrac{ \ve^i }{ \algbrac{e_i}{Z} }
    = -\algbrac{\algbrac{e_i}{\ve^i}}{Z} - \algbrac{e_i}{\ve^i(Z)}
\vspace{0.1em}
\end{equation*}
Breaking $Z \in \s{0}{\cL^*\tens TM}$ into $\cL^*$- and $TM$-factors, \thref{cor:ppg-alg-trace} and the symmetry of the bracket imply that $\algbrac{ \algbrac{e_i}{\ve^i} }{ Z } = -\ve^i(e_i) Z + \tfrac{1}{2}r(n+1) Z = -\tfrac{1}{2} (rn-r) Z$.  Since $\algbrac{e_i}{\ve^i(Z)} = -Z$ by Table \ref{tbl:ppg-alg-Z2gr}, the claim for $\cL^*\tens TM$ follows.  Finally, we easily compute that $\quab$ acts on $\cL^*$  by
\vspace{0.1em}
\begin{equation*}
  \quab \lambda
    = -\algbrac{\ve^i}{\lambda \ltens e_i}
    = \ve^i(e_i \ltens \lambda)
    = rn \lambda
\vspace{0.1em}
\end{equation*}
by the previous results and \thref{cor:ppg-alg-dims}. \end{proof}

\begin{prop} The BGG splitting operator associated to $\bW$ is given by
\begin{equation} \label{eq:ppg-bgg-metricsplit}
  \bggrepr{\bW} : h \mapsto \colvect{ h \\ Z^{\D} \\ \lambda^{\D} }
    \defeq \colvect{ h \\
                     \tfrac{2}{rn-r+2} \ve^i(\D_{e_i}h) \\
                     \tfrac{2}{rn(rn-r+2)} (\D^2_{e_i,e_j} h)(\ve^i,\ve^j)
                       - \tfrac{1}{rn} h(\nRic{e_i}{},\ve^i) }
    \colvectpunct{.}
\end{equation}
\end{prop}

\begin{proof} By definition $\bggrepr{\bW} \defeq \bggpi{\bW} \circ \mr{repr} : \s{0}{\liehom{0}{\cW}} \to \s{0}{\cW}$ where, since $\im(\mr{repr}) \subseteq \ker\liebdy$, we have $\bggpi{\bW}(s) = s - \quab_M^{-1} \liebdy (\d^{\bW} s)$; see \eqref{eq:para-bgg-projrepr} The exterior covariant derivative $\d^{\bW}$ coincides with the tractor connection \eqref{eq:ppg-bgg-metrictractor} on $\cW$, so that
\vspace{-0.2em}
\begin{equation*}
  \liebdy \big( \d^{\bW} \mr{repr}(h) \big)
    = \liebdy \! \colvect{ \D h \\
                          -h(\nRic{}{},\bdot) \\
                          0 } \\
    = \colvect{ 0 \\
                \algbrac{\ve^i}{\D_{e_i}h} \\
                -\algbrac{\ve^i}{h(\nRic{e_i}{},\bdot)} }
    = \colvect{ 0 \\
                -\ve^i(\D_{e_i}h) \\
                h(\nRic{e_i}{},\ve^i) }
\vspace{-0.2em}
\end{equation*}
\wrt\ any local frame of $M$.  To compute the action of $\quab_M^{-1}$, we use the Neumann series \eqref{eq:para-bgg-neumann}.  \Wrt\ any Weyl structure, we have
\vspace{-0.2em}
\begin{equation*}
  (\quab_M - \quab) \quab^{-1} \liebdy \! \colvect{ \D h \\ -h(\nRic{}{},\bdot) \\ 0 }
    = \colvect{ 0 \\
                0 \\
                \tfrac{2}{rn-r+2} \ve^i(\D_{e_i} \ve^j(\D_{e_j}h) ) }
\vspace{-0.2em}
\end{equation*}
by \eqref{eq:para-bgg-quabladiff} and \thref{lem:ppg-bgg-Wquabla}.  Therefore
\begin{align*}
  \quab_M^{-1}\liebdy \! \colvect{ \D h \\ -h(\nRic{}{},\bdot) \\ 0 }
  &= \left( \id - \quab^{-1}(\quab_M - \quab) \right) \quab^{-1} \liebdy \colvect{ \D h \\ -h(\nRic{}{},\bdot) \\ 0 } \\
  &= \colvect{ 0 \\
               -\tfrac{2}{rn-r+2} \ve^i(\D_{e_i}h) \\
               \tfrac{1}{rn} h(\nRic{e_i}{},\ve^i)
                 - \tfrac{2}{rn(rn-r+2)} (\D^2_{e_i,e_j} h)(\ve^i,\ve^j) }
    \colvectpunct{.}
\end{align*}
Using that $\bggpi{\bW} = \id - \quab_M^{-1}\liebdy \circ \d^{\bW}$, the result follows. \end{proof}

By the general theory, the first BGG operator is given by $\bgg{\bW} \defeq \mr{proj} \circ \bggpi[\bW]{1} \circ \D^{\bW} \circ \bggrepr{\bW}$.  The projection $\bggpi[\bW]{1}$ acts by the identity on homology, so the linear metric equation is the projection of \eqref{eq:ppg-bgg-metricsplit} onto the \non zero\ graded component of highest weight.  For $r=1,2,4$ and $8$ this gives
\vspace{-0.15em}
\begin{equation} \label{eq:ppg-bgg-metriceqn1}
  \bgg{\bW}(h) = \D h - \algbrac{Z^{\D}}{},
\vspace{-0.15em}
\end{equation}
where $Z^{\D} \defeq \tfrac{2}{rn-r+2} \, \ve^i(\D_{e_i}h) \in \s{0}{\cL^*\tens TM}$ is a \non zero\ multiple of $\liebdy(\D h)$.  Equivalently, we may write
\vspace{-0.15em}
\begin{equation*}
  \D h = 0 \mod \cL^*\tens TM,
\vspace{-0.15em}
\end{equation*}
where we identify $\cL^*\tens TM$ with its injective image in $\cL^*\tens\cB$ via the bracket $Z \mapsto \algbrac{Z}{}$.

Equation \eqref{eq:ppg-bgg-metriceqn1} is satisfied tautologically when $r=n$: there $\cL^*\tens\cB \isom \cL$ is $1$-dimensional, spanned by the (inverse) conformal metric $\conf$, and the Weyl connections are conformal in the sense that $\D \mathtt{c} = 0$ for all $\D \in \Dspace$.  Therefore if $h = f\conf$ we have $\algbrac{\ve^i(\D_{e_i} h)}{} = \algbrac{\conf(\d f,\bdot)}{} = \d f \tens \conf$ via the isomorphism $\cL^*\tens\cB \isom \cL$ provided by $\conf$.  Therefore the first BGG equation is given by
\vspace{0.25em}
\begin{equation} \label{eq:ppg-bgg-confmetric} \begin{aligned}
  \bgg{\bW^*}(h)
    &= \big( (\D^2_{e_i,\bdot} h)(\ve^i,\bdot) - h(\nRic{}{},\bdot) \big) \\
    &\qquad
      - \tfrac{1}{n} \big( (\D^2_{e_i,e_j} h)(\ve^i,\ve^j) - h(\nRic{e_i}{},\ve^i) \big) \id.
\end{aligned}
\vspace{0.25em}
\end{equation}
Identifying $\cL^*\tens\cB \isom \cL$ as above, \eqref{eq:ppg-bgg-confmetric} may be written as $\bgg{\bW^*}(\ell) = (\D^2 \ell - \ell \ltens \nRic{}{})_{\trfree}$, where the subscript ``$\trfree$'' denotes the \tracefree\ part \wrt\ $\conf$.  This is the \emph{\einstein\ scale equation} from conformal geometry, whose solutions parametrise \einstein\ metrics in the conformal class.  In particular, the conformal metric equation has more in common with the previous hessian-type equations; we shall explain this in Subsection \ref{ss:ppg-bgg-hess}.

The invariance of \eqref{eq:ppg-bgg-metriceqn1} and \eqref{eq:ppg-bgg-confmetric} of course follow from the fact that they are first BGG operators, so there is no need to do calculations as in \thref{prop:proj-bgg-metriceqn,prop:cproj-bgg-metriceqn,prop:qtn-bgg-metriceqn}.  The variations of the quantities $(h,Z^{\D},\lambda^{\D})$ \wrt\ a Weyl structure are also easily calculated: one finds that
\vspace{0.25em}
\begin{equation} \label{eq:ppg-bgg-metricweyld}
  \weyld{\gamma} Z^{\D} = h(\gamma,\bdot)
  \quad\text{and}\quad
  \weyld{\gamma} \lambda^{\D} = \gamma(Z^{\D}).
\vspace{0.25em}
\end{equation}
By the Taylor expansion \eqref{eq:para-calc-taylor}, it follows that $Z^{\D} \mapsto Z^{\D} + h(\gamma,\bdot)$ and $\lambda^{\D} \mapsto \lambda^{\D} + \gamma(Z^{\D}) + \tfrac{1}{2}h(\gamma,\gamma)$ under change of Weyl structure according to $\D \mapsto \D + \algbrac{}{\gamma}$.

\bgroup
\addtolength{\topsep}{0.5em}
\begin{cor} \thlabel{cor:ppg-bgg-lcconn} There is a bijection between \non degenerate\ solutions $h$ of \eqref{eq:ppg-bgg-metriceqn1} and metric connections in the class $\Dspace$ of Weyl connections. \end{cor}
\egroup

\begin{proof} If $h \in \s{0}{\cL^*\tens\cB}$ is a \non degenerate\ solution of \eqref{eq:ppg-bgg-metriceqn1} with $\D h = \algbrac{Z^{\D}}{}$, then $h^{-1}(Z^{\D},\bdot)$ is a $1$-form.  Hence $\D^g \defeq \D - \algbrac{}{h^{-1}(Z^{\D},\bdot)} \in \Dspace$ is independent of $\D \in \Dspace$ and satisfies $\D^g h = 0$, so that $\D^g$ is the \LC\ connection of the corresponding metric $g \defeq (\det h)^{1/r} \ltens h^{-1} \in \s{0}{\Symm{2}T^*M}$.  Conversely if $\D^g \in \Dspace$ is the \LC\ connection of a metric $g \in \s{0}{\cB^*}$, we have $\D^g h = 0$ for $h \defeq (\det g)^{1/r(n+1)} \ltens g^{-1} \in \s{0}{\cL^*\tens\cB}$, meaning that $h$ is a solution of \eqref{eq:ppg-bgg-metriceqn1}. \end{proof}

Since the linear metric equation is a first BGG operator, its solutions are in bijection with the parallel sections of a prolongation operator $\d^{\cW}$.  Thanks to Table \ref{tbl:ppg-alg-Z2gr}, we can compute this prolongation fairly explicitly.  By the differential identity from \itemref{prop:ppg-bgg-calc}{dbt}, it turns out that in fact $\d^{\cW}$ is a prolongation connection $\D^{\cW}$.

\begin{thm} \thlabel{thm:ppg-bgg-metricprol} There is a linear isomorphism between solutions of the linear metric equation and the parallel sections of the prolongation connection
\begin{equation} \label{eq:ppg-bgg-metricprol}
  \D^{\cW}_X \! \colvect{ h \\ Z \\ \lambda }
    = \colvect{ \D_X h - \algbrac{Z}{X} \\
                \D_X Z - \lambda \ltens X - h(\nRic{X}{},\bdot) \\
                \D_X \lambda - \nRic{X}{Z} }
    - \colvect{ 0 \\
                w \, \liebdy \algbrac{\Weyl{}{}}{h}_X \\
                \quab^{-1}\liebdy \algbrac{ \CY{}{} - w \, \quab\CY{}{} }{h}_X }
\end{equation}
on sections of $\cW \isom (\cL^*\tens\cB) \dsum (\cL^*\tens TM) \dsum \cL^*$, where $w \in \bR$ is the eigenvalue of $\quab^{-1}$ on $\liebdy \algbrac{\Weyl{}{}}{h}$.  The isomorphism is given explicitly by the splitting operator \eqref{eq:ppg-bgg-metricsplit}. \end{thm}

\begin{proof} We apply the general prolongation procedure.  Equation \eqref{eq:ppg-bgg-metriccurv} gives
\vspace{-1em}
\begin{equation*}
  \liebdy \! \left( \Curv[\bW]{}{ \colvect{h \\ Z \\ \lambda} } \right)
    = \colvect{ 0 \\
                \liebdy \algbrac{\Weyl{}{}}{h} \\
                \liebdy \algbrac{\Weyl{}{}}{Z} + \liebdy \algbrac{\CY{}{}}{h} }
    \colvectpunct{.}
\end{equation*}
The term $\liebdy \algbrac{\Weyl{}{}}{Z}$ vanishes, since
$
  \sum{i}{} \, \algbrac{ \ve^i }{ \Weyl{e_i,X}{Z} } = -\sum{i}{} \ve^i( \Weyl{e_i,X}{Z} ) = 0.
$
Next we have
\begin{align*}
  (\quab_M - \quab) \, \liebdy \! \left( \Curv{}{ \colvect{h \\ Z \\ \lambda} } \right)
    &= \ve^i \acts \left( (\D_{e_i} + \nRic{e_i}{})
      \colvect{ 0 \\
                \liebdy \algbrac{\Weyl{}{}}{h} \\
                \ast }
      \right)
    = \colvect{ 0 \\
                 0 \\
                 \liebdy( \D \liebdy \algbrac{\Weyl{}{}}{h} ) }
    \colvectpunct{,}
\end{align*}
where we do not care about the $\cL^*$-slot since it ``drops off the bottom'' upon application of $\ve^i$ anyway.  The Neumann series then gives
\vspace{0.2em}
\begin{equation*}
  \quab_M^{-1} \liebdy \! \left( \Curv{}{ \colvect{h \\ Z \\ \lambda} } \right)
    = \colvect{ 0 \\
                \quab^{-1} \liebdy \algbrac{\Weyl{}{}}{h} \\
                \quab^{-1} \liebdy \algbrac{\CY{}{}}{h}
                  - \liebdy( \D \quab^{-1}\liebdy \algbrac{\Weyl{}{}}{h}) }
    \colvectpunct{.}
\vspace{0.2em}
\end{equation*}
We next simplify the term $\liebdy(\D \quab^{-1} \liebdy \algbrac{\Weyl{}{}}{h})$.  Since $\quab^{-1}$ acts on $\liebdy \algbrac{\Weyl{}{}}{h}$ by a scalar $w \in \bR$, it suffices to consider $\liebdy(\D \liebdy \algbrac{\Weyl{}{}}{h})$, which equals
\begin{align} \label{eq:ppg-bgg-metricprol-1}
  &\liebdy(\D \liebdy \algbrac{\Weyl{}{}}{h})_X
    \notag \\
  &\quad
    = \sum{i,j}{} \, \ve^i \acts \big(
      \D_{e_i} ( \ve^j \acts \algbrac{\Weyl{e_j,X}{}}{h} )
       -  \ve^j \acts \algbrac{\Weyl{e_j,\D_{e_i}}{X}}{h} \big)
       \notag \\
  &\quad
    = -\sum{i,j}{} \, \ve^i \acts \big(
      \D_{e_i}( \Weyl{e_j,X}{h(\ve^j,\bdot)} ) - \Weyl{e_j,\D_{e_i}X}{h(\ve^j,\bdot)} \big)
      \notag \\
  &\quad
    = -\sum{i,j}{} \, \ve^i \acts \big(
      (\D_{e_i}\Weyl{}{})_{e_j,X} \acts h(\ve^j,\bdot)
      + \Weyl{e_j,X}{ (\D_{e_i}h)(\ve^j,\bdot) } \big)
      \notag \\
  &\quad
    = -\sum{i,j}{} \, \Big(
      \algbrac{ \ve^i }{ \algbrac{(\D_{e_i}\Weyl{}{})_{e_j,X}}{h(\ve^j,\bdot)} }
      + \algbrac{ \ve^i }{ \Weyl{e_j,X}{(\D_{e_i}h)(\ve^j,\bdot)} } \Big)
      \notag \\
  &\quad \begin{aligned}
    &\!= -\sum{i,j}{} \, \Big(
        \algbrac{ \algbrac{ \ve^i }{ (\D_{e_i}\Weyl{}{})_{e_j,X} } }
                { h(\ve^j,\bdot) }
      + \algbrac{ (\D_{e_i}\Weyl{}{})_{e_j,X} }
                { \algbrac{ \ve^i }{ h(\ve^j,\bdot) } } \\
    &\hspace{5.5em}
      + \algbrac{ \ve^i }
                { \Weyl{e_j,X}{ \algbrac{Z^{\D}}{e_i}(\ve^j,\bdot) }} \Big).
  \end{aligned}
\end{align}
The second term on the \rhs\ vanishes, since
\begin{equation*} \begin{aligned}
  \algbrac{ (\D_{e_i}\Weyl{}{})_{e_j,X} }
          { h(\ve^i,\ve^j) }
    &= \algbrac{ \D_{e_i}\Weyl{e_j,X}{} - \Weyl{e_j,\D_{e_i}X}{} }
               { h(\ve^i,\ve^j) } \\
    &= \D_{e_i} \algbrac{ \Weyl{e_j,X}{} }{ h(\ve^i,\ve^j) }
      - \algbrac{ \Weyl{e_j,X}{} }{ \D_{e_i} h(\ve^i,\ve^j) } \\
    & \qquad
      - \algbrac{ \Weyl{e_j,\D_{e_i}X}{} }{ h(\ve^i,\ve^j) } \\
    &= 0
\end{aligned}
\end{equation*}
as $\Weyl{X,Y}{}$ acts trivially on $\cL^*$ by \thref{cor:ppg-bgg-Wtrace}.  The third term in \eqref{eq:ppg-bgg-metricprol-1} equals
\begin{align} \label{eq:ppg-bgg-metricprol-2}
  &\algbrac{ \ve^i }
           { \algbrac{ \Weyl{e_j,X}{} }
                     { \algbrac{ \algbrac{Z^{\D}}{e_i} }
                               { \ve^j } } } \notag \\
  &\quad
     =  \algbrac{ \ve^i }
                { \algbrac{ \Weyl{e_j,X}{} }
                          { \ve^j(Z^{\D}) \ltens e_i
                          + \algbrac{ Z^{\D} }{ \algbrac{e_i}{\ve^j}} } } \notag \\
  &\quad
     = -\ve^j(Z^{\D}) (\tr \Weyl{e_j,X}{})
      + \algbrac{ \ve^i }
                { \algbrac{ \Weyl{e_j,X}{} }
                          { \algbrac{ Z^{\D} }{ \algbrac{e_i}{\ve^j} } } } \notag \\
  &\quad
     =  \algbrac{ \ve^i }
                { \algbrac{ \algbrac{ \Weyl{e_j,X}{} }
                                    { Z^{\D} } }
                          { \algbrac{ e_i }{ \ve^j } } }
      + \algbrac{ \ve^i }
                { \algbrac{ Z^{\D} }
                          { \algbrac{ \algbrac{ \Weyl{e_j,X}{} }
                                              { e_i } }
                                    { \ve^j } } } \notag \\
  &\quad\!
    \begin{aligned}
      &=  \algbrac{ \algbrac{ \ve^i }
                            { \algbrac{ \Weyl{e_j,X}{} }
                                      { Z^{\D} } } }
                  { \algbrac{e_i}{\ve^j} }
        + \algbrac{ \algbrac{ \Weyl{e_j,X}{} }
                            { Z^{\D} } }
                  { \algbrac{ \ve^i }
                            { \algbrac{e_i}{\ve^j} } } \\
      &\qquad
        + \algbrac{ \ve^i }
                  { \algbrac{ Z^{\D} }
                            { \algbrac{ \algbrac{ \Weyl{e_j,X}{} }
                                                { e_i } }
                                      { \ve^j } } }  
    \end{aligned}
\end{align}
The first term here is $-\ve^j(e_i) \ve^j(\Weyl{e_j,X}{Z^{\D}} = -(\liebdy \Weyl{}{})_X(Z^{\D}) = 0$ by \thref{lem:ppg-alg-Lbrac}, while the second term equals $\tfrac{1}{2}r(n+1) \algbrac{ \algbrac{ \Weyl{e_j,X}{} }{ Z^{\D} } }{ \ve^j } = \tfrac{1}{2} (\liebdy \Weyl{}{})_X(Z^{\D}) = 0$ by \thref{cor:ppg-alg-trace}.  Applying the Bianchi identity \eqref{eq:ppg-bgg-Wbianchi} to the third term and using \eqref{eq:ppg-bgg-metricprol-2} yields
\begin{align*}
  &\algbrac{ \ve^i }
           { \algbrac{ \Weyl{e_j,X}{} }
                     { \algbrac{ \algbrac{Z^{\D}}{e_i} }
                               { \ve^j } } } \\
  &\quad \begin{aligned}
    &= -\algbrac{ \ve^i }
                { \algbrac{ Z^{\D} }
                          { \algbrac{ \algbrac{ \Weyl{X,e_i}{} }
                                              { e_j }
                                    + \algbrac{ \Weyl{e_i,e_j}{} }
                                              { X } }
                                    { \ve^j } } } \\
    &=  \algbrac{ \algbrac{ Z^{\D} }
                          { \algbrac{ \algbrac{ \Weyl{X,e_i}{} }
                                              { e_j } }
                                    { \ve^j } } }
                { \ve^i }
      + \algbrac{ \algbrac{ Z^{\D} }
                          { \algbrac{ \algbrac{ \Weyl{e_i,e_j}{} }
                                              { X } }
                                    { \ve^j } } }
                { \ve^i } \\
    &=  \algbrac{ \ve^i(Z^{\D}) }
                { \algbrac{ \algbrac{ \Weyl{X,e_i}{} }
                                    { e_j } }
                          { \ve^j } }
      + \algbrac{ Z^{\D} }
                { \algbrac{ \algbrac{ \algbrac{ \Weyl{e_i,e_j}{} }
                                              { X } }
                                    { \ve^j } }
                          { \ve^i } } \\
    &\qquad
      + \algbrac{ \ve^i(Z^{\D}) }
                { \algbrac{ \algbrac{ \Weyl{e_i,e_j}{} }
                                    { X } }
                          { \ve^j } }
      + \algbrac{ Z^{\D} }
                { \algbrac{ \algbrac{ \algbrac{ \Weyl{e_i,e_j}{} }
                                              { X } }
                                    { \ve^j } }
                          { \ve^i } }
  \end{aligned}
\end{align*}
The first term on the \rhs\ equals $\ve^i(Z^{\D}) \tr(\Weyl{X,e_i}{}) = 0$ by \thref{cor:ppg-bgg-Wtrace}.  The third term is $-\ve^i(Z^{\D}) (\liebdy \Weyl{}{})_{e_i}(X) = 0$.  The last term is symmetric under exchange of $\ve^i$ and $\ve^j$, but skew-symmetric under exchange of $e_i$ and $e_j$; since we are summing over $i,j$, this term must vanish.  Consequently
\vspace{-0.2em}
\begin{equation*}
  \algbrac{ \ve^i }
          { \algbrac{ \Weyl{e_j,X}{} }
                    { \algbrac{ \algbrac{Z^{\D}}{e_i} }
                              { \ve^j } } }
    = \algbrac{ Z^{\D} }
              { \algbrac{ \algbrac{ \algbrac{ \Weyl{X,e_i}{} }
                                            { e_j } }
                                  { \ve^j } }
                         { \ve^i } }.
\vspace{-0.2em}
\end{equation*}
The inner algebraic bracket is a section of $T^*M$; contracting with a vector field $Y$ yields
\begin{align*}
  \Killing{ \algbrac{ \algbrac{ \algbrac{ \Weyl{X,e_i}{} }
                                        { e_j } }
                              { \ve^j } }
                    { \ve^i } }
          { Y }
    &=  \Killing{ \Weyl{X,e_i}{} }
                { \algbrac{ e_j }
                          { \algbrac{ \algbrac{Y}{\ve^i} }
                                    { \ve^j } } } \\
    &=  \Killing{ \Weyl{X,e_i}{} }
                { \algbrac{ Y }
                          { \algbrac{ \algbrac{e_j}{\ve^i} }
                                    { \ve^j } }
                + \algbrac{ \algbrac{Y}{\ve^i} }
                          { \algbrac{e_j}{\ve^j} } } \displaybreak \\
    &=  \tfrac{1}{2}r(n+1)
        \killing{ \Weyl{X,e_i}{} }
                { -\algbrac{Y}{\ve^i} + \algbrac{ \algbrac{Y}{\ve^i} }{ \id } } \\
    &= -\tfrac{1}{2}r(n+1) \ve^i(\Weyl{X,e_i}{Y}) \\
    &=  \tfrac{1}{2}r(n+1) (\liebdy \Weyl{}{})_X(Y)
     = 0
\end{align*}
Therefore the third term in \eqref{eq:ppg-bgg-metricprol-2} also vanishes, leaving
\begin{align*}
  \liebdy( \D \liebdy \algbrac{\Weyl{}{}}{h} )_X
    &= -\algbrac{ \algbrac{ \ve^i }
                          { (\D_{e_i} \Weyl{}{})_{e_j,X} } }
                { h(\ve^j, \bdot) } \\
    &= -\ve^i( (\D_{e_i} \Weyl{}{})_{e_j,X} h(\ve^j,\bdot) ) \\
    &=  h( (\quab \CY{}{})_{e_i,X}, \ve^i )
\end{align*}
by the differential Bianchi identity of \itemref{prop:ppg-bgg-calc}{dbt}.  Reintroducing the eigenvalue $w \in \bR$ of $\quab^{-1}$ on $\liebdy \algbrac{W}{h}$ now gives the desired formula for the curvature correction. \end{proof}

The eigenvalue $w$ of $\quab^{-1}$ on $\liebdy \algbrac{\Weyl{}{}}{h}$ empirically equals $rn/2$, but unfortunately the author could not find a general argument for this.

\medskip
\subsection{Normalised Ricci curvature and \einstein\ metrics} 
\label{ss:ppg-bgg-nric}

Having described the metrisability of \ppgs, our next goal is to better understand the normalised Ricci tensor $\nRic[g]{}{}$ of a compatible metric $g$.  As observed with the classical \proj\ structures, there is a close relationship between $\nRic[g]{}{}$ and the Ricci curvature $\sRic[g]{}{}$ of the \riem\ curvature tensor, which we describe in \thref{prop:ppg-bgg-nric}.  This allows us to understand the role of \einstein\ metrics in the theory, generalising \thref{prop:qtn-bgg-normal} in \qtn ic\ geometry.  First however, we must develop some algebraic tools for dealing with sections of $\cB^*$.

\begin{lem} \thlabel{lem:ppg-bgg-p0B} Let $g \in \s{0}{\cB^*}$ be \non degenerate\ and suppose that $A$ is a $g$-\self adjoint\ section of $\alg{gl}{TM}$.  Then $\b{g} \defeq g(A\bdot,\bdot)$ is a section of $\cB^*$ if and only if $A$ is a section of the subbundle $\fp^0_M \leq \alg{gl}{TM}$. \end{lem}

\begin{proof} Suppose first that $\b{g} \defeq g(A\bdot,\bdot)$ is a section of $\cB^*$ and let $\pi \defeq (\det g)^{-1/r(n+1)}$.  Then since $g$ is \non degenerate, $\pi$ is nowhere-vanishing and so we may write
\begin{equation*}
  AX = g^{-1}(\b{g}(X,\bdot), \bdot) = \pi^{-1} \ltens g^{-1}( \pi \ltens \b{g}(X,\bdot), \bdot )
\end{equation*}
for all $X \in \s{0}{TM}$.  Using Table \ref{tbl:ppg-alg-Z2gr} this gives
\begin{align*}
  AX
    = - \pi^{-1} \ltens g^{-1}( \algbrac{ \pi \ltens \b{g} }{ X }, \bdot )
    = \algbrac{ \pi^{-1} \ltens g^{-1} }
               { \algbrac{ \pi \ltens \b{g} }{ X } }
    = \algbrac{ \algbrac{ \pi^{-1} \ltens g^{-1} }
                         { \pi \ltens \b{g} } }
               { X },
\end{align*}
where the last equality holds since $\pi^{-1} \ltens g^{-1} \in \s{0}{\cL^*\tens\cB}$ has zero bracket with $X \in \s{0}{TM}$.  Here $\algbrac{ \pi^{-1} \ltens g^{-1} }{ \pi \ltens \b{g} }$ is a section of $\fp^0_M \dsum \liecenter{\fq^0}_M$.  But since $\liecenter{\fq^0}_M$ acts trivially on $TM$, it has trivial intersection with $\alg{gl}{TM}$; thus the component of $\algbrac{\pi^{-1} \ltens g^{-1}}{\pi \ltens \b{g}}$ in $\liecenter{\fq^0}_M$ vanishes.

Conversely suppose that $A \in \s{0}{\fp^0_M}$.  Then since $A$ is \self adjoint\ \wrt\ $g$, we may write
\begin{align*}
  \pi \ltens g(AX,Y)
    &= \tfrac{1}{2}( \pi \ltens g(AX,Y) + \pi \ltens g(X,AY) ) \\
    &= -\tfrac{1}{2} \pi \ltens (A \acts g) (X,Y) \\
    &= -\tfrac{1}{2} \algbrac{A}{\pi \ltens g}(X,Y) + \tfrac{1}{2} \algbrac{A}{\pi} g(X,Y) \\
    &= -\tfrac{1}{2} \algbrac{A}{\pi \ltens g}(X,Y) + \tfrac{\tr A}{r(n+1)} \pi \ltens g(X,Y)
\end{align*}
by \thref{cor:ppg-alg-trace}.  Since the action of $\fp^0_M$ preserves the subbundle $\cL \tens \cB^* \leq \cW^*$, it follows that $\pi \ltens g(A\bdot,\bdot)$ is a section of $\cL \tens \cB^*$, and hence $g(A\bdot,\bdot)$ is a section of $\cB^*$ since $\pi$ is nowhere-vanishing. \end{proof}

\thref{lem:ppg-bgg-p0B} frames an intuitive notion from \cproj\ and \qtn ic\ geometry in Lie-theoretic terms: the composition of a hermitian inner product with a \self adjoint\ endomorphism is hermitian if and only if the endomorphism is complex- or \qtn-linear.

\begin{cor} \thlabel{cor:ppg-bgg-powA} Suppose that $A \in \s{0}{\fp^0_M}$ is \self adjoint\ \wrt\ $g$.  Then $A^k \in \s{0}{\fp^0_M}$ for all $k \in \bN$. \end{cor}

\begin{proof} Suppose that $A^k \in \s{0}{\fp^0_M}$ for some $k\in\bN$.  Then $g_k \defeq g(A^k\bdot, \bdot)$ is a section of $\cB^*$ by \thref{lem:ppg-bgg-p0B}, so that $g_{k+1} \defeq g_k(A\bdot, \bdot) = g(A^{k+1}\bdot, \bdot)$ is a section of $\cB^*$ as well.  Then $A^{k+1} \in \s{0}{\fp^0_M}$ by \thref{lem:ppg-bgg-p0B} again; now apply induction. \end{proof}

We next show that $\nRic[g]{}{}$ is a section of $\cB^*$ whenever $g$ is a compatible metric.  Note that this is tautological for conformal geometries: there the metric equation is the \einstein\ scale equation, and a solution necessarily has $\nRic[g]{}{}$ a multiple of the metric $g$ in the conformal class.  We extract the following technical lemma for later use.

\begin{lem} \thlabel{lem:ppg-bgg-weylp0} Let $h$ be a \non degenerate\ linear metric with corresponding metric $g \defeq (\det h)^{1/r} \ltens h^{-1}$.  Then $X \mapsto (\det h)^{1/r} \liebdy \algbrac{\Weyl{}{}}{h}_X$ defines a $g$-\self adjoint\ section of $\fp^0_M$. \end{lem}

\begin{proof} Denote the given endomorphism by $\Phi : X \mapsto (\det h)^{1/r} \ltens \liebdy \algbrac{\Weyl{}{}}{h}_X$.  As an initial sanity check, we note that $\liebdy \algbrac{\Weyl{}{}}{h}$ is an $(\cL^*\tens TM)$-valued $1$-form, so that $\Phi(X)$ is a section of $TM$ as required.

We first prove that $\Phi$ is \self adjoint\ \wrt\ $g$.  Since $\Weyl{}{}$ acts trivially on $\cL$,
\vspace{-0.1em}
\begin{equation*}
  g(\Phi(X), Y)
    = (\det h)^{1/r} \ltens g(\algbrac{ \ve^i }{ \Weyl{e_i,X}{h} }, Y)
    = -g(\Weyl{e_i,X}{\ve^{i\sharp}}, Y)
\vspace{-0.1em}
\end{equation*}
since $\Weyl{e_i,X}{\ve^i} = 0$ by virtue of $\liebdy \Weyl{}{} = 0$.  Then by the Leibniz rule,
\vspace{-0.1em}
\begin{equation*} \begin{aligned}
  g(\Weyl{e_i,X}{\ve^{i\sharp}}, Y)
    &= -(\Weyl{e_i,X}{g})(\ve^{i\sharp}, Y) - g(\ve^{i\sharp}, \Weyl{e_i,X}{Y}) \\
    &= -(\Weyl{e_i,X}{g})(\ve^{i\sharp}, Y) - (\liebdy \Weyl{}{})_{X}(Y).
\end{aligned}
\vspace{-0.1em}
\end{equation*}
By \thref{cor:ppg-bgg-lcconn}, the \LC\ connection $\D^g$ of $g$ lies in $\Dspace$.  Then since $\Curv[g]{}{g} = 0$ we have $\Weyl{}{g} = \algbracw{\id}{\nRic[g]{}{}}{} \acts g$ by \eqref{eq:ppg-bgg-curvdecomp}, giving
\begin{align} \label{eq:ppg-bgg-weylp0-1}
  (\Weyl{e_i,X}{g})(\ve^{i\sharp}, Y)
    &=  g( \algbrac{ \algbracw{\id}{\nRic[g]{}{}}{e_i,X} }
                   { \ve^{i\sharp} },
           Y )
      + g( \ve^{i\sharp},
           \algbrac{ \algbracw{\id}{\nRic[g]{}{}}{e_i,X} }
                   { Y } ) \notag \\
    &=  g( \algbrac{ \algbrac{e_i}{\nRic[g]{X}{}} }
                   { \ve^{i\sharp} },
           Y )
      - g( \algbrac{ \algbrac{X}{\nRic[g]{e_i}{}} }
                   { \ve^{i\sharp} },
           Y ) \notag \\
    &\qquad
      + \ve^i( \algbrac{ \algbrac{e_i}{\nRic[g]{X}{}} }
                       { Y } )
      - \ve^i( \algbrac{ \algbrac{X}{\nRic[g]{e_i}{}} }
                       { Y } ) \notag \\
    &\! \begin{aligned}
     &= \killing{ \algbrac{ \algbrac{e_i}{\nRic[g]{X}{}} }
                          { \ve^{i\sharp} }^{\flat} }
                { Y }
       - g( \algbrac{ \algbrac{X}{\nRic[g]{e_i}{}} }
                    { \ve^{i\sharp} },
            Y ) \\
     &\qquad
       + \tfrac{1}{2}r(n+1) \nRic[g]{X}{Y}
       - \ve^i( \algbrac{ \algbrac{X}{\nRic[g]{e_i}{}} }
                        { Y } ),
    \end{aligned}
\end{align}
where we have evaluated the third term using \thref{cor:ppg-alg-trace}.  Using \threfit{lem:app-alg-musical}, the first term above equals
\vspace{0.1em}
\begin{equation*}
\begin{aligned}
  \killing{ \algbrac{ \algbrac{e_i}{\nRic[g]{X}{}} }
                    { \ve^{i\sharp} }^{\flat} }
          { Y }
    &= \killing{ \algbrac{ \algbrac{ e_i^{\flat} }
                                   { {\nRic[g]{X}{}}^{\sharp} } }
                         { \ve^i } }
               { Y } \\
    &= \killing{ \algbrac{ \algbrac{ {\nRic[g]{X}{}}^{\sharp} }
                                   { e_i^{\flat} } }
                         { Y } }
               { \ve^i } \\
    &= (2-r) \nRic[g]{X}{Y}
\end{aligned}
\vspace{0.1em}
\end{equation*}
by \threfit{lem:app-alg-flat}\ref{lem:app-alg-flat-noA}, which is also symmetric in $X,Y$ by \thref{cor:ppg-bgg-exactweyl}.  For the second term on the \rhs\ of \eqref{eq:ppg-bgg-weylp0-1}, we note that
\begin{equation} \label{eq:ppg-bgg-weylp0-2}
  g( \algbrac{ \algbrac{X}{\nRic[g]{e_i}{}} }
             { \ve^{i\sharp} },
      Y )
    =  \killing{ \algbrac{ \algbrac{X}{\nRic[g]{e_i}{}} }
                         { \ve^{i\sharp} }^{\flat} }
               { Y }
    =  \ve^i( \algbrac{ \algbrac{Y}{X^{\flat}} }
                      { {\nRic[g\,\sharp]{e_i}{}} } )
\end{equation}
by using \threfit{lem:app-alg-musical} again.  Writing $\rho : X \mapsto {\nRic[g\,\sharp]{X}{}}$, the symmetry of $\nRic[g]{}{}$ implies that $\rho$ is \self adjoint\ \wrt\ $g$.  Alternating \eqref{eq:ppg-bgg-weylp0-2} in $X,Y$ then yields
\vspace{0.2em}
\begin{equation*}
  \ve^i( \algbrac{ \algbrac{Y}{X^{\flat}} }
                 { {\nRic[g\,\sharp]{e_i}{}} }
       - \algbrac{ \algbrac{X}{Y^{\flat}} }
                 { {\nRic[g\,\sharp]{e_i}{}} } )
    = -\tr( \algbracw{\id}{g}{X,Y} \circ \rho ).
\vspace{0.2em}
\end{equation*}
However, $\algbracw{\id}{g}{X,Y}$ is skew-adjoint \wrt\ $g$ by \threfit{cor:app-alg-skew}, while $\rho$ is \self adjoint; therefore the trace vanishes identically and the term $g(\algbrac{ \algbrac{X}{\nRic[g]{e_i}{}} }{ \ve^{i\sharp} }, Y)$ is symmetric in $X,Y$.  Since the fourth term in \eqref{eq:ppg-bgg-weylp0-1} is evidently symmetric in $X,Y$, we conclude that $g(\Weyl{e_i,X}{ \ve^{i\sharp} }, Y)$ is symmetric in $X,Y$, \ie\ that $\Phi$ is $g$-\self adjoint.

The proof that $\Phi$ defines a section of $\fp^0_M$ unfortunately requires some case-by-case analysis using the classification of Subsection \ref{ss:ppg-class-real}.  For $r=1$ we have $\fp^0_M = \alg{gl}{TM}$, so there is nothing to prove.  For $r=8$ we have one of two \oct ic\ geometries, for which the Weyl curvature vanishes identically (see Figure \ref{fig:ppg-bgg-hasseE}); thus there is nothing to prove since $\Phi$ is identically zero.  For $r=n$ we have various conformal geometries, and a linear metric defines an \einstein\ metric $g$.  Then $\nRic[g]{}{}$ is proportional to $g$, giving $\algbracw{\id}{\nRic[g]{}{}}{} \acts g = 0$ and hence $\algbrac{\Weyl{}{}}{h} = 0$; thus $\Phi = 0$ again.

For the remaining geometries, we are going construct a subbundle $\cZ \leq \fp^0_M$ in which the Weyl curvature takes values, and which may be identified with the $\alg{gl}{TM}$-centraliser of a finite-dimensional subbundle $\cA \leq \fp^0_M$ spanned by endomorphisms $\{ A_i \}_i$ that satisfy $A_i^2 = \pm \id$, all with the same sign, and $g(A_i\bdot, A_i\bdot) = g(\bdot, \bdot)$.  Supposing that we have done this, for each $A \in \s{0}{\cA}$ we may decompose $\Weyl{}{}$ into $\pm A$-hermitian parts as
\vspace{0.3em}
\begin{equation*}
  \Weyl{X,Y}{} = \underbrace{ \tfrac{1}{2}(\Weyl{X,Y}{} + \Weyl{AX,AY}{}) }
                           _{ \eqdef \, \Weyl[A+]{X,Y}{} }
               + \underbrace{ \tfrac{1}{2}(\Weyl{X,Y}{} - \Weyl{AX,AY}{}) }
                           _{ \eqdef \, \Weyl[A-]{X,Y}{} },
\vspace{0.3em}
\end{equation*}
where $\Weyl[A\pm]{AX,AY}{} = \pm \Weyl[A\pm]{X,Y}{}$.  Henceforth we write ``$\pm$'' to mean the sign of $A^2 = \pm \id$.  Then since $A$ commutes with $\Weyl{X,Y}{}$, the Bianchi identity \eqref{eq:ppg-bgg-Wbianchi} gives
\begin{align*}
  \Weyl[A\pm]{AX,Y}{Z}
    &=  \tfrac{1}{2}( \Weyl{AX,Y}{Z} + \Weyl{X,AY}{Z} ) \\
    &= -\tfrac{1}{2}( \Weyl{Y,Z}{AX} + \Weyl{Z,AX}{Y} + \Weyl{AY,Z}{X} + \Weyl{Z,X}{AX} ) \\
    &= -\tfrac{1}{2}A( \Weyl{Y,Z}{X} \pm \Weyl{Z,AX}{AY} \pm \Weyl{AY,Z}{AX} + \Weyl{Z,X}{Y} ) \\
    &=  \tfrac{1}{2}A( \Weyl{X,Y}{Z} \pm \Weyl{AX,AY}{Z} ) \\
    &= A \Weyl[A\pm]{X,Y}{Z},
\end{align*}
\ie\ $\Weyl[A\pm]{AX,Y}{} = A \circ \Weyl[A\pm]{X,Y}{}$.  Since $A$ is invertible, $\{Ae_i\}_i$ is a local frame of $TM$ with dual coframe $\{ \mp A\ve^i \}$; indeed, $\mp A \ve^i(Ae_j) = \pm \ve^i (A^2 e_j) = +\delta_{ij}$ as required.  Moreover, since $g$ is $A$-hermitian, we have $AX^{\flat} = \mp (AX)^{\flat}$.  Then
\vspace{-0.1em}
\begin{equation*} \begin{aligned}
  \Weyl[A\pm]{e_i,X}{ \ve^{i\sharp} }
    &= \Weyl[A\pm]{Ae_i,X}{ A\ve^{i\sharp} } \\
    &= \mp A^2 \Weyl[A\pm]{e_i,X}{ \ve^{i\sharp} }
     = -\Weyl[A\pm]{e_i,X}{ \ve^{i\sharp} },
\end{aligned}
\vspace{-0.1em}
\end{equation*}
which consequently vanishes.  Thus the only component of $\Weyl{}{}$ which contributes to $\Weyl{e_i,X}{ \ve^{i\sharp} }$ is the component $\Weyl[\mp]{}{}$ which is $\mp A$-hermitian \wrt\ all $A$ in our chosen basis of $\cA$, giving $\Weyl{e_i,AX}{ \ve^{i\sharp} } = A \Weyl{e_i,X}{ \ve^{i\sharp} }$.  By construction of $\cA$, this implies that $\Phi : X \mapsto \Weyl{e_i,X}{\ve^{i\sharp}}$ takes values in the subbundle $\cZ \leq \fp^0_M$ in which $\Weyl{}{}$ takes values; in particular, $\Phi$ defines a section of $\fp^0_M$.%
\newcommand{\liectr}[1]{\liecenter[\alg{gl}{}]{#1}}  
It remains to construct the subbundles $\cZ$ and $\cA$, which we do on a case-by-case basis:
\begin{slimitemize}
  \item For the real form $\fg = \alg{sl}{n+1,\bC}$ in type \type{A}{2n+1}, we have $\fp^0 = \alg{gl}{n,\bC}$ and are dealing with \cproj\ geometry.  We may identify $\fp^0_M$ with the $\alg{gl}{TM}$-centraliser of the complex structure $J$, which preserves $g$ and satisfies $J^2 = -\id$.  Thus we may take $\cZ = \fp^0_M$ and $\cA = \linspan{J}{} \leq \fp^0_M$ in this case.
  
  \item For the split real form $\fg = \alg{sl}{n+1,\bR} \dsum \alg{sl}{n+1,\bR}$ in type \type{A}{2n+1} we have $\fp^0 = \alg{gl}{n,\bR} \dsum \alg{gl}{n,\bR}$, embedded into $\alg{gl}{2n,\bR}$ as block-diagonal matrices $\inlinematrix{ A & 0 \\ 0 & B }$.  A generic element of $\alg{gl}{2n,\bR}$ lies in $\fp^0$ if and only if it commutes with $I = \inlinematrix{ \id & 0 \\ 0 & -\id }$, which satisfies $I^2 = \id$, so that $\fp^0 = \liectr{I}$.  Moreover from Table \ref{tbl:app-tbl-Wd} we have
  \vspace{-0.1em}
  \begin{equation*}
    B^* = \dynkinAAp[0,0,0,1,-2, 0,0,0,1,-2]{2}{0}{3}{ppg-bgg-normal-BdAs}
        = \bR^{n*} \etens \conj{\bR^{n*}},
  \vspace{-0.1em}
  \end{equation*}
  where conjugate denotes representations of the second factor.  Since $I$ acts by $\id$ on $\bR^{n*}$ and by $-\id$ on $\conj{\bR^{n*}}$, we conclude that elements of $B^*$ are $I$-hermitian.  Thus we may take $\cZ = \fp^0_M$ again, and $\cA$ to be the span of the endomorphism of $TM$ induced pointwise by $I$.
  
  \item For the real form $\fg = \alg{sl}{n+1,\bH}$ in type \type{D}{2n+2}, we have $\fp^0 = \alg{gl}{n,\bH} \dsum \alg{sp}{1}$ and are dealing with \qtn ic\ geometry.  We may take $\cZ$ to be the \qtn-linear endomorphisms $\alg{gl}{TM,\cQ}$ of $TM$, and $\cA$ to be the \qtn ic\ bundle $\cQ$.  The Weyl curvature is \qtn-linear by \itemref{prop:qtn-para-calc}{Wprop}.
  
  \item The split real form $\fg = \alg{sl}{2n+2,\bR}$ in type \type{D}{2n+2} is the only tricky case.  We have $\fp^0 = \alg{s}{ \alg{gl}{2n,\bR} \dsum \alg{gl}{2,\bR} }$, giving \grassmannian\ geometry as studied in \cite[\S 4.1.3]{cs2009-parabolic1} and \cite{gs1999-qtntwistor}.  We embed $\fp^0$ into $\alg{gl}{4n,\bR}$ using the outer matrix product
  \begin{equation*}
    \left( A, \inlinematrix{ a & b \\ c & d } \right) \mapsto
      \begin{pmatrix}
        aA & bA \\
        cA & dA
      \end{pmatrix}
      \colvectpunct[-0.7em]{.}
  \end{equation*}
  The tangent bundle has highest weight
  \begin{equation*}
    TM = \dynkinApp[1,0,0,0,0,1]{3}{0}{3}{ppg-bgg-normal-gpDs}
      = \cE \etens \cH,
  \end{equation*}
  where $\cE, \cH$ are associated to the representations $E \defeq \bR^{2n}$ and $H \defeq \bR^2$, giving
  \begin{equation*}
    \alg{gl}{TM}
      = (\cE \tens \cE^*) \etens (\cH \tens \cH^*)
      = \alg{gl}{\cE} \etens \alg{sl}{\cH} \dsum \alg{gl}{\cE} \etens \bR;
  \end{equation*}
  \cf\ \eqref{eq:qtn-para-glEH}.  The summand $\alg{gl}{\cE} \etens \bR$ may be viewed as the space of all endomorphisms commuting with $\alg{sl}{\cH} \isom \alg{sl}{2,\bR}$, which is $3$-dimensional with basis
  \begin{equation} \label{eq:ppg-bgg-weylp0-grassbasis}
    \begin{bmatrix} \id & 0 \\
                    0   & -\id \end{bmatrix}, \quad
    \begin{bmatrix} 0   & \id \\
                    \id & 0    \end{bmatrix}
      \quad\text{and}\quad
    \begin{bmatrix} \tfrac{1}{\sqrt{2}}\id & \id \\
                    \tfrac{1}{2}\id        & -\tfrac{1}{\sqrt{2}}\id \end{bmatrix}
      \colvectpunct[-0.9em]{\, .}
  \end{equation}
  These all satisfy $A^2 = \id$.  Moreover
  \begin{equation*}
    \cB^* = \dynkinApp[0,0, 0,1,0,-2,0]{2}{0}{5}{ppg-bgg-normal-BdDs}
          = \Wedge{2}\cE^* \etens \Wedge{2}\cH^*
  \end{equation*}
  and $\alg{sl}{\cH}$ acts trivially on $\Wedge{2}\cH^* \isom \bR$, so that $\alg{sl}{\cH}$ preserves elements of $\cB^*$.  Thus we take $\cA = \alg{sl}{\cH}$ with the basis \eqref{eq:ppg-bgg-weylp0-grassbasis}.  By \cite[p.\ 377]{cs2009-parabolic1}, the Weyl curvature takes values in $\alg{sl}{\cE} \leq \cZ \defeq \alg{gl}{\cE} \etens \bR$.
\end{slimitemize}
This exhausts all the remaining real forms, so we are done. \end{proof}

\begin{thm} \thlabel{thm:ppg-bgg-nRicB} Let $h$ be a \non degenerate\ solution of the linear metric equation with \LC\ connection $\D \in \Dspace$.  Then $\nRic[g]{}{}$ is a section of $\cB^*$. \end{thm}

\begin{proof} Since $g \defeq (\det h)^{-1/r} \ltens h^{-1}$ is a section of $\cB^*$, by \thref{lem:ppg-bgg-p0B} it suffices to show that $X \mapsto (\nRic[g]{X}{})^{\sharp}$ is a section of $\fp^0_M$.  Evaluating the metric prolongation connection \eqref{eq:ppg-bgg-metricprol} \wrt\ $\D = \D^g$, we find that $Z^{\D} = 0$ and hence
\vspace{-0.15em}
\begin{equation*}
  g^{-1}(\nRic[g]{}{}, \bdot) = -\pi \ltens \lambda^g \ltens \id
    - \pi \ltens \quab^{-1} \liebdy \algbrac{\Weyl{}{}}{h}.
\vspace{-0.15em}
\end{equation*}
The first term on the \rhs\ clearly lies in $\fp^0_M$, and $X \mapsto \pi \ltens \liebdy \algbrac{\Weyl{}{}}{h}_X$ is a section of $\fp^0_M$ by \thref{lem:ppg-bgg-weylp0}.  Since $\fp^0_M$ acts on $\cW$ preserving the grading, $\algbrac{\Weyl{}{}}{h}$ lies in the irreducible subbundle $\cL^*\tens\cB$, which implies that $\liebdy \algbrac{\Weyl{}{}}{h}$ also lies in an irreducible subbundle of $T^*M \tens \cW$.  Therefore $\quab^{-1}$ acts by a scalar, so we see that $X \mapsto g^{-1}(\nRic[g]{X}{}, \bdot)$ is a section of $\fp^0_M$ as required. \end{proof}

For \proj\ differential geometry, $\cB^* = \Symm{2}T^*M$ and this statement is vacuous.  In \cproj\ and \qtn ic\ geometries, \thref{thm:ppg-bgg-nRicB} amounts to saying that $\nRic[g]{}{}$ is $J$- or Q-hermitian.  As noted above, this result is tautological for conformal geometries.

With the technical work completed, it is straightforward to relate $\nRic[g]{}{}$ and $\sRic[g]{}$.

\begin{prop} \thlabel{prop:ppg-bgg-nric} Let $g$ be a \non degenerate\ compatible metric with \LC\ connection $\D^g \in \Dspace$.  Then
\vspace{-0.15em}
\begin{equation} \label{eq:ppg-bgg-ricnric}
  \sRic[g]{} = -\tfrac{1}{2}(rn+3r-4) \nRic[g]{}{}.
\vspace{-0.15em}
\end{equation}
\end{prop}

\begin{proof} Since $g$ is \non degenerate and $\nRic[g]{}{}$ is symmetric, we may write $\nRic[g]{}{} = g(\rho\bdot,\bdot)$ for some endomorphism $\rho \in \s{0}{\alg{gl}{TM}}$ which is evidently \self adjoint\ \wrt\ $g$.  Then $\rho \in \s{0}{\fp^0_M}$ by \thref{lem:ppg-bgg-p0B,thm:ppg-bgg-nRicB}.  We have
\vspace{-0.15em}
\begin{equation*}
  (\rho X)^{\flat}(Y) = g(\rho X,Y) = g(X,\rho Y) = X^{\flat}(\rho Y) = -\algbrac{\rho}{X^{\flat}}(Y)
\vspace{-0.15em}
\end{equation*}
by \self adjointness, so the Ricci curvature $\sRic[g]{}$ is given by
\begin{align*}
  \sRic[g]{X,Y} \defeq \ve^i( \Curv[g]{e_i,X}{Y} ) &= \ve^i( \Weyl{e_i,X}{Y} ) - \ve^i( \algbrac{ \algbracw{\id}{\nRic[g]{}{}}{e_i,X} }{ Y }) \\
  &= -\ve^i( \algbrac{ \algbrac{e_i}{(\rho X)^{\flat}} }{ Y } - \algbrac{ \algbrac{X}{(\rho e_i)^{\flat}} }{ Y } )
\end{align*}
since $\liebdy\Weyl{}{} = 0$.  By \thref{cor:ppg-alg-trace}, the first term on the \rhs\ evaluates to
\begin{equation*}
  -\ve^i( \algbrac{ \algbrac{Y}{(\rho X)^{\flat}} }{ e_i } ) = -\tfrac{1}{2}r(n+1) (\rho X)^{\flat}(Y) = -\tfrac{1}{2}r(n+1) \nRic[g]{X}{Y}.
\end{equation*}
For the second term, the Jacobi identity yields
\begin{align*}
  \ve^i( \algbrac{ \algbrac{X}{(\rho e_i)^{\flat}} }{ Y } )
  &= -\ve^i( \algbrac{ \algbrac{X}{\algbrac{\rho}{e_i^{\flat}}} }{ Y } \\
  &= -\ve^i( \algbrac{ \algbrac{ \algbrac{X}{\rho} }{e_i^{\flat} } }{ Y } + \algbrac{ \algbrac{ \rho }{ \algbrac{X}{e_i^{\flat}} } }{ Y } ) \\
  &= -\ve^i( -\algbrac{ \algbrac{\rho X}{e_i^{\flat}} }{ Y } + \algbrac{ \algbrac{\rho}{Y} }{ \algbrac{X}{e_i^{\flat}} } + \algbrac{ \rho }{ \algbrac{ \algbrac{X}{e_i^{\flat}} }{ Y } } ) \\
  &= \ve^i( \algbrac{ \algbrac{\rho X}{e_i^{\flat}} }{ Y } ) + \ve^i( \algbrac{ \algbrac{X}{e_i^{\flat}} }{ \rho Y } - \ve^i( \rho \algbrac{ \algbrac{X}{e_i^{\flat}} }{ Y } ) \\
  &= (2-r)g(\rho X,Y) + (2-r)g(X,\rho Y) - (2-r)g(\rho X,Y) \\
  &= (2-r) \nRic[g]{X}{Y}
\end{align*}
by \threfit{lem:app-alg-flat}.  Therefore
\begin{align*}
  \sRic[g]{X,Y}
    &= -\tfrac{1}{2}r(n+1) \nRic[g]{X}{Y} + (2-r)\nRic[g]{X}{Y} \\
    &= -\tfrac{1}{2}(rn+3r-4) \nRic[g]{X}{Y}
\end{align*}
as claimed. \end{proof}

Note that \eqref{eq:ppg-bgg-ricnric} is much simpler than the formula relating $\sRic[\D]{}$ and $\nRic{}{}$ for an arbitrary Weyl connection $\D \in \Dspace$, which will contain, for example, skew-symmetric components.  Indeed, generally we only have that $\nRic{}{}$ is a section of $T^*M \tens T^*M$, which is not irreducible.
The reducible nature of $T^*M \tens T^*M$ explains why we see different scalings on the irreducible components of $\nRic{}{}$ in the classical cases.  Note also that $\nRic[g]{}{} \defeq -\quab_M^{-1} \sRic[g]{}$, so \thref{prop:ppg-bgg-nric} calculates the action of $\quab_M^{-1}$ on sections of $\cB^* \leq T^*M \tens T^*M$.

The classification of \thref{thm:ppg-class-cpxclass} implies that the constant $\tfrac{1}{2}(rn+3r-4)$ is always positive for $n>1$, which equals $\tfrac{1}{2}(n-1)$ when $r=1$; equals $n+1$ when $r=2$; and equals $2n+4$ when $r=4$.  Thus \thref{prop:ppg-bgg-nric} is consistent with Propositions \ref{prop:proj-para-nric}, \ref{prop:cproj-para-calc}\ref{prop:cproj-para-calc-nric} and \ref{prop:qtn-para-calc}\ref{prop:qtn-para-calc-nric}.  It is interesting to note that $\tfrac{1}{2}(rn+3r-4) = 2(r-1)$ when $n=1$, implying in particular that $\nRic[g]{}{}$ cannot be recovered from $\sRic[g]{} = 0$ for $1$-dimensional \proj\ structures.

Our description of normalised Ricci curvature allows us to understand the role of \einstein\ metrics in the theory: they are normal solutions of the linear metric equation.

\begin{defn} A solution $\sigma \in \s{0}{\liehom{0}{\cV}}$ of a first BGG equation is \emph{normal} if it is of the form $\sigma = \bggproj{\bV}(s)$ for some $\D^{\bV}$-parallel section $s \in \s{0}{\cV}$. \end{defn}

Equivalently, a solution $\sigma$ is normal if the curvature corrections of the prolongation connection act trivially on $\bggrepr{\bV}(\sigma)$.

\begin{thm} \thlabel{thm:ppg-bgg-normal} A \non degenerate\ linear metric $h \in \s{0}{\cL^*\tens\cB}$ is a normal solution if and only if the corresponding metric $g \defeq (\det h)^{-1/r} \ltens h^{-1}$ is \einstein. \end{thm}

\begin{proof} The proof is similar to the \qtn ic\ case of \thref{prop:qtn-bgg-normal}.  First suppose that $(h,Z^{\D},\lambda^{\D})$ is $\D^{\bW}$-parallel.  Then writing $\D^{\bW}$ \wrt\ the \LC\ connection $\D^g$ of $g \defeq (\det h)^{1/r} \ltens g^{-1}$, we have $Z^g = \tfrac{2}{rn-r+2} \liebdy (\D^g h) = 0$ and hence
\begin{equation*}
  \D^{\bW}_X \colvect{ h \\ Z^g \\ \lambda^g }
    = \colvect{ 0 \\
                -\lambda^g \ltens X - h(\nRic{X}{},\bdot) \\
                \D_X \lambda^g }
    = 0.
\end{equation*}
Therefore $\lambda^{\D}$ is a constant multiple of the global trivialisation $(\det h)^{-1/r}$ of $\cL^*$, say $\lambda^{\D} = c(\det h)^{-1/r}$.  Since $\sRic[g]{} = -\tfrac{1}{2} (rn+3r -4) \nRic[g]{}{}$ by \thref{prop:ppg-bgg-nric}, we obtain
\begin{equation*}
  \sRic[g]{X,\bdot} = \tfrac{1}{2}(rn+3r-4) h^{-1}(\lambda^{\D} \ltens X, \bdot)
    = \tfrac{1}{2}c(rn+3r-4) g(X,\bdot)
\end{equation*}
for all vector fields $X$, \ie\ that $g$ is \einstein.

Conversely, suppose that $g$ is an \einstein\ metric.  Then $\nRic[g]{}{} = cg$ for some constant $g$ by \thref{prop:ppg-bgg-nric}, a section of $\cB^*$.  We must show that the curvature corrections in \eqref{eq:ppg-bgg-metricprol} vanish identically.  Calculating $\D^{\cW}$ \wrt\ the \LC\ connection $\D^g$ of $g$, the Cotton--York tensor $\CY[g]{}{} \defeq \d^{\D^g} \nRic[g]{}{} = c (\d^{\D^g} g) = 0$ clearly vanishes, so it remains to show that $\liebdy \algbrac{\Weyl{}{}}{h} = 0$.  Since $\Weyl{}{}$ acts trivially on $\cL$, it suffices to show that $\Weyl{e_i,X}{\ve^{i\sharp}} = 0$ for all $X\in\s{0}{TM}$, for which
\begin{equation} \label{eq:ppg-bgg-normal-1}
  \Weyl{e_i,X}{\ve^{i\sharp}}
    = \Curv[g]{e_i,X}{\ve^{i\sharp}} + c \algbracw{\id}{g}{e_i,X} \acts \ve^{i\sharp}.
\end{equation}
For the first term on the \rhs\ of \eqref{eq:ppg-bgg-normal-1},
\begin{align*}
  g(\Curv[g]{e_i,X}{\ve^{i\sharp}},Y)
    &= -\ve^i(\Curv[g]{e_i,X}{Y}) \\
    &= -\sRic[g]{X,Y}
     = \tfrac{1}{2}(rn+3r-4) \nRic[g]{X}{Y}
\end{align*}
by \thref{prop:ppg-bgg-nric}, so that $\Curv[g]{e_i,X}{\ve^{i\sharp}} = \tfrac{1}{2}(rn+3r-4) cX$.  For the second term on the \rhs\ of \eqref{eq:ppg-bgg-normal-1},
\begin{align*}
  g( \algbracw{\id}{g}{e_i,X} \acts \ve^{i\sharp}, Y )
    &= -g( g^{-1}(\ve^i,\bdot), \algbracw{\id}{g}{e_i,X} \acts Y ) \\
    &= -\ve^i( \algbrac{ \algbrac{e_i}{X^{\flat}} }{ Y }
       - \algbrac{ \algbrac{X}{e_i^{\flat}} }{ Y } ) \\
    &= -\tfrac{1}{2}r(n+1) g(X,Y) + (2-r)g(X,Y) \\
    &= -\tfrac{1}{2}(rn+3r-4) g(X,Y)
\end{align*}
by \threfit{cor:app-alg-skew} and \threfit{lem:app-alg-flat}.  Therefore
\begin{equation*}
  \Weyl{e_i,X}{\ve^{i\sharp}}
    = \tfrac{1}{2}(rn+3r-4) \nRic[g]{X}{Y} - \tfrac{1}{2}(rn+3r-4) c g(X,Y) = 0
\end{equation*}
as required. \end{proof}

\subsection{The hessian equation} 
\label{ss:ppg-bgg-hess}

For \cproj\ and \qtn ic\ geometries, the first BGG operator associated to the dual $\fg$-representation $\bW^*$ was a second order hessian equation.  For the most part, the same is true in general.  As for the metric equation, we start by obtaining an explicit formula for the first BGG operator.  We continue to fix a Weyl structure and assume that $\Tor{} = 0$.

\begin{lem} \thlabel{lem:ppg-bgg-hesstractor} The tractor connection $\D^{\bW^*}$ and its curvature $\Curv[\bW^*]{}{}$ may be written
\vspace{0.2em}
\begin{equation*} \begin{aligned}
  \D^{\bW^*}_X \! \colvect{ \ell \\ \eta \\ \theta }
  &= \colvect{ \D_X \ell + \eta(X) \\
               \D_X \eta + \theta(X,\bdot) + \ell \ltens \nRic{X}{} \\
               \D_X \theta + \algbrac{\nRic{X}{}}{\eta} } \\[0.5em]
  \Curv[\bW^*]{X,Y}{ \colvect{ \ell \\ \eta \\ \theta } }
  &= \colvect{ 0 \\
               \Weyl{X,Y}{\eta} + \CY{X,Y}{} \ltens \ell \\
               \Weyl{X,Y}{\theta} + \algbrac{\CY{X,Y}{}}{\eta} }
\end{aligned}
\vspace{0.2em}
\end{equation*}
\wrt\ any Weyl structure. \end{lem}

\begin{proof} This follows from the general \formulae\ and the algebraic brackets. \end{proof}

Note the different sign conventions in \thref{lem:ppg-bgg-hesstractor} compared to the \cproj\ hessian \eqref{eq:cproj-bgg-hessprol} and \qtn ic\ hessian \eqref{eq:qtn-bgg-hessprol}, which we intentionally chose to avoid having minus signs throughout those derivations.  The signs above are those dictated by our conventions from Chapter \ref{c:para}.

\begin{lem} \thlabel{lem:ppg-bgg-Wdquabla} $\quab$ acts trivially on $\cL$, by the identity on $\cL\tens T^*M$, and by multiplication with $r$ on $\cL\tens\cB^*$. \end{lem}

\begin{proof} The $\ell$-slot is immediate, since $\cL \isom \liehom{0}{\cW^*}$.  For the $\eta$-slot, we have
\vspace{0.25em}
\begin{equation*}
  \quab \eta = \liebdy \liediff \eta
    = \algbrac{ \ve^i }{ \algbrac{e_i}{\eta} }
    = \algbrac{ \ve^i }{ \eta(e_i) }
    = \eta(e_i) \ltens \ve^i
    = \eta.
\vspace{0.25em}
\end{equation*}
For the $\theta$-slot, we have $\quab \theta = \sum{i}{} \, \algbrac{ \ve^i }{ \theta(e_i,\bdot) }$ and hence
\vspace{0.25em}
\begin{equation*} \begin{aligned}
  (\quab \theta)(X,Y)
  &= \algbrac{ \algbrac{ \algbrac{\ve^i}{\theta(e_i,\bdot)} }{ X } }{ Y } \\
  &= -\algbrac{ \algbrac{ \algbrac{X}{\ve^i} }{ Y } }
              { \theta(e_i,\bdot) }
    + \algbrac{ \algbrac{X}{\ve^i} }{ \theta(e_i,Y) }
    + \algbrac{ \algbrac{Y}{\ve^i} }{ \theta(e_i,X) } \\
  &= - \theta(e_i, \algbrac{ \algbrac{X}{\ve^i} }{ Y }) + 2\theta(X,Y)
\end{aligned}
\vspace{0.25em}
\end{equation*}
by the Jacobi identity.  Suppose first that $\theta$ is \non degenerate\ as an $\cL$-valued bilinear form on $TM$.  Then we may write $\theta = \ell \ltens g$ for a \non degenerate\ section $g$ of $\cB^*$, for which the first term equals
\begin{align*}
  -\ell \ltens g(e_i, \algbrac{ \algbrac{X}{\ve^i} }{ Y })
  &= -\ell \ltens \killing{ e_i^{\flat} }{ \algbrac{ \algbrac{X}{\ve^i} }{ Y } } \\
  &= -\ell \ltens \ve^i( \algbrac{ \algbrac{X}{e_i^{\flat}} }{ Y } ) \\
  &= (r-2) \theta(X,Y)
\end{align*}
by \threfit{lem:app-alg-flat}.  Combining this with the previous expression, we obtain $\quab \theta = r \theta$ in the case that $\theta$ is \non degenerate.  Since the \non degenerate\ bilinear forms constitute a dense open subset of $\cL\tens\cB^*$, the general result follows by continuity. \end{proof}

\begin{prop} The BGG splitting operator associated to $\bW^*$ is given by
\vspace{0.2em}
\begin{equation} \label{eq:ppg-bgg-hesssplit}
  \bggrepr{\bW^*} : \ell \mapsto
    \colvect{ \ell \\
              -\D\ell \\
              \tfrac{1}{r} \algbrac{ \ve^i }{ \D_{e_i}\D\ell - \ell \ltens \nRic{e_i}{} } }
    \colvectpunct{.}
\vspace{0.2em}
\end{equation}
\wrt\ any Weyl structure. \end{prop}

\begin{proof} We calculate the operator $\bggpi{\bW^*} : \s{0}{\cW^*} \to \s{1}{\cW^*}$, which is given by $\bggpi{\bW^*}(s) \defeq s - \quab_M^{-1}\liebdy (\D^{\bW^*}s)$ since $\liebdy s = 0$ for all $s \in \s{0}{\cW^*}$.  We have
\begin{equation*}
  \liebdy \big( \D^{\bW^*} \mr{repr}(\ell) \big)
    = \ve^i \acts \colvect{ \D_{e_i} \ell \\
                            \ell \ltens \nRic{e_i}{} \\
                            0 }
    = \colvect{ 0 \\
                (\D_{e_i}\ell) \ltens \ve^i \\
                \algbrac{ \ve^i }{ \ell \ltens \nRic{e_i}{} } }
        \colvectpunct{,}
\end{equation*}
where we note that $(\D_{e_i}\ell)\ve^i = \D\ell$.  To compute the action of $\quab_M^{-1}$, we note that
\begin{align*}
  &(\quab_M - \quab) \quab^{-1} \liebdy \big( \D^{\bW^*} \mr{repr}(\ell) \big) \\
  &\quad
    = \ve^j \acts
      \D_{e_j} \colvect{ 0 \\
                         \quab^{-1} \D\ell \\
                         \quab^{-1} \algbrac{ \ve^i }{ \ell \ltens \nRic{e_i}{} } }
    + \ve^i \acts \nRic{e_j}{} \acts
      \colvect{ 0 \\
                \quab^{-1} \D\ell \\
                \quab^{-1} \algbrac{ \ve^i }{ \ell \ltens \nRic{e_i}{} } }
    = \colvect{ 0 \\
                0 \\
                \algbrac{ \ve^i }{ \D_{e_i} \D\ell } }
    \colvectpunct{.}
\end{align*}
The Neumann series \eqref{eq:para-bgg-neumann} for $\quab_M^{-1}$ then gives
\vspace{-0.2em}
\begin{equation*} \begin{aligned}
  \quab_M^{-1}\liebdy \big( \D^{\bW^*} \mr{repr}(\ell) \big)
    &= \big( \id - \quab^{-1}(\quab_M - \quab) \big) \quab^{-1}
         \liebdy \big( \D^{\bW^*} \mr{repr}(\ell) \big) \\
  &= \colvect{ 0 \\
               \D\ell \\
               \tfrac{1}{r} \algbrac{ \ve^i }{ \ell \ltens \nRic{e_i}{} - \D_{e_i} \D \ell } }
\end{aligned}
\vspace{-0.2em}
\end{equation*}
by \thref{lem:ppg-bgg-Wdquabla}.  The claimed form of $\bggrepr{\bW^*}$ now follows. \end{proof}

It follows easily that the first BGG operator $\bgg{\bW^*}$ is given by the projection to the first homology of
\vspace{-0.2em}
\begin{equation} \label{eq:ppg-bgg-hessproj}
  \mr{proj}
    \colvect{ 0 \\
              -(\D^2 \ell - \ell \ltens \nRic{}{})
                + \tfrac{1}{r} \algbrac{ \ve^i }{ \D_{e_i}\D\ell - \ell\ltens\nRic{e_i}{} } \\
              \tfrac{1}{r} \D\algbrac{ \ve^i }{ \D_{e_i}\D\ell - \ell \ltens \nRic{e_i}{} } }
    \colvectpunct{.}
\vspace{-0.2em}
\end{equation}
Using \thref{cor:ppg-bgg-Wtrace} and the assumption $\Tor{} = 0$, it is easy to see that $\D^2 \ell - \ell \ltens \nRic[\D]{}{}$ defines a section of $\cL \tens \Symm{2}T^*M$.  The term $\tfrac{1}{r} \algbrac{ \ve^i }{ \D_{e_i} \D\ell - \ell\ltens\nRic{e_i}{} }$ should then be viewed as the projection to $\cL\tens\cB^*$ of $\D^2\ell - \ell\ltens\nRic{}{}$; this is made more precise by \thref{prop:ppg-bgg-hesseqn}.  If the second slot is \non zero, we then have
\begin{equation} \label{eq:ppg-bgg-hesseqn}
  \bgg{\bW^*}(\ell) = -(\D^2\ell - \ell\ltens\nRic{}{}) + \tfrac{1}{r} \algbrac{ \ve^i }{ \D_{e_i}\D\ell - \ell\ltens\nRic{e_i}{} },
\end{equation}
the projection away from $\cL \tens \cB^*$ inside $\cL \tens \Symm{2}T^*M$.   For $\fh$ simple and not of type \type{C}{n+1}, observe that $\cB^*$ is complementary to the Cartan square $\Cartan{2}T^*M$ inside $\Symm{2}T^*M$.

For type \type{C}{n+1}, we have $\cB^* = \Symm{2}T^*M$ and one can check that \eqref{eq:ppg-bgg-hesseqn} is trivial.  The first BGG equation is given by projection onto the third slot of \eqref{eq:ppg-bgg-hessproj}, \ie\ by the third-order ``Tanno equation'' \cite{m2010-gallottanno, mm2010-gtincomplete, t1978-someeqns}
\vspace{-0.2em}
\begin{equation*}
  \ell \mapsto \algbrac{ \ve^i }{ \D(\D_{e_i}\D\ell) - (\D\ell)\ltens\nRic{e_i}{}
    - \ell\ltens(\D\nRic{}{})_{e_i} } - \algbrac{ \nRic{}{} }{ \D\ell }.
\vspace{-0.2em}
\end{equation*}
This is precisely the third-order operator considered in Subsection \ref{ss:cproj-bgg-hess}.  This difference is already visible from the representation theory: for type \type{C}{n+1} we have
\vspace{-0.2em}
\begin{equation*}
  \liehom{0}{\bW^*} = \dynkinAp[0,0,0,0,2]{2}{0}{3}{ppg-bgg-hess-Chom0}
  \quad\text{and}\quad
  \liehom{1}{\bW^*} = \dynkinAp[0,0,0,3,-4]{2}{0}{3}{ppg-bgg-hess-Chom1},
\end{equation*}
which have weights $\tfrac{2n}{n+1}$ and $-\tfrac{n+3}{n+1}$ respectively, so that the first BGG operator has order $\tfrac{1}{n+1}(2n+(n+3)) = 3$.  There is in fact a second order operator defined on $\cL^{1/2}$, but we will not need this here; for $r=1$ we shall understand that the hessian is trivial.

\begin{defn} \thlabel{defn:ppg-bgg-hesseqn} Equation \eqref{eq:ppg-bgg-hesseqn} is called the \emph{hessian equation}. \end{defn}

\begin{rmk} In the conformal case ($r=n$), the representation $\bW$ is isomorphic to its dual via the conformal metric.  Therefore the conformal hessian and conformal metric equation coincide: they are both the \einstein\ scale equation
\smallskip\vspace{0.2em}
\begin{equation*}
  \bgg{\bW^*} : \dynkinBDp[1,0,0,0,0,0,0,0]{2}{0}{2}{ppg-bgg-hess-BDdom}
    \to \dynkinBDp[-3,2,0,0,0,0,0,0]{2}{0}{2}{ppg-bgg-hess-BDim}
\vspace{0.2em}
\end{equation*}
parametrising the family of \einstein\ scales. \end{rmk}

As for the \cproj\ and \qtn ic\ hessians, solutions of the hessian equation may be characterised using exact Weyl structures.

\begin{prop} \thlabel{prop:ppg-bgg-hesseqn} A nowhere-vanishing section $\ell \in \s{0}{\cL}$ satisfies $\bgg{\bW^*}(\ell) = 0$ if and only if the normalised Ricci tensor $\nRic[\D^{\ell}]{}{}$ of $\D^{\ell}$ is a section of $\cB^*$. \end{prop}

\begin{proof} Suppose first that $\bgg{\bW^*}(\ell) = 0$.  Then calculating \wrt\ $\D^{\ell}$ gives $\ell \ltens \nRic[\D^{\ell}]{}{} = \tfrac{1}{r} \algbrac{\ve^i}{\nRic[\D^{\ell}]{e_i}{}}$, where using \thref{cor:ppg-bgg-exactweyl} we view $\nRic[\D^{\ell}]{}{}$ as a section of $\Symm{2}T^*M$ on the \lhs, and as a $T^*M$-valued $1$-form on the \rhs.  In particular $\ell \ltens \nRic[\D^{\ell}]{}{}$ is a section of $\cL \tens \cB^*$ which, since $\ell$ is nowhere-vanishing, gives that $\nRic[\D^{\ell}]{}{}$ is a section of $\cB^*$.  On the other hand, suppose that $\nRic[\D^{\ell}]{}{}$ is a section of $\cB^*$.  Then
\vspace{0.2em}
\begin{equation*} \begin{aligned}
  \algbrac{\ve^i}{\ell \ltens \nRic[\D^{\ell}]{e_i}{}}(X,Y)
    &= \algbrac{ \algbrac{ \algbrac{\ve^i}{\ell \ltens \nRic[\D^{\ell}]{e_i}{}} }
                         { X } }
               { Y } \\
    &= \algbrac{ \algbrac{ \algbrac{ \ve^i }
                                   { \algbrac{e_i}{\ell \ltens \nRic[\D^{\ell}]{}{}} } }
                         { X } }
               { Y } \\
    &= \algbrac{ \algbrac{ \algbrac{ \algbrac{\ve^i}{e_i} }
                                   { \ell \ltens \nRic[\D^{\ell}]{}{} } }
                         { X } }
               { Y } \\
    &= - \algbrac{ \algbrac{ \algbrac{ \algbrac{e_i}{\ve^i} }
                                     { X } }
                           { \ell \ltens \nRic[\D^{\ell}]{}{} } }
                 { Y }
       - \algbrac{ \algbrac{ \algbrac{e_i}{\ve^i} }
                           { \ell \ltens \nRic[\D^{\ell}]{X}{} } }
                 { Y } \\
    &= - \tfrac{1}{2}r(n+1) \algbrac{ \algbrac{X}{\ell \ltens \nRic[\D^{\ell}]{}{}} }
                                    { Y }
       - (rn-\tfrac{1}{2}r(n+1)) \algbrac{ \ell \ltens \nRic[\D^{\ell}]{X}{} }
                                         { Y } \\
    &= \tfrac{1}{2}r(n+1) \ell \ltens \nRic[\D^{\ell}]{X}{Y}
       - (rn - \tfrac{1}{2}r(n+1)) \ell \ltens \nRic[\D^{\ell}]{X}{Y} \\
    &= r \, \ell \ltens \nRic[\D^{\ell}]{X}{Y},
\end{aligned}
\vspace{0.2em}
\end{equation*}
which implies that $\bgg{\bW^*}(\ell) = \ell \ltens \nRic[\D^{\ell}]{}{} - \tfrac{1}{r} r \ell \ltens \nRic[\D^{\ell}]{}{} = 0$ as required. \end{proof}

\begin{cor} \thlabel{cor:ppg-bgg-metrichess} Let $h \in \s{0}{\cL^*\tens\cB}$ be a \non degenerate\ solution of the linear metric equation.  Then $\ell \defeq (\det h)^{1/r} \in \s{0}{\cL}$ is a solution of the hessian. \end{cor}

\begin{proof} If $h$ is \non degenerate\ then $\pi \defeq (\det h)^{1/r}$ is a nowhere-vanishing section of $\cL$ which is parallel \wrt\ the \LC\ connection $\D^g$ of $h$, and $\nRic[g]{}{} \in \s{0}{\cB^*}$ by \thref{thm:ppg-bgg-nRicB}.  The result now follows from \thref{prop:ppg-bgg-hesseqn}. \end{proof}

Finally, the general theory gives an isomorphism between the solution space of the hessian equation and the space of parallel sections of a prolongation connection on $\cW^*$.

\begin{thm} \thlabel{thm:ppg-bgg-hessianprol} For $r>1$, there is a linear isomorphism between the solutions of the hessian equation \eqref{eq:ppg-bgg-hesseqn} and the parallel sections of the prolongation connection
\vspace{-0.3em}
\begin{equation} \label{eq:ppg-bgg-hessprol}
  \D^{\cW^*}_X \! \colvect{ \ell \\ \eta \\ \theta }
  = \colvect{ \D_X \ell   + \eta(X) \\
              \D_X \eta - \theta(X,\bdot)               + \ell \ltens \nRic{X}{} \\
              \D_X \theta  + \algbrac{\nRic{X}{}}{\eta}
            }
  - \quab^{-1} \!\!
    \colvect{ 0 \\
              0 \\
              \Weyl{e_i,X}{ \algbrac{\ve^i}{\eta} } + \algbrac{\ve^i}{ \ell \ltens \CY{e_i,X}{} }
            }
\vspace{-0.3em}
\end{equation}
on sections of $\cW^* \isom \cL \dsum (\cL\tens T^*M) \dsum (\cL\tens \cB^*)$.  The isomorphism is given explicitly by the splitting operator \eqref{eq:ppg-bgg-hesssplit}.  For $r=1$ the hessian equation is trivial, with solution space isomorphic to $\s{0}{\cL}$. \end{thm}

\begin{proof} The distinction between $r=1$ and $r>1$ is discussed above; for $r>1$ this is a straightforward application of the general theory.  By \thref{lem:ppg-bgg-hesstractor} we have
\vspace{-0.3em}
\begin{equation} \label{eq:ppg-bgg-hessprol-1}
  \liebdy \! \left( \Curv[\bW^*]{}{ \colvect{ \ell \\ \eta \\ \theta } } \right)
    = \colvect{ 0 \\
                0 \\
                \liebdy( \Weyl{}{\eta} + \ell \ltens \CY{}{} ) }
\vspace{-0.3em}
\end{equation}
Since the nilpotent differential operator $\quab_M - \quab$ appearing in the Neumann series \eqref{eq:para-bgg-neumann} necessarily lowers the weight, it acts trivially in \eqref{eq:ppg-bgg-hessprol-1}.  Therefore
\vspace{-0.3em}
\begin{equation*}
  \quab_M^{-1} \liebdy \! \left(
    \Curv[\bW^*]{}{ \colvect{ \ell \\ \eta \\ \theta } } \right)
    = \quab^{-1} \!\!
      \colvect{ 0 \\
                0 \\
                \liebdy( \Weyl{}{\eta} + \ell \ltens \CY{}{}) }
    \colvectpunct{.}
\vspace{-0.3em}
\end{equation*}
To obtain \eqref{eq:ppg-bgg-hessprol}, it remains to calculate the action of $\liebdy$.  This is immediate for the Cotton--York term, while for the Weyl term we have
\vspace{-0.2em}
\begin{equation*}
  \liebdy \algbrac{\Weyl{}{}}{\eta}_X
    = \sum{i}{} \, \algbrac{ \ve^i }{ \algbrac{\Weyl{e_i,X}{\eta}} }
    = \sum{i}{} \, \Weyl{e_i,X}{ \algbrac{\ve^i}{\eta} }
\vspace{-0.2em}
\end{equation*}
by the Jacobi identity, since $\liebdy \Weyl{}{} = 0$ implies that $\algbrac{\ve^i}{\Weyl{e_i,X}{}} = \ve^i \circ \Weyl{e_i,X}{} = 0$. \end{proof}

Unfortunately the author was unsuccessful in finding a general expression for the action of $\quab^{-1}$ in \eqref{eq:ppg-bgg-hessprol}.  Closed expressions are given by \thref{thm:cproj-bgg-hessprol,thm:qtn-bgg-hessprol} in \cproj\ and \qtn ic\ geometries.  For octonionic geometry, the harmonic curvature consists solely of the Cartan torsion $\Tor{}$; in particular the Weyl curvature vanishes identically.  The prolongation of the conformal hessian (\ie\ the \einstein\ scale equation) is discussed by Hammerl in \cite{h2008-confbgg, h2009-naturalprol}.

\chapter{\Ppgs\ of mobility two} 
\label{c:mob2}

\renewcommand{\algbracadornment}{}
\BufferDynkinLocaltrue
\renewcommand{\dynkinnameoffset}{-0.75}

The classical formulations of \proj, \cproj\ and \qtn ic\ geometries proceed by assuming the existence of a ``background'' metric and imposing an appropriate equivalence condition.  The linear metric equations then become the main equations \eqref{eq:proj-class-maineqn}, \eqref{eq:cproj-class-maineqn} and \eqref{eq:qtn-class-maineqn}, with equivalent metrics inducing a $1$-parameter family called a \emph{metrisability pencil}.  We describe a similar formulation for \ppg\ in Section \ref{s:mob2-pencils}, as well as generalising results of Matveev and Topalov \cite{mt1998-trajectory, tm2003-geodintegrability} on the geodesic flow.  A metrisability pencil also induces a $1$-parameter family of commuting vector fields, defined using the bilinear differential pairings from the BGG complex, which we study in Section \ref{s:mob2-vfs}.

Given a metrisability pencil, it is common to look at the eigenvalues of the corresponding solution $A$ of the main equation.  These turn out to be tightly controlled in the classical cases; as we shall see in Section \ref{s:mob2-evals}, similar results hold in the general case.  In \cproj\ geometry these eigenvalues give insight to a classification of pencils \cite{acg2006-ham2forms1}; unfortunately the author was unable to make any progress in this direction.

We fix, once and for all, a manifold $M$ supporting a \ppg\ with parameters $(r,n)$, where we assume that $n>0$ as per \thref{asmpt:ppg-alg-n>0}.  We also continue to assume that the Cartan torsion $\Tor{}$ vanishes.

\smallskip
\section{Metrisability pencils and integrability} 
\label{s:mob2-pencils}

The linearity of the metric equation \eqref{eq:ppg-bgg-metriceqn1} implies that any two linearly independent solutions $h,\b{h}$ span a two-dimensional real vector space $V$ of solutions.  More invariantly, we can avoid choosing a basis by supposing that a \ppg\ admits a two-dimensional family $\bs{h}$ of linear metrics parameterised by $V$.

\begin{defn} \thlabel{defn:mob2-pencils} A \emph{(metrisability) pencil} is a family $\bs{h} \in\s{0}{\cL^*\tens\cB\tens V^*}$ of linear metrics parameterised by a two-dimensional real vector space $V$. \end{defn}

Sometimes it will be fruitful to ``projectivise'' \thref{defn:mob2-pencils} as follows.  Let $\cO_V(1)$ denote the tautological line bundle over the projective line $\pr{V}$, whose fibre $\cO_V(1)_{[v]}$ over $[v]\in\pr{V}$ is the line spanned by $v$.  There is a natural map $\pr{V} \times V^* \to \cO_V(1)$ given by $([v],\alpha) \mapsto \alpha(v) \in\cO_V(1)_{[v]}$, which induces an isomorphism of $V^*$ with the space of algebraic functions $\pr{V} \to \cO_V(1)$ and therefore a \emph{natural lift} $\nat : \pr{V} \to V\tens \cO_V(1)$.  We may then identify a pencil $\bs{h}$, which is by definition a section of $\cL^*\tens\cB\tens V^*$ over $M$, with a section of $\cL^*\tens\cB \tens\cO_V(1)$ over $M\times \pr{V}$.  Indeed, the terminology ``pencil'' more typically refers to the \proj\ line $\pr{V}$.

While the previous definition is more natural, the concreteness obtained by choosing a basis $\{h,\b{h}\}$ of $V$ is often beneficial.  This picture may be projectivised by identifying $h,\b{h}$ with the points at $\infty,0$ respectively in an affine chart for $\pr{V}$ with affine parameter $t$; then any linear metric in the pencil, other than $h$, is proportional to $h_t \defeq \b{h}-th$ for some $t\in\bR$.  As a final ingredient, introduce the endomorphism $A = A(h,\b{h})$ of $T^*M$ defined by $\b{h} = h(A\bdot,\bdot)$.  Clearly $A$ is \self adjoint\ \wrt\ both $h,\b{h}$, and moreover defines a section of $\fp^0_M$ by \thref{lem:ppg-bgg-p0B}.  It follows that $A$ is \self adjoint\ \wrt\ any metric $h_t$ in the pencil; of course, we may write $h_t = h(A_t\bdot,\bdot)$ for $A_t \defeq A-t\id$.

\begin{rmk} For conformal geometries, $\cL^* \tens \cB \isom \cL$ is spanned by the (inverse) conformal metric $\conf$.  Since here $\fp^0_M = \alg{co}{TM} = \alg{so}{TM} \dsum (M\times\bR)$, any \self adjoint\ section of $\fp^0_M$ is necessarily a (functional) multiple of the identity.  However, note that this does not obstruct the existence of metrisability pencils: for conformal geometries the first-order linear metric equation \eqref{eq:ppg-bgg-metriceqn1} is trivial, and the second-order equation \eqref{eq:ppg-bgg-confmetric} coincides with the hessian equation \eqref{eq:ppg-bgg-hesseqn}.  Linear metrics are then \einstein\ metrics in the conformal class; pencils of \einstein\ metrics do exist, although they are rare \cite{b1924-confeinstein, b1925-confeachother}. \end{rmk}

\subsection{The pfaffian of a linear metric} 
\label{ss:mob2-pencils-pf}

Here we formalise a concept used previously.  Recall the $\fg$-representation $\bW \defeq \fh/\fq$, and consider its zeroth homology $\liehom{0}{\bW} \isom L^*\tens B$, where $B$ is a $\fp^0$-subrepresentation of $\Symm{2}(\fg/\fp)$.  Any element $h \in L^*\tens B$ may be viewed as a linear map $h : \fp^{\perp} \to L^*\tens \fg/\fp$, whose determinant is a linear map
\begin{equation*}
  \det h : L^{-r(n+1)/2} \to L^{-rn} \tens L^{r(n+1)/2}
\end{equation*}
since $\Wedge{rn}(\fg/\fp) \isom L^{r(n+1)/2}$ by \thref{prop:ppg-alg-Lwedge}.  By an abuse of notation we obtain an element $\det h \in L^r$ by dualising, which coincides with the square of the volume form of $h$.  As usual, $\det h$ vanishes if and only if $h$ is degenerate.

\begin{defn} \thlabel{defn:mob2-pencils-pf} The \emph{pfaffian} of $h$ is the element $\pf{h} \defeq (\det h)^{1/r}$ of $L$. \end{defn}

In terms of associated bundles, the determinant of a section $h$ of $\cL^*\tens\cB$ defines a section of $\cL^r$, yielding a pfaffian $\pf{h} \in \s{0}{\cL}$.  For a pencil $\bs{h}$ of linear metrics, $\bs{h}$ is, by definition, a linear functional on $V$, so that $\det \bs{h}$ is a homogeneous polynomial of degree $rn$ on $V$.  We will see in \thref{cor:mob2-pencils-pfpoly} that the pfaffian $\bs{\pi} \defeq \pf{\bs{h}}$ defines a section of $\cL \tens \Symm{n} V^*$ over $M$,  or equivalently a section of $\cL \tens \cO_{V}(n)$ over $M\times\pr{V}$, thus giving a polynomial of degree $n$ in any affine chart.

As before, the pfaffian allows us to transform more efficiently between a (\non degenerate) linear metric $h \in \s{0}{\cL^*\tens\cB}$ and a metric in the usual sense.  Indeed, $g\defeq (\pf{h})^{-1} \ltens h^{-1}$ is a section of $\Symm{2}T^*M$, and taking the top exterior power gives
\vspace{0.2em}
\begin{equation*}
  \det g = (\pf{h})^{-rn} \ltens (\det h)^{-1} = (\pf{h})^{-r(n+1)}
\vspace{0.2em}
\end{equation*}
so that $h=(\det g)^{1/r(n+1)} \ltens g^{-1}$.  As an immediate result, we may give a more concrete relation between the linear metric equation and the main equations \eqref{eq:proj-class-maineqn}, \eqref{eq:cproj-class-maineqn} and \eqref{eq:qtn-class-maineqn} of \proj, \cproj\ and \qtn ic\ geometries.  Consider the endomorphism $A$ associated to two linear metrics $h, \b{h}$.  Substituting $h = (\det g)^{1/r(n+1)} \ltens g^{-1}$ and $\b{h} = (\det\b{g})^{1/r(n+1)} \ltens \b{g}^{-1}$ into the defining equation $\b{h} = h(A\bdot,\bdot)$, we see that
\vspace{0.2em}
\begin{equation} \label{eq:mob2-pencils-endoA}
  A = \left( \frac{\det\b{g}}{\det g} \right)^{1/r(n+1)} \b{\sharp}\circ\flat
\vspace{0.2em}
\end{equation}
is just the familiar endomorphism featured in these main equations.  Moreover, we have the following characterisation of the section $\b{Z}^{\D}$ of $\cL^*\tens TM$ defined by \eqref{eq:ppg-bgg-metricsplit}.

\begin{prop} \thlabel{prop:mob2-pencils-Zgradient} Let $\bs{h}$ be a metrisability pencil with \non degenerate\ linear metrics $h,\b{h}$ at $\infty,0$ in an affine chart, and form the endomorphism $A=A(h,\b{h})$.  Let $\D$ be the \LC\ connection of $g \defeq (\pf{h})^{-1/r} \ltens h^{-1}$ and consider the section $\b{Z}^{\D}$ of $\cL^*\tens TM$ satisfying $\D_X\b{h} = \algbrac{\b{Z}^{\D}}{X}$.  Then $(\pf{h})\b{Z}^{\D} = \tfrac{1}{r}\grad[g]{(\tr A)}$. \end{prop}

\begin{proof} Differentiating the identity $\b{h} = h(A\bdot,\bdot)$ \wrt\ $\D$ yields $\D_X\b{h} = h((\D_X A)\bdot,\bdot)$, so that $\D_X A = h^{-1}\circ\algbrac{\b{Z}^{\D}}{X}$.  Taking a trace \wrt\ a local frame $\{e_i\}_i$ with dual coframe $\{\ve^i\}_i$, the \lhs\ gives $\left( (\D_X A)(\ve^i) \right)(e_i) = \d(\tr A)(X)$.  The \rhs\ gives
\begin{align*}
  h^{-1} \big( \algbrac{\b{Z}^{\D}}{X}(\ve^i,\bdot), e_i \big)
    &= \killing{ \algbrac{ e_i }{ h^{-1} } }
               { \algbrac{ \algbrac{\b{Z}^{\D}}{X} }{ \ve^i } } \\
    &= \killing{ h^{-1} }
               { \algbrac{ \algbrac{ \algbrac{\b{Z}^{\D}}{X} }
                                   { \ve^i } }
                         { e_i } } \\
    &= \killing{ h^{-1} }
               { \algbrac{ \algbrac{e_i}{\ve^i} }
                         { \algbrac{\b{Z}^{\D}}{X} } }
\end{align*}
by the invariance of $\killing{}{}$ and the Jacobi identity.  Since $\algbrac{ \algbrac{e_i}{\ve^i} }{ X } = \tfrac{1}{2}r(n+1) X$ for all vector fields $X$ by \thref{cor:ppg-alg-trace}, we obtain $\algbrac{ \algbrac{e_i}{\ve^i} }{ \b{Z}^{\D} } = -\tfrac{1}{2}r(n-1) \b{Z}^{\D}$ by the Leibniz rule and hence
\begin{align*}
  \Algbrac{ \algbrac{e_i}{\ve^i} }{ \algbrac{\b{Z}^{\D}}{X} }
    &= (-\ve^i(e_i) + \tfrac{1}{2}r(n+1) \algbrac{\b{Z}^{\D}}{X}
       + \tfrac{1}{2}r(n+1) \algbrac{\b{Z}^{\D}}{X} \\
    &= r \algbrac{\b{Z}^{\D}}{X}
\end{align*}
by the Jacobi identity.  Then $\tr \big( h^{-1}\circ \algbrac{\b{Z}^{\D}}{X} \big) = r \killing{ h^{-1} }{ \algbrac{\b{Z}^{\D}}{X} } = r h^{-1}(\b{Z}^{\D},X)$.  Thus $\d(\tr A) = r h^{-1}(\b{Z}^{\D}, \bdot)$, so that writing $h^{-1} = (\pf{h}) \ltens g$ and applying $\sharp$ to both sides completes the proof. \end{proof}

\vspace{0.1em}
\subsection{Adjugate tensors} 
\label{ss:mob2-pencils-adj}

In addition to the pfaffian, there is another operation that we can perform on linear metrics.  By definition the flat model of a \ppg\ has a projective embedding $G\acts\fp \injto \pr{\bW}$, so that Kostant's \thref{thm:lie-para-kostant} tells us that $G\acts\fp$ is an intersection of quadrics whose defining equations are given by projection away from the Cartan square in $\Symm{2}\bW$.  We let $\bU \defeq \Symm{2}\bW / \Cartan{2}\bW$, so that \eqref{eq:lie-para-kostant} identifies $\bU^* \subseteq \Symm{2}\bW^*$ with the quadratic defining equations of $G\acts\fp$.
	
\begin{prop} \thlabel{prop:mob2-pencils-repnU} $\liehom{0}{\bU^*} \isom L^2 \tens B^*$. \end{prop}

\begin{proof} It suffices to consider the case that $\fh$ is simple.  Then by \thref{prop:ppg-alg-gr}, the $\fg$-representation $\bW^*$ has associated graded representation $\bW^* \isom L \dsum (L\tens \fp^{\perp}) \dsum (L\tens B^*)$ \wrt\ any algebraic Weyl structure for $\fp$.  Therefore
\begin{equation} \label{eq:mob2-pencils-repnU-1} \begin{aligned}
  \Symm{2}\bW^*
  \isom L^2 \dsum (&L^2\tens \fp^{\perp}) \dsum (L^2\tens B^*)
    \dsum (L^2\tens \Symm{2}\fp^{\perp}) \\
  & \dsum (L^2\tens \fp^{\perp}\tens B^*) \dsum (L^2\tens \Symm{2}B^*),
\end{aligned} \end{equation}
so it suffices to exhibit $L^2 \tens B^*$ as the highest weight summand \emph{not} lying in $\Cartan{2} \bW^*$.  By \eqref{eq:ppg-alg-eigenvals}, the summands of \eqref{eq:mob2-pencils-repnU-1} are written in order of \non increasing\ weight.

Since $\Cartan{2}\bW^*$ is the highest weight $\fg$-subrepresentation of $\Symm{2}\bW^*$, its associated graded representation must contain the highest weight summand $L^2$ of $\Symm{2}\bW^*$.  Dualising, it follows also that $\Cartan{2}\bW^*$ contains the lowest weight irreducible summand, which is an irreducible $\fp^0$-subrepresentation of $L^2 \tens \Symm{2}B^*$.  Then the associated graded of $\Cartan{2}\bW^*$ must contain at least one irreducible summand of each weight between the weights of its highest and lowest components; \cf\ \cite[Eqn.\ (2.3)]{cds2005-ricci}.  The only summand in \eqref{eq:mob2-pencils-repnU-1} of the same weight as $\fp^{\perp} \acts L^2$ is $L^2 \tens \fp^{\perp}$, so this must lie in $\Cartan{2}\bW^*$.  In all cases except $r=1$, Table \ref{tbl:ppg-class-cpx} implies that $B^*$ is complementary to the Cartan square in $\Symm{2}\fp^{\perp}$; thus in these cases $L^2 \tens B^*$ cannot lie in $\Cartan{2}\bW^*$, so is the highest weight summand.  For $r=1$, we have $B^* = \Symm{2} \fp^{\perp}$ and thus have two copies of $L^2 \tens \Symm{2}\fp^{\perp}$.  The kernel of the anti-diagonal action of $\fp^{\perp}$ on $(L^2 \tens B^*) \dsum (L^2 \tens B^*)$ lies in $\Cartan{2}\bW^*$ by the results of \cite{e2005-cartanprod}, so that its complement forms a highest weight component of $\bU^*$ isomorphic to $L^2 \tens B$.  Thus in all cases, $L^2 \tens B^*$ is a highest weight summand in $\gr{(\bU^*)}$.  Since $B^*$ is irreducible and the remaining summands in \eqref{eq:mob2-pencils-repnU-1} have strictly lower weights, we conclude that $\liehom{0}{\bU^*} \isom L^2 \tens B^*$. \end{proof}

Let us calculate the complex representation $\bU^* \leq \Symm{2}\bW^*$ for each irreducible complex \ppg, as well as its graded components, using Table \ref{tbl:ppg-class-cpx}.  For notational convenience, we identify all representations and Lie algebras with their complexifications.  This information is summarised in Table \ref{tbl:tbl-app-Ud}.

\begin{typeblock}{C}{n+1} Here $\fg = \alg{sl}{n+1,\bC}$ and $\bW^* = \Symm{2}\bC^{n+1*}$, giving $\Symm{2}\bW^* \isom \Symm{4}\bC^{n+1*} \dsum \Cartan{2}\Wedge{2}\bC^{n+1*}$.  The Cartan square is $\Cartan{2}\bW^* = \Symm{4}\bC^{n+1*}$, so that
\begin{align*}
  \bU^*
    &= \dynkinSLR[0,0,0,2,0]{2}{0}{3}{mob2-pencils-adj-CUd}
       = \Cartan{2}\Wedge{2}\bC^{n+1*} \\
    &\isom \dynkinname{ \dynkinSLRp[0,0,0,2,0]{2}{0}{3}{mob2-pencils-adj-CUd1} }
                      { L^2 \tens B^* }
     \dsum \dynkinname{ \dynkinSLRp[0,0,0,1,1,-1]{2}{0}{4}{mob2-pencils-adj-CUd2} }
                      { L^2 \tens \fp^{\perp} \cartan \Wedge{2}\fp^{\perp} }
     \dsum \dynkinname{ \dynkinSLRp[0,0,0,2,0,-2]{2}{0}{4}{mob2-pencils-adj-CUd3} }
                      { L^2 \tens \Cartan{2}\Wedge{2}\fp^{\perp} }.
\end{align*}
\end{typeblock}

\begin{typeblock}{A}{2n+1} Here $\fg = \alg{sl}{n+1,\bC} \dsum \alg{sl}{n+1,\bC}$ and $\bW^* = \bC^{n+1*} \etens \conj{\bC^{n+1*}}$, giving $\Symm{2}\bW^* \isom (\Symm{2}\bC^{n+1*} \etens \Symm{2}\conj{\bC^{n+1*}}) \dsum (\Wedge{2}\bC^{n+1*} \etens \Wedge{2}\conj{\bC^{n+1*}})$.  The first summand here is the Cartan square, so that
\begin{align*}
  \bU^*
    &= \dynkinAA[0,0,0,1,0,0,0,0,1,0]{2}{0}{3}{mob2-pencils-adj-AUd}
        = \Wedge{2}\bC^{n+1*} \etens \Wedge{2}\conj{\bC^{n+1*}} \\
    &\isom \dynkinname[-0.5 ]{ \dynkinAAp[0,0,0,1,0,0,0,0,1,0]{2}{0}{3}
                                         {mob2-pencils-adj-AUd1} }
                             { L^2 \tens B^* }
     \dsum \dynkinname[-0.5 ]{ \left[ \dynkinAAp[0,0,0,1,0,-1,0,0,0,0,1,0]{2}{0}{4}
                                               {mob2-pencils-adj-AUd2}
                                      \dsum \cpxconj \right] }
                             { \zbox{(L^2 \tens \Wedge{1,0}M\cartan\Wedge{0,2}M)
                                     \dsum \cpxconj }}
     \dsum \dynkinname[-0.5 ]{ \dynkinAAp[0,0,0,1,0,-1,0,0,0,1,0,-1]{2}{0}{4}
                                        {mob2-pencils-adj-AUd3} }
                             { L^2 \tens \Wedge{1,1}M }
\end{align*}
where by $\cpxconj$ we mean the complex-conjugate representation. \end{typeblock}

\begin{typeblock}{D}{2n+2} Here $\fg = \alg{so}{4n+4,\bC}$ and $\bW^* = \Wedge{2}\bC^{2n+2*}$, thus giving $\Symm{2}\bW^* \isom \Cartan{2} \Wedge{2}\bC^{2n+2*} \dsum \Wedge{4}\bC^{2n+2*}$.  By \cite{e2005-cartanprod} the Cartan square consists of the \tracefree\ symmetric elements; therefore
\begin{align*}
  \bU^*
    &= \dynkinA[0,0,1,0,0,0]{1}{0}{5}{mob2-pencils-adj-DUd}
       = \Wedge{4}\bC^{2n+2*} \\
    &\isom \dynkinname{ \dynkinApp[0,0,1,0,0,0]{1}{0}{5}{mob2-pencils-adj-DUd1} }
                      { \Wedge{2}E^* \tens \Wedge{2}H }
     \dsum \dynkinname{ \dynkinApp[0,0,1,0,0,-1,1]{1}{0}{6}{mob2-pencils-adj-DUd2} }
                      { \Wedge{3}E^* \tens H }
     \dsum \dynkinname{ \dynkinApp[0,0,1,0,0,0,-1,0]{1}{0}{7}{mob2-pencils-adj-DUd3} }
                      { \Wedge{4}E^* }.
\end{align*}
The representations $E \isom \bC^{2n}$ and $H \isom \bC^2$ are defined in \eqref{eq:qtn-para-EH}.
\end{typeblock}

\begin{typeblock}{E}{7} Here $\fg = \alg[_6]{e}{\bC}$ and $\bW^*$ is dual to the $27$-dimensional representation $\erepn$ of $\alg[_6]{e}{\bC}$.  Calculating using \texttt{LiE},%
\footnote{An online service is available at \url{http://wwwmathlabo.univ-poitiers.fr/~maavl/LiE/form.html}.}
the symmetric square decomposes as
\vspace{0.2em}
\begin{equation*}
  \Symm{2} \! \left( \dynkinE[0,0,0,0,1,0]{6}{mob2-pencils-adj-ES2} \right)
    \isom \dynkinE[0,0,0,0,2,0]{6}{mob2-pencils-adj-EO2}
    \dsum \dynkinE[1,0,0,0,0,0]{6}{mob2-pencils-adj-EU}
    \colvectpunct[-0.7em]{\, ,}
\vspace{0.2em}
\end{equation*}
where the first summand is the Cartan square and the second is $\erepn[\bC]$.  Therefore
\vspace{0.1em}
\begin{equation*}
\begin{aligned}
  \bU^*
    &= \dynkinE[1,0,0,0,0,0]{6}{mob2-pencils-adj-EUd}
       = \erepn[\bC] \\
    &\isom \dynkinname[-0.15]{ \dynkinEp[1,0,0,0,0,0]{6}{mob2-pencils-adj-EUd1} }
                             { L^* \tens B }
     \dsum \dynkinname[-0.15]{ \dynkinEp[0,0,0,0,-1,1]{6}{mob2-pencils-adj-EUd2} }
                             { L^* \tens \fg/\fp }
     \dsum \dynkinname[-0.15]{ \dynkinEp[0,0,0,0,-1,0]{6}{mob2-pencils-adj-EUd3} }
                             { L^* }
     \isom \bW.
\end{aligned}
\vspace{0.1em}
\end{equation*}
Notice in particular that $\bU^* \isom \bW$. \end{typeblock}
  
\begin{typeblock}{BD}{n+4} Here $\fg = \alg{so}{n+2,\bC}$ and $\bW^*=\bC^{n+2*}$.  Since the Cartan product in the orthogonal algebras is the \tracefree\ symmetric product, $\Symm{2} \bW^* \isom \Symm[\trfree]{2} \bW^* \dsum \bC$.  Therefore $\bU^*$ is the trivial representation,
\vspace{0.5em}
\begin{equation*}
  \bU^*
    = \dynkinSO[0,0,0,0,0,0,0,0]{2}{0}{2}{mob2-pencils-adj-BDUd}
    = \bC
\vspace{0.2em}
\end{equation*}
This conforms with what we know already: in type \type{BD}{n+4} we have $\bW \isom \bW^*$ and hence $L^2 \tens B^* \isom L \tens L^* \isom \bC$ must be the trivial representation. \end{typeblock}

Recalling that $L^*\tens B = \liehom{0}{\bW}$ and $L = \liehom{0}{\bW^*}$, the pairing $(L^2\tens B^*) \times (L^*\tens B) \to L$ is really a pairing on Lie algebra homology, so is induced by a $\fp$-invariant pairing $\bU^*\times \bW \to \bW^*$ of $\fg$-representations.  Dualising then yields a map $\adj{} : L^*\tens B \to L^2\tens B^*$ such that the contraction $\killing{\adj h}{h}$ equals  $\pf{h}$ for all $h\in L^*\tens B$.

\begin{defn} \thlabel{defn:mob2-pencils-adj} The image of $h\in L^*\tens B$ under the map $\adj{} : L^*\tens B \to L^2\tens B^*$ is called the \emph{adjugate} of $h$, denoted by $h^*\defeq \adj{h}$. \end{defn}

Clearly if $h\in L^*\tens B$ is \non degenerate\ then the condition $\killing{h^*}{h} = \pf{h}$ implies that $h^* = (\pf{h}) \ltens h^{-1} = (\det h)^{1/r} \ltens h^{-1}$.

\begin{prop} \thlabel{prop:mob2-pencils-adjpoly} Given linearly independent elements $h, \b{h} \in L^*\tens B$, the adjugate of $h_t \defeq \b{h} - th$ is a polynomial of degree $n-1$ in $t$. \end{prop}
		
\begin{proof} We complexify $h_t$ and proceed on a case-by-case basis using \thref{thm:ppg-class-cpxclass}.  In each case, the construction is evidently polynomial of degree $n-1$ in $t$.

For $r=1$, we view $h_t$ as a linear map $\fp^{\perp} \to L^* \tens \fg/\fp$.  Via the \non degenerate\ pairing $\wedge: \fg/\fp \times \Wedge{n-1} \fg/\fp \to \Wedge{n} \fg/\fp \isom L^{r(n+1)/2}$, the $(n-1)$-fold wedge power $\Wedge{n-1} h_t : \Wedge{n-1} \fp^{\perp} \to L^{-n+1} \tens \Wedge{n-1} \fg/\fp$ may be viewed as the adjugate linear map $\fg/\fp \to L^2 \tens \fp^{\perp}$.  The contraction with $h_t$ is the $n$-fold wedge power, which is $\det h_t = \pf{h_t}$ by definition.

For $r=2$ we have $\cpxrepn{(L^* \tens B)} \isom L^{-1,-1} \tens \Symm[\bC]{1,1}(\fg/\fp)$, where we retain notation from Chapter \ref{c:cproj}.  Then $h_t$ defines a conjugate-linear map $h_t : \Wedge[\bC]{1,0} \to L^{-1,-1} \tens T^{0,1}$, whose $(n-1)$-wedge power yields a conjugate-linear adjugate map $T^{1,0} \to L^{2,2} \tens \Wedge[\bC]{0,1}$ via the pairings $\Wedge{1,0} \times \Wedge[\bC]{n-1} \Wedge{1,0} \to \Wedge[\bC]{n} \Wedge{1,0} \isom L^{-(n+1),0}$ and $T^{1,0} \times \Wedge[\bC]{n-1} T^{1,0} \to \Wedge[\bC]{n} T^{0,1} \isom L^{0,n+1}$.  This may equivalently be viewed as an element of $\cpxrepn{(L^2 \tens B^*)}$, whose contraction with $h$ is the usual matrix pfaffian $\pf{h} = (\det h)^{1/2}$.

For $r=4$, the complexification of $\fg/\fp$ decomposes as $\cpxrepn{(\fg/\fp)} = E \etens H$ for $E \isom \bC^{2n}$ and $H \isom \bC^2$, where $\Wedge[\bC]{2n} E \isom \Wedge[\bC]{2} H$.  The complexification of $h_t$ may be viewed as an element of $\Wedge[\bC]{2}E$, so that its $(n-1)$-fold wedge power defines an element of $\Wedge[\bC]{2n-2}E \isom \Wedge[\bC]{2}E^* \etens \Wedge[\bC]{2}H \isom \cpxrepn{(L^2 \tens B^*)}$ via the induced pairing $\Wedge[\bC]{2}E \times \Wedge[\bC]{2n-2} \to \Wedge[\bC]{2}H$.  The contraction with $h_t$ is an element of $\Wedge[\bC]{n}(\Wedge[\bC]{2}E) \isom \Wedge[\bC]{2}H$ which equals $\pf{h_t} = (\det h_t)^{1/4}$.

For $r=8$, we have $n=2$ and $h_t$ is an element of the standard representation of $\fp \isom \alg{co}{10,\bC}$, defining metrics on the $16$-dimensional half-spin representations.  These are dual, with contraction the conformal quadratic form on the standard representation $\bC^{10}$.  This must be a \non zero\ multiple of $\pf{h_t} = (\det h_t)^{1/8}$ by homogeneity, so that $h_t^* = h_t$ up to normalisation.  Finally, $L^2 \tens B^*$ is the trivial representation for $r=n$ and thus $h_t^* = 1$ up to normalisation. \end{proof}

\begin{cor} \thlabel{cor:mob2-pencils-pfpoly} The pfaffian $\pf{h_t}$ is a polynomial of degree $n$ in $t$. \end{cor}

\begin{proof} Since $h_t \defeq \b{h} - th$ is affine in $t$, the result follows immediately from \thref{prop:mob2-pencils-adjpoly} and the defining relation $\killing{h_t^*}{h_t} = \pf{h_t}$. \end{proof}

Moving to associated bundles, $h_t^*$ and $\pf{h_t}$ are polynomials of degrees $n-1$ and $n$ in $t$, whose coefficients are sections of $\cL^2 \tens \cB^*$ and $\cL$ respectively.

\vspace{-0.2em}
\subsection{Integrals of the geodesic flow} 
\label{ss:mob2-pencils-integrals}

Recall from Chapter \ref{c:proj} that the Weyl connections of \proj\ differential geometry are those linear connections with the same (unparameterised) geodesics.  If a \proj\ manifold admits a pencil of compatible metrics, it does not seem unreasonable to expect the dynamics of the (co)geodesic flow to have special properties.  In this subsection we describe the geodesic flow of a metric as a hamiltonian mechanical system and recall some related results from \proj\ and \cproj\ geometry.  

We first recall the basic elements of symplectic geometry and hamiltonian mechanics; a detailed introduction may be found in \cite{s2001-symplectic}.  A \emph{symplectic structure} on a manifold $M$ is the assignment of a closed \non degenerate\ $2$-form $\omega\in\s{2}{}$.  Then $\omega$ induces an inverse $2$-vector $\omega^{-1}\in\s{0}{\Wedge{2}TM}$ and hence a \emph{Poisson bracket}
\begin{align*}
  \poisson{}{} : \s{0}{} \times \s{0}{} &\to \s{0}{} \\
  \poisson{f}{g} &\defeq \omega^{-1}(\d f,\d g),
\end{align*}
\ie\ a skew-symmetric bilinear pairing on $\s{0}{}$ satisfying the product rule $\poisson{fg}{h} = \poisson{f}{h}g + f\poisson{g}{h}$; the Poisson bracket also equips $\s{0}{}$ with a Lie bracket.  Then $h \mapsto \poisson{f}{h}$ is a derivation of $\s{0}{}$ for all $f\in\s{0}{}$, giving a \emph{hamiltonian vector field} $X_f$ which satisfies $\p_{X_f}h = \poisson{f}{h}$. Given a symplectic form $\omega$, we may also write $X_f = \omega^{-1}(\d f,\bdot)$.  The assignment $f\mapsto X_f$ is linear and satisfies $\liebrac{X_f}{X_g} = X_{\poisson{f}{g}}$.

A choice of $H\in\s{0}{}$ makes the triple $(M,\omega,H)$ into a \emph{hamiltonian mechanical system}, where $H$ is the \emph{hamiltonian} (or energy function).  This yields a favoured vector field $X_H$, the hamiltonian vector field of $H$, whose integral curves provide a $1$-parameter family of transformations of $M$ which describe the evolution of the hamiltonian system according to Hamilton's equations \cite[Eqn.\ 22.11]{l2003-smoothmfds}.  In this picture, a function $f\in\s{0}{}$ is preserved by the system if and only if it is constant along the integral curves of $X_H$; such a function is called an \emph{integral} of the hamiltonian system.

\begin{prop} \thlabel{prop:mob2-pencils-integrals} $f\in\s{0}{}$ is an integral of $(M,\omega,H)$ if and only if $\poisson{f}{H}=0$. \noproof \end{prop}

In particular $\poisson{H}{H}=0$ by skew-symmetry, which corresponds to conservation of energy.  If $f_1,f_2 \in \s{0}{}$ are integrals, the Jacobi identity implies that $\poisson{f_1}{f_2}$ is again an integral, so that the space of integrals is a Lie subalgebra of $\s{0}{}$.  Two integrals are said to be in \emph{involution} if they also Poisson-commute, \ie\ $\poisson{f_1}{f_2}=0$.  If $\dim M = 2n$ then the linear space of hamiltonian vector fields corresponding to integrals is isotropic, so has dimension at most $n$.  The system is  \emph{(Liouville) integrable} if it admits $n$ linearly independent integrals in involution.

We can apply this to the geodesic flow of a \riem\ manifold $(M,g)$.  The cotangent bundle $\pi:T^*M\to M$ of $M$ admits a canonical $1$-form defined by
\begin{equation*}
  \Theta_{\alpha}(X) \defeq \alpha(\pi_{*}(X))
\end{equation*}
for all $\alpha\in T^*M$, viewed as maps $\alpha:TM\to\bR$, and all $X\in\s[T^*M]{0}{TT^*M}$.  Clearly $\Omega=\d\Theta$ defines a symplectic form on $T^*M$, called the \emph{canonical symplectic form}, so that $(M,\Omega)$ is a symplectic manifold.  The inverse metric $g^{-1}\in\s{0}{\Symm{2}TM}$ may be viewed as a homogeneous quadratic function on $T^*M$, so may be used as a hamiltonian function, and the integral curves of its hamiltonian vector field $X_{g^{-1}}\in\s[T^*M]{0}{TT^*M}$ project to geodesics on $M$ under $\pi$.  Pulling everything back to $TM$ using the musical isomorphisms of $g$, we obtain the \emph{geodesic flow} of $(M,g)$.  It is straightforward to show that the $g^{-1}$ is constant along geodesics.

\begin{defn} A smooth function $f:TM\to\bR$ is an \emph{integral of $g$} if $s \mapsto f(\gamma'(s))$ is constant for all affinely parameterised geodesics $\gamma$. \end{defn}

For homogeneous polynomial functions, we may reformulate being an integral of $g$ in terms of its \LC\ connection $\D^g$, see \cite[Prop.\ $1$]{t2003-geodesiccomp}.

\begin{lem} \thlabel{lem:mob2-pencils-quadint} A homogeneous polynomial function $Q:TM\to\bR$ is an integral of $g$ if and only if the symmetrisation ${\sym}{(\D^g Q)}$ of $\D^g Q$ vanishes. \noproof \end{lem}

Note that for homogeneous quadratic functions $Q \in \s{0}{\Symm{2}T^*M}$,
\vspace{-0.3em}
\begin{equation*}
  \sym(\D^g Q)(X,Y,Z)
    = (\D^g_X Q)(Y,Z) + (\D^g_Y Q)(Z,X) + (\D^g_Z Q)(X,Y)
\vspace{-0.3em}
\end{equation*}
for all $X,Y,Z \in TM$.

As remarked previously, one might expect a metrisability pencil on a \proj\ manifold (\ie\ $r=1$) to admit integrals of the metrics of the pencil.  Given metrics $g, \b{g}$ at $\infty,0$, a $1$-parameter family of homogeneous quadratic integrals defined by
\vspace{-0.3em}
\begin{equation*}
  I_t \defeq (\det(A-t\id)) g((A-t\id)^{-1}\bdot, \bdot)
\vspace{-0.3em}
\end{equation*}
were discovered by Matveev and Topalov \cite{mt1998-trajectory, tm2003-geodintegrability} and studied further by Topalov \cite{t2000-hierarchy, t2001-geodesichierarchies}, where $A$ is the usual endomorphism \eqref{eq:mob2-pencils-endoA}.  In particular, the $I_t$ are polynomial of degree $n-1$ in $t$ and mutually Poisson-commuting \wrt\ the Poisson bracket on $TM$ induced by $g$.  Similar integrals were discovered in the \cproj\ setting by Kiyohara and Topalov \cite{kt2011-kahlerliouville} (see also \cite{cemn2015-cproj}), where the Weyl connections are those linear connections with the same c-geodesics.  The \cproj\ integrals are also polynomial and mutually Poisson-commute.  In the next subsection, we will define an analogous $1$-parameter family of Poisson-commuting integrals of any metric of the pencil.

\vspace{-0.2em}
\subsection{Killing $2$-tensors} 
\label{ss:mob2-pencils-killing}

Consider a \non degenerate\ linear metric $h\in\s{0}{\cL^*\tens\cB}$ on $M$, with corresponding metric $g\defeq (\pf{h})^{-1} \ltens h^{-1}$.  Our next goal is to describe homogeneous polynomial integrals of $g$ independently of $\D\in \Dspace$.  This may be done by introducing weights to counteract the effect of changing the connection in $\Dspace$, and by using \thref{lem:mob2-pencils-quadint}.

\begin{prop} \thlabel{prop:mob2-pencils-killingeqn} For sections $k$ of $\cL^2 \tens \cB^*$, the equation
\vspace{-0.3em}
\begin{equation} \label{eq:mob2-pencils-killingeqn}
  (\D_X k)(Y,Z) + (\D_Y k)(Z,X) + (\D_Z k)(X,Y) = 0
\vspace{-0.3em}
\end{equation}
is independent of $\D\in \Dspace$. \end{prop}

\begin{proof} Since $\weyld{\gamma} \D_X k = \algbrac{X}{\gamma} \acts k$, the Leibniz rule gives
\begin{align} \label{eq:mob2-pencils-killingeqn-1}
  \weyld{\gamma}( (\D_X k)(Y,Z) )
    &= (\algbrac{X}{\gamma} \acts k)(Y,Z) \notag \\
    &= \algbrac{X}{\gamma} \acts k(Y,Z) - k(\algbrac{ \algbrac{X}{\gamma} }{ Y }, Z)
                                        - k(Y, \algbrac{ \algbrac{X}{\gamma} }{ Z }) \notag \\
    &= 2\gamma(X) k(Y,Z) - k(\algbrac{ \algbrac{X}{\gamma} }{ Y }, Z)
                         - k(Y, \algbrac{ \algbrac{X}{\gamma} }{ Z })
\end{align}
since $k(Y,Z)$ is a section of $\cL^2$.  To analyse the remaining algebraic bracket terms, we write $k = \ell \ltens \theta$ for some $\ell \in \s{0}{\cL}$ and $\theta \in \s{0}{\cL\tens\cB^*}$.  Then $k( \algbrac{ \algbrac{X}{\gamma} }{ Y }, Z ) = \ell \ltens \theta( \algbrac{ \algbrac{X}{\gamma} }{ Y }, Z )$, for which
\begin{align} \label{eq:mob2-pencils-killingeqn-2}
  \theta( \algbrac{ \algbrac{X}{\gamma} }{ Y }, Z )
    &=  \algbrac{ Z }
                { \algbrac{ \algbrac{ \algbrac{X}{\gamma} }
                                    { Y } }
                          { \theta } }
      \notag \\
    &=  \algbrac{ Z }
                { \algbrac{ \algbrac{ \algbrac{X}{\theta} }
                                    { \gamma } }
                          { Y } }
      + \algbrac{ Z }
                { \algbrac{ \algbrac{X}{\gamma} }
                          { \algbrac{Y}{\theta} } }
      \notag \\
    &=  \algbrac{ Z }
                { \algbrac{ \algbrac{ \theta(X,\bdot) }
                                    { \gamma } }
                          { Y } }
      + \algbrac{ Z }
                { \algbrac{ \algbrac{X}{\gamma} }
                          { \theta(Y,\bdot) } }
      \notag \\
    &=  \algbrac{ \algbrac{ \theta(X,Z) }
                          { \gamma } }
                { Y }
      + \algbrac{ \algbrac{ \theta(X,\bdot) }
                          { \algbrac{Z}{\gamma} } }
                { Y }
      \notag \\
    &\qquad
      + \algbrac{ \algbrac{ Z }
                          { \algbrac{X}{\gamma} } }
                { \theta(Y,\bdot) }
      + \algbrac{ \algbrac{X}{\gamma} }
                { \theta(Y,Z) }
      \notag \\
    &= -\algbrac{ \theta(X,Z) \ltens \gamma }
                { Y }
      - \algbrac{ \theta(X,Y) }
                { \algbrac{Z}{\gamma} }
      + \algbrac{ \theta(X,\bdot) }
                { \algbrac{ \algbrac{Z}{\gamma} }
                          { Y } }
      \notag \\
    &\qquad
      - \theta( Y, \algbrac{ \algbrac{X}{\gamma} }{ Z } )
      + \gamma(X) \theta(Y,Z)
      \notag \\
    &=  \gamma(Y) \theta(X,Z) + \gamma(Z) \theta(X,Y)
      - \theta( X, \algbrac{ \algbrac{Z}{\gamma} }{ Y } )
      \notag \\
    &\qquad
      - \theta( Y, \algbrac{ \algbrac{X}{\gamma} }{ Z } )
      + \gamma(X) \theta(Y,Z).
\end{align}
By symmetry of $\algbrac{ \algbrac{X}{\gamma} }{ Y }$ in $X,Y$ and \eqref{eq:mob2-pencils-killingeqn-1}, we have
\begin{align*}
  &\tfrac{1}{2} \weyld{\gamma} \left( (\D_X k)(Y,Z) + (\D_Y k)(Z,X) + (\D_Z k)(X,Y) \right) \\
    &\quad =
      \gamma(X) k(Y,Z) + \gamma(Y) k(Z,X) + \gamma(Z) k(X,Y) \\
    &\quad \qquad
      - k( \algbrac{ \algbrac{X}{\gamma} }{ Y }, Z )
      - k( \algbrac{ \algbrac{Y}{\gamma} }{ Z }, X )
      - k( \algbrac{ \algbrac{Z}{\gamma} }{ X }, Y ).
\end{align*}
Taking the tensor product of \eqref{eq:mob2-pencils-killingeqn-2} with $\ell$ and substituting the result on the \rhs, it is clear that all terms cancel.  Therefore \eqref{eq:mob2-pencils-killingeqn} is independent of $\D \in \Dspace$. \end{proof}

Calculating \wrt\ the \LC\ connection of $h$, \thref{lem:mob2-pencils-quadint} implies immediately that $k\in\s{0}{\cL^2\tens\cB^*}$ is a solution if and only if $(\pf{h})^{-2} k$, viewed as a section of $\Symm{2}T^*M$, is a homogeneous quadratic integral of $g \defeq (\pf h)^{-1} \ltens h^{-1}$.

\begin{defn} The equation \eqref{eq:mob2-pencils-killingeqn} is called the \emph{Killing $2$-tensor equation}, and its solutions $k\in\s{0}{\cL^2\tens\cB^*}$ are called \emph{Killing $2$-tensors}. \end{defn}

In \proj\ differential geometry, it is known \cite{mss2014-killingprol} that the Killing $2$-tensor equation is the first BGG operator associated to the representation $\bU^*$ from Subsection \ref{ss:mob2-pencils-adj}.  The same holds in \cproj\ geometry \cite[Cor.\ 5.8]{cemn2015-cproj}.  In these cases \eqref{eq:mob2-pencils-killingeqn} prolongs to a linear connection on the associated bundle $\cU^* \defeq \assocbdl{F^P}{P}{\bU^*}$, and the dimension of the solution space is bounded above by $\dim\bU^*$.  The Killing $2$-tensor equation is trivial for conformal geometries, since there $\bU^*$ is the trivial representation, and more naturally one studies the \emph{conformal Killing equation}; see \cite[\S 2.4]{gs2008-confkilling}.

\begin{prop} \thlabel{prop:mob2-pencils-killingmetric} The adjugate of any linear metric $h$ is a Killing $2$-tensor. \end{prop}

\begin{proof} If $h$ is \non degenerate\ then by \thref{cor:ppg-bgg-lcconn}, the \LC\ connection $\D^g$ of $h$ lies in $\Dspace$.  Since we may write $h^* = (\pf{h}) \ltens h^{-1}$, it follows that $h^*$ is $\D^g$-parallel and hence \eqref{eq:mob2-pencils-killingeqn} holds.  In the general case, note that the $1$-jet of $h^*$ is polynomial in the $1$-jet of $h$, so the result follows by continuity. \end{proof}

If $M$ admits a pencil $\bs{h}\in\s{0}{\cL^*\tens\cB\tens V^*}$ of linear metrics, then every linear metric $\bs{h}_v$ defines a Killing $2$-tensor $\bs{h}_v^*$.  Applying \thref{prop:mob2-pencils-adjpoly} to the metrisability pencil $\bs{h}$ we may view its adjugate $\bs{h}^*$ as a section of $\cL^2\tens\cB^* \tens \cO_V(n-1)$ over $M\times\pr{V}$, \ie\ as a homogeneous polynomial of degree $n-1$ in the parameters of the pencil.  Choosing an affine chart, the $1$-parameter family $h_t^*$ of Killing $2$-tensors is polynomial of degree $n-1$ in $t$, whose coefficients are also Killing $2$-tensors.  Trivialising $\cL$ \wrt\ the pfaffian $\pf{\bs{h}_u} \in\s{0}{\cL}$ of some metric $\bs{h}_u$ in the pencil, each Killing $2$-tensor $\bs{h}_v^*$ yields a quadratic integral of the geodesic flow of $\bs{g}_u = (\pf{\bs{h}_u})^{-1} \ltens \bs{h}_u^{-1}$.

Viewing $\bs{h}$ instead as a section of $\cL^*\tens\cB\tens\cO_V(1)$ over $M\times\pr{V}$ and choosing an affine chart with metrics $h,\b{h}$ at $\infty,0$ respectively, every metric in the pencil is proportional to a metric in the $1$-parameter family $h_t = \b{h}-th$.  This yields a $1$-parameter family $h_t^*$ of Killing $2$-tensors, and it is straightforward to check that we recover the quadratic integrals of the geodesic flow of $g=(\pf{h})^{-1} \ltens h^{-1}$ by trivialising $\cL$ \wrt\ the metric at $\infty$.  Indeed, in this trivialisation
\begin{equation} \label{eq:mob2-pencils-killingpencil} \begin{split}
  (\pf{h})^{-2} \ltens h_t^*
  &= (\pf{h_t}) \ltens (\pf{h})^{-1}\left[
    (\pf{h_t})^{-1} \ltens h^{-1}(A_t^{-1}\bdot,\bdot) \right] \\
  &= (\pf{h_t}) \ltens (\pf{h})^{-1} g(A_t^{-1}\bdot,\bdot) \\
  &= (\det A_t)^{1/r} g(A_t^{-1}\bdot,\bdot),
\end{split} \end{equation}
where $A_t\in\s{0}{\fp^0_M}$ denotes both the endomorphism satisfying $h_t=h(A_t\bdot,\bdot)$ and its transpose.  It is clear from \eqref{eq:mob2-pencils-killingpencil} that we recover Topalov and Matveev's integrals \cite{mt1998-trajectory, t2000-hierarchy, t2003-geodesiccomp} of the geodesic flow in the projective case,  and Kiyohara and Topalov's integrals \cite[Prop.\ 2.1]{kt2011-kahlerliouville} in the \cproj\ case.  It is well-known that the corresponding integrals of the geodesic flow mutually commute \wrt\ the Poisson bracket induced on $TM$ by any metric in the pencil.  The same result holds in general.

\begin{thm} \thlabel{thm:mob2-pencils-integralscomm} Let $\bs{h}$ be a metrisability pencil.  Then the components of $\bs{h}^*$ Poisson-commute \wrt\ the Poisson bracket on $TM$ induced by any metric in $\bs{h}$. \end{thm}

To prove \thref{thm:mob2-pencils-integralscomm}, we shall first need to describe the Poisson bracket on $TM$ induced by a linear metric in more detail.  We restrict attention to homogeneous polynomial functions on $T^*M$ of degree $k$, which may equivalently be viewed as sections of $\Symm{k}TM$ over $M$.  In this picture the canonical Poisson bracket becomes the \emph{Schouten--Nijenhuis bracket} on symmetric multivectors, given by%
\footnote{This may be more familiar in Penrose's abstract index notation \cite{pr1987-spinors1}.  Writing $Q,R$ with $k,\ell$ contravariant indicies, employing the summation convention and letting round brackets denote symmetrisation, the Schouten--Nijenhuis bracket reads $\poisson{Q}{R}^{a\cdots e} = k\, Q^{f(a\cdots b} \D_f R^{cd\cdots e)} - \ell\, \D_f^{(a\cdots bc} R^{d \cdots e)f}$.}
\vspace{0.2em}
\begin{equation*} \begin{gathered}
  \poisson{}{} : \s{0}{\Symm{k}TM} \times \s{0}{\Symm{\ell}TM}
            \to \s{0}{\Symm{k+\ell-1}TM}, \\
  \poisson{Q}{R}
    \defeq \sym{\left( k(\ve^i\intprod Q)\tens \D_{e_i}R
                - \ell(\ve^i\intprod R)\tens \D_{e_i}Q \right) }
\end{gathered}
\vspace{0.2em}
\end{equation*}
for any local frame $\{e_i\}_i$ of $TM$ with dual coframe $\{\ve^i\}_i$.  In particular, the Schouten--Nijenhuis bracket between homogeneous quadratic functions is given by
\vspace{0.2em}
\begin{equation} \label{eq:mob2-pencils-snquad}
  \poisson{Q}{R} = 2\sym{\left( Q(\ve^i,\bdot)\tens \D_{e_i}R
                              - R(\ve^i,\bdot)\tens \D_{e_i}Q \right)},
\vspace{0.2em}
\end{equation}
a section of $\Symm{3}TM$, where evaluating the symmetrisation of $Q(\ve^i,\bdot)\tens(\D_{e_i}R)$, a section of $TM\tens\Symm{2}TM$, on $\alpha,\beta,\gamma\in\s{1}{}$ entails taking the cyclic sum over $\alpha,\beta,\gamma$.

\begin{lem} \thlabel{lem:mob2-pencils-integralscomm} Let $h$ be a \non degenerate\ linear metric, with pfaffian $\pi\defeq \pf{h}$ and adjugate $h^* \defeq \adj{h}$.  Then
\vspace{0.2em}
\begin{equation} \label{eq:mob2-pencils-integralderiv}
  \D_X h^* = h^{-1}(Z^{\D},X)\, h^*
    + \pi \tens \algbrac{ h^{-1}(X,\bdot) }{ h^{-1}(Z^{\D},\bdot) }
\vspace{0.2em}
\end{equation}
for all $\D\in\Dspace$, where $Z^{\D}\in\s{0}{\cL^*\tens TM}$ satisfies $\D_X h = \algbrac{Z^{\D}}{X}$. \end{lem}

\begin{proof} By \thref{cor:ppg-bgg-lcconn} there is $\D^h\in \Dspace$ with $\D^h h=0$; therefore also $\D^h\pi = 0$ and $\D^h h^*=0$.  For any other $\D\in \Dspace$ we can write $\D = \D^h + \algbrac{}{\gamma}$ for some $\gamma\in\s{1}{}$, and it follows that $\D_X h = \algbrac{ \algbrac{X}{\gamma} }{ h } = \algbrac{ h(\gamma,\bdot) }{ X }$ by the Jacobi identity.  Applying $\liebdy$ then gives $\gamma = h^{-1}(Z^{\D},\bdot)$.

By \non degeneracy\ of $h$, we can write $h^* = \pi \ltens h^{-1}$.  Then since $\D h^* = \algbrac{}{\gamma} \acts h^*$, the Leibniz rule yields
\begin{align*}
  \D_X h^*
  &= \algbrac{ \algbrac{X}{\gamma} }{ \pi } \tens h^{-1}
    + \pi \tens \algbrac{ \algbrac{X}{\gamma} }{ h^{-1} } \\
  &= \gamma(X)\pi \ltens h^{-1} + \pi \ltens \left(
    \algbrac{\algbrac{X}{h^{-1}}}{\gamma} + \algbrac{X}{\algbrac{\gamma}{h^{-1}}} \right) \\
  &= \gamma(X)h^* + \pi \ltens \algbrac{h^{-1}(X,\bdot)}{\gamma},
\end{align*}
Substituting $\gamma=h^{-1}(Z^{\D},\bdot)$ from above gives the result. \end{proof}

\begin{proofof}{thm:mob2-pencils-integralscomm} Choose an affine chart for $\bs{h}$ with linear metrics $h,\b{h}$ at $\infty,0$ respectively.  It then suffices to consider the $1$-parameter family $h_t \defeq \b{h} - th$ of linear metrics, for which $\D_X h_t = \algbrac{Z_t^{\D}}{X}$ with $Z_t^{\D}=\b{Z}^{\D}-tZ^{\D}$ for all $t\in\bR$.

The $1$-parameter family $h_t^*$ of Killing $2$-tensors yields a $1$-parameter family $Q_t\in\s{0}{\Symm{2}TM}$ of quadratic integrals of $g \defeq (\pf{h})^{-1} \ltens h^{-1}$, which may be written as
\vspace{0.1em}
\begin{equation*}
  Q_t(\alpha,\beta) = h_t^*( h(\alpha,\bdot), h(\beta,\bdot) )
\vspace{0.1em}
\end{equation*}
for $\alpha,\beta\in\s{1}{}$.  We must show that the Schouten--Nijenhuis bracket \eqref{eq:mob2-pencils-snquad} vanishes for all $s,t\in\bR$.  Since $\poisson{Q_s}{Q_t}=0$ trivially when $s=t$, we assume that $s\neq t$. It suffices also to assume that $h_s,h_t$ are \non degenerate; the general case will then follow by continuity.

Let $\D \in \Dspace$ be the \LC\ connection of $h$.  Then since $\D h=0$,
\vspace{0.1em}
\begin{equation*}
  (\D_X Q_t)\left( \alpha,\beta \right)
    = (\D_X h_t^*)\left( h(\alpha),h(\beta) \right)
\vspace{0.1em}
\end{equation*}
for all $\alpha,\beta\in\s{1}{}$, where we write $h(\alpha) \defeq h(\alpha,\bdot)$ for notational convenience.  Since $h$ is \non degenerate\ and $\pi \defeq \pf{h}$ is $\D$-parallel, it suffices to take $1$-forms of the form $\pi^{-1} \ltens h^{-1}(X,\bdot)$ and prove that
\vspace{0.1em}
\begin{equation} \label{eq:mob2-pencils-integralscomm-1}
  \sym{\left[ h_s^*(h(\ve^i),\bdot) \tens \D_{e_i}h_t^*
            - h_t^*(h(\ve^i),\bdot) \tens \D_{e_i}h_s^* \right]} = 0
\vspace{0.1em}
\end{equation}
as a section of $\cL^3\tens\Symm{3}T^*M$.  Notice that we may write
\begin{equation*} \begin{aligned}
  h &= \tfrac{1}{s-t} \big( (\b{h}-th) - (\b{h}-sh) \big) \\
    &= \tfrac{1}{s-t}( h_t - h_s )
\end{aligned} \end{equation*}
for any $s\neq t$, so that \eqref{eq:mob2-pencils-integralscomm-1} holds if and only if
\vspace{0.05em}
\begin{equation} \label{eq:mob2-pencils-integralscomm-2}
  \sym{\left[
  \begin{aligned}
    &   h_s^*(h_t(\ve^i)\bdot) \tens \D_{e_i}h_t^*
      - h_s^*(h_s(\ve^i)\bdot) \tens \D_{e_i}h_t^* \\
    & \quad - h_t^*(h_t(\ve^i)\bdot) \tens \D_{e_i}h_s^*
            + h_t^*(h_s(\ve^i)\bdot) \tens \D_{e_i}h_s^*
  \end{aligned}
  \right]} = 0.
\vspace{0.05em}
\end{equation}
Let us analyse the terms inside the symmetrisation one at a time.  Writing $\pi(t) \defeq \pf{h_t}$, for the third term above we have
\begin{align*}
  h_t^*(h_t(\ve^i,\bdot),X)\, (\D_{e_i}h_s^*)(Y,Z)
  &= \pi(t) h_t^{-1}(h_t(\ve^i),X)\, (\D_{e_i}h_s^*)(Y,Z) \\
  &= \pi(t) \ve^i(X)\, (\D_{e_i}h_s^*)(Y,Z) \\
  &= \pi(t) \, (\D_X h_s^*)(Y,Z),
\end{align*}
whose symmetrisation in $X,Y,Z$ vanishes because $h_s^*$ is a Killing $2$-tensor.  Switching the roles of $s,t$, the second term in \eqref{eq:mob2-pencils-integralscomm-2} also vanishes.

Observe now that, since we are computing \wrt\ the \LC\ connection $\D$ of $h$, $Z^{\D}=0$ and therefore $Z^{\D}_t = \b{Z}^{\D}-tZ^{\D} = \b{Z}^{\D}$ for all $t\in\bR$.  By Lemma \ref{lem:mob2-pencils-integralscomm} it follows that the first term in \eqref{eq:mob2-pencils-integralscomm-2} is given by
\begin{align*}
  & h_s^*(h_t(\ve^i,\bdot),X)\tens \D_{e_i}h_t^* \\
  &\quad = \pi(s) \ve^i\left( h_t(h_s^*(X,\bdot),\bdot) \right) \big(
    h_t^{-1}(\b{Z}^{\D},e_i)h_t^* + \pi(t)\tens
    \algbrac{h_t^{-1}(e_i,\bdot)}{h_t^{-1}(\b{Z}^{\D},\bdot)} \big) \\
  &\quad = \pi(s)\pi(t) \big( h_s^{-1}(\b{Z}^{\D},X)h_t^{-1}
    + \algbrac{h_s^{-1}(X,\bdot)}{h_t^{-1}(\b{Z}^{\D},\bdot)} \big).
\end{align*}
We note that $h_s^{-1}(X,\bdot)\in\s{0}{\cL\tens T^*M}$ and $h_t^{-1}(Z_t^{\D},\bdot)\in\s{1}{}$, so that the algebraic bracket term is a section of $\cL\tens\cB^*$ as required.  Writing $\alpha_s \defeq h_s^{-1}(\b{Z}^{\D},\bdot) \in \s{1}{}$, the previous two calculations imply that \eqref{eq:mob2-pencils-integralscomm-2} is equivalent to
\begin{equation} \label{eq:mob2-pencils-integralscomm-3}
\vspace{0.1em}
  \pi(s)\pi(t) \sym[X,Y,Z]{ \! \left[
  \begin{aligned}
    & \alpha_s(X)h_t^{-1}(Y,Z) + \algbrac{h_s^{-1}(X,\bdot)}{\alpha_t}(Y,Z) \\
    & \quad + \alpha_t(X)h_s^{-1}(Y,Z) + \algbrac{h_t^{-1}(X,\bdot)}{\alpha_s}(Y,Z) \,
  \end{aligned}
    \right]} = 0,
\vspace{0.1em}
\end{equation}
where by $\sym[X,Y,Z]{}$ we mean the cyclic sum over $X,Y,Z$.  Since $\pi(t) \in \s{0}{\cL}$, we have
\begin{equation} \label{eq:mob2-pencils-integralscomm-4}
  \pi(t) \alpha_s(X) h_t^{-1}(Y,Z)
    = \algbrac{\algbrac{X}{\alpha_s}}{\pi(t)}\tens h_t^{-1}(Y,Z),
\end{equation}
a section of $\cL^3$.  The Jacobi identity then yields
\begin{equation*}
  \algbrac{ h_t^{-1}(X,\bdot) }{ \alpha_s }
  = \algbrac{ \algbrac{X}{h_t^{-1}} }{ \alpha_s }
  = \algbrac{ \algbrac{X}{\alpha_s} }{ h_t^{-1} }
\end{equation*}
since the bracket between $T^*M$ and $\cL\tens\cB^*$ is zero.  Therefore, combining this with \eqref{eq:mob2-pencils-integralscomm-4} and omitting the evaluation on vector fields $Y,Z$, the first and last terms inside the symmetrisation in \eqref{eq:mob2-pencils-integralscomm-3} equal
\begin{align} \label{eq:mob2-pencils-integralscomm-5}
  &\pi(s)\tens\left(
    \algbrac{\algbrac{X}{\alpha_s}}{\pi(t)} \tens h_t^{-1}
    + \pi(t) \tens \algbrac{\algbrac{X}{\alpha_s}}{h_t^{-1}}
    \right) \notag \\
  &\quad = \pi(s) \tens
    \left( \algbrac{X}{\alpha_s} \acts \pi(t) \ltens h_t^{-1} \right) \notag \\
  &\quad = \pi(s) \tens \left( \weyld{\alpha_s} \D_X h_t^* \right),
\end{align}
where $\weyld{\alpha_s}$ is the linear variation of Weyl structure \wrt\ $\alpha_s$.  The symmetrisation of \eqref{eq:mob2-pencils-integralscomm-5} is zero by \thref{prop:mob2-pencils-killingeqn}.  Therefore, upon symmetrisation, the first and last terms in \eqref{eq:mob2-pencils-integralscomm-3} sum to zero; exchanging $s,t$, so do the middle two terms.  Since \eqref{eq:mob2-pencils-integralscomm-3} was equivalent to commutativity of the $h_t^*$, we are done. \end{proofof}

\section{Canonical vector fields} 
\label{s:mob2-vfs}

In \cproj\ geometry, the role of the canonically defined Killing vector fields associated to a metrisability pencil has been emphasised by a number of authors, most notably in \cite{acg2006-ham2forms1,acgt2004-ham2forms2,bmr2015-localcproj,mr2012-yanoobata}.  These Killing fields are defined as symplectic gradients $J\grad[g]{\sigma_i}$ for certain functions $\sigma_i$ associated to the pencil, and mutually commute.  More recently, analogous gradient fields $\grad[g]{\sigma_i}$ have been discovered in \proj\ geometry by Eastwood \cite{e2014-commutingvfs} and developed further in the \cproj\ setting in \cite{bmr2015-localcproj, cemn2015-cproj}.  While these gradient fields mutually commute, they are not Killing in general.  The present author believes that these gradient fields $X(t)$ play the more fundamental role, and it is merely a happy accident that the symplectic gradients of \cproj\ geometry are Killing fields.  Indeed, in \qtn ic\ geometry the vector fields $J_aX(t)$ are not Killing unless the manifold is hypercomplex.

In this section, we develop the theory of the canonically defined gradient vector fields associated to a pencil.  Subsection \ref{ss:mob2-vfs-pairing} is devoted to an invariant description of these vector fields, using the bilinear pairing provided by the general BGG calculus, while Subsection \ref{ss:mob2-vfs-comm} undertakes a somewhat gruesome proof that they mutually commute.

\medskip
\subsection{Definition using BGG pairings} 
\label{ss:mob2-vfs-pairing}

By \thref{prop:ppg-alg-gr}, an algebraic Weyl structure for $\fp$ induces decompositions
\begin{equation} \label{eq:mob2-vfs-Wgradings}
\begin{aligned}
  \bW   &\isom (L^*\tens B) \dsum (L^*\tens \fg/\fp) \dsum L^* \\
  \bW^* &\isom L \dsum (L\tens \fp^{\perp}) \dsum (L\tens B^*)
\end{aligned}
\end{equation}
of $\bW$ and $\bW^*$, where $B \leq \Symm{2}(\fg/\fp)$ is a $\fp^0$-subrepresentation and, by \thref{prop:ppg-alg-Lwedge}, $L \isom (\Wedge{rn} \fg/\fp)^{2/r(n+1)}$.  Since also $\fg \isom \fg/\fp \dsum \fp^0 \dsum \fp^{\perp}$, taking Lie brackets of appropriate elements in \eqref{eq:mob2-vfs-Wgradings} yields a bilinear pairing
\begin{equation} \label{eq:mob2-vfs-hompairing}
  \bW\times\bW^* \to \fg.
\end{equation}
The $(\fg/\fp)$-valued component of \eqref{eq:mob2-vfs-hompairing} is $\fp$-invariant, so descends to a pairing $(L^*\tens B)\times L \to \fg/\fp$ on zeroth homology modules.  The theory of Section \ref{ss:para-bgg-bgg} then provides a bilinear differential pairing $(\cL^*\tens \cB) \times \cL \to TM$ between associated bundles of the form \eqref{eq:para-calc-bggcup}.  In terms of the BGG splitting operators \eqref{eq:ppg-bgg-metricsplit} and \eqref{eq:ppg-bgg-hesssplit} this is
\begin{equation} \label{eq:mob2-vfs-pairing}
  (h,\pi) \mapsto X(h,\pi) \defeq h(\eta^{\D},\bdot) - \pi \ltens Z^{\D}
\end{equation}
for any $\D \in \Dspace$, where $\eta^{\D} \defeq \D\pi$ and $Z^{\D} \defeq \tfrac{2}{rn-r+2} \liebdy (\D h)$ as in \eqref{eq:ppg-bgg-metricprol} and \eqref{eq:ppg-bgg-hessprol}.

\begin{lem} \thlabel{lem:mob2-vfs-Xinv} $X(h,\pi)$ is independent of $\D \in \Dspace$. \end{lem}

\begin{proof} From Section \ref{s:ppg-bgg} we have $\weyld{\gamma} Z^{\D} = h(\gamma,\bdot)$ and $\weyld{\gamma} \eta^{\D} = \algbrac{ \algbrac{}{\gamma} }{ \pi } = \pi \ltens \gamma$, so that $\weyld{\gamma} X(h,\pi) = h(\pi \ltens \gamma, \bdot) - \pi \ltens h(\gamma, \bdot) = 0$. \end{proof}

If $h$ is \non degenerate\ then $X(h,\pi)$ is a gradient vector field \wrt\ $g \defeq (\pf{h})^{-1} \ltens h^{-1}$: we may trivialise the line bundle $\cL$ \wrt\ $\pf{h}$ and write $X(h,\pi)$ \wrt\ the \LC\ connection $\D$ of $h$, yielding
\vspace{-0.3em}
\begin{equation} \label{eq:mob2-vfs-gradient}
  X(h,\pi) = h(\D\pi,\bdot) = g^{-1}( \d((\pf{h})^{-1}\pi), \bdot )
    = \grad[g]( (\pf{h})^{-1}\pi ).
\vspace{-0.3em}
\end{equation}
Since this simple fact shall be important in the sequel, we record the following property of gradient vector fields.

\begin{lem} \thlabel{lem:mob2-vfs-gradient} Let $(M,g)$ be a (\pseudo)\riem\ manifold with \LC\ connection $\D$.  Then $\D X$ is \self adjoint\ \wrt\ $g$ for all gradient vector fields $X$. \end{lem}

\begin{proof} Suppose that $X = \grad[g]{f}$.  Then for all vector fields $Y,Z$, we have $g(\D_Y X,Z) = g(\D_Y \d f^{\sharp},Z) = (\D_Y \d f)(Z) = \p_Y \p_Z f - \d f(\D_Y Z)$.  Since $\D$ is \torsionfree, we obtain $g(\D_Y X,Z) - g(Y, \D_Z X) = \p_{\liebrac{Y}{Z}}f - \d f(\liebrac{Y}{Z}) = 0$ as required. \end{proof}

Observe also that $X(h,\pf{h})=0$ for all $h\in\s{0}{\cL^*\tens B}$: either $h$ is degenerate, in which case $\pf{h}=0$ and $\eta^{\D}=0$; or else $h$ is \non degenerate\ and we may evaluate $X(h,\pf{h})$ \wrt\ the \LC\ connection $\D$ of $h$, giving $\eta^{\D} = 0$ and $Z^{\D} = 0$.

Suppose now that we have a pencil $\bs{h} \in\s{0}{\cL^*\tens\cB\tens V^*}$ of compatible linear metrics.  The pfaffian $\bs{\pi} \defeq \pf{\bs{h}}$ is a section of $\cL \tens \Symm{n}V^*$, so may be viewed as a homogeneous polynomial of degree $n$ in the parameters of the pencil.  We would like to define the vector fields $X(\bs{h},\bs{\pi})$, but since $X(h,\pf{h})=0$ some care is needed to avoid making a trivial definition.  There are two equivalent approaches:
\begin{slimitemize}
  \item We may choose an affine chart on $\pr{V}$ with linear metrics $h,\b{h}$ at $\infty,0$ respectively.  Then each metric in the pencil is proportional to some $h_t \defeq \b{h}-th$ and $\pi(t) \defeq \pf{h_t}$ is a polynomial of degree $n$ in $t$ by \thref{cor:mob2-pencils-pfpoly}.  Since \eqref{eq:mob2-vfs-pairing} is linear in $h$,
  \vspace{-0.3em}
  \begin{equation} \label{eq:mob2-vfs-chart} \begin{split}
    X(h_s,\pi(t))
    &= X(h_s,\pi(t)) - X(h_t,\pi(t)) \\
    &= X(\b{h}-sh,\pi(t)) - X(\b{h}-th,\pi(t)) = (t-s)X(h,\pi(t))
  \end{split}
  \vspace{-0.3em}
  \end{equation}
  for each $s,t\in\bR$, where $X(t) \defeq X(h,\pi(t))$ is independent of $s$ and polynomial of degree $n-1$ in $t$.
  
  \item More invariantly, we view $\bs{h}$ and $\bs{\pi}$ respectively as sections of $\cL^*\tens\cB\tens V^*$ and $\cL \tens \Symm{n}V^*$ over $M$, so that $X(\bs{h},\bs{\pi})$ is a section of $TM \tens V^* \tens \Symm{n}V^*$.  The Clebsch--Gordan formula \cite{fh1991-repntheory} gives
  \vspace{-0.3em}
  \begin{equation*}
    V^* \tens \Symm{n}V^* \isom \Symm{n+1}V^* \dsum (\Symm{n-1}V^* \tens \Wedge{2}V^*)
  \vspace{-0.3em}
  \end{equation*}
  and, since $X(h,\pf{h}) = 0$ for each $h\in\s{0}{\cL^*\tens\cB}$, the $\Symm{n+1}V^*$-component must vanish identically.  It follows that we may write $X(\bs{h},\bs{\pi}) = \bs{X} \symm \ve$ for a section $\bs{X}$ of $TM \tens \Symm{n-1}V^*$ and a constant area form $\ve \in \Wedge{2}V^*$.  Of course, $\bs{X}$ is identified with $X(t)$ from above in an affine chart, while the evaluation of $\ve$ on $(1,s),(1,t)\in V$ is (proportional to) the coefficient $(t-s)$ in \eqref{eq:mob2-vfs-chart}.
\end{slimitemize}

The pencil $\bs{h}$ therefore induces a homogeneous family $\bs{X}$ of gradient vector fields.  For concreteness, we shall often work with the vector fields $X(t)$ in an affine chart.

\begin{defn} \label{defn:mob2-vfs-vfs} The homogeneous family $\bs{X} \in\s{0}{TM\tens\Symm{n-1}V^*}$ will be called the \emph{family of canonical vector fields} associated to the pencil $\bs{h}$. \end{defn}

\begin{rmk} In an affine chart, $X(t) = \sum{i=1}{n-1} (-1)^i t^i X_i$ is polynomial of degree $n-1$ in $t$.  Identifying $X(t)$ with the gradient field $\grad[g]{( (\pf{h})^{-1}\pi(t) )}$, we have $X_i = \grad[g]{\sigma_i}$ for $\sigma_i$ the elementary symmetric function of the $n$ roots of $\pi(t)$; these are Eastwood's gradient fields in \proj\ geometry \cite{e2014-commutingvfs}.  In \cproj\ geometry, since $J$ is parallel and $X_i$ is gradient, the symplectic gradients $K_i \defeq J\grad[g]{\sigma_i}$ are Killing fields. \end{rmk}

\begin{prop} \thlabel{prop:mob2-vfs-DXt} Choose an affine chart for $\pr{V}$ with \non degenerate\ linear metrics $h,\b{h}$ at $\infty, 0$ respectively, and let $\D^g$ be the \LC\ connection of $g \defeq (\pf h)^{-1/r} \ltens h^{-1}$.  Then $\D^g X(t)$, is a $g$-\self adjoint\ section of $\fp^0_M$ for all $t \in \bR$. \end{prop}

\begin{proof} The fact that $\D X(t)$ is \self adjoint\ \wrt\ $g$ follows from \eqref{eq:mob2-vfs-gradient} and \thref{lem:mob2-vfs-gradient}.  Differentiating \eqref{eq:mob2-vfs-pairing} using the prolongation \eqref{eq:ppg-bgg-metricprol} gives
\begin{align*}
  \D_Y X(t)
    &=  \algbrac{Z^{\D}}{Y}(\D\pi(t), \bdot) + h(\D_Y \D\pi(t), \bdot) \\
    &\qquad - (\D_Y Z^{\D}) \ltens \pi(t) - (\D_Y \pi(t)) \ltens Z^{\D} \\
    &=  \algbrac{ \algbrac{Z^{\D}}{\D\pi(t)} }{ Y }
      + \algbrac{Z^{\D}}{\D_Y \pi(t)}
      + h(\D_Y \D\pi(t) - \pi(t) \ltens \nRic{Y}{}, \bdot) \\
    &\qquad - (\D_Y \pi(t)) \ltens Z^{\D}
      - \lambda^{\D} \ltens \pi(t) \ltens Y
      - \pi(t) \ltens \quab^{-1} \liebdy \algbrac{\Weyl{}{}}{h}_Y \\
    &= \algbrac{ \algbrac{Z^{\D}}{\D\pi(t)} }{ Y }
      + h(\D_Y \D\pi(t) - \pi(t) \ltens \nRic{Y}{}, \bdot) \\
    &\qquad - \lambda^{\D} \ltens \pi(t) \ltens Y
      - \pi(t) \ltens \quab^{-1} \liebdy \algbrac{\Weyl{}{}}{h}_Y
\end{align*}
for all vector fields $Y$.  The bracket $\algbrac{Z^{\D}}{\D\pi(t)}$ lies in $\fp^0_M \dsum \liecenter{\fq^0}_M$; however $\alg{gl}{TM}$ has zero intersection with $\liecenter{\fq^0}_M$, so $\algbrac{Z^{\D}}{\D\pi(t)}$ is a section of $\fp^0_M$.  Since $h_t$ is a linear metric for each $t\in\bR$, its pfaffian $\pi(t)$ defines a solution of the hessian equation by \thref{cor:ppg-bgg-metrichess}.  Thus $\D^2\pi(t) - \pi(t) \ltens \nRic{}{}$ is a section of $\cL\tens\cB^*$, implying that $Y \mapsto h(\D_Y \D\pi(t) - \pi(t) \ltens \nRic{Y}{}, \bdot)$ is a section of $\fp^0_M$ by \thref{lem:ppg-bgg-p0B}.  Since $\lambda^{\D} \ltens \pi(t)$ is just a smooth function, the third term is proportional to the identity map and hence a section of $\fp^0_M$.  Finally, the Weyl curvature term defines a section of $\fp^0_M$ by \thref{lem:ppg-bgg-weylp0}.  Therefore the \rhs\ above is the image of $Y$ under a section of $\fp^0_M$. \end{proof}

By \thref{prop:mob2-pencils-Zgradient}, the vector field $\pi \ltens \b{Z}^{\D}$ may be identified with the vector field $\Lambda \defeq \mu^{\sharp}$ appearing in the main equation \eqref{eq:cproj-class-maineqn} of \cproj\ geometry.  Thus \thref{prop:mob2-vfs-DXt} generalises the well-known fact that $\Lambda$ is holomorphic; see \cite{cmr2015-cprojmob, mr2012-yanoobata}.  In \qtn ic\ geometry, one finds that in fact $\D X(t)$ is \qtn-linear, implying that $X(t)$ is a gradient \qtn ic\ vector field by \cite[Lem.\ 1(2)]{am1994-qtnlaplacian}.

\subsection{Commutativity} 
\label{ss:mob2-vfs-comm}

The objective of this subsection is to prove that the canonical vector fields mutually commute.  This is well-known in \cproj\ geometry \cite{acg2006-ham2forms1, acgt2004-ham2forms2, bmr2015-localcproj, mr2012-yanoobata}, where one usually works with the canonical Killing fields $K(t) \defeq JX(t)$.

\vspace{-0.15em}
\begin{thm} \thlabel{thm:mob2-vfs-vfs} The components of $X(\bs{h},\bs{\pi})$ mutually commute.
\vspace{-0.15em}
\end{thm}

The proof of \thref{thm:mob2-vfs-vfs} will require some preliminary technical work.  Note that is suffices to proof that the components of $\bs{X} \in \s{0}{TM\tens\Symm{n-1}V^*}$ from above commute; equivalently, we may choose an affine chart on $\pr{V}$ with metrics $h,\b{h}$ at $\infty,0$ respectively and show that $\liebrac{X(s)}{X(t)}=0$ for all $s,t\in\bR$.

\vspace{-0.2em}
\begin{lem} \thlabel{lem:mob2-vfs-bracs} Let $g$ be a \non degenerate\ section of $\cB^*$ and suppose that $A \in \s{0}{\fp^0_M}$ is \self adjoint\ \wrt\ $g$.  Then for all vector fields $X$, we have:
\begin{enumerate}
  \item \label{lem:mob2-vfs-bracs-1}
  $\algbrac{ \algbrac{X}{X^\flat} }{ X^{\flat} } = -g(X,X)X^{\flat}$; and
      
  \item \label{lem:mob2-vfs-bracs-2}
  $\algbrac{ \algbrac{AX}{X^{\flat}} }{ X^{\flat} } = -g(AX,X)X^{\flat}$,
\end{enumerate}
where $\flat : TM \to T^*M$ is the musical isomorphism $X \mapsto g(X,\bdot)$ of $g$.
\vspace{-0.2em}
\end{lem}

\begin{proof} \proofref{lem:mob2-vfs-bracs}{1}
Let $h \defeq (\det g)^{1/r(n+1)} \ltens g^{-1}$ be the corresponding section of $\cL^*\tens\cB$.  For an arbitrary vector field $Y$, we may write
\begin{align} \label{eq:mob2-vfs-bracs-1}
  \algbrac{ \algbrac{X}{X^{\flat}} }
                 { X^{\flat} } (Y)
    &= h \big( \algbrac{ \algbrac{X}{X^{\flat}} }
                       { X^{\flat} }, h^{-1}(Y,\bdot) \big) \notag \\
    &= \killing{ h^{-1}(Y,\bdot) }
               { h( \algbrac{ \algbrac{X}{X^{\flat}}}
                            { X^{\flat} }, \bdot) } \notag \\
    &= \killing{ \algbrac{ Y }{ h^{-1} } }
               { \algbrac{ h }
                         { \algbrac{ \algbrac{X}{X^{\flat}} }
                                   { X^{\flat} } } } \notag \\
    &= \killing{ Y }
               { \algbrac{ h^{-1} }
                           { \algbrac{ h }
                                     { \algbrac{ \algbrac{X}{X^{\flat}} }
                                               { X^{\flat} } } } }
\end{align}
by the Jacobi identity and invariance of the Killing form on $\fh_M$.  Using Table \ref{tbl:ppg-alg-Z2gr}, successive applications of the Jacobi identity to last display yield
\begin{align} \label{eq:mob2-vfs-bracs-2}
  &\algbrac{ h^{-1} }
           { \algbrac{ h }
                     { \algbrac{ \algbrac{X}{X^{\flat}} }
                               { X^{\flat} } } }
    \notag \\ &\quad
  = \algbrac{ h^{-1} }
            { \algbrac{ \algbrac{ X }
                                { \algbrac{h}{X^{\flat}} } }
                      { X^{\flat} } }
    + \algbrac{ h^{-1} }
              { \algbrac{ \algbrac{X}{X^{\flat}} }
                        { \algbrac{h}{X^{\flat}} } }
    \notag \\ &\quad
  = \algbrac{ h^{-1} }
            { \algbrac{ \algbrac{ X }
                                { h(X^{\flat},\bdot) } }
                      { X^{\flat} } }
    + \algbrac{ h^{-1} }
              { \algbrac{ \algbrac{X}{X^{\flat}} }
                        { h(X^{\flat},\bdot) } }
    \notag \\ &\quad
  = \algbrac{ \algbrac{ \algbrac{h^{-1}}{X} }
                      { h(X^{\flat},\bdot) } }
            { X^{\flat} }
    + \algbrac{ \algbrac{X}{ \algbrac{ h^{-1} }
                                     { h(X^{\flat},\bdot) } } }
              { X^{\flat} }
    \notag \\ &\qquad\qquad
    + \algbrac{ \algbrac{ \algbrac{h^{-1}}{X} }
                        { X^{\flat} } }
              { h(X^{\flat},\bdot) }
    + \algbrac{ \algbrac{X}{X^{\flat}} }
              { \algbrac{ h^{-1} }
                        { h(X^{\flat},\bdot) } }
    \notag \\ &\quad\!
  \begin{aligned}
  &= -\algbrac{ \algbrac{ h^{-1}(X,\bdot) }
                       { h(X^{\flat},\bdot) } }
             { X^{\flat} }
    - \algbrac{ \algbrac{ h^{-1}(X,\bdot) }
                        { X^{\flat} } }
              { h(X^{\flat},\bdot) }
    \\ &\qquad\quad
    + 2\algbrac{ \algbrac{X}{X^{\flat}} }
               { X^{\flat} }.
  \end{aligned}
\end{align}
Writing $\pi \defeq \pf{h} \neq 0$, we have $h(X^{\flat},\bdot) = \pi^{-1} \ltens X$ and $h^{-1}(X,\bdot) = \pi X^{\flat}$.  Making this substitution in \eqref{eq:mob2-vfs-bracs-2} and using the Jacobi identity on the second term, we obtain
\begin{align} \label{eq:mob2-vfs-bracs-3}
  &\algbrac{ h^{-1} }
          { \algbrac{ h }
                    { \algbrac{ \algbrac{X}{X^{\flat}} }
                              { X^{\flat} } } }
    \notag \\ &\quad
  = -\algbrac{ \algbrac{ \pi X^{\flat} }
                        { \pi^{-1}X } }
              { X^{\flat} }
    - \algbrac{ \algbrac{ \pi X^{\flat} }
                        { X^{\flat} } }
              { \pi^{-1}X }
   + 2\algbrac{ \algbrac{X}{X^{\flat}} }
              { X^{\flat} }
    \notag \\ &\quad
  = -2\algbrac{ \algbrac{ \pi X^{\flat} }
                        { \pi^{-1}X } }
              { X^{\flat} }
    - \algbrac{ \pi X^{\flat} }
              { \algbrac{ X^{\flat} }
                        { \pi^{-1}X } }
   + 2\algbrac{ \algbrac{X}{X^{\flat}} }
              { X^{\flat} }.
\end{align}
Consider the first term on the \rhs\ of \eqref{eq:mob2-vfs-bracs-3}.  Since we may write $\pi^{-1} \ltens X = \algbrac{\pi^{-1}}{X}$, the Jacobi identity yields
\begin{align} \label{eq:mob2-vfs-bracs-4}
  \algbrac{ \algbrac{ \pi X^{\flat} }
                    { \algbrac{ \pi^{-1} }
                              { X } } }
           { X^{\flat} }
    \notag
  &= \algbrac{ \algbrac{ \algbrac{ \pi X^{\flat} }
                                 { \pi^{-1} } }
                       { X } }
             { X^{\flat} }
   + \algbrac{ \algbrac{ \pi^{-1} }
                       { \algbrac{ \pi X^{\flat} }
                                 { X } } }
             { X^{\flat} }
    \notag \\
  &= \algbrac{ \algbrac{ X^{\flat} }
                       { X } }
             { X^{\flat} }
   - \algbrac{ \algbrac{ \pi^{-1} }
                       { \pi g(X,X) } }
             { X^{\flat} },
\end{align}
where we have used that the inner-most brackets are the contractions $\algbrac{\pi X^{\flat}}{\pi^{-1}} = X^{\flat}$ and $\algbrac{\pi X^{\flat}}{X} = -\pi g(X,X)$.  The bracket between $\pi^{-1}\in\s{0}{\cL^*}$ and $X^{\flat}\in\s{1}{}$ vanishes, so applying the Jacobi identity to the second term in \eqref{eq:mob2-vfs-bracs-4} gives
\begin{align} \label{eq:mob2-vfs-bracs-5}
  \algbrac{ \algbrac{ \pi X^{\flat} }
                    { \algbrac{ \pi^{-1} }
                              { X } } }
           { X^{\flat} }
    \notag
  &= -\algbrac{ \algbrac{X}{X^{\flat}} }
              { X^{\flat} }
    - \algbrac{ \pi^{-1} }
              { \algbrac{ \pi g(X,X) }
                        { X^{\flat} } }
    \notag \\
  &= -\algbrac{ \algbrac{X}{X^{\flat}} }
              { X^{\flat} }
    + \algbrac{ \pi^{-1} }
              { \pi g(X,X)X^{\flat} }
    \notag \\
  &= -\algbrac{ \algbrac{X}{X^{\flat}} }
              { X^{\flat} }
    - g(X,X)X^{\flat}.
\end{align}
Consider now the second term in \eqref{eq:mob2-vfs-bracs-3}.  Since the inner algebraic bracket is (minus) the contraction, we have
\begin{equation} \label{eq:mob2-vfs-bracs-6}
  -\algbrac{ \pi X^{\flat} }
           { \algbrac{ X^{\flat} }
                     { \pi^{-1}X } }
    = \algbrac{ \pi X^{\flat} }
              { \pi^{-1} g(X,X) }
    = g(X,X)X^{\flat}.
\end{equation}
Now substituting \eqref{eq:mob2-vfs-bracs-4}, \eqref{eq:mob2-vfs-bracs-5} and \eqref{eq:mob2-vfs-bracs-6} into \eqref{eq:mob2-vfs-bracs-1}, we obtain
\begin{equation*}
  \killing{ Y }
          { \algbrac{ \algbrac{X}{X^{\flat}} }
                    { X^{\flat} } }
  = \killing{ Y }
            { 4\algbrac{ \algbrac{X}{X^{\flat}} }
                       { X^{\flat} }
              + 3g(X,X)X^{\flat} }.
\end{equation*}
and therefore
\begin{equation*}
  \killing{ Y }
          { \algbrac{ \algbrac{X}{X^{\flat}} }
                    { X^{\flat} }
          + g(X,X)X^{\flat} }
    = 0
\end{equation*}
by rearranging.  Since this holds for all vector fields $Y$, the result follows by \non degeneracy\ of the Killing form.

\smallskip

\proofref{lem:mob2-vfs-bracs}{2}
Following the steps of \ref{lem:mob2-vfs-bracs-1} with $\algbrac{\algbrac{X}{X^{\flat}}}{X^{\flat}}$ replaced by $\algbrac{\algbrac{AX}{X^{\flat}}}{X^{\flat}}$ until we reach equations analogous to \eqref{eq:mob2-vfs-bracs-4} and \eqref{eq:mob2-vfs-bracs-5}, we find that
\begin{align} \label{eq:mob2-vfs-bracs-7}
  \killing{ Y }
          { \algbrac{ \algbrac{AX}{X^{\flat}} }
                    { X^{\flat} } }
  &= \killing{ Y }
             { -2\algbrac{ \algbrac{(AX)^{\flat}}{X} }
                         { X^{\flat} }
               -2\algbrac{ \algbrac{ \pi^{-1} }
                                   { \algbrac{ \pi(AX)^{\flat} }
                                             { X^{\flat} } } }
                         { X^{\flat} } } \notag \\
  & \begin{aligned}
      = \big\langle
        \vphantom{\algbrac{\algbrac{(AX)^{\flat}}{X}}{X^{\flat}}} 
          Y ,& -2\algbrac{ \algbrac{(AX)^{\flat}}{X} }
                                 { X^{\flat} }
                      + 2g(AX,X)X^{\flat}  \\[-0.2em]
       &\qquad + g(X,X)(AX)^{\flat}
                      + 2\algbrac{ \algbrac{AX}{X^{\flat}} }
                                 { X^{\flat} }
       \big\rangle.
     \end{aligned}
\end{align}
Consider the first term above.  Since $A$ is \self adjoint\ \wrt\ $g$, we have $(AX)^{\flat}(Y) = g(X,AY) = (X^{\flat}\circ A)(Y) = \algbrac{X^{\flat}}{A}(Y)$, and thus $(AX)^{\flat} = \algbrac{X^{\flat}}{A}$.  By the Jacobi identity, the first term on the \rhs\ above equals
\begin{align} \label{eq:mob2-vfs-bracs-8}
  \algbrac{ \algbrac{(AX)^{\flat}}{X} }
          {X^{\flat} }
  &= \algbrac{ \algbrac{ \algbrac{X^{\flat}}{A} }
                       { X } }
             { X^{\flat} }
    \notag \\
  &= \algbrac{ \algbrac{ \algbrac{X^{\flat}}{X} }
                       { A } }
             { X^{\flat} }
   + \algbrac{ \algbrac{X^{\flat}}{AX} }
             { X^{\flat} }
    \notag \\
  &\! \begin{aligned}
    &= \algbrac{ \algbrac{ \algbrac{X^{\flat}}{X} }
                         { X^{\flat} } }
               { A }
     - \algbrac{ \algbrac{X^{\flat}}{X} }
               { (AX)^{\flat} } \\
    &\qquad
     + \algbrac{ \algbrac{X^{\flat}}{AX} }
               { X^{\flat} }.
  \end{aligned}
\end{align}
The first term on the \rhs\ equals $-\algbrac{g(X,X)X^{\flat}}{A} = -g(X,X)(AX)^{\flat}$ by part \ref{lem:mob2-vfs-bracs-1}, while the second term equals $-\algbrac{\algbrac{(AX)^{\flat}}{X}}{X^{\flat}}$ by symmetry in $X^{\flat}$ and $(AX)^{\flat}$.  Since this term equals the negative of the \lhs, rearranging \eqref{eq:mob2-vfs-bracs-8} gives
\vspace{0.15em}
\begin{equation} \label{eq:mob2-vfs-bracs-9}
  \algbrac{ \algbrac{(AX)^{\flat}}{X} }
          { X^{\flat} }
  = -\tfrac{1}{2} \big(
    g(X,X)(AX)^{\flat} + \algbrac{ \algbrac{AX}{X^{\flat}} }
                                 { X^{\flat} }
    \big).
\vspace{0.15em}
\end{equation}
Substituting \eqref{eq:mob2-vfs-bracs-9} into \eqref{eq:mob2-vfs-bracs-7}, the $g(X,X)(AX)^{\flat}$ terms cancel and we obtain
\begin{equation} \label{eq:mob2-vfs-bracs-10}
  \killing{ Y }
          { \algbrac{ \algbrac{AX}{X^{\flat}} }
                    { X^{\flat} } }
  = \killing{ Y }
            { 3\algbrac{ \algbrac{AX}{X^{\flat}} }
                       { X^{\flat} }
            + 2g(AX,X)X^{\flat} }.
\end{equation}
Noting that the first term on the \rhs\ of \eqref{eq:mob2-vfs-bracs-10} equals three times the \lhs, rearranging and dividing by two gives
\begin{equation*}
  \killing{ Y }
          { \algbrac{ \algbrac{AX}{X^{\flat}} }
                    { X^{\flat} }
          + g(AX,X)X^{\flat} }
    = 0
\end{equation*}
and hence the result by the \non degeneracy\ of the Killing form on $\fg_M$. \end{proof}

\begin{rmk} \thref{lem:mob2-vfs-bracs} is trivial to verify when working with a particular geometry, since the algebraic bracket $\algbrac{}{} : TM\times T^*M \to \fp^0_M$ is known explicitly.  However, care must be taken in the conformal case: the bracket in \itemref{lem:mob2-vfs-bracs}{2} equals
\begin{equation*}
  \algbrac{ \algbrac{AX}{X^{\flat}} }{ X^{\flat} }
     = -2g(AX,X)X^{\flat} + g(X,X)(AX)^{\flat},
\end{equation*}
so one must use that the $g$-\self adjoint\ sections of $\fp^0_M$ are multiples of the identity. \end{rmk}

It will be convenient henceforth to work with the \LC\ connection $\D$ of $h$, so that $X(t) = h(\D\pi(t),\bdot)$.  We can also trivialise the line bundles $\cL$ and $\cL^*$ using the pfaffian $\pi\defeq \pf{h} \in\s{0}{\cL}$ and its inverse $\pi^{-1} \in\s{0}{\cL^*}$.  As discussed in Subsection \ref{ss:mob2-vfs-pairing}, if $h_t$ is degenerate then $\pi(t)=0$ and hence $X(t)=0$ identically.

Form the usual endomorphism $A = \algbrac{h^{-1}}{\b{h}}$ and denote its transpose map by the same symbol.  The following result is well-known in the \proj\ \cite{bt2014-newholonomy, bkm2009-fubini} and \cproj\ \cite{bmr2015-localcproj, cemn2015-cproj} cases.  However, we pursue a different proof to these sources.

\begin{prop} \thlabel{prop:mob2-vfs-ADLcommute} Let $\D,\b{\D}$ be the \LC\ connections of $h,\b{h}$ respectively and suppose that $\b{Z}^{\D} \in\s{0}{\cL^*\tens TM}$ satisfies $\D\b{h} = \algbrac{\b{Z}^{\D}}{}$.  Then the endomorphisms $A, \D(\pi\b{Z}^{\D}) \in \s{0}{\alg{gl}{TM}}$ commute. \end{prop}

\begin{proof} From the prolongation \eqref{eq:ppg-bgg-metricprol} of the linear metric equation, we have
\begin{equation} \label{eq:mob2-vfs-ADLcommute-1}
  \D_X \b{Z}^{\D}
  = \b{\lambda}^{\D}X + \b{h}(\nRic{X}{},\bdot)
    + \quab^{-1} \liebdy \algbrac{\Weyl{}{}}{\b{h}}_X.
\end{equation}
Since we are assuming that $\D$ is \torsionfree, the Weyl curvature $\Weyl{}{}$ and hence also the final Weyl curvature term $\quab^{-1}\p \algbrac{\Weyl{}{}}{\b{h}}$ in \eqref{eq:mob2-vfs-ADLcommute-1} are \proj ly\ invariant, and therefore may be computed \wrt\ either $\D$ or $\b{\D}$; in the latter case we obtain
\begin{equation*}
  \quab^{-1}\p \algbrac{\Weyl{}{}} {\b{h}}_X
    = -\b{\lambda}^{\b{\D}}X - \b{h}(\nRic[\b{\D}]{X}{},\bdot)
\end{equation*}
since $\b{Z}^{\b{\D}}=0$.  Substituting this into \eqref{eq:mob2-vfs-ADLcommute-1}, we obtain
\begin{equation*}
  \D_X \b{Z}^{\D}
    = (\b{\lambda}^{\D} - \b{\lambda}^{\b{\D}})X
      - \b{h}(\nRic{X}{} - \nRic[\b{\D}]{X}{},\bdot).
\end{equation*}
Using \thref{cor:ppg-bgg-exactweyl} and that $A$ is \self adjoint\ \wrt\ $h^{-1}$, it follows that
\begin{align} \label{eq:mob2-vfs-ADLcommute-2}
  h^{-1}( \D_{AX} \b{Z}^{\D}, Y )
    &= (\b{\lambda}^{\D} - \b{\lambda}^{\b{\D}})h^{-1}(AX,Y)
      - h^{-1}( \b{h}(\nRic{AX}{} - \nRic[\b{\D}]{AX}{}, \bdot), Y) \notag \\
    &= (\b{\lambda}^{\D} - \b{\lambda}^{\b{\D}})h^{-1}(AY,X)
      - (\nRic{AY}{} - \nRic[\b{\D}]{AY}{})(AX) \notag \\
    &= h^{-1}(\D_{AY} \b{Z}^{\D}, X)
\end{align}
for all vector fields $X,Y$.  Since $h^{-1} = \pi \ltens g$, \eqref{eq:mob2-vfs-ADLcommute-2} now implies that $\D(\pi \b{Z}^{\D}) \circ A$ is \self adjoint\ \wrt\ $g$.  But by \thref{prop:mob2-pencils-Zgradient}, $\pi \ltens \b{Z}^{\D}$ is a gradient vector field \wrt\ $g$, so is $g$-\self adjoint\ by \thref{lem:mob2-vfs-gradient}.  Then since $A$, $\D(\pi \ltens \b{Z}^{\D})$ and $\D(\pi \ltens \b{Z}^{\D}) \circ A$ are all \self adjoint, it follows easily that $A$ commutes with $\D(\pi \ltens \b{Z}^{\D})$. \end{proof}

\vspace{0.1em}
\begin{lem} \thlabel{lem:mob2-vfs-derivs} Choose $t\in\bR$ such that $h_t$ is \non degenerate.  Then for all vector fields $X$ and all $\D \in \Dspace$, the following identities hold:
\begin{enumerate}
  \item \label{lem:mob2-vfs-derivs-1}
  $\D_X\pi(t) = h_t^{-1}(\b{Z}^{\D},X) \ltens \pi(t)$; and
  
  \item \label{lem:mob2-vfs-derivs-2}
  $\D_X h_t^{-1} = \algbrac{ h_t^{-1}(X,\bdot) }{ h_t^{-1}(\b{Z}^{\D},\bdot) }$.
\end{enumerate}
\end{lem}

\begin{proof} Since $h_t$ is \non degenerate, there is a connection $\D^t\in\Dspace$ with $\D^t h_t=0$.  Writing $\D=\D^t + \algbrac{}{\gamma_t}$ for some $\gamma_t \in \s{1}{}$, we then have $\D_X h_t = \algbrac{ \algbrac{X}{\gamma_t} }{ h_t } = \algbrac{h_t(\gamma_t,\bdot)}{X}$.  But also $\D_X h_t = \algbrac{Z_t^{\D}}{X}$ for a section $\b{Z}^{\D}$ of $\cL^* \tens TM$, which is independent of $t$; applying $\liebdy$ to both expressions gives $\gamma_t = h_t^{-1}(\b{Z}^{\D}, \bdot)$.

\smallskip

\proofref{lem:mob2-vfs-derivs}{1}
Since $\D^t h_t=0$ we also have $\D^t\pi(t)=0$, giving $\D_X\pi(t) = \algbrac{ \algbrac{X}{h_t^{-1}(\b{Z}^{\D},\bdot)} }{ \pi(t) } = h_t^{-1}(\b{Z}^{\D},X) \pi(t)$ as claimed.

\smallskip

\proofref{lem:mob2-vfs-derivs}{2}
Also $\D^th_t^{-1}=0$, so that $\D_X h_t^{-1} = \algbrac{ \algbrac{X}{h_t^{-1}(\b{Z}^{\D},\bdot)} }{ h_t^{-1} }$.  Applying the Jacobi identity and observing that the bracket between $h_t^{-1}(\b{Z}^{\D},\bdot)$ and $h_t^{-1}$ is trivial, we immediately obtain
\begin{equation*}
  \D_X h_t^{-1} = \algbrac{ \algbrac{X}{h_t^{-1}} }{ h_t^{-1}(\b{Z}^{\D}, \bdot) }
    = \algbrac{ h_t^{-1}(X,\bdot) }{ h_t^{-1}(\b{Z}^{\D},\bdot) }
\end{equation*}
as required. \end{proof}

We are finally in a position to prove that the vector fields $X(t)$ mutually commute.

\begin{proofof}{thm:mob2-vfs-vfs} The result is trivial if $s=t$, so assume that $s\neq t$.  Moreover if $h_t$ is degenerate then $X(t)=0$ identically, so we may assume that $h_s,h_t$ are both \non degenerate.

By \itemref{lem:mob2-vfs-derivs}{1} we have $\D\pi(t) = \pi(t) \ltens h_t^{-1}(\b{Z}^{\D},\bdot)$, so that writing $X(t)$ \wrt\ the \LC\ connection $\D$ of $h$ yields
\begin{equation} \label{eq:mob2-vfs-vfs-1}
  X(t) = \pi(t) \ltens h\big( h_t^{-1}(\b{Z}^{\D},\bdot), \bdot \big).
\end{equation}
Using both parts of \thref{lem:mob2-vfs-derivs}, we find that the covariant derivative of $X(t)$ \wrt\ $\D$ is given by
\vspace{0.05em}
\begin{equation} \label{eq:mob2-vfs-vfs-2} \begin{split}
  \D_Y X(t)
  &= \pi(t)h_t^{-1}(\b{Z}^{\D},Y) h(h_t^{-1}(\b{Z}^{\D}))
    + \pi(t) h\big( \algbrac{h_t^{-1}(Y)}{h_t^{-1}(\b{Z}^{\D})}(\b{Z}^{\D}) \big) \\
  &\qquad + \pi(t) h( h_t^{-1}(\D_Y\b{Z}^{\D})),
\end{split}
\vspace{0.05em}
\end{equation}
where we have written $h(\alpha) \defeq h(\alpha,\bdot)$ and $h_t^{-1}(X) \defeq h_t^{-1}(X,\bdot)$ for notational convenience.  Since the algebraic bracket between $\pi(t)h_t^{-1}(\b{Z}^{\D},Y) \in\s{0}{\cL}$ and $h(h_t^{-1}(\b{Z}^{\D})) \in\s{0}{\cL^*\tens TM}$ is just the contraction, the first term on the \rhs\ of \eqref{eq:mob2-vfs-vfs-2} equals
$\algbrac{h_t^{-1}(\pi(t)\b{Z}^{\D},Y) }
         { h(h_t^{-1}(\b{Z}^{\D})) } $.
Writing $h(h_t^{-1}(\b{Z}^{\D})) = \algbrac{h}{h_t^{-1}(\b{Z}^{\D})}$ and noting that the bracket between $\algbrac{\pi(t)\b{Z}^{\D}}{h_t^{-1}(Y)} \in\s{0}{\cL}$ and $h \in\s{0}{\cL^*\tens\cB}$ is trivial, applying the Jacobi identity to this term yields
\begin{align} \label{eq:mob2-vfs-vfs-3}
  &\algbrac{h_t^{-1}(\pi(t)\b{Z}^{\D},Y) }
           { h(h_t^{-1}(\b{Z}^{\D})) }
    \notag \\ &\quad
  = \Algbrac{ \algbrac{ \pi(t)\b{Z}^{\D} }
                      { h_t^{-1}(Y) } }
            { \algbrac{ h }
                      { h_t^{-1}(\b{Z}^{\D}) } }
    \notag \\ &\quad
  = \Algbrac{ h }
            { \algbrac{ \algbrac{ \pi(t)\b{Z}^{\D} }
                                { h_t^{-1}(Y) } }
                      { h_t^{-1}(\b{Z}^{\D}) } }
    \notag \\ &\quad\! \begin{aligned}
    & = \Algbrac{ h }
                { \algbrac{ \algbrac{ \pi(t)\b{Z}^{\D} }
                                    { h_t^{-1}(\b{Z}^{\D}) } }
                          { h_t^{-1}(Y) } } \\
    &\qquad
      + \Algbrac{ h }
                { \algbrac{ \pi(t)\b{Z}^{\D} }
                          { \algbrac{ h_t^{-1}(Y) }
                                    { h_t^{-1}(\b{Z}^{\D}) } } }
    \end{aligned}
    \notag \\ &\quad\! \begin{aligned}
    & = -h\big( \algbrac{ \algbrac{ \pi(t)\b{Z}^{\D} }
                                  { h_t^{-1}(\b{Z}^{\D}) } }
                       { h_t^{-1}(Y) }, \bdot \big) \\
    &\qquad
       - h\big( \algbrac{ h_t^{-1}(Y) }
                        { h_t^{-1}(\b{Z}^{\D}) } (\pi(t)\b{Z}^{\D},\bdot)
         , \bdot \big).
  \end{aligned}
\end{align}
The second term on the \rhs\ of \eqref{eq:mob2-vfs-vfs-3} cancels with the second term on the \rhs\ of \eqref{eq:mob2-vfs-vfs-2}, so that
\begin{equation} \label{eq:mob2-vfs-vfs-4}
  \D_Y X(t)
    = -h\big( \algbrac{ \algbrac{ \pi(t)\b{Z}^{\D} }
                                 { h_t^{-1}(\b{Z}^{\D}) } }
                       { h_t^{-1}(Y) }, \bdot \big)
    + \pi(t) h(h_t^{-1}(\D_Y\b{Z}^{\D})).
\end{equation}

Now let $Y=X(s) = \pi(s)h(h_s^{-1}(\b{Z}^{\D}))$.  We first deal with the second term on the \rhs\ above.  Since we may write $h_t = h(A_t\bdot,\bdot)$ for an invertible \self adjoint\ endomorphism $A_t: T^*M\to T^*M$, we have $h(h_t^{-1}(Y)) = A_t^{-1}Y$, where we use the same symbol to denote the transpose endomorphism $A_t:TM\to TM$.  Trivialising $\cL$ \wrt\ $\pi\defeq \pf{h}$ and writing $\pi(t) = p(t)\pi$ for all $t\in\bR$, \eqref{eq:mob2-vfs-vfs-1} becomes $X(t) = p(t)A_t^{-1}(\pi\b{Z}^{\D})$.  Therefore, the term under scrutiny becomes
\begin{equation} \label{eq:mob2-vfs-vfs-5}
  \pi(t) h(h_t^{-1}(\D_{X(s)} \b{Z}^{\D}))
    = p(s)p(t) A_t^{-1}\D_{A_s^{-1}(\pi \b{Z}^{\D})} (\pi\b{Z}^{\D}).
\end{equation}
Since $A=\b{h}\circ h^{-1} : TM \to TM$ and $\D(\pi\b{Z}^{\D})$ commute by \thref{prop:mob2-vfs-ADLcommute}, so also do $A_t^{-1}$ and $\D(\pi\b{Z}^{\D})$; combining this with the fact that $A_s^{-1}$ and $A_t^{-1}$ commute for all $s,t\in\bR$, we see that \eqref{eq:mob2-vfs-vfs-5} is symmetric in $s,t$.  Thus, it suffices to show that the first term alone in \eqref{eq:mob2-vfs-vfs-4} (with $Y=X(s)$) vanishes upon alternation in $s,t$.

Since we may write $h=\tfrac{1}{s-t}(h_t-h_s)$, we have
\begin{equation*}
  h_t^{-1}(X(s),\bdot)
    = \tfrac{1}{s-t} \pi(s) \big( h_s^{-1}(\b{Z}^{\D}) - h_t^{-1}(\b{Z}^{\D}) \big),
\end{equation*}
and the first term on the \rhs\ of \eqref{eq:mob2-vfs-vfs-4} becomes
\begin{equation*}
  -\tfrac{1}{s-t}p(s)p(t) h \big(
    \algbrac{ \algbrac{ \pi\b{Z}^{\D} }
                      { h_t^{-1}(\b{Z}^{\D}) } }
            { \pi h_s^{-1}(\b{Z}^{\D}) - \pi h_t^{-1}(\b{Z}^{\D}) }, \, \bdot
  \big).
\end{equation*}
Since both $\tfrac{1}{s-t}p(s)p(t)$ and $\pi h_s^{-1}(\b{Z}^{\D}) - \pi h_t^{-1}(\b{Z}^{\D})$ are skew-symmetric in $s,t$, alternation yields
\begin{equation} \label{eq:mob2-vfs-vfs-6}
  \tfrac{1}{s-t}p(s)p(t) h \big(
    \algbrac{ \algbrac{ \pi\b{Z}^{\D} }
                      { h_s^{-1}(\b{Z}^{\D}) - h_t^{-1}(\b{Z}^{\D}) } }
            { \pi h_s^{-1}(\b{Z}^{\D}) - \pi h_t^{-1}(\b{Z}^{\D}) }, \, \bdot
  \big).
\end{equation}
Considering the metric $g \defeq \pi^{-1}h^{-1}$ corresponding to the linear metric $h$, we have $h_s^{-1}(\b{Z}^{\D}) - h_t^{-1}(\b{Z}^{\D}) = g \big( (A_s^{-1}-A_t^{-1})(\pi\b{Z}^{\D}), \bdot \big)$.  It is straightforward to see that $A_s^{-1}-A_t^{-1} = (s-t)A_s^{-1}A_t^{-1}$, so that the previous expression yields
\begin{equation*}
  \pi \b{Z}^{\D} = \tfrac{1}{s-t}A_s A_t \big( h_s^{-1}(\b{Z}^{\D}) - h_t^{-1}(\b{Z}^{\D}) \big)^{\sharp}.
\end{equation*}
Then \eqref{eq:mob2-vfs-vfs-6} is of the form $\algbrac{ \algbrac{A_{s,t} X}{X^{\flat}} }{ \pi X^{\flat} }$, where
\begin{align*}
  A_{s,t} &\defeq \tfrac{1}{(s-t)^2}p(s)p(t) A_s A_t\\
  \text{and}\quad
  X &\defeq h \big( h_s^{-1}(\b{Z}^{\D}) - h_t^{-1}(\b{Z}^{\D}), \bdot \big).
\end{align*}
Here $A_{s,t}$ is a $g$-\self adjoint\ element of $\fp^0_M$ by \thref{lem:ppg-bgg-p0B}, so the Leibniz rule and \itemref{lem:mob2-vfs-bracs}{2} yield
\begin{align*}
  \algbrac{ \algbrac{A_{s,t}X}{X^{\flat}} }{ \pi X^{\flat} }
  &= \algbrac{ \algbrac{A_{s,t}X}{X^{\flat}} }{ \pi } \tens X^{\flat}
   + \pi \tens \algbrac{ \algbrac{A_{s,t}X}{X^{\flat}} }{ X^{\flat} } \\
  &= g(X,A_{s,t}X)X^{\flat} - g(A_{s,t}X,X)X^{\flat} = 0.
\end{align*}
We conclude that $\D_{X(s)}X(t)$ is symmetric in $s,t$; since we assume that $\D$ is \torsionfree, it follows that $[X(s),X(t)] = \D_{X(s)}X(t) - \D_{X(t)}X(s) = 0$ for all $s,t \in \bR$. \end{proofof}

\begin{rmk} \thlabel{rmk:mob2-vfs-rmks} The author suspects a more conceptually satisfactory proof of \thref{thm:mob2-vfs-vfs} is available, which exploits the facts that $\bs{h}$ and $\bs{\pi}$ are in the kernels of the appropriate BGG operators, and that $X(\bs{h},\bs{\pi})$ is the result of the BGG bilinear pairing operator $\cW\times\cW^* \to \fg$.  One might expect strong restrictions on the pairings when there is a $2$-dimensional family of solutions of each BGG operator.  Unfortunately, the author was unable to make progress in this direction. \end{rmk}

\section{Relative eigenvalues and the order of a pencil} 
\label{s:mob2-evals}

Suppose that $M$ admits a pencil $\bs{h}$ of linear metrics containing a \non degenerate\ metric, which we view as a section of $\cL^*\tens\cB \tens \cO_V(1)$ over $M\times \pr{V}$.  The pfaffian $\bs{\pi} \defeq \pf{\bs{h}}$ lies in the kernel of the hessian \eqref{eq:ppg-bgg-hesseqn} by \thref{cor:ppg-bgg-metrichess}.  Choosing an affine chart on $\pr{V}$ with linear metrics $h,\b{h}$ at $\infty,0$, \thref{cor:mob2-pencils-pfpoly} allows us to view $\bs{\pi}$ as a polynomial $\pi(t)$ of degree $n$ in $t$.  Writing $\b{h} = h(A\bdot,\bdot)$ as before,
\begin{equation*}
  \pi(t) \defeq \pf{h_t} = (\det (A-t\id))^{1/r}(\pf{h})
\end{equation*}
and therefore the (possibly complex-valued) roots $\xi:M\to\bC$ of $\pi(t)$ are precisely the eigenvalues of $A$.  The following is immediate from \thref{cor:mob2-pencils-pfpoly}.

\begin{cor} \thlabel{cor:mob2-evals-mult} The eigenvalues of $A$ have algebraic multiplicity divisible by $r$. \noproof \end{cor}

Forgetting the affine chart on $\pr{V}$ and viewing $\bs{\pi}$ as a homogeneous polynomial of degree $n$ in the parameters of the pencil, the roots of $\bs{\pi}$ are functions $\xi : M \to \pr[\bC]{V\tens\bC}$.

\begin{defn} \thlabel{defn:mob2-evals} The \emph{(relative) eigenvalues} of $\bs{h}$ are the roots $\xi:M\to \pr[\bC]{V\tens\bC}$ of the pfaffian $\bs{\pi}=\pf{\bs{h}}$.  In an affine chart, we identify these with functions $\xi : M \to \bC$. \end{defn}

Note that since $A$ is \self adjoint\ \wrt\ all metrics in the pencil, its eigenvalues are either real-valued or occur in complex-conjugate pairs.  In particular, if $\bs{h}$ contains a positive definite linear metric then the eigenvalues are necessarily real-valued.  Since we will not assume this, it will be necessary to deal with the Jordan normal form of $\bs{h}$, whose details we recall in Subsection \ref{ss:mob2-evals-jordan}.  The eigenvalues of $\bs{h}$ are not \apriori\ smooth functions $\xi:M\to\bC$, so we must deal with the set of points on which the eigenvalues are ``well-behaved''.  There are two relevant notions of well-behaved here; we study both in Subsection \ref{ss:mob2-evals-regular} and show them to be equivalent, following Topalov's treatment in the \proj\ case \cite{t2000-hierarchy} and in $PQ^{\epsilon}$-\proj\ geometry \cite{t2003-geodesiccomp}.

In Subsection \ref{ss:mob2-evals-order} we define an integer invariant of a pencil, called its \emph{order}.  This definition is well-known (and more readily accessible) in \riem\ \cproj\ geometry, where the order equals the maximal dimension of the pointwise span of the canonical Killing fields \cite{acg2006-ham2forms1}.  Recently this definition has been extended to \pseudo\riem\ \cproj\ geometry \cite[\S 5.6]{cemn2015-cproj}.
Finally, we contemplate pencils containing a \riem\ metric in Subsection \ref{ss:mob2-evals-riem}.

\vspace{-0.2em}
\subsection{Review of Jordan normal forms} 
\label{ss:mob2-evals-jordan}
\vspace{-0.1em}

The theory of Jordan normal forms over the complex numbers is well-known, but less common over the real numbers.  We review this theory here for completeness; a readable introduction may be found in \cite{w2009-jordan}.

Let $V$ be a complex vector space and $A\in\alg{gl}{V}$ an endomorphism.  The minimal polynomial of $A$ splits as $m_A(t) = \prod{i=1}{k} (t-\xi_i)^{m_i}$ for $\xi_i$ the eigenvalues of $A$.  If an eigenvalue $\xi$ has geometric multiplicity $m$, the Jordan decomposition theorem states that there is a basis of $V$ in which the restriction of $A$ to the \emph{generalised eigenspace} $G_A(\xi) \defeq \ker(A-\xi\id)^m$ of $\xi$ is represented by a sum of the \emph{Jordan blocks} of the form
\vspace{-0.2em}
\begin{equation} \label{eq:mob2-evals-jordan}
  J_k(\xi) =
  \begin{bmatrix}
    \xi    & 1      & \cdots & 0      & 0      \\
    0      & \xi    & \cdots & 0      & 0      \\
    \vdots & \vdots & \ddots & \vdots & \vdots \\
    0      & 0      & \cdots & \xi    & 1      \\
    0      & 0      & \cdots & 0      & \xi
  \end{bmatrix}
  \colvectpunct[-2.3em]{,}
\vspace{-0.2em}
\end{equation}
with $A$ represented by a direct sum of such blocks.  For later use, we note that if $\xi\neq 0$ then the Jordan block $J_k(\xi)$ is invertible, with inverse
\vspace{-0.2em}
\begin{equation} \label{eq:mob2-evals-invjordan}
  J_k(\xi)^{-1} =
  \begin{bmatrix}
    \xi^{-1} & -\xi^{-2} & \cdots & (-1)^{k+1}\xi^{-k} \\
    0        & \xi^{-1}  & \cdots & (-1)^k\xi^{-k+1}   \\
    \vdots   & \vdots    & \ddots & \vdots             \\
    0        & 0         & \cdots & \xi^{-1}
  \end{bmatrix}
  \colvectpunct[-2.4em]{.}
\vspace{-0.2em}
\end{equation}
The part $v_1,\ldots,v_k$ of the Jordan basis spanning $G_A(\xi)$ is called the \emph{Jordan chain}, and satisfies $Av_i = \xi v_i + v_{i-1}$.   The following lemma will be useful later.

\begin{lem} \thlabel{lem:mob2-evals-jordanli} Let $v_1,\ldots,v_k$ be generalised eigenvectors of an endomorphism $A$ corresponding to distinct eigenvalues $\xi_1,\ldots,\xi_k$.  Then $v_1,\ldots,v_k$ are linearly independent. \end{lem}

\begin{proof} By definition, there are $m_i\in\bN$ such that $(A-\xi_i\id)^{m_i}v_i=0$ for each $i=1,\ldots,k$.  Suppose first that $v\in \ker(A-\xi_i)^{m_i} \cap \ker(A-\xi_j)$ for $i\neq j$.  Then $Av = \xi_j v$ yields $(A-\xi_i\id)^{m_i}v = (\xi_j-\xi_i)^{m_i}v = 0$, so that $\xi_i\neq \xi_j$ implies $v=0$ and hence
\begin{equation*}
  \ker(A-\xi_i\id)^{m_i} \cap \ker(A-\xi_j\id) = 0.
\end{equation*}
Now suppose that $v\in \ker(A-\xi_i\id)^{m_i} \cap \ker(A-\xi_j\id)^{m_j}$ for $i\neq j$ and $m_j>1$.  Then $(A-\xi_j\id)^{m_j-1}v \in \ker(A-\xi_i\id)^{m_i} \cap \ker(A-\xi_j\id) = 0$ by the above, so that $v \in \ker(A-\xi_i\id)^{m_i} \cap \ker(A-\xi_j\id)^{m_j-1}$.  By induction we conclude that $v=0$.

We prove the lemma by induction on $k$.  It is obvious for $k=1$, so assume it holds for some $k\geq 1$ and let $v_i \in \ker(A-\xi_i\id)^{m_i}$ for $i=1,\ldots,k+1$.  Suppose that $\sum{i=1}{k+1}a_iv_i = 0$.  Applying $(A-\xi_{k+1}\id)^{m_{k+1}}$ to this linear combination, we obtain
\begin{equation*}
  \sum{i=1}{k+1} a_i(A-\xi_{k+1}\id)^{m_{k+1}}v_i
    = \sum{i=1}{k} a_i(A-\xi_{k+1}\id)^{m_{k+1}}v_i = 0.
\end{equation*}
By the first paragraph we must have $(A-\xi_{k+1}\id)^{m_{k+1}}v_i \neq 0$ for $i\neq k+1$, so that $(A-\xi_{k+1}\id)^{m_{k+1}}v_i \in \ker(A-\xi_i\id)^{m_i}$ implies that $a_i=0$ for $i=1,\ldots,k$ by the induction hypothesis.  Since $v_{k+1}\neq 0$, it follows also that $a_{k+1}=0$. \end{proof}

The Jordan normal form over a real vector space $V$ may be described as follows.  The inclusion $\alg{gl}{n,\bC} \injto \alg{gl}{2n,\bR}$ identifies each complex entry $z=a+\cpx{i}b\in\bC$ of a matrix in $\alg{gl}{n,\bC}$ with the $2\by 2$ real matrix $\inlinematrix{a & -b \\ b & \phantom{-}a} \in \alg{gl}{2,\bR}$.  Complexifying $V$ and writing each complex eigenvalue of $A$ in real and complex parts, this substitution in the Jordan normal form \eqref{eq:mob2-evals-jordan} yields we arrive at the real Jordan normal form of $A$.

If $v_1,\ldots,v_k \in V\tens\bC$ comprise a Jordan chain for an eigenvalue $\xi\in \bC$, we can write $\xi = a+\cpx{i}b$ and $v_i=x_i+\cpx{i}y_i$ for $x_i,y_i\in V$.  Comparing real and imaginary parts, the identity $Av_i = \xi v_i + v_{i-1}$ gives $Ax_i = ax_i - by_i + x_{i-1}$ and $Ay_i = ay_i + bx_i + y_{i-1}$.  However, it is often convenient to complexify and avoid the real Jordan normal form altogether.

\subsection{Stable and regular points} 
\label{ss:mob2-evals-regular}

This subsection acts mainly as a technical stepping stone for defining the order of a pencil in Subsection \ref{ss:mob2-evals-order}.  Since we are not assuming our metrisability pencil $\bs{h}$ admits \riem\ members, we must entertain the possibility that $\bs{h}$ admits \non trivial\ Jordan blocks.  Choosing an affine chart with linear metrics $h,\b{h}$ at $\infty,0$ respectively, this means that the minimal and characteristic polynomials of the usual endomorphism $A=h^{-1}\circ \b{h}$ of $T^*M$ need not coincide.  The minimal polynomial yields a section $\tilde{\pi}(t)$ of $\cL$, which is polynomial in $t$ of degree at most $n = \deg\pi(t)$ and whose roots are the eigenvalues of $\bs{h}$.  Denote by $\rk{\bs{h}\at{x}}$ the degree of $\tilde{\pi}(t)$ at $x\in M$, and write $\rk{\bs{h}} = \max[x\in M]{\rk{\bs{h}\at{x}}}$.

\begin{defn} \thlabel{defn:mob2-evals-regular} Let $x\in M$.
\begin{enumerate}
  \item \label{defn:mob2-evals-regular-st}
  Call $x$ a \emph{stable point} if there is an open neighbourhood $U\ni x$ on which $\rk{\bs{h}\at{y}}$ is constant for all $y\in U(x)$.  Denote the set of stable points in $M$ by $M^0$.  A point is \emph{singular} if it is not regular.

  \item \label{defn:mob2-evals-regular-reg}
  Call $x$ a \emph{regular point} if the rank $\rk{\bs{h}\at{x}}$ is maximal among all points in $M$, \ie\ $\rk{\bs{h}\at{x}} = \rk{\bs{h}}$.  Denote the set of regular points by $M^{\mr{reg}}$.
\end{enumerate}
\end{defn}

The goal of this subsection is to show that, if $M$ is connected, the stable set $M^0$ and regular set $M^{\mr{reg}}$ coincide.  We adapt the proof given by Topalov in $PQ^{\epsilon}$-\proj\ geometry \cite{t2003-geodesiccomp}, which is largely algebraic and depends critically on the quadratic integrals of the geodesic flow constructed in Subsection \ref{ss:mob2-pencils-killing}.

\begin{lem} \thlabel{lem:mob2-evals-dense} Fix $x\in M$ and suppose that $U\subset\bR$ is an open interval containing no eigenvalues of $\bs{h}\at{x}$.  Then:
\begin{enumerate}
  \item \label{lem:mob2-evals-dense-1}
  $\dim\linspan{h_t^*\at{x}}{t\in U} = \rk{\bs{h}\at{x}}$;
  
  \item \label{lem:mob2-evals-dense-2}
  $\dim\linspan{h_t^*\at{x}(X,\bdot)}{t\in U} \leq \rk{\bs{h}\at{x}}$ for every $X\in T_xM$, with equality for $X$ lying in a dense open subset of $T_xM$; and
  
  \item \label{lem:mob2-evals-dense-3}
  if $t_1,\ldots,t_k\in U$ are $k \defeq \rk{\bs{h}\at{x}}$ distinct real numbers, the pointwise Killing $2$-tensors $h_{t_1}^*\at{x}, \ldots,$ $h_{t_k}^*\at{x}$ are linearly independent elements of $\cL^2_x\tens\cB^*_x$.
\end{enumerate}
\end{lem}

\begin{proof} For notational simplicity, we omit the evaluation at $x$ from all tensors in this proof.  Since by \thref{prop:mob2-pencils-adjpoly} the Killing $2$-tensor $h_t^*$ is a polynomial of degree $n-1$ in $t$, we may write $h_t^* = \sum{i=0}{n-1} (-1)^i k_i\,t^i$, where the $k_i$ are Killing $2$-tensors by linearity of the Killing equation \eqref{eq:mob2-pencils-killingeqn}.  Since this is a Vandermonde system (see \cite[App.\ B]{acg2006-ham2forms1}) and therefore invertible, we have
\begin{equation} \label{eq:mob2-evals-denselem-1}
\begin{split}
  \linspan{h_t^*}{t\in U}
    &= \linspan{k_0,\ldots,k_{n-1}}{} \\
  \text{and}\quad
  \linspan{h_t^*(X,\bdot)}{t\in U}
    &= \linspan{k_0(X,\bdot),\ldots,k_{n-1}(X,\bdot)}{}
\end{split}
\end{equation}
for each $X\in T_xM$.  Let $\ell\in\bN$ denote the maximal rank of the $h_i^*(X,\bdot)$, \ie\
\begin{equation*}
  \ell \defeq \max[X\in T_xM \,]{ \dim{\linspan{h_0^*(X,\bdot),\ldots,h_{n-1}^*(X,\bdot)}{}} }.
\end{equation*}
The remainder of the proof is more straightforward if we assume without loss of generality that the linear metric $h$ at $\infty$ is in Jordan normal form, for which we should complexify our current picture.  Denoting the set of eigenvalues of $\bs{h}\at{x}$ by $S\subset\bC$, complexifying yields complex Killing $2$-tensors $h_z^*$ for $z\in \bC\setminus S$.  As above, we may write $h_z^* = \sum{i=0}{n-1} (-1)^i k_i^{\bC} \, z^i$ for complex Killing $2$-tensors $k_i^{\bC}\in (\cL^2_x\tens\cB^*_x)\tens\bC$.  Clearly then
\begin{equation} \label{eq:mob2-evals-denselem-3}
  \dim\linspan{h_0^*, \ldots, h_{n-1}^*}{}
    = \dim[\bC]\linspan[\bC]{k_0^{\bC}, \ldots, k_{n-1}^{\bC}}{},
\end{equation}
and arguing using a Vandermonde system as above yields a complex analogue of \eqref{eq:mob2-evals-denselem-1}:
\begin{equation} \label{eq:mob2-evals-denselem-4}
\begin{split}
  \linspan[\bC]{h_z^*}{z\in\bC\setminus S}
    &= \linspan[\bC]{k_0^{\bC}, \ldots, k_{n-1}^{\bC}}{} \\
  \text{and}\quad
  \linspan[\bC]{h_z^*(X,\bdot)}{z\in\bC\setminus S}
    &= \linspan[\bC]{k_0^{\bC}(X,\bdot), \ldots, k_{n-1}^{\bC}(X,\bdot)}{}
\end{split}
\end{equation}
for each $X\in \bC T_xM$.  Combining \eqref{eq:mob2-evals-denselem-1}, \eqref{eq:mob2-evals-denselem-3} and \eqref{eq:mob2-evals-denselem-4} we obtain
\begin{equation} \label{eq:mob2-evals-denselem-5}
\begin{split}
  \dim\linspan{h_t^*}{t\in U} &= \dim[\bC]\linspan[\bC]{h_z^*}{z\in \bC\setminus S} \\
  \text{and}\quad
  \ell &= \max[X\in \bC T_xM \,]{ \dim[\bC]\linspan[\bC]{h_z^*(X,\bdot)}{z\in \bC\setminus S} }.
\end{split}
\end{equation}
With these preliminaries in hand, we may continue with the proof of the lemma.

\smallskip

\proofref{lem:mob2-evals-dense}{1}
Continue to work with the complexification $h^{\bC}$ of $h$, given in Jordan normal form.  By the formula $h_z^* = (\det[\bC] h_z)^{1/r}(h_z^{\bC})^{-1}$ for $z\in \bC\setminus S$, the direct sum decomposition of $h^{\bC}$ contains a Jordan block $J_k(\xi)$ if and only if the direct sum decomposition of $h_z^*$ contains a Jordan block $\pi(t) J_k(\xi-z)^{-1}$.  Since the size $k$ of this block is at most the multiplicity of $\xi$ as a root of $\pi(z)$, it is clear that $\rk{\bs{h}^{\bC}\at{x}} = \dim[\bC] \linspan[\bC]{h_z^*}{z\in \bC\setminus S}$ by \eqref{eq:mob2-evals-invjordan}.  Combining this with \eqref{eq:mob2-evals-denselem-5} completes the proof.

\smallskip

\proofref{lem:mob2-evals-dense}{2}
We have $\dim[\bC]\linspan[\bC]{h_z^*(X,\bdot)}{z\in\bC\setminus S} \leq \ell = \rk{\bs{h}\at{x}}$ by \eqref{eq:mob2-evals-denselem-3} and the second equation of \eqref{eq:mob2-evals-denselem-5}.  By definition of $\ell$, there is an $X\in T_xM$ and a multi-index $I = (i_1,\ldots,i_{\ell})$ for $i_j \in \setof{0,\ldots,n-1}{}$ such that the weighted $\ell$-covector
\begin{equation*}
  h_I^*(X,\bdot) \defeq
  h_{i_1}^*(X,\bdot) \wedge\cdots\wedge h_{i_{\ell}}^*(X,\bdot)
  \in \cL^{2\ell}_x \tens \Wedge{\ell} T^*_xM
\end{equation*}
is \non zero, with $h_J^*(X,\bdot)=0$ for multi-indices $J$ of length greater than $\ell$.  Choosing a basis in $T_xM$, the equation $h_I^*(X,\bdot)=0$ is a homogeneous polynomial of degree $\ell$ in the components of $X$.  It follows that the set of tangent vectors $X$ for which $h_I^*(X,\bdot)\neq 0$ is open and dense in $T_xM$.

\smallskip

\proofref{lem:mob2-evals-dense}{3}
Writing $h_t^* = \sum{i=1}{n-1} (-1)^i k_i \, t^i$ as above, evidently the $k_i$ are linearly independent and independent of $t$.  If $t_1,\ldots,t_k$ are $k \defeq \rk{\bs{h}\at{x}} = \dim\linspan{k_0,\ldots,k_{n-1}}{}$ distinct real numbers, a linear combination of the $h_{t_i}^*$ has the form
\begin{equation} \label{eq:mob2-evals-denselem-6} \begin{aligned}
  \sum{i=1}{k} a_i h_{t_i}^*
    &= \sum{i=1}{k} \sum{j=1}{n-1} (-1)^j a_i k_j t_i^j \\
    &= (-1)^{n-1}k_{n-1} \sum{i=1}{k} a_i t_i^{n-1} + \cdots + k_0 \sum{i=1}{k} a_i.
\end{aligned} \end{equation}
This is a $k\by k$ Vandermonde system relating the $t_i$ and $a_i$, which is invertible because the $t_i$ are distinct.  Thus, a linear combination \eqref{eq:mob2-evals-denselem-6} is zero if and only if the $a_i$ all vanish, \ie\ $h_{t_1}^*, \ldots, h_{t_k}^*$ are linearly independent. \end{proof}

Our next task is to show that a generic point of $M$ is stable.

\begin{prop} \thlabel{prop:mob2-evals-dense} The set $M^0 \subseteq M$ of stable points in $M$ is open and dense. \end{prop}

\begin{proof} Clearly $M^0 \subseteq M$ is open by construction.  We claim that if $x\in M$ is singular, then every open neighbourhood $U \ni x$ contains a point $y \in U$ for which $\rk{\bs{h}\at{y}} > \rk{\bs{h}\at{x}}$.

To prove the claim, choose $k \defeq \rk{\bs{h}\at{x}}$ distinct real numbers $t_1, \ldots, t_k$ which are not eigenvalues of $A\at{x}$.  Using \itemref{lem:mob2-evals-dense}{3}, the Killing $2$-tensors $h_{t_1}^*\at{x}, \ldots, h_{t_k}^*\at{x}$ are linearly independent elements of $\cL^2_x \tens \cB^*_x$, and by the smoothness of $h_t^*$ and the implicit function theorem there is an open neighbourhood $U \ni x$ on which this linear independence holds.  By \itemref{lem:mob2-evals-dense}{1},
\begin{equation*}
  \dim\linspan{h_t^*\at{y}}{t\in\bR} = \rk{\bs{h}\at{y}}
\end{equation*}
for all $y \in U$, so that the fact that $h_{t_1}^*\at{y}, \ldots, h_{t_k}^*\at{y}$ are linearly independent implies that $\rk{\bs{h}\at{y}} \geq k = \rk{\bs{h}\at{x}}$ for all $y\in U$.  However $x$ is singular, so that $\rk{\bs{h}}$ is \non constant\ on every open neighbourhood of $x$.  In particular, every open neighbourhood of $x$ contains a point $y$ for which $\rk{\bs{h}\at{y}} > \rk{\bs{h}\at{x}}$.

It remains to show that $M^0 \subseteq M$ is dense.  Suppose that $S \subset M$ is an open set consisting of singular points, and take $x_1 \in S$ with $\rk{\bs{h}\at{x_1}} = k_1$.  Applying the previous claim to $x=x_1$ yields a singular point $x_2 \in S$ with rank $\rk{\bs{h}\at{x_2}} = k_2 > k_1$.  Iterating now yields a singular point $x_{\ell} \in S$ of maximal rank, any neighbourhood of which contains a point of strictly larger rank.  This contradicts the maximality of $\rk{\bs{h}\at{x_{\ell}}}$, so that the singular open subset $S$ cannot exist.  It follows that the singular subset $S \defeq M\setminus M^0$ has empty interior; equivalently, the closure of $M^0 = M \setminus (M \setminus S)$ equals the whole of $M$, so that $M^0$ is dense in $M$. \end{proof}

The inclusion $M^{\mr{reg}} \subseteq M^0$ was proved in the final paragraph of \thref{prop:mob2-evals-dense}.

\begin{cor} \thlabel{cor:mob2-evals-dense} Every regular point is stable, \ie\ $M^{\mr{reg}} \subseteq M^0$.  \qed \end{cor}

\begin{thm} \thlabel{thm:mob2-evals-regular} If $M$ is connected, every stable point is regular; thus $M^{\mr{reg}}=M^0$. \end{thm}

\begin{proof} Suppose that $x\in M^0$ with $\rk{\bs{h}\at{x}}=k$.  By \itemref{lem:mob2-evals-dense}{3} there are distinct real numbers $t_1,\ldots,t_k \in\bR$ such that the pointwise Killing $2$-tensors $h_{t_1}^*\at{x}, \ldots, h_{t_k}^*\at{x}$ are linearly independent.  By smoothness this also holds in a neighbourhood $U\ni x$ and, shrinking $U$ if necessary, we may assume that $\bs{h}$ has constant rank $k$ on $U$.  It follows that $h_t^* = \sum{i=1}{k} a_i h_{t_i}^*$ on $U$ for smooth functions $a_i\in\s[U]{0}{\bR}$.  Trivialising $\cL^2$ \wrt\ the metric $h$ at $\infty$, both $(\pf{h})^{-2}h_t^*$ and the $(\pf{h})^{-2}h_{t_i}^*$ are integrals of the geodesic flow of $g\defeq(\pf{h})^{-1}h^{-1}$, so by \thref{prop:mob2-pencils-integrals} they Poisson-commute with the hamiltonian of $g$.  By linearity of the Poisson bracket, it follows that the functions $a_i$ are also integrals of $g$, so are constant on $U$.

Viewing $h_t^*\at{x}$ as a quadratic function on $T_xM$, we have $\d_X h_t^*\at{x} = 2h_t^*\at{x}(X,\bdot)$ for all $X\in T_xM$ once we make the canonical identification $T_X T_xM \isom T_xM$.  By \itemref{lem:mob2-evals-dense}{2}, there is a dense open subset $V\subset T_xM$ for which $h_{t_1}^*\at{x}(X,\bdot), \ldots,$ $h_{t_k}^*\at{x}(X,\bdot)$ are linearly independent, so that $\d_X h_{t_1}^*\at{x}, \ldots, \d_X h_{t_k}^*\at{x}$ are also linearly independent on $V$.  By the implicit function theorem, their linear independence holds in a small neighbourhood of $x$ which, shrinking $U$ if necessary, coincides with $U$.

Now let $x_0\in M^0$ and $x\in M^{\mr{reg}}$ have ranks $k_0\leq k$ respectively, and suppose that $k_0 < k$.  By the above, there are $k$ distinct real numbers $t_1,\ldots,t_k$ and a neighbourhood $U\ni x$ on which $h_{t_1}^*\at{x}, \ldots, h_{t_k}^*\at{x}$ are functionally independent.  Suppose that $x,x_0$ may be connected by a geodesic $\gamma$, with $\gamma(0)=x$ and $\gamma(1)=x_0$.  If $\gamma'(0)\in T_xM$ does not lie in the subset $V\subset T_xM$ from above, replace $\gamma$ with a geodesic with $\gamma(0)=x$ and $\gamma'(0)\in V$.  Since $V$ is dense in $T_xM$ and the exponential map is smooth, $\gamma'(0)\in V$ can be chosen so that $\gamma(1)=\exp \gamma'(0)$ lies in the neighbourhood of $x_0$ on which $\rk{\bs{h}\at{x_0}}=k_0$ is constant.  Replacing $x_0$ with $\gamma(1)\in M^0$, we may assume without loss of generality that $\gamma'(0)\in V$.  By the above, $\d_{\gamma'(0)} h_{t_1}^* \at{x}, \ldots,$ $\d_{\gamma'(0)} h_{t_k}^* \at{x}$ are linearly independent at $x$.  Moreover since $\d_X h_t^*\at{x} = 2h_t^*\at{X}(X,\bdot)$ from above, the $(\pf{h})^{-2}\at{x} \d_X h_{t_i}^* \at{X}$ are also integrals of $g$.  Noting that $\d_{\gamma'(1)} h_{t_i}^* \at{x_0}$ is the image of $\d_{\gamma'(0)} h_{t_i}^*\at{x}$ under the $1$-parameter family of diffeomorphisms of $TM$ induced by the geodesic flow of $g$, we conclude that $\d_{\gamma'(1)} h_{t_1}^*\at{x_0}, \ldots, \d_{\gamma'(1)} h_{t_k}^* \at{x_0}$ are linearly independent at $x_0$.

On the other hand, by \itemref{lem:mob2-evals-dense}{2} we know that the $h_{t_i}^*\at{x_0}$ have span at most $k_0 < k$, and therefore by the previous paragraph we must have $h_{t_{k_0+1}}^*\at{x_0} = \sum{i=1}{k_0} a_i h_{t_i}^*\at{x_0}$ for constants $a_i$.  Differentiating gives $\d_{\gamma'(1)} h_{t_{k_0+1}}^*\at{x_0} = \sum{i=1}{k_0} a_i \d_{\gamma'(1)}h_{t_i}^*\at{x_0}$, contradicting the linear independence of $\d_{\gamma'(1)} h_{t_1}^*\at{x_0}, \ldots, \d_{\gamma'(1)} h_{t_k}^*\at{x_0}$ deduced in the previous paragraph.  Thus it is impossible that $k_0<k$, giving $k_0=k$ by maximality of $k$, and therefore all stable points which are geodesically connected to a regular point are themselves regular.  Finally, if $M$ is connected, any two points may be connected by a finite sequence of geodesic segments, so that $M^{\mr{reg}} \subseteq M^0$. \end{proof}

Recall that a set of functions are \emph{functionally independent} if their differentials are linearly independent on a dense open subset.  The following is immediate.

\begin{cor} \thlabel{cor:mob2-evals-regular} If $t_1,\ldots,t_k$ are $k\defeq \rk{\bs{h}}$ distinct real numbers which are not eigenvalues of $\bs{h}$, the Killing $2$-tensors $h_{t_1}^*, \ldots, h_{t_k}^*$ are functionally independent on $TM$.  \qed \end{cor}

When $M$ is connected, we shall dispense with the notation $M^{\mr{reg}}$ and simply write $M^0$ for the subset of regular points, which necessarily have maximal rank $\rk{\bs{h}}$.

\medskip
\subsection{The order of a pencil} 
\label{ss:mob2-evals-order}
\vspace{0.1em}

The order of a metrisability pencil $\bs{h}$ will be a canonical integer $\ell\in\bZ$ associated to the eigenvalues of $\bs{h}$.  To define it, we will split the study of the eigenvalues of $\bs{h}$ into the study of the constant and \non constant\ eigenvalues.  We begin by analysing eigenvalues of higher geometric multiplicity, following the methods of \cite[Prop.\ 1]{bm2013-priemlocalform}.

\begin{thm} \thlabel{thm:mob2-evals-gm2} Let $U\subset M^0$ be an open subset of regular points of $M$, and suppose that $\xi:U\to\bC$ is a smooth eigenvalue of $\bs{h}$ with geometric multiplicity at least two at all points of $U$.  Then:
\begin{enumerate}
  \item \label{thm:mob2-evals-gm2-const}
  $\xi$ is a constant function; and
  
  \item \label{thm:mob2-evals-gm2-two}
  if $M$ is connected, $\xi$ is an eigenvalue of $\bs{h}$ of geometric multiplicity at least two at all points of $M$.
\end{enumerate}
\end{thm}

\begin{proof} \proofref{thm:mob2-evals-gm2}{const}
Choose an affine chart for $\bs{h}$ with \non degenerate\ linear metrics at $\infty, 0$, and form the usual endomorphism $A$ defined by $\b{h} = h(A\bdot,\bdot)$.  Since $A$ is self-adjoint \wrt\ any metric in the pencil, there are two cases to consider: either
\begin{enumerate*}[label=(\alph*)]
  \item \label{item:mob2-evals-gm2-const-a} $\xi:U\to\bR$ is real-valued; or
  \item \label{item:mob2-evals-gm2-const-b} $\xi, \conj{\xi} :  U\to\bC$ is a pair of complex-conjugate eigenvalues.
\end{enumerate*}
We treat these two cases individually.

\smallskip

\proofref{item:mob2-evals-gm2}{const-a}
Suppose that $\xi:U\to\bR$ is real-valued and \non constant.  Then there is a point $x_0\in U$ such that $\d\xi\at{x_0}\neq 0$, and, shrinking $U$ if necessary, it follows that $x_0$ is a regular value of $\xi$. Then by Sard's theorem, the level set
\vspace{-0.2em}
\begin{equation*}
  U_0 \defeq \setof{x\in U}{\xi(x)=\xi(x_0)} \subseteq U
\vspace{-0.2em}
\end{equation*}
is a regular submanifold of $U$ of codimension one.  By the assumption that $\xi$ is \non constant, we can choose a point $y\in U$ with $y\notin U_0$.  Let $\gamma_{x,y}$ be a geodesic joining a point $x\in U_0$ to $y$ with $\gamma_{x,y}(0)=x$ and $\gamma_{x,y}(1)=y$, and consider all such geodesics as $x\in U_0$ varies.  Choose $y\notin U_0$ in a small neighbourhood of $x_0$ onto which the exponential map is a diffeomorphism from a neighbourhood of zero in $T_{x_0}M$.  Since $U_0$ has codimension one, it follows that for such $y$ the subset
\vspace{-0.2em}
\begin{equation} \label{eq:mob2-evals-gm2-1}
  V_{x_0,y} \defeq \bigcup_{x \in U_0} \linspan{\gamma_{x,y}'(1)}{} \subseteq T_yM
\vspace{-0.4em}
\end{equation}
contains a \non empty\ open subset of $T_yM$.  Since $\xi$ is also a root of $\pi(t)\defeq(\det h_t)^{1/r}$, the linear metric $h_{\xi(x_0)} \at{x}$ is degenerate and hence $(\pf{h_{\xi(x_0)}}) \ltens h_{\xi(x_0)}^* \at{x} = 0$ for all $x\in U_0$.  In particular, $(\pf{h_{\xi(x_0)}}) \ltens h_{\xi(x_0)}^*\at{x} (\gamma_{x,y}'(0),\gamma_{x,y}'(0))=0$, yielding $(\pf{h_{\xi(x_0)}}) \ltens h_{\xi(x_0)}^* \at{y}(v,v)=0$ for all $v\in V_{x_0,y}$.  Since $(\pf{h_{\xi(x_0)}}) \ltens h_{\xi(x_0)}^*\at{y}$ is bilinear and $V_{x_0,y}$ contains a \non empty\ open subset of $T_yM$, it follows that $(\pf{h_{\xi(x_0)}}) \ltens h_{\xi(x_0)}^*\at{y} = 0$ and hence $\pf{h_{\xi(x_0)}} \at{y} = 0$.  Thus $\xi(x_0)$ is an eigenvalue of $\bs{h}$ at $y$, so that $\xi$ is constant in a neighbourhood of $x_0$.  This contradicts our assumption that $\d\xi\at{x_0}\neq 0$.

\smallskip

\proofref{item:mob2-evals-gm2}{const-b}
Suppose now that $\xi, \conj{\xi}:U\to\bC$ are a pair of complex-conjugate eigenvalues.  For each $z\in\bC$, $h_z^*$ is a complex Killing $2$-tensor, so both its real and imaginary parts are real Killing $2$-tensors.  Suppose that $x_0\in U$ is a point with $\d\xi\at{x_0}\neq 0$, and let $U_0 \defeq \setof{ x \in M }{ \xi(x)=\xi(x_0) }$ as in \ref{item:mob2-evals-gm2-const-a}.

Suppose first that the differentials of $\Re{\xi}$ and $\Im{\xi}$ are proportional.  Then, shrinking $U$ if necessary, $U_0$ is a regular submanifold of codimension one in $U$.  Arguing as in \ref{item:mob2-evals-gm2-const-a}, we conclude that both $\Re{\xi}, \Im{\xi}$ and hence $\xi$ are constant in a neighbourhood of $x_0$, contradicting our assumption that $\d\xi\at{x_0}\neq 0$.

Suppose instead that the differentials of $\Re{\xi}$ and $\Im{\xi}$ are not proportional.  Then, shrinking $U$ if necessary, $U_0$ is regular submanifold of $U$ of codimension two.  Let $y\in U$ and, as before, consider all geodesics $\gamma_{x,y}$ connecting points $x\in U_0$ to $y$; we suppose that $\gamma_{x,y}(0)=x$ and $\gamma_{x,y}(1)=y$.  Arguing as in \ref{item:mob2-evals-gm2-const-a}, if $y$ lies in a sufficiently small open subset of $x_0$ then the subset $V_{x_0,y}$ defined by \eqref{eq:mob2-evals-gm2-1} contains a \non empty\ submanifold of $T_yM$ of codimension one.  Then, since
\vspace{-0.15em}
\begin{equation*}
  h_{\xi(x_0)}^* \at{y}(X,X)
  = h_{\xi(x_0)}^*\at{x} \left( \gamma_{x,y}'(0), \gamma_{x,y}'(0) \right) = 0
\vspace{-0.15em}
\end{equation*}
for all $X\in V_{x_0,y}$, we conclude that $\mathopsl{Re}\! \big[ h_{\xi(x_0)}^* \at{y}(X,X) \big]$ and $\mathopsl{Im} \! \left[ h_{\xi(x_0)}^* \at{y}(X,X) \right]$ are proportional for all $X\in T_yM$.  Since these components are Killing $2$-tensors, the function of proportionality must be constant, so that there is a complex number $0\neq a+b\cpx{i} \in\bC$ satisfying
\vspace{-0.1em}
\begin{equation} \label{eq:mob2-evals-gm2-2}
  (a+b\cpx{i}) h_{\xi(x_0)}^* \at{y}(X,X) = (a-b\cpx{i}) h_{\conj{\xi}(x_0)}^* \at{y}(X,X)
\vspace{-0.1em}
\end{equation}
for all $X\in T_yM$.  At points $y\in U$ with $y\notin U_0$, $\xi(y)\neq \xi(x_0)$ and thus $h_{\xi(x_0)}\at{y}$ is \non degenerate.  Taking a trace in \eqref{eq:mob2-evals-gm2-2} \wrt\ the linear metric $h$ at $\infty$ yields that $h_{\xi(x_0)}^*\at{y}$ and $h_{\conj{\xi}(x_0)}^*\at{y}$ are proportional, implying that $\xi(x_0) = \conj{\xi}(x_0)$.  This contradicts our assumption that $\xi$ is not real.

\smallskip

\proofref{thm:mob2-evals-gm2}{two}
Suppose now that $\xi:U\to\bC$ is a constant (real- or complex-valued) eigenvalue of $\bs{h}$ of geometric multiplicity at least two.  Choose a point $y\in M\setminus U$ which can be connected to a point $x\in U$ by a geodesic $\gamma$, with $\gamma(0)=y$ and $\gamma(1)=x\in U$; thus $\gamma=\gamma_{y,\gamma'(0)}$.  Consider also the geodesics $\gamma_{y,Y}$ with $\gamma_{y,Y}(0)=y$ and $\gamma_{y,Y}'(0)=Y\in T_yM$.  Since the exponential map is smooth, there is an open neighbourhood $V\subset T_yM$ of $\gamma'(0)$ such that $\gamma_{y,Y}(1) = \exp_y Y$ lies in $U$ for all $Y\in V$.  Now, because $\xi$ is an eigenvalue of geometry multiplicity at least two, we have $h_{\xi}^* \at{x} (\gamma_{y,Y}'(1), \gamma_{y,Y}'(1)) = 0$ for all $Y\in V$.  Since $(\pf{h})^{-2}h_{\xi}^*$ is an integral of $g\defeq (\pf{h})^{-2}h^{-1}$, this yields $h_{\xi}^*\at{y} (\gamma_{y,v}'(0), \gamma_{y,v}'(0)) = h_{\xi}^*\at{y} (Y,Y)=0$ for all $Y\in V$ also.  Then $h_{\xi}^* \at{y}$ vanishes on a \non empty\ open subset of $T_yM$, so vanishes identically on the whole of $T_yM$ by bilinearity, and it follows that $\xi$ is an eigenvalue of geometric multiplicity at least two at $y$.

Any point in a neighbourhood $U_0 \ni y$ is also geodesically connected to a point in $U$, so $\xi$ is an eigenvalue of geometric multiplicity at least two in $U_0$.  Finally, if $M$ is connected, then any two points may be connected by a piecewise geodesic curve. \end{proof}

In particular, the contrapositive of \itemref{thm:mob2-evals-gm2}{const} implies that there is a single Jordan block for each \non constant\ eigenvalue of $\bs{h}$.  Note however that the opposite implication does not hold in general: constant eigenvalues may still have geometric multiplicity one.

Next we would like to analyse the derivatives of the eigenvalues of the pencil $\bs{h}$.  In \cproj\ geometry, Calderbank \etal\ have recently shown \cite[Cor.\ 5.17]{cemn2015-cproj} that $\d\xi$ is an eigenform with eigenvalue $\xi$; see also \cite{bmr2015-localcproj}.  This is not quite true for a general \ppg.  We will need the following lemma; this is presumably well-known, but the author could not find a reference.

\begin{lem} \thlabel{lem:mob2-evals-gmproj} Suppose that $A : V \to V$ is a linear endomorphism of a vector space $V$, and let $\xi$ of be an eigenvalue of $A$ with generalised eigenspace $G_A(\xi)$.  Then the linear map given by projection onto $G_A(\xi)$ is polynomial in $A$. \end{lem}

\begin{proof} Let $m$ denote the multiplicity of $\xi$ as a root of the minimal polynomial $m_A(t)$ of $A$.  Then $q(t) \defeq m_A(t) / (t-\xi)^m$ is polynomial in $t$ which is coprime to $(t-\xi)^m$, so by Euclid's algorithm there are polynomials $a(t), b(t)$ such that
\begin{equation} \label{eq:mob2-evals-gmproj-1}
  a(t) q(t) = 1 - b(t) (t-\xi)^m.
\end{equation}
Substituting the linear operator $A$ in place of $t$, we have
\begin{align*}
  \big( a(A) q(A) \big)^2
    &= a(A) q(A) \big( 1 - b(A) (A-\xi\id)^m \big) \\
    &= a(A) q(A) - a(A) b(A) q(A) (A-\xi\id)^m \\
    &= a(A) q(A) - a(A) b(A) m_A(A) \\
    &= a(A) q(A)
\end{align*}
since $m_A(A) = 0$ by definition, so that $a(A) q(A)$ is a projection.  Moreover, since $G_A(\xi)$ coincides with the kernel of $(A - \xi\id)^m$, the identity $(A - \xi\id)^m q(A) = m_A(A) = 0$ implies that $q(A)$ takes values in $G_A(\xi)$.  Since $A$ preserves $G_A(\xi)$, it follows that $a(A) q(A)$ also takes values in $G_A(\xi)$.  Finally, \eqref{eq:mob2-evals-gmproj-1} implies that $a(A) q(A)$ restricts to the identity on $G_A(\xi)$.  Therefore $a(A) q(A)$ is the projection to $G_A(\xi)$, which is polynomial in $A$ as required. \end{proof}

\begin{thm} \thlabel{thm:mob2-evals-dxi} Let $\xi : M \to \bC$ be an eigenvalue of the pencil $\bs{h}$ in a neighbourhood $U$ of regular points.  Choose an affine chart for $\bs{h}$ with \non degenerate\ metrics $h, \b{h}$ at $\infty, 0$, and form the usual endomorphism $A$ defined by $\b{h} = h(A\bdot, \bdot)$.  Then $\d\xi$ is a generalised eigenform of $A$ with eigenvalue $\xi$ on $U$. \end{thm}

\begin{proof} Choose an affine chart of $\bs{h}$ with metrics $h, \b{h}$ at $\infty, 0$, with $h$ \non degenerate, and form the usual endomorphism $A$ defined by $\b{h} = h(A\bdot, \bdot)$.  Then $\xi$ is an eigenvalue of $A$, and in a neighbourhood of a point $x \in U$, we can choose a local frame \wrt\ which $A$ is in Jordan normal form.
If $\xi$ is constant then $\d\xi=0$ and there is nothing to prove. Therefore we assume that $\d\xi\at{x} \neq 0$ for at least one point $x \in U$.

Suppose that the generalised eigenspace $G_A(\xi)$ has dimension $k$.  Then by the contrapositive of \itemref{thm:mob2-evals-gm2}{const}, $G_A(\xi)$ is spanned by a single Jordan chain $\alpha_k, \ldots, \alpha_1$ of generalised eigenforms, which satisfy $A\alpha_i = \xi\alpha_i + \alpha_{i-1}$ with $\alpha_0 \defeq 0$.  Let $X_k, \ldots, X_1$ be the dual coframe.  Then the transpose endomorphism $A^{\top} : TM \to TM$ satisfies
\vspace{-0.1em}
\begin{equation} \label{eq:mob2-evals-dxi-1}
\begin{aligned}
  \alpha_i(A^{\top}X_i)
    &= (A\alpha_i)(X_j) \\
    &= (\xi\alpha_i + \alpha_{i-1})(X_j)
     = \alpha_i(\xi X_j + X_{j+1}),
\end{aligned}
\vspace{-0.1em}
\end{equation}
so that the restriction of $A^{\top}$ to $G_A(\xi)$ is also in Jordan normal \wrt\ the $X_i$, providing we reverse their order.  In particular $A^{\top}X_i = \xi X_i + X_{i+1}$, where $X_{k+1} \defeq 0$.

Let $\D \in \Dspace$ be the \LC\ connection of $h$ provided by \thref{cor:ppg-bgg-lcconn}.  Then $\D_X \b{h} = \algbrac{\b{Z}^{\D}}{X}$ for some $\b{Z}^{\D} \in \s{0}{\cL^* \tens TM}$ by \eqref{eq:ppg-bgg-metricprol}, so that the defining relation $\b{h} = h(A\bdot, \bdot)$ gives $(\D_X A)\alpha = \algbrac{ h^{-1} }{ \algbrac{ \algbrac{\b{Z}^{\D}}{X} }{ \alpha } }$ for all $\alpha \in \s{1}{}$.  For all $i=1,\ldots,k$, differentiating the identity $A\alpha_i = \xi\alpha_i + \alpha_{i-1}$ then gives
\vspace{-0.1em}
\begin{equation} \label{eq:mob2-evals-dxi-2}
  \algbrac{ h^{-1} }
          { \algbrac{ \algbrac{\b{Z}^{\D}}{X} }
                    { \alpha_i } }
    + A(\D_X \alpha_i)
    = \d\xi(X)\alpha_i + \xi(\D_X \alpha_i) + \D_X \alpha_{i-1}.
\vspace{-0.1em}
\end{equation}
Now \eqref{eq:mob2-evals-dxi-1} implies that
\vspace{-0.1em}
\begin{equation*}
  (A \D_X \alpha_i)(X_i)
    = (\D_X \alpha_i)(A^{\top} X_i)
    = (\D_X \alpha_i)(\xi X_i + X_{i+1}),
\vspace{-0.1em}
\end{equation*}
so that by contracting \eqref{eq:mob2-evals-dxi-2} with $X_i$ and rearranging we obtain
\begin{align*}
  \d\xi(X)
    &=  \killing{ \algbrac{ h^{-1} }
                          { \algbrac{ \algbrac{\b{Z}^{\D}}{X} }
                                    { \alpha_i } } }
                { X_i } \\
    &\qquad
      + (\D_X \alpha_i)(\xi X_i + X_{i+1})
      - \xi(\D_X \alpha_i)(X_i)
      - (\D_X \alpha_{i-1})(X_i) \\
    &=  \killing{ \algbrac{ h^{-1} }
                          { \algbrac{ \algbrac{\b{Z}^{\D}}{X} }
                                    { \alpha_i } } }
                { X_i } \\
    &\qquad
       + (\D_X \alpha_i)(X_{i+1}) - (\D_X \alpha_{i-1})(X_i).
\end{align*}
Summing over $i=1,\ldots,k$, the two terms on the last line telescopically cancel, leaving
\begin{equation} \label{eq:mob2-evals-dxi-3}
  k \, \d\xi(X) = \Sum{i=1}{k} \, 
    \killing{ \algbrac{ h^{-1} }
                      { \algbrac{ \algbrac{\b{Z}^{\D}}{X} }
                                { \alpha_i } } }
            { X_i }.
\end{equation}
For each $i$, invariance of the Killing form gives
\begin{equation*}
  \killing{ \algbrac{ h^{-1} }
                    { \algbrac{ \algbrac{\b{Z}^{\D}}{X} }
                              { \alpha_i } } }
          { X_i }
    =  \killing{ \algbrac{ \algbrac{ \alpha_i }
                                   { \algbrac{X_i}{h^{-1}} } }
                         { \b{Z}^{\D} } }
               { X }.
\end{equation*}
On the other hand $\d\xi(X) = \killing{\d\xi}{X}$, thus giving
\begin{equation} \label{eq:mob2-evals-dxi-5}
  k \, \d\xi = \Sum{i=1}{k} \,
    \algbrac{ \algbrac{ \alpha_i }
                      { \algbrac{X_i}{h^{-1}} } }
            { \b{Z}^{\D} }
\end{equation}
by \non degeneracy\ of the Killing form.  For each summand, the Jacobi identity yields
\begin{align} \label{eq:mob2-evals-dxi-4}
  \algbrac{ \algbrac{ \alpha_i }
                    { \algbrac{X_i}{h^{-1}} } }
          { \b{Z}^{\D} }
    &=  \algbrac{ \algbrac{ \algbrac{\alpha_i}{X_i} }
                          { h^{-1} } }
                { \b{Z}^{\D} } \notag \\
    &=  \algbrac{ \algbrac{ \algbrac{\alpha_i}{X_i} }
                          { \b{Z}^{\D} } }
                { h^{-1} }
      + \algbrac{ \algbrac{\alpha_i}{X_i} }
                { \algbrac{h^{-1}}{\b{Z}^{\D}} } \notag \\
    &=  h^{-1} \big( \algbrac{ \algbrac{X_i}{\alpha_i} }
                             { \b{Z}^{\D} },
                     \bdot \big)
      - \algbrac{ \algbrac{X_i}{\alpha_i} }
                { h^{-1}(\b{Z}^{\D}, \bdot) }.
\end{align}
For notational convenience, write $\pi \defeq \pf{h}$ and suppose that $g \defeq \pi^{-1} \ltens h^{-1}$ is the metric associated to $h$.  Then we can write $\b{Z}^{\D} = \pi^{-1} \ltens \Lambda$ for $\Lambda \defeq \pi \ltens \b{Z}^{\D} \in \s{0}{TM}$,%
\footnote{By \thref{prop:mob2-pencils-Zgradient}, $\Lambda$ coincides with the vector field of the same name appearing in \cite{bmr2015-localcproj, cmr2015-cprojmob, mr2012-yanoobata}.}
and the first term on the \rhs\ of \eqref{eq:mob2-evals-dxi-4} evaluates to
\vspace{-0.1em}
\begin{equation*} \begin{aligned}
  h^{-1} \big( \algbrac{ \algbrac{X_i}{\alpha_i} }
                       { \pi^{-1} \ltens \Lambda }, \bdot \big)
    &= -h^{-1}( \alpha_i(X_i) \pi^{-1} \ltens \Lambda, \bdot)
      + h^{-1}( \pi^{-1} \ltens \algbrac{ \algbrac{X_i}{\alpha_i} }
                                        { \Lambda }, \bdot) \\
    &= -\Lambda^{\flat}
      + \algbrac{ \algbrac{X_i}{\alpha_i} }
                { \Lambda }^{\flat}
\end{aligned}
\vspace{-0.1em}
\end{equation*}
It follows from \eqref{eq:mob2-evals-dxi-4} and \threfit{lem:app-alg-musical} that
\begin{align*}
  \algbrac{ \algbrac{ \alpha_i }
                    { h^{-1}(X_i,\bdot) } }
          { \b{Z}^{\D} }
    &= -\Lambda^{\flat}
      + \algbrac{ \algbrac{X_i}{\alpha_i} }
                { \Lambda }^{\flat}
      - \algbrac{ \algbrac{X_i}{\alpha_i} }
                { \Lambda^{\flat} } \\
    &= -\Lambda^{\flat}
      + \algbrac{ \algbrac{X_i}{\alpha_i} }
                { \Lambda }^{\flat}
      + \algbrac{ \algbrac{\alpha_i^{\sharp}}{X_i} }
                { \Lambda }^{\flat} \\
    &= -\Lambda^{\flat}
      + \algbrac{ \algbraco{\id}{g}{X_i,\alpha_i^{\sharp}} }
                { \Lambda }^{\flat}.
\end{align*}
Substituting the last display into \eqref{eq:mob2-evals-dxi-5}, applying $\sharp$ to both sides and dividing by $k \neq 0$, we obtain
\begin{equation} \label{eq:mob2-evals-dxi-6}
  \d\xi^{\sharp}
    = -\Lambda
     + \tfrac{1}{k} \Sum{i=1}{k} \,
       \algbrac{ \algbraco{\id}{g}{X_i,\alpha_i^{\sharp}} }
               { \Lambda }.
\end{equation}

We next evaluate the symmetrised algebraic bracket.  First observe that since $A$ is $g$-\self adjoint, the musical isomorphisms $\flat = \sharp^{-1} : TM \to T^*M$ restrict to isomorphisms of the $\xi$-generalised eigenspaces $G_{A^{\top}}$ of $A^{\top}$ and $G_A(\xi)$ of $A$.  Since $\{X_i\}_{i=1}^k$ is a basis of $G_{A^{\top}}(\xi)$ with dual basis $\{\alpha_i\}_{i=1}^k$, for any $Y \in \s{0}{TM}$ and any $\beta \in \s{1}{}$ we have
\vspace{0.2em}
\begin{equation} \label{eq:mob2-evals-dxi-7}
  \Sum{i=1}{k} \beta( \algbrac{ \algbrac{X_i}{\alpha_i} }{ Y } )
    = \Sum{i=1}{k} \alpha_i( \algbrac{ \algbrac{Y}{\beta} }{ X_i } )
    = \tr( \algbrac{Y}{\beta} \at{G_{A^{\top}}(\xi)} ),
\vspace{0.2em}
\end{equation}
where $\bdot\at{G_{A^{\top}}(\xi)}$ denotes the restriction to $G_{A^{\top}}(\xi)$.  If $\Pi_{\xi} : \s{0}{TM} \to \s{0}{TM}$ is the projection onto $G_{A^{\top}}(\xi)$ then it is easy to see that $\Pi_{\xi}$ is \self adjoint\ \wrt\ $g$.  Moreover $\Pi_{\xi}$ is polynomial in $A^{\top} \in \s{0}{\fp^0_M}$ by \thref{lem:mob2-evals-gmproj}, hence itself a section of $\fp^0_M$ by \thref{cor:ppg-bgg-powA}.  Therefore we may apply \threfit{lem:app-alg-trA} to \eqref{eq:mob2-evals-dxi-7}, yielding
\begin{align*}
  \tr( \algbrac{Y}{\beta} \at{G_A(\xi)} )
    = \tr( \algbrac{Y}{\beta} \circ \Pi_{\xi} )
    = \tfrac{1}{2} \big( (\tr \Pi_{\xi})\beta(Y) + r\beta(\Pi_{\xi}Y) \big).
\end{align*}
Since $G_A(\xi)$ has dimension $k$ by assumption, the projection $\Pi_{\xi}$ has trace $k$.  Moreover if $Y$ lies in a generalised eigenspace for an eigenvalue $\mu \neq \xi$, we have $\Pi_{\xi} Y = 0$ and hence $\tr( \algbrac{Y}{\beta} \at{G_{A^{\top}}(\xi)} ) = \tfrac{k}{2} \beta(Y)$.  Since $\beta \in \s{1}{}$ is arbitrary, we conclude from \eqref{eq:mob2-evals-dxi-7} that
\begin{equation*}
  \Sum{i=1}{k} \, \algbrac{X_i}{\alpha_i}
    = \tfrac{k}{2}\id \mod G_{A^{\top}}(\xi).
\end{equation*}
Since $\{X_i^{\flat}\}_{i=1}^k$ and $\{\alpha_i^{\sharp}\}_{i=1}^k$ are also dual bases of $G_A(\xi)$ and $G_{A^{\top}}(\xi)$, we similarly obtain $\Sum{i=1}{k} \algbrac{X_i^{\flat}}{\alpha_i^{\sharp}} = \tfrac{k}{2}\id \bmod G_{A^{\top}}(\xi)$, thus giving
\vspace{-0.1em}
\begin{equation} \label{eq:mob2-evals-dxi-8}
  \Sum{i=1}{k} \, \algbraco{\id}{g}{X_i,\alpha_i^{\sharp}}
    = k \, \id \mod G_{A^{\top}}(\xi).
\vspace{-0.1em}
\end{equation}

Now write $\Lambda = \Sum{\mu}{} \Lambda_{\mu}$ according to the generalised eigenspace decomposition of $A^{\top}$, where the sum runs over all distinct eigenvalues of $A$, and $\Lambda_{\mu} \in G_{A^{\top}}(\mu)$.  Applying \eqref{eq:mob2-evals-dxi-8} to the expression \eqref{eq:mob2-evals-dxi-6} for $\d\xi^{\sharp}$, we conclude that
\vspace{-0.1em}
\begin{equation} \label{eq:mob2-evals-dxi-9}
  \d\xi^{\sharp} = - \Lambda_{\xi}
    + \Sum{i=1}{k} \, \algbrac{ \algbraco{\id}{g}{X_i,\alpha_i^{\sharp}} }
                              { \Lambda_{\xi} }.
\vspace{-0.1em}
\end{equation}
It remains to see that $\algbrac{X_i}{\alpha_i}$ preserves $G_{A^{\top}}(\xi)$.  If $\xi$ has multiplicity $m$ as a root of the minimal polynomial of $A$ then $G_{A^{\top}}(\xi)$ coincides with the kernel of $A_{\xi}^m \defeq (A^{\top} - \xi\id)^m$, which is a section of $\fp^0_M$ by \thref{cor:ppg-bgg-powA}.  Then for any $Y \in G_{A^{\top}}(\xi)$, we have
\begin{align*}
  A_{\xi}^m \big( \algbrac{ \algbrac{X_i}{\alpha_i} }
                          { Y } \big)
    &=  \algbrac{ \algbrac{ \algbrac{A_{\xi}^m}{X_i} }
                          { \alpha_i } }
                { Y }
      + \algbrac{ \algbrac{ X_i }
                          { \algbrac{A_{\xi}^m}{\alpha_i} } }
                { Y } \\
    &\qquad
      + \algbrac{ \algbrac{X_i}{\alpha_i} }
                { \algbrac{A_{\xi}^m}{Y} } \\
    &=  0,
\end{align*}
whence $\algbrac{ \algbrac{X_i}{\alpha_i} }{ Y } \in G_{A^{\top}}(\xi)$ also.  Similarly $\algbrac{\alpha_i^{\sharp}}{X_i}$ preserves $G_{A^{\top}}(\xi)$, so that \eqref{eq:mob2-evals-dxi-9} implies that $\d\xi^{\sharp} \in G_{A^{\top}}(\xi)$.  Finally, we recall that $\flat$ restricts to an isomorphism between $G_{A^{\top}}(\xi)$ and $G_A(\xi)$, giving $\d\xi \in G_A(\xi)$ as required. \end{proof}

\begin{cor} \thlabel{cor:mob2-evals-order} There is an integer $\ell\in\bN$, equal to the number of distinct \non constant\ eigenvalues of $\bs{h}$, such that the span of the $\d\xi_i\at{x}$ has dimension $\ell$ for all $x\in M^0$.  In particular, the $\xi_i$ are functionally independent on $M$. \end{cor}

\begin{proof} Suppose $x\in M^0$ is a regular point.  Then, by the results of Subsection \ref{ss:mob2-evals-regular}, the minimal polynomial of $\bs{h}$ has constant and maximal degree in a neighbourhood $U\ni x$; equivalently, the number of distinct eigenvalues of $\bs{h}\at{y}$ is constant and maximal for all $y\in U$.  It follows that if $\xi_j(x)=\xi_k(x)$ for $j\neq k$ then $x\notin M^0$, so that the $\xi_i$ are pointwise distinct on $M^0$; \cf\ \cite[Rmk.\ 5.4]{cemn2015-cproj}.  Applying \thref{thm:mob2-evals-dxi}, the $\d\xi_i\at{x}$ lie in distinct generalised eigenspaces of $\bs{h}\at{x}$ and are therefore linearly independent at $x$ by \thref{lem:mob2-evals-jordanli}.  Since $M^0\subset M$ is open and dense by \thref{prop:mob2-evals-dense}, it also follows that the $\xi_i$ are functionally independent. \end{proof}

The integer $\ell\in\bN$ from \thref{cor:mob2-evals-order} is then an invariant of the metrisability pencil $\bs{h}$, prompting the following definition; \cf\ \cite{acg2006-ham2forms1, acgt2004-ham2forms2, bmr2015-localcproj}.

\begin{defn} \thlabel{defn:mob2-evals-order} The integer $\ell \defeq \dim[\bC]\linspan[\bC]{\d\xi_1\at{M_0}, \ldots, \d\xi_n\at{M_0}}{}$ satisfying $0 \leq \ell \leq n$ from \thref{cor:mob2-evals-order} is called the \emph{order} of the metrisability pencil $\bs{h}$. \end{defn}

In the case $\ell = n$, all roots of $\pi(t)$ are distinct and \non constant.  Then the minimal and characteristic polynomials of $\bs{h}$ coincide, \ie\ $\tilde{\pi}(t)=\pi(t)$, and each eigenvalue has both characteristic and geometric multiplicities equal to one.  In particular, $\bs{h}$ is semisimple.  On the other hand, if $\ell=0$ then $\bs{h}$ has no \non constant\ eigenvalues.  Applying \thref{prop:mob2-pencils-Zgradient} to an affine chart we conclude that $\b{Z}^{\D}=0$, so that all linear metrics in $\bs{h}$ are affinely equivalent.

\begin{rmk} \thlabel{rmk:mob2-evals-amgm} In \cproj\ geometry (and also hypercomplex geometry), the existence of Killing vector fields allows us to replace ``geometric multiplicity'' with ''algebraic multiplicity'' in the previous results.  Indeed, the Killing fields provide linear integrals of the geodesic flow, and let us argue using the characteristic polynomial $\pi(t)$ of $A$ rather than the minimal polynomial $\tilde{\pi}(t)$; see \cite[Lem.\ 5.16]{cemn2015-cproj} for details. \end{rmk}

\vspace{0.5em}
\subsection{Special features of \riem\ pencils} 
\label{ss:mob2-evals-riem}
\vspace{0.3em}

Let $\bs{h}$ be a \non degenerate\ pencil of linear metrics.  Choosing an affine chart with linear metrics $h,\b{h}$ at $\infty,0$ respectively, we remarked above that the eigenvalues of a pencil $\bs{h}$ are either real-valued or come in complex-conjugate pairs.  In particular, if it is possible to choose $h,\b{h}$ to be positive definite, the eigenvalues of $\bs{h}$ are real-valued.  This is the situation most commonly encountered in the literature.

\begin{defn} \thlabel{defn:mob2-evals-riempencil} A \emph{\riem\ pencil} is a pencil admitting an affine chart with linear metrics of \riem\ signature at $\infty,0$. \end{defn}

Note that it is sufficient for the pencil to admit a single linear metric $\b{h}$ of \riem\ signature.  Indeed, putting this metric at $0$ and any other linear metric $h$ at $\infty$ defines an affine chart, and $h_t \defeq \b{h}-th$ will have \riem\ signature for $t$ close to zero.

In projective geometry, it is known that the (necessarily real-valued) eigenvalues of such a pencil have a global ordering \cite{bm2003-geombenenti, m2005-globalordering, tm2003-geodintegrability}, and a similar result is implicit in \cproj\ by the results of \cite{acgt2004-ham2forms2}.  Adapting the proof given by Bolsinov and Matveev in \cite[Thm.\ 3]{bm2003-geombenenti}, we generalise this result to all \ppgs\ admitting a \riem\ pencil $\bs{h}$.  We fix an affine chart with \riem\ linear metrics $h,\b{h}$ at $\infty,0$ and denote the (possibly indistinct) $n$ roots of $\pi(t)\defeq \pf{h_t}$ by $\xi_1,\ldots,\xi_n$.

\begin{lem} \thlabel{lem:mob2-evals-ordering} Let $x\in M$, and consider the roots $t_1(x,X) \leq\cdots\leq t_{n-1}(x,X)$ of the degree $n-1$ polynomial $h_t^*\at{x}(X,X)$ at $x$.  Then:
\begin{enumerate}
  \item \label{lem:mob2-evals-ordering-1}
  $t_i(x,X) \in\bR$ and satisfies $\xi_i(x) \leq t_i(x,X) \leq \xi_{i+1}(x)$ for all $X\in T_xM$ and every $i=1,\ldots,n-1$; and
  
  \item \label{lem:mob2-evals-ordering-2}
  if $\xi_i(x)<\xi_{i+1}(x)$ then for every $c\in\bR$ the set
  \begin{equation*}
    V_c \defeq \setof{ X\in T_xM }{ t_i(x,X) = c } \subseteq T_xM
  \end{equation*}
  has zero Lebesgue measure.
\end{enumerate}
\end{lem}

\begin{proof} \proofref{lem:mob2-evals-ordering}{1}
Let $g\defeq (\pf{h})^{-1} \ltens h^{-1}$ be the metric associated to $h$.  Since $h,\b{h}$ are \riem\ and the $\xi_i$ are the eigenvalues of the usual endomorphism $A$, we can decompose $T_xM = \Dsum{i=1}{n}E_A(\xi_i)\at{x}$ into the eigenspaces of $A\at{x} \in\alg{gl}{T_xM}$.  In this decomposition, $h_t^*\at{x} \in \cL^2_x\tens\cB^*_x$ takes the diagonal block form
\begin{equation} \label{eq:mob2-evals-ordering-1}
  h_t^*\at{x}
    = (\pf{h}\at{x})^2 \Sum{i=1}{n} \left( \prod{j\neq i}{} (\xi_j(x)-t) \right) g_i\at{x},
\end{equation}
where $g_i \defeq g\at{E_A(\xi_i)}$ is the restriction of $g$ to the $\xi_i$-eigendistribution $E_A(\xi_i)$.  Let $X\in T_xM$ and write $X = X_1 + \cdots + X_n$ according to the eigenspace decomposition of $T_xM$.  Since the coefficients of $h_t^*\at{x}(X,X)$ depend continuously on the $\xi_i(x)$ and $X$, it suffices to prove \ref{lem:mob2-evals-ordering-1} when all $\xi_i(x)$ are distinct and $X\neq 0$.

Write $c_i \defeq \xi_i(x)$ for notational convenience.  Reordering if necessary, we may suppose that the $\xi_i$ are labelled so that $c_1 \leq \cdots \leq c_n$.  By \eqref{eq:mob2-evals-ordering-1} we obtain
\begin{align*}
  h_{c_i}^*\at{x}(X,X)
  &= (\pf{h}\at{x})^2 \Prod{j\neq i}{} (c_j-c_i) g_i(X,X) \\
  &= (\pf{h}\at{x})^2
    \underbrace{(c_1-c_i) \cdots (c_{i-1}-c_i)}_{\text{sign }(-1)^{i-1}} \,
    \underbrace{(c_{i+1}-c_i) \cdots (c_n-c_i)}_{> \, 0} \,
    \underbrace{g(X_i,X_i)}_{\geq \, 0}.
\end{align*}
If $g(X_i,X_i)=0$ then $c_i=t_i(x,X)$ is a root of $h_t^*\at{x}$, and likewise for $c_{i+1}$; otherwise $h_{c_i}^*\at{x}$ and $h_{c_{i+1}}^*\at{x}$ have opposite signs, so that there is a root $t_i(x,X)$ of $h_t^*\at{x}$ in the interval $(c_i,c_{i+1)}$.  Since $h_t^*(v,v)$ has degree $n-1$ and there $n-1$ such roots, we have accounted for all possible roots and \ref{lem:mob2-evals-ordering-1} follows.

\smallskip

\proofref{lem:mob2-evals-ordering}{2}
Suppose first that $\xi_i(x) < c < \xi_{i+1}(x)$.  Then $h_c^*\at{x}$ is not identically zero, so that $V_c = \setof{X\in T_xM}{h_c^*\at{x}(X,X)=0}$ is a \non zero\ quadric in $T_xM$.  Any such quadric has zero Lebesgue measure.

Now suppose that $c=\xi_i(x)$ or $c=\xi_{i+1}(x)$; without loss of generality we assume the former.  Let $k_i$ be the multiplicity of $\xi_i(x)$ as a root of $\pi(t)$.  Then every direct summand of \eqref{eq:mob2-evals-ordering-1} has a factor $(\xi_i(x)-t)^{k_i-1}$, so that
\begin{equation*}
  \hat{h}_t^*\at{x} = (\xi_i(x)-t)^{-k_i+1} h_t^*\at{x}
\end{equation*}
is a \non zero\ quadratic form which is polynomial of degree $n-k_i$ in $t$.  Moreover $\hat{h}_t^*\at{x}$ vanishes on $V_c$, so that $V_c$ is contained in a \non trivial\ quadric in $T_xM$.  As above, this quadric has zero Lebesgue measure. \end{proof}

\begin{thm} \thlabel{thm:mob2-evals-ordering} Suppose that $M$ is connected and consider a \riem\ pencil with linear metrics $h,\b{h}$ at $\infty,0$.  Let $\xi_1, \ldots, \xi_n$ denote the roots of $\pi(t) \defeq \pf{h_t}$.  Then:
\begin{enumerate}
  \item \label{thm:mob2-evals-ordering-1}
  $\xi_i(x) \leq \xi_{i+1}(y)$ for all $x,y\in M$; and
  
  \item \label{thm:mob2-evals-ordering-2}
  if $\xi_i(x)<\xi_{i+1}(x)$ at some $x\in M$, then $\xi_i(y)<\xi_{i+1}(y)$ at almost every $y\in M$.
\end{enumerate} \end{thm}

\begin{proof} \proofref{thm:mob2-evals-ordering}{1}
Take $x,y\in M$ and join them by a piecewise geodesic $\gamma:[0,1]\to M$ with $\gamma(0)=x$ and $\gamma(1)=y$.  As in \thref{lem:mob2-evals-ordering}, let $X\in T_xM$ and consider the $1$-parameter family $(\pf{h})^{-2}h_t^*(X,X)$ of integrals of the geodesic flow of $g\defeq (\pf{h})^{-1}h^{-1}$ and its ordered roots $t_1(x,X)\leq\cdots\leq t_{n-1}(x,X)$.  Since $(\pf{h})^{-2}h_t^*(X,X)$ is an integral of the geodesic flow of $g$, so are the $t_i(x,X)$.  It follows that $t_i(\gamma(0),\gamma'(0)) = t_i(\gamma(1),\gamma'(1))$, so that by \itemref{lem:mob2-evals-ordering}{1} we have
\begin{equation*}
  \xi_i(x) = \xi_i(\gamma(0)) \leq t_i(\gamma(0),\gamma'(0))
    = t_i(\gamma(1),\gamma'(1)) \leq \xi_{i+1}(\gamma(1)) = \xi_{i+1}(y)
\end{equation*}
as required.

\smallskip

\proofref{thm:mob2-evals-ordering}{2}
Suppose that $\xi_i(y) = \xi_{i+1}(y)$ holds for every $y$ in some open subset $U\subset M$.  Then item \ref{thm:mob2-evals-ordering-1} implies that
\begin{equation*}
  \xi_i(y_1) \leq \xi_{i+1}(y_2) = \xi_i(y_2) \leq \xi_{i+1}(y_1) = \xi_i(y_1)
\end{equation*}
for every $y_1,y_2\in U$, so that $\xi_i$ equals a constant $c\in\bR$ on $U$.

Now consider all possible geodesics $\gamma_{x,y}:[0,1]\to M$ joining $x$ to some $y\in U$.  Let $V_c \subseteq T_xM$ denote the set of their initial velocity vectors $\gamma_{x,y}'(0)$.  Then $t_i(\gamma(0),\gamma'(0)) = t_i(\gamma(1),\gamma'(1)) = c$ since $\gamma(1)\in U$, so that $V_c \subseteq \setof{X\in T_xM}{t_i(x,X)=0}$.  By \itemref{lem:mob2-evals-ordering}{2} this set has zero Lebesgue measure, and since $U\subset \exp V_c$ is follows that $U$ also has zero Lebesgue measure. \end{proof}

\appendix
\renewcommand{\chaptername}{appendix}
\renewcommand{\thesection}{\Alph{section}}
\renewcommand{\thetable}{\Alph{chapter}.\arabic{table}}

\chapter{Some algebraic identities} 
\label{c:app-alg}

\BufferDynkinLocaltrue
\renewcommand{\dynkinnameoffset}{-0.75}

This appendix contains some algebraic identities whose proofs would have disrupted the flow of the main thesis.  They are all formal consequences of the Jacobi identity and the brackets provided by Table \ref{tbl:ppg-alg-Z2gr}, and may be phrased either in terms of elements of the $\bZ^2$-graded Lie algebra $\fh$ or its associated graded bundle.  We assume throughout that a \ppg\ with parameters $(r,n)$ has been fixed, together with (algebraic) Weyl structures for both parabolics.

\begin{lem} \thlabel{lem:app-alg-musical} $\liebrac{\liebrac{X}{\alpha}}{Y}^{\flat} = \liebrac{\liebrac{X^{\flat}}{\alpha^{\sharp}}}{Y^{\flat}}$. \end{lem}

\begin{proof} Write $\pi \defeq (\det g)^{-1/r(n+1)} \in L$ and $h \defeq \pi^{-1} \ltens g^{-1}$, so that $X^{\flat} = -\liebrac{\pi^{-1}}{\liebrac{X}{h^{-1}}}$ and $\alpha^{\sharp} = \liebrac{\pi}{\liebrac{h}{\alpha}}$.  Then by the Jacobi identity we have
\begin{align} \label{eq:app-alg-musical-1}
  \liebrac{ X^{\flat} }{ \alpha^{\sharp} }
  &= \liebrac{ X^{\flat} }
             { \liebrac{\pi}{h(\alpha,\bdot)} } \notag \\
  &= \liebrac{ \pi X^{\flat} }
             { h(\alpha,\bdot) }
   + \liebrac{ \pi }
             { \liebrac{ X^{\flat} }
                       { h(\alpha,\bdot) } } \notag \\
  &= \liebrac{ h^{-1}(X,\bdot) }
             { h(\alpha,\bdot) }
     - \liebrac{ \pi }
               { \liebrac{ \liebrac{\pi^{-1}}{h^{-1}(X,\bdot)} }
                         { h(\alpha,\bdot) } } \notag \\
  &= \liebrac{ h^{-1}(X,\bdot) }
             { h(\alpha,\bdot) }
     - \liebrac{ \pi }
               { \liebrac{ \pi^{-1} }
                         { \liebrac{h^{-1}(X,\bdot)}{h(\alpha,\bdot)} } }.
\end{align}
Applying the Jacobi identity to the term $\liebrac{h^{-1}(X,\bdot)}{h(\alpha,\bdot)}$ in \eqref{eq:app-alg-musical-1} gives
\vspace{-0.1em}
\begin{equation*}
\begin{aligned}
  \liebrac{ h^{-1}(X,\bdot) }{ h(\alpha,\bdot) }
  &= \liebrac{ \liebrac{X}{h^{-1}} }{ \liebrac{h}{\alpha} } \\
  &= \liebrac{ \liebrac{ \liebrac{X}{h^{-1}} }{ h } }{ \alpha }
   + \liebrac{ h }{ \liebrac{ \liebrac{X}{h^{-1}} }{ \alpha } } \\
  &= \liebrac{X}{\alpha} + \liebrac{ h }{ \liebrac{ \liebrac{X}{\alpha} }{ h^{-1} } }.
\end{aligned}
\vspace{-0.1em}
\end{equation*}
Substituting the last display into \eqref{eq:app-alg-musical-1} then yields
\begin{align} \label{eq:app-alg-musical-2}
  \liebrac{ X^{\flat} }{ \alpha^{\sharp} }
  &= \liebrac{X}{\alpha} + \liebrac{ h }{ \liebrac{\liebrac{X}{\alpha}}{h^{-1}} } \notag \\
  &\qquad
     - \liebrac{ \pi }{ \liebrac{ \pi^{-1} }
                                { \liebrac{X}{\alpha} } }
     - \liebrac{ \pi }{ \liebrac{ \pi^{-1} }
                                { \liebrac{ h }
                                          { \liebrac{ \liebrac{X}{\alpha} }
                                                    { h^{-1} } } } } \notag \\
  &= \liebrac{X}{\alpha} + \liebrac{h}{\liebrac{\liebrac{X}{\alpha}}{h^{-1}}}
     - \liebrac{\pi}{\alpha(X)\pi^{-1}},
\end{align}
where the last term vanishes by virtue of the Jacobi identity and the fact that $\pi^{-1} \in L^*$ has trivial bracket with both $h \in L^*\tens B$ and $\liebrac{\liebrac{X}{\alpha}}{h^{-1}} \in L\tens B^*$.  Taking the Lie bracket with $Y^{\flat}$, equation \eqref{eq:app-alg-musical-2} yields
\begin{equation} \label{eq:app-alg-musical-3} \begin{split}
  \liebrac{ \liebrac{X^{\flat}}{\alpha^{\sharp}} }{ Y^{\flat} }
  &= \liebrac{ \liebrac{X}{\alpha} }{ Y^{\flat} }
     + \liebrac{ \liebrac{h}{\liebrac{\liebrac{X}{\alpha}}{h^{-1}}} }
               { Y^{\flat} } \\
  &\qquad
     - \alpha(X)\liebrac{ \liebrac{\pi}{\pi^{-1}} }
                        { Y^{\flat} }.
\end{split} \end{equation}
Writing $Y^{\flat} = \pi^{-1} \ltens h^{-1}(Y,\bdot)$ and using the Jacobi identity, the second term on the \rhs\ of \eqref{eq:app-alg-musical-3} equals
\begin{align*}
  \liebrac{ \liebrac{ h }
                    { \liebrac{ \liebrac{X}{\alpha} }
                              { h^{-1} } } }
          { \pi^{-1}h^{-1}(Y,\bdot) }
  &= \liebrac{ \liebrac{h}{\pi^{-1}h^{-1}(Y,\bdot)} }
             { \liebrac{ \liebrac{X}{\alpha} }
                       { h^{-1} } } \\
  &= \liebrac{ \pi^{-1}Y }
             { \liebrac{ \liebrac{X}{\alpha} }
                       { h^{-1} } } \\
  &= \liebrac{ \liebrac{\pi^{-1}}{Y} }
             { \liebrac{ \liebrac{X}{\alpha} }
                       { h^{-1} } }.
\end{align*}
Noting that $\liebrac{ \liebrac{X}{\alpha} }{ h^{-1} } \in L\tens B^*$ brackets trivially with $\pi^{-1} \in L^*$, this equals
\begin{align} \label{eq:app-alg-musical-4}
  & \liebrac{ \liebrac{\pi^{-1}}{Y} }
            { \liebrac{ \liebrac{X}{\alpha} }
                      { h^{-1} } } \notag \\
  &\quad
    = \liebrac{ \pi^{-1} }
              { \liebrac{ \liebrac{ Y }
                                  { \liebrac{X}{\alpha} } }
                        { h^{-1} } }
      + \liebrac{ \pi^{-1} }
                { \liebrac{ \liebrac{X}{\alpha} }
                          { \liebrac{Y}{h^{-1}} } } \notag \\
  &\quad
    = - \liebrac{ \pi^{-1} }
                { \liebrac{ \liebrac{\liebrac{X}{\alpha}}{Y} }
                          { h^{-1} } }
      + \liebrac{ \pi^{-1} }
                { \liebrac{ \liebrac{X}{\alpha} }
                          { h^{-1}(Y,\bdot) } } \notag \\
  &\quad
    = - \liebrac{ \pi^{-1} }
                { h^{-1}(\liebrac{\liebrac{X}{\alpha}}{Y},\bdot) }
      + \liebrac{ \alpha(X)\pi^{-1} }
                { h^{-1}(Y,\bdot) } \notag \\
  &\qquad\quad
     + \liebrac{ \liebrac{X}{\alpha} }
               { \liebrac{ \pi^{-1} }
                         { h^{-1}(Y,\bdot) } } \notag \\
  &\quad
    = \liebrac{ \liebrac{X}{\alpha} }{ Y }^{\flat} - \alpha(X)Y^{\flat}
      - \liebrac{ \liebrac{X}{\alpha} }{ Y^{\flat} }.
\end{align}
On the other hand, the third term on the \rhs\ of \eqref{eq:app-alg-musical-3} equals
\begin{equation*}
  -\alpha(X)\liebrac{\liebrac{\pi}{\pi^{-1}}}{Y^{\flat}}
    = -\alpha(X)\liebrac{\liebrac{\pi}{Y^{\flat}}}{\pi^{-1}}
    = \alpha(X)Y^{\flat}.
\end{equation*}
Substituting these expressions into the \rhs\ of \eqref{eq:app-alg-musical-3}, the first and third terms of \eqref{eq:app-alg-musical-3} cancel with the third and second terms of \eqref{eq:app-alg-musical-4} respectively, leaving only the desired term $\liebrac{ \liebrac{X}{\alpha} }{ Y }^{\flat}$. \end{proof}

\vfill
\begin{cor} \thlabel{cor:app-alg-skew} $\liebracw{\id}{g}_{X,Y} \in \fp^0$ is skew-adjoint \wrt\ $g$. \end{cor}

\begin{proof} Since the Killing form between $\fg/\fp$ and $\fp^{\perp}$ is simply contraction,
\begin{align*}
  g(\liebrac{\liebracw{\id}{g}_{X,Y}}{Z}, W)
    &= \killing{ \liebrac{\liebrac{Z}{Y^{\flat}}}{X} - \liebrac{\liebrac{Z}{X^{\flat}}}{Y} }
               { W^{\flat} } \\
    &= \killing{ Z }
               {   \liebrac{\liebrac{Y^{\flat}}{X}}{W^{\flat}} 
                 - \liebrac{\liebrac{X^{\flat}}{Y}}{W^{\flat}} }
\end{align*}
by invariance of $\killing{}{}$.  But $\liebrac{ \liebrac{Y^{\flat}}{X} }{ W^{\flat} } = \liebrac{ \liebrac{Y}{X^{\flat}} }{ W }^{\flat}$ by \threfit{lem:app-alg-musical}, from which the result easily follows. \end{proof}

\begin{lem} \thlabel{lem:app-alg-flat} Suppose that $A \in \fp^0$ is \self adjoint\ \wrt\ a \non degenerate\ element $g \in B^*$, and let $\flat = \sharp^{-1} : \fg/\fp \to \fp^{\perp}$ be the musical isomorphisms of $g$.  Then:
\begin{enumerate}
  \item \label{lem:app-alg-flat-noA}
  $\ve^i( \liebrac{ \liebrac{X}{e_i^{\flat}} }{ Y } ) = (2-r)g(X,Y)$; and
  
  \item \label{lem:app-alg-flat-A}
  $\ve^i( A \liebrac{ \liebrac{X}{e_i^{\flat}} }{ Y } ) = (2-r)g(AX,Y)$
\end{enumerate}
for all $X,Y \in \fg/\fp$. \end{lem}

\medskip
\begin{proof} \proofref{lem:app-alg-flat}{noA} Let $h \defeq \pi^{-1} \ltens \ltens g^{-1}$ be the corresponding element of $L^*\tens B$, where $\pi \defeq (\det g)^{-1/r(n+1)} \in L$.  Then we may write $e_i^{\flat} = h^{-1}(\pi^{-1} \ltens e_i, \bdot) = \liebrac{ h^{-1} }{ \liebrac{\pi^{-1}}{e_i} }$, for which Table \ref{tbl:ppg-alg-Z2gr} and successive applications of the Jacobi identity yield
\begin{equation*} \begin{aligned}
  \liebrac{ \liebrac{X}{e_i^{\flat}} }{ Y }
    &= \liebrac{ \liebrac{ X }
                         { \liebrac{h^{-1}}{\pi^{-1} \ltens e_i} } }
               { Y } \\
    &= \liebrac{ \liebrac{ h^{-1}(X,\bdot) }
                         { \pi^{-1} \ltens e_i } }
               { Y }
     + \liebrac{ \liebrac{ h^{-1} }
                         { \liebrac{X}{\pi^{-1} \ltens e_i} } }
               { Y } \\
    &= -\liebrac{ h^{-1}(X,Y) }{ \pi^{-1} \ltens e_i }
     + \liebrac{ h^{-1}(X,\bdot) }
               { \liebrac{\pi^{-1} \ltens e_i}{Y} } \\
    &\qquad
     - \liebrac{ h^{-1}(Y,\bdot) }
               { \liebrac{e_i}{\pi^{-1} \ltens X} }
\end{aligned}
\end{equation*}
Now $\liebrac{ \pi^{-1} \ltens e_i }{ Y } = \liebrac{ \liebrac{\pi^{-1}}{e_i} }{ Y } = \liebrac{ \pi^{-1} \ltens Y }{ e_i }$ since $\fg/\fp$ is abelian, so we obtain
\begin{align*}
  &\liebrac{ \liebrac{X}{e_i^{\flat}} }{ Y } \\
  &\quad \begin{aligned}
    &= -\pi^{-1} \ltens h^{-1}(X,Y)e_i
     + \liebrac{ h^{-1}(X,\bdot) }
               { \liebrac{\pi^{-1} \ltens Y}{e_i} }
     - \liebrac{ h^{-1}(Y,\bdot) }
               { \liebrac{e_i}{\pi^{-1} \ltens X} } \\
    &= -g(X,Y)e_i \\
    &\qquad
     + \liebrac{ \liebrac{ \liebrac{h^{-1}(X,\bdot)}{\pi^{-1}} }
                         { Y } }
               { e_i }
     - \liebrac{ \liebrac{ \pi^{-1} }
                         { h^{-1}(X,Y) } }
               { e_i }
     - \liebrac{ \pi^{-1} \ltens Y }
               { h^{-1}(X,e_i) } \\
    &\qquad
     + \liebrac{ \liebrac{ \liebrac{h^{-1}(Y,\bdot)}{\pi^{-1}} }
                         { X } }
               { e_i }
     - \liebrac{ \liebrac{ \pi^{-1} }
                         { h^{-1}(Y,X) } }
               { e_i }
     - \liebrac{ \pi^{-1} \ltens X }
               { h^{-1}(Y,e_i) } \\
    &= -g(X,Y)e_i
     + \liebrac{ \liebrac{X^{\flat}}{Y} }
               { e_i }
     - \liebrac{ \pi^{-1} \ltens e_i }
               { h^{-1}(X,Y) }
     + g(X,e_i)Y \\
    &\qquad
     + \liebrac{ \liebrac{Y^{\flat}}{X} }
               { e_i }
     - \liebrac{ \pi^{-1} \ltens e_i }
               { h^{-1}(Y,X) }
     + g(Y,e_i)X \\
    &= -g(X,Y)e_i
     - \liebrac{ \liebrac{Y}{X^{\flat}} }
               { e_i }
     + g(X,Y)e_i + g(X,e_i)Y \\
    &\qquad
     - \liebrac{ \liebrac{X}{Y^{\flat}} }
               { e_i }
     + g(Y,X)e_i + g(X,e_i)Y \\
    &= g(X,Y)e_i + g(X,e_i)Y + g(Y,e_i)X
     - \liebrac{ \liebrac{X}{Y^{\flat}} }{ e_i }
     - \liebrac{ \liebrac{Y}{X^{\flat}} }{ e_i }.
\end{aligned} \end{align*}
Using \thref{cor:ppg-alg-trace}, evaluating on $\ve^i$ and summing over $i$ now yields
\begin{equation*} \begin{aligned}
  \ve^i \big( \liebrac{ \liebrac{X}{e_i^{\flat}} }{ Y } \big)
    &= -rn g(X,Y) + g(X,Y) + g(Y,X) \\
    &\qquad
      - \tfrac{1}{2}r(n+1) g(Y,X) - \tfrac{1}{2}r(n+1) g(X,Y) \\
    &= (2-r) g(X,Y)
\end{aligned}
\end{equation*}
as claimed.

\smallskip

\proofref{lem:app-alg-flat}{A} Using the Gram-Schimdt algorithm, we may suppose that $\{e_i\}_i$ is a $g$-orthonormal frame.  Then by the Jacobi identity and \ref{lem:app-alg-flat-noA},
\begin{align} \label{eq:app-alg-flat-A-1}
  \ve^i \big( A\liebrac{ \liebrac{X}{e_i^{\flat}} }{ Y } \big)
    &= \ve^i \big( \liebrac{ \liebrac{AX}{e_i^{\flat}} }{ Y }
     + \liebrac{ \liebrac{X}{e_i^{\flat}} }
               { AY } \notag \\
    &\qquad
     + \liebrac{ \liebrac{X}{\liebrac{A}{e_i^{\flat}}} }
               { Y } \big) \notag \\
  &= 2(2-r) g(AX,Y)
   + \killing{ \ve^i }
             { \liebrac{ \liebrac{ X }
                                 { \liebrac{A}{e_i^{\flat}} } }
                       { Y } }.
\end{align}
However $g$-orthonormality of the $e_i$ implies that $e_i^{\flat} = \ve^i$, so this becomes
\begin{align*}
  \killing{ \ve^i }
          { \liebrac{ \liebrac{ X }
                              { \liebrac{A}{e_i^{\flat}} } }
                    { Y } }
    &= \killing{ e_i^{\flat} }
               { \liebrac{ \liebrac{ X }
                                   { \liebrac{A}{\ve^i} } }
                         { Y } } \\
    &= \killing{ \liebrac{A}{\ve^i} }
               { \liebrac{ \liebrac{ X }
                                   { \liebrac{e_i^{\flat}} } }
                         { Y } } \\
    &= -\ve^i \big( A \liebrac{ \liebrac{X}{e_i^{\flat}} }
                              { Y } \big),
\end{align*}
which is minus the \lhs\ of \eqref{eq:app-alg-flat-A-1}; rearranging now completes the proof. \end{proof}

\vspace{0.6em}
\begin{lem} \thlabel{lem:app-alg-trA} Let $A \in \fp^0$ be \self adjoint\ \wrt\ a \non degenerate\ element $g \in B^*$.  Then $\liebrac{Ae_i}{\ve^i} = \tfrac{1}{2}((\tr A)\id + rA)$. \end{lem}

\begin{proof} Let $h \defeq \pi^{-1} \ltens g^{-1}$ be the corresponding element of $L^* \tens B$, where we write $\pi \defeq (\det g)^{1/r(n+1)}$ for brevity.  Then $g(\liebrac{ \liebrac{Ae_i}{\ve^i} }{ X }, Y) = h^{-1}( \liebrac{\liebrac{Ae_i}{\ve^i}}{X}, \pi^{-1}\ltens Y)$, for which the Jacobi identity and Table \ref{tbl:ppg-alg-Z2gr} yield
\begin{align} \label{eq:app-alg-trA-1}
  &h^{-1}( \liebrac{\liebrac{Ae_i}{\ve^i}}{X}, \pi^{-1}\ltens Y)
      \notag \\
    &\quad =  \killing{ \liebrac{ \liebrac{ \liebrac{Ae_i}{\ve^i} }
                                { X } }
                      { h^{-1} } }
            { \liebrac{\pi^{-1}}{Y} }
      \notag \\
    &\quad = -\killing{ \liebrac{ \pi^{-1} }
                                { \liebrac{ \liebrac{ \liebrac{Ae_i}{\ve^i} }
                                          { X } }
                                { h^{-1} } } }
                      { Y }
      \notag \\
    &\quad = -\killing{ \liebrac{ \liebrac{ \pi^{-1} \ltens Ae_i }
                                { \ve^i } }
                      { X }
            + \liebrac{ \liebrac{Ae_i}{\ve^i} }
                      { \pi^{-1} \ltens X } }
            { \liebrac{ h^{-1} }{ Y } }
      \notag \\
    &\quad =  \killing{ \liebrac{ \ve^i(\pi^{-1} \ltens Ae_i) }
                      { X }
            + \liebrac{ \liebrac{Ae_i}{\pi^{-1} \ltens X} }
                      { \ve^i }
            + \liebrac{ Ae_i }
                      { \liebrac{\ve^i}{\pi^{-1}\ltens X} } }
            { h^{-1}(Y,\bdot) }
      \notag \\
    &\quad =  \killing{ \ve^i(\pi^{-1} \ltens Ae_i)X + \ve^i(\pi^{-1} \ltens X)Ae_i }
                      { h^{-1}(Y,\bdot) }
      \notag \\
    &\quad\qquad
      + \killing{ \liebrac{ \liebrac{Ae_i}{\pi^{-1} \ltens X} }
                          { \ve^i } }
                { h^{-1}(Y,\bdot) }
      \notag \\
    &\quad =  (\tr A)g(X,Y) + g(AX,Y)
            + \killing{ \liebrac{ \liebrac{Ae_i}{\pi^{-1} \ltens X} }
                                { \ve^i } }
                      { h^{-1}(Y,\bdot) }.
\end{align}
We evaluate the final term on the \rhs\ of \eqref{eq:app-alg-trA-1} separately.  We have
\begin{align} \label{eq:app-alg-trA-2}
  &\killing{ \liebrac{ \liebrac{Ae_i}{\pi^{-1} \ltens X} }
                    { \ve^i } }
          { h^{-1}(Y,\bdot) }
      \notag \\
    &\quad =  \killing{ \liebrac{Ae_i}{\pi^{-1} \ltens X} }
                      { \liebrac{ \liebrac{\ve^i}{Y} }
                                { h^{-1} } }
      \notag \\
    &\quad = -\killing{ \liebrac{ h^{-1} }
                                { \liebrac{Ae_i}{\pi^{-1} \ltens X} } }
                      { \liebrac{Y}{\ve^i} }
      \notag \\
    &\quad =  \killing{ \liebrac{ h^{-1}(Ae_i,\bdot) }
                                { \pi^{-1} \ltens X } }
                      { \liebrac{Y}{\ve^i} }
            - \killing{ \liebrac{Ae_i}{X^{\flat}} }
                      { \liebrac{Y}{\ve^i} }
      \notag \\
    &\quad = -\killing{ \liebrac{ \pi^{-1} }
                                { h^{-1}(Ae_i,X) }
                      - \liebrac{(Ae_i)^{\flat}}{X} }
                      { \liebrac{Y}{\ve^i} }
            - \killing{ \liebrac{ X^{\flat} }
                                { \liebrac{Ae_i}{\ve^i} } }
                      { Y }
      \notag \displaybreak \\
    &\quad =  \killing{ h^{-1}(Ae_i,X) }
                      { \ve^i(Y) \ltens \pi^{-1} }
            - \killing{ \liebrac{(Ae_i)^{\flat}}{X} }
                      { \liebrac{Y}{\ve^i} }
            - \killing{ \liebrac{ X^{\flat} }
                                { \liebrac{Ae_i}{\ve^i} } }
                      { Y }
      \notag \\
    &\quad =  g(AX,Y)
            - \killing{ \liebrac{(Ae_i)^{\flat}}{X} }
                      { \liebrac{Y}{\ve^i} }
            - \killing{ \liebrac{ X^{\flat} }
                                { \liebrac{Ae_i}{\ve^i} } }
                      { Y }.
\end{align}
We have $(Ae_i)^{\flat} = -\liebrac{A}{e_i^{\flat}}$ since $A$ is $g$-\self adjoint, so the second term on the \rhs\ of \eqref{eq:app-alg-trA-2} equals
\begin{align} \label{eq:app-alg-trA-3}
  -\killing{ \liebrac{(Ae_i)^{\flat}}{X} }
           { \liebrac{Y}{\ve^i} }
    &=  \killing{ \liebrac{ \liebrac{A}{e_i^{\flat}} }
                          { X } }
                { \liebrac{Y}{\ve^i} }
      \notag \\
    &=  \killing{ \liebrac{ A }
                          { \liebrac{e_i^{\flat}}{X} }
                - \liebrac{AX}{e_i^{\flat}} }
                { \liebrac{Y}{\ve^i} }
      \notag \\
    &= -\killing{ \liebrac{e_i^{\flat}}{X} }
                { \liebrac{AY}{\ve^i}
                + \liebrac{ Y }
                          { \liebrac{A}{\ve^i} } }
      - \killing{ \liebrac{ \liebrac{AX}{e_i^{\flat}} }
                          { Y } }
                { \ve^i }
      \notag \\
    &=  \ve^i( \liebrac{ \liebrac{X}{e_i^{\flat}} }
                       { AY } )
      - \killing{ \liebrac{e_i^{\flat}}{X} }
                { \liebrac{ Y }
                          { \liebrac{A}{\ve^i} } }
      - \ve^i( \liebrac{ \liebrac{AX}{e_i^{\flat}} }
                       { Y } )
      \notag \\
    &=  (2-r)g(X,AY)
      - \killing{ \liebrac{e_i^{\flat}}{X} }
                { \liebrac{ Y }
                          { \liebrac{A}{\ve^i} } }
      - (2-r)g(AX,Y)
      \notag \\
    &= -\ve^i( A \liebrac{ \liebrac{X}{e_i^{\flat}} }
                         { Y } )
      \notag \\
    &= -(2-r)g(AX,Y)
\end{align}
by \threfit{lem:app-alg-flat}\ref{lem:app-alg-flat-A}.  The third term on the \rhs\ of \eqref{eq:app-alg-trA-2} equals
\begin{align*}
  \killing{ \liebrac{ X^{\flat} }
                    { \liebrac{Ae_i}{\ve^i} } }
          { Y }
    &=  \killing{ \liebrac{ X^{\flat} }
                          { \liebrac{Ae_i}{\ve^i} }^{\sharp} }
                { Y^{\flat} } \\
    &=  \killing{ \liebrac{ X }
                          { \liebrac{(Ae_i)^{\flat}}{\ve^{i\sharp}} } }
                { Y^{\flat} }
\end{align*}
by \threfit{lem:app-alg-flat}.  We may as well assume that $\{e_i\}_i$ is a $g$-orthonormal frame, so that $e_i^{\flat} = \ve^i$.  Then the last display becomes
\begin{align} \label{eq:app-alg-trA-4}
  \killing{ \liebrac{ X^{\flat} }
                    { \liebrac{Ae_i}{\ve^i} } }
          { Y }
    &= -\killing{ \liebrac{ X }
                          { \liebrac{ \liebrac{A}{\ve^i} }
                                    { e_i } } }
                { Y^{\flat} }
      \notag \\
    &= -\killing{ \liebrac{ A }
                          { \liebrac{ X }
                                    { \liebrac{\ve^i}{e_i} } }
                - \liebrac{ AX }
                          { \liebrac{\ve^i}{e_i} }
                - \liebrac{ X }
                          { \liebrac{\ve^i}{Ae_i} } }
                { Y^{\flat} }
      \notag \\
    &= -\killing{ \tfrac{1}{2}r(n+1)AX - \tfrac{1}{2}r(n+1)AX
                + \liebrac{ \liebrac{Ae_i}{\ve^i} }
                          { X } }
                { Y^{\flat} }
      \notag \\
    &= g( \liebrac{ \liebrac{Ae_i}{\ve^i} }
                  { X }, Y ).
\end{align}
Substituting \eqref{eq:app-alg-trA-2}, \eqref{eq:app-alg-trA-3} and \eqref{eq:app-alg-trA-4} into \eqref{eq:app-alg-trA-1} now yields
\begin{align*}
  g(\liebrac{ \liebrac{Ae_i}{\ve^i} }{ X }, Y)
    &=  (\tr A)g(X,Y) + g(AX,Y) + g(AX,Y) \\
    &\qquad
      - (2-r)g(AX,Y) - g(\liebrac{ \liebrac{Ae_i}{\ve^i} }{ X }, Y) \\
    &=  (\tr A)g(X,Y) + rg(AX,Y),
\end{align*}
from which the result follows. \end{proof}

\chapter{Tables} 
\label{c:app-tbl}

\BufferDynkinLocaltrue
\renewcommand{\dynkinnameoffset}{-0.75}

For the convenience of the reader, we collect together some representation-theoretic information about the flat models of both the ``big'' and ``small'' \Rspaces\ of a \ppg.

\vfill

\begin{table}[h]
\newcommand{\e}[3][0em]{ \begin{gathered} #2 \\[#1] \leq#3 \end{gathered} }
\newcommand{\tpad}[1][-1.3ex]{ && \\[#1] }
\newcommand{\bpad}[1][-1.1ex]{ \\[#1] && }
  \centering
  \begin{tabular}{| c | *{2}{>{$}c<{$}|} }
    \hline \tpad[-2ex]
    \bfseries{Type}
      & H\acts\fq
      & G\acts\fp
    \bpad[-2.5ex] \\ \hline \tpad[-2.8ex]
    \type{C}{n+1}
      & \e{ \dynkinCp{2}{0}{3}{tables-rspace-Cq} }
          { \alg{sp}{2n+2,\bC} }
      & \e{ \dynkinAp{2}{0}{2}{tables-rspace-Cp} }
          { \alg{sl}{n+1,\bC} }
    \bpad \\ \hline \tpad[-2.8ex]
    \type{A}{2n+1}
      & \e{ \dynkinApk{2}{2}{tables-rspace-Aq} }
          { \alg{sl}{2n+2,\bC} }
      & \e[-1.6ex]{ \dynkinAAp{2}{0}{2}{tales-rspace-Ap} }
                  { \alg{sl}{n+1,\bC}\dsum\alg{sl}{n+1,\bC} }
    \bpad \\ \hline \tpad
    \type{D}{2n+2}
      & \e{ \dynkinDp{3}{0}{2}{tables-rspace-Dq} }
          { \alg{so}{4n+4,\bC} }
      & \e{ \dynkinApp{3}{0}{3}{tables-rspace-Dp} }
          { \alg{sl}{2n+2,\bC} }
    \bpad \\ \hline \tpad[-2.8ex]
    \type{E}{7}
      & \e{ \dynkinEp{7}{tables-rspace-Eq} }
          { \alg[_7]{e}{\bC} }
      & \e{ \dynkinEp{6}{tables-rspace-Ep} }
          { \alg[_6]{e}{\bC} }
    \bpad \\ \hline \tpad
    \type{BD}{n+4}
      & \e[-1.0ex]{ \dynkinBDp{3}{0}{2}{tables-rspace-BDq} }
                  { \alg{so}{n+4,\bC} }
      & \e[-1.0ex]{ \dynkinBDp{2}{0}{2}{tables-rspace-BDp} }
                  { \alg{so}{n+2,\bC} }
    \bpad \\ \hline
  \end{tabular}
  \caption[The symmetric \Rspaces\ $H\acts\fq$ and $G\acts\fp$]
          {The symmetric \Rspace s\ $H\acts\fq$ and $G\acts\fp$ associated to the irreducible complex \ppg.}
  \label{tbl:app-tbl-rspace}
\end{table}
\pagebreak

\begin{table}[p]
\newcommand{\e}[3][0em]{ \begin{gathered} #2 \\[#1] \leq #3 \end{gathered} }
\newcommand{\tpad}[1][-1.3ex]{ &&& \\[#1] }
\newcommand{\bpad}[1][-1.1ex]{ \\[#1] &&& }
  \centering
  \begin{tabular}{| c | *{2}{>{$}c<{$}|} m{6.25em} |}
    \hline \tpad[-2.0ex] 
    \bf{Type}
      & \fq\leq\fh
      & \fp\leq\fg
      & \bf{Geometry}
    \bpad[-2.5ex] \\ \hline \tpad[-3.0ex]
    \type{C}{n+1} &
      \e{ \dynkinCp{2}{0}{3}{tables-real-Ch} }
        { \alg{sp}{2n+2,\bR} } &
      \e{ \dynkinSLRp{2}{0}{2}{tables-real-Cg} }
        { \alg{sl}{n+1,\bR} } &
      \Proj\ \newline geometry
    \bpad \\ \hline \tpad[-3.0ex]
    \type{A}{2n+1} &
      \e[-1.6ex]{ \dynkinASUp{2}{0}{2}{tables-real-A1h} }
                { \alg{sl}{2n+2,\bR} } &
      \e[-1.6ex]{ \dynkinAAp{2}{0}{2}{tables-real-A1g} }
                { \alg{sl}{n+1,\bR}\dsum\alg{sl}{n+1,\bR} } &
      ???
    \bpad \\ \cline{2-4} \tpad[-3.0ex]
    ~ &
      \e[-1.6ex]{ \dynkinSUp{2}{0}{2}{tables-real-A2h} }
                { \alg{su}{n+1,n+1} } &
      \e[-1.6ex]{ \dynkinSLCp{2}{0}{2}{tables-real-A2g} }
                { \alg{sl}{n+1,\bC} } &
      \Cproj\ \newline geometry
    \bpad \\ \hline \tpad
    \type{D}{2n+2} &
      \e{ \dynkinDp{3}{0}{2}{tables-real-D1h} }
        { \alg{so}{4n+4,\bR} } &
      \e{ \dynkinApp{3}{0}{3}{tables-real-D1g} }
        { \alg{sl}{2n+2,\bR} } &
      \Grassmannian\ \newline of $2$-planes
    \bpad \\ \cline{2-4} \tpad
    ~ &
      \e{ \dynkinSOsp{3}{0}{2}{tables-real-D2h} }
        { \alg[^*]{so}{4n+4} } &
      \e{ \dynkinSLHp{3}{0}{3}{tables-real-D2g} }
        { \alg{sl}{2n+2,\bR} } &
      Almost \newline \qtn ic\ \newline geometry
    \bpad \\ \hline \tpad[-3.0ex]
    \type{E}{7} &
      \e{ \dynkinEp{7}{tables-real-E1h} }
        { \alg[_{7(7)}]{e}{} = \erealform{V} } &
      \e{ \dynkinEp{6}{tables-real-E1g} }
        { \alg[_{6(6)}]{e}{} = \erealform{I} } &
      ???
    \bpad \\ \cline{2-4} \tpad[-3.0ex]
    ~ &
      \e{ \dynkinEVIIp{tables-real-E2h} }
        { \alg[_{7(-25)}]{e}{} = \erealform{VII} } &
      \e{ \dynkinEIVp{tables-real-E2g} }
        { \alg[_{6(-26)}]{e}{} = \erealform{IV} } &
      Cayley plane
    \bpad \\ \hline \tpad
    \type{BD}{n+4} &
      \e[-1.0ex]{ \dynkinSOpqshortp{3}{2}{tables-real-BDh} }
                { \alg{so}{p+2,q+2} } &
      \e[-1.0ex]{ \dynkinSOpqshortp{2}{2}{tables-real-BDg} }
                { \alg{so}{p+1,q+1} } &
      Conformal \newline geometry \newline of signature \newline $(p,q)$
    \bpad \\ \hline
  \end{tabular}
  \caption[Real forms of \ppgs]
          {Real forms of the irreducible \ppgs.}
  \label{tbl:app-tbl-real}
\end{table}

\begin{landscape}
  \begin{table}[p]
  \newcommand{\e}[2]{ \begin{array}{c} #1 \ifx \relax#2\relax \\ ~ \else \\ =#2 \fi \end{array} }
  \newcommand{\tpad}[1][-1.1ex]{ &&&& \\[#1] }
  \newcommand{\bpad}[1][-0.9ex]{ \\[#1] &&&& }
    \centering
    \begin{tabular}{| c | *{4}{>{$}c<{$}|} }
      \hline \tpad[-2ex]
      \text{\bf Type}
        & \bW
        & L^*\tens B
        & L^*\tens\fg/\fp
        & L^*
      \bpad[-2.6ex] \\ \hline \tpad
      \type{C}{n+1}
        & \e{ \dynkinA[2,0,0,0]{2}{0}{2}{tables-W-C} }
            { \Symm[\bC]{2}\bC^{n+1} }
        & \e{ \dynkinAp[2,0,0,0]{2}{0}{2}{tables-W-C1} }
            {}
        & \e{ \dynkinAp[1,0,0,-1]{2}{0}{2}{tables-W-C2} }
            {}
        & \e{ \dynkinAp[0,0,0,-2]{2}{0}{2}{tables-W-C3} }
            {}
      \bpad \\ \hline \tpad
      \type{A}{2n+1}
        & \e{ \dynkinAA[1,0,0,0, 1,0,0,0]{2}{0}{2}{tables-W-A} }
            { \bC^{n+1}\etens\conj{\bC^{n+1}} }
        & \e{ \dynkinAAp[1,0,0,0, 1,0,0,0]{2}{0}{2}{tables-W-A1} }
            {}
        & \e{ \dynkinAAp[1,0,0,0, 0,0,0,-1]{2}{0}{2}{tables-W-A2} \dsum \cpxconj }
            {}
        & \e{ \dynkinAAp[0,0,0,-1, 0,0,0,-1]{2}{0}{2}{tables-W-A3} }
            {}
      \bpad \\ \hline \tpad
      \type{D}{2n+2}
        & \e{ \dynkinA[0,1,0,0,0,0]{3}{0}{3}{tables-W-D} }
            { \Wedge[\bC]{2}{\bC^{2n+2}} }
        & \e{ \dynkinApp[0,1,0,0,0,0]{3}{0}{3}{tables-W-D1} }
            {}
        & \e{ \dynkinApp[1,0,0,0,-1,1]{3}{0}{3}{tables-W-D2} }
            {}
        & \e{ \dynkinApp[0,0,0,0,-1,0]{3}{0}{3}{tables-W-D3} }
            {}
      \bpad \\ \hline \tpad
      \type{E}{7}
        & \e{ \dynkinE[1,0,0,0,0,0]{6}{tables-W-E} }
            { \erepn[\bC] }
        & \e{ \dynkinEp[1,0,0,0,0,0]{6}{tables-W-E1} }
            {}
        & \e{ \dynkinEp[0,0,0,0,-1,1]{6}{tables-W-E2} }
            {}
        & \e{ \dynkinEp[0,0,0,0,-1,0]{6}{tables-W-E3} }
            {}
      \bpad \\ \hline \tpad
      \type{BD}{n+4}
        & \e{ \dynkinBD[1,0,0,0,0,0,0,0]{2}{0}{2}{tables-W-BD} }
            { \bC^{n+2} }
        & \e{ \dynkinBDp[1,0,0,0,0,0,0,0]{2}{0}{2}{tables-W-BD1} }
            {}
        & \e{ \dynkinBDp[-1,1,0,0,0,0,0,0,0]{2}{0}{3}{tables-W-BD2} }
            {}
        & \e{ \dynkinBDp[-1,0,0,0,0,0,0,0]{2}{0}{2}{tables-W-BD3} }
            {}
      \bpad \\ \hline
    \end{tabular}
    \caption[$\bW$ and its graded components]
            {The infinitesimal isotropy representation $\bW \defeq \fh/\fq$ and its graded components.}
    \label{tbl:app-tbl-W}
  \end{table}
\end{landscape}

\begin{landscape}
  \begin{table}[p]
  \newcommand{\e}[2]{ \begin{array}{c} #1 \ifx \relax#2\relax \\ ~ \else \\ =#2 \fi \end{array} }
  \newcommand{\tpad}[1][-1.1ex]{ &&&& \\[#1] }
  \newcommand{\bpad}[1][-0.9ex]{ \\[#1] &&&& }
    \centering
    \begin{tabular}{| c | *{4}{>{$}c<{$}|} }
      \hline \tpad[-2ex]
      \text{\bf Type}
        & \bW^*
        & L
        & L\tens\fp^{\perp}
        & L\tens B^*
      \bpad[-2.6ex] \\ \hline \tpad
      \type{C}{n+1}
        & \e{ \dynkinA[0,0,0,2]{2}{0}{2}{tables-Wd-C} }
            { \Symm[\bC]{2}\bC^{n+1*} }
        & \e{ \dynkinAp[0,0,0,2]{2}{0}{2}{tables-Wd-C1} }
            {}
        & \e{ \dynkinAp[0,0,0,1,0]{2}{0}{3}{tables-Wd-C2} }
            {}
        & \e{ \dynkinAp[0,0,0,2,-2]{2}{0}{3}{tables-Wd-C3} }
            {}
      \bpad \\ \hline \tpad
      \type{A}{2n+1}
        & \e{ \dynkinAA[0,0,0,1, 0,0,0,1]{2}{0}{2}{tables-Wd-A} }
            { \bC^{n+1*}\etens\conj{\bC^{n+1*}} }
        & \e{ \dynkinAAp[0,0,0,1, 0,0,0,1]{2}{0}{2}{tables-Wd-A1} }
            {}
        & \e{ \dynkinAAp[0,0,0,1,-1, 0,0,0,0,-1]{2}{0}{3}{tables-Wd-A2} \dsum \cpxconj }
            {}
        & \e{ \dynkinAAp[0,0,0,1,-1, 0,0,0,1,-1]{2}{0}{3}{tables-Wd-A3} }
            {}
      \bpad \\ \hline \tpad
      \type{D}{2n+2}
        & \e{ \dynkinA[0,0,0,0,1,0]{3}{0}{3}{tables-Wd-D} }
            { \Wedge[\bC]{2}{\bC^{2n+2*}} }
        & \e{ \dynkinApp[0,0,0,0,1,0]{3}{0}{3}{tables-Wd-D1} }
            {}
        & \e{ \dynkinApp[0,0,0,1,-1,1]{2}{0}{4}{tables-Wd-D2} }
            {}
        & \e{ \dynkinApp[0,0,0,1,0,-1,0]{2}{0}{5}{tables-Wd-D3} }
            {}
      \bpad \\ \hline \tpad
      \type{E}{7}
        & \e{ \dynkinE[0,0,0,0,1,0]{6}{tables-Wd-E} }
            { \erepn[\bC^*] }
        & \e{ \dynkinEp[0,0,0,0,1,0]{6}{tables-Wd-E1} }
            {}
        & \e{ \dynkinEp[0,0,0,0,-1,1]{6}{tables-Wd-E2} }
            {}
        & \e{ \dynkinEp[1,0,0,0,-1,0]{6}{tables-Wd-E3} }
            {}
      \bpad \\ \hline \tpad
      \type{BD}{n+4}
        & \e{ \dynkinBD[1,0,0,0,0,0,0,0]{2}{0}{2}{tables-Wd-BD} }
            { \bC^{n+2} }
        & \e{ \dynkinBDp[1,0,0,0,0,0,0,0]{2}{0}{2}{tables-Wd-BD1} }
            {}
        & \e{ \dynkinBDp[-1,1,0,0,0,0,0,0,0]{3}{0}{2}{tables-Wd-BD2} }
            {}
        & \e{ \dynkinBDp[-1,0,0,0,0,0,0,0]{2}{0}{2}{tables-Wd-BD3} }
            {}
      \bpad \\ \hline
    \end{tabular}
    \caption[$\bW^*$ and its graded components]
            {The representation $\bW^* \isom \fq^{\perp}$ and its graded components.}
    \label{tbl:app-tbl-Wd}
  \end{table}
\end{landscape}

\begin{landscape}
  \begin{table}[p]
  \newcommand{\e}[2]{ \begin{array}{c} #1 \ifx \relax#2\relax \\ ~ \else \\ =#2 \fi \end{array} }
  \newcommand{\tpad}[1][-1.1ex]{ &&&& \\[#1] }
  \newcommand{\bpad}[1][-0.9ex]{ \\[#1] &&&& }
    \centering
    \begin{tabular}{| c | *{4}{>{$}c<{$}|} }
      \hline \tpad[-2ex]
      \text{\bf Type}
        & \bU^*
        & L^2\tens B^*
        & L^2\tens\fp^{\perp}\wedge B^*
        & L^2\tens B^*\wedge B^*
      \bpad[-2.6ex] \\ \hline \tpad
      \type{C}{n+1}
        & \e{ \dynkinA[0,0,0,2,0]{2}{0}{3}{tables-Ud-C} }
            { \Cartan{2}\Wedge{2}\bC^{n+1*} }
        & \e{ \dynkinAp[0,0,0,2,0]{2}{0}{3}{tables-Ud-C1} }
            {}
        & \e{ \dynkinAp[0,0,0,1,1,-1]{2}{0}{4}{tables-Ud-C2} }
            {}
        & \e{ \dynkinAp[0,0,0,2,0,-2]{2}{0}{4}{tables-Ud-C3} }
            {}
      \bpad \\ \hline \tpad
      \type{A}{2n+1}
        & \e{ \dynkinAA[0,0,0,1,0, 0,0,0,1,0]{2}{0}{3}{tables-Ud-A} }
            { \Wedge{2}\bC^{n+1*}\etens\Wedge{2}\conj{\bC^{n+1*}} }
        & \e{ \dynkinAAp[0,0,0,1,0, 0,0,0,1,0]{2}{0}{3}{tables-Ud-A1} }
            {}
        & \e{ \dynkinAAp[0,0,0,1,0,-1, 0,0,0,0,1,0]{2}{0}{4}{tables-Ud-A2} \dsum \cpxconj }
            {}
        & \e{ \dynkinAAp[0,0,0,1,0,-1, 0,0,0,1,0,-1]{2}{0}{4}{tables-Ud-A3} }
            {}
      \bpad \\ \hline \tpad
      \type{D}{2n+2}
        & \e{ \dynkinA[0,0,1,0,0,0,0]{1}{0}{5}{tables-Ud-D} }
            { \Wedge[\bC]{4}{\bC^{2n+2*}} }
        & \e{ \dynkinApp[0,0,1,0,0,0,0]{1}{0}{5}{tables-Ud-D1} }
            {}
        & \e{ \dynkinApp[0,0,1,0,0,-1,1]{1}{0}{6}{tables-Ud-D2} }
            {}
        & \e{ \dynkinApp[0,0,1,0,0,0,-1,0]{1}{0}{7}{tables-Ud-D3} }
            {}
      \bpad \\ \hline \tpad
      \type{E}{7}
        & \e{ \dynkinE[1,0,0,0,0,0]{6}{tables-Ud-E} }
            { \erepn[\bC] }
        & \e{\dynkinEp[1,0,0,0,0,0]{6}{tables-Ud-E1} }
            {}
        & \e{ \dynkinEp[0,0,0,0,-1,1]{6}{tables-Ud-E2} }
            {}
        & \e{ \dynkinEp[0,0,0,0,-1,0]{6}{tables-Ud-E3} }
            {}
      \bpad \\ \hline \tpad
      \type{BD}{n+4}
        & \e{ \dynkinBD[0,0,0,0,0,0,0,0]{2}{0}{2}{tables-Ud-BD} }
            { \bC }
        & \e{ \dynkinBDp[0,0,0,0,0,0,0,0]{2}{0}{2}{tables-Ud-BD1} }
            {}
        & \e{ - }
            {}
        & \e{ - }
            {}
      \bpad[-2.3ex] \\ \hline
    \end{tabular}
    \caption[$\bU^*$ and its graded components]
            {The representation $\bU^* \defeq \Symm{2}\bW^* / \Cartan{2}\bW^*$ defining the adjugate, together with its graded components.}
    \label{tbl:tbl-app-Ud}
  \end{table}
\end{landscape}

\addcontentsline{toc}{chapter}{Bibliography}
\bibliography{bib}

\chapter*{Index of notation} 
\label{c:app-symb}

\BufferDynkinLocaltrue
\renewcommand{\dynkinnameoffset}{-0.75}

\addcontentsline{toc}{chapter}{Index of notation}

\begin{multicols*}{2} \small
\newcommand{\fletter}[1]{\medskip \item[\textbf{#1}]}
\newcommand{\letter}[1]{\bigskip}
\newcommand{\e}[1]{\item[{$#1$}]}

\begin{itemize}[
  itemindent=0em,
  align=left,
  leftmargin=5.2em,
  rightmargin=0.8em,
  labelwidth=4.7em,
  itemsep=-0.3em
]

  \letter{A}
\e{\AD}                             Adjoint action of $G$ on $G$
\e{\Ad}                             Adjoint action of $G$ on $\fg$
\e{\ad}                             Adjoint action of $\fg$ on $\fg$
\e{A(g,\b{g})}                      Solution of main equation

  \letter{B} 
\e{B}                               Inverse metric repn.\ of a PPG

  \letter{C}
\e{\CP}                             Complex \proj\ space
\e{\cpxbdl{E}}                      Complexification of $E$
\e{\CY{}{}}                         Cotton--York tensor of $\D$
\e{\erepn[\bC]}                     $27$-dimensional repn.\ of $\alg[_6]{e}{\bC}$
\e{\eerepn[\bC]}                    $56$-dimensional repn.\ of $\alg[_7]{e}{\bC}$

  \letter{D}
\e{\d}                              Exterior derivative
\e{\d^{\D}}                         Exterior covariant derivative induced by $\D$
\e{\bgg{\bV}}                       First BGG operator of $\bV$
\e{\p}                              Lie homology boundary map
\e{\p^*}                            Lie cohomology differential
\e{\weyld{\gamma}}                  Change in Weyl structure wrt.\ $\gamma$
\e{\div[g]{}}                       Divergence wrt.\ $g$

  \letter{E}
\e{\{e_i\}_i}                       Local frame of $M$
\e{\{\ve^i\}_i}                     Coframe dual to $\{e_i\}_i$
\e{\exp}                            Exponential map

  \letter{F}
\e{F^P}                             Cartan $P$-frame bundle

  \letter{G}
\e{G\acts\fp}                       Flat model of a PPG
\e{\alg{gl}{V}}                     Endomorphisms of $V$
\e{\grad[g]{}}                      Gradient wrt.\ $g$
\e{\Grass{k}{V}}                    \Grassmannian\ of $k$-planes in $V$%

  \letter{H}
\e{\bH}                             \Qtn\ algebra
\e{\HP}                             \Qtn ic\ \proj\ space
\e{H\acts\fq}                       Big \Rspace\ of a PPG
\e{\liehom[\fg]{k}{V}}              $k$th $V$-valued Lie homology of $\fg$
\e{\liecohom[\fg]{k}{V}}            $k$th $V$-valued Lie algebra cohomology of $\fg$

  \letter{I}
\e{\cpx{i},\cpx{i_a}}               Imaginary units in $\bC$ or $\bH$

  \letter{J}
\e{J,J_a}                           Almost complex structures
\e{\jetbdl{k}{E}}                   $k$th jet bundle of $E$


  \letter{L}
\e{L}                               Bundle of scales of a PPG
\e{\lie{X}}                         Lie derivative wrt.\ $X$


  \letter{N}
\e{n}                               \Proj\ dimension of a PPG
\e{\bN}                             Natural numbers $\bZ_{\geq 0}$
\e{\Nijen{}{}}                      Nijenhuis torsion of $J$

  \letter{O}
\e{\bO}                             \Oct\ algebra
\e{\OP}                             \Oct ic\ \proj\ plane

  \letter{P}
\e{\pr{V}}                          Projectivisation of $V$
\e{\pf{}}                           Pfaffian

  \letter{Q}
\e{\cQ}                             Almost \qtn ic\ structure

  \letter{R}
\e{r}                               Scalar dimension of a PPG
\e{\RP}                             Real \proj\ space
\e{\nRic{}{}}                       Normalised Ricci curvature of $\D$
\e{\Curv{}{}}                       Curvature of $\D$
\e{\rk{}}                           Rank
\e{\sRic{}{}}                       Ricci curvature of $\D$

  \letter{S}
\e{\Sph[n]}                         Unit sphere in $\bR^{n+1}$
\e{\spinrepn, \spinrepn^{\pm}}      (Half)-spin repn.\ of $\alg{so}{n,\bC}$
\e{\Scal[g]}                        Scalar curvature of $g$

  \letter{T}
\e{\bT}                             Standard repn.\ of a PPG
\e{\Tor[\D]{}{}}                    Torsion of $\D$
\e{\tr}                             Trace

  \letter{U}
\e{\bU^*}                           Defining quadric of a PPG

  \letter{V}
\e{\realrepn{V}}                    Real repn.\ underlying $V$
\e{\cpxrepn{V}}                     Complexification of $V$
\e{\vol[g]}                         Volume form of $g$

  \letter{W}
\e{\bW}                             Isotropy representation of $\fh/\fq$
\e{\Weyl[\D]{}{}}                   Weyl curvature of $\D$
\e{\omega}                          \Kahler\ form
\e{\Omega^Q}                        (Fundamental) Kraines form
\e{\s{k}{E}}                        Sections of $\Wedge{k}T^*M\tens E$ over $M$


  \letter{Z}
\e{\liecenter{\fg},\grpcenter{G}}   Centre of $\fg,G$

\bigskip\medskip
  \fletter{Miscellaneous symbols}

\e{\tens}                           Tensor product
\e{\etens}                          External tensor product
\e{\wedge}                          Wedge product
\e{\symm}                           Symmetric product
\e{\cartan}                         Cartan product
\e{\liebrac{}{}}                    Lie bracket
\e{\liebracw{}{}}                   Skew-symmetrisation of $\liebrac{}{}$
\e{\killing{}{}}                    Killing form
\e{\algbrac{}{}}                    Algebraic bracket
\e{\algbracw{}{}{}}                 Skew-symmetrisation of $\algbrac{}{}$
\e{\algbraco{}{}{}}                 Symmetrisation of $\algbrac{}{}$
\e{\poisson{}{}}                    Poisson bracket

\e{\D}                              (Weyl) connection
\e{\Dspace}                         Space of Weyl connections
\e{\D^{\bV}}                        Tractor connection induced by $\bV$
\e{\D^{\cV}}                        Prolongation connection on $\cV$
\e{\quab}                           Algebraic laplacian
\e{\quab_M}                         First-order laplacian
\e{\intprod}                        Contraction

\e{\bdot\at{x}}                     Evaluation at a point $x$
\e{\Injto}                          Injective map
\e{\Surjto}                         Surjective map
\e{\flat,\sharp}                    Musical isomorphisms

\bigskip\medskip
  \fletter{Generic notations}

\e{M}                               Smooth manifold
\e{m}                               Typical dimension of a manifold
\e{X,Y,Z}                           Vector fields
\e{\alpha,\beta,\gamma}             Differential forms
\e{g}                               (\Pseudo)\riem\ metric

\e{G,H}                             Semisimple Lie groups
\e{\fg,\fh}                         Semisimple Lie algebras
\e{P,Q}                             Parabolic subgroups
\e{\fp,\fq}                         Parabolic subalgebras
\e{\fp^{\perp}}                     Killing polar of $\fp$
\e{\fp^0}                           Reductive Levi factor $\fp/\fp^{\perp}$ of $\fp$
\e{\opp{\fp}}                       Parabolic opposite to $\fp$

\e{V,W}                             Lie algebra repns.\
\e{\bV,\bW}                         Lie algebra repns., restricted to a parabolic subalgebra
\e{\cV,\cW}                         Bundles associated to either $V,W$ or $\bV,\bW$

\e{h}                               Element of $L^*\tens B$
\e{Z}                               Element of $L^*\tens\fg/\fp$
\e{\lambda}                         Element of $L^*$
\e{\ell,\pi}                        Elements of $L$
\e{\zeta,\eta}                      Elements of $L\tens\fp^{\perp}$
\e{\theta}                          Element of $L\tens B^*$

\end{itemize}
\end{multicols*}

\end{document}